\documentclass[9pt]{amsart}
\usepackage{amsmath}
\usepackage{amsxtra}
\usepackage{amscd}
\usepackage{amsthm}
\usepackage{amsfonts}
\usepackage{amssymb}
\usepackage{eucal}
\usepackage[all]{xy}
\usepackage{graphicx}
\usepackage[usenames]{color}
\usepackage[hidelinks]{hyperref}

\usepackage{tikz-cd}

\newtheorem{cor}[subsubsection]{Corollary}

\newtheorem{lem}[subsubsection]{Lemma}
\newtheorem{prop}[subsubsection]{Proposition}

\newtheorem{conj}[subsubsection]{Conjecture}

\newtheorem{thm}[subsubsection]{Theorem}
\newtheorem{mainthm}[subsubsection]{Main Theorem}

\theoremstyle{remark}
\newtheorem{rem}[subsubsection]{Remark}

\theoremstyle{definition}

\theoremstyle{remark}

\newcommand{\thmref}[1]{Theorem~\ref{#1}}

\newcommand{\secref}[1]{Sect.~\ref{#1}}
\newcommand{\lemref}[1]{Lemma~\ref{#1}}
\newcommand{\propref}[1]{Proposition~\ref{#1}}
\newcommand{\corref}[1]{Corollary~\ref{#1}}
\newcommand{\conjref}[1]{Conjecture~\ref{#1}}

\numberwithin{equation}{section}

\newcommand{\nc}{\newcommand}
\nc{\renc}{\renewcommand}
\nc{\ssec}{\subsection}
\nc{\sssec}{\subsubsection}
\nc{\on}{\operatorname}

\nc{\ips}{{\iota_P^{(S)}}}
\nc{\ipms}{{\iota_{P^-}^{(S)}}}
\nc{\sfpps}{{\sfp_P^{(S)}}}
\nc{\sfppms}{{\sfp_{P^-}^{(S)}}}

\nc\ol{\overline}
\nc\ul{\underline}
\nc\wt{\widetilde}
\nc\tboxtimes{\wt{\boxtimes}}
\nc\tstar{\wt{\star}}
\nc{\alp}{\alpha}

\nc{\ZZ}{{\mathbb Z}}
\nc{\NN}{{\mathbb N}}
\nc{\OO}{{\mathbb O}}
\renc{\SS}{{\mathbb S}}
\nc{\DD}{{\mathbb D}}
\nc{\GG}{{\mathbb G}}

\nc{\Fq}{{\mathbb F}_q}
\nc{\Fqb}{\ol{\mathbb F}_q}
\nc{\Ql}{{\mathbb Q}_\ell}
\nc{\Qlb}{{\ol{\mathbb Q}_\ell}}
\nc{\id}{\text{id}}
\nc\X{\mathcal X}

\nc{\red}{\on{red}}
\nc{\Ho}{\on{Ho}}
\nc{\Hom}{\on{Hom}}
\nc{\coHom}{\ul{\on{coHom}}}
\nc{\coMaps}{{\bf{coMaps}}}
\nc{\coef}{\on{coef}}
\nc{\Lie}{\on{Lie}}
\nc{\Loc}{\on{Loc}}
\nc{\coLoc}{\on{coLoc}}
\nc{\Pic}{\on{Pic}}
\nc{\Bun}{\on{Bun}}
\nc{\IC}{\on{IC}}
\nc{\Aut}{\on{Aut}}
\nc{\rk}{\on{rk}}
\nc{\Sh}{\on{Sh}}
\nc{\Perv}{\on{Perv}}
\nc{\pos}{{\on{pos}}}
\nc{\Conv}{\on{Conv}}
\nc{\Sph}{\on{Sph}}
\nc{\Sym}{\on{Sym}}
\nc{\BunBb}{\overline{\Bun}_B}
\nc{\BunNb}{\overline{\Bun}_{N,\rho(\omega_X)}}
\nc{\BunTb}{\overline{\Bun}_T}
\nc{\BunBbm}{\overline{\Bun}_{B^-}}
\nc{\BunBbel}{\overline{\Bun}_{B,el}}
\nc{\BunBbmel}{\overline{\Bun}_{B^-,el}}
\nc{\Buno}{\overset{o}{\Bun}}
\nc{\BunPb}{{\overline{\Bun}_P}}
\nc{\BunBM}{\Bun_{B(M)}}
\nc{\BunBMb}{\overline{\Bun}_{B(M)}}
\nc{\BunPbw}{{\widetilde{\Bun}_P}}
\nc{\BunBP}{\widetilde{\Bun}_{B,P}}
\nc{\GUb}{\overline{G/U}}
\nc{\GUPb}{\overline{G/U(P)}}

\nc{\Hhom}{\underline{\on{Hom}}}
\nc\syminfty{\on{Sym}^{\infty}}
\nc\lal{\ol{\lambda}}
\nc\xl{\ol{x}}
\nc\thl{\ol{\theta}}
\nc\nul{\ol{\nu}}
\nc\mul{\ol{\mu}}
\nc\Sum\Sigma
\nc{\oX}{\overset{o}{X}{}}
\nc{\hl}{\overset{\leftarrow}h{}}
\nc{\hr}{\overset{\rightarrow}h{}}
\nc{\M}{{\mathcal M}}
\nc{\N}{{\mathcal N}}
\nc{\F}{{\mathcal F}}
\nc{\D}{{\mathcal D}}
\nc{\Q}{{\mathcal Q}}
\nc{\Y}{{\mathcal Y}}
\nc{\G}{{\mathcal G}}
\nc{\E}{{\mathcal E}}
\nc{\CalC}{{\mathcal C}}
\nc\Dh{\widehat{\D}}

\nc{\C}{{\mathcal C}}
\nc{\K}{{\mathcal K}}
\renewcommand{\H}{{\mathcal H}}

\nc{\T}{{\mathcal T}}
\nc{\V}{{\mathcal V}}
\renc{\P}{{\mathcal P}}
\nc{\A}{{\mathcal A}}
\nc{\B}{{\mathcal B}}
\nc{\U}{{\mathcal U}}

\nc{\Gr}{{\on{Gr}}}

\nc{\frn}{{\check{\mathfrak u}(P)}}

\nc{\fC}{\mathfrak C}
\nc{\fT}{\mathfrak T}
\nc{\p}{\mathfrak p}
\nc{\q}{\mathfrak q}
\nc\f{{\mathfrak f}}

\nc{\qo}{{\mathfrak q}}
\nc{\po}{{\mathfrak p}}
\nc{\s}{{\mathfrak s}}
\nc\w{\text{w}}

\renewcommand{\mod}{{\on{-mod}}}

\nc\mathi\iota
\nc\Spec{\on{Spec}}
\nc\Proj{\on{Proj}}
\nc\Mod{\on{Mod}}
\nc{\tw}{\widetilde{\mathfrak t}}
\nc{\pw}{\widetilde{\mathfrak p}}
\nc{\qw}{\widetilde{\mathfrak q}}
\nc{\jw}{\widetilde j}

\nc{\grb}{\overline{\Gr}}
\nc{\I}{\mathcal I}

\nc{\lambdach}{{\check\lambda}}
\nc{\Lambdach}{{\check\Lambda}{}}
\nc{\much}{{\check\mu}}
\nc{\omegach}{{\check\omega}}
\nc{\nuch}{{\check\nu}}
\nc{\etach}{{\check\eta}}
\nc{\alphach}{{\check\alpha}}
\nc{\oblvtach}{{\check\oblvta}}
\nc{\rhoch}{{\check\rho}}
\nc{\ch}{{\check h}}

\nc{\Hb}{\overline{\H}}

\nc{\BA}{{\mathbb{A}}}
\nc{\BC}{{\mathbb{C}}}
\nc{\BE}{{\mathbb{E}}}
\nc{\BF}{{\mathbb{F}}}
\nc{\BG}{{\mathbb{G}}}
\nc{\BL}{{\mathbb{L}}}
\nc{\BM}{{\mathbb{M}}}
\nc{\BO}{{\mathbb{O}}}
\nc{\BD}{{\mathbb{D}}}
\nc{\BN}{{\mathbb{N}}}
\nc{\BP}{{\mathbb{P}}}
\nc{\BQ}{{\mathbb{Q}}}
\nc{\BR}{{\mathbb{R}}}
\nc{\BV}{{\mathbb{V}}}
\nc{\BW}{{\mathbb{W}}}
\nc{\BZ}{{\mathbb{Z}}}
\nc{\BS}{{\mathbb{S}}}
\nc{\Deep}{{\bf{deep}}}
\nc{\deep}{deep}

\nc{\CA}{{\mathcal{A}}}
\nc{\CB}{{\mathcal{B}}}

\nc{\CE}{{\mathcal{E}}}
\nc{\CF}{{\mathcal{F}}}
\nc{\CH}{{\mathcal{H}}}

\nc{\CL}{{\mathcal{L}}}
\nc{\CC}{{\mathcal{C}}}
\nc{\CG}{{\mathcal{G}}}
\nc{\CalD}{{\mathcal{D}}}
\nc{\CM}{{\mathcal{M}}}
\nc{\CN}{{\mathcal{N}}}
\nc{\CK}{{\mathcal{K}}}
\nc{\CO}{{\mathcal{O}}}
\nc{\CP}{{\mathcal{P}}}
\nc{\CQ}{{\mathcal{Q}}}
\nc{\CR}{{\mathcal{R}}}
\nc{\CS}{{\mathcal{S}}}
\nc{\CT}{{\mathcal{T}}}
\nc{\CU}{{\mathcal{U}}}
\nc{\CV}{{\mathcal{V}}}
\nc{\CW}{{\mathcal{W}}}
\nc{\CX}{{\mathcal{X}}}
\nc{\CY}{{\mathcal{Y}}}
\nc{\CZ}{{\mathcal{Z}}}
\nc{\CI}{{\mathcal{I}}}

\nc{\csM}{{\check{\mathcal A}}{}}
\nc{\oM}{{\overset{\circ}{\mathcal M}}{}}
\nc{\obM}{{\overset{\circ}{\mathbf M}}{}}
\nc{\oCA}{{\overset{\circ}{\mathcal A}}{}}
\nc{\obA}{{\overset{\circ}{\mathbf A}}{}}
\nc{\ooM}{{\overset{\circ}{M}}{}}
\nc{\osM}{{\overset{\circ}{\mathsf M}}{}}
\nc{\vM}{{\overset{\bullet}{\mathcal M}}{}}
\nc{\nM}{{\underset{\bullet}{\mathcal M}}{}}
\nc{\oD}{{\overset{\circ}{\mathcal D}}{}}
\nc{\obD}{{\overset{\circ}{\mathbf D}}{}}
\nc{\oA}{{\overset{\circ}{A}}{}}
\nc{\op}{{\overset{\bullet}{\mathbf p}}{}}
\nc{\oU}{{\overset{\bullet}{\mathcal U}}{}}
\nc{\oZ}{{\overset{\circ}{\mathcal Z}}{}}
\nc{\ofZ}{{\overset{\circ}{\mathfrak Z}}{}}
\nc{\oF}{{\overset{\circ}{\fF}}}

\nc{\fa}{{\mathfrak{a}}}
\nc{\ofa}{\overset{\circ}{\mathfrak{a}}}
\nc{\fb}{{\mathfrak{b}}}
\nc{\fd}{{\mathfrak{d}}}
\nc{\ff}{{\mathfrak{f}}}
\nc{\fg}{{\mathfrak{g}}}
\nc{\fgl}{{\mathfrak{gl}}}
\nc{\fh}{{\mathfrak{h}}}
\nc{\fj}{{\mathfrak{j}}}
\nc{\fl}{{\mathfrak{l}}}
\nc{\fm}{{\mathfrak{m}}}
\nc{\ofm}{\overset{\circ}{\mathfrak{m}}}
\nc{\fn}{{\mathfrak{n}}}
\nc{\fu}{{\mathfrak{u}}}
\nc{\fp}{{\mathfrak{p}}}
\nc{\fr}{{\mathfrak{r}}}
\nc{\fs}{{\mathfrak{s}}}
\nc{\ft}{{\mathfrak{t}}}
\nc{\oft}{\overset{\circ}{\mathfrak{t}}}
\nc{\fz}{{\mathfrak{z}}}
\nc{\fsl}{{\mathfrak{sl}}}
\nc{\hsl}{{\widehat{\mathfrak{sl}}}}
\nc{\hgl}{{\widehat{\mathfrak{gl}}}}
\nc{\hg}{{\widehat{\mathfrak{g}}}}
\nc{\hp}{{\widehat{\mathfrak{p}}}}
\nc{\hm}{{\widehat{\mathfrak{m}}}}
\nc{\chg}{{\widehat{\mathfrak{g}}}{}^\vee}
\nc{\hn}{{\widehat{\mathfrak{n}}}}
\nc{\chn}{{\widehat{\mathfrak{n}}}{}^\vee}

\nc{\fA}{{\mathfrak{A}}}
\nc{\fB}{{\mathfrak{B}}}
\nc{\fD}{{\mathfrak{D}}}
\nc{\fE}{{\mathfrak{E}}}
\nc{\fF}{{\mathfrak{F}}}
\nc{\fG}{{\mathfrak{G}}}
\nc{\fK}{{\mathfrak{K}}}
\nc{\fL}{{\mathfrak{L}}}
\nc{\fM}{{\mathfrak{M}}}
\nc{\fN}{{\mathfrak{N}}}
\nc{\fP}{{\mathfrak{P}}}
\nc{\fU}{{\mathfrak{U}}}
\nc{\fV}{{\mathfrak{V}}}
\nc{\fZ}{{\mathfrak{Z}}}

\nc{\ba}{{\mathbf{a}}}
\nc{\bb}{{\mathbf{b}}}
\nc{\bc}{{\mathbf{c}}}
\nc{\bd}{{\mathbf{d}}}
\nc{\bbf}{{\mathbf{f}}}
\nc{\be}{{\mathbf{e}}}
\nc{\bi}{{\mathbf{i}}}
\nc{\bj}{{\mathbf{j}}}
\nc{\bh}{{\mathbf{h}}}
\nc{\bm}{{\mathbf{m}}}
\nc{\bn}{{\mathbf{n}}}
\nc{\bo}{{\mathbf{o}}}
\nc{\bp}{{\mathbf{p}}}
\nc{\bq}{{\mathbf{q}}}
\nc{\bu}{{\mathbf{u}}}
\nc{\bv}{{\mathbf{v}}}
\nc{\bx}{{\mathbf{x}}}
\nc{\bs}{{\mathbf{s}}}
\nc{\by}{{\mathbf{y}}}
\nc{\bw}{{\mathbf{w}}}
\nc{\bA}{{\mathbf{A}}}
\nc{\bK}{{\mathbf{K}}}
\nc{\bB}{{\mathbf{B}}}
\nc{\bC}{{\mathbf{C}}}
\nc{\bG}{{\mathbf{G}}}
\nc{\bD}{{\mathbf{D}}}
\nc{\bE}{{\mathbf{E}}}
\nc{\bH}{{{\mathbf{H}}}}
\nc{\bL}{{\mathbf{L}}}
\nc{\bM}{{\mathbf{M}}}
\nc{\bN}{{\mathbf{N}}}
\nc{\bO}{{\mathbf{O}}}
\nc{\bQ}{{\mathbf{Q}}}
\nc{\bV}{{\mathbf{V}}}
\nc{\bW}{{\mathbf{W}}}
\nc{\bX}{{\mathbf{X}}}
\nc{\bZ}{{\mathbf{Z}}}
\nc{\bS}{{\mathbf{S}}}

\nc{\sA}{{\mathsf{A}}}
\nc{\sB}{{\mathsf{B}}}
\nc{\sC}{{\mathsf{C}}}
\nc{\sD}{{\mathsf{D}}}
\nc{\sF}{{\mathsf{F}}}
\nc{\sG}{{\mathsf{G}}}
\nc{\sH}{{\mathsf{H}}}
\nc{\sK}{{\mathsf{K}}}
\nc{\sM}{{\mathsf{M}}}
\nc{\sN}{{\mathsf{N}}}
\nc{\sO}{{\mathsf{O}}}
\nc{\sW}{{\mathsf{W}}}
\nc{\sQ}{{\mathsf{Q}}}
\nc{\sP}{{\mathsf{P}}}
\nc{\sR}{{\mathsf{R}}}
\nc{\sT}{{\mathsf{T}}}
\nc{\sZ}{{\mathsf{Z}}}
\nc{\sfi}{{\mathsf{i}}}
\nc{\sfj}{{\mathsf{j}}}
\nc{\sfp}{{\mathsf{p}}}
\nc{\sfq}{{\mathsf{q}}}
\nc{\sfs}{{\mathsf{s}}}
\nc{\sft}{{\mathsf{t}}}
\nc{\sr}{{\mathsf{r}}}
\nc{\sfk}{{\mathsf{k}}}
\nc{\sa}{{\mathsf{s}}}
\nc{\sg}{{\mathsf{g}}}
\nc{\sn}{{\mathsf{n}}}
\nc{\sh}{{\mathsf{h}}}
\nc{\sff}{{\mathsf{f}}}
\nc{\sfb}{{\mathsf{b}}}
\nc{\sfc}{{\mathsf{c}}}
\nc{\sfe}{{\mathsf{e}}}
\nc{\sd}{{\mathsf{d}}}

\nc{\BK}{{\bar{K}}}

\nc{\tA}{{\widetilde{\mathbf{A}}}}
\nc{\tB}{{\widetilde{\mathcal{B}}}}
\nc{\tg}{{\widetilde{\mathfrak{g}}}}
\nc{\tG}{{\widetilde{G}}}
\nc{\TM}{{\widetilde{\mathbb{M}}}{}}
\nc{\tO}{{\widetilde{\mathsf{O}}}{}}
\nc{\tU}{{\widetilde{\mathfrak{U}}}{}}
\nc{\TZ}{{\tilde{Z}}}
\nc{\tx}{{\tilde{x}}}
\nc{\tbv}{{\tilde{\bv}}}
\nc{\tfP}{{\widetilde{\mathfrak{P}}}{}}
\nc{\tz}{{\tilde{\zeta}}}
\nc{\tmu}{{\tilde{\mu}}}

\nc{\urho}{\underline{\rho}}
\nc{\uB}{\underline{B}}
\nc{\uC}{{\underline{\mathbb{C}}}}
\nc{\ui}{\underline{i}}
\nc{\uj}{\underline{j}}
\nc{\ofP}{{\overline{\mathfrak{P}}}}
\nc{\oB}{{\overline{\mathcal{B}}}}
\nc{\og}{{\overline{\mathfrak{g}}}}
\nc{\oI}{{\overline{I}}}

\nc{\eps}{\varepsilon}
\nc{\hrho}{{\hat{\rho}}}

\nc{\one}{{\mathbf{1}}}
\nc{\two}{{\mathbf{t}}}

\nc{\Rep}{{\mathop{\operatorname{\rm Rep}}}}
\nc{\Tot}{{\mathop{\operatorname{\rm Tot}}}}
\nc{\Ker}{{\mathop{\operatorname{\rm Ker}}}}
\nc{\im}{{\mathop{\operatorname{\rm Im}}}}
\nc{\Hilb}{{\mathop{\operatorname{\rm Hilb}}}}
\nc{\End}{{\mathop{\operatorname{\rm End}}}}
\nc{\Ext}{{\mathop{\operatorname{\rm Ext}}}}
\nc{\CHom}{{\mathop{\operatorname{{\mathcal{H}}\it om}}}}
\nc{\CEnd}{{\mathop{\operatorname{{\mathcal{E}}\it nd}}}}
\nc{\GL}{{\mathop{\operatorname{\rm GL}}}}
\nc{\gr}{{\mathop{\operatorname{\rm gr}}}}
\nc{\HN}{{\mathop{\operatorname{\rm HN}}}}
\nc{\Id}{{\mathop{\operatorname{\rm Id}}}}
\nc{\de}{{\mathop{\operatorname{\rm def}}}}
\nc{\length}{{\mathop{\operatorname{\rm length}}}}
\nc{\supp}{{\mathop{\operatorname{\rm supp}}}}

\nc{\Cliff}{{\mathsf{Cliff}}}
\nc{\Fl}{\on{Fl}}
\nc{\Fib}{{\mathsf{Fib}}}
\nc{\Coh}{{\on{Coh}}}
\nc{\QCoh}{{\on{QCoh}}}
\nc{\IndCoh}{{\on{IndCoh}}}
\nc{\FCoh}{{\mathsf{FCoh}}}

\nc{\reg}{{\text{\rm reg}}}

\nc{\cplus}{{\mathbf{C}_+}}
\nc{\cminus}{{\mathbf{C}_-}}
\nc{\cthree}{{\mathbf{C}_\bullet}}
\nc{\Qbar}{{\bar{Q}}}
\nc\Eis{\on{Eis}}
\nc\Eisb{\ol\Eis{}}
\nc\Eisr{\on{Eis}^{rat}{}}
\nc\wh{\widehat}
\nc{\Def}{\on{Def_{\check{\fb}}(E)}}
\nc{\barZ}{\overline{Z}{}}
\nc{\barbarZ}{\overline{\barZ}{}}
\nc{\barpi}{\overline\pi}
\nc{\barbarpi}{\overline\barpi}
\nc{\barpip}{\overline\pi{}^+}
\nc{\barpim}{\overline\pi{}^-}

\nc{\fq}{\mathfrak q}

\nc{\fqb}{\ol{\sfq}{}}
\nc{\fpb}{\ol{\sfp}{}}
\nc{\fpr}{{\mathsf{pair}^{rat}}{}}
\nc{\fqr}{{\sfq^{rat}}{}}

\nc{\hattimes}{\wh\otimes}

\nc{\bOmega}{{\overline{\Omega(\check \fn)}}}

\nc{\seq}[1]{\stackrel{#1}{\sim}}

\nc{\cT}{{\check{T}}}
\nc{\cG}{{\check{G}}}
\nc{\cM}{{\check{M}}}
\nc{\cB}{{\check{B}}}
\nc{\cP}{{\check{P}}}
\nc{\cN}{{\check{N}}}

\nc{\ct}{{\check{\mathfrak t}}}
\nc{\cg}{{\check{\fg}}}
\nc{\cb}{{\check{\fb}}}
\nc{\cp}{{\check{\fp}}}
\nc{\cn}{{\check{\fn}}}
\nc{\cm}{{\check{\fm}}}

\nc{\cLambda}{{\check\Lambda}}

\nc{\cla}{{\check\lambda}}
\nc{\cmu}{{\check\mu}}
\nc{\cnu}{{\check\nu}}
\nc{\ceta}{{\check\eta}}

\nc{\DefbE}{{\on{Def}_{\cB}(E_\cT)}}

\nc{\imathb}{{\ol{\imath}}}
\nc{\rlr}{\overset{\longrightarrow}{\underset{\longrightarrow}\longleftarrow}}

\nc{\oBun}{\overset{\circ}\Bun}
\nc{\LS}{\on{LS}}
\nc{\BunBbb}{\ol{\ol{Bun}}_B}
\nc{\BunBr}{\Bun_B^{rat}}
\nc{\BunBrsg}{\Bun_B^{rat,\on{s.g.}}}
\nc{\BunBrp}{\Bun_B^{rat,polar}}
\nc{\BunBrpbg}{\Bun_B^{rat,polar,\on{b.g.}}}
\nc{\BunBrpsg}{\Bun_B^{rat,polar,\on{s.g.}}}
\nc{\BunTrp}{\Bun_T^{rat,polar}}
\nc{\BunTrpbg}{\Bun_T^{rat,polar,\on{b.g.}}}
\nc{\BunTrpsg}{\Bun_T^{rat,polar,\on{s.g.}}}
\nc{\BunNr}{\Bun_N^{rat}}
\nc{\BunNre}{\Bun_N^{enh,rat}}
\nc{\BunTr}{\Bun_T^{rat}}
\nc{\Vect}{\on{Vect}}
\nc{\Whit}{\on{Whit}}
\nc{\CTb}{\ol{\on{CT}}}
\nc{\Ran}{{\on{Ran}}}
\nc{\CTr}{\on{CT}^{rat}{}}
\nc\jmathr{\jmath^{rat}{}}
\nc{\ux}{\underline{x}}
\nc{\clambda}{{\check\lambda}}
\nc{\calpha}{{\check\alpha}}
\nc{\ind}{{\mathbf{ind}}}
\nc{\coinv}{{\mathbf{coinv}}}
\nc{\oblv}{{\mathbf{oblv}}}
\nc{\free}{{\mathbf{free}}}
\nc{\ox}{{\overline{x}}}
\nc{\cLa}{\check{\Lambda}}
\nc{\StinftyCat}{\on{DGCat}}
\nc{\inftyCat}{\infty\on{-Cat}}
\nc{\inftygroup}{\infty\on{-Grpd}}
\nc{\Dmod}{\on{D-mod}}
\nc{\CMaps}{{\mathcal Maps}}
\nc{\Maps}{\on{Maps}}
\nc{\affSch}{\on{Sch}^{\on{aff}}}
\nc{\dr}{{\on{dR}}}
\nc{\oCF}{\overset{\circ}\CF}
\nc{\oCY}{\overset{\circ}\CY}
\nc{\opi}{\overset{\circ}\pi}
\nc{\leqG}{\underset{G}\leq}
\nc{\leqM}{\underset{M}\leq}
\nc{\leqGad}{\underset{G_{ad}}\leq}
\nc{\leqMad}{\underset{M_{ad}}\leq}
\nc{\Tr}{\on{Tr}}
\nc{\Frob}{{\on{Frob}}}
\nc{\DGCat}{\on{DGCat}}
\nc{\tDGCat}{2\on{-DGCat}_{\on{u.g.}}}
\nc{\ev}{\on{ev}}
\nc{\mmod}{\on{-}\mathbf{mod}}
\nc{\sotimes}{\overset{!}\otimes}
\nc{\Shv}{\on{Shv}}
\nc{\Spc}{\on{Spc}}
\nc{\Res}{\on{Res}}
\nc{\bDelta}{{\mathbf{\Delta}}}
\nc{\bMaps}{{\mathbf{Maps}}}
\nc{\cD}{\mathcal D}
\nc{\ocD}{\overset{\circ}\cD}
\nc{\ppart}{(\!(t)\!)}
\nc{\qqart}{[\![t]\!]}
\nc{\oCU}{\overset{\circ}{\CU}}
\nc{\Exc}{{\mathcal{E}xc}}
\nc{\Sht}{\on{Sht}}
\nc{\Nilp}{{\on{Nilp}}}
\nc{\Drinf}{\on{Drinf}}
\nc{\Sing}{\on{Sing}}
\nc{\IndLisse}{\Lisse}
\nc{\Shvl}{\on{Shv}_{\on{lisse}}} 
\nc{\Lisse}{\on{Lisse}}
\nc{\Mir}{\on{Mir}}
\nc{\fSet}{\on{fSet}}
\nc{\qLisse}{\on{QLisse}}
\nc{\Ev}{\on{Ev}}
\nc{\Sat}{\on{Sat}}
\nc{\Se}{\on{Se}}
\nc{\coSht}{\on{co-Sht}}
\nc{\coCK}{\on{co-}\!\CK}
\nc{\FLE}{\on{FLE}}
\nc{\BRST}{\on{BRST}}
\nc{\KL}{\on{KL}}
\nc{\crit}{{\on{crit}}}
\nc{\Op}{{\on{Op}}}
\nc{\MOp}{\on{MOp}}
\nc{\Wak}{\on{Wak}}
\nc{\Av}{\on{Av}}
\nc{\semiinf}{{\frac{\infty}{2}}}
\nc{\DS}{\on{DS}}
\nc{\dR}{{\on{dR}}}
\nc{\Poinc}{{\on{Poinc}}}
\renc{\det}{\on{det}}
\nc{\oG}{\overset{\circ}{G}}
\nc{\mer}{{\on{mer}}}
\nc{\mf}{{\on{mon-free}}}
\nc{\Mmf}{{\on{mon-free}_\cM}}
\nc{\Gmf}{{\on{mon-free}_\cG}}
\nc{\GMmf}{{\on{mon-free}_{\cG,\cM}}}
\nc{\Pmf}{{\on{mon-free}_{\cP^-}}}
\nc{\Accs}{{\mathsf{Accs}}} 
\nc{\semiinfAccs}{\semiinf\!\on{-}\!{\mathsf{Accs}}} 
\nc{\BunPmt}{{\widetilde{\Bun}_{P^-}}}
\nc{\BunPmtpR}{{\widetilde{\Bun}_{P^-,\Ran}}}
\nc{\BunPmtpZ}{{\widetilde{\Bun}_{P^-,\CZ}}}
\nc{\BunPmtpZp}{{\widetilde{\Bun}_{P^-,\CZ;}}}
\nc{\BunPmtpZsub}{{\widetilde{\Bun}_{P^-,\CZ^{\subseteq}}}}
\nc{\BunPmtpS}{{\widetilde{\Bun}_{P^-,S}}}
\nc{\BunPmtpx}{{\widetilde{\Bun}_{P^-,\ul{x}}}}
\nc{\BunNtpx}{{\widetilde{\Bun}_{N_P^-,\ul{x}}}}
\nc{\BunNbR}{{\ol{\Bun}_{N,\rho(\omega_X),\Ran}}}
\nc{\BunNbg}{{\ol{\Bun}_{N,\rho(\omega_X),\on{gen}}}}
\nc{\BunNbx}{{\ol{\Bun}_{N,\rho(\omega_X),\ul{x}}}}
\nc{\BunNbZ}{{\ol{\Bun}_{N,\rho(\omega_X),\CZ}}}
\nc{\BunNbZsub}{{\ol{\Bun}_{N,\rho(\omega_X),\CZ^\subseteq}}}
\nc{\BunNbS}{{\ol{\Bun}_{N,\rho(\omega_X),S}}}
\nc{\BunNbMg}{{\ol{\Bun}_{N(M),\rho(\omega_X),\on{gen}}}}
\nc{\BunNbMZ}{{\ol{\Bun}_{N(M),\rho(\omega_X),\CZ}}}
\nc{\BunNbMZsub}{{\ol{\Bun}_{N(M),\rho(\omega_X),\CZ^\subseteq}}}
\nc{\BunNbMR}{{\ol{\Bun}_{N(M),\rho(\omega_X),\Ran}}}
\nc{\BunNbMgM}{{\ol{\Bun}_{N(M),\rho_M(\omega_X),\on{gen}}}}
\nc{\BunNbMZM}{{\ol{\Bun}_{N(M),\rho_M(\omega_X),\CZ}}}
\nc{\BunNbMZsubM}{{\ol{\Bun}_{N(M),\rho_M(\omega_X),\CZ^\subseteq}}}
\nc{\BunNbMRM}{{\ol{\Bun}_{N(M),\rho_M(\omega_X),\Ran}}}

\nc{\Zas}{{\mathsf{Zas}}}

\begin{document}


\vskip1cm

\title[Proof of the geometric Langlands conjecture III]{Proof of the geometric Langlands conjecture III: \\
compatibility with parabolic induction}

\author{Justin Campbell, Lin Chen, Dennis Gaitsgory and Sam Raskin}\footnote{With an appendix coauthored additionally by J.~Faergeman, K.~Lin and N.~Rozenblyum.}

\dedicatory{To B.~Feigin}

%
\begin{abstract}
We establish the compatibility of the Langlands functor $\BL_G$ with the operations of Eisenstein series constant term,
and deduce that $\BL_G$ induces an equivalence on Eisenstein-generated subcategories. 
\end{abstract} 

\date{\today}

\maketitle

\bigskip

\bigskip


\tableofcontents

\section*{Introduction}

\ssec{What is done in this paper?}

This paper is the third in the series of five, in the course of which a proof of the geometric Langlands
conjecture (stated as \cite[Conjecture 1.6.7]{GLC1}) will be given. 

\sssec{}

As far as the program of proving the geometric Langlands conjecture is concerned, in this paper 
the following two steps toward the proof are performed: 

\begin{itemize}

\item It is shown that the geometric Langlands functor 
\begin{equation} \label{e:Langlands functor intro}
\BL_G:\Dmod_{\frac{1}{2}}(\Bun_G)\to \IndCoh_\Nilp(\LS_\cG)
\end{equation}
is compatible with the functors of \emph{Eisenstein series} and \emph{constant term},
see Theorems \ref{t:L and Eis} and \ref{t:CT compat}, respectively. 

\medskip

\item Assuming the geometric Langlands conjecture for proper Levi subgroups, 
it is shown (\thmref{t:main}) that $\BL_G$ induces an equivalence on Eisenstein-generated subcategories
$$\Dmod_{\frac{1}{2}}(\Bun_G)_{\Eis}\overset{\sim}\to \IndCoh_\Nilp(\LS_\cG)_{\on{red}}.$$

\end{itemize}

\medskip

Another result that concerns the global geometric Langlands program, proved here, and which is of independent
interest, is \thmref{t:left adjoint as dual}, which says that:

\medskip

\begin{itemize}

\item The left adjoint functor of $\BL_G$ can be obtained from the functor \emph{dual} to $\BL_G$, by composing
with the \emph{Miraculous functor} and Chevalley involution.

\end{itemize}

\sssec{}

Other results established in this paper concern the \emph{local Langlands theory}. Apart from being of independent interest,
these results provide local ingredients for the proofs of global theorems mentioned above.

\medskip

The main local result is \thmref{t:local Jacquet}, which establishes the compatibility of the critical FLE with 
the Jacquet functors on the geometric and spectral sides. 

%
%
%
%
%
%

%
%
%
%
%

\ssec{The logical structure: compatibility with Eisenstein series and constant terms}

We will now describe the logical structure of the paper from the point of view of the geometric Langlands conjecture.

\sssec{}

The Langlands functor $\BL_G$ is constructed so that it makes the following diagram commute\footnote{Up to a cohomlogical shift, which we omit 
in the Introduction.} 
\begin{equation} \label{e:coeff comp L Intro}
\CD
\Whit^!(\Gr_{G,\Ran}) @>{\on{CS}_G}>> \Rep(\cG)_\Ran \\
@A{\on{coeff}_G}AA @AA{\Gamma^{\IndCoh,\on{spec}}_\cG}A \\
\Dmod_{\frac{1}{2}}(\Bun_G) @>{\BL_G}>> \IndCoh_\Nilp(\LS_\cG). 
\endCD
\end{equation} 

Here:

\begin{itemize}

\item $\on{CS}_G$ is the geometric Casselman-Shalika equivalence;

\item $\on{coeff}_G$ is the functor of Whittaker coefficient(s);

\medskip

\item $\Gamma^{\IndCoh,\on{spec}}_\cG$ is the functor right adjoint to the localization functor 
$$\Rep(\cG)_\Ran\overset{\Loc_\cG^{\on{spec}}}\longrightarrow \QCoh(\LS_\cG)\hookrightarrow \IndCoh_\Nilp(\LS_\cG).$$

\end{itemize} 

\medskip

The geometric Langlands conjecture (\cite[Conjecture 1.6.7]{GLC2}) says that the functor $\BL_G$ is an equivalence. 

\sssec{}

The compatibility of the geometric Casselman-Shalika equivalence with Jacquet functors implies that $\BL_G$ 
is compatible with the \emph{Eisenstein functors},
i.e., for a standard (negative) parabolic $P^-$ with Levi quotient $M$, the following diagram commutes (again, up to a cohomological shift):

\begin{equation} \label{e:L and Eis Intro}
\CD
\Dmod_{\frac{1}{2}}(\Bun_M) @>{\BL_M}>> \IndCoh_\Nilp(\LS_\cM) \\
@V{\Eis^-_{!,\rho_P(\omega_X)}}VV @VV{\Eis^{-,\on{spec}}}V \\
\Dmod_{\frac{1}{2}}(\Bun_G) @>{\BL_G}>> \IndCoh_\Nilp(\LS_\cG). 
\endCD
\end{equation} 

In the above formula, $\Eis^-_{!,\rho_P(\omega_X)}$ is the \emph{translated} Eisenstein series functor,
see \secref{ss:rho shift Eis}. 

\sssec{}

Given diagram \eqref{e:L and Eis Intro}, by passing to right adjoint functors along the vertical arrows, we obtain a diagram 
\begin{equation} \label{e:top lid nat trans Intro}
\xy
(0,-20)*+{\Dmod_{\frac{1}{2}}(\Bun_G)}="X";
(0,0)*+{\Dmod_{\frac{1}{2}}(\Bun_M)}="Y";
(50,-20)*+{ \IndCoh_\Nilp(\LS_\cG).}="X'";
(50,0)*+{\IndCoh_\Nilp(\LS_\cM),}="Y'";
{\ar@{->}^{\on{CT}^-_{*,\rho_P(\omega_X)}} "X";"Y"};
{\ar@{->}_{\on{CT}^{-,\on{spec}}} "X'";"Y'"};
{\ar@{->}_{\BL_G} "X";"X'"};
{\ar@{->}^{\BL_M} "Y";"Y'"};
{\ar@{=>} "Y";"X'"}
\endxy
\end{equation} 
that commutes \emph{up to a natural transformation}. 

\medskip

However, it is entirely not obvious that the natural transformation in \eqref{e:top lid nat trans Intro} is an isomorphism.
Ultimately, we establish that it is an isomorphism (\thmref{t:L and CT}), but this comes after we prove our main result, \thmref{t:main}. 

\sssec{}

What we prove \emph{a priori} is that there exists \emph{some} natural isomorphism that makes the diagram 
\begin{equation} \label{e:L and CT Intro}
\CD
\Dmod_{\frac{1}{2}}(\Bun_M) @>{\BL_M}>> \IndCoh_\Nilp(\LS_\cM) \\
@A{\on{CT}^-_{*,\rho_P(\omega_X)}}AA @AA{\on{CT}^{-,\on{spec}}}A \\
\Dmod_{\frac{1}{2}}(\Bun_G) @>{\BL_G}>> \IndCoh_\Nilp(\LS_\cG) 
\endCD
\end{equation} 
commute\footnote{We do not know, and are not sure that it is true, that the natural isomorphism in \eqref{e:L and CT Intro} equals
one of \eqref{e:top lid nat trans Intro}. One can show, however, that the two differ by a (non-zero) scalar.}. 

\medskip

The existence of the commutative diagram \eqref{e:L and CT Intro} is one of the main results of this 
paper (\thmref{t:CT compat}). It uses the commutative 
square\footnote{Up to a graded line, omitted in the Introduction.}
\begin{equation} \label{e:Langlands critical compat Intro}
\CD
\Dmod_{\frac{1}{2}}(\Bun_G) @>{\BL_G}>> \IndCoh_\Nilp(\LS_\cG) \\
@A{\Loc_G}AA @AA{\Poinc^{\on{spec}}_\cG}A \\
\KL(G)_{\crit,\Ran} @>{\FLE_{G,\crit}}>>  \IndCoh^*(\Op^{\on{mon-free}}_\cG(\cD^\times))_\Ran,
\endCD
\end{equation}
constructed in \cite[Theorem 18.5.2]{GLC2}.

\sssec{}

We prove the existence of \eqref{e:L and CT Intro} by constructing a commutative cube
(see diagram \eqref{e:CT cube}) that relates the diagram \eqref{e:Langlands critical compat Intro}
for $G$ with a similar diagram for $M$, with the crucial ingredient being the compatibility
of the critical FLE with Jacquet functors, given by \thmref{t:local Jacquet} mentioned above. 

\sssec{}

Of course, the compatibility of $\BL_G$ with the critical localization functor, expressed by diagram \eqref{e:Langlands critical compat Intro}
plays a much bigger role in this project than just proving the existence of \eqref{e:L and CT Intro}. 

\medskip

In the next paper, it will be used to show that the functor $\BL_G$ is \emph{ambidextrous} (at least on the cuspidal part),
which would be another crucial step towards the proof of the geometric Langlands conjecture.

\begin{rem}
One can say that our approach to the proof of the geometric Langlands conjecture consists of playing diagrams
\eqref{e:coeff comp L Intro} and \eqref{e:Langlands critical compat Intro} one against the other.

\medskip

Note that the approach to geometric Langlands via \eqref{e:coeff comp L Intro} was essentially the idea behind
Drinfeld's founding work \cite{Dr} (taken up by \cite{Laum} and developed further in \cite{FGV}), and the approach 
via \eqref{e:Langlands critical compat Intro} was the idea behind the work of Beilinson-Drinfeld \cite{BD1}.

\end{rem} 

\ssec{The logical structure: equivalence on Eisenstein parts}

\sssec{}

The existence of diagram \eqref{e:L and CT Intro} immediately implies that the functor $\BL_G$ admits a left
adjoint, to be denoted $\BL_G^L$, which makes the diagram 
\begin{equation} \label{e:L L and Eis Intro}
\CD
\Dmod_{\frac{1}{2}}(\Bun_M) @<{\BL^L_M}<< \IndCoh_\Nilp(\LS_\cM) \\
@V{\Eis^-_{!,\rho_P(\omega_X)}}VV @VV{\Eis^{-,\on{spec}}}V \\
\Dmod_{\frac{1}{2}}(\Bun_G) @<{\BL^L_G}<< \IndCoh_\Nilp(\LS_\cG). 
\endCD
\end{equation} 
commute. 

\medskip

Having both diagrams \eqref{e:L and Eis Intro} and \eqref{e:L L and Eis Intro} implies that the functors
$\BL_G$ and $\BL_G^L$ send the subcategories
$$\Dmod_{\frac{1}{2}}(\Bun_G)_{\Eis}\subset \Dmod_{\frac{1}{2}}(\Bun_G)$$
and
$$\IndCoh_\Nilp(\LS_\cG)_{\on{red}}\subset \IndCoh_\Nilp(\LS_\cG)$$
to one another. 

\medskip

The main result of this paper, \thmref{t:main}, says that the resulting adjoint functors
\begin{equation} \label{e:L on Eis Intro}
\BL_G:\Dmod_{\frac{1}{2}}(\Bun_G)_{\Eis}\leftrightarrows \IndCoh_\Nilp(\LS_\cG)_{\on{red}}:\BL_G^L
\end{equation} 
are mutually inverse equivalences, provided that we know that the geometric Langlands conjecture holds
for all proper Levi subgroups of $G$. 

\medskip

We will now explain the logic of the proof of this theorem.

\sssec{}

Given diagrams \eqref{e:L and Eis Intro} and \eqref{e:L L and Eis Intro}, and using the inductive hypothesis
that GLC holds for Levi quotients of proper parabolic subgroups, it suffices to show that $\BL_G^L$ is
fully faithful. 

\sssec{}

Consider the composition
$$\BL_G\circ \BL_G^L,$$
viewed as a monad acting on $\IndCoh_\Nilp(\LS_\cG)$.

\medskip

We show (\thmref{t:AG}) that this monad is given by the action of an associative algebra object
$$\CA_G\in \QCoh(\LS_\cG)$$
(we view $\IndCoh_\Nilp(\LS_\cG)$ as a module category over $\QCoh(\LS_\cG)$). 

\sssec{}

The assertion that the functor $\BL_G^L$ is fully faithful is equivalent to the assertion that the unit
\begin{equation} \label{e:unit AG Intro}
\CO_{\LS_\cG}\to \CA_G
\end{equation}
is an isomorphism.\footnote{In fact, combining the main theorem of \cite{FR} with the commutation of 
\eqref{e:L and CT Intro} and the inductive hypothesis that GLC holds for proper Levis, we obtain that the fully-faithfulness 
of $\BL_G^L$ is equivalent to the geometric Langlands conjecture.} 

\medskip

The assertion that $\BL_G^L|_{\IndCoh_\Nilp(\LS_\cG)_{\on{red}}}$ is fully faithful is equivalent to the assertion that the
map 
\begin{equation} \label{e:unit AG red Intro}
\CO_{\LS^{\on{red}}_\cG}\to \CA_G|_{\LS^{\on{red}}_\cG},
\end{equation}
induced by \eqref{e:unit AG Intro}, is an isomorphism, where $\LS^{\on{red}}_\cG\subset \LS_\cG$ is any
closed substack whose underlying subset consists of reducible local systems.

\sssec{}

The latter assertion is equivalent to the map
\begin{equation} \label{e:unit AG par Intro}
\CO_{\LS_{\cP^-}}\to (\sfq^{-,\on{spec}})^*(\CA_G)
\end{equation}
being an isomorphism for any \emph{proper} standard (negative) parabolic $\cP^-\subset \cG$, where
$$\sfq^{-,\on{spec}}:\LS_{\cP^-}\to \LS_\cG$$
is the natural projection. 

\medskip

Now, a simple but crucial observation is given by \lemref{l:line bundle is enough}, which says that 
in order to prove that \eqref{e:unit AG par Intro} is an isomorphism, it is enough to show that the object
$$(\sfq^{-,\on{spec}})^*(\CA_G) \in \QCoh(\LS_{\cP^-})$$
is a \emph{line bundle}.

\medskip

Further, we show that in order to prove that $(\sfq^{-,\on{spec}})^*(\CA_G)$ is a line bundle, it suffices to show
that its \emph{restriction} along the map
$$\LS_\cM\to \LS_{\cP^-}$$
is a line bundle on $\LS_\cM$.

\sssec{}

Now, the functor of pullback along
$$\LS_\cM\to \LS_{\cP^-}\to \LS_\cG,$$
is (up to tensoring by a line bundle) the functor of ``compactified constant term",\footnote{This is so on the subcategory 
$\QCoh(\LS_\cG)\subset \IndCoh(\LS_\cG)$, which is what is relevant for us here. The actual 
definition of $\on{CT}^{-,\on{spec}}_{!*}$ is more involved, see \secref{ss:comp spec CT}.}
denoted $\on{CT}^{-,\on{spec}}_{!*}$, introduced in \secref{ss:comp spec CT}. 

\medskip

Note that
$$\CA_G\simeq \BL_G\circ \BL_G^L(\CO_{\LS_\cG}).$$

We will prove that there exists a canonical isomorphism
\begin{equation} \label{e:!* back and forth Intro}
\on{CT}^{-,\on{spec}}_{!*}\circ \BL_G\circ \BL_G^L\simeq \on{CT}^{-,\on{spec}}_{!*},
\end{equation} 
see \propref{p:!* back and forth}, so that
$$\on{CT}^{-,\on{spec}}_{!*}(\CA_G)\simeq \on{CT}^{-,\on{spec}}_{!*}(\CO_{\LS_\cG}).$$

\medskip

This will establish the desired property of $(\sfq^{-,\on{spec}})^*(\CA_G)$, since this property clearly holds for
$\CO_{\LS_\cG}$. 

\sssec{}

In order to prove \eqref{e:!* back and forth Intro}, we break it into two isomorphisms
\begin{equation} \label{e:!* back Intro}
\on{CT}^{-,\on{spec}}_{!*}\circ \BL_G\simeq \BL_M\circ \on{CT}^-_{!*,\rho_P(\omega_X)}
\end{equation}
and 
\begin{equation} \label{e:!* forth Intro}
\on{CT}^-_{!*,\rho_P(\omega_X)}\circ \BL^L_G \simeq  \BL^L_M\circ \on{CT}^{-,\on{spec}}_{!*},
\end{equation}
for a suitably defined functor
$$\on{CT}^-_{!*,\rho_P(\omega_X)}:\Dmod_{\frac{1}{2}}(\Bun_G)\to \Dmod_{\frac{1}{2}}(\Bun_M).$$

Then \eqref{e:!* back and forth Intro} would follow by the induction hypothesis that the geometric
Langlands conjecture holds for $M$, and hence $\BL_M\circ \BL_M^L\simeq \on{Id}$.

\sssec{} \label{sss:compactified compat Intro}

The upshot of the above discussion is that in addition to the commutative diagram \eqref{e:L and CT Intro},
one needs to have a commutative diagram that involves the pair
$(\on{CT}^-_{!*,\rho_P(\omega_X)},\on{CT}^{-,\on{spec}}_{!*})$.
This would imply the isomorphism \eqref{e:!* back Intro}.

\medskip

It will turn out by duality considerations that \eqref{e:!* forth Intro} is equivalent to a cousin of \eqref{e:L and Eis Intro},
where instead of $(\Eis^-_{!,\rho_P(\omega_X)},\Eis^{-,\on{spec}})$ one uses suitably defined functors
$(\Eis^-_{!*,\rho_P(\omega_X)},\Eis^{-,\on{spec}}_{!*})$.  

\medskip

This brings us to the idea of \emph{enhancement}. 

\ssec{All this enhancement business}

\sssec{}

Informally speaking, the enhancement procedure(s)
$$\Eis^-_!\rightsquigarrow \Eis^-_{!*}, \quad \Eis^{-,\on{spec}}\rightsquigarrow \Eis^{-,\on{spec}}_{!*}$$
and
$$\on{CT}^-_*\rightsquigarrow \on{CT}^-_{!*}, \quad \on{CT}^{-,\on{spec}}\rightsquigarrow \on{CT}^{-,\on{spec}}_{!*}$$
is a mechanism to encode the (pre)compositions of the original functors with the Hecke functors for $G$ and $M$
(resp., $\cG$ and $\cM$). 

\sssec{}

Technically, this idea is realized as follows. We introduce the \emph{local semi-infinite categories}, denoted 
\begin{equation} \label{e:semiinf cat Intro}
\on{I}(G,P^-)^{\on{loc}} \text{ and } \on{I}(\cG,\cP^-)^{\on{spec,loc}},
\end{equation} 
respectively.

\medskip

The category on the geometric side, i.e., $\on{I}(G,P^-)^{\on{loc}}$, is the category of (half-twisted) D-modules on  
$\Gr_G$, satisfying the parabolic semi-infinite equivariance condition, i.e., equivariance 
with respect to $\fL(N^-_P)\cdot \fL^+(M)$, up to a small renormalization procedure.  

\medskip

On the spectral side, the category $\on{I}(\cG,\cP^-)^{\on{spec,loc}}$ is that of (appropriately defined) 
ind-coherent sheaves on the fiber product
\begin{equation} \label{e:loc spe Intro}
\LS_\cM(\cD)\underset{\LS_\cM(\cD^\times)}\times \LS_{\cP^-}(\cD^\times)\underset{\LS_\cG(\cD^\times)}\times \LS_\cG(\cD),
\end{equation} 
where $\cD$ (resp., $\cD^\times$) is the formal (resp., formal punctured) disc. 

\medskip

The categories \eqref{e:semiinf cat Intro} are acted on by
$$\Sph_G\otimes \Sph_M \text{ and } \Sph^{\on{spec}}_\cG\otimes \Sph^{\on{spec}}_\cM,$$ respectively.
So they give rise to functors between
$$\Sph_G\mmod \leftrightarrow \Sph_M\mmod \text{ and } \Sph^{\on{spec}}_\cG\mmod\leftrightarrow \Sph^{\on{spec}}_\cM\mmod,$$
respectively. 

\medskip

One forms the global enhanced categories\footnote{The tensor products below should be spread over the Ran space,
a point which we elide in the Introduction.} 
$$\Dmod_{\frac{1}{2}}(\Bun_G)^{-,\on{enh}}:=
\Dmod_{\frac{1}{2}}(\Bun_G)\underset{\Sph_G}\otimes \on{I}(G,P^-)^{\on{loc}},$$
$$\Dmod_{\frac{1}{2}}(\Bun_M)^{-,\on{enh}}:=
\Dmod_{\frac{1}{2}}(\Bun_M)\underset{\Sph_M}\otimes \on{I}(G,P^-)^{\on{loc}}$$
and 
$$\IndCoh_\Nilp(\LS_\cG)^{-,\on{enh}}:=\IndCoh_\Nilp(\LS_\cG)\underset{\Sph^{\on{spec}}_\cG}\otimes \on{I}(\cG,\cP^-)^{\on{spec,loc}},$$
$$\IndCoh_\Nilp(\LS_\cM)^{-,\on{enh}}:=\IndCoh_\Nilp(\LS_\cM)\underset{\Sph^{\on{spec}}_\cM}\otimes \on{I}(\cG,\cP^-)^{\on{spec,loc}},$$
and we show that the enhanced Eisenstein and constant term functors naturally extend to functors between the enhanced categories.

\sssec{}

We show that the (derived) geometric Satake equivalences
$$\Sph_G\overset{\Sat_G}\simeq \Sph_\cG^{\on{spec}} \text{ and } 
\Sph_M\overset{\Sat_M}\simeq \Sph_\cM^{\on{spec}}$$
extend to an equivalence
$$ \on{I}(G,P^-)^{\on{loc}}\overset{\Sat^{-,\semiinf}}\simeq \on{I}(\cG,\cP^-)^{\on{spec,loc}},$$
and we show that the enhanced Langlands functor(s)
$$\BL_G^{-,\on{enh}}:\Dmod_{\frac{1}{2}}(\Bun_G)^{-,\on{enh}}\to \IndCoh_\Nilp(\LS_\cG)^{-,\on{enh}}$$
and 
$$\BL_M^{-,\on{enh}}:\Dmod_{\frac{1}{2}}(\Bun_M)^{-,\on{enh}}\to \IndCoh_\Nilp(\LS_\cM)^{-,\on{enh}}$$
are compatible with the enhanced versions of the Eisenstein and constant term functors. 

\sssec{}

Now, the point of all this is that one can express the ``compactified" functors 
$$\Eis^-_{!*},\,\, \Eis^{-,\on{spec}}_{!*},\,\, \on{CT}^-_{!*},\, \on{CT}^{-,\on{spec}}_{!*}$$
in terms\footnote{This expression is provided using the \emph{semi-infinite IC sheaf} constructed in \cite{Ga4,Ga5}.} 
of the enhanced functors
$$\Eis^{-,\on{enh}}_!,\,\, \Eis^{-,\on{spec,enh}},\,\, \on{CT}^{-,\on{enh}}_*,\,\, \on{CT}^{-,\on{spec,enh}},$$
respectively,
and the enhanced compatibility of the Langlands functor with the Eisenstein and constant term functors
implies the compatibility with the ``compactified" functors, alluded to in \secref{sss:compactified compat Intro}. 

\sssec{}

The above explains why we need the enhancements for the proof of \thmref{t:main}. 

\medskip

However, there is one
more aspect that makes these enhancements necessary: they are used in the proof in the main local result,
\thmref{t:local Jacquet}, which expresses the compatibility of the FLE with Jacquet functors, see \secref{ss:local Intro}
below. 

\sssec{}

We conclude this subsection with the following remark. The enhancement 
$$\Eis^-_!\rightsquigarrow \Eis^{-,\on{enh}}_!$$
is in fact something natural when one traces back the path from sheaves to functions. 

\medskip

Recall that in the classical theory of automorphic functions, the Eisenstein operator is defined as pull-push along the following diagram
\begin{equation} \label{e:Eis classical}
\vcenter
{\xy
(0,0)*+{G(\BA)/G(K)}="X";
(40,0)*+{G(\BA)/N^-_P(\BA)\cdot M(K).}="Y";
(20,20)*+{G(\BA)/P^-(K)}="Z";
{\ar@{->} "Z";"X"};
{\ar@{->} "Z";"Y"};
\endxy}
\end{equation}

This is a map of $G(\BA)$-representations. When one takes spherical vectors everywhere, one obtains the diagram
\begin{equation} \label{e:Eis classical sph}
\vcenter
{\xy
(0,0)*+{G(\BO)\backslash G(\BA)/G(K)}="X";
(40,0)*+{G(\BO)\backslash G(\BA)/N^-_P(\BA)\cdot M(K).}="Y";
(20,20)*+{G(\BO)\backslash G(\BA)/P^-(K)}="Z";
{\ar@{->} "Z";"X"};
{\ar@{->} "Z";"Y"};
\endxy}
\end{equation}

Now, since for a local field $\CK$ with local ring $\CO$, we have
\begin{equation} \label{e:Iwasawa}
G(\CO)\backslash G(\CK)\simeq P^-(\CO)\backslash P^-(\CK),
\end{equation} 
we have
$$G(\BO)\backslash G(\BA)\simeq P^-(\BO)\backslash P^-(\BA),$$
and hence we can rewrite \eqref{e:Eis classical sph} as 
\begin{equation} \label{e:Eis classical naive}
\vcenter
{\xy
(0,0)*+{G(\BO)\backslash G(\BA)/G(K)}="X";
(40,0)*+{M(\BO)\backslash M(\BA)/M(K).}="Y";
(20,20)*+{P^-(\BO)\backslash P^-(\BA)/P^-(K)}="Z";
{\ar@{->} "Z";"X"};
{\ar@{->} "Z";"Y"};
\endxy}
\end{equation}

Now, if we take the geometric counterpart of \eqref{e:Eis classical naive}, we obtain the diagram
$$
\vcenter
{\xy
(0,0)*+{\Bun_G}="X";
(40,0)*+{\Bun_M,}="Y";
(20,20)*+{\Bun_{P^-}}="Z";
{\ar@{->} "Z";"X"};
{\ar@{->} "Z";"Y"};
\endxy}
$$
which is the diagram that defines the functor $\Eis^-_!$.

\medskip

By contrast, the geometric counterpart of \eqref{e:Eis classical sph} gives rise to the functor $\Eis^{-,\on{enh}}_!$. 

\medskip

The point is that the geometric counterpart of the set-theoretic isomorphism
$$G(\CO)\backslash G(\CK)/N^-_P(\CK)\cdot M(\CO)\simeq 
P^-(\CO)\backslash P^-(\CK)/N^-_P(\CK)\cdot M(\CO)= 
M(\CO)\backslash M(\CK)/M(\CO)$$
is a \emph{stratification} of $\Gr_G$ by orbits of $\fL(N^-_P)\cdot \fL^+(M)$, so that 
the category $\on{I}(G,P^-)^{\on{loc}}$ is related to $\Sph_M$ by a monadic adjunction,
which is not at all an equivalence. 

\ssec{The local FLE/Jacquet compatibility} \label{ss:local Intro}

In this subsection we describe the main local result, \thmref{t:local Jacquet}, whose proof occupies 
most of Part I of the paper.

\sssec{}

In \cite{GLC2} we constructed the factorization category $\IndCoh^*(\Op^\mf_\cG)$, and we established
the \emph{Fundamental Local Equivalence} (the FLE) 
$$\FLE_{G,\crit}:\KL(G)_\crit \simeq \IndCoh^*(\Op^\mf_\cG).$$

There is a similar equivalence for $M$, and also one that incorporates a ``$\rho$"-shift (see \eqref{e:FLE M crit}): 
\begin{equation} \label{e:FLE M crit Intro}
\FLE_{M,\crit_M+\rhoch_P}:\KL(M)_{\crit_M+\rhoch_P} \simeq \IndCoh^*(\Op_{\cM,\rhoch_P}^{\on{mon-free}}).
\end{equation} 

\sssec{}

In \secref{sss:BRST summary} we construct a Jacquet functor
$$J^{-,\Sph}_{\on{KM},\rho_P(\omega_X)}:\KL(G)_{\crit_G}\to \KL(M)_{\crit_M+\rhoch_P}.$$

\medskip

In \secref{sss:spectral Jacquet Theta} we construct a Jacquet functor
$$J^{-,\on{spec}}_{\Op,\Theta}:\IndCoh^*(\Op_\cG^{\on{mon-free}})\to \IndCoh^*(\Op_{\cM,\rhoch_P}^{\on{mon-free}}).$$

\medskip

\thmref{t:local Jacquet} says that the diagram of factorization categories 
\begin{equation} \label{e:local Jacquet Intro}
\CD
\KL(M)_{\crit_M+\rhoch_P} @>{\on{FLE}_{M,\crit_G+\rhoch_P}}>> \IndCoh^*(\Op_{\cM,\rhoch_P}^{\on{mon-free}}) \\
@A{J^{-,\Sph}_{\on{KM},\rho_P(\omega_X)}}AA @AA{J^{-,\on{spec}}_{\Op,\Theta}}A \\
\KL(G)_{\crit_G} @>{\on{FLE}_{G,\crit_G}}>>  \IndCoh^*(\Op_\cG^{\on{mon-free}}).
\endCD
\end{equation} 
commutes. 

\sssec{}

\thmref{t:local Jacquet} is the most difficult among the local results in this paper and its prequel (\cite{GLC2})
and it digs deep into the original works of Feigin and Frenkel on 
Wakimoto modules (with some developments of the latter in \cite{FG1}).

\medskip

In order to prove \thmref{t:local Jacquet}, we replace \eqref{e:local Jacquet Intro} by a dual diagram
\begin{equation}  \label{e:dual Jacquet diagram Intro}
\CD
\KL(M)_{\crit_G-\rhoch_P} @>{\FLE_{M,\crit_G+\rhoch_P}\circ \tau_M}>> \IndCoh^*(\Op_{\cM,\rhoch_P}^{\on{mon-free}}) \\
@V{\Wak^{-,\Sph}_{\rho_P(\omega_X)}}VV  @VV{\on{co}\!J_\Op^{-,\on{spec}}}V \\ 
\KL(G)_{\crit_G}   @>{\on{FLE}_{G,\crit_G}\circ \tau_G}>>   \IndCoh^*(\Op_\cG^{\on{mon-free}}),
\endCD
\end{equation} 
where $\Wak^{-,\Sph}_{\rho_P(\omega_X)}$ is the \emph{spherical} Wakimoto functor. 

\medskip

Further, in order to establish the commutation of \eqref{e:dual Jacquet diagram Intro}, we form its enhancement

\medskip

\begin{equation}  \label{e:dual Jacquet enh diagram Intro}
\CD
\KL(M)_{\crit_G-\rhoch_P} @>{\FLE_{M,\crit_G+\rhoch_P}\circ \tau_M}>> \IndCoh^*(\Op_{\cM,\rhoch_P}^{\on{mon-free}}) \\
@V{\Wak^{-,\on{enh}}_{\rho_P(\omega_X)}}VV @VV{\on{co}\!J^{-,\on{spec},\on{enh}}_\Op}V \\
\KL(G)_{\crit_G}\underset{\Sph_G}\otimes \on{I}(G,P^-)^{\on{loc}}_{\rho_P(\omega_X)} @>{\FLE^{-,\on{enh}}_{G,\crit}}>>
\IndCoh^*(\Op_\cG^{\on{mon-free}}) \underset{\Sph_\cG^{\on{spec}}}\otimes \on{I}(\cG,\cP^-)^{\on{spec,loc}},
\endCD
\end{equation}
where the lower horizontal arrow is the \emph{enhancement} of the functor $\FLE_{G,\crit}$. This is the second essential 
usage of the enhancement procedure in this paper.\footnote{We also need (the dual of) the enhanced diagram \eqref{e:dual Jacquet enh diagram Intro}
to prove the compatibility of the \emph{enhanced} Langlands functor with constant terms.}  

\medskip

In turns out that, although the terms of \eqref{e:dual Jacquet enh diagram Intro} look more complicated than those 
in \eqref{e:dual Jacquet diagram Intro}, this diagram is more manageable since the functor $\Wak^{-,\on{enh}}_{\rho_P(\omega_X)}$
is more directly related to the classical Wakimoto functor. 

\sssec{}

The proof of the commutation of \eqref{e:dual Jacquet enh diagram Intro} uses techniques from the theory of factorization algebras.

\medskip

We map the lower right corner, i.e., the category $\IndCoh^*(\Op_\cG^{\on{mon-free}}) \underset{\Sph_\cG^{\on{spec}}}\otimes \on{I}(\cG,\cP^-)^{\on{spec,loc}}$
to the category of factorization modules over a certain factorization algebra 
$$\Omega(R)^\Op \in \IndCoh^*(\Op^\mer_\cG)\otimes \Rep(\cM).$$

This functor turns out to be \emph{almost an equivalence} (technically, an equivalence on eventually coconnective subcategories). The bulk
of the proof of the commutation of \eqref{e:dual Jacquet enh diagram Intro} consists of showing that the two circuits in 
\eqref{e:dual Jacquet enh diagram Intro} become isomorphic when composed with the above functor. 

\medskip

At the end of the day, this comes down to the computation of the Drinfeld-Sokolov reduction of Wakimoto modules, 
together with an action of the Feigin-Frenkel center (which is given by Miura transform). 

\sssec{}

One can say that the above calculation is the computational core of this paper, hence the dedication to B.~Feigin. 

\ssec{Description of the actual contents}

This paper consists of three Parts. 

\sssec{Contents of Part I}

This Part is devoted to the study of Jacquet functors (and their duals) on the two sides
of Langlands correspondence. And on each side there are two pairs of such functors:
Whittaker-style and Kac-Moody style. 

\medskip

In \secref{s:semiinf} we introduce and study the \emph{local semi-infinite category} on the geometric side, $\on{I}(G,P^-)^{\on{loc}}$.
It serves as a connector between categories associated with $G$ and $M$. 

\medskip

In \secref{s:spectral semiinf} we introduce the spectral counterpart of $\on{I}(G,P^-)^{\on{loc}}$, denoted 
$\on{I}(\cG,\cP^-)^{\on{spec,loc}}$. 

\medskip

Over a fixed formal disc, this is the category of ind-coherent sheaves
on \eqref{e:loc spe Intro}; it is well-defined because this fiber product is a (derived) algebraic stack locally
of finite type. However, when we work in families over the Ran space, we no longer have the locally
of finite type property, and special care is needed in order to define the category of ind-coherent sheaves
(in fact, there are two such categories that one can consider: $\IndCoh^!(-)$ and $\IndCoh^*(-)$). These issues are dealt with
in \secref{s:semiinf spec}.

\medskip

Coming back to \secref{s:spectral semiinf}, we state the existence of an equivalence \eqref{e:semiinf cat Intro}
(\thmref{t:semiinf geom Satake}), and a closely related assertion, \thmref{t:semiinf CS}. The proofs of these
theorems is delegated to the Appendix, Sects. \ref{s:Sith}-\ref{s:semiinf Sat}. 

\medskip

In \secref{s:KM} we study various versions of the Jacquet functors between
the categories of Kac-Moody modules for $\fg$ and $\fm$. These Jacquet functors are obtained essentially by decorating 
the functor of BRST reduction with respect to $\fL(\fn^-_P)$. We also study the dual functors, which yield various 
versions of the Wakimoto module construction. 

\medskip

In \secref{s:Op} we study Jacquet functors between the various categories of ind-coherent sheaves on opers 
for $\cG$ and $\cM$. The central tool is a geometric object called \emph{parabolic Miura opers}, which appeared
prominently in the works of Feigin and Frenkel. 

\medskip

In \secref{s:FLE and Jacquet} we state our main local result, \thmref{t:local Jacquet}. We reformulate it in dual 
terms, and also state the enhanced version, \thmref{t:local Jacquet dual enh}. The proof of 
\thmref{t:local Jacquet dual enh} spans Sects. \ref{s:proof local Jacquet dual enh} and \ref{s:ff in plus}. 

\medskip

In \secref{s:proof local Jacquet dual enh} we supply the computational/representation-theoretic part of the proof, 
with one statement, \thmref{t:FG1}, delegated to the Appendix, \secref{s:proof of Wak}. 

\medskip

In \secref{s:ff in plus} we deal with convergence issues. 

\sssec{Contents of Part II}

This Part is thematically close to Part II of \cite{GLC2} in that it deals 
with the interactions between Eisenstein and constant term functors (on either side of 
Langlands correspondence) with local-to-global functors studied in {\it loc. cit.} (there are two such on each side). 

\medskip

We note
that \emph{none} of the results from Part II cross the Langlands bridge, i.e., we work on the geometric and
spectral sides separately. We first deal with the geometric side (Sects. \ref{s:Eis}-\ref{s:Eis and Loc})
and then with the spectral side (Sects. \ref{s:spec CT}-\ref{s:Eis and Op}). 

\medskip

In \secref{s:Eis} we introduce the geometric Eisenstein and constant term functors, as well as their
enhancements.

\medskip 

In \secref{s:Eis glob} we interpret the enhanced Eisenstein and constant term functors via
Drinfeld's compactifications. This interpretation is needed for the next section, where we compute 
constant terms of \emph{Poincar\'e functors}. 

\medskip 

In \secref{s:Eis and Whit} we perform the first calculation of the type: ``what do we get when we compose 
a local-to-global functor with the constant term functor". Namely, we compute the composition of the Poincar\'e functor
with the constant term functor. The result of \thmref{t:CT co of Poinc} says that the resulting functor equals the
composition of a local Jacquet functor (from the Whittaker category of $G$ to that of $M$), followed by the Poincar\'e functor
for $M$. 

\medskip

However, there is a caveat: we need to precompose the local functor with the functor of ``insertion of the unit
along the marked Ran space."  This feature will be common to all four instances of such calculations, performed
in Part II. 

\medskip

Additionally, in \secref{s:Eis and Whit} we perform a similar calculation, but for the enhanced functors. 

\medskip 

In \secref{s:Eis and Loc} we calculate the composition of the localization functor $$\Loc_G:\KL(G)_{\crit,\Ran}\to \Dmod_{\frac{1}{2}}(\Bun_G)$$
with the constant term functor. The answer is given again as a composition of a local Jacquet functor 
(from $\KL(G)_\crit$ to $\KL(M)_\crit$), followed by the localization functor for $M$, modulo the same caveat as was mentioned above. 
We also perform an analogous computation for the enhanced version of the functors involved. 

\medskip

In \secref{s:spec CT} we introduce the spectral Eisenstein and constant term functors, as well as their enhancements. 
We calculate the composition of the spectral localization functor $$\Loc^{\IndCoh,\on{spec}}_\cG:\Rep(\cG)_\Ran\to \IndCoh(\LS_\cG)$$
with the spectral constant term functor (and also for the enhanced versions of these functors). 

\medskip

In \secref{s:Eis and Op} we calculate the composition of the \emph{spectral Poincar\'e functor} and the spectral constant 
term functor. 

\sssec{Contents of Part III}

In this Part, we finally turn to problems of Langlands duality in the global setting, by combining the results of Parts I and II.

\sssec{}

In \secref{s:L and Eis} we prove that the Langlands functor is compatible with the Eisenstein functors, i.e., 
we establish the commutation of the diagram \eqref{e:L and Eis Intro}, as well as its enhanced version.

\medskip

In \secref{s:L and CT} we prove that the Langlands functor is compatible with the constant term functors, i.e., 
we establish the commutation of the diagram \eqref{e:L and CT Intro}, as well as its enhanced version.

\medskip

In \secref{s:left adj}, we prove the following three results: 

\smallskip

\begin{itemize}

\item The Langlands functor $\BL_G$ admits a left adjoint (denoted $\BL_G^L$; moreover, the functor $\BL_G^L$
is compatible with the Eisenstein functors);

\smallskip

\item The left adjoint of the Langlands functor is isomorphic to the dual of the Langlands functor, up to Miraculous
duality and the Chevalley involution;

\smallskip

\item The composition $\BL_G\circ \BL_G^L$, viewed as an endofunctor of $\IndCoh_\Nilp(\LS_\cG)$
is given by tensor product with an object\footnote{A priori, this could have been an object of 
$\IndCoh_{\Nilp\times \Nilp}(\LS_\cG\times \LS_\cG)\underset{\QCoh(\LS_\cG)\otimes \QCoh(\LS_\cG)}\otimes \QCoh(\LS_\cG)$.}
of $\QCoh(\LS_\cG)$.

\end{itemize}

\medskip

Finally, in \secref{s:Eis equiv} we state and prove the main result of this paper, \thmref{t:main}, which says that the Langlands
functor $\BL_G$ induces an equivalence between the Eisenstein-generated subcategories of the two sides. 

\ssec{Conventions and notation}

\sssec{}

Notations and conventions in this paper follow those adopted in \cite{GLC2}. 
We only need to introduce those that have to do with parabolic subgroups of $G$.

\medskip

To simplify the discussion, we choose a Cartan subgroup $T\subset G$.

\sssec{}

For a standard parabolic subgroup $P$ we denote by $P^-$ its negative, so that $M:=P\cap P^-$ is a Levi
subgroup in both. We let $N_P$ (resp., $N^-_P$) the unipotent radical of $P$ (resp., $P^-$). 

\medskip

We will identify the Levi \emph{quotients} of $P$ and $P^-$ by means of 
$$P/N_P \simeq P\cap P^- \simeq P^-/N^-_P.$$

\sssec{}

We let $\Lambda$ denote the \emph{coweight} lattice of $G$. We denote by $\Lambda_G^+\subset \Lambda$
the submonoid of dominant coweights. 

\medskip

We let $\cLambda$ denote the weight lattice of $G$, and we let $\cLambda^+_G$ the submonoid of dominant weights.

\sssec{}

We let $\cP$ (resp., $\cP^-$) the positive (resp., negative) parabolic in $\cG$, corresponding to $P$ (resp., $P^-$),
so that $\cM$ is the Langlands dual of $M$. 

\medskip

We let $\cN_P$ (resp., $\cN^-_P$) denote the unipotent radical of $\cP$ (resp., $\cP^-$). We denote by
$\cn_P$ and $\cn^-_P$ the corresponding Lie algebras. 

\medskip

We let $B(M),N(M),B^-(M),N^-(M)$ (resp., $\cB(M),\cN(M),\cB^-(M),\cN^-(M)$) the corresponding subgroups 
of $M$ (resp., $\cM$). 

\sssec{Why this minus sign?}

The reason that in most of the paper we work with the \emph{negative} parabolic $P^-$ is that when consider its 
interaction with the Whittaker model (which is taken with respect to the \emph{positive} $N$), the open Bruhat
cell corresponds to the unit element of the Weyl group. 

\medskip

In other words, if we worked with $P$, we would have
the longest element $w_0\in W$ enter a lot of formulas.  

\ssec{Acknowledgements}

First, we wish to thank our collaborators on the project of proving GLC, D.~Arinkin and D.~Beraldo. 

\medskip

In addition, we wish to acknowledge the contributions of the following mathematicians
(in the order in which their ideas appear in the paper): 

\medskip

The Langlands-type equivalences of local semi-infinite categories given by Theorems 
\ref{t:semiinf geom Satake} and \ref{t:semiinf CS}, at the pointwise level, reproduce results
from Bezrukavnikov's theory. His results was the reason that we thought these equivalences
were possible in the first place. 

\medskip 

The theory of Wakimoto and BRST functors, as well as the notion of Miura opers, 
which play a crucial role in Part I of the paper, 
was pioneered and developed by B.~Feigin and E.~Frenkel.

\medskip

The study of Eisenstein series in the geometric context was pioneered by G.~Laumon. 

\medskip

The constant term/Whittaker calculations are crucially based on the idea of Zastava spaces,
introduced by M.~Finkelberg and I.~Mirkovi\'c. 

\medskip

Drinfeld's compactifications, which are extensively used in Sects. \ref{s:Eis and Whit} 
were (obviously) introduced by V.~Drinfeld. 

\medskip

The idea that factorization homology appears in local-to-global calculations for the localization functor 
goes back to A.~Beilinson. 

\medskip

The idea of enhancement of Eisenstein series originates in the joint work of the third-named author with
A.~Braverman, and is strongly influenced by ideas taught to him by J.~Bernstein. 

\medskip

The third- and fourth-named authors wish to thank IHES, where a significant part of this paper was written. 

\medskip

The work of D.G. was supported by NSF grant DMS-2005475. 
The work of S.R. was supported by a Sloan Research Fellowship and NSF grants DMS-2101984 and DMS-2416129 
while this work was in preparation.

\newpage 

\centerline{\bf Part I: Jacquet functors in the local theory}

\bigskip



The key players in this Part are the \emph{local semi-infinite category}, denoted $\on{I}(G,P^-)^{\on{loc}}_{\rho_P(\omega_X)}$,
and its spectral counterpart, denoted $\on{I}(\cG,\cP^-)^{\on{spec,loc}}$. 

\medskip

Objects of $\on{I}(G,P^-)^{\on{loc}}_{\rho_P(\omega_X)}$ can be used as \emph{generalized Jacquet functors}
$$\Whit^!(G)\to \Whit^!(M)$$
and objects of $\on{I}(\cG,\cP^-)^{\on{spec,loc}}$ can be used as \emph{generalized spectral Jacquet functors}
$$\Rep(\cG)\to \Rep(\cM).$$

The first main result of this Part is \thmref{t:semiinf geom Satake}, which says that there is a canonically defined
equivalence of factorization categories
$$\Sat^{-,\semiinf}:\on{I}(G,P^-)^{\on{loc}}_{\rho_P(\omega_X)} \simeq \on{I}(\cG,\cP^-)^{\on{spec,loc}},$$
essentially determined by the fact that it intertwines the Jacquet functors on the two sides via
the Casselman-Shalika equivalences for $G$ and $M$, respectively. The proof
of this theorem is delegated to the Appendices, Sects. \ref{s:Sith}-\ref{s:semiinf Sat}. 

\medskip

In addition, objects of (the dual of) $\on{I}(G,P^-)^{\on{loc}}_{\rho_P(\omega_X)}$ can be used as  \emph{generalized Jacquet functors}
$$\KL(G)_{\crit}\to \KL(M)_\crit$$
and objects of (the dual of) $\on{I}(\cG,\cP^-)^{\on{spec,loc}}$ can be used as \emph{generalized Jacquet functors} for opers
$$\IndCoh^*(\Op^\mf_\cG)\to \IndCoh^*(\Op^\mf_\cM).$$

The second main result of this Part, \thmref{t:local Jacquet enh}, says that the equivalence $\Sat^{-,\semiinf}$ intertwines
the corresponding Jacquet functors via the FLE equivalences for $G$ and $M$, respectively. 

\medskip

Both of the above compatibilities would be used in a crucial way in Part II in order to establish the compatibility of the
Langlands functor with the functors of Eisenstein series and constant term. 

\bigskip

\section{The local semi-infinite category} \label{s:semiinf}

In this section we study the category that ultimately allows one to connect categories 
(on the local geometric side) associated with the group $G$ and corresponding 
categories for its Levi subgroups.

\medskip

The category in question is the \emph{local spherical semi-infinite category}, which is a \emph{twisted, translated and renormalized version} of 
\begin{equation} \label{e:semiinf-naive}
\Dmod(\Gr_G)^{-,\semiinf}:= \Dmod(\Gr_G)^{\fL(N^-_P)\cdot \fL^+(M)}.
\end{equation} 

\medskip

The twist and the translation referred to above are familiar from \cite{GLC2}. The twist has to do with the fact that we deal not 
with plain D-modules, but with critically twisted ones. The translation has to do with the $T$-bundle $\rho(\omega_X)$, see 
\cite[Sect. 1.2.4]{GLC2}. The renormalization has also appeared in \cite[Sect. 1.5.2]{GLC2}, but will be reviewed in more detail
in \secref{ss:renorm}.



%

\ssec{The corrected Jacquet functor}

In addition to the twisting, translation and renormalization mentioned above, there is one tweak that is 
inherent to the passage between $G$ and $M$ on the geometric side: a certain cohomological 
shift built into Jacquet functors. The need for this shift appears via the requirement of compatibility 
with factorization. 

\medskip

This shift will be introduced in this subsection. 

\sssec{}

Let $P^-$ be the (negative) standard parabolic of $G$ with Levi quotient $M$. Consider the restrictions
of the line bundles 
$$\det_{\Gr_G} \text{ and } \det_{\Gr_M}$$
along the maps
\begin{equation} \label{e:Jacquet diagram Gr}
\Gr_G \overset{\sfp^-}\leftarrow \Gr_{P^-}\overset{\sfq^-}\to \Gr_M,
\end{equation} 
respectively.

\medskip

Denote by $\det_{\Gr_{G,M}}$ their ratio:
$$\det_{\Gr_{G,M}}=\det_{\Gr_G}|_{\Gr_{P^-}}\otimes (\det_{\Gr_M}|_{\Gr_{P^-}})^{\otimes -1};$$
it naturally descends (as a factorization line bundle) to $\Gr_M$. By a slight abuse of notation,
we will denote the resulting line bundle on $\Gr_M$ by the same symbol $\det_{\Gr_{G,M}}$. 

\sssec{}

We consider $\det_{\Gr_{G,M}}$ as an (evenly) graded line bundle on $\Gr_{M}$, so that its portion over 
the connected component $\Gr_M^\lambda$ has grading
$$2\langle \lambda,2\rhoch_P\rangle,$$
where $2\rhoch_P$ is the character of $M$ equal to the determinant of its action on $\fn(P)$.

\medskip

It was shown in \cite[Sect. 5.2]{GLys} that $\det_{\Gr_{G,M}}$ admits a canonical square root,\footnote{Which depends on the choice of a square root
of $\omega_X$.}
to be denoted $\on{det}^{\otimes \frac{1}{2}}_{\Gr_{G,M}}$, 
viewed as a $\BZ$-graded and hence $\BZ/2\BZ$-graded (=super) factorization line bundle, so that its portion over 
the connected component $\Gr_M^\lambda$ has grading
$$\langle \lambda,2\rhoch_P\rangle.$$



\sssec{} \label{sss:triv ratio gerbe}

Ignoring the grading, 
The line bundle $\det^{\otimes \frac{1}{2}}_{\Gr_{G,M}}$ gives rise to an identification of the pullbacks of the gerbes
$$\det^{\frac{1}{2}}_{\Gr_G} \text{ and } \det^{\frac{1}{2}}_{\Gr_M}$$
to $\Gr_{P^-}$. 

\medskip

Due to the above identification of gerbes, we have a well-defined functor 
$$J_{\Gr}^{-,*,\on{nv}}:\Dmod_{\frac{1}{2}}(\Gr_G)\to \Dmod_{\frac{1}{2}}(\Gr_M),$$
given by !-pull and *-push along \eqref{e:Jacquet diagram Gr}.

\medskip

We will refer to $J_{\Gr}^{-,*,\on{nv}}$ as the ``naive" Jacquet functor. 

\sssec{} \label{sss:naive Jacquet}

However, due to the fact that $\det^{\otimes \frac{1}{2}}_{\Gr_{G,M}}$ does not factorize as a line bundle, but only 
as a super line bundle, the above identification of gerbes is \emph{incompatible} with factorization.

\medskip

Hence, the functor $J_{\Gr}^{-,*,\on{nv}}$ is \emph{not} compatible with factorization either.

\sssec{} \label{sss:corrected Jacquet}

We introduce the corrected Jacquet functor 
$$J_{\Gr}^{-,*}:\Dmod_{\frac{1}{2}}(\Gr_G)\to \Dmod_{\frac{1}{2}}(\Gr_M),$$
as follows: 

\medskip

Over a connected component $\Gr_M^\lambda$, we set
$$J_{\Gr}^{-,*}(-):=J_{\Gr}^{-,*,\on{nv}}(-)[\langle \lambda,2\rhoch_P\rangle].$$

The functor $J_{\Gr}^{-,*}$ is naturally compatible with the factorization structures, due to the sign rules. 

\ssec{Digression: renormalized equivariant categories} \label{ss:renorm}

\sssec{}

Let $\bC$ be a category acted on by an algebraic group $H$. Assume that the category 
$\bC^H$ is compactly generated. 

\medskip

We let $\bC^{H,\on{ren}}$ denote the following renormalized version of the category $\bC^H$:

\medskip

By definition, $\bC^{H,\on{ren}}$ is the ind-completion of the full subcategory of $\bC^H$,
consisting of objects, whose images under the forgetful functor
$$\bC^H\to \bC$$
are compact. 

\medskip

The ind-extension of the inclusion 
$$(\bC^{H,\on{ren}})^c\hookrightarrow \bC^H$$
is a functor
\begin{equation} \label{e:unren}
\bC^{H,\on{ren}}\to  \bC^H.
\end{equation}

\sssec{}

Note that when $H$ is unipotent, an object of $\bC^H$ is compact \emph{if} (and only if) if its image in $\bC$
is compact, so in this case, \eqref{e:unren} is an equivalence. 

\sssec{}

We have a tautological inclusion 
$$(\bC^H)^c\subset (\bC^{H,\on{ren}})^c,$$
whose ind-extension is a (fully faithful) functor
$$\bC^H\to \bC^{H,\on{ren}},$$
left adjoint to \eqref{e:unren}. 

\medskip

I.e., we obtain $\bC^H$ is a colocalization of $\bC^{H,\on{ren}}$:
\begin{equation} \label{e:ren coloc}
\bC^H\rightleftarrows \bC^{H,\on{ren}}.
\end{equation}. 

\sssec{Example}

Take $\bC=\Vect$, equipped with the trivial action of $H$. Note that $$\Vect^H\simeq \Dmod(\on{pt}/H).$$

\medskip

Then $\Vect^{H,\on{ren}}$ is compactly generated by $k\in \bC^H$, and so
$$\Vect^{H,\on{ren}}\simeq \on{C}^\cdot(\on{pt}/H)\mod.$$

The full subcategory
$$\Vect^H\subset \Vect^{H,\on{ren}}$$
corresponds to the full subcategory
$$\on{C}^\cdot(\on{pt}/H)\mod_{0}\subset \on{C}^\cdot(\on{pt}/H)\mod,$$
generated by the augmentation module. 

\sssec{}

The category $\bC^H$ is acted on tautologically by the symmetric monoidal category
$\Vect^H$.

\medskip

Similarly, $\Vect^{H,\on{ren}}$ acquires a symmetric monoidal structure, and $\bC^{H,\on{ren}}$
is a module over it.

\medskip

Moreover, the functor \eqref{e:unren} is symmetric monoidal, and we have
$$\bC^H\simeq \Vect^H\underset{\Vect^{H,\on{ren}}}\otimes \bC^{H,\on{ren}}.$$

\sssec{Example} \label{e:non-ren rep}

Take $\bC=\fh\mod$, so that 
$$\bC^H\simeq \Rep(H).$$

In this case, it is easy to see that the adjoint functors \eqref{e:ren coloc}
are mutually inverse equivalences. 

\sssec{}

The renormalization construction works also for group-schemes. We will apply it in the case when
$\ul\bC$ is a sheaf of categories over $\Ran$, with the group-scheme in question being $\fL^+(H)$
for an algebraic group $H$. We refer the reader to \cite[Sect. 6.11]{CR} for details.

\ssec{The local spherical semi-infinite category}

\sssec{}

Consider the factorization category
$$\on{I}(G,P^-)^{\on{loc}}:=
\Dmod_{\frac{1}{2}}(\Gr_G)^{\fL(N^-_P)\cdot \fL^+(M),\on{ren}},$$
where the renormalization procedure is taken with respect to $\fL^+(M)$, viewed as acting on the category 
$\Dmod_{\frac{1}{2}}(\Gr_G)^{\fL(N^-_P)}$ (the corresponding compact generation property is established
in \cite[Sect. 1.5]{Ga5} for $P=B$ and in \cite{FH} for an arbitrary parabolic). 
 
\medskip

Convolution equips it with a natural action of $\Sph_G$.
(Recall that according to the convention in \cite[Sect. 1.5.4]{GLC2}, we identify 
left and right actions of $\Sph_G$ via the involution $\sigma$.) 

\sssec{} \label{sss:corrected Sph action}

The trivialization of the gerbe in \secref{sss:triv ratio gerbe} gives rise to an action of $\Sph_M$ on $\on{I}(G,P^-)^{\on{loc}}$.
However, as in \secref{sss:naive Jacquet}, this action is \emph{incompatible} with factorization. 

\medskip

We introduce a corrected version of $\Sph_M$-action on $\on{I}(G,P^-)^{\on{loc}}$ by precomposing the one above
with the automorphism of $\Sph_M$ that acts as cohomological shift $[\langle \lambda,2\rhoch_P\rangle]$
on the connected component $\Gr_M^\lambda$. 

\medskip

From now on, unless explicitly mentioned otherwise, we will only consider this corrected version of the 
$\Sph_M$-action on $\on{I}(G,P^-)^{\on{loc}}$. 

\medskip

The resulting $\Sph_M$-action on $\on{I}(G,P^-)^{\on{loc}}$ is compatible with factorization, and commutes with the
$\Sph_G$-action, making $\on{I}(G,P^-)^{\on{loc}}$ into a $(\Sph_G,\Sph_M)$-bimodule. 

\sssec{}

By the same logic, we have a factorization functor 
\begin{multline} \label{e:IGP to M}
\on{I}(G,P^-)^{\on{loc}}=\Dmod_{\frac{1}{2}}(\Gr_G)^{\fL(N^-_P)\cdot \fL^+(M),\on{ren}}
\to  \Dmod_{\frac{1}{2}}(\Gr_{P^-})^{\fL(N^-_P)\cdot \fL^+(M),\on{ren}}\simeq \\
\simeq \Dmod_{\frac{1}{2}}(\Gr_M)^{\fL^+(M),\on{ren}}=:\Sph_M \overset{[\text{shift}]}\to \Sph_M,
\end{multline} 
compatible with the $\Sph_M$-actions, where:

\smallskip

\begin{itemize}

\item The second arrow is $(\sfp^-)^!$;

\medskip

\item The third arrow is the equivalence, induced by 
$$\Dmod_{\frac{1}{2}}(\Gr_M) \overset{(\sfq^-)^!}\simeq \Dmod_{\frac{1}{2}}(\Gr_{P^-})^{\fL(N^-_P)};$$

\item The the last arrow is the cohomological shift $[\langle \lambda,2\rhoch_P\rangle]$
on $\Gr_M^\lambda$;

\end{itemize} 

\medskip

We denote this functor by $\oblv_{\semiinf\to\Sph_M}$. 

\sssec{}

According to \cite[Proposition 1.5.3]{Ga5}, the functor \eqref{e:IGP to M} admits a (factorization) 
left adjoint.\footnote{In \cite{Ga5} only the case $P=B$ is considered. The case of an arbitrary parabolic 
is considered in the forthcoming paper \cite{FH}.}

\medskip

Explicitly, this left adjoint is given by
\begin{multline} \label{e:M to IGP}
\Sph_M\overset{[-\text{shift}]}\to 
\Sph_M:=\Dmod_{\frac{1}{2}}(\Gr_M)^{\fL^+(M),\on{ren}}\overset{(\sfq^-)^!}\simeq \\
\simeq 
\Dmod_{\frac{1}{2}}(\Gr_{P^-})^{\fL(N^-_P)\cdot \fL^+(M),\on{ren}}\overset{(\sfp^-)_!}\longrightarrow \Dmod_{\frac{1}{2}}(\Gr_G)^{\fL(N^-_P)\cdot \fL^+(M),\on{ren}}
=\on{I}(G,P^-)^{\on{loc}}.
\end{multline}

In the above formula, $(\sfp^-)_!$ is the functor, induced by the same-named functor
$$(\sfp^-)_!:\Dmod_{\frac{1}{2}}(\Gr_{P^-})^{\fL(N^-_P)}\to \Dmod_{\frac{1}{2}}(\Gr_G)^{\fL(N^-_P)},$$
left adjoint to $(\sfp^-)^!$.

\medskip

The content of \cite[Proposition 1.5.3]{Ga5} is that the partially defined functor $(\sfp^-)_!$, left adjoint to
$$(\sfp^-)^!:\Dmod_{\frac{1}{2}}(\Gr_G)\to \Dmod_{\frac{1}{2}}(\Gr_{P^-}),$$
is defined on 
$$\Dmod_{\frac{1}{2}}(\Gr_{P^-})^{\fL(N^-_P)}\subset \Dmod_{\frac{1}{2}}(\Gr_{P^-}),$$
and has a natural factorization structure.

\sssec{} \label{sss:Delta expl}

We will denote the above left adjoint by $\ind_{\Sph_M\to\semiinf}$. 

\medskip

Denote
$$\Delta^{-,\semiinf}:=\ind_{\Sph_M\to\semiinf}(\delta_{1,\Gr_M}).$$

Explicitly, $\Delta^{-,\semiinf}$ is given by the !-extension of the dualizing sheaf on the 
$\fL(N^-_P)$-orbit through the origin in $\Gr_G$. (Again, it is the assertion of \cite[Proposition 1.5.3]{Ga5} that 
this !-extension is well-defined over the Ran space and its formation commutes with !-pullbacks.)

\sssec{}

Since $\Sph_M$ is rigid, the functor $\ind_{\Sph_M\to\semiinf}$
is also compatible with the $\Sph_M$-actions.\footnote{Again, we do not distinguish between left and right
actions of $\Sph_M$.} In particular, we have:
$$\ind_{\Sph_M\to\semiinf}(-)\simeq (-)\underset{M}\star \Delta^{-,\semiinf},$$
where $\underset{M}\star$ refers to the $\Sph_M$-action on $\on{I}(G,P^-)^{\on{loc}}$.

\medskip

Since the functor \eqref{e:IGP to M}
is conservative, we obtain a monadic adjunction
\begin{equation} \label{e:IGP to M adj}
\ind_{\Sph_M\to\semiinf}:\Sph_M\rightleftarrows \on{I}(G,P^-)^{\on{loc}}:\oblv_{\semiinf\to\Sph_M}
\end{equation} 
as $\Sph_M$-module categories.

\sssec{} \label{sss:Omega tilde}

Consider the above adjunction \eqref{e:IGP to M adj}. Due to the monadicity, there exists a canonically defined (factorization) associative algebra object 
$$\wt\Omega \in \Sph_M$$
(where $\Sph_M$ is viewed as a monoidal factorization category) so that
\begin{equation} \label{e:IGP as Omega}
\on{I}(G,P^-)^{\on{loc}}\simeq \wt\Omega\mod^r(\Sph_M)
\end{equation} 
as $\Sph_M$-module categories, 
and the adjunction \eqref{e:IGP to M adj} identifies with 
$$\ind_{\wt\Omega}:\Sph_M\rightleftarrows \wt\Omega\mod^r(\Sph_M):\oblv_{\wt\Omega}.$$

\sssec{}

The factorization category $\on{I}(G,P^-)^{\on{loc}}$ carries a natural \emph{unital structure} (see \cite[Sect. C.11]{GLC2} for what this means),
for which the functor $\ind_{\Sph_M\to\semiinf}$ is (strictly) unital.

\medskip

The factorization unit of $\one_{\on{I}(G,P^-)^{\on{loc}}}\in \on{I}(G,P^-)^{\on{loc}}$ is the above object 
$\Delta^{-,\semiinf}$. 

\sssec{}

Note that we also have: 
$$\Delta^{-,\semiinf}\simeq \Av_!^{\fL(N^-_P)}(\delta_{1,\Gr_G}),$$
where $\Av_!^{\fL(N^-_P)}$ is the left adjoint\footnote{This left adjoint is defined by an argument similar to one in \cite[Lemma 2.3.4(2,4)]{Ch1}.} 
to the forgetful functor
$$\Dmod_{\frac{1}{2}}(\Gr_G)^{\fL(N^-_P)\cdot \fL^+(M),\on{ren}}\hookrightarrow \Dmod_{\frac{1}{2}}(\Gr_G)^{\fL^+(M),\on{ren}}.$$

\medskip

Hence, we obtain an adjunction of factorization functors
\begin{equation} \label{e:adj semiinf Sph}
\Av_!^{\fL(N^-_P)}:\Sph_G \rightleftarrows \on{I}(G,P^-)^{\on{loc}}:\Av^{\fL^+(G)/\fL^+(M)}_*.
\end{equation}
with the left adjoint being unital. 

\medskip

The functors in \eqref{e:adj semiinf Sph} are compatible with the $\Sph_G$-actions. Hence, we have
$$\Av_!^{\fL(N^-_P)}(-)\simeq (-)\underset{G}\star \Delta^{-,\semiinf},$$
where $\underset{G}\star$ refers to the $\Sph_G$-action on $\on{I}(G,P^-)^{\on{loc}}$.

\sssec{}

In addition to the functor $\ind_{\Sph_M\to\semiinf}$, there exists a more simply-minded functor, denoted $\ind^*_{\Sph_M\to\semiinf}$, given by 
\begin{multline} \label{e:M to IGP *}
\Sph_M\overset{[-\text{shift}]}\to 
\Sph_M=\Dmod_{\frac{1}{2}}(\Gr_M)^{\fL^+(M),\on{ren}}\overset{(\sfq^-)^!}\simeq \\
\simeq 
\Dmod_{\frac{1}{2}}(\Gr_{P^-})^{\fL(N^-_P)\cdot \fL^+(M),\on{ren}}\overset{(\sfp^-)_*}\longrightarrow 
\Dmod_{\frac{1}{2}}(\Gr_G)^{\fL(N^-_P)\cdot \fL^+(M),\on{ren}}
=\on{I}(G,P^-)^{\on{loc}}.
\end{multline}

I.e., the difference between $\ind^*_{\Sph_M\to\semiinf}$ and $\ind_{\Sph_M\to\semiinf}$ is that we use $(\sfp^-)_*$
instead of $(\sfp^-)_!$.

\medskip

Denote by $\nabla^{-,\semiinf}$ the image of $\delta_{1,\Gr_M}=\one_{\Sph_M}$ under \eqref{e:M to IGP *}. 

\medskip

Note that the functor $\ind^*_{\Sph_M\to\semiinf}$ is also compatible with the actions of $\Sph_M$. Hence, we have: 
$$\ind^*_{\Sph_M\to\semiinf} \simeq (-)\star \nabla^{-,\semiinf}.$$

\medskip

Note also that the composition
$$\oblv_{\semiinf\to\Sph_M}\circ \ind^*_{\Sph_M\to\semiinf}$$
is canonically isomorphic to the identity endofunctor of $\Sph_M$. 

\medskip

In particular,
\begin{equation} \label{e:oblv of nabla}
\oblv_{\semiinf\to\Sph_M}(\nabla^{-,\semiinf})\simeq \delta_{1,\Gr_M}.
\end{equation} 

\sssec{} \label{sss:Omega tilde neg}

Note that the category $\Sph_M$ is naturally graded by
$$\Lambda_M/\Lambda_{[M,M]_{\on{sc}}}\simeq \pi_{1,\on{alg}}(M)$$ 
and in terms of this grading, the algebra $\wt\Omega$ 
is graded by the sub-monoid 
$$\Lambda^{\on{pos}}_{G,P} \subset \Lambda_M/\Lambda_{[M,M]_{\on{sc}}}$$ spanned by the images of the simple coroots,
with the $0$-weight component being $\one_{\Sph_M}$. 

\medskip

In particular, the algebra $\wt\Omega$ admits a unique augmentation, and hence we have a well-defined object
\begin{equation} \label{e:aug Sph}
\delta_{1,\Gr_M}=\one_{\Sph_M}\in \wt\Omega\mod^r(\Sph_M).
\end{equation} 

\medskip

It follows from \eqref{e:oblv of nabla} that under the equivalence \eqref{e:IGP as Omega}
the object $\nabla^{-,\semiinf}\in \on{I}(G,P^-)^{\on{loc}}$ corresponds to \eqref{e:aug Sph}. 

\sssec{}

We now recall that in addition to the factorization algebras 
$$\Delta^{-,\semiinf} \text{ and } \nabla^{-,\semiinf},$$
there exists another factorization algebra in $\on{I}(G,P^-)^{\on{loc}}$, denoted 
$$\IC^{-,\semiinf},$$
constructed in \cite{Ga4,Ga5}.\footnote{In {\it loc.cit.} this object was defined for $P=B$. The case of a general parabolic
is similar, and will be considered in detail in \cite{FH}.}

\medskip

It is unital and is equipped with homomorphisms of unital factorization algebras 
$$\Delta^{-,\semiinf} \to \IC^{-,\semiinf}\to \nabla^{-,\semiinf}.$$

\ssec{Digression: \texorpdfstring{$\Sph_G$}{Sph} as a rigid category}

\sssec{}

Let $\bA$ be a rigid monoidal category (see \cite[Chapter 1, Definition 9.1.2]{GaRo3} for what this means). Let 
$\on{mult}:\bA\otimes \bA\to \bA$ denote the binary operation, and let 
$\on{mult}^R$ denote its right adjoint. Let $\one_\bA\in \bA$ denote the monoidal unit. 

\medskip

Recall  (see \cite[Chapter 1, Sect. 9.2.1]{GaRo3}) that the object
$$\on{mult}^R(\one_\bA)\in \bA\otimes \bA$$
provides the unit for a self-duality on $\bA$.

\sssec{}

Note, however, that the are \emph{two} such self-dualities: 

\medskip

On the one hand, we have the identification
$$\bA \overset{\psi^l}\simeq \bA^\vee,$$
for which 
$$(\psi^l\otimes \on{id})(\on{mult}^R(\one_\bA))\in \bA^\vee\otimes \bA$$
is the unit.

\medskip

On the other hand, we have the identification
$$\bA \overset{\psi^r}\simeq \bA^\vee,$$
for which 
$$(\on{id}\otimes \psi^r)(\on{mult}^R(\one_\bA))\in \bA\otimes \bA^\vee$$
is the unit.

\sssec{} \label{sss:pivotal}

Both $\psi^l$ and $\psi^r$ are monoidal equivalences, where we equip $\bA^\vee$ with a monoidal structure
dual to the comonoidal structure on $\bA$, where the latter is given by passing to right adjoints from the
original monoidal structure.  

\medskip

Let $\psi^{l\to r}$ be the endofunctor of $\bA$, given by
$$(\psi^r)^{-1}\circ \psi^l.$$

It is a monoidal self-equivalence of $\bA$. 

\medskip

A \emph{pivotal structure} on $\bA$ is a trivialization of $\psi^{l\to r}$
(as a monoidal self-equivalence). 

\medskip

Given a pivotal structure on $\bA$, we will use notation
$$\psi^l=:\psi:=\psi^r.$$

\sssec{Example}
Suppose that $\bA$ is compactly generated. Then $\psi^l$ (resp., $\psi^r$) 
is given on compact objects by the operation of left (resp., right) dual.

\medskip

Hence, $\psi^{l\to r}$ is given on compact objects as the square of functor
of left duality.

\sssec{} \label{sss:pivot modules}

Let $\bM^l$ (resp., $\bM^r$) be a left (resp., right) $\bA$-module category. Consider the 
action functors
$$\bA\otimes \bM^l \overset{\on{act}}\to \bM^l \text{ and } \bM^r\otimes \bA\overset{\on{act}}\to \bM^r.$$

Consider their duals 
$$\bM^l\overset{\on{act}^\vee}\to \bA^\vee\otimes \bM^l \text{ and } \bM^r\overset{\on{act}^\vee}\to \bM^r\otimes \bA^\vee.$$

\medskip

We have the commutative diagrams
$$
\CD
\bM^l @>{\on{act}^\vee}>> \bA^\vee\otimes \bM^l  \\
@A{\on{Id}}AA @A{\sim}A{\psi^l\otimes \on{Id}}A \\
\bM^l @>{\on{act}^R}>> \bA \otimes \bM^l 
\endCD
$$
and 
$$
\CD
\bM^r @>{\on{act}^\vee}>>  \bM^r\otimes \bA^\vee \\
@A{\on{Id}}AA @A{\sim}A{\on{Id}\otimes \psi^r}A \\
\bM^r @>{\on{act}^R}>>  \bM^r\otimes \bA.
\endCD
$$

In particular, we have a commutative diagram
$$
\CD
\bM^r @>{\on{act}^\vee}>>  \bM^r\otimes \bA^\vee \\
@A{\on{Id}}AA @A{\sim}A{\on{Id}\otimes \psi^l}A \\
\bM^r_{\psi^{l\to r}} @>{\on{act}^R}>>  \bM^r_{\psi^{l\to r}}\otimes \bA,
\endCD
$$
where $\bM^r_{\psi^{l\to r}}$ the right $\bA$-module category, obtained by precomposing 
the action functor with the automorphism $\psi^{l\to r}$ of $\bA$.

\medskip

From here it follows that we have a canonical identification

\begin{equation} \label{e:tensor vs cotensor}
\bM^l\overset{\bA}\otimes \bM^r \simeq \bM^r_{\psi^{l\to r}}\underset{\bA}\otimes \bM^l,
\end{equation}
where
$$\bM^l\overset{\bA}\otimes \bM^r:=\on{Funct}_{\bA\otimes \bA^o}(\bA,\bM^l\otimes \bM^r),$$
see\footnote{The notation $\bA\rightsquigarrow \bA^o$ means reversing the monoidal structure.} 
\cite[Chapter 1, Proposition 9.5.6]{GaRo3}. 

\medskip

In particular, we obtain that if $\bA$ is equipped with a pivotal structure, we have a canonical identification
\begin{equation} \label{e:tensor vs cotensor level}
\bM^r\overset{\bA}\otimes \bM^l \simeq \bM^r\underset{\bA}\otimes \bM^l.
\end{equation}

\medskip

Further, we obtain that if $\bA$ is equipped with a pivotal structure, and $\bM^l$ (resp., $\bM^r$) is dualizable as a plain DG category, then its 
dualizable as an $\bA$-module category, and its $\bA$-module dual identifies with $(\bM^l)^\vee$ (resp., $(\bM^r)^\vee$)
as a right (resp., left) $\bA$-module category. 

\sssec{}

The above discussion applies verbatim to the case when $\bA$ is a monoidal factorization category.

\sssec{}

We take $\bA=\Sph_G$. We claim that $\bA$ is equipped with a naturally defined pivotal structure.

\medskip

Indeed, in this case we have a well-defined \emph{Verdier duality} identification
\begin{equation} \label{e:Verdier Sph}
\Sph_G \simeq \Sph^\vee_G
\end{equation}
as monoidal factorization categories (see \secref{sss:pivotal}), and $\psi^l$ (resp., $\psi^r$)
is obtained from \eqref{e:Verdier Sph} by composing (resp., precomposing) with the anti-involution
$\sigma$ of $\Sph_G$.

\medskip

The identification $\psi^l\simeq \psi^r$ follows now from the fact that $\sigma$ commutes
with Verdier duality on $\Sph_G$.

\sssec{} \label{sss:duality over Sph} 

Thus, thanks to \secref{sss:pivot modules}, we will freely identity duals of $\Sph_G$-modules 
as plain categories with their duals as $\Sph_G$-modules.

\medskip

Furthermore, given a left and right module categories $\bM^l$ and $\bM^r$ over $\Sph_G$ we will freely identify
\begin{equation} \label{e:tensor vs cotensor Sph}
\bM^r\overset{\Sph_G}\otimes \bM^l \simeq \bM^r\underset{\Sph_G}\otimes \bM^l.
\end{equation}

\sssec{}

We now consider the case of $\bA=\Sph_\cG^{\on{spec}}$. 

\medskip

On the one hand, the geometric Satake equivalence
$$\on{Sat}_G:\Sph_G\simeq \Sph_\cG^{\on{spec}}$$
allows us to transport the pivotal structure from $\Sph_G$ to $\Sph_\cG^{\on{spec}}$

\medskip

However, on the other hand, one would like to understand the geometric meaning of this 
pivotal structure on $\Sph_\cG^{\on{spec}}$ purely spectral side. This will be done in a future
publication, and in the meantime we will avoid using it in main the narrative in this paper. This 
means that some categories that could ultimately be identified will \emph{not} be identified
in this paper. 

\medskip

Thus, given a left and right module categories $\bM^l$ and $\bM^r$ over $\Sph_G$, we have 
\begin{equation} \label{e:tensor vs cotensor Sph spec}
\bM^r\overset{\Sph^{\on{spec}}_\cG}\otimes \bM^l \simeq \bM^r_{\psi_G^{l\to r}}\underset{\Sph^{\on{spec}}_\cG}\otimes \bM^l,
\end{equation}
where $\psi_G^{l\to r}$ is the corresponding monoidal self-equivalence of $\Sph_\cG^{\on{spec}}$. 

\sssec{An observation} \label{sss:tensor vs cotensor}

Let $B$ be an coassociative coalgebra in a symmetric monoidal category $\bO$. 
and let $M^l$ and $M^r$ be left and right $B$-comodules in $\bO$, respectively. In this case
we can consider the object
$$M^r\overset{B}\otimes M^l\in \bO.$$

Let now $A$ be an associative algebra in $\bO$, dualizable as an object. Let $M^l$ and $M^r$ be 
left and right $A$-modules in $\bO$, respectively. Set $B:=A^\vee$. Passing to the duals 
of the action maps, we can regard $M^l$ and $M^r$ as $B$-comodules.

\medskip

Unwinding the definitions, we obtain a canonical isomorphism
$$M^r\overset{B}\otimes M^l\simeq M^l\overset{A}\otimes M^r.$$

We will mostly apply this for $\bO=\DGCat$. 

\ssec{The dual of the spherical semi-infinite category}

\sssec{} 

Consider the category 
$$\Dmod(\Gr_G)_{\fL(N^-_P)}.$$

It is proved in \cite{Ch1} that it is dualizable as a factorization category. Once the dualizability is established,
it follows formally that Verdier duality on $\Gr_G$ gives rise to an identification
$$(\Dmod(\Gr_G)^{\fL(N^-_P)})^\vee \simeq \Dmod(\Gr_G)_{\fL(N^-_P)}.$$

\sssec{} \label{sss:IGP co}

Denote
$$\on{I}(G,P^-)^{\on{loc}}_{\on{co}}:=\Dmod_{\frac{1}{2}}(\Gr_G)_{\fL(N^-_P)}^{\fL^+(M),\on{ren}}.$$

By similar logic, we have a canonical identification
\begin{equation} \label{e:dual as co}
\left(\on{I}(G,P^-)^{\on{loc}}\right)^\vee \simeq \on{I}(G,P^-)^{\on{loc}}_{\on{co}}.
\end{equation} 

\sssec{} \label{sss:which action on dual}

Note that the identification \eqref{e:dual as co} is compatible with the $\Sph_M$-actions in the 
following sense:

\begin{itemize}

\item The action on $\left(\on{I}(G,P^-)^{\on{loc}}\right)^\vee$ is induced by the 
$\Sph_M$-action on $\on{I}(G,P^-)^{\on{loc}}$ specified in \secref{sss:corrected Sph action}
(as always, we pass from right to left $\Sph_M$-modules using $\sigma$). 

\medskip

\item The action on $\on{I}(G,P^-)^{\on{loc}}_{\on{co}}$ is obtained from the natural geometric
action, by applying the inverse cohomological shift to the one from \secref{sss:corrected Sph action}.

\end{itemize}

\medskip

Since $\Sph_M$ is rigid, the duality \eqref{e:dual as co} realizes $\on{I}(G,P^-)^{\on{loc}}_{\on{co}}$
as a $\Sph_M$-module dual of $\on{I}(G,P^-)^{\on{loc}}$. 

\sssec{}

Note that the functor dual to 
$$\oblv_{\semiinf\to\Sph_M}:\on{I}(G,P^-)^{\on{loc}}\to \Sph_M$$
is a functor, to be denoted 
\begin{equation} \label{e:M to IGP co}
\ind_{\Sph_M\to\semiinf,\on{co}}:\Sph_M\to \on{I}(G,P^-)^{\on{loc}}_{\on{co}},
\end{equation} 
given by
\begin{multline} \label{e:M to IGP co expl}
\Sph_M\overset{[\text{shift}]}\to 
\Sph_M=\Dmod_{\frac{1}{2}}(\Gr_M)^{\fL^+(M),\on{ren}}\simeq \\
\simeq 
\Dmod_{\frac{1}{2}}(\Gr_{P^-})^{\fL^+(M),\on{ren}}_{\fL(N^-_P)}\overset{(\sfp^-)_*}\longrightarrow \Dmod_{\frac{1}{2}}(\Gr_G)_{\fL(N^-_P)}^{\fL^+(M),\on{ren}}
=\on{I}(G,P^-)^{\on{loc}}_{\on{co}},
\end{multline}
where:

\begin{itemize} 

\item The first arrow is the cohomological shift $[\langle \lambda,2\rhoch_P\rangle]$
on $\Gr_M^\lambda$;

\medskip

\item The third arrow is induced by the equivalence
$$\Dmod_{\frac{1}{2}}(\Gr_{P^-})_{\fL(N^-_P)}\overset{(\sfq^-)_*}\longrightarrow \Dmod_{\frac{1}{2}}(\Gr_M).$$

\end{itemize} 

\sssec{}

The factorization category $\on{I}(G,P^-)^{\on{loc}}_{\on{co}}$ has a natural unital structure, and the functor 
$\ind_{\Sph_M\to\semiinf,\on{co}}$ is (strictly) unital.

\medskip

The factorization unit in $\on{I}(G,P^-)^{\on{loc}}_{\on{co}}$ is the object, to be denoted $\Delta^{-,\semiinf}_{\on{co}}$, 
equal to the image of $\delta_{1,\Gr}\in \Dmod_{\frac{1}{2}}(\Gr_G)^{\fL^+(P^-),\on{ren}}$ under the natural projection
$$\Dmod_{\frac{1}{2}}(\Gr_G)^{\fL^+(P^-),\on{ren}}\to \Dmod_{\frac{1}{2}}(\Gr_G)_{\fL(N^-_P)}^{\fL^+(M),\on{ren}}=
\Delta^{-,\semiinf}_{\on{co}}.$$   

\sssec{}

By duality, we obtain that the functor \eqref{e:M to IGP co} admits a right adjoint, 
to be denoted $\oblv_{\semiinf\to\Sph_M,\on{co}}$, namely, the dual of $\ind_{\Sph_M\to\semiinf,}$, 
so that we have a monadic adjunction
\begin{equation} \label{e:IGP to M adj co}
\ind_{\Sph_M\to\semiinf,\on{co}}:\Sph_M\rightleftarrows \on{I}(G,P^-)^{\on{loc}}_{\on{co}}:\oblv_{\semiinf\to\Sph_M,\on{co}}.
\end{equation} 

\sssec{}

Note also that we can also consider the functor
\begin{multline} \label{e:M to IGP co naive}
\on{I}(G,P^-)^{\on{loc}}_{\on{co}}=\Dmod_{\frac{1}{2}}(\Gr_G)_{\fL(N^-_P)}^{\fL^+(M),\on{ren}} \overset{(\sfp^-)^!}\longrightarrow
\Dmod_{\frac{1}{2}}(\Gr_{P^-})_{\fL(N^-_P)}^{\fL^+(M),\on{ren}}\simeq \\
\simeq \Dmod_{\frac{1}{2}}(\Gr_M)^{\fL^+(M),\on{ren}}=\Sph_M \overset{[\text{shift}]}\simeq \Sph_M.
\end{multline} 

Denote this functor by $\oblv^*_{\semiinf\to\Sph_M,\on{co}}$. By construction, we have
$$\oblv^*_{\semiinf\to\Sph_M,\on{co}}\simeq (\ind^*_{\Sph_M\to\semiinf})^\vee.$$

\sssec{}

Let $\bA$ be a monoidal category and $A\in \bA$ an associative algebra. We consider
$A\mod(\bA)$ as a right $\bA$-module category, and $A\mod^r(\bA)$ as a left $\bA$-module category.

\medskip

Tautologically, $A\mod(\bA)$ and $A\mod^r(\bA)$ are each other's duals (as right and left $\bA$-module categories,
respectively). 

\medskip

Note also that 
$$A\mod^r(\bA)\simeq A^{o}\mod(\bA^{o}),$$
where:

\begin{itemize}

\item For a monoidal category $\bA$, we denoted by $\bA^{o}$ denotes the monoidal category category obtained 
by reversing the monoidal operation;

\medskip

\item For an associative algebra $A$ in a monoidal category $\bA$,
we denote by $A^{o}$ the corresponding associative algebra
in $\bA^{o}$. 

\end{itemize}

\sssec{} \label{sss:algebra for dual category}

Let $\wt\Omega_{\on{co}}$ be the associative (factorization) algebra object in $\Sph_M$,
so that \eqref{e:IGP to M adj co} identifies with
$$\ind_{\wt\Omega_{\on{co}}}:\Sph_M\rightleftarrows \wt\Omega_{\on{co}}\mod^r(\Sph_M):\oblv_{\wt\Omega_{\on{co}}}.$$

\medskip

From the duality between
$$\on{I}(G,P^-)^{\on{loc}} \text{ and } \on{I}(G,P^-)^{\on{loc}}_{\on{co}}$$
as $\Sph_M$-module categories, we obtain an identification
\begin{equation} \label{e:op and sigma}
\wt\Omega_{\on{co}} \simeq \sigma((\wt\Omega)^{o}).
\end{equation}

\sssec{}

We now quote the following result of \cite{Ch1}:

\begin{thm} \label{t:geom intertwiner}
The composite functor 
$$\Dmod(\Gr_G)^{\fL(N(P))} \hookrightarrow 
\Dmod(\Gr_G)\twoheadrightarrow \Dmod(\Gr_G)_{\fL(N^-_P)}$$
is an equivalence (as factorization categories).
\end{thm}

\medskip

From \thmref{t:geom intertwiner} we formally obtain:

\begin{cor}  \label{c:geom intertwiner}
The functor 
\begin{multline} \label{e:long intertwiner}
\on{I}(G,P)^{\on{loc}}:=\Dmod_{\frac{1}{2}}(\Gr_G)^{\fL(N(P))\cdot \fL^+(M),\on{ren}} \hookrightarrow \\
\hookrightarrow \Dmod_{\frac{1}{2}}(\Gr_G)^{\fL^+(M),\on{ren}}\twoheadrightarrow 
\Dmod_{\frac{1}{2}}(\Gr_G)_{\fL(N^-_P)}^{\fL^+(M),\on{ren}}=\on{I}(G,P^-)^{\on{loc}}_{\on{co}}
\end{multline}
is an equivalence (as factorization categories).
\end{cor}

In what follows we will denote by $\Theta_{\on{I}(G,P^-)^{\on{loc}}}$ the functor inverse to the equivalence of \eqref{e:long intertwiner}. 

\sssec{}

Applying the Chevalley involution on $G$ (normalized so that it swaps $P$ and $P^-$, and thus is compatible
with the Chevalley involution $\tau_M$ on $M$), we obtain an equivalence
$$\on{I}(G,P)^{\on{loc}} \overset{\tau_G}\simeq \on{I}(G,P^-)^{\on{loc}}.$$

Thus, composing \eqref{e:long intertwiner} with \eqref{e:dual as co} and $\tau_G$, we obtain a self-duality
\begin{equation} \label{e:self-duality IGP}
(\on{I}(G,P^-)^{\on{loc}})^\vee \simeq \on{I}(G,P^-)^{\on{loc}}. 
\end{equation}

By construction, the equivalence \eqref{e:self-duality IGP} is compatible with the actions of $\Sph_M$, where:

\begin{itemize}

\item The action on $(\on{I}(G,P^-)^{\on{loc}})^\vee$ is the one specified in \secref{sss:which action on dual};

\medskip

\item The action on $\on{I}(G,P^-)^{\on{loc}}$ is precomposition of the action 
specified in \secref{sss:corrected Sph action} with the automorphism $\tau_M$ of $\Sph_M$. 

\end{itemize}

\sssec{}

As in \secref{sss:algebra for dual category}, from the equivalence \eqref{e:self-duality IGP} we obtain an isomorphism 
\begin{equation} \label{e:Omega sigma and tau}
\tau_M(\wt\Omega) \simeq \sigma((\wt\Omega)^{o}), 
\end{equation} 
as associative (factorization) algebras in $\Sph_M$. 

\ssec{A twisting construction: recollections}  \label{ss:twist by G-bundle}

In this and the next subsections we will (re)introduce the translation by $\rho(\omega_X)$,
mentioned in the preamble to this section. 

%

\medskip

The reader may choose to skip this material on the first pass and returned to when necessary. 

\sssec{} \label{sss:twist by M-bundle}

Recall the twisting construction from \cite[Sect. 1.2]{GLC2}. Let $\CP_M$ be an $M$-torsor on $X$. We can consider $\CP_M$-twisted versions of all objects in sight, i.e.,
$$\Gr_{M,\CP_M},\,\, \Gr_{G,\CP_M},\,\, \fL(N^-_P)_{\CP_M},$$
see \cite[Sect. 1.2]{GLC2}.

\medskip

We will denote by subscript $\CP_M$ the categories associated with the corresponding twisted geometric
objects, i.e.,
$$\Sph_{M,\CP_M},\,\, \Sph_{G,\CP_M},\,\, \on{I}(G,P^-)^{\on{loc}}_{\CP_M}.$$

In particular, we have a monadic adjunction 
$$\ind_{\Sph_M\to\semiinf}:\Sph_{M,\CP_M}\rightleftarrows \on{I}(G,P^-)^{\on{loc}}_{\CP_M}:\oblv_{\semiinf\to\Sph_M},$$
and the corresponding associative (factorization) algebra object
$$\wt\Omega_{\CP_M}\in \Sph_{M,\CP_M}.$$

We will denote by the same symbol $\Delta^{-,\semiinf}$ the factorization unit in $\on{I}(G,P^-)^{\on{loc}}_{\CP_M}$,
and by
$$\IC^{-,\semiinf}\in \on{I}(G,P^-)^{\on{loc}}_{\CP_M}$$
the corresponding ``semi-infinite IC sheaf". 

\sssec{} \label{sss:remove twist by M-bundle}

Note, however, that the local Hecke stacks for $M$ (or $G$), i.e.,
$$\on{Hecke}^{\on{loc}}_M:=\fL^+(M)\backslash \fL(M)/\fL^+(M) \text{ and } \on{Hecke}^{\on{loc}}_G:=\fL^+(G)\backslash \fL(G)/\fL^+(G)$$ 
are canonically isomorphic to their twisted versions, see \cite[Equation (1.3)]{GLC2}.

\medskip

So we have canonical identifications of monoidal (factorization) categories
$$\Sph_M\overset{\alpha_{\CP_M,\on{taut}}}\simeq \Sph_{M,\CP_M}\text{ and } 
\Sph_G\overset{\alpha_{\CP_M,\on{taut}}}\simeq \Sph_{G,\CP_M}.$$

Similarly, we have a canonical equivalence
\begin{equation} \label{e:IGP local no twist}
\alpha_{\CP_M,\on{taut}}:\on{I}(G,P^-)^{\on{loc}}\simeq \on{I}(G,P^-)^{\on{loc}}_{\CP_M}.
\end{equation} 

Tautologically, we have:
$$\alpha_{\CP_M,\on{taut}}(\wt\Omega)\simeq \wt\Omega_{\CP_M}$$
and 
$$\alpha_{\CP_M,\on{taut}}(\Delta^{-,\semiinf})\simeq \Delta^{-,\semiinf}.$$

\sssec{} \label{sss:IGP co twist}

As in \secref{sss:IGP co}, we can consider
the category $\on{I}(G,P^-)^{\on{loc}}_{\on{co},\CP_M}$. 

\medskip

We still have an equivalence
\begin{equation} \label{e:dual as co twisted}
\left(\on{I}(G,P^-)^{\on{loc}}_{\CP_M}\right)^\vee \simeq \on{I}(G,P^-)^{\on{loc}}_{\on{co},\CP_M},
\end{equation} 
compatible with $\Sph_M$-actions.

\medskip

Similarly, we can consider the associative (factorization) algebra object 
$$\wt\Omega_{\on{co},\CP_M}\in \Sph_M,$$
and as in \eqref{e:op and sigma} we have
\begin{equation} \label{e:op and sigma twisted}
\wt\Omega_{\on{co},\CP_M} \simeq \sigma((\wt\Omega_{\CP_M})^o).
\end{equation}

\sssec{} \label{sss:twisted self-duality IGP}

The assertion of \thmref{t:geom intertwiner} renders automatically to the twisted context, so that the functor 
\begin{equation} \label{e:long intertwiner twisted}
\on{I}(G,P)^{\on{loc}}_{\CP_M}\to \on{I}(G,P^-)^{\on{loc}}_{\on{co},\CP_M}
\end{equation} 
is an equivalence. 

\medskip

The Chevalley involution induces an equivalence
$$\on{I}(G,P)^{\on{loc}}_{\CP_M} \overset{\tau_G}\simeq \on{I}(G,P^-)^{\on{loc}}_{\tau_M(\CP_M)}.$$

One can compose it with the equivalence
$$\on{I}(G,P^-)^{\on{loc}}_{\tau_M(\CP_M)} \overset{(\alpha_{\CP_M,\on{taut}})\circ (\alpha_{\tau_M(\CP_M),\on{taut}})^{-1}}
\simeq \on{I}(G,P)^{\on{loc}}_{\CP_M},$$
and thus obtain again a self-duality 
\begin{equation} \label{e:self-duality IGP twisted}
(\on{I}(G,P^-)^{\on{loc}}_{\CP_M})^\vee \simeq  \on{I}(G,P^-)^{\on{loc}}_{\CP_M}.
\end{equation} 

\sssec{}

Similarly to \eqref{e:Omega sigma and tau}, the equivalence \eqref{e:self-duality IGP twisted} yields an isomorphism of associative 
(factorization) algebras in $\Sph_M$:
\begin{equation} \label{e:Omega sigma and tau twisted}
(\alpha_{\tau_M(\CP_M),\on{taut}})^{-1}(\tau_M(\wt\Omega_{\CP_M}))\simeq (\alpha_{\CP_M,\on{taut}})^{-1} \circ \sigma((\wt\Omega_{\CP_M})^o).
\end{equation} 

\ssec{Twisting with respect to the center}  \label{ss:central twist}

This subsection is a continuation of the preceding one. Here we will assume that the $M$-bundle 
$\CP_M$ is induced from a $Z_M$-bundle $\CP_{Z_M}$. 

\medskip

By contrast with the previous section, in which ``not much was going on", here there is something to pay
attention to: the non-triviality of the automorphism \eqref{e:autom Sph} of $\Sph_M$. 

\medskip

That said, this subsection can also be skipped on the first pass, and returned to when necessary. 

\sssec{} \label{sss:twist by central bundle} 

Assume now that in the context of \cite[Sect. 1.2]{GLC2}, the group  
$H$ is abelian and the action of $\fL^+(H)$ on $\CY$ is trivial. In this case, we have canonical isomorphisms
$$\CY_{\CP_H}\simeq \CY \text{ and } \fL(H)_{\CP_H}\simeq \fL(H).$$

In particular, we obtain an equivalence
$$\alpha_{\CP_H,\on{cent}}: \Dmod(\CY)_{\CP_H}\overset{\sim}\to \Dmod(\CY)$$
and hence \emph{another} identification
$$\alpha_{\CP_H,\on{cent}}: (\Dmod(\CY)^{\fL^+(H)})_{\CP_H}\overset{\sim}\to \Dmod(\CY)^{\fL^+(H)}.$$

We will denote by 
$$(\on{transl}_{\CP_{Z_H}})^*:= \alpha_{\CP_H,\on{cent}}\circ \alpha_{\CP_H,\on{taut}}$$
the resulting auto-equivalence of $\Dmod(\CY)^{\fL^+(H)}$. 

\medskip

The functor $(\on{transl}_{\CP_{Z_H}})^*$ is the pullback along the automorphism, denoted $\on{transl}_{\CP_{Z_H}}$,
of the stack $\CY/\fL^+(H)$ given by translation by the point $\CP_H\in \on{pt}/\fL^+(H)$ via the action map
$$\on{pt}/\fL^+(H)\times \CY/\fL^+(H)\to \CY/\fL^+(H).$$

\sssec{} \label{sss:twist by central bundle Gr}

A typical example of the situation of \secref{sss:twist by central bundle} is when 
$\CY=\Gr_M$, so that $\CY/\fL^+(M)$ is the local Hecke stack
$$\on{Hecke}^{\on{loc}}_M:=\fL^+(M)\backslash \fL(M)/\fL^+(M).$$

\medskip

Let $H$ map to the center of $M$. In this case, the action of $\fL^+(H)$ on $\Gr_M$ is trivial.

\medskip

The above automorphism $\on{transl}_{\CP_{Z_H}}$ of 
$$\on{Hecke}^{\on{loc}}_M=\{\CP'_M,\CP''_M, \CP'_M|_{\cD^\times}\sim \CP''_M|_{\cD^\times}\}$$
is given by tensoring the $M$-bundles involved by $\CP_H$, using the canonical map
$$\on{pt}/H\times \on{pt}/M\to \on{pt}/M.$$

\sssec{} \label{sss:twist by central M-bundle}

Assume now that in the context of \secref{sss:twist by M-bundle}, the $M$-torsor $\CP_M$ is induced by a $Z_M$-torsor $\CP_{Z_M}$. Hence,
we obtain \emph{another} identification 
\begin{equation} \label{e:alpha on Sph}
\Sph_{M,\CP_{Z_M}}\overset{\alpha_{\CP_{Z_M},\on{cent}}}\simeq \Sph_M
\end{equation} 
as monoidal (factorization) categories.

\medskip

The composite 
\begin{equation} \label{e:autom Sph}
(\on{transl}_{\CP_{Z_M}})^*:=\alpha_{\CP_{Z_M},\on{cent}} \circ \alpha_{\CP_{Z_M},\on{taut}}
\end{equation}
is a monoidal (factorization) automorphism of $\Sph_M$.

\begin{rem} \label{r:non-trivial transl}

The automorphism $(\on{transl}_{\CP_{Z_M}})^*$ of $\Sph_M$ is \emph{non-trivial}. However, its composition with 
the naive Satake functor
$$\on{Sat}_M^{\on{nv}}:\Rep(\cM)\to \Sph_M.$$
\emph{is canonically} trivial. 

\end{rem} 

\sssec{Convention} \label{sss:convent Sph-action on IGP}

Henceforth, unless explicitly specified otherwise, when we consider $\on{I}(G,P^-)^{\on{loc}}_{\CP_{Z_M}}$ as acted on by $\Sph_M$,
we will do so using the identification $\alpha_{\CP_{Z_M},\on{cent}}$ of \eqref{e:alpha on Sph}. 

\medskip

In particular, we will view $\wt\Omega_{\CP_{Z_M}}$ as an associative (factorization) algebra object
in $\Sph_M$ via $\alpha_{\CP_{Z_M},\on{cent}}$.

%
%
%
%
%

\medskip

Tautologically, we have
$$\wt\Omega_{\CP_{Z_M}} \simeq (\on{transl}_{\CP_{Z_M}})^*(\wt\Omega),$$
as associative (factorization) algebras in $\Sph_M$. 

\sssec{}  \label{sss:self-duality twisted geom cen}

Note that the equivalence 
\begin{equation} \label{e:self-duality IGP twisted center}
(\on{I}(G,P^-)^{\on{loc}}_{\CP_{Z_M}})^\vee \simeq  \on{I}(G,P^-)^{\on{loc}}_{\CP_{Z_M}}
\end{equation} 
of \eqref{e:self-duality IGP twisted} is compatible with the $\Sph_M$-actions, where

\begin{itemize}

\item The action on $\left(\on{I}(G,P^-)^{\on{loc}}_{\CP_{Z_M}}\right)^\vee$ is induced by the 
$\Sph_M$-action on $\on{I}(G,P^-)^{\on{loc}}_{\CP_{Z_M}}$ specified in \secref{sss:convent Sph-action on IGP}
(as always, we pass from right to left $\Sph_M$-modules using $\sigma$);

\medskip

\item The action on $\on{I}(G,P^-)^{\on{loc}}_{\CP_{Z_M}}$ is precomposition of the action 
specified in \secref{sss:corrected Sph action} with the automorphisms $\tau_M$ and 
$(\on{transl}_{\CP_{Z_M}\otimes \tau(\CP^{\otimes -1}_{Z_M})})_*$ of $\Sph_M$. 

\end{itemize} 

\sssec{}

The equivalence \eqref{e:self-duality IGP twisted center} translates into the following isomorphism of 
associative (factorization) algebras in $\Sph_M$:

\begin{equation} \label{e:Omega sigma and tau twisted center}
\tau_M(\wt\Omega_{\CP_{Z_M}}) \simeq 
(\on{transl}_{\CP^{\otimes -1}_{Z_M}\otimes \tau(\CP_{Z_M})})^*
\left(\sigma((\wt\Omega_{\CP_{Z_M}})^{o})\right).
\end{equation} 

Note, however, that \eqref{e:Omega sigma and tau twisted center} can be equivalently obtained by applying
$(\on{transl}_{\CP^{\otimes -1}_{Z_M}})^*$ to \eqref{e:Omega sigma and tau}. 

\sssec{}

In practice we will take 
$$\CP_{Z_M}:=\rho_P(\omega_X):=2\rho_P(\omega^{\otimes \frac{1}{2}}_X).$$

\sssec{}

Note that
$$\tau(\rho_P(\omega_X))=-\rho_P(\omega_X).$$

Hence, in this case, \eqref{e:Omega sigma and tau twisted center} specializes to 
\begin{equation} \label{e:Omega sigma and tau rho}
\tau_M(\wt\Omega_{\rho_P(\omega_X)}) \simeq 
(\on{transl}_{-2\rho_P(\omega_X)})^*\left(\sigma((\wt\Omega_{\rho_P(\omega_X)})^{o})\right).
\end{equation} 

\ssec{The enhanced Jacquet functor} \label{ss:enh Jacquet} 

In this subsection we introduce the \emph{raison d'\^etre} of $\on{I}(G,P^-)^{\on{loc}}$, namely, how it can be used 
to define generalized Jacquet functors. 

\sssec{}
 
Note that we can identify 
$$\Dmod(\Gr_G)^{\fL(N^-_P)\cdot \fL^+(M)}  \simeq  \Bigl(\Dmod(\Gr_G)\otimes \Dmod(\Gr_M)\Bigr)^{\fL(P^-)}.$$

This observation allows us to view objects of $\Dmod(\Gr_G)^{\fL(N^-_P)\cdot \fL^+(M)}$ 
\emph{kernels of functors} between $\Dmod(\Gr_G)$ and $\Dmod(\Gr_M)$.

\medskip

Namely, we consider the factorization functor
\begin{multline} \label{e:gen Jacquet functor non-ren}
\Dmod(\Gr_G)\otimes \Dmod(\Gr_G)^{\fL(N^-_P)\cdot \fL^+(M)} \simeq \\
\simeq \Dmod(\Gr_G)\otimes \Bigl(\Dmod(\Gr_G)\otimes \Dmod(\Gr_M)\Bigr)^{\fL(P^-)}\to \\
\to \Dmod(\Gr_G)\otimes \Dmod(\Gr_G)\otimes \Dmod(\Gr_M) \overset{(\sotimes) \otimes (\on{Id})}\longrightarrow 
\Dmod(\Gr_G)\otimes \Dmod(\Gr_M) 
\overset{\Gamma_\dr(\Gr_G,-)\otimes \on{Id}}\longrightarrow \\
\to \Dmod(\Gr_M).
\end{multline}

This functor is equivariant with respect to the $\fL(P^-)$-actions on $\Dmod(\Gr_G)$ (via $\fL(P^-)\to \fL(G)$)
and on $\Dmod(\Gr_M)$ (via $\fL(P^-)\to \fL(M)$), respectively. 

\sssec{}

The functor \eqref{e:gen Jacquet functor non-ren} has a variant, given by 
\begin{multline} \label{e:gen Jacquet functor}
\Dmod_{\frac{1}{2}}(\Gr_G)\otimes \on{I}(G,P^-)^{\on{loc}}\to 
\Dmod_{\frac{1}{2}}(\Gr_G)\otimes \Dmod_{\frac{1}{2}}(\Gr_G)^{\fL(N^-_P)\cdot \fL^+(M)} \simeq \\
\simeq \Dmod_{\frac{1}{2}}(\Gr_G)\otimes \left(\Dmod_{\frac{1}{2}}(\Gr_G)\otimes \Dmod_{\frac{1}{2}}(\Gr_M)\right)^{\fL(P^-)}\to \\
\to \Dmod_{\frac{1}{2}}(\Gr_G)\otimes \Dmod_{\frac{1}{2}}(\Gr_G)\otimes \Dmod_{\frac{1}{2}}(\Gr_M) 
\overset{(\sotimes) \otimes (\on{Id})}\longrightarrow \\
\to \Dmod(\Gr_G)\otimes \Dmod_{\frac{1}{2}}(\Gr_M)
\overset{\Gamma_\dr(\Gr_G,-)\otimes \on{Id}}\longrightarrow \Dmod_{\frac{1}{2}}(\Gr_M)\overset{[\text{shift}]}\to \Dmod_{\frac{1}{2}}(\Gr_M),
\end{multline}
where:

\begin{itemize}

\item The category $\left(\Dmod_{\frac{1}{2}}(\Gr_G)\otimes \Dmod_{\frac{1}{2}}(\Gr_M)\right)^{\fL(P^-)}$ makes sense 
due to the identification of the two multiplicative $\mu_2$-gerbes on $\fL(P^-)$ (one obtained by restriction from
$\fL(G)$, and another from $\fL(M)$) that results from the existence of the square root $\on{det}^{\otimes \frac{1}{2}}_{\Gr_{G,M}}$;

\item The $\sotimes$-tensor product functor 
$$\Dmod_{\frac{1}{2}}(\Gr_G)\otimes \Dmod_{\frac{1}{2}}(\Gr_G)\to \Dmod(\Gr_G)$$
makes sense due to the fact that the square of the gerbe $\det_{\Gr_G}^{\frac{1}{2}}$ is canonically trivial. 

\medskip

\item The last arrow is the cohomological shift by $\langle \lambda,2\rhoch_P\rangle$ on 
$\Gr_M^\lambda$.

\end{itemize}

\sssec{}

Denote the functor \eqref{e:gen Jacquet functor} by $J_\Gr^{-,\on{pre-enh}}$. It has the following properties:

\medskip

\begin{enumerate}

\item It has an equivariance property with respect to $\fL(P^-)$, 
similar to that of \eqref{e:gen Jacquet functor non-ren}; 

\smallskip

\item It carries a natural factorization structure;

\smallskip

\item It is compatible with the $\Sph_M$-actions.

\end{enumerate}

\medskip

Furthermore, $J_\Gr^{-,\on{pre-enh}}$ factors as 
$$\Dmod_{\frac{1}{2}}(\Gr_G)\otimes \on{I}(G,P^-)^{\on{loc}} \to
\Dmod_{\frac{1}{2}}(\Gr_G)\underset{\Sph_G}\otimes \on{I}(G,P^-)^{\on{loc}}\to \Dmod_{\frac{1}{2}}(\Gr_M),$$
and the resulting functor
\begin{equation} \label{e:J Gr untwisted}
\Dmod_{\frac{1}{2}}(\Gr_G)\underset{\Sph_G}\otimes \on{I}(G,P^-)^{\on{loc}}\to \Dmod_{\frac{1}{2}}(\Gr_M)
\end{equation}
also has the properties (1)-(3) above. 

\medskip

We denote the functor \eqref{e:J Gr untwisted} by $J_\Gr^{-,\on{enh}}$. 

\sssec{} 

We will denote by
$$J_{\Gr}^{-,!}:\Dmod_{\frac{1}{2}}(\Gr_G)\to \Dmod_{\frac{1}{2}}(\Gr_M)$$
the functor obtained from $J_\Gr^{-,\on{pre-enh}}$ by 
inserting $\Delta^{-,\semiinf}$ along the second factor.

\medskip

Note that the functor $\Dmod_{\frac{1}{2}}(\Gr_G)\to \Dmod_{\frac{1}{2}}(\Gr_M)$, obtained from $J_\Gr^{-,\on{pre-enh}}$
by inserting $\nabla^{-,\semiinf}$ along the second factor, identifies with the functor $J_{\Gr}^{-,*}$ from
\secref{sss:corrected Jacquet}. 

\medskip

We will denote the functor $\Dmod_{\frac{1}{2}}(\Gr_G)\to \Dmod_{\frac{1}{2}}(\Gr_M)$, obtained from 
$J_\Gr^{-,\on{pre-enh}}$ by 
inserting $\IC^{-,\semiinf}$ along the second factor by $J_{\Gr}^{-,!*}$. 

\sssec{}  \label{sss:tricky twisted enh Jacq}

Let $\CP_M$ be an $M$-torsor on $X$. Note also that the functor \eqref{e:gen Jacquet functor} admits a twisted version
\begin{equation} \label{e:gen Jacquet functor twisted}
\Dmod_{\frac{1}{2}}(\Gr_{G,\CP_M})\otimes \on{I}(G,P^-)^{\on{loc}}_{\CP_M}\to \Dmod_{\frac{1}{2}}(\Gr_{M,\CP_M}).
\end{equation}

\medskip

Let us now 
be given an $M$-bundle $\CP'_M$ and a $Z_M$-bundle $\CP''_{Z_M}$ on $X$, so that 
$$\CP_M\simeq \CP'_M\otimes \CP''_{Z_M}.$$

\medskip

We have the following version of the functor \eqref{e:gen Jacquet functor twisted}:
\begin{multline} \label{e:gen Jacquet functor twisted center}
\Dmod_{\frac{1}{2}}(\Gr_{G,\CP_M})\otimes \on{I}(G,P^-)^{\on{loc}}_{\CP''_{Z_M}}
\overset{\on{Id}\otimes \alpha_{\CP'_M,\on{taut}}}\longrightarrow 
\Dmod_{\frac{1}{2}}(\Gr_{G,\CP_M})\otimes \on{I}(G,P^-)^{\on{loc}}_{\CP_M} \overset{\text{\eqref{e:gen Jacquet functor twisted}}}\longrightarrow \\
\to \Dmod_{\frac{1}{2}}(\Gr_{M,\CP_M})\overset{\alpha_{\CP''_{Z_M},\on{cent}}}\longrightarrow \Dmod_{\frac{1}{2}}(\Gr_{M,\CP'_M}).
\end{multline} 

\sssec{} \label{sss:prop enh Jacq}

The functor \eqref{e:gen Jacquet functor twisted center} inherits the following properties:

\medskip

\begin{enumerate}

\item It is equivariant with respect to the natural action of $\fL(P^-)_{\CP_M}$ on 
$\Dmod_{\frac{1}{2}}(\Gr_{G,\CP_M})$ and the action of 
$$\fL(P^-)_{\CP_M}\twoheadrightarrow \fL(M)_{\CP_M}\overset{\alpha_{\CP''_{Z_M},\on{cent}}}\simeq \fL(M)_{\CP'_M}$$
on $\Dmod_{\frac{1}{2}}(\Gr_{M,\CP'_M})$.

\medskip

\item It carries a natural factorization structure;

\bigskip

\item It is compatible with the $\Sph_M$-action on $\on{I}(G,P^-)^{\on{loc}}_{\CP''_{Z_M}}$ 
from \secref{sss:convent Sph-action on IGP}, and the action of $\Sph_M$ on $\Dmod_{\frac{1}{2}}(\Gr_{M,\CP'_M})$
obtained from the identification
$$\Sph_M\overset{\alpha_{\CP'_M,\on{taut}}}\simeq  \Sph_{M,\CP'_M},$$
(see \secref{sss:remove twist by M-bundle}). 

\end{enumerate}

\sssec{}

In practice, we will take 
$$\CP'_M:=\rho_M(\omega_X),\,\, \CP'_{Z_M}:=\rho_P(\omega_X),$$
so that 
$$\CP_M=\rho(\omega_X).$$

We will consider the resulting functor
\begin{equation} \label{e:gen Jacquet functor twisted omega}
\Dmod_{\frac{1}{2}}(\Gr_{G,\rho(\omega_X)})\otimes \on{I}(G,P^-)^{\on{loc}}_{\rho_P(\omega_X)}\to \Dmod_{\frac{1}{2}}(\Gr_{M,\rho_M(\omega_X)}).
\end{equation}

We will denote the functor \eqref{e:gen Jacquet functor twisted omega} by $J_{\Gr}^{-,\on{pre-enh}}$, i.e., the same symbol as in the untwisted situation.

\sssec{}

The functor $J_{\Gr}^{-,\on{pre-enh}}$ naturally factors as
\begin{multline*} 
\Dmod_{\frac{1}{2}}(\Gr_{G,\rho(\omega_X)})\otimes \on{I}(G,P^-)^{\on{loc}}_{\rho_P(\omega_X)}\to \\
\to \Dmod_{\frac{1}{2}}(\Gr_{G,\rho(\omega_X)})\underset{\Sph_G}\otimes \on{I}(G,P^-)^{\on{loc}}_{\rho_P(\omega_X)}\to 
\Dmod_{\frac{1}{2}}(\Gr_{M,\rho_M(\omega_X)}).
\end{multline*}

We will denote the resulting functor 
$$\Dmod_{\frac{1}{2}}(\Gr_{G,\rho(\omega_X)})\underset{\Sph_G}\otimes \on{I}(G,P^-)^{\on{loc}}_{\rho_P(\omega_X)}\to 
\Dmod_{\frac{1}{2}}(\Gr_{M,\rho_M(\omega_X)})$$
by $J_{\Gr}^{-,\on{enh}}$.  It inherits properties (1)-(3) from \secref{sss:prop enh Jacq}. We will refer to it as the \emph{enhanced Jacquet functor}. 

%

\sssec{}

We will denote by the same symbols 
\begin{equation} \label{e:Jacquet Whit functors}
J_{\Gr}^{-,!},\,\, J_{\Gr}^{-,!*} \text{ and } J_{\Gr}^{-,*}
\end{equation} 
the functors
$$\Dmod_{\frac{1}{2}}(\Gr_{G,\rho(\omega_X)})\to \Dmod_{\frac{1}{2}}(\Gr_{M,\rho_M(\omega_X)})$$
obtained from $J_{\Gr}^{-,\on{pre-enh}}$ by inserting the objects
$$\Delta^{-,\semiinf},\,\, \IC^{-,\semiinf} \text{ and } \nabla^{-,\semiinf}\in \on{I}(G,P^-)^{\on{loc}}_{\rho_P(\omega_X)},$$
respectively, along the second factor.

\sssec{} \label{sss:Jacqet to Whit}

The equivariance property of \eqref{e:gen Jacquet functor twisted omega} with respect\footnote{Here $N(M)$ is the maximal unipotent
subgroup of the Levi $M$.} to $\fL(N(M))\subset \fL(P^-)$ implies that the functors 
$J_{\Gr}^{-,\on{pre-enh}}$ and $J_{\Gr}^{-,\on{enh}}$ restrict to functors
\begin{equation} \label{e:Whit semiinf pair}
\Whit^!(G)\otimes \on{I}(G,P^-)^{\on{loc}}_{\rho_P(\omega_X)}\to \Whit^!(M) \text{ and }
\Whit^!(G)\underset{\Sph_G}\otimes \on{I}(G,P^-)^{\on{loc}}_{\rho_P(\omega_X)}\to \Whit^!(M).
\end{equation} 
respectively. 

\medskip

We will denote the resulting functors \eqref{e:Whit semiinf pair} by 
\begin{equation} \label{e:Whit semiinf ult}
J_{\Whit}^{-,\on{pre-enh}} \text{ and } J_{\Whit}^{-,\on{enh}},
\end{equation} 
respectively. 

\medskip

We will denote by
$$J_{\Whit}^{-,!},\,\, J_{\Whit}^{-,!*} \text{ and } J_{\Whit}^{-,*}$$
the functors
$$\Whit^!(G)\to \Whit^!(M)$$
obtained from the functors \eqref{e:Jacquet Whit functors} by restriction. 

\sssec{}

Since the functor $J_{\Gr}^{-,\on{enh}}$ respects the $\Sph_M$-actions on the two sides, by duality (see \secref{sss:duality over Sph}), it gives rise to a 
$\Sph_G$-linear functor
\begin{equation} \label{e:Jacquet Whit co}
\Dmod_{\frac{1}{2}}(\Gr_{G,\rho(\omega_X)})\to \Dmod_{\frac{1}{2}}(\Gr_{M,\rho_M(\omega_X)})
\underset{\Sph_M}\otimes \on{I}(G,P^-)^{\on{loc}}_{\on{co},\rho_P(\omega_X)},
\end{equation} 
which we will denote by $J_{\Gr}^{-,\on{enh}_{\on{co}}}$. 
 
\medskip 

The functor \eqref{e:Jacquet Whit co} gives rise to a $\Sph_G$-linear functor
$$J_{\Whit}^{-,\on{enh}_{\on{co}}}:\Whit^!(G)\to \Whit^!(M) \underset{\Sph_M}\otimes \on{I}(G,P^-)^{\on{loc}}_{\on{co},\rho_P(\omega_X)}.$$

\sssec{}

In what follows we will use the notations
$$\Whit^!(M)^{-,\on{enh}}:=\Whit^!(M) \underset{\Sph_M}\otimes \on{I}(G,P^-)^{\on{loc}}_{\rho_P(\omega_X)}$$
and 
$$\Whit^!(M)^{-,\on{enh}_{\on{co}}}:=\Whit^!(M) \underset{\Sph_M}\otimes \on{I}(G,P^-)^{\on{loc}}_{\on{co},\rho_P(\omega_X)}.$$

So, $J_{\Whit}^{-,\on{enh}_{\on{co}}}$ is a functor
$$\Whit^!(G)\to \Whit^!(M)^{-,\on{enh}_{\on{co}}}.$$ 

We recover the functor $J_{\Whit}^{-,!}$ from  $J_{\Whit}^{-,\on{enh}_{\on{co}}}$ as the composite
$$\Whit^!(G)\overset{J_{\Whit}^{-,\on{enh}_{\on{co}}}} \longrightarrow \Whit^!(M)^{-,\on{enh}_{\on{co}}}\overset{\on{Id}\otimes \oblv_{\semiinf\to\Sph_M,\on{co}}}
\longrightarrow\Whit^!(M) \underset{\Sph_M}\otimes \Sph_M=\Whit^!(M).$$

\section{Spectral semi-infinite category/ies}  \label{s:spectral semiinf}

In this section we introduce the spectral counterpart of $\on{I}(G,P^-)^{\on{loc}}_{\rho_P(\omega_X)}$, which is the category 
$$\on{I}(\cG,\cP^-)^{\on{spec,loc}},$$
thought of as ind-coherent sheaves on the (factorization) space 
\begin{equation} \label{e:triple product}
\LS^\reg_\cG\underset{\LS^\mer_\cG}\times \LS^\mer_{\cP^-}\underset{\LS^\mer_\cM}\times \LS^\reg_\cM.
\end{equation} 

When doing so, we will have to overcome difficulties having to do with the fact that the Ran version of \eqref{e:triple product}
is not locally (almost) of finite type. This discussion is mostly delegated to the Appendix, \secref{s:semiinf spec}. 

\medskip

The main result stated in this section is \thmref{t:semiinf geom Satake}, which says that there is a canonical equivalence
$$\on{I}(G,P^-)^{\on{loc}}_{\rho_P(\omega_X)}\overset{\Sat^{-,\semiinf}}\simeq \on{I}(\cG,\cP^-)^{\on{spec,loc}};$$
the proof of this theorem is also delegated to the Appendix and occupies Sects. \ref{s:Sith}-\ref{s:semiinf Sat}.



\ssec{The spherical spectral semi-infinite category}

\sssec{} \label{sss:Hecke spec GP}

Denote by
$$\on{Hecke}^{\on{spec,loc}}_{\cG,\cP^-}$$
the factorization space \eqref{e:triple product}. 

\sssec{}

For a fixed $\ul{x}\in \Ran$, the corresponding space
$$\on{Hecke}^{\on{spec,loc}}_{\cG,\cP^-,\ul{x}}$$ is an algebraic stack locally almost of finite type. Indeed, for $\ul{x}$
being a singleton $x$, we have\footnote{The isomorphism below follows from the fact that the formal completion $(\LS^\mer_{\sH,x})^\wedge_\reg$ of $\LS^\mer_{\sH,x}$ along
$\LS^\reg_{\sH,x}$ identifies with $\sh^\wedge/\on{Ad}(\sH)$, where $\sh^\wedge$ is the formal completion of $\sh$ at $0$.}
$$\on{Hecke}^{\on{spec,loc}}_{\cG,\cP^-,x}\simeq \on{pt}/\cG\underset{\cg/\on{Ad}(\cG)}\times \cp^-/\on{Ad}(\cP^-)\underset{\cm/\on{Ad}(\cM)}\times \on{pt}/\cM.$$

\medskip

So the category $\IndCoh(\on{Hecke}^{\on{spec,loc}}_{\cG,\cP^-,\ul{x}})$ is well-defined.

\medskip

However, when $\ul{x}$ moves in families over $\Ran$, as an algebro-geometric object, $\on{Hecke}^{\on{spec,loc}}_{\cG,\cP^-}$ does not
fit into the framework developed in \cite[Appendex A.5]{GLC2}, in which the category $\IndCoh^*(-)$ is defined. 

\medskip

Rather, we refer the reader to
\secref{ss:semiinf spec}, where we define
$$\IndCoh^*(\on{Hecke}^{\on{spec,loc}}_{\cG,\cP^-})$$
as a (unital) factorization category. 

\medskip

In addition, in \secref{sss:dual semiinf spec}, we will define the factorization category $$\IndCoh^!(\on{Hecke}^{\on{spec,loc}}_{\cG,\cP^-}),$$
which will turn out to be the dual of $\IndCoh^*(\on{Hecke}^{\on{spec,loc}}_{\cG,\cP^-})$. 

\sssec{}

We will use the notations:
$$\on{I}(\cG,\cP^-)^{\on{spec,loc}}:=\IndCoh^*(\on{Hecke}^{\on{spec,loc}}_{\cG,\cP^-})$$ 
and 
$$\on{I}(\cG,\cP^-)_{\on{co}}^{\on{spec,loc}}:=\IndCoh^!(\on{Hecke}^{\on{spec,loc}}_{\cG,\cP^-}).$$ 

\sssec{} \label{sss:spec Hecke action on semiinf}

At a fixed point $\ul{x}\in \Ran$, the category $\on{I}(\cG,\cP^-)^{\on{spec,loc}}_{\ul{x}}=\IndCoh^*(\on{Hecke}^{\on{spec,loc}}_{\cG,\cP^-,\ul{x}})$ 
carries a natural structure of bimodule category with respect to 
$$\Sph^{\on{spec}}_{\cM,\ul{x}}:=\IndCoh^*(\on{Hecke}^{\on{spec,loc}}_{\cM,\ul{x}}) \text{ and }
\Sph^{\on{spec}}_{\cG,\ul{x}}:=\IndCoh^*(\on{Hecke}^{\on{spec,loc}}_{\cG,\ul{x}}).$$

In \secref{sss:semiinf Hecke}, we will construct a bimodule structure on $\on{I}(\cG,\cP^-)^{\on{spec,loc}}$
with respect to 
$$\Sph^{\on{spec}}_\cM \text{ and } \Sph^{\on{spec}}_\cG$$
in the factorization setting. 
  
\sssec{}

We have a correspondence
$$
\CD
\LS^\reg_{\cP^-} @>{\sfp^{-,\on{spec}}\times \iota}>> \LS^\reg_\cG\underset{\LS^\mer_\cG}\times \LS^\mer_{\cP^-}\\
@V{\sfq^{-,\on{spec}}}VV \\
\LS^\reg_\cM
\endCD
$$
and its base change
\begin{equation} \label{e:monadic diag IGP spec}
\CD
\LS^\reg_{\cP^-}\underset{\LS^\mer_\cM}\times \LS^\reg_\cM @>>>  \on{Hecke}^{\on{spec,loc}}_{\cG,\cP^-} \\
@VVV \\
\on{Hecke}^{\on{spec,loc}}_\cM. 
\endCD
\end{equation}

\sssec{} \label{sss:Omega tilde spec ptws}

For a fixed $\ul{x}\in \Ran$, the functors of !-pull and *-push along \eqref{e:monadic diag IGP spec} define a forgetful functor 
$$\oblv^{\on{spec}}_{\semiinf\to\Sph_\cM}:\on{I}(\cG,\cP^-)^{\on{spec,loc}}_{\ul{x}}\to \Sph_{\cM,\ul{x}}^{\on{spec}},$$
which admits a left adjoint, denoted $\ind^{\on{spec}}_{\Sph_\cM\to\semiinf}$, given by *-pull followed by *-push. 

\medskip

Thus, we obtain an adjoint pair
\begin{equation} \label{e:IGP to M adj spec x}
\ind^{\on{spec}}_{\Sph_\cM\to\semiinf}:\Sph_{\cM,\ul{x}}^{\on{spec}}\rightleftarrows \on{I}(\cG,\cP^-)_{\ul{x}}^{\on{spec,loc}}:\oblv^{\on{spec}}_{\semiinf\to\Sph_\cM}.
\end{equation}

Moreover, it is easy to see that the adjunction \eqref{e:IGP to M adj spec} is monadic
and respects the $\Sph_{\cM,\ul{x}}^{\on{spec}}$-actions.

\sssec{} \label{sss:Omega tilde spec Ran}

In \secref{sss:semiinf spec oblv}, we will construct the functor 
$$\oblv^{\on{spec}}_{\semiinf\to\Sph_\cM}:\on{I}(\cG,\cP^-)^{\on{spec,loc}}\to \Sph_\cM^{\on{spec}}$$
in the factorization setting.

\medskip

Moreover, we will show that $\oblv^{\on{spec}}_{\semiinf\to\Sph_M}$ admits a left adjoint, denoted $\ind^{\on{spec}}_{\Sph_\cM\to\semiinf}$,
which is a (strictly unital) factorization functor. 

\medskip

As in \secref{sss:Omega tilde}, we can view the adjunction 
\begin{equation} \label{e:IGP to M adj spec}
\ind^{\on{spec}}_{\Sph_\cM\to\semiinf}:\Sph_\cM^{\on{spec}}\rightleftarrows \on{I}(\cG,\cP^-)^{\on{spec,loc}}:\oblv^{\on{spec}}_{\semiinf\to\Sph_\cM}
\end{equation}
as between 
$\Sph_\cM^{\on{spec}}$-module categories.

\sssec{}

The factorization unit 
$$\one_{\on{I}(\cG,\cP^-)^{\on{spec,loc}}}\in \on{I}(\cG,\cP^-)^{\on{spec,loc}}$$
is the $(\IndCoh,*)$-direct image of 
$$\CO_{\LS^\reg_{\cP^-}}\in \QCoh(\LS^\reg_{\cP^-}) \simeq \IndCoh^*(\LS^\reg_{\cP^-})$$
along 
\begin{equation} \label{e:iota basic}
\LS^\reg_{\cP^-}\overset{\iota}\to \LS^\reg_\cG\underset{\LS^\mer_\cG}\times \LS^\reg_{\cP^-}\underset{\LS^\mer_\cM}\times \LS^\reg_\cM =
\on{Hecke}^{\on{spec,loc}}_{\cG,\cP^-}.
\end{equation} 

We will also use the notation
$$\Delta^{-,\on{spec},\semiinf}:=\one_{\on{I}(\cG,\cP^-)^{\on{spec,loc}}}.$$

Since the functor $\ind^{\on{spec}}_{\Sph_\cM\to\semiinf}$ is $\Sph_\cM^{\on{spec}}$-linear, we have
$$\ind^{\on{spec}}_{\Sph_\cM\to\semiinf}(-)\simeq (-)\underset{\cM}\star \Delta^{-,\on{spec},\semiinf},$$
where $\underset{\cM}\star$ denotes the $\Sph_\cM^{\on{spec}}$-action on $\on{I}(\cG,\cP^-)^{\on{spec,loc}}$.

\sssec{}

The $\Sph^{\on{spec}}_\cG$-action on $\Delta^{-,\on{spec},\semiinf}$ is a $\Sph^{\on{spec}}_\cG$-linear (strictly unital) factorization functor
\begin{equation} \label{e:Sph G spec to semiing}
\Sph^{\on{spec}}_\cG\to \on{I}(\cG,\cP^-)^{\on{spec,loc}}, \quad (-)\underset{\cG}\star \Delta^{-,\on{spec},\semiinf},
\end{equation} 
where $\underset{\cG}\star$ denotes the $\Sph_\cM^{\on{spec}}$-action on $\on{I}(\cG,\cP^-)^{\on{spec,loc}}$.

\medskip

Explicitly, the functor \eqref{e:Sph G spec to semiing} is given by $(\IndCoh,*)$-pullback, followed by $(\IndCoh,*)$-direct image along the diagram
$$
\CD
\LS^\reg_\cG\underset{\LS^\mer_\cG}\times \LS^\reg_{\cP^-} @>>> \LS^\reg_\cG\underset{\LS^\mer_\cG}\times \LS^\mer_{\cP^-}
\underset{\LS^\mer_\cM}\times \LS^\reg_\cM \\
@VVV \\
\LS^\reg_\cG\underset{\LS^\mer_\cG}\times \LS^\reg_\cG.
\endCD 
$$

The functor \eqref{e:Sph G spec to semiing} admits a factorization right adjoint, given by !-pullback, followed by $(\IndCoh,*)$-direct image
\emph{in the opposite direction} along the above diagram. 

\sssec{}

In addition to the factorization functor $\ind^{\on{spec}}_{\Sph_\cM\to\semiinf}$, one can consider another functor, denoted 
$$\ind^{\on{spec},*}_{\Sph_\cM\to\semiinf}:\Sph_\cM^{\on{spec}}\to \on{I}(\cG,\cP^-)^{\on{spec,loc}},$$
defined as !-pullback, followed by $(\IndCoh,*)$-direct image along the diagram\footnote{Note the appearance of the \emph{positive} parabolic $\cP$
in this diagram.}
$$
\CD
& & & & \on{Hecke}^{\on{spec,loc}}_{\cG,\cP^-} \\
& & & & @VV{=}V \\
\LS^\reg_\cP\underset{\LS^\mer_\cP}\times \LS^\reg_\cM @>{=}>>
\LS^\reg_\cP\underset{\LS^\mer_\cP}\times \LS^\mer_{\cM}\underset{\LS^\mer_\cM}\times \LS^\reg_\cM @>>>
\LS^\reg_\cG\underset{\LS^\mer_\cG}\times \LS^\mer_{\cP^-}\underset{\LS^\mer_\cM}\times \LS^\reg_\cM  \\
@VVV \\
\LS^\reg_\cM\underset{\LS^\mer_\cM}\times \LS^\reg_\cM.
\endCD
$$

Note that the composition
$$\oblv^{\on{spec}}_{\semiinf\to\Sph_\cM}\circ \ind^{\on{spec},*}_{\Sph_\cM\to\semiinf}$$
is isomorphic to the identity endofunctor of $\Sph_\cM^{\on{spec}}$. This follows by base change from the fact that the diagram
$$
\CD
\LS^\reg_\cM @>>> \LS^\reg_{\cP^-} \\
@VVV @VVV \\
\LS^\reg_\cP\underset{\LS^\mer_\cP}\times \LS^\mer_{\cM} @>>> \LS^\reg_\cG\underset{\LS^\mer_\cG}\times \LS^\mer_{\cP^-}
\endCD
$$
is Cartesian. 

\medskip

We denote
$$\nabla^{-,\on{spec},\semiinf}:=\ind^{\on{spec},*}_{\Sph_\cM\to\semiinf}(\one_{\Sph_\cM^{\on{spec}}})\in \on{I}(\cG,\cP^-)^{\on{spec,loc}}.$$

We have
\begin{equation} \label{e:oblv nabla spec}
\oblv^{\on{spec}}_{\semiinf\to\Sph_\cM}(\nabla^{-,\on{spec},\semiinf})\simeq \one_{\Sph_\cM^{\on{spec}}}.
\end{equation}

\sssec{}

In addition to the factorization algebras 
$$\Delta^{-,\on{spec},\semiinf} \text{ and } \nabla^{-,\on{spec},\semiinf},$$
we will consider the following factorization algebra object in $\on{I}(\cG,\cP^-)^{\on{spec,loc}}$,
denoted
$$\IC^{-,\on{spec},\semiinf}.$$

\medskip

Namely, it is is defined as the mage of $\CO_{\LS^\reg_{\cM}}\in \QCoh(\LS^\reg_\cM)$ along
$$\QCoh(\LS^\reg_\cM)\to \QCoh(\LS^\reg_{\cP^-})\simeq \IndCoh^*(\LS^\reg_{\cP^-})\overset{\iota^\IndCoh_*}\longrightarrow
\on{Hecke}^{\on{spec,loc}}_{\cG,\cP^-},$$
where the first arrow is \emph{direct image} along the map
\begin{equation} \label{e:LS M to LS P}
\LS^\reg_\cM\to \LS^\reg_{\cP^-},
\end{equation} 
corresponding to the canonical map\footnote{A choice of a simpitting $\cM\to \cP^-$ gives rise to a map
$\on{pt}/\cM\to \on{pt}/\cP^-$; every two splittings differ by conjugation by a a \emph{uniquely} defined element of
$\cN^-_P$, and the resulting two maps $\on{pt}/\cM\rightrightarrows \on{pt}/\cP^-$ coincide.} $\on{pt}/\cM\to \on{pt}/\cP^-$.
 
\ssec{The (factorization) algebras \texorpdfstring{$\wt\Omega^{\on{spec}}$}{wtOmega} 
and \texorpdfstring{$\Omega^{\on{spec}}$}{Omega}} \label{ss:wt Omega spec}

The definition of $\on{I}(\cG,\cP^-)^{\on{spec,loc}}$ involves some complicated algebraic geometry in infinite type. Yet, the datum
of this category can be equivalently recorded by a more elementary object: an associative factorization algebra in $\Sph_\cM^{\on{spec}}$,
denoted $\wt\Omega^{\on{spec}}$. 

\medskip

We introduce it in the subsection. 

\sssec{} \label{sss:wt Omega spec}

Let $$\wt\Omega^{\on{spec}}\in \Sph_\cM^{\on{spec}}$$
denote the associative (factorization) algebra object 
so that \eqref{e:IGP to M adj spec} identifies with the adjunction
\begin{equation} \label{e:IGP to M adj spec as Omega}
\ind_{\wt\Omega^{\on{spec}}}:\Sph_\cM^{\on{spec}}\rightleftarrows \wt\Omega^{\on{spec}}\mod^r(\Sph^{\on{spec}}_\cM):\oblv_{\wt\Omega^{\on{spec}}}.
\end{equation}

\sssec{}

In terms of the identification
\begin{equation} \label{e:IGP spec via Omega spec}
\on{I}(\cG,\cP^-)^{\on{spec,loc}}\simeq \wt\Omega^{\on{spec}}\mod^r(\Sph^{\on{spec}}_\cM),
\end{equation} 
the object
$$\Delta^{-,\on{spec},\semiinf}\in \on{I}(\cG,\cP^-)^{\on{spec,loc}}$$ 
corresponds to
$$\wt\Omega^{\on{spec}}\in \wt\Omega^{\on{spec}}\mod^r(\Sph^{\on{spec}}_\cM).$$

\sssec{} 

Note that the adjunction
$$(\sfq^{-,\on{spec}})^*:\QCoh(\LS^\reg_\cM)\rightleftarrows \QCoh(\LS^\reg_{\cP^-}):(\sfq^{-,\on{spec}})_*$$
identifies with
\begin{equation} \label{e:P to M adj spec}
\Res^\cM_{\cP^-}:\Rep(\cM)\rightleftarrows \Rep(\cP^-):\on{inv}_{\cN^-_P},
\end{equation}
where $\on{inv}_{\cN^-_P}$ is the functor of $\cn^-_P$-invariants (a.k.a. cohomological
Chevalley complex). 

\sssec{} \label{sss:Omega spec}

Let $$\Omega^{\on{spec}}\in \on{FactAlg}^{\on{untl}}(\Rep(\cM))$$ denote the (commutative) factorization algebra 
corresponding (via \cite[Corollary C.8.8]{GLC2}) to 
$$\on{inv}_{\cN^-_P}(k)\in \on{ComAlg}(\Rep(\cM));$$
note that the latter is the same as 
the cohomological Chevalley complex 
$$\on{C}_{\on{chev}}^\cdot(\cn^-_P),$$
of $\cn^-_P$ with coefficients in the trivial module. 

\medskip

The (monadic) adjunction \eqref{e:P to M adj spec} can therefore be rewritten as
\begin{equation} \label{e:P to M adj spec as Omega}
\ind_{\Omega^{\on{spec}}}:\Rep(\cM) \rightleftarrows \Omega^{\on{spec}}\mod^r(\Rep(\cM)):\oblv_{\Omega^{\on{spec}}}.
\end{equation}

\sssec{}

By a slight abuse of notation, we will denote by the same symbol $\Omega^{\on{spec}}$ its image along the
functor
$$\Rep(\cM) \overset{\on{nv}}\to \Sph_\cM^{\on{spec}}.$$

\medskip

From diagram \eqref{e:monadic diag IGP spec} we obtain that there exists a canonically defined map of associative algebra objects
\begin{equation} \label{e:Omega spec to Omega tilde}
\Omega^{\on{spec}}\to \wt\Omega^{\on{spec}}, 
\end{equation}
so that the adjunction \eqref{e:IGP to M adj spec as Omega} factors as
a composition of 
\begin{equation} \label{e:P to M adj spec as Omega der}
\ind_{\Omega^{\on{spec}}}: \Sph_\cM^{\on{spec}} \rightleftarrows \Omega^{\on{spec}}\mod^r(\Sph_\cM^{\on{spec}}):\oblv_{\Omega^{\on{spec}}}.
\end{equation}
and 
\begin{equation} \label{e:P to M adj spec as Omega and tilde}
\ind_{\Omega^{\on{spec}}\to \wt\Omega^{\on{spec}}}:
\Omega^{\on{spec}}\mod^r(\Sph_\cM^{\on{spec}} ) \rightleftarrows \wt\Omega^{\on{spec}}\mod^r(\Sph^{\on{spec}}_\cM):
\oblv_{\wt\Omega^{\on{spec}}\to \Omega^{\on{spec}}}.
\end{equation}

\sssec{} \label{sss:Omega tilde spec augm}

Note that the category $\Sph^{\on{spec}}_\cM$ is naturally graded by
$$\Lambda_M/\Lambda_{[M,M]_{\on{sc}}}\simeq \pi_{1,\on{alg}}(M)\simeq \Hom(Z(\cM),\BG_m),$$ 
and in terms of this grading, the algebras
$$\wt\Omega^{\on{spec}} \text{ and } \Omega^{\on{spec}}$$
are graded by the sub-monoid 
$$\Lambda^{\on{pos}}_{G,P}\subset \Lambda_M/\Lambda_{[M,M]_{\on{sc}}}$$ spanned by the images of the positive simple coroots,
with the $0$-weight component being $\one_{\Sph^{\on{spec}}_\cM}$. 

\medskip

In particular, both $\wt\Omega^{\on{spec}}$ and $\Omega^{\on{spec}}$ are naturally augmented. In particular, we have the 
well-defined objects
$$\one_{\Sph^{\on{spec}}_\cM}\in \wt\Omega^{\on{spec}}\mod^r(\Sph^{\on{spec}}_\cM) \text{ and }
\one_{\Sph^{\on{spec}}_\cM}\in \Omega^{\on{spec}}\mod^r(\Sph^{\on{spec}}_\cM).$$

\sssec{}

By \eqref{e:oblv nabla spec}, in terms of the equivalence \eqref{e:IGP spec via Omega spec}, the object
$$\nabla^{-,\on{spec},\semiinf}\in \on{I}(\cG,\cP^-)^{\on{spec,loc}}$$ 
corresponds to
$$\one_{\Sph^{\on{spec}}_\cM}\in \wt\Omega^{\on{spec}}\mod^r(\Sph^{\on{spec}}_\cM).$$

\sssec{}

By construction, the object 
$$\IC^{-,\on{spec},\semiinf}\in \on{I}(\cG,\cP^-)^{\on{spec,loc}}$$ 
corresponds to
$$\ind_{\Omega^{\on{spec}}\to \wt\Omega^{\on{spec}}}(\one_{\Sph^{\on{spec}}_\cM})\in \wt\Omega^{\on{spec}}\mod^r(\Sph^{\on{spec}}_\cM).$$

\ssec{Spectral Jacquet functors}

This subsection is a spectral counterpart of \secref{ss:enh Jacquet}: we will show how objects of $\on{I}(\cG,\cP^-)^{\on{spec,loc}}$
give rise to generalized Jacquet functors $\Rep(\cG)\to \Rep(\cM)$. 

\sssec{}

Consider the diagram 
\begin{equation} \label{e:gen Jacquet functor spec diagram}
\vcenter
{\xy
(0,0)*+{\LS^\reg_\cG}="A";
(40,0)*+{\LS^\reg_\cM.}="B";
(20,20)*+{\LS^\reg_\cG\underset{\LS^\mer_\cG}\times \LS^\mer_{\cP^-}\underset{\LS^\mer_\cM}\times \LS^\reg_\cM}="C";
{\ar@{->}_{p_1} "C";"A"};
{\ar@{->}^{p_2} "C";"B"};
\endxy}
\end{equation}

\medskip

At a fixed point $\ul{x}\in \Ran$, the corresponding diagram
$$
\vcenter
{\xy
(0,0)*+{\LS^\reg_{\cG,\ul{x}}}="A";
(40,0)*+{\LS^\reg_{\cM,\ul{x}}}="B";
(20,20)*+{\LS^\reg_{\cG,\ul{x}}\underset{\LS^\mer_{\cG,\ul{x}}}\times \LS^\mer_{\cP^-,\ul{x}}\underset{\LS^\mer_{\cM,\ul{x}}}\times \LS^\reg_{\cM,\ul{x}}}="C";
{\ar@{->}_{p_1} "C";"A"};
{\ar@{->}^{p_2} "C";"B"};
\endxy}
$$
consists of stacks locally almost of finite type. Hence, we have a well-defined functor
\begin{multline*} 
\Rep(\cG)_{\ul{x}}\otimes \on{I}(\cG,\cP^-)^{\on{spec,loc}}_{\ul{x}}=
\QCoh(\LS^\reg_{\cG,\ul{x}})\otimes 
\IndCoh^*(\LS^\reg_{\cG,\ul{x}}\underset{\LS^\mer_{\cG,\ul{x}}}\times \LS^\mer_{\cP^-,\ul{x}}\underset{\LS^\mer_{\cM,\ul{x}}}\times \LS^\reg_{\cM,\ul{x}})
\overset{p_1^*(-)\otimes (-)}\longrightarrow  \\
\to \IndCoh^*(\LS^\reg_{\cG,\ul{x}}\underset{\LS^\mer_{\cG,\ul{x}}}\times \LS^\mer_{\cP^-,\ul{x}}\underset{\LS^\mer_{\cM,\ul{x}}}\times \LS^\reg_{\cM,\ul{x}})
\overset{(p_2)^\IndCoh_*}\longrightarrow \IndCoh^*(\LS^\reg_{\cM,\ul{x}}) \simeq \Rep(\cM)_{\ul{x}}.
\end{multline*} 

\sssec{}

In \secref{sss:gen Jacquet} we will show that the functor $(p_2)^\IndCoh_*$ in the composition
\begin{multline} \label{e:gen Jacquet functor spec} 
\Rep(\cG)\otimes \on{I}(\cG,\cP^-)^{\on{spec,loc}}=
\QCoh(\LS^\reg_\cG)\otimes 
\IndCoh^*(\LS^\reg_\cG\underset{\LS^\mer_\cG}\times \LS^\mer_{\cP^-}\underset{\LS^\mer_\cM}\times \LS^\reg_\cM)
\overset{p_1^*(-)\otimes (-)}\longrightarrow  \\
\to \IndCoh^*(\LS^\reg_\cG\underset{\LS^\mer_\cG}\times \LS^\mer_{\cP^-}\underset{\LS^\mer_\cM}\times \LS^\reg_\cM)
\overset{(p_2)^\IndCoh_*}\longrightarrow \IndCoh^*(\LS^\reg_\cM) \simeq \Rep(\cM)
\end{multline} 
is well-defined also in the factorization setting. 

\medskip

We will denote the functor \eqref{e:gen Jacquet functor spec} by $J^{-,\on{spec,pre-enh}}$.

\sssec{} For a fixed point $\ul{x}\in \Ran$, the functor $J^{-,\on{spec,pre-enh}}$ naturally factors as 
$$\Rep(\cG)_{\ul{x}}\otimes \on{I}(\cG,\cP^-)^{\on{spec,loc}}_{\ul{x}}\to
\Rep(\cG)_{\ul{x}}\underset{\Sph^{\on{spec}}_{\cG,\ul{x}}}\otimes \on{I}(\cG,\cP^-)^{\on{spec,loc}}_{\ul{x}}\to \Rep(\cM)_{\ul{x}}.$$

In \secref{sss:gen Jacquet Sph} we will show that a similar factorization takes place when $\ul{x}$ is allowed to move
in families over $\Ran$, thereby giving rise to a factorization functor
$$\Rep(\cG)\underset{\Sph^{\on{spec}}_\cG}\otimes \on{I}(\cG,\cP^-)^{\on{spec,loc}}\to \Rep(\cM),$$
which we denote by $J^{-,\on{spec,enh}}$. We will refer to it as the \emph{spectral enhanced Jacquet functor}. 

\sssec{}

We will denote by 
$$J^{-,\on{spec},!},\,\, J^{-,\on{spec},!*} \text{ and } J^{-,\on{spec},*}$$
the functors
$$\Rep(\cG)\to \Rep(\cM),$$
obtained from $J^{-,\on{spec,pre-enh}}$ by inserting the objects
$$\Delta^{-,\on{spec},\semiinf}, \IC^{-,\on{spec},\semiinf} \text{ and } \nabla^{-,\on{spec},\semiinf}\in \on{I}(\cG,\cP^-)^{\on{spec,loc}},$$
respectively, along the second factor.

\sssec{}

Unwinding the definitions, we obtain that the functor
$J^{-,\on{spec},!}$ identifies canonically with the $(\IndCoh,*)$-pullback, followed by 
$(\IndCoh,*)$-pushforward along the diagram
\begin{equation} \label{e:gen Jacquet functor spec diag reg}
\xy
(0,0)*+{\LS^\reg_\cG}="A";
(40,0)*+{\LS^\reg_\cM.}="B";
(20,20)*+{\LS^\reg_{\cP^-}}="C";
{\ar@{->}_{\sfp^{-,\on{spec}}} "C";"A"};
{\ar@{->}^{\sfq^{-,\on{spec}}} "C";"B"};
\endxy
\end{equation}

I.e., $J^{-,\on{spec},!}$ is the functor
$$\Rep(\cG) \overset{\Res^\cG_{\cP^-}}\longrightarrow \Rep(\cP^-) \overset{\on{inv}_{\cN^-_P}}\longrightarrow \Rep(\cM).$$

\sssec{}

Similarly, the functor $J^{-,\on{spec},!*}$ is simply the restriction functor
$$\Rep(\cG) \overset{\Res^\cG_\cM}\longrightarrow \Rep(\cM)$$
along $\cM\to \cG$.

\sssec{}

Finally, the functor $J^{-,\on{spec},*}$ is described as follows. Consider the fiber product
$$\CY:=\LS^\reg_\cM\underset{\LS^\reg_\cM\underset{\LS^\mer_\cM}\times \LS^\reg_\cM}\times (\LS^\reg_\cP\underset{\LS^\mer_\cP}\times \LS^\reg_\cM)
\simeq \LS^\reg_\cP\underset{\LS^\mer_\cP\underset{\LS^\mer_\cM}\times \LS^\reg_\cM}\times \LS^\reg_\cM.$$

It is equipped with maps 
$$p_1:\CY\to \LS^\reg_\cG \text{ and } p_2:\CY\to \LS^\reg_\cM.$$

Set
$$\CF_\CY:=p_2^!(\CO_{\LS^\reg_\cM})\in \IndCoh^*(\CY).$$

Then $J^{-,\on{spec},*}$ is given by 
$$\Rep(\cG)\simeq \QCoh(\LS^\reg_\cG)\overset{p_1^*(-)\otimes \CF_\CY}\to  \IndCoh^*(\CY)\overset{(p_2)^\IndCoh_*}\longrightarrow 
\IndCoh^*(\LS^\reg_\cM)\simeq \Rep(\cM).$$

\sssec{}

The functor $J^{-,\on{spec,enh}}$ is $\Sph^{\on{spec}}_\cM$-linear, and hence, by duality, it gives rise to a $\Sph^{\on{spec}}_\cG$-linear
functor
$$\Rep(\cG)\to  \Rep(\cM) \underset{\Sph^{\on{spec}}_\cM}\otimes \on{I}(\cG,\cP^-)_{\on{co}}^{\on{spec,loc}},$$
which we denote by $J^{-,\on{spec,enh}_{\on{co}}}$. 

\ssec{Semi-infinite geometric Satake} 

In this subsection we state one of the two main theorems in Part I of the paper, namely, \thmref{t:semiinf geom Satake}. 

\sssec{} \label{sss:curse}

Recall (see \cite[Sect. 1.8]{GLC2}) that the anti-involutions $\sigma$ and $\sigma^{\on{spec}}$
of $\Sph_G$ and $\Sph_\cG^{\on{spec}}$, respectively, are compatible with the geometric Satake
equivalence 
\begin{equation}  \label{e:Sat G}
\Sph_G \overset{\on{Sat}_G}\simeq \Sph^{\on{spec}}_{\cG}
\end{equation}
\emph{up to} the Chevalley involution $\tau_G$. 

\medskip

Hence, $\tau_G$ comes up every time when we compare the passage between left and right modules categories
on the two sides of geometric Satake. 

\sssec{}

Denote
$$\on{CS}_{G,\tau}:=\on{CS}_G\circ \tau_G \text{ and } \on{CS}_{M,\tau}:=\on{CS}_M\circ \tau_M.$$

\medskip

We are now ready to state the \emph{semi-infinite geometric Satake} equivalence:

\begin{thm} \label{t:semiinf geom Satake}
There exists an equivalence of unital factorization categories
$$\Sat^{-,\semiinf}:\on{I}(G,P^-)^{\on{loc}}_{\rho_P(\omega_X)}\to \on{I}(\cG,\cP^-)^{\on{spec,loc}},$$
compatible with the actions of 
\begin{equation}  \label{e:Sat G tau}
\Sph_G \overset{\on{Sat}_G}\simeq \Sph^{\on{spec}}_{\cG},
\end{equation}
which makes the diagram
\begin{equation} \label{e:semiinf geom Satake}
\CD
\Whit^!(G)\underset{\Sph_G}\otimes \on{I}(G,P^-)^{\on{loc}}_{\rho_P(\omega_X)} @>{(\on{CS}_{G,\tau})\otimes \Sat^{-,\semiinf}}>>  
\Rep(\cG)\underset{\Sph^{\on{spec}}_\cG}\otimes \on{I}(\cG,\cP^-)^{\on{spec,loc}} \\
@V{J^{-,\on{enh}}_{\Whit}}VV @VV{J^{-,\on{spec,enh}}}V \\
\Whit^!(M) @>{\on{CS}_{M,\tau}}>> \Rep(\cM)
\endCD
\end{equation} 
commute. 
Moreover, the functor $\Sat^{-,\semiinf}$ and the functors in \eqref{e:semiinf geom Satake} are compatible 
with the actions of
\begin{equation}  \label{e:Sat M}
\Sph_M \overset{\on{Sat}_{M,\tau}}\simeq \Sph^{\on{spec}}_{\cM}.
\end{equation}
\end{thm} 

We will supply a proof in \secref{s:semiinf Sat}.

%
%

\begin{rem}
The existence of \emph{an} equivalence 
$$\on{I}(G,P^-)^{\on{loc}}_{\rho_P(\omega_X)}\simeq \on{I}(\cG,\cP^-)^{\on{spec,loc}}$$
was established in \cite{CR}.\footnote{At the pointwise level, the existence of an equivalence 
between $\on{I}(G,P^-)^{\on{loc}}_{\rho_P(\omega_X)}$ and $\on{I}(\cG,\cP^-)^{\on{spec,loc}}$ had been established in \cite{ABG},
which was the reason \thmref{t:semiinf geom Satake} was conjectured.}
However, instead of the commutativity 
of \eqref{e:semiinf geom Satake}, the construction in {\it loc. cit.} makes a different diagram commute (see Remark \ref{r:Ups}
below). 

\medskip

So, the assertion of \thmref{t:semiinf geom Satake} does not quite follow from \cite{CR}.
Our proof of \thmref{t:semiinf geom Satake}
follows ideas similar to those of \cite{CR}, but is formally 
independent of it.\footnote{Except that we use \cite[Theorem 4.6.1]{CR} 
as a key technical point, see the proof of \propref{p:Omega almost ULA}.}

\medskip

That said, one can deduce \thmref{t:semiinf geom Satake} 
from the results of \cite{CR} by an argument involving duality; this will be addressed in a future publication. 

\end{rem}

\ssec{Consequences for Jacquet functors}

\sssec{}

Being an equivalence of unital factorization categories, the functor $\Sat^{-,\semiinf}$ satisfies
$$\Sat^{-,\semiinf}(\Delta^{-,\semiinf})\simeq \Delta^{-,\on{spec},\semiinf}.$$

From here, we obtain:

\begin{cor} \label{c:!-Jacquet on Whit}
The following diagram of factorization categories and functors commutes:
\begin{equation} \label{e:!-Jacquet on Whit}
\CD
\Whit^!(G) @>{\on{CS}_{G,\tau}}>>  \Rep(\cG) \\
@V{J^{-,!}_{\Whit}}VV @VV{\on{inv}_{\cN^-_P}=J^{-,\on{spec},!}}V \\
\Whit^!(M) @>{\on{CS}_{M,\tau}}>>  \Rep(\cM).
\endCD
\end{equation} 
\end{cor}

\begin{rem} \label{r:! Jaq}

Although we stated the commutativity of \eqref{e:!-Jacquet on Whit} as a corollary of \thmref{t:semiinf geom Satake}, when it comes 
to proofs, the actual logic of the proofs is the other way around: one establishes \corref{c:!-Jacquet on Whit} first, and this is 
as one of the initial steps in the proof of \thmref{t:semiinf geom Satake}. 

%

\end{rem}

\begin{rem}  \label{r:Ups}

What follows from \cite{CR} is that the diagram
\begin{equation} \label{e:J* Whit}
\CD
\Whit^!(G) @>{\on{CS}_{G,\tau}}>>  \Rep(\cG) \\
@V{J^{-,*}_{\Whit}}VV @VV{J^{-,\on{spec},*}}V \\
\Whit^!(M) @>{\on{CS}_{M,\tau}}>>  \Rep(\cM)
\endCD
\end{equation}
commutes.

\end{rem}

\sssec{}

For future reference, we introduce the functor
\begin{equation} \label{e:J Theta Whit}
J^{-,!}_{\Whit,\Theta}:\Whit_*(G)\to \Whit_*(M)
\end{equation}
so that the diagram
$$
\CD
\Whit_*(G) @>{J^{-,!}_{\Whit,\Theta}}>> \Whit_*(M) \\
@V{\Theta_{\Whit(G)}}V{\sim}V @V{\sim}V{\Theta_{\Whit(M)}}V \\
\Whit^!(G) @>{J^{-,!}_{\Whit}}>> \Whit^!(M)
\endCD
$$
commutes, where the vertical arrows are as in \cite[Sect. 1.3.12]{GLC2}. 

\medskip

From \corref{c:!-Jacquet on Whit} we obtain a commutative diagram
\begin{equation} \label{e:!-Jacquet on Whit Theta}
\CD
\Whit^*(G) @<{\FLE_{\cG,\infty}}<<  \Rep(\cG) \\
@V{J^{-,!}_{\Whit,\Theta}}VV @VV{\on{inv}_{\cN^-_P}}V \\
\Whit^*(M) @<{\FLE_{\cM,\infty}}<<  \Rep(\cM).
\endCD
\end{equation} 

\sssec{}

The compatibility of the functor $\Sat^{-,\semiinf}$ with \eqref{e:Sat M} yields
a commutative diagram 
$$
\CD
\Sph_M @>{\Sat_{M,\tau}}>>  \Sph_\cM^{\on{spec}} \\
@V{\ind_{\Sph_M\to\semiinf}}VV @VV{\ind_{\Sph^{\on{spec}}_\cM\to\semiinf}}V \\
\on{I}(G,P^-)^{\on{loc}}_{\rho_P(\omega_X)} @>{\Sat^{-,\semiinf}}>>  \on{I}(\cG,\cP^-)^{\on{spec,loc}}.
\endCD
$$

Passing to right adjoints, we obtain a commutative diagram 
$$
\CD
\Sph_M @>{\Sat_M}>>  \Sph_\cM^{\on{spec}} \\
@A{\oblv_{\semiinf\to\Sph_M}}AA @AA{\oblv_{\semiinf\to\Sph^{\on{spec}}_\cM}}A \\
\on{I}(G,P^-)^{\on{loc}}_{\rho_P(\omega_X)} @>{\Sat^{-,\semiinf}}>>  \on{I}(\cG,\cP^-)^{\on{spec,loc}}. 
\endCD
$$

\medskip

Hence, we can view the equivalence $\Sat^{-,\semiinf}$ as the statement that we have an isomorphism of 
associative (factorization) algebra objects 
\begin{equation} \label{e:semiinf Sat via Omega}
\Sat_{M,\tau}(\wt\Omega_{\rho_P(\omega_X)})\simeq \wt\Omega^{\on{spec}}.
\end{equation} 

\sssec{}

Note that for reasons of augmentation, we have
$$\Sat^{-,\semiinf}(\nabla^{-,\semiinf})\simeq \nabla^{-,\on{spec},\semiinf}.$$

Note that this isomorphism immediately implies the commutativity of \eqref{e:J* Whit}. 

\sssec{}

The next assertion will be proved in \secref{ss:semiinf Sat IC}:

\begin{prop} \label{p:semiinf Sat IC}
We have a canonical isomorphism of unital factorization algebras
$$\Sat^{-,\semiinf}(\IC^{-,\semiinf})\simeq \IC^{-,\on{spec},\semiinf}.$$
\end{prop}

\begin{cor} \label{c:semiinf Sat IC}
The following diagram of factorization categories and functors commutes:
\begin{equation} \label{e:IC Jacquet on Whit and Rep}
\CD
\Whit^!(G) @>{\on{CS}_{G,\tau}}>>  \Rep(\cG) \\
@V{J^{-,!*}_{\Whit}}VV @VV{\on{Res}^\cG_\cM\simeq J^{-,\on{spec},!*}}V \\
\Whit^!(M) @>{\on{CS}_{M,\tau}}>>  \Rep(\cM).
\endCD
\end{equation} 
\end{cor}

\begin{rem}
The commutativity of \eqref{e:IC Jacquet on Whit and Rep} for $P=B$ (essentially) follows
from \cite[Theorem 4.15.1]{Ra4}. A direct proof of \corref{c:semiinf Sat IC} (without referring to
\propref{p:semiinf Sat IC}) for an arbitrary parabolic will be given in \cite{FH}. 



\end{rem}

\ssec{Self-duality on the spectral side}

\sssec{}

Recall that according to \eqref{e:self-duality IGP twisted}, we have a canonical equivalence
$$(\on{I}(G,P^-)^{\on{loc}}_{\rho_P(\omega_X)})^\vee \simeq \on{I}(G,P^-)^{\on{loc}}_{\rho_P(\omega_X)}.$$

Combining with \thmref{t:semiinf geom Satake}, 
we obtain an equivalence
\begin{equation} \label{e:self-duality I G P spec}
\on{I}(\cG,\cP^-)^{\on{spec,loc}}_{\on{co}} \simeq \on{I}(\cG,\cP^-)^{\on{spec,loc}}.
\end{equation}

\begin{rem}  \label{r:IGP spec duality}

Note according to \secref{sss:self-duality twisted geom cen}, the identification \eqref{e:self-duality I G P spec}
is compatible with:

\begin{itemize}

\item The natural actions of $\Sph_\cG^{\on{spec}}$ on the two sides;

\medskip

\item The natural actions of $\Sph_\cM^{\on{spec}}$ on the two sides, \emph{up to}
the automorphism of $\Sph_\cM^{\on{spec}}$ that corresponds under $\Sat_M$ to
the automorphism $\on{transl}^*_{2\rho_P(\omega_X)}$ of $\Sph_M$. 

\end{itemize} 

\end{rem}

\sssec{}

One can describe the equivalence \eqref{e:self-duality I G P spec} purely in spectral terms. However, this
is non-trivial, and will be addressed in a separate publication. 

\medskip

For example, at a fixed $\ul{x}\in \Ran$, in which case $\on{Hecke}^{\on{spec,loc}}_{\cG,\cP^-,\ul{x}}$ is a stack
locally almost of finite type and so
$$\on{I}(\cG,\cP^-)^{\on{spec,loc}}_{\on{co}} \simeq \IndCoh(\on{Hecke}^{\on{spec,loc}}_{\cG,\cP^-,\ul{x}})\simeq
\on{I}(\cG,\cP^-)^{\on{spec,loc}},$$
the identification \eqref{e:self-duality I G P spec} corresponds to the usual (i.e., Serre) self-duality of 
$\IndCoh(\on{Hecke}^{\on{spec,loc}}_{\cG,\cP^-,\ul{x}})$

\medskip

For the purposes of this paper, we will need a particular property of \eqref{e:self-duality I G P spec}, stated as \propref{p:semiinf spec IC dual}
below. 

\sssec{} \label{sss:Delta spec co to Delta spec}

Let $\Delta^{-,\on{spec},\semiinf}_{\on{co}}$ be the factorization unit of $\on{I}(\cG,\cP^-)^{\on{spec,loc}}_{\on{co}}$. 
Explicitly, 
$$\Delta^{-,\on{spec},\semiinf}_{\on{co}}=\iota^\IndCoh_*(\omega_{\LS^\reg_{\cP^-}}),$$
where:

\begin{itemize}

\item The morphism $\iota$ is as in \eqref{e:iota basic};

\medskip

\item The functor 
$$\IndCoh^!(\LS^\reg_{\cP^-})\overset{\iota^\IndCoh_*}\longrightarrow \IndCoh^!(\on{Hecke}^{\on{spec,loc}}_{\cG,\cP^-})=
\on{I}(\cG,\cP^-)^{\on{spec,loc}}_{\on{co}}$$
is the dual of the functor $\iota^!$ from \eqref{e:arcs to loops semiinf}, 
and so is the left adjoint of the tautologically defined functor
$$\on{I}(\cG,\cP^-)^{\on{spec,loc}}_{\on{co}}=\IndCoh^!(\on{Hecke}^{\on{spec,loc}}_{\cG,\cP^-})\overset{\iota^!}\to
\IndCoh^!(\LS^\reg_{\cP^-}).$$

\end{itemize}

Since \eqref{e:self-duality I G P spec} is an equivalence of unital factorization categories, we obtain that the functor
\eqref{e:self-duality I G P spec} send 
$$\Delta^{-,\on{spec},\semiinf}_{\on{co}}\in \on{I}(\cG,\cP^-)^{\on{spec,loc}}_{\on{co}}\,\, \mapsto
\Delta^{-,\on{spec},\semiinf}\in \on{I}(\cG,\cP^-)^{\on{spec,loc}}.$$

\sssec{} \label{sss:semiinf IC spec co}

Let 
$$\IC^{-,\on{spec},\semiinf}_{\on{co}}\in \on{I}(\cG,\cP^-)^{\on{spec,loc}}_{\on{co}}$$
be the factorization algebra object equal to the image of 
$$\omega_{\LS^\reg_\cM}\in \IndCoh^!(\LS^\reg_\cM)$$ 
along the functor
$$\IndCoh^!(\LS^\reg_\cM)\to \IndCoh^!(\LS^\reg_{\cP^-})\overset{\iota^\IndCoh_*}\longrightarrow \IndCoh^!(\on{Hecke}^{\on{spec,loc}}_{\cG,\cP^-}),$$
where:

\begin{itemize}

\item The first arrow is the dual of
$$\IndCoh^*(\LS^\reg_{\cP^-})\simeq \QCoh(\LS^\reg_{\cP^-})\to \QCoh(\LS^\reg_\cM)\simeq \IndCoh^*(\LS^\reg_\cM),$$
with the middle arrow being the functor of *-pullback along \eqref{e:LS M to LS P};

\medskip

\item The second arrow is the functor from \secref{sss:Delta spec co to Delta spec} above. 

\end{itemize} 

\sssec{}

Note that the (factorization) functor 
$$\on{I}(\cG,\cP^-)^{\on{spec,loc}}\to \Vect$$
corresponding to the (factorization algebra) object $\IC^{-,\on{spec},\semiinf}_{\on{co}}$ is the composition
\begin{multline*}
\on{I}(\cG,\cP^-)^{\on{spec,loc}}=\IndCoh^*(\on{Hecke}^{\on{spec,loc}}_{\cG,\cP^-})\overset{\iota^!}\to
\IndCoh^*(\LS^\reg_{\cP^-})\simeq \\
\simeq \QCoh(\LS^\reg_{\cP^-})\to \QCoh(\LS^\reg_\cM)\overset{\Gamma(\LS^\reg_\cM,-)}\longrightarrow \Vect,
\end{multline*}
where the functor $\QCoh(\LS^\reg_{\cP^-})\to \QCoh(\LS^\reg_\cM)$ is *-pullback with respect to \eqref{e:LS M to LS P}. 

\sssec{}

We claim:

\begin{prop} \label{p:semiinf spec IC dual}
The factorization algebra $\IC^{-,\on{spec},\semiinf}_{\on{co}}$ corresponds under the identification \eqref{e:self-duality I G P spec} to the 
factorization algebra $\IC^{-,\on{spec},\semiinf}$.
\end{prop}

The proof will be given in \secref{sss:semiinf spec IC dual}. 

\ssec{The semi-infinite Casselman-Shalika equivalence}  \label{ss:semiinf CS}

Just like the proof of the (derived) geometric Satake equivalence 
$$\Sph_G\overset{\Sat_G}\simeq \Sph^{\on{spec}}_\cG$$
in \cite{CR} is bootstrapped from the geometric Casselman-Shalika equivalence 
$$\on{CS}_G:\Whit^!(G) \simeq \Rep(\cG),$$
the proof of \thmref{t:semiinf geom Satake} is based on another local Langlands-type equivalence,
\thmref{t:semiinf CS}. 
 
\medskip

We introduce the actors in the latter equivalence in this subsection, and the proof is delegated to
\secref{s:Sith}. 

\sssec{}

Consider the factorization prestack 
\begin{equation} \label{e:ambient spec}
\LS^\Mmf_{\cP^-}:=\LS^\mer_{\cP^-}\underset{\LS^\mer_\cM}\times \LS^\reg_\cM.
\end{equation} 

When we work over a fixed point $\ul{x}\in \Ran$, we have
\begin{equation} \label{e:ambient spec at x}
\LS^\Mmf_{\cP^-,\ul{x}}\simeq
(\cn^-_P/\on{Ad}(\cP^-))^{\times |\ul{x}|}.
\end{equation} 

Hence, the prestack \eqref{e:ambient spec at x} is an algebraic stack locally of finite type, and thus the category
$$\IndCoh(\LS^\Mmf_{\cP^-,\ul{x}})$$
is well-defined. Moreover, we have map
\begin{equation}  \label{e:iota into ambient x}
\iota:\LS^\reg_{\cP^-}\to \LS^\Mmf_{\cP^-},
\end{equation} 
which is a \emph{closed embedding} over $\ul{x}$, and which gives rise to a pair of adjoint functors
\begin{equation}  \label{e:iota into ambient x adj}
\iota^\IndCoh_*:\IndCoh(\LS^\reg_{\cP^-,\ul{x}}) \rightleftarrows \IndCoh(\LS^\Mmf_{\cP^-,\ul{x}}):\iota^!.
\end{equation} 

\medskip

However, we do not know how to make sense of $\IndCoh^*(-)$ or $\IndCoh^!(-)$ of \eqref{e:ambient spec} as a factorization
category. 

\begin{rem} \label{r:no LS N mer}
That said, one can consider the factorization category
\begin{equation} \label{e:ambient spec QCo co}
\QCoh_{\on{co}}(\LS^\Mmf_{\cP^-}),
\end{equation}
see \cite[Sect. B.13.8]{GLC2}.

\medskip

Over a fixed point $\ul{x}\in \Ran$, the functor
$$\Psi_{\LS^\Mmf_{\cP^-,\ul{x}}}:
\IndCoh(\LS^\Mmf_{\cP^-,\ul{x}})\to
\QCoh_{\on{co}}(\LS^\Mmf_{\cP^-,\ul{x}})$$
is an equivalence, since \eqref{e:ambient spec at x} is smooth. 

\medskip

So, it may be that \eqref{e:ambient spec QCo co} is a reasonable candidate for 
$$\IndCoh^*(\LS^\Mmf_{\cP^-}).$$

\medskip

However, for example, the basic thing we do not know is whether with this definition we have an adjunction
$$\iota^\IndCoh_*:\IndCoh(\LS^\reg_{\cP^-}) \rightleftarrows \IndCoh(\LS^\mer_{\cP^-}\underset{\LS^\mer_{\cM}}\times \LS^\reg_\cM):\iota^!$$
of factorization functors, extending the pointwise adjunction \eqref{e:iota into ambient x adj}. 

\end{rem}

\sssec{}  \label{sss:ambient spec semiinf init}

Let
$$(\LS^\Mmf_{\cP^-})^\wedge_\mf$$
be the formal completion of \eqref{e:ambient spec} along \eqref{e:iota into ambient x}. 

\medskip

In \secref{ss:ambient spec semiinf} we will show how to construct 
$$\IndCoh^*((\LS^\Mmf_{\cP^-})^\wedge_\mf)$$
as a factorization category.

\sssec{} \label{sss:functor to Rep M}

The category $\IndCoh^*((\LS^\Mmf_{\cP^-})^\wedge_\mf)$ will come equipped with the following 
identification
\begin{equation} \label{e:Whit IGP conv spec}
\Rep(\cG)\underset{\Sph_\cG^{\on{spec}}}\otimes \on{I}(\cG,\cP^-)^{\on{spec,loc}}\simeq  \IndCoh^*((\LS^\Mmf_{\cP^-})^\wedge_\mf),
\end{equation}
see \secref{sss:ambient via IGP}.

\medskip

With respect to this equivalence, the functor 
$$J^{-,\on{spec,enh}}:\Rep(\cG)\underset{\Sph_\cG^{\on{spec}}}\otimes \on{I}(\cG,\cP^-)^{\on{spec,loc}}\to \Rep(\cM)$$
corresponds to the functor
$$\IndCoh^*((\LS^\Mmf_{\cP^-})^\wedge_\mf)\overset{((\sfq^{-,\on{spec}})^\wedge_\mf)^\IndCoh_*}\longrightarrow \IndCoh^*(\LS_\cM^\reg)\simeq \Rep(\cM),$$
where $(\sfq^{-,\on{spec}})^\wedge_\mf$ denotes the map
$$(\LS^\Mmf_{\cP^-})^\wedge_\mf\to \LS_\cM^\reg.$$

\sssec{}

Combing \eqref{e:Whit IGP conv spec} with \thmref{t:semiinf geom Satake} we obtain:

\begin{thm} \label{t:semiinf CS}
There exists an equivalence of unital factorization categories
\begin{equation} \label{e:ambient}
\on{CS}^{-,\semiinf}_\tau:
\Whit^!(G)\underset{\Sph_G}\otimes \on{I}(G,P^-)^{\on{loc}}_{\rho_P(\omega_X)} \simeq \IndCoh^*((\LS^\Mmf_{\cP^-})^\wedge_\mf),
\end{equation}
which makes the following diagram commute
\begin{equation} \label{e:semiinf CS diag}
\CD
\Whit^!(G)\underset{\Sph_G}\otimes \on{I}(G,P^-)^{\on{loc}}_{\rho_P(\omega_X)}  @>{\on{CS}^{-,\semiinf}_\tau}>> \IndCoh^*((\LS^\Mmf_{\cP^-})^\wedge_\mf) \\
@V{J^{-,\on{enh}}_{\Whit}}VV @VV{((\sfq^{-,\on{spec}})^\wedge_\mf)^\IndCoh_*}V \\
\Whit^!(M) @>{\on{CS}_{M,\tau}}>> \Rep(\cM).
\endCD
\end{equation}
Furthermore, the functor $\on{CS}^{-,\semiinf}_\tau$ is compatible with the actions of
$$\Sph_M \overset{\Sat_{M,\tau}}\simeq \Sph_\cM^{\on{spec}}$$
on the two sides. 
\end{thm}

We view \thmref{t:semiinf CS} as a semi-infinite version of the Casselman-Shalika equivalence. 

\begin{rem}
As far as the actual logic of the proofs is concerned, we will prove \thmref{t:semiinf CS} first, and then use
it to deduce \thmref{t:semiinf geom Satake}.
\end{rem}

%
%
%
%
%

\section{Jacquet and Wakimoto functors for Kac-Moody modules} \label{s:KM}

In this subsection we introduce generalized Jacquet functors between categories of
Kac-Moody representations for $G$ and $M$. The input for the construction is
the functor of BRST reduction with respect to the loops Lie algebra into $\fn^-_P$. 

\medskip

We also study the \emph{dual} functors, which go from Kac-Moody representations of $M$
to those of $G$, called the \emph{Wakimoto} functors. 

\medskip

We import the definitions pertaining to Kac-Moody modules and the Kazhdan-Lusztig category
from \cite[Sect. 2]{GLC2}.

\medskip

For most of this section we will work at an arbitrary level $\kappa$, but at some point we will
specialize to the case $\kappa=\crit_G$. 

\ssec{Twisting by \texorpdfstring{$Z^0_\cG$}{ZcG}-torsors}  \label{ss:twisting by dual torsors}

In this subsection we will introduce a new kind of twist, specific to Kac-Moody extensions:
it turns out that one can twist the Kac-Moody extension for $G$ by a torsor with respect to the
the (connected) center of $\cG$. 

\medskip

In the bulk of the paper, the observations from this subsection will be applied when the reductive group
in question is the Levi subgroup $M$ of the original $G$. The specific twist that we will need is the Langlands
dual of the \emph{Miura shift} for opers, see \secref{sss:Miura shift}. 

\sssec{} \label{sss:twist by ZcG}

Let $[\fg,\fg]$ be the Lie algebra of the derived group of $G$, so that $\fg_{\on{ab}}:=\fg/[\fg,\fg]$ is the cocenter of $\fg$. 

\medskip

Consider the vector space $\fL(\fg_{\on{ab}})/\fL^+(\fg_{\on{ab}})$. Its dual, viewed as a group, acts by automorphisms 
of the Kac-Moody extension, and hence also by automorphisms of the categories $\hg\mod_\kappa$ and $\KL(G)_\kappa$. 

\medskip

In particular, given a $(\fL(\fg_{\on{ab}})/\fL^+(\fg_{\on{ab}}))^*$-torsor $\CP$, we can form a twist $\hg_{\kappa+\CP}$
of $\hg$ and consider the twisted versions of
$\hg\mod_\kappa$ and $\KL(G)_\kappa$, denoted 
$$\hg\mod_{\kappa+\CP} \text{ and } \KL(G)_{\kappa+\CP},$$
respectively. 

\sssec{}

Note that we have a canonical duality
$$\fg_{\on{ab}}\simeq (z_\cg)^*.$$

In particular, we can identify 
$$(\fL(\fg_{\on{ab}})/\fL^+(\fg_{\on{ab}}))^*\simeq \fL^+(z_\cg\otimes \omega_X),$$
where the right-hand side is the space of sections of $z_\cg$-valued 1-forms on the formal
disc. 

\sssec{} \label{sss:create g ab tors}

Let $Z^0_\cG$ denote the neutral connected component of the center of $\cG$. Consider the homomorphism
$$\on{dlog}:\fL^+(Z^0_\cG)\to \fL^+(z_\cg\otimes \omega_X).$$

Thus, starting from a $Z^0_\cG$-torsor $\CP_{Z^0_\cG}$ on $X$, we can:

\begin{itemize}

\item Produce a $\fL^+(Z^0_\cG)$-torsor (by restricting to the formal disc);

\item Induce a $\fL^+(z_\cg\otimes \omega_X)$-torsor using $\on{dlog}$;

\item Think of the latter a $(\fL(\fg_{\on{ab}})/\fL^+(\fg_{\on{ab}}))^*$-torsor, which by a slight abuse of notation we denote by
$$\on{dlog}(\CP_{Z^0_\cG}).$$

\end{itemize}

The above construction is naturally compatible with factorization.

\medskip

We denote the resulting (factorization) categories of Kac-Moody modules by
$$\hg\mod_{\kappa+\on{dlog}(\CP_{Z^0_\cG})} \text{ and } \KL(G)_{\kappa+\on{dlog}(\CP_{Z^0_\cG})},$$
respectively. 

\sssec{} 

Denote $G_{\on{ab}}=G/[G,G]$, so that it is a torus dual to $Z^0_\cG$.
The Contou-Carr\`ere symbol defines a bilinear pairing
\begin{equation} \label{e:ContCar}
\langle-,-\rangle_{\on{CtCr}}:\fL(Z^0_\cG)\times \fL(G_{\on{ab}})\to \BG_m
\end{equation} 

We normalize \eqref{e:ContCar} so that for an element $g\in \fL(Z^0_\cG)$, the resulting homomorphism
$$\langle g,-\rangle_{\on{CtCr}}: \fL(G_{\on{ab}})\to \BG_m$$
is such that its differential, which is an element of
$$\fL(\fg_{\on{ab}})^*\simeq \fL(z_\cg\otimes \omega_X)$$
is given by $\on{dlog}(g)\in \fL(z_\cg\otimes \omega_X)$. 

\sssec{} \label{sss:ZGc torsor to line bundle}

Restricting $\langle-,-\rangle_{\on{CtCr}}$ along 
$$\fL^+(Z^0_\cG)\times \fL(G)\to \fL^+(Z^0_\cG)\times \fL(G_{\on{ab}}),$$
we obtain that a $Z^0_\cG$-torsor $\CP_{Z^0_\cG}$ on $X$ gives rise to a multiplicative 
line bundle $\CL^{\on{loc}}_{\fL(G),\CP_{Z^0_\cG}}$ on $\fL(G)$, i.e., 
a central extension by means of $\BG_m$, compatible with
the factorization structure.

\medskip

The induced central extension of $\fL(\fg)$ is $\hg\mod_{\on{dlog}(\CP_{Z^0_\cG})}$.

\sssec{} \label{sss:change level}

Hence, tensoring by $\CL^{\on{loc}}_{\fL(G),\CP_{Z^0_\cG}}$ defines a bijection 
between categories acted on (strongly) by $\fL(G)$ at levels $\kappa$ and $\kappa+\on{dlog}(\CP_{Z^0_\cG})$,
respectively.

\medskip

In particular, we obtain an action of $\fL(G)$ at level $\kappa$ on the category
$\hg\mod_{\kappa+\on{dlog}(\CP_{Z^0_\cG})}$, and we have the convolution functor
$$\Dmod_\kappa(\Gr_G)\otimes \KL(G)_{\kappa+\on{dlog}(\CP_{Z^0_\cG})}\overset{\star}\to \hg\mod_{\kappa+\on{dlog}(\CP_{Z^0_\cG})}.$$

\sssec{Notational convention} \label{sss:lambdach twist}

Assume for a moment that the above $Z^0_\cG$-torsor is of the form 
$\lambdach(\omega_X)$
where $\lambdach:\BG_m\to Z^0_\cG$.

\medskip

In this case, we will use a shorthandnotation
$$\hg\mod_{\kappa+\lambdach}:=\hg\mod_{\kappa+\on{dlog}(\lambdach(\omega_X))}
\text{ and } \KL(G)_{\kappa+\lambdach}:=\KL(G)_{\kappa+\on{dlog}(\lambdach(\omega_X))}$$
for the corresponding twisted categories. 

\medskip

Note that the assignment
$$\lambdach\rightsquigarrow \on{dlog}(\lambdach(\omega_X))$$
is linear. So the resulting torsor with respect to 
$$\fL^+(z_\cg\otimes \omega_X)\simeq (\fL(\fg_{\on{ab}})/\fL^+(\fg_{\on{ab}}))^*$$
makes sense for any $\lambdach\in z_\cg$. 

\medskip

Thus, the twisted categories 
$$\hg\mod_{\kappa+\lambdach} \text{ and } \KL(G)_{\kappa+\lambdach}$$
are defined for any $\lambdach\in z_\cg$. 

\ssec{Jacquet functor(s) for Kac-Moody representations} \label{ss:BRST}

In this subsection we recall the construction of the (usual) BRST reduction functor from Kac-Moody 
representations of $G$ to those of $M$. 

\sssec{} \label{sss:level for m}

Note that that $\kappa|_{\fp^-}$ factors through $\fp^-\to \fm$. We will use the same character $\kappa$ to denote the resulting
level for $\fm$. 

\medskip

Note that the central extensions of $\fL(\fp^-)$ induced from $\hg_\kappa$ and $\hm_\kappa$ via the maps
$$\fL(\fg) \leftarrow \fL(\fp^-)\to \fL(\fm),$$
respectively, are canonically isomorphic.

\medskip

We denote the resulting central extension of $\fL(\fp^-)$ by $\hp^-_\kappa$. 

\sssec{} \label{sss:level for p}

It is a basic fact (see \cite[Proof of Theorem 2.8.17, point (c) on page 142]{BD2}) that the Tate extension of $\fL(\fp^-)$ identifies with 
$$\hp^-_{-\crit_G-\crit_M+\rhoch_P},$$
where:

\begin{itemize}

\item The subscripts $\crit_G$ and $\crit_M$ are the critical levels for $G$ and $M$, respectively.
 
\medskip

\item The subscript $\rhoch_P$ is as in \secref{sss:lambdach twist}.

\end{itemize}

\medskip

In particular, for any $\kappa$ we have a canonical duality
\begin{equation} \label{e:duality KL P}
(\hp^-\mod_{\crit_G+\kappa})^\vee \simeq \hp^-\mod_{\crit_M-\kappa-\rhoch_P}
\end{equation} 
as unital factorization categories, compatible with the $\fL(G)$-action, where:

\smallskip

\begin{itemize}

\item On the left-hand side, $\fL(P^-)$ acts at level $-(\crit_G+\kappa)$ restricted along $P^-\to G$;

\smallskip

\item On the right-hand side, $\fL(P^-)$ acts at level $\crit_M-\kappa-\rhoch_P$ restricted along $P^-\to M$;

\smallskip

\item We pass between actions at these levels as in \secref{sss:change level}.

\end{itemize}

\sssec{}

The projection
$$\hp^-\mod_{\crit_M-\kappa-\rhoch_P}\to \hm\mod_{\crit_M-\kappa-\rhoch_P}$$
gives rise to a (strictly) unital factorization functor
$$\Res^{\hm}_{\hp^-}:\hm\mod_{\crit_M-\kappa-\rhoch_P}\to \hp^-\mod_{\crit_M-\kappa-\rhoch_P}.$$

Its dual, which is a functor
\begin{equation} \label{e:BRST from P}
\BRST^-:\hp^-\mod_{\crit_G+\kappa}\to \hm\mod_{\crit_M+\kappa+\rhoch_P},
\end{equation} 
is called the functor of BRST reduction with respect to $\fL(\fn^-_P)$. 

\medskip

Being the dual of a (strictly) unital factorization functor, the functor $\BRST^-$ acquires
a lax unital factorization structure.

\sssec{}

Note that the functor $\Res^{\hm}_{\hp^-}$ respects the actions of $\fL(P^-)$ on the two sides. Hence, the functor
$\BRST^-$ acquires a similar structure, where: 

\begin{itemize}

\item On the left-hand side, $\fL(P^-)$ acts at level $\crit_G+\kappa$ restricted along $P^-\to G$;

\item On the left-hand side, $\fL(P^-)$ acts at level $\crit_M+\kappa+\rhoch_P$ restricted along $P^-\to M$;

\item We pass between actions at these levels as in \secref{sss:change level}.

\end{itemize}

\medskip

In particular, the functor $\BRST^-$ is $\fL(N^-_P)$-invariant, and hence canonically factors as
\begin{equation} \label{e:BRST from P coinv}
\hp^-\mod_{\crit_G+\kappa} \twoheadrightarrow (\hp^-\mod_{\crit_G+\kappa})_{\fL(N^-_P)}\overset{\ol\BRST^-}\longrightarrow
\hm\mod_{\crit_M+\kappa+\rhoch_P},
\end{equation}
where $\ol\BRST^-$ acquires a $\fL(M)$-equivariant structure. 

\medskip

In particular, the functor $\BRST^-$ induces a functor
$$\hp^-\mod_{\crit_G+\kappa}^{\fL^+(P^-)}\to \hm\mod^{\fL^+(P^-)}_{\crit_M+\kappa+\rhoch_P}\simeq \hm\mod^{\fL^+(M)}_{\crit_M+\kappa+\rhoch_P}=
\KL(M)_{\crit_M+\kappa+\rhoch_P},$$
denoted by the same symbol $\BRST^-$,
and the functor $\ol\BRST^-$ induces a functor
$$(\hp^-\mod_{\crit_G+\kappa})^{\fL^+(M)}_{\fL(N^-_P)}\to \hm\mod^{\fL^+(M)}_{\crit_M+\kappa+\rhoch_P}=
\KL(M)_{\crit_M+\kappa+\rhoch_P},$$
denoted by the same symbol $\ol\BRST^-$. 

\sssec{}

We let 
$$J^-_{\on{KM}}:\hg\mod_{\crit_G+\kappa}\to \hm\mod_{\crit_M+\kappa+\rhoch_P}$$
denote the composition 
$$\hg\mod_{\crit_G+\kappa}\to \hp^-\mod_{\crit_G+\kappa} \overset{\BRST^-}\to \hm\mod_{\crit_M+\kappa+\rhoch_P}.$$

\medskip

We will denote by the same symbol $J^-_{\on{KM}}$ the induced functor
$$\hg\mod_{\crit_G+\kappa}^{\fL^+(P^-)}\to \hm\mod^{\fL^+(P^-)}_{\crit_M+\kappa+\rhoch_P}\simeq \hm\mod^{\fL^+(M)}_{\crit_M+\kappa+\rhoch_P}=
\KL(M)_{\crit_M+\kappa+\rhoch_P}.$$

\sssec{}

Let $$\ol{J}^-_{\on{KM}}:(\hg\mod_{\crit_G+\kappa})_{\fL(N^-_P)}\to \hm\mod_{\crit_M+\kappa+\rhoch_P}$$
denote the composition 
$$(\hg\mod_{\crit_G+\kappa})_{\fL(N^-_P)}\to (\hp^-\mod_{\crit_G+\kappa})_{\fL(N^-_P)} \overset{\ol\BRST^-}\to \hm\mod_{\crit_M+\kappa+\rhoch_P}.$$

We will denote by the same symbol $\ol{J}^-_{\on{KM}}$ the induced functor
$$(\hg\mod_{\crit_G+\kappa})^{\fL^+(M)}_{\fL(N^-_P)}\to \hm\mod^{\fL^+(M)}_{\crit_M+\kappa+\rhoch_P}=
\KL(M)_{\crit_M+\kappa+\rhoch_P}.$$

\sssec{}

We will denote by $J_{\on{KM}}^{-,\Sph}$ the composition
\begin{equation} \label{e:J KM Sph}
\KL(G)_{\crit_G+\kappa}=\hg\mod_{\crit_G+\kappa}^{\fL^+(G)}\overset{\oblv_{\fL^+(G)\to \fL^+(P^-)}}\longrightarrow 
\hg\mod_{\crit_G+\kappa}^{\fL^+(P^-)} \overset{J^-_{\on{KM}}}\to \KL(M)_{\crit_M+\kappa+\rhoch_P}.
\end{equation} 

\ssec{The enhanced Jacquet functor for Kac-Moody representations}

In this subsection we couple the Jacquet functor $\ol{J}^-_{\on{KM}}$ from the previous subsection with objects
$\on{I}(G,P^-)^{\on{loc}}_{\on{co}}$ to obtain generalized Jacquet functors that map $\KL(G)_{\crit_G}$ to
$\KL(M)_{\crit_M}$.

\sssec{} \label{sss:J enh}

We now specialize to the critical level. In this case we will 
introduce one more version of the Jacquet functor, which will play a fundamental role in the sequel. 

\medskip

Note that we have a fully faithful functor
$$\KL(G)_{\crit_G}\underset{\Sph_G}\otimes \Dmod_{\frac{1}{2}}(\Gr_G)_{\fL(N^-_P)}^{\fL^+(M)}\hookrightarrow 
(\hg\mod_{\crit_G})^{\fL^+(M)}_{\fL(N^-_P)}.$$

Pre-composing with 
$$\on{I}(G,P^-)^{\on{loc}}_{\on{co}}:=
\Dmod_{\frac{1}{2}}(\Gr_G)_{\fL(N^-_P)}^{\fL^+(M),\on{ren}}\to 
\Dmod_{\frac{1}{2}}(\Gr_G)_{\fL(N^-_P)}^{\fL^+(M)}$$
and composing with $\ol{J}^-_{\on{KM}}$, 
we obtain a functor
\begin{equation} \label{e:Jacquet KM enh co}
\KL(G)^{-,\on{enh}_{\on{co}}}_{\crit_G}:=
\KL(G)_{\crit_G}\underset{\Sph_G}\otimes \on{I}(G,P^-)^{\on{loc}}_{\on{co}}\to \KL(M)_{\crit_M+\rhoch_P},
\end{equation}
which we denote by
$$J^{-,\on{enh}_{\on{co}}}_{\on{KM}}.$$

\sssec{}

The functor $J^{-,\on{enh}_{\on{co}}}_{\on{KM}}$ is $\Sph_M$-linear. Hence, by \secref{sss:duality over Sph}, 
the datum of $J^{-,\on{enh}_{\on{co}}}_{\on{KM}}$ is equivalent to the datum of
a $\Sph_G$-linear functor
\begin{equation} \label{e:Jacquet KM enh}
\KL(G)_{\crit_G}\to \KL(M)_{\crit_M+\rhoch_P}\underset{\Sph_M}\otimes \on{I}(G,P^-)^{\on{loc}}=:\KL(M)^{-,\on{enh}}_{\crit_M+\rhoch_P}.
\end{equation}

We denote the functor \eqref{e:Jacquet KM enh} by
$$J^{-,\on{enh}}_{\on{KM}}.$$

\sssec{}

Note that the functor $J_{\on{KM}}^{-,\on{Sph}}$ of \eqref{e:J KM Sph} is recovered from $J^{-,\on{enh}}_{\on{KM}}$ as
$$\KL(G)_{\crit_G} \overset{J^{-,\on{enh}}_{\on{KM}}}\longrightarrow 
\KL(M)^{-,\on{enh}}_{\crit_M+\rhoch_P}\overset{\oblv_{\on{enh}}}\longrightarrow \KL(M)_{\crit_M+\rhoch_P},$$
where $\oblv_{\on{enh}}$ is the functor 
$$\KL(M)_{\crit_M+\rhoch_P}\underset{\Sph_M}\otimes \on{I}(G,P^-)^{\on{loc}} \overset{\on{Id}\otimes \oblv_{\semiinf\to \Sph}}\longrightarrow \\
\to \KL(M)_{\crit_M+\rhoch_P}\underset{\Sph_M}\otimes \Sph_M= \KL(M)_{\crit_M+\rhoch_P}.$$

\ssec{The Wakimoto functor(s)}

In this subsection we define the \emph{Wakimoto} functor, which is a functor dual to $J^-_{\on{KM}}$. 

\medskip

Note, however,
that our definition of the Wakimoto functor (or the values of this functor on particular objects) is different from what is 
called \emph{Wakimoto modules} in most of the literature: our Wakimoto modules are semi-infinite in nature. The
usual Wakimoto modules can be obtained from ours by a subsequent averaging procedure, see \secref{ss:usual Wak}. 

\sssec{} \label{sss:semiinf KM}

Let 
\begin{equation} \label{e:original Wak}
\Wak^{-,\semiinf}:\hm\mod_{\crit_M-\kappa-\rhoch_P}\to \hg\mod_{\crit_G-\kappa}
\end{equation}
be the functor dual to the functor $J^-_{\on{KM}}$. 

\medskip

Since $J^-_{\on{KM}}$ factors as 
$$\hg\mod_{\crit_G+\kappa}\to (\hg\mod_{\crit_G+\kappa})_{\fL(N^-_P)}\overset{\ol{J}^-_{\on{KM}}}\longrightarrow \hm\mod_{\crit_M+\kappa+\rhoch_P},$$
the functor $\Wak^{-,\semiinf}$ takes values in 
$$(\hg\mod_{\crit_G-\kappa})^{\fL(N^-_P)}\subset \hg\mod_{\crit_G-\kappa}.$$

\medskip

The resulting functor 
$$\hm\mod_{\crit_M-\kappa-\rhoch_P}\to (\hg\mod_{\crit_G-\kappa})^{\fL(N^-_P)}$$
is the dual of the functor $\ol{J}^-_{\on{KM}}$.

\medskip

We will keep the same symbol $\Wak^{-,\semiinf}$ for the functor 
\begin{equation} \label{e:semiinf Wak}
\KL(M)_{\crit_M-\kappa-\rhoch_P}:=\hm\mod_{\crit_M-\kappa-\rhoch_P}^{\fL^+(M)} \to (\hg\mod_{\crit_G-\kappa})^{\fL(N^-_P)\cdot \fL^+(M)}
=:\hg\mod_{\crit_G-\kappa}^{-,\semiinf}
\end{equation}

\begin{rem}

Recall (see \cite[Sect. B.14]{GLC2}) that $\hg\mod_{\crit_G+\kappa}$ has a natural structure of (unital) factorization 
category, and this structure in inherited by $(\hg\mod_{\crit_G+\kappa})_{\fL(N^-_P)}$. 

\medskip

However, we do not know whether the operation of $\fL(N^-_P)$-\emph{invariants} commutes with colimits or
tensor products. Hence, a priori 
$$\hg\mod_{\crit_G-\kappa}^{\fL^+(N^-_P)} \text{ and } \hg\mod_{\crit_G-\kappa}^{-,\semiinf}$$
are only \emph{lax} factorization categories (see \cite[Sect. B.11.12]{GLC2} for what this means), equipped with a unital
structure.\footnote{For a unital factorization category $\bC$, equipped with an action of $\fL(G)$, the category $\bC^{\fL^+(N^-_P)}$
acquires a unital structure, in which the factorization unit is $\Delta^{-,\semiinf}\star \one_{\bC}$.} 

\medskip

The functor $\Wak^{-,\semiinf}$ (in both variants) has a natural structure of (lax unital) factorization functor. 

\end{rem} 

\begin{rem}

The version of the Wakimoto functor given by \eqref{e:semiinf Wak} may seem somewhat exotic. For example, it produces objects
that are not seen by the t-structure on $\hg\mod_{\crit_G-\kappa}$ (i.e., all of their cohomologies are $0$). 

\medskip

One recovers from it the usual Wakimoto functor by composing $\Wak^{-,\semiinf}$
with the functor
$$\hg\mod_{\crit_G-\kappa}\overset{\Av_*^{\fL^+(N)}}\longrightarrow \hg\mod_{\crit_G-\kappa}^{\fL^+(N)}\hookrightarrow
\hg\mod_{\crit_G-\kappa},$$
see \secref{sss:usual Wak}.

\end{rem}

\sssec{}

Denote:
\begin{equation} \label{e:Sph Wak}
\Wak^{-,\on{Sph}}:=\on{Av}_*^{\fL^+(G)/\fL^+(M)}\circ \Wak^{-,\semiinf}, \quad
\KL(M)_{\crit_M-\kappa-\rhoch_P}\to \KL(G)_{\crit_G-\kappa}.
\end{equation}

This is the functor dual to $J_{\on{KM}}^{-,\on{Sph}}$. It has a natural (lax unital) factorization structure. 

\sssec{} \label{sss:Wak Sph via CDO}

The functor \eqref{e:Sph Wak} can be explicitly described in terms of the duality
\begin{equation} \label{e:duality shifted KLM}
(\KL(M)_{\crit_M-\kappa-\rhoch_P})^\vee\simeq \KL(M)_{\crit_G+\kappa-\rhoch_P}.
\end{equation} 

Namely, the corresponding object of
$$(\KL(M)_{\crit_M-\kappa-\rhoch_P})^\vee\otimes \KL(G)_{\crit_G-\kappa}\simeq
\KL(M)_{\crit_G+\kappa-\rhoch_P} \otimes \KL(G)_{\crit_G-\kappa}$$
is given by
$$(J^{-,\Sph}_{\on{KM}}\otimes \on{Id})(\mathfrak{CDO}(G)_{\crit_G+\kappa,\crit_G-\kappa}),$$
where 
$$\mathfrak{CDO}(G)_{\crit_G+\kappa,\crit_G-\kappa}\in  \KL(G)_{\crit_G+\kappa}\otimes  \KL(G)_{\crit_G-\kappa}$$
is as in \cite[Sect. 2.2.3]{GLC2}. 

\ssec{Restriction of Wakimoto modules}

\sssec{} 

We now consider the relationship between the Wakimoto construction and
restriction to the opposite Borel. 

\medskip 

\newcommand{\hb}{\widehat{\fb}}

Let $\fb(M) = \fp^- \cap \fb$ (resp., $\fn(M) = \fp^- \cap \fn$) be the Lie algebra of the positive
Borel (resp., its unipotent radical) in $M$. 

\medskip 

Note that as in \secref{sss:level for p}, the Tate extension for $\fL(\fb(M))$ is 
$-\crit_M-\rhoch_M$ while the Tate extension for $\fL(\fb)$ itself is
$-\crit_G-\rhoch_G$, so the difference between these two Tate extensions
is the same as the difference between the Tate extensions for 
$\fL(\fp)$ and $\fL(\fg)$. Therefore, the following result makes sense:

\begin{prop}\label{p:wak-explicit, almost}

The following diagram of (lax unital) factorization functors commutes:
\[
\begin{tikzcd}[row sep = large, column sep = scriptsize]
\hp^-\mod_{\crit_M+\kappa-\rhoch_P}
\arrow[rr,shorten <= 1ex, shorten >= 1ex, "\oblv^{\vee}"]
\arrow[d,"\oblv"]
&&
\hg\mod_{\crit_G+\kappa}
\arrow[d,"\oblv"]
\\
\hb(M)\mod_{\crit_M+\kappa-\rhoch_P}
\arrow[rr,shorten <= 1ex, shorten >= 1ex,"\oblv^{\vee}"]
&&
\hb\mod_{\crit_G+\kappa}
\end{tikzcd}
\]

\noindent Here the functors $\oblv^{\vee}$ are dual to forgetful functors 
$\oblv$ under the equivalences
\[
\begin{gathered}
(\hp^-\mod_{\crit_M+\kappa-\rhoch_P})^{\vee} \simeq 
\hp^-\mod_{\crit_G-\kappa}, \quad 
(\hg\mod_{\crit_G+\kappa})^{\vee} \simeq 
\hg\mod_{\crit_G-\kappa},
\\
(\hb(M)\mod_{\crit_M+\kappa-\rhoch_P})^{\vee} \simeq
\hb(M)\mod_{-\kappa+\rhoch_G}, \quad 
(\hb\mod_{\crit_G+\kappa})^{\vee} \simeq 
\hb\mod_{-\kappa+\rhoch_G}
\end{gathered}
\]

\noindent given by semi-infinite cohomology.

\end{prop}

\begin{rem}
Unwinding the abstraction, \propref{p:wak-explicit, almost} amounts
to the identity
\begin{equation} \label{e:triang g}
\Gamma^{\IndCoh}(\fL(G),\delta_1) \simeq 
\Gamma^{\IndCoh}(\fL(P^-),\delta_1) \underset{\hb(M)}{\overset{\semiinf}{\otimes}}
\Gamma^{\IndCoh}(\fL(B),\delta_1),
\end{equation}
where $\Gamma^{\IndCoh}(\fL(H),\delta_1)$ is the unit of the duality between $\wh\fh\mod_{?}$ and $\wh\fh\mod_{??}$
for $H=G,P^-,B,B(M)$ and $?$ and $??$ are appropriate levels. 

\medskip

So, \eqref{e:triang g} expresses the triangular decomposition of $\fg$.

\end{rem}

\begin{rem} \label{r:semiinf restr gen}

A similar result holds for any parabolic $\fp_2$ in place of the Borel,
or more generally, any subalgebra $\fh \subset \fg$ with 
the property that $\fp^-+\fh = \fg$. In particular, we can apply it to $\fh=\fn$, see \lemref{l:wak restr to n}
below. 

\end{rem}

\sssec{}

As a formal consequence of \propref{p:wak-explicit, almost}, we obtain:

\begin{cor}\label{c:wak-explicit}

The following diagram of (lax unital) factorization functors commutes:
\[
\begin{tikzcd}[row sep = large, column sep = scriptsize]
\hm\mod_{\crit_M+\kappa-\rhoch_P}
\arrow[rr,shorten <= 1ex, shorten >= 1ex, "\Wak^{-,\semiinf}"]
\arrow[d,"\oblv"]
&&
\hg\mod_{\crit_G+\kappa}
\arrow[d,"\oblv"]
\\
\hb(M)\mod_{\crit_M+\kappa-\rhoch_P}
\arrow[rr,shorten <= 1ex, shorten >= 1ex,"\oblv^{\vee}"]
&&
\hb\mod_{\crit_G+\kappa}
\end{tikzcd}
\]

\end{cor}

\sssec{}

For future reference, we note (see Remark \ref{r:semiinf restr gen}): 

\begin{lem} \label{l:wak restr to n} 
The following diagram of (lax unital) factorization commutes
\[
\begin{tikzcd}[row sep = large, column sep = scriptsize]
\hp^-\mod_{\crit_M+\kappa-\rhoch_P}
\arrow[rr,shorten <= 1ex, shorten >= 1ex, "\oblv^{\vee}"]
\arrow[d,"\oblv"]
&&
\hg\mod_{\crit_G+\kappa}
\arrow[d,"\oblv"]
\\
\fL(\fn(M))\mod
\arrow[rr,shorten <= 1ex, shorten >= 1ex,"\oblv^{\vee}"]
&&
\fL(\fn)\mod,
\end{tikzcd}
\]
where the lower horizontal arrow is the dual of the natural restriction functor. 
\end{lem}

\begin{cor} \label{c:wak restr to n} 
The following diagram of (lax unital) factorization functors commutes:
\[
\begin{tikzcd}[row sep = large, column sep = scriptsize]
\hm\mod_{\crit_M+\kappa-\rhoch_P}
\arrow[rr,shorten <= 1ex, shorten >= 1ex, "\Wak^{-,\semiinf}"]
\arrow[d,"\oblv"]
&&
\hg\mod_{\crit_G+\kappa}
\arrow[d,"\oblv"]
\\
\fL(\fn(M))\mod
\arrow[rr,shorten <= 1ex, shorten >= 1ex,"\oblv^{\vee}"]
&&
\fL(\fn)\mod,
\end{tikzcd}
\]
where the lower horizontal arrow is the dual of the natural restriction functor. 
\end{cor}

\ssec{The ``usual" Wakimoto functor} \label{ss:usual Wak}

\sssec{} \label{sss:usual Wak}

The ``usual" Wakimoto functor, to be denoted
$$\Wak:\hm\mod_{\crit_M+\kappa-\rhoch_P}\to \hg\mod_{\crit_G+\kappa}$$
is defined as the composition
$$\hm\mod_{\crit_M+\kappa-\rhoch_P}\overset{\Wak^{-,\semiinf}}\longrightarrow (\hg\mod_{\crit_G+\kappa})^{\fL(N^-_P)}\to
\hg\mod_{\crit_G+\kappa}\overset{\Av^{\fL^+(N_P)}_*}\longrightarrow (\hg\mod_{\crit_G+\kappa})^{\fL^+(N_P)}.$$

This is the version of the Wakimoto functor that was studied in the early papers of Feigin and Frenkel (see \cite[Sect. 2.4]{Ga2}).

\medskip

The functor $\Wak$ induces a functor
\begin{equation} \label{e:Wak on KLM}
\KL(M)_{\crit_M+\kappa-\rhoch_P}\to (\hg\mod_{\crit_G+\kappa})^{\fL^+(P)},
\end{equation} 
which we will denote by the same character $\Wak$. 

\sssec{} \label{sss:usual Wak via CDO}

Here is an explicit description of the functor $\Wak$ in the spirit of \secref{sss:Wak Sph via CDO}. 

\medskip

Let $\oG:=P\cdot P^-\subset G$ be the big Bruhat cell. We can consider the factorization algebra
$$\mathfrak{CDO}(\oG)_{\crit_G+\kappa,\crit_G-\kappa}\in  \hg\mod_{\crit_G+\kappa}\otimes  \hg\mod_{\crit_G-\kappa}.$$

Consider the factorization algebra 
\begin{equation} \label{e:ker for Wak}
(\on{Id}\otimes J^-_{\on{KM}})(\mathfrak{CDO}(\oG)_{\crit_G+\kappa,\crit_G-\kappa})\in (\hg\mod_{\crit_G+\kappa})^{\fL^+(P)}\otimes  \KL(M)_{\crit_G-\kappa+\rhoch_P}.
\end{equation}

Then the object \eqref{e:ker for Wak} corresponds to the functor $\Wak$ of \eqref{e:Wak on KLM} in terms of the duality \eqref{e:duality shifted KLM}.

\sssec{} \label{sss:BW}

We will denote by $\BW_{\fg,\crit_G+\kappa}$ the image of $\Wak(\on{Vac}(M)_{\crit_M+\kappa})$ under the forgetful functor
$$(\hg\mod_{\crit_G+\kappa})^{\fL^+(P)}\to \hg\mod_{\crit_G+\kappa}\to \Vect.$$

This is a \emph{classical} factorization algebra. The unital structure on $\Wak^{-,\semiinf}$ (and hence on $\Wak$)
equips $\BW_{\fg,\crit_G+\kappa}$ with a homomorphism
$$\BV_{\fg,\crit_G+\kappa}\to \BW_{\fg,\crit_G+\kappa}.$$

\sssec{}

The following is immediate from \corref{c:wak-explicit}:

\begin{cor}\label{c:wak-usual-explicit prel}

The following diagram of (lax unital) factorization functors commutes:
\[
\begin{tikzcd}[row sep = large, column sep = large]
\hm\mod_{\crit_M+\kappa-\rhoch_P}
\arrow[rr,shorten <= 1ex, shorten >= 1ex, "\Wak"]
\arrow[d,"\oblv"]
&&
\hg\mod_{\crit_G+\kappa}^{\fL^+(N_P)}
\arrow[d,"\oblv"]
\\
\hb(M)\mod_{\crit_M+\kappa-\rhoch_P}
\arrow[rr]
&&
\hb\mod_{\crit_G+\kappa}^{\fL^+(N_P)},
\end{tikzcd}
\]
where the bottom horizontal arrow is the dual of the forgetful functor. 

\end{cor}

\sssec{}

From \corref{c:wak-usual-explicit prel}, we obtain:

\begin{cor}\label{c:wak-usual-explicit}

The following diagram of (lax unital) factorization functors commutes:
\[
\begin{tikzcd}[row sep = large, column sep = large]
\hm\mod_{\crit_M+\kappa-\rhoch_P}
\arrow[rr,shorten <= 1ex, shorten >= 1ex, "\Wak"]
\arrow[d,"\oblv"]
&&
\hg\mod_{\crit_G+\kappa}^{\fL^+(N_P)}
\arrow[d,"\oblv"]
\\
\hb(M)\mod_{\crit_M+\kappa-\rhoch_P}
\arrow[rr,shorten <= 1ex, shorten >= 1ex,"\oblv(-) \otimes 
\mathfrak{CDO}(N_P)"]
&&
\fL(N_P)\mod^{\fL^+(N_P)},
\end{tikzcd}
\]
where in the bottom horizontal arrow, the functor $\oblv$ is the forgetful functor
$$\hb(M)\mod_{\crit_M+\kappa-\rhoch_P}\to \Vect.$$

\end{cor}

\begin{rem}
\corref{c:wak-usual-explicit} 
expresses the following idea: the functor $\Wak$ produces 
objects that are ``semi-infinite free" with respect to $\fL(\fn_P)$, see \cite[Sect. 2.2]{Ga2}.
\end{rem}

\subsubsection{}

We have the following key result.

\begin{prop}\label{p:wak-exact}

At the pointwise level, the functor 
$\Wak$ is t-exact.

\end{prop}

\begin{rem} \label{r:wak-exact}

Note that \corref{c:wak-usual-explicit} implies that the t-exactness property in \propref{p:wak-exact}
holds after applying the forgetful functor 
$\oblv:\hg\mod_{\crit_G+\kappa}^{\fL^+(N_P)} \to \Vect$.

\medskip

However, since the t-structure on $\hg\mod_{\crit_G+\kappa}^{\fL^+(N_P)}$ is not left-separated,
this is not enough to deduce the assertion of proposition. We fill in the 
details in \secref{ss:wak-exact}.

\end{rem}

\ssec{The enhanced Wakimoto functor}

\sssec{}

We now specialize the critical level. In this case, we will encounter one more version of the Wakimoto functor.

\medskip

Namely, let 
$$\Wak^{-,\on{enh}}:\KL(M)_{\crit_M-\rhoch_P}\to \KL(G)_{\crit_G}\underset{\Sph_G}\otimes \on{I}(G,P^-)^{\on{loc}}$$
be the functor dual to the functor $J^{-,\on{enh}_{\on{co}}}_{\on{KM}}$ of \eqref{e:Jacquet KM enh co},
see \secref{sss:duality over Sph}. 

\sssec{}

Explicitly, the functor $\Wak^{-,\on{enh}}$ is the composition
$$\KL(M)_{\crit_M-\rhoch_P} \overset{\Wak^{-,\semiinf}}\to 
\hg\mod_{\crit_G}^{-,\semiinf}\to \KL(G)_{\crit_G}\underset{\Sph_G}\otimes \on{I}(G,P^-)^{\on{loc}},$$
where the second arrow, i.e.,
\begin{equation} \label{e:KM to Sph gen}
\hg\mod_{\crit_G}^{-,\semiinf}\to \KL(G)_{\crit_G}\underset{\Sph_G}\otimes \on{I}(G,P^-)^{\on{loc}}
\end{equation}
is the functor dual to
\begin{multline} \label{e:Sph gen KM co}
\KL(G)_{\crit_G}\underset{\Sph_G}\otimes \on{I}(G,P^-)^{\on{loc}}_{\on{co}}=
\KL(G)_{\crit_G}\underset{\Sph_G}\otimes \Dmod_{\frac{1}{2}}(\Gr_G)_{\fL(N^-_P)}^{\fL^+(M),\on{ren}}\to  \\
\to \KL(G)_{\crit_G}\underset{\Sph_G}\otimes \Dmod_{\frac{1}{2}}(\Gr_G)_{\fL(N^-_P)}^{\fL^+(M)}\hookrightarrow 
(\hg\mod_{\crit_G})^{\fL^+(M)}_{\fL(N^-_P)}.
\end{multline}

\sssec{}

We record:

\begin{lem} \label{l:KM to Sph gen}
The functors
$$\KL(G)_{\crit_G}\underset{\Sph_G}\otimes \Dmod_{\frac{1}{2}}(\Gr_G)_{\fL(N^-_P)}^{\fL^+(M),\on{ren}}\to  
 \KL(G)_{\crit_G}\underset{\Sph_G}\otimes \Dmod_{\frac{1}{2}}(\Gr_G)_{\fL(N^-_P)}^{\fL^+(M)}$$
and 
$$\KL(G)_{\crit_G}\underset{\Sph_G}\otimes \Dmod_{\frac{1}{2}}(\Gr_G)^{\fL(N^-_P)\cdot \fL^+(M),\on{ren}}\to  
 \KL(G)_{\crit_G}\underset{\Sph_G}\otimes \Dmod_{\frac{1}{2}}(\Gr_G)^{\fL(N^-_P)\cdot \fL^+(M)}$$
are equivalences. In particular, the functors \eqref{e:Sph gen KM co} and 
\begin{multline} \label{e:Sph gen KM}
\KL(G)_{\crit_G}\underset{\Sph_G}\otimes \on{I}(G,P^-)^{\on{loc}}=
\KL(G)_{\crit_G}\underset{\Sph_G}\otimes \Dmod_{\frac{1}{2}}(\Gr_G)^{\fL(N^-_P)\cdot \fL^+(M),\on{ren}}\to  \\
\to \KL(G)_{\crit_G}\underset{\Sph_G}\otimes \Dmod_{\frac{1}{2}}(\Gr_G)^{\fL(N^-_P)\cdot \fL^+(M)}\hookrightarrow 
(\hg\mod_{\crit_G})^{\fL^+(M)\cdot \fL(N^-_P)}=:\hg\mod_{\crit_G}^{-,\semiinf}
\end{multline}
are fully faithful, and the functor \eqref{e:KM to Sph gen} identifies with the right adjoint of \eqref{e:Sph gen KM}.
\end{lem}

\begin{proof}

It is enough to establish the corresponding pointwise assertion with respect to $\Ran$. 
In this case, by (the parabolic version of) \cite[Theorem 6.2.1 and Corollary 6.2.3]{Ra4}, at the pointwise level, we can replace both categories
$$\Dmod_{\frac{1}{2}}(\Gr_G)_{\fL(N^-_P)}^{\fL^+(M)} \text{ and }\Dmod_{\frac{1}{2}}(\Gr_G)^{\fL(N^-_P)\cdot \fL^+(M)}$$
by 
$$\Dmod_{\frac{1}{2}}(\Gr_G)^{I^-_P},$$
where $I^-_P$ is  the corresponding (negative)
parahoric subgroup.

\medskip

Similarly, 
$$\Dmod_{\frac{1}{2}}(\Gr_G)_{\fL(N^-_P)} \text{ and } \Dmod_{\frac{1}{2}}(\Gr_G)^{\fL(N^-_P)}$$
can both be identified with 
$$\Dmod_{\frac{1}{2}}(\Gr_G)^{I^{-,0}_P},$$
where $I^{-,0}_P$ is the unipotent radical of $I^-_P$. 

\medskip

We can rewrite $(\Dmod_{\frac{1}{2}}(\Gr_G)^{I^{-,0}_P})^{M,\on{ren}}$ as
$$\Dmod_{\frac{1}{2}}(I^{-,0}_P\backslash \fL(G)/K)^{M\times G,\on{ren}},$$
where $K\subset \fL^+(G)$ is the first congruence subgroup, and the renormalization construction
is taken with respect to the action of the group $M\times G$. 

\medskip

We can further rewrite this as
$$\Dmod_{\frac{1}{2}}(I^-_P\backslash \fL(G)/K)^{G,\on{ren}},$$
and where the renormalization is taken with respect to $G$. 

\medskip

Let now be $\bC_1,\bC_2$ be an arbitrary pair of categories, acted on by $\fL(G)$. We claim that the functor
\begin{equation} \label{e:non-ren abs}
\bC_1^{\fL^+(G)}\underset{\Sph_G}\otimes ((\bC_2)^K)^{G,\on{ren}}\to \bC_1^{\fL^+(G)}\underset{\Sph_G}\otimes ((\bC_2)^K)^{G}
\end{equation} 
is an equivalence. 

\medskip

Note that 
$$\on{C}^\cdot(\on{pt}/G)\mod \simeq \Dmod(\on{pt}/G)\simeq \Dmod(\on{pt}/\fL^+(G)))$$ maps to the center of $\Sph_G$. It suffices to show that the 
action of $\on{C}^\cdot(\on{pt}/G)$ on the left-hand side of \eqref{e:non-ren abs} (via the $((\bC_2)^K)^{G,\on{ren}}$ term) factors
via an action of $\on{C}^\cdot(\on{pt}/G)\mod_0$. 

\medskip

However, the above action equals the action via the $\bC_1^{\fL^+(G)}$ term, and the required factorization is manifest.

\end{proof} 

\sssec{}

Will denote the (fully fathful) functor 
$$\KL(G)_{\crit_G}\underset{\Sph_G}\otimes \on{I}(G,P^-)^{\on{loc}}\to \hg\mod_{\crit_G}^{-,\semiinf}$$
of \eqref{e:Sph gen KM} by $\Accs$. 

\medskip

Will denote its essential image by
$$\hg\mod_{\crit_G}^{-,\semiinfAccs}\subset \hg\mod_{\crit_G}^{-,\semiinf},$$
and we will refer it it as \emph{the spherically accessible} subcategory of $\hg\mod_{\crit_G}^{-,\semiinf}$.

\sssec{}

From \lemref{l:KM to Sph gen} we obtain:

\begin{cor} \label{c:KM to Sph gen}
The functor \eqref{e:KM to Sph gen} identifies with the right adjoint of $\Accs$. 
\end{cor} 

\sssec{}

From \corref{c:KM to Sph gen} we obatin: 

\begin{cor} \label{c:Wak Sph via semiinf}
The functor $\Wak^{-,\on{Sph}}$ is recovered from $\Wak^{-,\on{enh}}$ as
\begin{multline*}
\KL(M)_{\crit_M-\rhoch_P} \overset{\Wak^{-,\on{enh}}}\longrightarrow \KL(G)_{\crit_G}\underset{\Sph_G}\otimes \on{I}(G,P^-)^{\on{loc}}
\overset{\on{Id}\otimes ((-)\underset{G}\star \Delta^{-,\semiinf})^R}\longrightarrow \\
\to \KL(G)_{\crit_G}\underset{\Sph_G}\otimes \Sph_G=\KL(G)_{\crit_G}.
\end{multline*}
\end{cor}

Note that the functor $(\Delta^{-,\semiinf})^R$ appearing in the above formula is just 
\begin{multline*}\on{I}(G,P^-)^{\on{loc}}\simeq \Dmod_{\frac{1}{2}}(\Gr_G)^{\fL(N^-_P)\cdot \fL^+(M),\on{ren}}\to \\
\to \Dmod_{\frac{1}{2}}(\Gr_G)^{\fL^+(M),\on{ren}}\overset{\on{Av}_*^{\fL^+(G)/\fL^+(M)}}\longrightarrow 
\Dmod_{\frac{1}{2}}(\Gr_G)^{\fL^+(G),\on{ren}}=:\Sph_G.
\end{multline*}

\ssec{Twisting by \texorpdfstring{$Z_G$}{ZG}-torsors}

In this subsection we introduce the twists by central torsors (see \secref{ss:central twist}) in the context of Kac-Moody representations. 

\medskip

As in \secref{ss:twisting by dual torsors}, in the bulk of the paper, the observations from this subsection will be applied when the reductive group
in question is not the original $G$, but rather one of its Levi subgroups. 

\sssec{}

Let $\CP_G$ be a $G$-bundle on $X$. As in \cite[Sect. 1.2]{GLC2}, we can consider the $\CP_G$-twists of 
$\hg\mod_\kappa$ and $\KL(G)_\kappa$, denoted  
$$\hg\mod_{\kappa,\CP_G} \text{ and } \KL(G)_{\kappa,\CP_G},$$
respectively.

\sssec{} \label{sss:no twist on KL}

Note, however, that as in \cite[Sect. 1.2.3]{GLC2}, we have a canonical equivalence
\begin{equation} \label{e:alpha taut KL}
\alpha_{\CP_G,\on{taut}}:\KL(G)_\kappa\to \KL(G)_{\kappa,\CP_G}. 
\end{equation}

This equivalence fits into the commutative diagram
$$
\CD
\KL(G)_{\kappa} @>{\alpha_{\CP_G,\on{taut}}}>> \KL(G)_{\kappa,\CP_G} \\
@V{\Res^{(\hg,\fL^+(G))_\kappa}_{\fL^+(G)}}VV @VV{\oblv^{(\hg,\fL^+(G))_\kappa}_{\fL^+(G)}}V \\
\Rep(\fL^+(G)) @>{\alpha_{\CP_G,\on{taut}}}>> \Rep(\fL^+(G)_{\CP_G}).
\endCD
$$

\sssec{} \label{sss:kappa dlog}

Assume now that $\CP_G$ is induced from a $Z_G$-torsor $\CP_{Z_G}$. In this case, we have a canonical
identification
\begin{equation} \label{e:central twist of KM}
\alpha_{\CP_{Z_G},\on{cent}}:\hg\mod_{\kappa,\CP_{Z_G}}\overset{\sim}\to \hg\mod_{\kappa+\kappa(\on{dlog}(\CP_{Z_G}),-)},
\end{equation}
where:

\begin{itemize}

\item $\on{dlog}(\CP_{Z_G})$ is the $\fL^+(z_\fg\otimes \omega_X)$-torsor, 
induced by means of the map
$$\on{dlog}:\fL^+(Z_G)\to \fL^+(z_\fg\otimes \omega_X).$$

\item The level $\kappa(-,-)$ is viewed as defining a map 
$$\kappa(-,-):z_\fg\to z_\cg,$$
and so we can view 
$$\kappa(\on{dlog}(\CP_{Z_G}),-)$$
as a torsor with respect to 
$$\fL^+(z_\cg\otimes \omega_X)\simeq (\fL(\fg_{\on{ab}})/\fL^+(\fg_{\on{ab}}))^*.$$

\end{itemize}

\sssec{}

The equivalence \eqref{e:central twist of KM} induces an equivalence
\begin{equation} \label{e:central twist of KL}
\alpha_{\CP_{Z_G},\on{cent}}:\KL(G)_{\kappa,\CP_{Z_G}}\overset{\sim}\to\KL(G)_{\kappa+\kappa(\on{dlog}(\CP_{Z_G}),-)}.
\end{equation}

Composing with \eqref{e:alpha taut KL}, we obtain an equivalence, denoted 
\begin{equation} \label{e:central un-twist of KL}
\KL(G)_{\kappa}\overset{(\on{transl}_{\CP_{Z_G}})^*}\longrightarrow \KL(G)_{\kappa+\kappa(\on{dlog}(\CP_{Z_G}),-)}.
\end{equation}

This equivalence fits into the commutative diagram
$$
\CD
\KL(G)_{\kappa} @>{(\on{transl}_{\CP_{Z_G}})^*}>> \KL(G)_{\kappa+\kappa(\on{dlog}(\CP_{Z_G}),-)} \\
@V{\Res^{(\hg,\fL^+(G))}_{\fL^+(G)}}VV @VV{\Res^{(\hg,\fL^+(G))}_{\fL^+(G)}}V \\
\Rep(\fL^+(G)) & & \Rep(\fL^+(G)) \\
@V{=}VV @VV{=}V \\
\QCoh(\on{pt}/\fL^+(G)) @>{(\on{transl}_{\CP_{Z_G}})^*}>>  \QCoh(\on{pt}/\fL^+(G)) 
\endCD
$$
where the bottom horizontal arrow is the functor of pullback with respect to the automorphism $\on{transl}_{\CP_{Z_G}}$
of $\on{pt}/\fL^+(G)$, given by tensoring with $\CP_{Z_G}|_\cD$. 

\sssec{} \label{sss:trans for crit}

Note that when $\kappa|_{z_\fg}=0$, (e.g., when we are at the critical level), we have
$$\kappa(\on{dlog}(\CP_{Z_G}),-)=0.$$

So in this case, $(\on{transl}_{\CP_{Z_G}})^*$ is an endofunctor of $\KL(G)_\kappa$. It can be explicitly described as follows:

\medskip

For every trivialization of $\CP_{Z_G}$, we have
$$(\on{transl}_{\CP_{Z_G}})^*\simeq \on{Id}.$$

A change of trivialization given by an element $g\in \fL^+(Z_G)$ corresponds to the automorphism of the 
identity endofunctor of $\KL(G)_\kappa$ given by the action of $g$ on modules. 

\sssec{Notational convention} \label{sss:lambda twist}

Assume for a moment that the above $Z_G$-torsor is of the form $\lambda(\omega_X)$, 
where $\lambda:\BG_m\to Z_G$.

\medskip

In this case, we will use a shorthandnotation
$$\hg\mod_{\kappa+\kappa(\lambda,-)}:=\hg\mod_{\kappa+\kappa(\on{dlog}(\lambda(\omega_X)),-)} \text{ and }
\KL(G)_{\kappa+\kappa(\lambda,-)}:=\KL(G)_{\kappa+\kappa(\on{dlog}(\lambda(\omega_X)),-)}.$$ 

Note that this agrees with the notation introduced in \secref{sss:lambdach twist}, where we regard
$\kappa(\lambda,-)$ as an element of $z_\cg$. 

\medskip

As in {\it loc.cit.}, the notations 
$$\hg\mod_{\kappa+\kappa(\lambda,-)} \text{ and } \KL(G)_{\kappa+\kappa(\lambda,-)}$$ make sense for an arbitrary 
element $\lambda\in z_\fg$. 

\sssec{}

Similar constructions apply if instead of $\KL(G)_\kappa$ (resp., $\hg\mod_\kappa$), we start with the category of the form
$\KL(G)_{\kappa+\on{dlog}(\CP_{Z^0_\cG})}$ (resp., $\hg\mod_{\kappa+\on{dlog}(\CP_{Z^0_\cG})}$) for a $Z^0_\cG$-torsor
$\CP_{Z^0_\cG}$.

\medskip

We denote the resulting categories by 
$$\KL(G)_{\kappa+\on{dlog}(\CP_{Z^0_\cG})+\kappa(\on{dlog}(\CP_{Z_G}),-)} \text{ and }  
\hg\mod_{\kappa+\on{dlog}(\CP_{Z^0_\cG})+\kappa(\on{dlog}(\CP_{Z_G}),-)}$$
or
$$\KL(G)_{\kappa+\lambdach+\kappa(\lambda,-)}  \text{ and }  \hg\mod_{\kappa+\lambdach+\kappa(\lambda,-)},$$
respectively.

\ssec{The twisted Jacquet and Wakimoto functors}

In this subsection we will adapt the (generalized) Jacquet and Wakimoto functors to the twisted situation. 

\sssec{}

Let $\CP_M$ be an $M$-torsor. Consider the $\CP_M$-twists of the objects
appearing in \secref{ss:BRST}, see \secref{ss:twist by G-bundle}.

\medskip

In particular, we obtain the functor
$$J^-_{\on{KM},\CP_M}:\hg\mod_{\crit_G+\kappa,\CP_M} \to \hm\mod_{\crit_M+\kappa+\rhoch_P,\CP_M}.$$ 

\sssec{}

Note, however, that by \secref{sss:no twist on KL}, we have the equivalences
$$\alpha_{\CP_M,\on{taut}}:\hg\mod_{\crit_G+\kappa}^{\fL^+(M)}\to \hg\mod_{\crit_G+\kappa,\CP_M}^{\fL^+(M)}$$  and 
$$\alpha_{\CP_M,\on{taut}}:\KL(M)_{\crit_M+\kappa+\rhoch_P}\to \KL(M)_{\crit_M+\kappa+\rhoch_P,\CP_M}$$
so that the diagram
\begin{equation} \label{e:BRST twisted CD}
\CD
\hg\mod_{\crit_G+\kappa}^{\fL^+(M)} @>{\alpha_{\CP_M,\on{taut}}}>> \hg\mod_{\crit_G+\kappa,\CP_M}^{\fL^+(M)} \\
@V{J^-_{\on{KM}}}VV @VV{J^-_{\on{KM},\CP_M}}V \\
\KL(M)_{\crit_M+\kappa+\rhoch_P} @>{\alpha_{\CP_M,\on{taut}}}>> \KL(M)_{\crit_M+\kappa+\rhoch_P,\CP_M}
\endCD
\end{equation} 
commutes. 

\sssec{} \label{sss:BRST shift conv}

Assume now that $\CP_M$ is induced by a $Z_M$-bundle $\CP_{Z_M}$. In this case, we can identify
$$\hm\mod_{\crit_M+\kappa+\rhoch_P,\CP_{Z_M}}\overset{\alpha_{\CP_{Z_M},\on{cent}}}\simeq \hm\mod_{\crit_M+\kappa+\rhoch_P+\kappa(\on{dlog}(\CP_{Z_M}),-)}$$ and  
$$\KL(M)_{\crit_M+\kappa+\rhoch_P,\CP_{Z_M}}\overset{\alpha_{\CP_{Z_M},\on{cent}}}\simeq \KL(M)_{\crit_M+\kappa+\rhoch_P+\kappa(\on{dlog}(\CP_{Z_M}),-)}$$
see \eqref{e:central twist of KM} and \eqref{e:central twist of KL}, respectively. 
 
\medskip

Using the above identofications, we will view $J^-_{\on{KM},\CP_{Z_M}}$ as a functor
$$\hg\mod_{\crit_G+\kappa}^{\fL^+(M)} \to \KL(M)_{\crit_M+\kappa+\rhoch_P+\kappa(\on{dlog}(\CP_{Z_M}),-)}$$
equal to
\begin{multline} \label{e:BRST twisted}
\hg\mod_{\crit_G+\kappa}^{\fL^+(M)} \overset{\alpha_{\CP_{Z_M},\on{taut}}}\longrightarrow 
\hg\mod_{\crit_G+\kappa,\CP_{Z_M}}^{\fL^+(M)} \overset{J^-_{\on{KM},\CP_{Z_M}}}\longrightarrow \\
\to \KL(M)_{\crit_M+\kappa+\rhoch_P,\CP_{Z_M}}\overset{\alpha_{\CP_{Z_M},\on{cent}}}\simeq \KL(M)_{\crit_M+\kappa+\rhoch_P+\kappa(\on{dlog}(\CP_{Z_M}),-)},
\end{multline} 
or, which is the same,
\begin{equation} \label{e:BRST twisted bis}
\hg\mod_{\crit_G+\kappa}^{\fL^+(M)}  \overset{J^-_{\on{KM}}}\longrightarrow  \KL(M)_{\crit_M+\kappa+\rhoch_P}
\overset{(\on{transl}_{\CP_{Z_M}})^*}\simeq \KL(M)_{\crit_M+\kappa+\rhoch_P+\kappa(\on{dlog}(\CP_{Z_M}),-)}.
\end{equation}

\sssec{Warning}

Note that when $\kappa=0$ (i.e., we are at the critical level), both functors 
$$J^-_{\on{KM}} \text{ and } J^-_{\on{KM},\CP_{Z_M}}$$
can be viewed as mapping $\hg\mod_{\crit_G}^{\fL^+(M)}$ to the \emph{same} category 
$\KL(M)_{\crit_M+\rhoch_P}$.

\medskip

Yet, these two functors are different: namely, they differ by the automorphism of $\KL(M)_{\crit_M+\rhoch_P}$ given by
$(\on{transl}_{\CP_{Z_M}})^*$, see \secref{sss:trans for crit}. 

\medskip

In particular, they have different compatibilities with respect to the actions of $\Sph_M$. Indeed, in the twisted
version, while $\Sph_M$ acts naturally on $\hg\mod_{\crit_G}^{\fL^+(M)}$, its action on $\KL(M)_{\crit_M+\rhoch_P}$
is twisted by the automorphism $(\on{transl}_{\CP_{Z_M}})^*$ of $\Sph_M$, see Remark \ref{r:non-trivial transl}

\sssec{}

We will denote by $J^{-,\Sph}_{\on{KM},\CP_{Z_M}}$ the restriction of \eqref{e:BRST twisted} along
$$\KL(G)_{\crit_G+\kappa}:=\hg\mod_{\crit_G+\kappa}^{\fL^+(G)}\to \hg\mod_{\crit_G+\kappa}^{\fL^+(M)},$$
so this is a functor
\begin{equation} \label{e:BRST Sph} 
J^{-,\Sph}_{\on{KM},\CP_{Z_M}}:\KL(G)_{\crit_G+\kappa}\to \KL(M)_{\crit_M+\kappa+\rhoch_P-\kappa(\on{dlog}(\CP_{Z_M}),-)}.
\end{equation} 

From \eqref{e:BRST twisted CD} we obtain a commutative diagram
\begin{equation} \label{e:BRST twisted CD sph}
\CD
\KL(G)_{\crit_G+\kappa} @>{\on{Id}}>> \KL(G)_{\crit_G+\kappa}  \\
@V{J^{-,\Sph}_{\on{KM}}}VV @VV{J^{-,\Sph}_{\on{KM},\CP_{Z_M}}}V \\
\KL(M)_{\crit_M+\kappa+\rhoch_P} @>{(\on{transl}_{\CP_{Z_M}})^*}>> \KL(M)_{\crit_M+\kappa+\rhoch_P+\kappa(\on{dlog}(\CP_{Z_M}),-)}.
\endCD
\end{equation}

\sssec{}

Similar conventions apply to the Wakimoto functors. In particular, we obtain the functors
$$\Wak^{-,\semiinf}_{\CP_{Z_M}}:\KL(M)_{\crit_M-\kappa-\rhoch_P-\kappa(\on{dlog}(\CP_{Z_M}),-)}\to 
\hg\mod_{\crit_G-\kappa}^{-,\semiinf}$$
and
$$\Wak^{-,\Sph}_{\CP_{Z_M}}:\KL(M)_{\crit_M-\kappa-\rhoch_P-\kappa(\on{dlog}(\CP_{Z_M}),-)}\to 
\KL(G)_{\crit_G-\kappa}.$$

\sssec{} \label{sss:BRST summary}

We will apply the above constructions mostly in the case when $\CP_{Z_M}=\rho_P(\omega_X)$.
So we obtain the functors
$$\ol{J}^-_{\on{KM},\rho_P(\omega_X)}:
(\hg\mod_{\crit_G+\kappa})^{\fL^+(M)}_{\fL(N^-_P)}\to \KL(M)_{\crit_M+\kappa+\rhoch_P+\kappa(\rho_P,-)},$$
$$\Wak^{-,\semiinf}_{\rho_P(\omega_X)}:\KL(M)_{\crit_M-\kappa-\rhoch_P-\kappa(\rho_P,-)}\to \hg\mod_{\crit_G-\kappa}^{-,\semiinf},$$
$$\Wak_{\rho_P(\omega_X)}:\KL(M)_{\crit_M-\kappa-\rhoch_P-\kappa(\rho_P,-)}\to (\hg\mod_{\crit_G-\kappa})^{\fL^+(P)},$$
as well as
$$J^{-,\Sph}_{\rho_P(\omega_X)}:\KL(G)_{\crit_G+\kappa}\to \KL(M)_{\crit_M+\kappa+\rhoch_P+\kappa(\rho_P,-)}$$
and
$$\Wak^{-,\Sph}_{\rho_P(\omega_X)}:\KL(M)_{\crit_M-\kappa-\rhoch_P-\kappa(\rho_P,-)}\to \KL(G)_{\crit_G-\kappa}.$$

\sssec{}

We now specialize to the critical level. In this case the functors in \secref{sss:BRST summary} map
$$\ol{J}^-_{\on{KM},\rho_P(\omega_X)}:
(\hg\mod_{\crit_G})^{\fL^+(M)}_{\fL(N^-_P)}\to \KL(M)_{\crit_M+\rhoch_P},$$
$$\Wak^{-,\semiinf}_{\rho_P(\omega_X)}:\KL(M)_{\crit_M-\rhoch_P}\to \hg\mod_{\crit_G}^{-,\semiinf},$$
$$\Wak_{\rho_P(\omega_X)}:\KL(M)_{\crit_M-\rhoch_P}\to (\hg\mod_{\crit_G})^{\fL^+(P)},$$
as well as
$$J^{-,\Sph}_{\rho_P(\omega_X)}:\KL(G)_{\crit_G}\to \KL(M)_{\crit_M+\rhoch_P}$$
and
$$\Wak^{-,\Sph}_{\rho_P(\omega_X)}:\KL(M)_{\crit_M-\rhoch_P}\to \KL(G)_{\crit_G}.$$

\sssec{} \label{sss:KM Jacquet summary 1}

In what follows we will also need twisted versions of the functors $J^{-,\on{enh}_{\on{co}}}_{\on{KM}}$, $J^{-,\on{enh}}_{\on{KM}}$ and $\Wak^{-,\on{enh}}$,
respectively.

\medskip

We let 
$$J^{-,\on{enh}_{\on{co}}}_{\on{KM},\rho_P(\omega_X)}:
\KL(G)_{\crit_G}\underset{\Sph_G}\otimes \on{I}(G,P^-)^{\on{loc}}_{\on{co},\rho_P(\omega_X)}\to \KL(M)_{\crit_M+\rhoch_P}$$
be the functor
\begin{multline*}
\KL(G)_{\crit_G}\underset{\Sph_G}\otimes \on{I}(G,P^-)^{\on{loc}}_{\on{co},\rho_P(\omega_X)} \overset{\on{Id}\otimes \alpha^{-1}_{\on{taut},\rho_P(\omega_X)}}\simeq
\KL(G)_{\crit_G}\underset{\Sph_G}\otimes \on{I}(G,P^-)^{\on{loc}}_{\on{co}} \overset{J^{-,\on{enh}_{\on{co}}}_{\on{KM}}}\longrightarrow \\
\to \KL(M)_{\crit_M+\rhoch_P}\overset{(\on{transl}_{\rho_P(\omega_X)})^*}\longrightarrow \KL(M)_{\crit_M+\rhoch_P},
\end{multline*}
which is the same as
\begin{multline*}
\KL(G)_{\crit_G}\underset{\Sph_G}\otimes \on{I}(G,P^-)^{\on{loc}}_{\on{co},\rho_P(\omega_X)} 
\overset{\alpha_{\on{taut},\rho_P(\omega_X)}\otimes \on{Id}}\simeq 
\KL(G)_{\crit_G,\rho_P(\omega_X)}\underset{\Sph_G}\otimes \on{I}(G,P^-)^{\on{loc}}_{\on{co},\rho_P(\omega_X)} \to \\
\to (\hg\mod_{\crit_G,\rho_P(\omega_X)})^{\fL^+(M)_{\rho_P(\omega_X)}}_{\fL(N^-_P)\rho_P(\omega_X)} 
\overset{\ol{J}^-_{\on{KM},\rho_P(\omega_X)}}\longrightarrow  \KL(M)_{\crit_M+\rhoch_P,\rho_P(\omega_X)}\to \\
\overset{\alpha_{\on{cent},\rho_P(\omega_X)}}\longrightarrow \KL(M)_{\crit_M+\rhoch_P}.
\end{multline*}

This functor intertwines the natural actions of $\Sph_M$ on the two sides (see \secref{sss:convent Sph-action on IGP} 
for the conventions regarding the action of $\Sph_M$
on $\on{I}(G,P^-)^{\on{loc}}_{\rho_P(\omega_X)}$). 

\medskip 

The composition
$$
\KL(G)_{\crit_G} \overset{\on{Id}\otimes \one_{\on{I}(G,P^-)^{\on{loc}}_{\on{co},\rho_P(\omega_X)} }}\longrightarrow 
\KL(G)_{\crit_G}\underset{\Sph_G}\otimes \on{I}(G,P^-)^{\on{loc}}_{\on{co},\rho_P(\omega_X)}
\overset{J^{-,\on{enh}_{\on{co}}}_{\on{KM},\rho_P(\omega_X)}}\longrightarrow \KL(M)_{\crit_M+\rhoch_P}
$$
recovers the functor $J^{-,\Sph}_{\rho_P(\omega_X)}$. 

\sssec{} \label{sss:KM Jacquet summary 2}

We let 
$$J^{-,\on{enh}}_{\on{KM},\rho_P(\omega_X)}:\KL(G)_{\crit_G}\to 
\KL(M)_{\crit_M+\rhoch_P} \underset{\Sph_M}\otimes \on{I}(G,P^-)^{\on{loc}}_{\rho_P(\omega_X)}$$
be the functor obtained from $J^{-,\on{enh}_{\on{co}}}_{\on{KM},\rho_P(\omega_X)}$ by duality. 

\medskip

This functor intertwines the natural actions of $\Sph_G$ on the two sides. 

\medskip

The composition
$$\KL(G)_{\crit_G}\overset{J^{-,\on{enh}}_{\on{KM},\rho_P(\omega_X)}}\longrightarrow 
\KL(M)_{\crit_M+\rhoch_P} \underset{\Sph_M}\otimes \on{I}(G,P^-)^{\on{loc}}_{\rho_P(\omega_X)} 
\overset{\oblv_{\on{enh}}}\longrightarrow \KL(M)_{\crit_M+\rhoch_P},$$
recovers the functor $J^{-,\Sph}_{\rho_P(\omega_X)}$, where
$\oblv_{\on{enh}}$ denotes the composition 
$$\KL(M)_{\crit_M+\rhoch_P} \underset{\Sph_M}\otimes \on{I}(G,P^-)^{\on{loc}}_{\rho_P(\omega_X)} 
\overset{\on{Id}\otimes \oblv_{\semiinf\to \Sph}}\longrightarrow 
\KL(M)_{\crit_M+\rhoch_P} \underset{\Sph_M}\otimes \Sph_M= \KL(M)_{\crit_M+\rhoch_P},$$
 
\sssec{} \label{sss:Jacquet KM summary}

In the sequel we will also use the notations
$$\KL(G)^{-,\on{enh}_{\on{co}}}_{\crit_G,\rho_P(\omega_X)}:=
\KL(G)_{\crit_G}\underset{\Sph_G}\otimes \on{I}(G,P^-)^{\on{loc}}_{\on{co},\rho_P(\omega_X)}$$
and
$$\KL(M)_{\crit_M+\rhoch_P,\rho_P(\omega_X)}^{-,\on{enh}}:=
\KL(M)_{\crit_M+\rhoch_P} \underset{\Sph_G}\otimes \on{I}(G,P^-)^{\on{loc}}_{\rho_P(\omega_X)},$$
so that we can view $J^{-,\on{enh}_{\on{co}}}_{\on{KM},\rho_P(\omega_X)}$ as a functor
$$\KL(G)^{-,\on{enh}_{\on{co}}}_{\crit_G,\rho_P(\omega_X)}\to \KL(M)_{\crit_M+\rhoch_P}$$
and $J^{-,\on{enh}}_{\on{KM},\rho_P(\omega_X)}$ as a functor
$$\KL(G)_{\crit_G}\to \KL(M)_{\crit_M+\rhoch_P,\rho_P(\omega_X)}^{-,\on{enh}}.$$

\sssec{}

We let 
$$\Wak^{-,\on{enh}}_{\rho_P(\omega_X)}:
\KL(M)_{\crit_M-\rhoch_P}\to \KL(G)_{\crit_G}\underset{\Sph_G}\otimes \on{I}(G,P^-)^{\on{loc}}_{\rho_P(\omega_X)}$$
be the functor dual to $J^{-,\on{enh}_{\on{co}}}_{\on{KM},\rho_P(\omega_X)}$.

\medskip

I.e., $\Wak^{-,\on{enh}}_{\rho_P(\omega_X)}$ is the composition
\begin{multline*}
\KL(M)_{\crit_M-\rhoch_P} \overset{(\on{transl}_{-\rho_P(\omega_X)})^*}\longrightarrow \KL(M)_{\crit_M-\rhoch_P} \overset{\Wak^{-,\on{enh}}}\longrightarrow \\
\to \KL(G)_{\crit_G}\underset{\Sph_G}\otimes \on{I}(G,P^-)^{\on{loc}} \overset{\on{Id}\otimes \alpha_{\on{taut},\rho_P(\omega_X)}}\longrightarrow 
\KL(G)_{\crit_G}\underset{\Sph_G}\otimes \on{I}(G,P^-)^{\on{loc}}_{\rho_P(\omega_X)}:=\KL(G)_{\crit_G,\rho_P(\omega_X)}^{-,\on{enh}}.
\end{multline*}

\ssec{Drinfeld-Sokolov reduction of Wakimoto modules} \label{ss:DS of Wak}

In this subsection we compute the composition of the Drinfeld-Sokolov reduction and Wakimoto functors.
This computation lies at the core of the proof of \thmref{t:local Jacquet dual enh}. 

\medskip

At its essence, this computation goes back to Feigin and Frenkel, and is lies also in the core of their work on 
duality for W-algebras.

\sssec{}

We have the following basic assertion, which follows by duality from  \corref{c:wak restr to n}:

\begin{lem} \label{l:DS of Wak}
The following diagram of (lax unital) factorization functors commutes:
\begin{equation} \label{e:DS of Wak}
\CD
\hm\mod_{\crit_M+\kappa-\rhoch_P,\rho(\omega_X)} @>{\alpha_{\rho_P(\omega),\on{cent}}}>> \hm\mod_{\crit_M+\kappa-\rhoch_P+\kappa(\rho_P,-),\rho_M(\omega_X)}  \\
@V{\Wak^{-,\semiinf}_{\rho(\omega_X)}}VV @VV{\DS_M}V \\
\hg\mod_{\crit_G+\kappa,\rho(\omega_X)}  @>{\DS_G}>>  \Vect.
\endCD
\end{equation}
\end{lem}

\sssec{}

From \lemref{l:DS of Wak} we obtain:

\begin{cor} \label{c:DS of Wak}
The following diagram of (lax unital) factorization functors commutes:
$$
\CD
\KL(M)_{\crit_M+\kappa-\rhoch_P+\kappa(\rho_P,-)} & @>{\alpha_{\rho_M(\omega),\on{taut}}}>>&  \KL(M)_{\crit_M+\kappa-\rhoch_P+\kappa(\rho_P,-),\rho_M(\omega_X)}  \\
@V{\Wak^{-,\semiinf}_{\rho_P(\omega_X)}}VV & & @VV{\DS_M}V \\
\hg\mod^{\fL^+(M)}_{\crit_G+\kappa}  @>{\alpha_{\rho(\omega),\on{taut}}}>> \hg\mod^{\fL^+(M)}_{\crit_G+\kappa,\rho(\omega_X)}  @>{\DS_G}>>  \Vect.
\endCD
$$
\end{cor}

\sssec{} \label{sss:DS of usual Wak}

For future reference we observe that the functors
$$\hg\mod_{\crit_G,\rho(\omega_X)} \overset{\DS_G}\to \Vect$$
and
$$\hg\mod_{\crit_G,\rho(\omega_X)}  \overset{\Av^{\fL^+(N_P)_{\rho(\omega_X)}}_*}\longrightarrow 
\hg\mod^{\fL^+(N_P)_{\rho(\omega_X)}}_{\crit_G,\rho(\omega_X)} \to \hg\mod_{\crit_G,\rho(\omega_X)} \overset{\DS_G}\to \Vect$$
are canonically isomorphic.

\medskip

In particular, from \corref{c:DS of Wak} we obtain:
\begin{cor} \label{c:DS of Wak usual}
The following diagram of (lax unital) factorization functors commutes:
$$
\CD
\KL(M)_{\crit_M+\kappa-\rhoch_P+\kappa(\rho_P,-)} & @>{\alpha_{\rho_M(\omega),\on{taut}}}>>&  \KL(M)_{\crit_M+\kappa-\rhoch_P+\kappa(\rho_P,-),\rho_M(\omega_X)}  \\
@V{\Wak_{\rho_P(\omega_X)}}VV & & @VV{\DS_M}V \\
\hg\mod^{\fL^+(P)}_{\crit_G+\kappa}  @>{\alpha_{\rho(\omega),\on{taut}}}>> \hg\mod^{\fL^+(P)}_{\crit_G+\kappa,\rho(\omega_X)}  @>{\DS_G}>>  \Vect.
\endCD
$$
\end{cor}

\sssec{}

Specializing to the critical level, from \corref{c:DS of Wak} we obtain:

\begin{cor} \label{c:DS of Wak crit}
The following diagrams of (lax unital) factorization functors commute:
$$
\CD
\KL(M)_{\crit_M-\rhoch_P} & @>{\alpha_{\rho_M(\omega),\on{taut}}}>>&  \KL(M)_{\crit_M-\rhoch_P,\rho_M(\omega_X)}  \\
@V{\Wak^{-,\semiinf}_{\rho_P(\omega_X)}}VV & & @VV{\DS_M}V \\
\hg\mod^{\fL^+(M)}_{\crit_G}  @>{\alpha_{\rho(\omega),\on{taut}}}>> \hg\mod^{\fL^+(M)}_{\crit_G,\rho(\omega_X)}  @>{\DS_G}>>  \Vect
\endCD
$$
and
$$
\CD
\KL(M)_{\crit_M-\rhoch_P} & @>{\alpha_{\rho_M(\omega),\on{taut}}}>>&  \KL(M)_{\crit_M-\rhoch_P,\rho_M(\omega_X)}  \\
@V{\Wak_{\rho_P(\omega_X)}}VV & & @VV{\DS_M}V \\
\hg\mod^{\fL^+(M)}_{\crit_G}  @>{\alpha_{\rho(\omega),\on{taut}}}>> \hg\mod^{\fL^+(M)}_{\crit_G,\rho(\omega_X)}  @>{\DS_G}>>  \Vect.
\endCD
$$
\end{cor}

\ssec{Proof of \propref{p:wak-exact}}\label{ss:wak-exact}

\sssec{}\label{sss:ignore-cnst}

As was explained in Remark \ref{r:wak-exact}, it is enough to prove 
that $\Wak$ is left t-exact.

\medskip 

Suppose $$\CM \in \hm\mod_{\crit_M+\kappa-\rhoch_P}^{\geq 0}$$ is coconnective.
We wish to show that $\Wak(\CM)$ is coconnective as well.
It suffices to do this in the case that 
$$\CM \in \hm\mod_{\crit_M+\kappa-\rhoch_P}^{\fL_n^+(M)}$$
for $\fL_n^+(M) \subset \fL^+(M)$ the $n$th congruence subgroup as
every coconnective $\CM$ is a filtered colimit of coconnective
$\fL_n^+(M)$-integrable modules.
Therefore, we assume $\CM$ is $\fL_n^+(M)$-integrable to start.

\medskip 

Now for $n$ fixed, it suffices to show that
$\Wak(\CM)$ is eventually coconnective, since \corref{c:wak-usual-explicit}
will then force it to be coconnective as is.

\sssec{}

Note that $\Wak(\CM) \in \hg\mod_{\crit_G+\kappa}^{\fL^+_n(G)\cdot\fL^+(N_P)}$
under our assumption. 

\medskip

Therefore, showing that $\Wak(\CM)$ is 
eventually coconnective is equivalent to showing that
\begin{equation}\label{eq:hom-vacs-to-wak}
\CHom(\BV_{\crit_G+\kappa,n},\Wak(\CM)) \in \Vect
\end{equation}

\noindent is eventually coconnective, 
where $\BV_{\crit_G+\kappa,n}$ is induced from 
the trivial module for $\fL_n^+(\fg) = t^n\fg[[t]]$.

\sssec{}

We can rewrite \eqref{eq:hom-vacs-to-wak} as
\begin{equation}\label{eq:vacs-wak-pairing}
\langle\BV_{\crit_G-\kappa,n}[n\cdot\dim\fg],\Wak(\CM)\rangle_{\hg}
\end{equation}

\noindent where $\langle-,\rangle_{\hg}$ is the pairing
\[
\hg\mod_{\crit_G-\kappa} \otimes \hg\mod_{\crit_G+\kappa} \to \Vect.
\]

\medskip 

By the definition of $\Wak$, we can rewrite \eqref{eq:vacs-wak-pairing} up to
a shift -- which we can ignore by \secref{sss:ignore-cnst} -- as
\begin{equation}\label{eq:bdd-jac}
\langle J^-_{\on{KM}}(\BV_{\crit_G-\kappa,n}),\CM\rangle_{\hm}.
\end{equation}

\sssec{}

We now make the following elementary observation:

\begin{lem}

Suppose $\fh$ is a Tate Lie algebra, $\mathfrak{k} \subset \fh$ is a subalgebra
that is a lattice, and $\fh_0 \subset \fh$ is a Tate Lie subalgebra.
Then for any $\CN \in \mathfrak{k}\mod$, $\on{ind}_{\mathfrak{k}}^{\fh}(\CN)$ has a canonical
$\bZ^{\geq 0}$-indexed filtration 
$$\on{fil}_{\bullet}(\on{ind}_{\mathfrak{k}}^{\fh}(\CN))$$
as an $\fh_0$-module so that
\[
\gr_i (\on{ind}_{\mathfrak{k}}^{\fh}(\CN))
\simeq \on{ind}_{\mathfrak{k}\cap \fh_0}^{\fh_0}
\Big(\Sym^i\big(\fh/(\mathfrak{k}+\fh_0)\big) \otimes \CN\Big) 
\in \fh_0\mod.
\]

\end{lem}

\sssec{}

We apply the lemma to $\fh = \hg$, $\fh_0 = \hp^-$, $\mathfrak{k} = \fL_n^+(\fg)$, 
and $\CN$ being the trivial module. It follows that $\BV_{\crit_G-\kappa,n}$
carries a filtration as a $\hp^-$-module so that
\[
\gr_{\bullet}(\BV_{\crit_G-\kappa,n}) = 
\on{ind}_{\mathfrak{\fL_n^+(\fp^-)}}^{\hp_{\crit_G-\kappa}}
\Big(\Sym \big(\fL(\fg)/(\fL_n^+(\fg)+\fL(\fp^-))\big)\Big).
\]

Therefore, it suffices to prove that 
\begin{equation}\label{eq:bdd-jac-2}
\big\langle \on{BRST}^-
\big(\on{ind}_{\mathfrak{\fL_n^+(\fp^-)}}^{\hp_{\crit_G-\kappa}}
\big(\Sym \big(\fL(\fg)/(\fL_n^+(\fg)+\fL(\fp^-))\big)\big)\big),
\CM\big\rangle_{\hm}
\end{equation}

\noindent is eventually coconnective.

\sssec{}

By general properties of semi-infinite cohomology, 
we can compute \eqref{eq:bdd-jac-2} as
\[
\on{inv}_{\fL^+_n(\fp^-)}
\Big(\Sym \big(\fL(\fg)/(\fL_n^+(\fg)+\fL(\fp^-))\big) \otimes \CM 
\otimes \det(\fL^+(\fp^-)/\fL_n^+(\fp^-))^{\vee}\Big)
[-\dim \fL^+(\fp^-)/\fL_n^+(\fp^-)].
\]

Because invariants is left t-exact, this complex is evidently bounded from 
below.

\section{Jacquet functor(s) for opers} \label{s:Op}

In the present section we introduce generalized Jacquet functors (as well as their duals) that connect the $\IndCoh^*(-)$ categories
on the spaces of monodromy-free opers for $\cG$ and $\cM$ to one another. Note, however, that the corresponding space of opers 
for $\cM$ is translated by the \emph{Miura shift}, see \secref{sss:Miura shift}. 

\medskip

The definition of these functors uses another geometric object that we call (parabolic) \emph{Miura opers}. 

\ssec{(Parabolic) Miura opers}

\sssec{Translated opers} \label{sss:transl oper} 

Let $\cG$ be a reductive group, and let $\CP_{Z^0_\cG}$ be a $Z^0_\cG$-torsor on $X$. We let
$$\Op_{\cG,\CP_{Z^0_\cG}}$$
be the following variant of the D-scheme $\Op_\cG$:

\medskip

In the definition of opers (see \cite[Sect. 3.1.1]{GLC2}), instead of requiring that the induced $\cT$-bundle be $\rhoch(\omega_X)$,
we require that it be
\begin{equation} \label{e:pre-Miura}
\rhoch(\omega_X)\otimes \CP_{Z^0_\cG}.
\end{equation} 

\medskip

We let 
$$\Op^\reg_{\cG,\CP_{Z^0_\cG}}\subset \Op^{\on{mon-free}}_{\cG,\CP_{Z^0_\cG}}\subset \Op^\mer_{\cG,\CP_{Z^0_\cG}}$$
denote the corresponding factorization spaces.\footnote{As in \cite[Sect. 3.1.3]{GLC2}, we distinguish notationally between the D-scheme $\Op_{\cG,\CP_{Z^0_\cG}}$
and the corresponding factorization space $\Op^\reg_{\cG,\CP_{Z^0_\cG}}$.}

\medskip

The material from \cite[Sect. 3.1]{GLC2} transfers verbatim to the present context.

\sssec{} \label{sss:Miura shift}

We now take the reductive group in question to be $\cM$, the Levi subgroup of a standard parabolic $\cP$.
We take
$$\CP_{Z^0_\cM}:=\rhoch_P(\omega_X).$$

Note that in this case, the $\cT$-bundle \eqref{e:pre-Miura} is
$$\rhoch_M(\omega_X)\otimes \rhoch_P(\omega_X)=\rhoch(\omega_X),$$
where $\rhoch$ in the right-hand side is the $\rhoch$ for $\cG$. 

\medskip

We will use a shorthandnotation
$$\Op_{\cM,\rhoch_P}:=\Op_{\cM,\rhoch_P(\omega_X)}.$$

\sssec{Example}

When $P=B$, we obtain that the scheme 
$$\Op_{\cT,\rhoch}$$
that classifies connections on the $\cT$-bundle $\rhoch(\omega_X)$.

\sssec{} \label{sss:Miura}

Let $\MOp_{\cG,\cP^-}$ denote the D-scheme of $\cP^-$-Miura opers, i.e., 
$$\MOp_{\cG,\cP^-}:=(\Op_\cG\underset{\LS_\cG}\times \LS_{\cP^-})^{\on{trans}},$$
where the superscript ``trans" refers to the condition that the $\cB$-reduction of the
$\cG$-bundles involved in the oper structure is transversal to the $\cP^-$-reduction.\footnote{Elsewhere in the literature, what we call ``Miura opers"
is sometimes referred to as ``generic Miura opers".}

\medskip

We have the natural forgetful map
\begin{equation} \label{e:MOp to Op}
\sfp^{\on{Miu}}:\MOp_{\cG,\cP^-}\to \Op_\cG
\end{equation} 

\sssec{}

Note also that we have a map 
\begin{equation} \label{e:MOp to Op M}
\sfq^{\on{Miu}}:\MOp_{\cG,\cP^-}\to \Op_{\cM,\rhoch_P},
\end{equation} 
constructed as follows:

\medskip

The $\cM$-bundle and the connection are induced from the $\cP^-$-bundle with its connection. The reduction
to $\cB(M)$ comes from the reduction of the original $\cG$-bundle to $\cB$.

\sssec{}

The following assertion is fundamental, albeit immediate:

\begin{lem} \label{l:Miura}
The map \eqref{e:MOp to Op M} is an isomorphism.
\end{lem} 

\begin{rem}
The composition 
$$\sfp^{\on{Miu}}\circ (\sfq^{\on{Miu}})^{-1}$$
is a map
$$\Op_{\cM,\rhoch_P}\to  \Op_\cG,$$
called the (parabolic) Miura transform.
\end{rem}

\begin{rem}

The assertion of \lemref{l:Miura} has the following classical counterpart: let
$$\fb =  \on{Fil}_0(\fg) \subset \on{Fil}_{-1}(\fg)\subset \fg$$
be the corresponding piece of the principal filtration. Then
$$(\on{Fil}_{-1}(\fg)/\on{Ad}(B)\underset{\fg/\on{Ad}(G)}\times \fp^-/\on{Ad}(P^-))
\underset{\on{pt}/B\underset{\on{pt}/G}\times \on{pt}/P^-}\times (\on{pt}/B\underset{\on{pt}/G}\times \on{pt}/P^-)^{\on{trans}},$$ 
projects isomorphically to
$$\on{Fil}_{-1}(\fm)/\on{Ad}(B_M),$$
where
$$(\on{pt}/B\underset{\on{pt}/G}\times \on{pt}/P^-)^{\on{trans}}\subset \on{pt}/B\underset{\on{pt}/G}\times \on{pt}/P^-$$
is the transversal locus, i.e., one corresponding to the open $B\times P^-$-orbit acting on $G$.

\end{rem} 

\ssec{Definition of the Jacquet functor for opers} \label{ss:spectral Jacquet}

In this subsection we introduce the unehnanced version of the Jacquet functor for opers. The definition will
be simple in that it does not require $\IndCoh^*(-)$ of algebro-geometric objects of semi-infinite nature. The 
latter will be appear, after some preparations, in the definition of the enhanced version, in \secref{ss:semiinf op}. 

\sssec{} \label{sss:Miura mf}

Consider the fiber product 
\begin{equation}  \label{e:P-mf MOp}
\MOp^\Pmf_{\cG,\cP^-}:=\MOp^\mer_{\cG,\cP^-}\underset{\LS^\mer_{\cP^-}}\times \LS^\reg_{\cP^-}.
\end{equation} 

Note that the maps
$$\Op^\mer_\cG\leftarrow \MOp^\Pmf_{\cG,\cP^-}
\to \Op^\mer_{\cM,\rhoch_P},$$
induced by $\sfp^{\on{Miu}}$ and $\sfq^{\on{Miu}}$, respectively, naturally factor via maps
\begin{equation} \label{e:Jacquet spec diag}
\Op^{\on{mon-free}}_\cG\overset{\sfp^{\on{Miu,mon-free}}}\longleftarrow 
\MOp^\Pmf_{\cG,\cP^-}
\overset{\sfq^{\on{Miu,mon-free}}}\longrightarrow \Op^{\on{mon-free}}_{\cM,\rhoch_P}.
\end{equation} 

\sssec{} \label{sss:Miura map mf}

Recall that according to \cite[Sect. 3.2.6]{GLC2}, 
$\MOp^\mer_{\cG,\cP^-}\simeq \Op^\mer_{\cM,\rhoch_P}$
is an ind-placid ind-affine factorization ind-scheme. 

\medskip

In addition, the map
$$\LS^\reg_{\cP^-}\to \LS^\mer_{\cP^-}$$
is ind-affine, locally almost of finite presentation (see \cite[Lemma B.7.13]{GLC2}). 
Hence, $\MOp^\Pmf_{\cG,\cP^-}$ is also an ind-placid ind-affine factorization ind-scheme. 

\medskip

In particular, we have well-defined mutually dual unital factorization categories
$$\IndCoh^!(\MOp^\Pmf_{\cG,\cP^-}) \text{ and } \IndCoh^*(\MOp^\Pmf_{\cG,\cP^-}),$$
see \cite[Sect. B.13.16-21]{GLC2}. And we have a pair of mutually dual (lax unital) factorization functors
$$(\sfp^{\on{Miu,mon-free}})^!:\IndCoh^!(\Op^{\on{mon-free}}_\cG)\to \IndCoh^!(\MOp^\Pmf_{\cG,\cP^-})$$ and 
$$(\sfp^{\on{Miu,mon-free}})^\IndCoh_*:\IndCoh^*(\MOp^\Pmf_{\cG,\cP^-}) \to \IndCoh^*(\Op^{\on{mon-free}}_\cG).$$

\medskip

Furthermore, the morphism $\sfq^{\on{Miu,mon-free}}$ is an ind-closed embedding,\footnote{Namely, the morphism 
in \cite[Lemma B.7.13]{GLC2} is an ind-closed embedding if $\sH$ is unipotent} locally almost of finite presentation. 
Hence, by \cite[Sect. A.10.11]{GLC2}, we have a pair of mutually dual (lax unital) factorization functors
$$(\sfq^{\on{Miu,mon-free}})^\IndCoh_*:\IndCoh^*(\MOp^\Pmf_{\cG,\cP^-}) \to \IndCoh^*(\Op^{\on{mon-free}}_{\cM,\rhoch_P})$$ and 
$$(\sfq^{\on{Miu,mon-free}})^!:\IndCoh^!(\Op^{\on{mon-free}}_{\cM,\rhoch_P})\to \IndCoh^!(\MOp^\Pmf_{\cG,\cP^-}).$$

\sssec{} \label{sss:J opers !}

Let 
$$J^{-,\on{spec}}_\Op:\IndCoh^!(\Op_\cG^{\on{mon-free}})\to \IndCoh^!(\Op_{\cM,\rhoch_P}^{\on{mon-free}})$$
denote the functor 
$$(\sfq^{\on{Miu,mon-free}})_*\circ (\sfp^{\on{Miu,mon-free}})^!.$$

\sssec{} \label{sss:spectral Jacquet Theta}

Denote by $J^{-,\on{spec}}_{\Op,\Theta}$ the functor
$$\IndCoh^*(\Op_\cG^{\on{mon-free}})\to \IndCoh^*(\Op_{\cM,\rhoch_P}^{\on{mon-free}}),$$
so that we have a commutative diagram
$$
\CD
\IndCoh^!(\Op_\cG^{\on{mon-free}}) @>{J^{-,\on{spec}}_{\Op}}>> \IndCoh^!(\Op_{\cM,\rhoch_P}^{\on{mon-free}}) \\
@V{\Theta_{\Op^\mf_\cG}}V{\sim}V @V{\sim}V{\Theta_{\Op^\mf_\cM}}V \\
\IndCoh^*(\Op_\cG^{\on{mon-free}}) @>{J^{-,\on{spec}}_{\Op,\Theta}}>> \IndCoh^*(\Op_{\cM,\rhoch_P}^{\on{mon-free}}),
\endCD
$$
where $\Theta_{\Op^\mf_\cG}$ and $\Theta_{\Op^\mf_\cM}$ are the identifications of \cite[Sect. 3.7.7]{GLC2}. 

\sssec{}

We will refer to $J^{-,\on{spec}}_{\Op}$ (and also $J^{-,\on{spec}}_{\Op,\Theta}$)
as the \emph{spectral Jacquet functor}. The main theorem in Part I of the paper will
establish its relationship with the functor $J^-_{\on{KM}}$ at the critical level. 

\sssec{}

We let
$$\on{co}\!J_\Op^{-,\on{spec}}:\IndCoh^*(\Op_{\cM,\rhoch_P}^{\on{mon-free}})\to \IndCoh^*(\Op_\cG^{\on{mon-free}})$$
the functor dual to $J^{-,\on{spec}}_{\Op}$. 

\medskip

Explicitly, it is given by
$$\on{co}\!J_\Op^{-,\on{spec}}=(\sfp^{\on{Miu,mon-free}})^\IndCoh_*\circ (\sfq^{\on{Miu,mon-free}})^!.$$


\ssec{The \emph{semi-infinite} Jacquet functor for opers, pointwise version}

In this subsection we start exploring the semi-infinite Jacquet functor for opers. 

\medskip

Ultimately, we will need a factorization version of this functor. However, by way of motivation, in this subsection we will work over a fixed point $\ul{x}\in \Ran$,
and we will be able to define a version of the Jacquet functor more general than in the factorization situation.

\sssec{}

Consider the (non-affine) D-scheme 
$$\Op_{\cG,\cP^-}:=\Op_{\cG}\underset{\on{pt}/\cG}\times \on{pt}/\cP^-$$
over $X$.

\medskip

Consider the corresponding factorization spaces:
$$\Op^\reg_{\cG,\cP^-}:=\Op^\reg_{\cG}\underset{\LS^\reg_{\cG}}\times \LS^\reg_{\cP^-};$$
$$\Op^\mer_{\cG,\cP^-}:=\Op^\mer_{\cG}\underset{\LS^\mer_{\cG}}\times \LS^\mer_{\cP^-};$$
$$\Op^\Pmf_{\cG,\cP^-}:=\Op^\mer_{\cG,\cP^-}\underset{\LS^\mer_{\cP^-}}\times \LS^\reg_{\cP^-};$$
$$\Op^\Mmf_{\cG,\cP^-}:=\Op^\mer_{\cG,\cP^-}\underset{\LS^\mer_{\cM}}\times \LS^\reg_{\cM};$$
$$\Op^\Gmf_{\cG,\cP^-}:=\LS^\reg_{\cG}\underset{\LS^\mer_{\cG}}\times \Op^\mer_{\cG,\cP^-}$$
and
$$\Op^\GMmf_{\cG,\cP^-}:=\LS^\reg_{\cG}\underset{\LS^\mer_{\cG}}\times \Op^\mer_{\cG,\cP^-}\underset{\LS^\mer_{\cM}}\times \LS^\reg_{\cM}.$$

Note that the latter space can also be written as
$$\Op^\mf_\cG\underset{\LS^\reg_\cG}\times \on{Hecke}^{\on{spec,loc}}_{\cG,\cP^-}.$$

\sssec{}

Consider the following maps:

\bigskip

\begin{itemize}

\item 
Let $\jmath$ denote the map
$$\MOp^\mer_{\cG,\cP^-}\to  \Op^\mer_{\cG,\cP^-},$$
and also the maps
$$\Op^\mf_{\cM,\rhoch_P}\overset{\sfq^{\on{Miu}}}\simeq \MOp^{\Mmf}_{\cG,\cP^-}:=\MOp^\mer_{\cG,\cP^-}\underset{\LS^\mer_{\cM}}\times \LS^\reg_{\cM}\to 
\Op^\Mmf_{\cG,\cP^-}$$
and 
$$\MOp^\Pmf_{\cG,\cP^-}:=\MOp^\mer_{\cG,\cP^-}\underset{\LS^\mer_{\cP^-}}\times \LS^\reg_{\cP^-}\to 
\Op^\Pmf_{\cG,\cP^-}$$
obtained by base change;

\medskip

\item
Let $\sfp^{-,\on{spec}}$ denote the map
$$\LS^\reg_{\cP^-} \to \LS^\reg_{\cG},$$
and also its base change my means of $\Op^{\on{mon-free}}_{\cG}\to \LS^\reg_{\cG}$, which is the map 
$$\Op^\Pmf_{\cG,\cP^-}\to \Op^{\on{mon-free}}_{\cG};$$

\medskip

\item $\iota$ denote the map
$$\LS^\reg_{\cP^-}\to \LS^\mer_{\cP^-}\underset{\LS^\mer_{\cM}}\times \LS^\reg_{\cM}=:\LS^\Mmf_{\cP^-},$$
and also the maps
\begin{equation} \label{e:iota Op}
\Op^\Pmf_{\cG,\cP^-}\to \Op^\Mmf_{\cG,\cP^-}.
\end{equation} 
and
$$\MOp^\Pmf_{\cG,\cP^-} \to \MOp^\Mmf_{\cG,\cP^-},$$
obtained by base change. 
\end{itemize}

\begin{rem} \label{r:j not open}
Note that although the map
$$\MOp^\reg_{\cG,\cP^-}\to \Op^\reg_{\cG,\cP^-}$$
is an open embedding, this is \emph{no longer the case} for the three maps denoted $\jmath$ above. 
\end{rem}

\sssec{}

Note that we have a Cartesian diagram 
\begin{equation} \label{e:MOpG OpG diag}
\CD
\MOp^\Pmf_{\cG,\cP^-} @>{\imath}>> \MOp^\Mmf_{\cG,\cP^-} @>{\sfq^{\on{Miu,mon-free}}}>{\sim}> \Op^\mf_{\cM,\rhoch_P}  \\
@V{\jmath}VV @VV{\jmath}V \\
\Op^\Pmf_{\cG,\cP^-}
@>{\iota}>> \Op^\Mmf_{\cG,\cP^-} \\
@V{\sfp^{-,\on{spec}}}VV \\
\Op^{\on{mon-free}}_{\cG}. 
\endCD
\end{equation} 

The composite vertical arrow in \eqref{e:MOpG OpG diag} is the map $\sfp^{\on{Miu,mon-free}}$. 

\sssec{}

Since
$$\LS^\Mmf_{\cP^-,\ul{x}}\simeq
(\cn^-_P/\on{Ad}(\cP^-))^{\times |\ul{x}|},$$
it is easy to see that $\Op^\Mmf_{\cG,\cP^-,\ul{x}}$
is an ind-scheme that admits a finite open cover
by ind-affine ind-placid ind-schemes.

\medskip

Hence, by \cite[Sect. A.10.12]{GLC2}, the functor
$$\iota^\IndCoh_*:\IndCoh^!(\Op^\Pmf_{\cG,\cP^-,\ul{x}})\to
\IndCoh^!(\Op^\Mmf_{\cG,\cP^-,\ul{x}})$$
is well-defined, and we have an isomorphism
\begin{equation} \label{e:base change GOP}
(\sfq^{\on{Miu,mon-free}})^\IndCoh_*\circ \jmath^!\simeq \jmath^!\circ \iota^\IndCoh_*,
\end{equation}
as functors
$$\IndCoh^!(\Op^\Pmf_{\cG,\cP^-,\ul{x}})\rightrightarrows
\IndCoh^!(\Op^\Mmf_{\cG,\cP^-,\ul{x}}).$$

\sssec{} \label{sss:semiinf Jaq opers}

We let 
$$J^{-,\on{spec},\semiinf}_\Op:\IndCoh^!(\Op^\Mmf_{\cG,\cP^-,\ul{x}}) \to \IndCoh^!(\Op^\mf_{\cM,\rhoch_P,\ul{x}})$$
be the functor
$$\IndCoh^!(\Op^\Mmf_{\cG,\cP^-,\ul{x}})\overset{\jmath^!}\longrightarrow 
\IndCoh^!(\Op^\mf_{\cM,\rhoch_P,\ul{x}}) \overset{(\sfq^{\on{Miu,mon-free}})^\IndCoh_*}\simeq \IndCoh^!(\Op^\mf_{\cM,\rhoch_P,\ul{x}}),$$
and refer to it as the \emph{semi-infinite} Jacquet functor for opers.

\medskip

By \eqref{e:base change GOP}, we have:
\begin{equation} \label{e:Jacquet opers via semiinf}
J^{-,\on{spec}}_\Op\simeq J^{-,\on{spec},\semiinf}_\Op \circ \iota^\IndCoh_*\circ (\sfp^{-,\on{spec}})^!.
\end{equation}

\begin{rem} 

Generalizing \cite{FG2}, we conjecture\footnote{In fact, we believe that this conjecture can be proved by the methods of {\it loc. cit.}}
that the category $\IndCoh^!(\Op^\Mmf_{\cG,\cP^-,\ul{x}})$
is equivalent to the category 
$$((\hg\mod_{\crit_G})_{\fL(N^-_P)}^{\fL^+(M)})_{\ul{x}},$$
so that the following diagram commutes:
$$
\CD
((\hg\mod_{\crit_G})_{\fL(N^-_P)}^{\fL^+(M)})_{\ul{x}} @>>> \IndCoh^!(\Op^\Mmf_{\cG,\cP^-,\ul{x}}) \\
@V{\ol{J}^-_{\on{KM}}}VV @VV{J^{-,\on{spec},\semiinf}_\Op}V \\
 \KL(M)_{\crit_M+\rhoch_P,\ul{x}} @>{\Theta^{-1}_{\Op^\mf_\cM}\circ \FLE_{M,\crit_M}}>> \IndCoh^!(\Op^\mf_{\cM,\rhoch_P,\ul{x}}). 
\endCD
$$

\end{rem} 

\sssec{} \label{sss:no OpMmf}

However, we do \emph{not} know how to make sense of 
$$\IndCoh^!(\Op^\Mmf_{\cG,\cP^-})$$
in the factorization setting. 

\medskip

What we will do instead is replace $\Op^\Mmf_{\cG,\cP^-}$ by its formal completion
along the map $\iota$ of \eqref{e:iota Op}. 
This will place us in a context where we can make sense of the categories involved in the factorization setting.  

\ssec{The semi-infinite category for opers} \label{ss:form coml Op GP}

Let 
$$(\Op^\Mmf_{\cG,\cP^-})^\wedge_\mf$$
be the formal completion of $\Op^\Mmf_{\cG,\cP^-}$
along the map $\iota$ of \eqref{e:iota Op}. 

\medskip

The goal of this subsection is to define unital factorization categories 
$$\IndCoh^*((\Op^\Mmf_{\cG,\cP^-})^\wedge_\mf) \text{ and } 
\IndCoh^!((\Op^\Mmf_{\cG,\cP^-})^\wedge_\mf).$$

We will do so by mimicking the procedure in \secref{ss:ambient spec semiinf}. 

\sssec{} \label{sss:Op mf triv}

Consider the fiber product
\begin{equation} \label{e:Op mf triv}
\on{pt}\underset{\LS^\reg_\cG}\times\Op_\cG^\mf.
\end{equation} 

This is a factorization ind-affine ind-scheme. It is \emph{not} ind-placid. Yet, we claim that the factorization
categories 
\begin{equation} \label{e:Op mf triv IndCoh}
\IndCoh^*(\on{pt}\underset{\LS^\reg_\cG}\times\Op_\cG^\mf) \text{ and } \IndCoh^!(\on{pt}\underset{\LS^\reg_\cG}\times\Op_\cG^\mf)
\end{equation} 
are well-defined and mutually dual. 

\medskip

Indeed, this follows by the same principle as in \cite[Sect. E.2.11]{GLC2}. Namely, by construction, \eqref{e:Op mf triv} carries an action
of $\fL^+_\nabla(\cG)$, so that
$$\fL^+_\nabla(\cG)\backslash (\on{pt}\underset{\LS^\reg_\cG}\times\Op_\cG^\mf) \simeq \Op_\cG^\mf,$$
while 
$$\IndCoh^*(\Op_\cG^\mf) \text{ and } \IndCoh^!(\Op_\cG^\mf)$$
are well-defined and mutually dual, since $\Op_\cG^\mf$ is placid. 

\sssec{}

Note now that we can rewrite \eqref{e:Op mf triv} also as
$$\on{pt}\underset{\LS^\mer_\cG}\times\Op_\cG^\mer.$$

\medskip

Hence, we obtain that the above $\fL^+_\nabla(\cG)$-action on \eqref{e:Op mf triv} extends to an action of 
$\fL_\nabla(\cG)$. From here we obtain an action of $\IndCoh^*(\fL_\nabla(\cG))$ on $\IndCoh^*(\on{pt}\underset{\LS^\reg_\cG}\times\Op_\cG^\mf)$
and a coaction of $\IndCoh^!(\fL_\nabla(\cG))$ on $\IndCoh^!(\on{pt}\underset{\LS^\reg_\cG}\times\Op_\cG^\mf)$. 

\medskip

In this way we obtain an action of $$\Sph_\cG^{\on{spec}}:=\IndCoh^*(\fL^+_\nabla(\cG)\backslash \fL_\nabla(\cG)/\fL^+_\nabla(\cG))$$
on $\IndCoh^*(\Op_\cG^\mf)$ and a coaction of $$(\Sph_\cG^{\on{spec}})^\vee\simeq \IndCoh^!(\fL^+_\nabla(\cG)\backslash \fL_\nabla(\cG)/\fL^+_\nabla(\cG))$$
on $\IndCoh^!(\Op_\cG^\mf)$; this is what the construction from \cite[Sect. E.8]{GLC2} amounts to.

\medskip

The above structures on the $\IndCoh^*(-)$ and $\IndCoh^!(-)$ versions are mutually dual in the natural sense.

\sssec{}

Recall the factorization space 
$$\Op^\GMmf_{\cG,\cP^-}=\LS^\reg_\cG\underset{\LS^\mer_\cG}\times \Op^\Mmf_{\cG,\cP^-},$$
and consider the fiber product
$$\on{pt}\underset{\LS^\reg_\cG}\times \Op^\GMmf_{\cG,\cP^-}\simeq
\on{pt}\underset{\LS^\mer_\cG}\times \Op^\Mmf_{\cG,\cP^-};$$
this is a factorization space acted on by $\fL_\nabla(\cG)$.

\medskip

Note that we can identify
$$\on{pt}\underset{\LS^\mer_\cG}\times \Op^\Mmf_{\cG,\cP^-} \simeq 
\Bigl(\on{pt}\underset{\LS^\reg_\cG}\times\Op_\cG^\mf\Bigr) \times \Bigl(\fL_\nabla(\cG/\cN^-_P)/\fL^+_\nabla(\cM)\Bigr).$$

Hence, we can identify 
$(\Op^\Mmf_{\cG,\cP^-})^\wedge_\mf$ with the \'etale quotient 
$$\left(\Bigl(\on{pt}\underset{\LS^\reg_\cG}\times\Op_\cG^\mf\Bigr) \times \Bigl(\fL_\nabla(\cG/\cN^-_P)/\fL^+_\nabla(\cM)\Bigr)\right)/\fL_\nabla(\cG).$$
We will use this presentation to access the corresponding $\IndCoh^*(-)$ and $\IndCoh^!(-)$ categories, cf. \secref{sss:ambient spec semiinf}

\sssec{}

By the above, we have a well-defined (unital) factorization category
\begin{multline*}
\IndCoh^*\left(\Bigl(\on{pt}\underset{\LS^\reg_\cG}\times\Op_\cG^\mf\Bigr) \times \Bigl(\fL_\nabla(\cG/\cN^-_P)/\fL^+_\nabla(\cM)\Bigr)\right):= \\
:=\IndCoh^*\left(\Bigl(\on{pt}\underset{\LS^\reg_\cG}\times\Op_\cG^\mf\Bigr)\right) \otimes
\IndCoh^*\left(\Bigl(\fL_\nabla(\cG/\cN^-_P)/\fL^+_\nabla(\cM)\Bigr)\right),
\end{multline*}
equipped with the (diagonal) action of the monoidal factorization category $\IndCoh^*(\fL_\nabla(\cG))$.

\medskip

We set
\begin{multline} \label{e:Op semminf * init}
\IndCoh^*((\Op^\Mmf_{\cG,\cP^-})^\wedge_\mf):=\\
:=\IndCoh^*\left(\Bigl(\on{pt}\underset{\LS^\reg_\cG}\times\Op_\cG^\mf\Bigr) \times \Bigl(\fL_\nabla(\cG/\cN^-_P)/\fL^+_\nabla(\cM)\Bigr)\right)
\underset{\IndCoh^*(\fL_\nabla(\cG))}\otimes \Vect.
\end{multline}

Similarly, we have a well-defined (unital factorization) category
\begin{multline*}
\IndCoh^!\left(\Bigl(\on{pt}\underset{\LS^\reg_\cG}\times\Op_\cG^\mf\Bigr) \times \Bigl(\fL_\nabla(\cG/\cN^-_P)/\fL^+_\nabla(\cM)\Bigr)\right):=\\
:=\IndCoh^!\left(\Bigl(\on{pt}\underset{\LS^\reg_\cG}\times\Op_\cG^\mf\Bigr)\right) \otimes
\IndCoh^!\left(\Bigl(\fL_\nabla(\cG/\cN^-_P)/\fL^+_\nabla(\cM)\Bigr)\right),
\end{multline*}
equipped with a coaction of the comonoidal factorization category $\IndCoh^!(\fL_\nabla(\cG))$.

\medskip

We set 
\begin{multline} \label{e:Op semminf ! init}
\IndCoh^!((\Op^\Mmf_{\cG,\cP^-})^\wedge_\mf):=\\
:=\IndCoh^!\left(\Bigl(\on{pt}\underset{\LS^\reg_\cG}\times\Op_\cG^\mf\Bigr) \times \Bigl(\fL_\nabla(\cG/\cN^-_P)/\fL^+_\nabla(\cM)\Bigr)\right)
\overset{\IndCoh^!(\fL_\nabla(\cG))}\otimes \Vect
\end{multline}
(see \secref{sss:tensor vs cotensor} for the $\overset{?}\otimes$ notation). 

\medskip

By construction, $\IndCoh^*((\Op^\Mmf_{\cG,\cP^-})^\wedge_\mf)$ and $\IndCoh^!((\Op^\Mmf_{\cG,\cP^-})^\wedge_\mf)$ are mutually dual unital factorization categories.

\sssec{}

Note that we can also write 
\begin{multline} \label{e:Op semminf *}
\IndCoh^*((\Op^\Mmf_{\cG,\cP^-})^\wedge_\mf)\simeq  \\
\simeq \IndCoh^*(\on{pt}\underset{\LS^\reg_\cG}\times\Op_\cG^\mf)\underset{\IndCoh^*(\fL_\nabla(\cG))}\otimes \IndCoh^*(\fL_\nabla(\cG/\cN^-_P))_{\fL^+_\nabla(\cM)}
\end{multline}
and
\begin{multline} \label{e:Op semminf !}
\IndCoh^!((\Op^\Mmf_{\cG,\cP^-})^\wedge_\mf)\simeq  \\
\simeq \IndCoh^!(\on{pt}\underset{\LS^\reg_\cG}\times\Op_\cG^\mf)\overset{\IndCoh^!(\fL_\nabla(\cG))}\otimes \IndCoh^!(\fL_\nabla(\cG/\cN^-_P))^{\fL^+_\nabla(\cM)}.
\end{multline}

And we can rewrite further: 

\begin{equation} \label{e:Op semiinf * bis}
\IndCoh^*((\Op^\Mmf_{\cG,\cP^-})^\wedge_\mf)
\simeq \IndCoh^*(\Op_\cG^\mf)\underset{\Sph^{\on{spec}}_\cG}\otimes \on{I}(\cG,\cP^-)^{\on{spec,loc}}
\end{equation} 
and 
\begin{equation} \label{e:Op semiinf ! bis}
\IndCoh^!((\Op^\Mmf_{\cG,\cP^-})^\wedge_\mf)
\simeq \IndCoh^!(\Op_\cG^\mf)\overset{(\Sph^{\on{spec}}_\cG)^\vee}\otimes \on{I}(\cG,\cP^-)^{\on{spec,loc}}_{\on{co}}, 
\end{equation} 
see \secref{sss:ambient via IGP}. 

\sssec{}

In particular, we have the (strictly unital) factorization functor
\begin{multline} \label{e:Op to Op semiinf *}
\IndCoh^*(\Op_\cG^\mf)\overset{\on{Id}\otimes \one_{\on{I}(\cG,\cP^-)^{\on{spec,loc}}}}\longrightarrow
\IndCoh^*(\Op_\cG^\mf)\otimes \on{I}(\cG,\cP^-)^{\on{spec,loc}} \to \\
\to \IndCoh^*(\Op_\cG^\mf)\underset{\Sph^{\on{spec}}_\cG}\otimes \on{I}(\cG,\cP^-)^{\on{spec,loc}}  
\simeq \IndCoh^*((\Op^\Mmf_{\cG,\cP^-})^\wedge_\mf)
\end{multline} 
which admits a right adjoint as a factorization functor, namely
\begin{multline} \label{Op semiinf * to Op}
\IndCoh^*((\Op^\Mmf_{\cG,\cP^-})^\wedge_\mf)= \\
=\IndCoh^*(\Op_\cG^\mf)\underset{\Sph^{\on{spec}}_\cG}\otimes \on{I}(\cG,\cP^-)^{\on{spec,loc}}  
\overset{\on{Id}\otimes (-\underset{\cG}\star \Delta^{-,\on{spec},\semiinf})^R}\longrightarrow \\
\to \IndCoh^*(\Op_\cG^\mf)\underset{\Sph^{\on{spec}}_\cG}\otimes \Sph^{\on{spec}}_\cG \simeq \IndCoh^*(\Op_\cG^\mf),
\end{multline} 
where
$(-\underset{\cG}\star \Delta^{-,\on{spec},\semiinf})^R$ is the right adjoint of the functor
$$\Sph_\cG^{\on{spec}}\overset{-\underset{\cG}\star \Delta^{-,\on{spec},\semiinf}}\longrightarrow \on{I}(\cG,\cP^-)^{\on{spec,loc}}.$$

\sssec{}

The dual of \eqref{Op semiinf * to Op} is a (strictly) unital factorization functor
\begin{equation} \label{e:Op to Op semiinf !}
\IndCoh^!(\Op_\cG^\mf)\to \IndCoh^!((\Op^\Mmf_{\cG,\cP^-})^\wedge_\mf).
\end{equation}

Note that the functor \eqref{e:Op to Op semiinf !} can also be interpreted as
\begin{multline} \label{e:Op to Op semiinf ! bis}
\IndCoh^!(\Op_\cG^\mf)\overset{\on{Id}\otimes \one_{\on{I}(\cG,\cP^-)^{\on{spec,loc}}_{\on{co}}}}\longrightarrow
\IndCoh^!(\Op_\cG^\mf)\otimes \on{I}(\cG,\cP^-)^{\on{spec,loc}}_{\on{co}}\to \\
\to \IndCoh^!(\Op_\cG^\mf)\overset{(\Sph^{\on{spec}}_\cG)^\vee}\otimes \on{I}(\cG,\cP^-)^{\on{spec,loc}}_{\on{co}}  
\simeq \IndCoh^!((\Op^\Mmf_{\cG,\cP^-})^\wedge_\mf),
\end{multline} 
where the second arrow is the left adjoint of the forgetful functor
$$\IndCoh^!(\Op_\cG^\mf)\overset{(\Sph^{\on{spec}}_\cG)^\vee}\otimes \on{I}(\cG,\cP^-)^{\on{spec,loc}}_{\on{co}}  \to
\IndCoh^!(\Op_\cG^\mf)\otimes \on{I}(\cG,\cP^-)^{\on{spec,loc}}_{\on{co}},$$
which is well-defined as a factorization functor, since $\Sph^{\on{spec}}_\cG$ is rigid. 

\sssec{}

Consider now the fiber product
$$\Op^\Pmf_{\cG,\cP^-}\simeq
\Op^\mf_\cG\underset{\LS^\reg_\cG}\times \LS^\reg_{\cP^-}.$$

By the same principle as in \secref{sss:Op mf triv}, we have well-defined unital factorization categories 
$$\IndCoh^*(\Op^\Pmf_{\cG,\cP^-}) \text{ and } 
\IndCoh^!(\Op^\Pmf_{\cG,\cP^-}).$$

\medskip

Namely,
$$\IndCoh^*(\Op^\Pmf_{\cG,\cP^-}) :=
\IndCoh^*(\on{pt}\underset{\LS^\reg_\cG}\times \Op^\mf_\cG)_{\fL^+_\nabla(\cP^-)}$$
and
$$\IndCoh^!(\Op^\Pmf_{\cG,\cP^-}) :=
\IndCoh^!(\on{pt}\underset{\LS^\reg_\cG}\times \Op^\mf_\cG)^{\fL^+_\nabla(\cP^-)}$$

\medskip

We have the naturally defined factorization functors
$$(\sfp^{-,\on{spec}})^\IndCoh_*:\IndCoh^*(\Op^\Pmf_{\cG,\cP^-})\to \IndCoh^*(\Op^\mf_\cG)$$
and 
$$(\sfp^{-,\on{spec}})^!:\IndCoh^!(\Op^\mf_\cG)\to \IndCoh^!(\Op^\Pmf_{\cG,\cP^-}).$$

\medskip

Moreover, the functor $(\sfp^{-,\on{spec}})^\IndCoh_*$ admits a left adjoint as a factorization functor
$$(\sfp^{-,\on{spec}})^*:\IndCoh^*(\Op^\mf_\cG)\to \IndCoh^*(\Op^\Pmf_{\cG,\cP^-}).$$

The above functors carry lax unital structures, and for $(\sfp^{-,\on{spec}})^*$ and $(\sfp^{-,\on{spec}})^!$
this structure is strict. 

\sssec{}

Let ${}'\!\iota$ denote the map
$$\Op^\Pmf_{\cG,\cP^-}\to (\Op^\Mmf_{\cG,\cP^-})^\wedge_\mf.$$

\medskip

We claim that we have adjoint pairs
\begin{equation} \label{e:iota wedge *}
({}'\!\iota)^\IndCoh_*:
\IndCoh^*(\Op^\Pmf_{\cG,\cP^-})\rightleftarrows 
\IndCoh^*((\Op^\Mmf_{\cG,\cP^-})^\wedge_\mf):({}'\!\iota)^!
\end{equation} 
and 
\begin{equation} \label{e:iota wedge !}
({}'\!\iota)^\IndCoh_*:
\IndCoh^!(\Op^\Pmf_{\cG,\cP^-})\rightleftarrows 
\IndCoh^!((\Op^\Mmf_{\cG,\cP^-})^\wedge_\mf):({}'\!\iota)^!.
\end{equation} 

Namely, we interpret the functor $({}'\!\iota)^\IndCoh_*$ in \eqref{e:iota wedge *} as 
\begin{multline*}
\IndCoh^*(\on{pt}\underset{\LS^\reg_\cG}\times\Op^\mf_\cG)\underset{\IndCoh^*(\fL^+_\nabla(\cP^-))}\otimes \Vect \simeq \\
\simeq \IndCoh^*(\on{pt}\underset{\LS^\reg_\cG}\times\Op^\mf_\cG)\underset{\IndCoh^*(\fL^+_\nabla(\cG))}\otimes
\IndCoh^*(\fL^+_\nabla(G/\cN^-_P))_{\fL^+_\nabla(\cM)}\to \\
\to \IndCoh^*(\on{pt}\underset{\LS^\reg_\cG}\times\Op^\mf_\cG)\underset{\IndCoh^*(\fL^+_\nabla(\cG))}\otimes
\IndCoh^*(\fL_\nabla(G/\cN^-_P))_{\fL^+_\nabla(\cM)}\to  \\
\to \IndCoh^*(\on{pt}\underset{\LS^\reg_\cG}\times\Op^\mf_\cG)\underset{\IndCoh^*(\fL_\nabla(\cG))}\otimes
\IndCoh^*(\fL_\nabla(G/\cN^-_P))_{\fL^+_\nabla(\cM)},
\end{multline*}
and it admits a right adjoint as a factorization functor by \secref{sss:right adj coinv}. 

\medskip

The functors in \eqref{e:iota wedge !} are set to be the duals of ones in \eqref{e:iota wedge *}.

\sssec{} \label{sss:factor Jacquet}

Unwinding the definitions, we obtain that the functor \eqref{e:Op to Op semiinf *} identifies with
$$\IndCoh^*(\Op_\cG^\mf) \overset{(\sfp^{-,\on{spec}})^*}\longrightarrow 
\IndCoh^*(\Op^\Pmf_{\cG,\cP^-}) \overset{({}'\!\iota)^\IndCoh_*}\longrightarrow \IndCoh^*((\Op^\Mmf_{\cG,\cP^-})^\wedge_\mf),$$
and hence its right adjoint, i.e., the functor \eqref{Op semiinf * to Op}, identifies with
$$\IndCoh^*((\Op^\Mmf_{\cG,\cP^-})^\wedge_\mf) 
\overset{({}'\!\iota)^!}\longrightarrow 
\IndCoh^*(\Op^\Pmf_{\cG,\cP^-})  \overset{(\sfp^{-,\on{spec}})_*}\longrightarrow \IndCoh^*(\Op_\cG^\mf).$$

Furthermore, the functor \eqref{e:Op to Op semiinf !} identifies with 
$$\IndCoh^!(\Op_\cG^\mf) \overset{(\sfp^{-,\on{spec}})^!}\longrightarrow 
\IndCoh^!(\Op^\Pmf_{\cG,\cP^-}) \overset{({}'\!\iota)^\IndCoh_*}\longrightarrow \IndCoh^!((\Op^\Mmf_{\cG,\cP^-})^\wedge_\mf).$$

\ssec{A formal completion of Miura opers}

The semi-infinite Jacquet functor for opers goes via the category of ind-coherent sheaves on the formal completion of
the space of Miura opers along the monodromy-free locus. 

\medskip

We study this category in the present subsection. 

\sssec{} \label{sss:form compl Miura}

Let
$$(\MOp^\Mmf_{\cG,\cP^-})^\wedge_\mf$$
denote the formal completion of 
$$\Op^\mf_{\cM,\rhoch_P}\simeq \MOp^\Mmf_{\cG,\cP^-}:=\MOp^\mer_{\cG,\cP^-}\underset{\LS_\cM^\mer}\times \LS_\cM^\reg$$
along the map 
\begin{equation} \label{e:G-mf MOp}
\MOp^\Pmf_{\cG,\cP^-} \overset{\sfq^{\on{Miu,mon-free}}}\longrightarrow \Op^\mf_{\cM,\rhoch_P},
\end{equation} 
i.e., the map 
$$\MOp^\Pmf_{\cG,\cP^-}\to
\MOp^\Mmf_{\cG,\cP^-}.$$

\medskip

Note that the image of the above map coincides set-theoretically with the preimage of
\begin{equation} \label{e:LS G reg ro mer}
\LS_\cG^\reg\subset \LS_\cG^\mer
\end{equation} 
along the map 
\begin{equation} \label{e:Op M to LS G}
\MOp^\Mmf_{\cG,\cP^-} \to \MOp^\mer_{\cG,\cP^-}\overset{\sfp^{\on{Miu}}}\to \Op^\mer_\cG\to \LS_\cG^\mer.
\end{equation} 

\sssec{}

Let $(\sfq^{\on{Miu,mon-free}})^\wedge$ denote the ind-closed embedding
$$(\MOp^\Mmf_{\cG,\cP^-})^\wedge_\mf\to \Op^\mf_{\cM,\rhoch_P}.$$

\medskip

Since \eqref{e:LS G reg ro mer} is locally almost of finite presentation (see \cite[Lemma B.7.13]{GLC2}), we obtain 
that the map $(\sfq^{\on{Miu,mon-free}})^\wedge$ is also locally almost of finite presentation. 

\medskip

In particular, since $\Op^\mf_{\cM,\rhoch_P}$ is ind-placid, we obtain that $(\MOp^\Mmf_{\cG,\cP^-})^\wedge_\mf$ is
also ind-placid (see \cite[Corollary A.9.10]{GLC2}). Hence, we have well-defined unital factorization categories
\begin{equation} \label{e:IndCoh Miu compl}
\IndCoh^*\left((\MOp^\Mmf_{\cG,\cP^-})^\wedge_\mf\right) \text{ and } \IndCoh^!\left((\MOp^\Mmf_{\cG,\cP^-})^\wedge_\mf\right)
\end{equation}
and pairs of adjoint functors:
\begin{multline} \label{e:q Miu adj *}
((\sfq^{\on{Miu,mon-free}})^\wedge)^\IndCoh_*:\IndCoh^*\left((\MOp^\Mmf_{\cG,\cP^-})^\wedge_\mf\right)\rightleftarrows \\
\rightleftarrows \IndCoh^*(\Op^\mf_{\cM,\rhoch_P}):((\sfq^{\on{Miu,mon-free}})^\wedge)^!
\end{multline}
and 
\begin{multline} \label{e:q Miu adj !}
((\sfq^{\on{Miu,mon-free}})^\wedge)^\IndCoh_*:\IndCoh^!\left((\MOp^\Mmf_{\cG,\cP^-})^\wedge_\mf\right)\rightleftarrows \\
\rightleftarrows \IndCoh^!(\Op^\mf_{\cM,\rhoch_P}):((\sfq^{\on{Miu,mon-free}})^\wedge)^!
\end{multline}
with the left adjoints being strictly unital. 

\sssec{}

Let $'\!\iota$ denote the map
$$\MOp^\Pmf_{\cG,\cP^-}\to (\MOp^\Mmf_{\cG,\cP^-})^\wedge_\mf.$$

We obtain that the adjunctions
$$(\sfq^{\on{Miu,mon-free}})^\IndCoh_*:\IndCoh^*(\MOp^\Pmf_{\cG,\cP^-})\rightleftarrows \IndCoh^*(\Op^\mf_{\cM,\rhoch_P}):
(\sfq^{\on{Miu,mon-free}})^!$$
and
$$(\sfq^{\on{Miu,mon-free}})^\IndCoh_*:\IndCoh^!(\MOp^\Pmf_{\cG,\cP^-})\rightleftarrows \IndCoh^!(\Op^\mf_{\cM,\rhoch_P}):
(\sfq^{\on{Miu,mon-free}})^!$$
factor as compositions of 
\begin{equation} \label{e:iota' MOp}
'\!\iota^\IndCoh_*:\IndCoh^*(\MOp^\Pmf_{\cG,\cP^-})\rightleftarrows \IndCoh^*((\MOp^\Mmf_{\cG,\cP^-})^\wedge_\mf):
{}'\!\iota^!
\end{equation} 
and
\begin{equation} \label{e:iota' MOp !}
'\!\iota^\IndCoh_*:\IndCoh^!(\MOp^\Pmf_{\cG,\cP^-})\rightleftarrows \IndCoh^!((\MOp^\Mmf_{\cG,\cP^-})^\wedge_\mf):
{}'\!\iota^!
\end{equation} 
with the adjunctions \eqref{e:q Miu adj *} and \eqref{e:q Miu adj !}, respectively. 

\sssec{}

We now claim that the categories \eqref{e:IndCoh Miu compl} admit descriptions similar to those in \eqref{e:Op semminf *} and 
\eqref{e:Op semminf !}, respectively. 
 
\medskip

Namely, recall the map \eqref{e:Op M to LS G}, and consider the fiber products 
$$\LS^\reg_\cG\underset{\LS^\mer_\cG}\times\MOp^\Mmf_{\cG,\cP^-}\text{ and }\on{pt}\underset{\LS^\mer_\cG}\times\MOp^\Mmf_{\cG,\cP^-}.$$

Both are ind-affine factorization ind-schemes, with $\LS^\reg_\cG\underset{\LS^\mer_\cG}\times\MOp^\Mmf_{\cG,\cP^-}$ being placid, so that
$$\IndCoh^*(\LS^\reg_\cG\underset{\LS^\mer_\cG}\times\MOp^\Mmf_{\cG,\cP^-}) \text{ and }
\IndCoh^!(\LS^\reg_\cG\underset{\LS^\mer_\cG}\times\MOp^\Mmf_{\cG,\cP^-})$$
are well-defined as (mutually dual) unital factorization categories.

\medskip

The ind-scheme $\on{pt}\underset{\LS^\mer_\cG}\times\MOp^\Mmf_{\cG,\cP^-}$ carries an action of $\fL_\nabla(\cG)$ and
$$\fL^+_\nabla(\cG)\backslash (\on{pt}\underset{\LS^\mer_\cG}\times\MOp^\Mmf_{\cG,\cP^-})\simeq
\LS^\reg_\cG\underset{\LS^\mer_\cG}\times\MOp^\Mmf_{\cG,\cP^-}.$$

Hence, by the same mechanism as in \secref{sss:Op mf triv}, the categories
$$\IndCoh^*(\on{pt}\underset{\LS^\mer_\cG}\times\MOp^\Mmf_{\cG,\cP^-}) \text{ and }
\IndCoh^!(\on{pt}\underset{\LS^\mer_\cG}\times\MOp^\Mmf_{\cG,\cP^-})$$
are well-defined as (mutually dual) unital factorization categories. 

\medskip

Moreover, $\IndCoh^*(\on{pt}\underset{\LS^\mer_\cG}\times\MOp^\Mmf_{\cG,\cP^-})$ carries an action of the monoidal factorization
category $\IndCoh^*(\fL_\nabla(\cG))$,
and $\IndCoh^!(\on{pt}\underset{\LS^\mer_\cG}\times\MOp^\Mmf_{\cG,\cP^-})$ carries a coaction of the comonoidal 
factorization category $\IndCoh^!(\fL_\nabla(\cG))$. 

\sssec{}

Let $''\!\iota$ denote the map
$$\MOp^\GMmf_{\cG,\cP^-}:=\LS^\reg_\cG\underset{\LS^\mer_\cG}\times\MOp^\Mmf_{\cG,\cP^-}\to (\MOp^\Mmf_{\cG,\cP^-})^\wedge_\mf.$$

\medskip

We have the adjunctions
$$({}''\!\iota)^\IndCoh_*:\IndCoh^*(\MOp^\GMmf_{\cG,\cP^-}) \rightleftarrows 
\IndCoh^*((\MOp^\Mmf_{\cG,\cP^-})^\wedge_\mf):({}''\!\iota)^!$$
and
$$({}''\!\iota)^\IndCoh_*:\IndCoh^!(\MOp^\GMmf_{\cG,\cP^-}) \rightleftarrows 
\IndCoh^!((\MOp^\Mmf_{\cG,\cP^-})^\wedge_\mf):({}''\!\iota)^!.$$

Note now that the functor
\begin{multline}  \label{e:MOp form comp pre}
\IndCoh^*(\on{pt}\underset{\LS^\mer_\cG}\times\MOp^\Mmf_{\cG,\cP^-}) \to
\IndCoh^*(\on{pt}\underset{\LS^\mer_\cG}\times\MOp^\Mmf_{\cG,\cP^-}) \underset{\IndCoh^*(\fL^+_\nabla(\cG))}\otimes \Vect\simeq \\
\simeq \IndCoh^*(\MOp^\GMmf_{\cG,\cP^-}) \overset{({}''\!\iota)^\IndCoh_*}\longrightarrow
\IndCoh^*((\MOp^\Mmf_{\cG,\cP^-})^\wedge_\mf)
\end{multline}
naturally factors as
\begin{multline*}
\IndCoh^*(\on{pt}\underset{\LS^\mer_\cG}\times\MOp^\Mmf_{\cG,\cP^-}) \to \\
\to \IndCoh^*(\on{pt}\underset{\LS^\mer_\cG}\times\MOp^\Mmf_{\cG,\cP^-}) \underset{\IndCoh^*(\fL_\nabla(\cG))}\otimes \Vect \to
\IndCoh^*((\MOp^\Mmf_{\cG,\cP^-})^\wedge_\mf).
\end{multline*}

\sssec{}

We claim:

\begin{prop} \label{p:MOp form comp}
The  functor
\begin{equation} \label{e:MOp form comp *}
\IndCoh^*(\on{pt}\underset{\LS^\mer_\cG}\times\MOp^\Mmf_{\cG,\cP^-}) \underset{\IndCoh^*(\fL_\nabla(\cG))}\otimes \Vect \to
\IndCoh^*((\MOp^\Mmf_{\cG,\cP^-})^\wedge_\mf)
\end{equation} 
is an equivalence.
\end{prop} 

\begin{proof}

First, we claim that the functor \eqref{e:MOp form comp *} preserves compactness. To prove this, we note that its precomposition
with 
\begin{multline} \label{e:MOp form comp 1}
\IndCoh^*(\MOp^\GMmf_{\cG,\cP^-}) \simeq
\IndCoh^*(\on{pt}\underset{\LS^\mer_\cG}\times\MOp^\Mmf_{\cG,\cP^-}) \underset{\IndCoh^*(\fL^+_\nabla(\cG))}\otimes \Vect\to \\
\to \IndCoh^*(\on{pt}\underset{\LS^\mer_\cG}\times\MOp^\Mmf_{\cG,\cP^-}) \underset{\IndCoh^*(\fL_\nabla(\cG))}\otimes \Vect,
\end{multline} 
which is the functor 
$$({}''\!\iota)^\IndCoh_*:\IndCoh^*(\MOp^\GMmf_{\cG,\cP^-}) \to
\IndCoh^*((\MOp^\Mmf_{\cG,\cP^-})^\wedge_\mf)$$
preserves compactness. 

\medskip

However, the functor \eqref{e:MOp form comp 1} itself also preserves compactness, and 
 its essential image generates the target. 

\medskip

Hence, by \cite[Proposition 6.2.6]{GLC2}, it suffices to show that the functor \eqref{e:MOp form comp *} is a pointwise equivalence
(i.e., induces an equivalence over each $\ul{x}\in \Ran$). 

\medskip

We rewrite the resulting functor 
$$\IndCoh^*(\Op^\mf_{\cM,\rhoch_P,\ul{x}}\underset{\LS^\mer_{\cG,\ul{x}}}\times \on{pt}) 
\underset{\IndCoh^*(\fL_\nabla(\cG)_{\ul{x}})}\otimes \Vect \to
\IndCoh^*((\MOp_{\cG,\cP^-,\ul{x}})^\wedge_\mf)$$
as
$$\IndCoh^*(\Op^\mf_{\cM,\rhoch_P,\ul{x}}\underset{\LS^\mer_{\cG,\ul{x}}}\times \LS^\reg_{\cG,\ul{x}}) 
\underset{\Sph^{\on{spec}}_{\cG,\ul{x}}}\otimes \Rep(\cG)_{\ul{x}}\to \IndCoh^*((\MOp_{\cG,\cP^-,\ul{x}})^\wedge_\mf),$$
and the required assertion follows by the same mechanism as \cite[Proposition 3.6.5]{GLC2}. 

\end{proof} 

\sssec{}

Similarly to \eqref{e:MOp form comp pre}, the functor
\begin{multline*} 
({}''\!\iota)^!:\IndCoh^!((\MOp^\Mmf_{\cG,\cP^-})^\wedge_\mf)\overset{({}''\!\iota)^!}\longrightarrow 
\IndCoh^!(\MOp^\GMmf_{\cG,\cP^-}) \simeq \\
\simeq 
\IndCoh^!(\on{pt}\underset{\LS^\mer_\cG}\times\MOp^\Mmf_{\cG,\cP^-}) \overset{\IndCoh^!(\fL^+_\nabla(\cG))}\otimes \Vect\to
\IndCoh^!(\on{pt}\underset{\LS^\mer_\cG}\times\MOp^\Mmf_{\cG,\cP^-})
\end{multline*}
naturally factors as
\begin{multline*} 
\IndCoh^!((\MOp^\Mmf_{\cG,\cP^-})^\wedge_\mf)\to 
\IndCoh^!(\on{pt}\underset{\LS^\mer_\cG}\times\MOp^\Mmf_{\cG,\cP^-}) \overset{\IndCoh^!(\fL_\nabla(\cG))}\otimes \Vect\to \\
\to \IndCoh^!(\on{pt}\underset{\LS^\mer_\cG}\times\MOp^\Mmf_{\cG,\cP^-}).
\end{multline*}

As a formal consequence of \propref{p:MOp form comp}, we obtain:

\begin{cor} \label{c:MOp form comp}
The  functor
\begin{equation} \label{e:MOp form comp !}
\IndCoh^!((\MOp^\Mmf_{\cG,\cP^-})^\wedge_\mf)\to 
\IndCoh^!(\on{pt}\underset{\LS^\mer_\cG}\times\MOp^\Mmf_{\cG,\cP^-}) \overset{\IndCoh^!(\fL_\nabla(\cG))}\otimes \Vect
\end{equation} 
is an equivalence.
\end{cor} 

\ssec{The semi-infinite Jacquet functor for opers, factorization version} \label{ss:semiinf op}

Recall the functor $J^{-,\on{spec},\semiinf}_\Op$ from \secref{sss:semiinf Jaq opers}. The goal of this subsection
is to introduce its variant, denoted 
$$J^{-,\on{spec},\semiinf_\mf}_\Op,$$
which will be a (lax unital) factorization functor 
$$\IndCoh^!((\Op^\Mmf_{\cG,\cP^-})^\wedge_\mf)\to \IndCoh^!(\Op^\mf_{\cM,\rhoch_P}),$$
as well as its dual, denoted
$$\on{co}\!J^{-,\on{spec},\semiinf_\mf}_\Op,$$
which will be (lax unital) factorization functor
$$\IndCoh^*(\Op^\mf_{\cM,\rhoch_P})\to  
\IndCoh^*((\Op^\Mmf_{\cG,\cP^-})^\wedge_\mf).$$

\medskip

Intuitively, the functor $J^{-,\on{spec},\semiinf_\mf}_\Op$ is pullback along the map
$$(\MOp^\Mmf_{\cG,\cP^-})^\wedge_\mf \to (\Op^\Mmf_{\cG,\cP^-})^\wedge_\mf,$$
induced by the map 
\begin{equation} \label{e:MOp to Op again}
\jmath:\MOp^\Mmf_{\cG,\cP^-} \to \Op^\Mmf_{\cG,\cP^-},
\end{equation} 
followed by pushforward along 
$$(\MOp^\Mmf_{\cG,\cP^-})^\wedge_\mf \to \MOp^\Mmf_{\cG,\cP^-}\overset{\sfq^{\on{Miu}}}\simeq \Op^\mf_{\cM,\rhoch_P}.$$

\sssec{}

Note that the map $\jmath$ of \eqref{e:MOp to Op again} induces a map, denoted by the same symbol $\jmath$, 
\begin{equation} \label{e:J G-mf}
\jmath:\on{pt}\underset{\LS^\mer_\cG}\times\MOp^\Mmf_{\cG,\cP^-}\to 
\on{pt}\underset{\LS^\mer_\cG}\times\Op^\Mmf_{\cG,\cP^-}
\end{equation}
as factorization prestacks, compatible with the actions of $\fL_\nabla(\cG)$.  

\begin{rem}
Note that the map \eqref{e:J G-mf} is \emph{not} an open embedding, see Remark \ref{r:j not open} above.
\end{rem}

\sssec{}

We have a pair of mutually dual (lax unital) factorization functors 
$$\jmath_*:\IndCoh^*(\on{pt}\underset{\LS^\mer_\cG}\times\MOp^\Mmf_{\cG,\cP^-})\to 
\IndCoh^*(\on{pt}\underset{\LS^\mer_\cG}\times\Op^\Mmf_{\cG,\cP^-})$$
and 
$$\jmath^!:\IndCoh^!(\on{pt}\underset{\LS^\mer_\cG}\times\Op^\Mmf_{\cG,\cP^-})\to
\IndCoh^!(\on{pt}\underset{\LS^\mer_\cG}\times\MOp^\Mmf_{\cG,\cP^-}).$$

They induce a pair of (lax unital) factorization functors
\begin{multline} \label{e:from Op to MOp *}
\jmath_*:
\IndCoh^*(\on{pt}\underset{\LS^\mer_\cG}\times\MOp^\Mmf_{\cG,\cP^-})\underset{\IndCoh^*(\fL_\nabla(\cG))}\otimes \Vect \to \\
\to \IndCoh^*(\on{pt}\underset{\LS^\mer_\cG}\times\Op^\Mmf_{\cG,\cP^-})\underset{\IndCoh^*(\fL_\nabla(\cG))}\otimes \Vect
\end{multline} 
\begin{multline} \label{e:from MOp to Op !}
\jmath^!:\IndCoh^!(\on{pt}\underset{\LS^\mer_\cG}\times\Op^\Mmf_{\cG,\cP^-})\overset{\IndCoh^!(\fL_\nabla(\cG))}\otimes \Vect\to \\
\to \IndCoh^!(\on{pt}\underset{\LS^\mer_\cG}\times\MOp^\Mmf_{\cG,\cP^-})\overset{\IndCoh^!(\fL_\nabla(\cG))}\otimes \Vect,
\end{multline} 
respectively.

\sssec{}

Using the identifications
$$\IndCoh^*(\on{pt}\underset{\LS^\mer_\cG}\times\Op^\Mmf_{\cG,\cP^-})\underset{\IndCoh^*(\fL_\nabla(\cG))}\otimes \Vect =:
\IndCoh^*((\Op^\Mmf_{\cG,\cP^-})^\wedge_\mf)$$
and 
$$\IndCoh^!(\on{pt}\underset{\LS^\mer_\cG}\times\Op^\Mmf_{\cG,\cP^-})\overset{\IndCoh^!(\fL_\nabla(\cG))}\otimes \Vect =:
\IndCoh^!((\Op^\Mmf_{\cG,\cP^-})^\wedge_\mf)$$
and 
$$\IndCoh^*(\on{pt}\underset{\LS^\mer_\cG}\times\MOp^\Mmf_{\cG,\cP^-})\underset{\IndCoh^*(\fL_\nabla(\cG))}\otimes \Vect 
\overset{\text{\propref{p:MOp form comp}}}\simeq 
\IndCoh^*((\MOp^\Mmf_{\cG,\cP^-})^\wedge_\mf)$$
and 
$$\IndCoh^!(\on{pt}\underset{\LS^\mer_\cG}\times\MOp^\Mmf_{\cG,\cP^-})\overset{\IndCoh^!(\fL_\nabla(\cG))}\otimes \Vect
\overset{\text{\corref{c:MOp form comp}}}\simeq 
\IndCoh^!((\MOp^\Mmf_{\cG,\cP^-})^\wedge_\mf),$$
we can interpret the functors \eqref{e:from Op to MOp *} and \eqref{e:from MOp to Op !} as a pair of mutually dual 
(lax unital) factorization functors 
\begin{equation} \label{e:from Op to MOp * bis}
\jmath_*:\IndCoh^*((\MOp^\Mmf_{\cG,\cP^-})^\wedge_\mf)\to \IndCoh^*((\Op^\Mmf_{\cG,\cP^-})^\wedge_\mf),
\end{equation} 
and 
\begin{equation} \label{e:from Op to MOp ! bis}
\jmath^!:\IndCoh^!((\Op^\Mmf_{\cG,\cP^-})^\wedge_\mf)\to \IndCoh^!((\MOp^\Mmf_{\cG,\cP^-})^\wedge_\mf),
\end{equation}
respectively.

\sssec{}

We define the functor 
$$J^{-,\on{spec},\semiinf_\mf}_\Op:\IndCoh^!((\Op^\Mmf_{\cG,\cP^-})^\wedge_\mf)\to \IndCoh^!(\Op^\mf_{\cM,\rhoch_P})$$
to be the composition
\begin{multline} \label{e:defn J semiinf}
\IndCoh^!((\Op^\Mmf_{\cG,\cP^-})^\wedge_\mf)\overset{\jmath^!}\to \\
\to \IndCoh^!((\MOp^\Mmf_{\cG,\cP^-})^\wedge_\mf)
\overset{((\sfq^{\on{Miu,mon-free}})^\wedge)^\IndCoh_*}\longrightarrow \IndCoh^!(\Op^\mf_{\cM,\rhoch_P}).
\end{multline} 

We define the functor 
$$\on{co}\!J^{-,\on{spec},\semiinf_\mf}_\Op:\IndCoh^*(\Op^\mf_{\cM,\rhoch_P})\to \IndCoh^*((\Op^\Mmf_{\cG,\cP^-})^\wedge_\mf)$$
to be 
\begin{multline} \label{e:defn coJ semiinf}
\IndCoh^*(\Op^\mf_{\cM,\rhoch_P})\overset{((\sfq^{\on{Miu,mon-free}})^\wedge)^!}\longrightarrow  \\
\to \IndCoh^*((\MOp^\Mmf_{\cG,\cP^-})^\wedge_\mf)\overset{\jmath^\IndCoh_*}\to \IndCoh^*((\Op^\Mmf_{\cG,\cP^-})^\wedge_\mf).
\end{multline} 

\sssec{}

We have a commutative diagram of functors
\begin{equation} \label{e:BC Miura ! pre}
\CD
\IndCoh^!(\MOp^\Pmf_{\cG,\cP^-}) @<{'\!\iota^!}<< \IndCoh^!((\MOp^\Mmf_{\cG,\cP^-})^\wedge_\mf) \\
@A{\jmath^!}AA @AA{\jmath^!}A  \\
\IndCoh^!(\Op^\Pmf_{\cG,\cP^-}) @<{'\!\iota^!}<< \IndCoh^!((\Op^\Mmf_{\cG,\cP^-})^\wedge_\mf).
\endCD
\end{equation}

By adjunction, it gives rise to a natural transformation
\begin{equation} \label{e:BC Miura !}
({}'\!\iota)^\IndCoh_*\circ \jmath^!\to  \jmath^!\circ ({}'\!\iota)^\IndCoh_*
\end{equation} 
as functors 
$$\IndCoh^!(\Op^\Pmf_{\cG,\cP^-})\rightrightarrows \IndCoh^!((\MOp^\Mmf_{\cG,\cP^-})^\wedge_\mf).$$

Similarly, we have a commutative diagram of functors 
\begin{equation} \label{e:BC Miura * pre}
\CD
\IndCoh^*(\MOp^\Pmf_{\cG,\cP^-}) @>{({}'\!\iota)^\IndCoh_*}>> \IndCoh^*((\MOp^\Mmf_{\cG,\cP^-})^\wedge_\mf) \\
@V{\jmath^\IndCoh_*}VV  @V{\jmath^\IndCoh_*}VV  \\
\IndCoh^*(\Op^\Pmf_{\cG,\cP^-}) @>{({}'\!\iota)^\IndCoh_*}>> \IndCoh^*((\Op^\Mmf_{\cG,\cP^-})^\wedge_\mf).
\endCD
\end{equation}

By adjunction, it gives rise to a natural transformation
\begin{equation} \label{e:BC Miura *}
\jmath^\IndCoh_* \circ ({}'\!\iota)^!\to ({}'\!\iota)^!\circ \jmath^\IndCoh_*
\end{equation}
as functors
$$\IndCoh^*((\MOp^\Mmf_{\cG,\cP^-})^\wedge_\mf)\rightrightarrows \IndCoh^*(\Op^\Pmf_{\cG,\cP^-}).$$ 

\medskip

We claim:

\begin{lem} \label{l:BC Miura}
The natural transformations \eqref{e:BC Miura !} and \eqref{e:BC Miura *} are isomorphisms.
\end{lem}

\begin{proof}

We will prove the assertion for \eqref{e:BC Miura *}; the one for \eqref{e:BC Miura !} would then follow by duality.

\medskip

Consider the fiber products
$$\on{pt}\underset{\LS^\reg_{\cP^-}}\times \Op^\Pmf_{\cG,\cP^-}\simeq \on{pt}\underset{\LS^\mer_{\cP^-}}\times \Op^\mer_{\cG,\cP^-}$$
and 
$$\on{pt}\underset{\LS^\reg_{\cP^-}}\times \MOp^\Pmf_{\cG,\cP^-}\simeq \on{pt}\underset{\LS^\mer_{\cP^-}}\times \MOp^\mer_{\cG,\cP^-}.$$

Both are factorization ind-affine ind-schemes equipped with an action of $\fL_\nabla(\cP^-)$. We have a commutative diagram
\begin{equation} \label{e:BC Miura ! pre 1}
\CD
(\on{pt}\underset{\LS^\mer_{\cP^-}}\times \MOp^\mer_{\cG,\cP^-})\times \fL_\nabla(\cP^-)/\fL^+_\nabla(\cP^-) 
@>{\on{id}\times \sfq^{-,\on{spec}}}>> (\on{pt}\underset{\LS^\mer_{\cP^-}}\times \MOp^\mer_{\cG,\cP^-}) \times \fL_\nabla(\cM)/\fL^+_\nabla(\cM)  \\
@V{\jmath\times \on{id}}VV @VV{\jmath \times \on{id}}V \\
(\on{pt}\underset{\LS^\mer_{\cP^-}}\times \Op^\mer_{\cG,\cP^-})\times \fL_\nabla(\cP^-)/\fL^+_\nabla(\cP^-) 
@>{\on{id}\times \sfq}>> (\on{pt}\underset{\LS^\mer_{\cP^-}}\times \Op^\mer_{\cG,\cP^-}) \times \fL_\nabla(\cM)/\fL^+_\nabla(\cM),
\endCD
\end{equation} 
where $\sfq^{-,\on{spec}}$ denotes the map $\fL_\nabla(\cP^-)/\fL^+_\nabla(\cP^-) \to \fL_\nabla(\cM)/\fL^+_\nabla(\cM)$.  

\medskip

We can interpret \eqref{e:BC Miura * pre} as obtained by the operation 
$$-\underset{\IndCoh^*(\fL_\nabla(\cP^-))}\otimes \Vect$$
from the diagram
$$
\CD
\IndCoh^*((\on{pt}\underset{\LS^\mer_{\cP^-}}\times \MOp^\mer_{\cG,\cP^-})\times \fL_\nabla(\cP^-)/\fL^+_\nabla(\cP^-))
@>{(\on{id}\times \sfq^{-,\on{spec}})^\IndCoh_*}>> \IndCoh^*(\on{pt}\underset{\LS^\mer_{\cP^-}}\times \MOp^\mer_{\cG,\cP^-}) \times \fL_\nabla(\cM)/\fL^+_\nabla(\cM)) \\
@V{(\jmath\times \on{id})^\IndCoh_*}VV @VV{(\jmath\times \on{id})^\IndCoh_*}V \\
\IndCoh^*((\on{pt}\underset{\LS^\mer_{\cP^-}}\times \Op^\mer_{\cG,\cP^-})\times \fL_\nabla(\cP^-)/\fL^+_\nabla(\cP^-))
@>{(\on{id}\times \sfq^{-,\on{spec}})^\IndCoh_*}>> \IndCoh^*((\on{pt}\underset{\LS^\mer_{\cP^-}}\times \MOp^\mer_{\cG,\cP^-})  \times \fL_\nabla(\cM)/\fL^+_\nabla(\cM)).
\endCD
$$

Hence, it is enough to show that the natural transformation
$$(\jmath\times \on{id})^\IndCoh_* \circ (\on{id}\times \sfq^{-,\on{spec}})^!\to (\on{id}\times \sfq^{-,\on{spec}})^!\circ \jmath^\IndCoh_*$$
is an isomorphism. 

\medskip

However, the latter diagram splits as a tensor product of
$$\jmath^\IndCoh_*:\IndCoh^*((\on{pt}\underset{\LS^\mer_{\cP^-}}\times \MOp^\mer_{\cG,\cP^-}))\to \IndCoh^*((\on{pt}\underset{\LS^\mer_{\cP^-}}\times \Op^\mer_{\cG,\cP^-}))$$
along the vertical arrows and 
$$(\sfq^{-,\on{spec}})^\IndCoh_*:\IndCoh(\fL_\nabla(\cP^-)/\fL^+_\nabla(\cP^-))\to \IndCoh(\fL_\nabla(\cM)/\fL^+_\nabla(\cM))$$
along the horizontal arrows, and the assertion follows. 


\end{proof}

\sssec{}

Combining \lemref{l:BC Miura} with \secref{sss:factor Jacquet} we obtain:

\begin{cor} \label{c:factor Jacquet} \hfill

\smallskip

\noindent{\em(a)} 
The functor $J^{-,\on{spec}}_\Op$ identifies canonically with
$$\IndCoh^!(\Op^\mf_\cG) \overset{\text{\eqref{e:Op to Op semiinf !}}}\longrightarrow
\IndCoh^!((\Op^\Mmf_{\cG,\cP^-})^\wedge_\mf)\overset{J^{-,\on{spec},\semiinf_\mf}_\Op}\longrightarrow \IndCoh^!(\Op^\mf_{\cM,\rhoch_P}).$$

\smallskip

\noindent{\em(b)} 
The functor $\on{co}\!J^{-,\on{spec}}_\Op$ identifies canonically with
$$\IndCoh^*(\Op^\mf_{\cM,\rhoch_P}) \overset{\on{co}\!J^{-,\on{spec},\semiinf_\mf}_\Op}\longrightarrow 
\IndCoh^*((\Op^\Mmf_{\cG,\cP^-})^\wedge_\mf)\overset{\text{\eqref{Op semiinf * to Op}}}\longrightarrow \IndCoh^*(\Op^\mf_\cG).$$

\end{cor} 

\ssec{The semi-infinite Jacquet functor for opers, enhanced version}

In this subsection we will give slightly different interpretations of the functors $J^{-,\on{spec},\semiinf_\mf}_\Op$ and $\on{co}\!J^{-,\on{spec},\semiinf_\mf}_\Op$.

\medskip

These interpretations play a crucial role in the statements of Theorems \ref{t:local Jacquet enh} and \ref{t:local Jacquet dual enh}.  

\sssec{}

Recall the identification 
$$\IndCoh^*((\Op^\Mmf_{\cG,\cP^-})^\wedge_\mf)
\simeq \IndCoh^*(\Op_\cG^\mf)\underset{\Sph^{\on{spec}}_\cG}\otimes \on{I}(\cG,\cP^-)^{\on{spec,loc}}$$
of \eqref{e:Op semiinf * bis}.

\medskip

Thus, we can interpret the functor $\on{co}\!J^{-,\on{spec},\semiinf_\mf}_\Op$ as a functor
$$\IndCoh^*(\Op^\mf_{\cM,\rhoch_P}) \to \IndCoh^*(\Op_\cG^\mf)\underset{\Sph^{\on{spec}}_\cG}\otimes \on{I}(\cG,\cP^-)^{\on{spec,loc}}$$

When viewed as such, we will denote it by
$$\on{co}\!J^{-,\on{spec},\on{enh}}_\Op.$$

\medskip

Unwinding the constructions, we obtain that the functor $\on{co}\!J^{-,\on{spec},\on{enh}}_\Op$ 
respects the actions of $\Sph^{\on{spec}}_\cM$ on the two sides. 

\medskip

Note that \corref{c:factor Jacquet}(b) says that the original functor $\on{co}\!J^{-,\on{spec}}_\Op$ identifies with
\begin{multline*}
\IndCoh^*(\Op^\mf_{\cM,\rhoch_P}) \overset{\on{co}\!J^{-,\on{spec},\on{enh}}_\Op}\longrightarrow 
\IndCoh^*(\Op_\cG^\mf)\underset{\Sph^{\on{spec}}_\cG}\otimes \on{I}(\cG,\cP^-)^{\on{spec,loc}}
\overset{\on{Id}\otimes (-\star \Delta^{-,\on{spec},\semiinf})^R}\longrightarrow  \\
\to \IndCoh^*(\Op_\cG^\mf)\underset{\Sph^{\on{spec}}_\cG}\otimes \Sph^{\on{spec}}_\cG\simeq \IndCoh^*(\Op_\cG^\mf).
\end{multline*}

\medskip

We will also use the notation
$$\IndCoh^*(\Op_\cG^\mf)^{-,\on{enh}}:=
\IndCoh^*(\Op_\cG^\mf)\underset{\Sph^{\on{spec}}_\cG}\otimes \on{I}(\cG,\cP^-)^{\on{spec,loc}}.$$

\sssec{} \label{sss:J Op ult}

Recall now the identification
$$\IndCoh^!((\Op^\Mmf_{\cG,\cP^-})^\wedge_\mf)
\simeq \IndCoh^!(\Op_\cG^\mf)\overset{(\Sph^{\on{spec}}_\cG)^\vee}\otimes \on{I}(\cG,\cP^-)^{\on{spec,loc}}_{\on{co}}$$
of \eqref{e:Op semiinf ! bis}.

\medskip

Thus, we can view the functor $J^{-,\on{spec},\semiinf_\mf}_\Op$ as a functor
$$\IndCoh^!(\Op_\cG^\mf)\overset{(\Sph^{\on{spec}}_\cG)^\vee}\otimes \on{I}(\cG,\cP^-)^{\on{spec,loc}}_{\on{co}}\to \IndCoh^!(\Op^\mf_{\cM,\rhoch_P}).$$

When viewed as such, we will denote it by
$$J^{-,\on{spec},\on{enh}_{\on{co}}}_\Op.$$

\medskip

Unwinding the construction, we obtain that the functor $J^{-,\on{spec},\on{enh}_{\on{co}}}_\Op$ respects the coactions of
$(\Sph^{\on{spec}}_\cM)^\vee$ on the two sides. 

\medskip

In what follows we will also use the notation
$$\IndCoh^!(\Op_\cG^\mf)^{-,\on{enh}_{\on{co}}}:=
\IndCoh^!(\Op_\cG^\mf)\overset{(\Sph^{\on{spec}}_\cG)^\vee}\otimes \on{I}(\cG,\cP^-)^{\on{spec,loc}}_{\on{co}}.$$

\sssec{}

It follows from \corref{c:factor Jacquet}(a) that the original functor 
$J^{-,\on{spec}}_\Op$ identifies with
\begin{multline*} 
\IndCoh^!(\Op^\mf_\cG)  \overset{\on{Id}\otimes \one_{\on{I}(\cG,\cP^-)^{\on{spec,loc}}_{\on{co}}}}\longrightarrow 
\IndCoh^!(\Op^\mf_\cG) \otimes \on{I}(\cG,\cP^-)^{\on{spec,loc}}_{\on{co}}\to \\
\to \IndCoh^!(\Op_\cG^\mf)\overset{(\Sph^{\on{spec}}_\cG)^\vee}\otimes \on{I}(\cG,\cP^-)^{\on{spec,loc}}_{\on{co}} 
\overset{J^{-,\on{spec},\on{enh}_{\on{co}}}_\Op}\longrightarrow \IndCoh^!(\Op^\mf_{\cM,\rhoch_P}),
\end{multline*} 
where the second arrow is the left adjoint of the forgetful functor
$$\IndCoh^!(\Op_\cG^\mf)\overset{(\Sph^{\on{spec}}_\cG)^\vee}\otimes \on{I}(\cG,\cP^-)^{\on{spec,loc}}_{\on{co}} \to
\IndCoh^!(\Op^\mf_\cG) \otimes \on{I}(\cG,\cP^-)^{\on{spec,loc}}_{\on{co}},$$
which exists thanks to rigidity.

%
%
%
%

\section{Compatibility of the FLE with the Jacquet functors} \label{s:FLE and Jacquet}

In this section we first formulate the theorem that expresses the compatibility of the critical FLE with
the Jacquet functors (\thmref{t:local Jacquet}), as well as its enhanced version
(\thmref{t:local Jacquet enh}). 

\medskip

However, in order to prove both these theorems, we will reformulate them in dual terms. Thus, we will
formulate Theorems \ref{t:local Jacquet dual} and \ref{t:local Jacquet dual enh}, 
which are equivalent to Theorems \ref{t:local Jacquet} and \ref{t:local Jacquet enh}, respectively. 

\medskip

A feature of the present situation is that although \thmref{t:local Jacquet dual} looks simpler
than its enhanced version, namely, \thmref{t:local Jacquet dual enh}, we will have to prove
the latter in order to prove the former. I.e., the enhanced statement ends up being more 
accessible than the unenhanced one. The reason for this lies in the computational core of the
proof, namely, our knowledge of the Drinfeld-Sokolov reduction of Wakimoto modules, see
\secref{ss:DS of Wak}. 

\ssec{The twisted critical FLE}

In this subsection we will adapt the setting the critical FLE to the situation twisted by a $Z^0_\cG$-torsor. 

\sssec{}

Let us return to the setting of \secref{ss:twisting by dual torsors}. Let $\CP_{Z^0_\cG}$
be a $Z^0_\cG$-torsor on $X$.

\medskip

Recall the \emph{critical FLE}, i.e., the equivalence
\begin{equation} \label{e:critical FLE}
\FLE_{G,\crit_G}:\KL(G)_{\crit_G} \to \IndCoh^*(\Op_\cG^\mf)
\end{equation}
of \cite[Theorem 6.1.4]{GLC2}.

\medskip

We will consider twists of the two sides 
of \eqref{e:critical FLE} by $\CP_{Z^0_\cG}$. 

\sssec{}

On the Kac-Moody side, we consider the factorization categories
$$\hg\mod_{\crit_G+\on{dlog}(\CP_{Z^0_\cG})} \text{ and } \KL(G)_{\crit_G+\on{dlog}(\CP_{Z^0_\cG})},$$
see \secref{sss:twist by ZcG}.

\medskip

On the oper side, we consider the factorization category 
$$\IndCoh^*(\Op_{\cG,\CP_{Z^0_\cG}}^{\on{mon-free}}),$$
see \secref{sss:transl oper}. 

\medskip

We claim that the equivalence \eqref{e:critical FLE} admits a twisted version, namely,
\begin{equation} \label{e:twisted FLE}
\on{FLE}_{G,\crit_G+\on{dlog}(\CP_{Z^0_\cG})}:\KL(G)_{\crit_G+\on{dlog}(\CP_{Z^0_\cG})}\to 
\IndCoh^*(\Op_{\cG,\CP_{Z^0_\cG}}^{\on{mon-free}}).
\end{equation}

\sssec{}

Indeed, in order to construct \eqref{e:twisted FLE}, one repeats the construction from \cite[Sect. 6.1]{GLC2}, using the following version of the
Feigin-Frenkel isomorphism \cite[Theorem 5.1.2]{GLC2}:

\begin{equation} \label{e:twisted FF}
\fz_{\fg,\CP_{Z^0_\cG}}\overset{\on{FF}_G}\simeq \CO_{\Op^\reg_{\cG,\CP_{Z^0_\cG}}},
\end{equation} 
where $\fz_{\fg,\CP_{Z^0_\cG}}$ is the classical center of $\BV_{\fg,\crit_G,\CP_{Z^0_\cG}}$. 

\medskip

The identification \eqref{e:twisted FF} is obtained as the tensor product of the untwisted Feigin-Frenkel isomorphism for $\fg_{\on{der}}:=[\fg,\fg]$
with the twisted version for the torus $\fg_{\on{ab}}$, see \cite[Sect. 5.1.3]{GLC2}. 

\sssec{} \label{sss:twisted FLE compat}

The following commutative diagram is the twisted version of \cite[Theorem 6.4.5]{GLC2}: 

\medskip

\begin{equation} \label{e:twisted FLE compat}
\CD
\Whit_*(G)\underset{\Sph_G}\otimes \KL(G)_{\crit_G+\on{dlog}(\CP_{Z^0_\cG})} @>{\on{Id}\otimes \alpha_{\rho(\omega_X),\on{taut}}}>> 
\Whit_*(G)\underset{\Sph_G}\otimes \KL(G)_{\crit_G+\on{dlog}(\CP_{Z^0_\cG}),\rho(\omega_X)}  \\
@VV{\FLE^{-1}_{\cG,\infty}\otimes \FLE_{G,\crit_G+\on{dlog}(\CP_{Z^0_\cG})}}V @VVV  \\
\Rep(\cG)\underset{\Sph^{\on{spec}}_\cG}\otimes \IndCoh^*(\Op_{\cG,\CP_{Z^0_\cG}}^{\on{mon-free}}) & & \Whit_*(\hg\mod_{\crit_G+\on{dlog}(\CP_{Z^0_\cG}),\rho(\omega_X)}) \\
@VVV @VV{\ol\DS^{\on{enh,rfnd}}}V \\
\IndCoh^*(\Op^\mer_{\cG,\CP_{Z^0_\cG}}) @<{\on{FF}_G}<< \IndCoh^*(``\Spec"(\fZ_{\fg,\CP_{Z^0_\cG}})), 
\endCD
\end{equation} 

(see \cite[Sect. 4.8]{GLC2} where the notation $\ol\DS^{\on{enh,rfnd}}$ is explained). 

\sssec{} \label{sss:duality KM twisted}

Note the dualities of \cite[Equations (2.3) and (2.4)]{GLC2} induce dualities
\begin{equation} \label{e:KM self-dual twisted}
(\hg\mod_{\crit_G-\on{dlog}(\CP_{Z^0_\cG})})^\vee \simeq \hg\mod_{\crit_G+\on{dlog}(\CP_{Z^0_\cG})}
\end{equation}
and 
\begin{equation} \label{e:KL self-dual twisted}
(\KL(G)_{\crit_G,-\on{dlog}(\CP_{Z^0_\cG})})^\vee \simeq \KL(G)_{\crit_G+\on{dlog}(\CP_{Z^0_\cG})}.
\end{equation}

Note also that the Chevalley involution $\tau_G$ induces equivalences 
$$\hg\mod_{\crit_G+\on{dlog}(\CP_{Z^0_\cG})} \simeq \hg\mod_{\crit_G-\on{dlog}(\CP_{Z^0_\cG})}$$
and
$$\KL(G)_{\crit_G+\on{dlog}(\CP_{Z^0_\cG})} \simeq \KL(G)_{\crit_G-\on{dlog}(\CP_{Z^0_\cG})}.$$

\sssec{}

The equivalence of \cite[Equation (3.21)]{GLC2} induces an equivalence
$$\Theta_{\Op_\cG^\mf}: \IndCoh^!(\Op^\mf_{\cG,\CP_{Z^0_\cG}}))\to \IndCoh^*(\Op^\mf_{\cG,\CP_{Z^0_\cG}}).$$

In particular, we obtain a self-duality
\begin{equation} \label{e:Op self-dual twisted}
\left(\IndCoh^*(\Op_{\cG,\CP_{Z^0_\cG}}^{\on{mon-free}})\right)^\vee\simeq 
\IndCoh^!(\Op_{\cG,\CP_{Z^0_\cG}}^{\on{mon-free}})\overset{\Theta_{\Op^\mf_\cG}}\simeq 
\IndCoh^*(\Op_{\cG,\CP_{Z^0_\cG}}^{\on{mon-free}}).
\end{equation}

It follows formally from \cite[Theorem 8.1.4]{GLC2} that we have a commutative diagram
\begin{equation} \label{e:FLE duality twisted}
\CD
\left(\IndCoh^*(\Op_{\cG,\CP_{Z^0_\cG}}^{\on{mon-free}})\right)^\vee @>{\text{\eqref{e:Op self-dual twisted}}}>{\sim}>  
\IndCoh^*(\Op^\mf_{\cG,\CP_{Z^0_\cG}}) \\
@V{(\FLE_{G,\crit_G+\on{dlog}(\CP_{Z^0_\cG})})^\vee}VV  @AA{\FLE_{G,\crit_G+\on{dlog}(\CP_{Z^0_\cG})}\circ \tau_G}A \\
(\KL(G)_{\crit_G+\on{dlog}(\CP_{Z^0_\cG})})^\vee @>{\sim}>{\text{\eqref{e:KL self-dual twisted}}}> \KL(G)_{\crit_G-\on{dlog}(\CP_{Z^0_\cG})}. 
\endCD
\end{equation}

%
%

\sssec{}

In practice, we will take the reductive group in question to be the Levi subgroup $M$ of a standard parabolic
$P$, and $\CP_{Z^0_\cM}:=\rhoch_P(\omega_X)$.

\medskip

So in this case, the equivalence \eqref{e:twisted FLE} specializes to 
\begin{equation} \label{e:FLE M crit}
\FLE_{M,\crit_M+\rhoch_P}:\KL(M)_{\crit_M+\rhoch_P} \simeq \IndCoh^*(\Op_{\cM,\rhoch_P}^{\on{mon-free}}),
\end{equation} 
see \secref{sss:lambdach twist} for the notational conventions.

\ssec{The initial formulation}

\sssec{}

Recall the functor
$$J^{-,\Sph}_{\on{KM},\rho_P(\omega_X)}:\KL(G)_{\crit_G}\to \KL(M)_{\crit_M-\rhoch_P},$$
see \secref{sss:BRST summary}.

\medskip 

Recall also the functor 
$$J^{-,\on{spec}}_{\Op,\Theta}:\IndCoh^*(\Op_\cG^{\on{mon-free}})\to \IndCoh^*(\Op_{\cM,\rhoch_P}^{\on{mon-free}}),$$
see \secref{sss:spectral Jacquet Theta}.

\sssec{}

The following theorem, which is one of the main results of this paper, 
 expresses the compatibility of the FLE with the Jacquet functors: 

\begin{thm} \label{t:local Jacquet}
The following diagram of (lax unital) factorization functors commutes
\begin{equation} \label{e:local Jacquet}
\CD
\KL(M)_{\crit_M+\rhoch_P} @>{\on{FLE}_{M,\crit_G+\rhoch_P}}>> \IndCoh^*(\Op_{\cM,\rhoch_P}^{\on{mon-free}}) \\
@A{J^{-,\Sph}_{\on{KM},\rho_P(\omega_X)}}AA @AA{J^{-,\on{spec}}_{\Op,\Theta}}A \\
\KL(G)_{\crit_G} @>{\on{FLE}_{G,\crit_G}}>>  \IndCoh^*(\Op_\cG^{\on{mon-free}}).
\endCD
\end{equation} 
\end{thm} 

\begin{rem}
Note that the statement of \thmref{t:local Jacquet} bears a similarity with that of \corref{c:!-Jacquet on Whit}.
\end{rem}

\ssec{The enhanced version}

\sssec{}

Recall that the equivalence
$$\FLE_{G,\crit_G}:\KL(G)_{\crit_G}\to \IndCoh^*(\Op_\cG^{\on{mon-free}})$$
is compatible with the actions of
$$\Sph_G\overset{\Sat_G}\simeq \Sph_\cG^{\on{spec}},$$
(see \cite[Theorem 6.4.5(a)]{GLC2}).

\medskip

Consider the equivalence
\begin{equation} \label{e:Chev dual FLE G}
\KL(G)_{\crit_G} \overset{\FLE_{G,\crit_G}}\simeq \IndCoh^*(\Op_\cG^{\on{mon-free}}) \overset{\Theta^{-1}_{\Op^\mf_\cG}}\simeq
\IndCoh^!(\Op_\cG^{\on{mon-free}}).
\end{equation}

This equivalence is compatible with the actions of 
$$\Sph_G \overset{\Sat_G}\simeq \Sph_\cG^{\on{spec}},$$
and hence with the coactions of
\begin{equation} \label{e:dual Sat}
\Sph^\vee_G \overset{(\Sat^\vee_G)^{-1}}\simeq (\Sph_\cG^{\on{spec}})^\vee.
\end{equation} 

Consider also a similar equivalence for $M$:
\begin{equation} \label{e:Chev dual FLE M}
\KL(M)_{\crit_M+\rhoch_P} \overset{\FLE_{M,\crit_G+\rhoch_P}}\simeq \IndCoh^*(\Op_{\cM,\rhoch_P}^{\on{mon-free}}) 
\overset{\Theta^{-1}_{\Op^\mf_\cG}}\simeq
\IndCoh^!(\Op_{\cM,\rhoch_P}^{\on{mon-free}}).
\end{equation}

\sssec{}

Consider the factorization category $\on{I}(G,P^-)^{\on{loc}}_{\on{co},\rho_P(\omega_X)}$
as a comodule over the comonoidal  factorization category $\Sph^\vee_G$. 

\medskip

Consider the factorization category 
$\on{I}(\cG,\cP^-)^{\on{spec,loc}}_{\on{co}}$
as a comodule over the comonoidal  factorization category $(\Sph_\cG^{\on{spec}})^\vee$. 

\medskip

The equivalence $\Sat^{-,\semiinf}$ induces an equivalence
\begin{equation} \label{e:dual semiinf Sat}
\on{I}(G,P^-)^{\on{loc}}_{\on{co},\rho_P(\omega_X)} \overset{((\Sat^{-,\semiinf})^\vee)^{-1}}\simeq 
 \on{I}(\cG,\cP^-)^{\on{spec,loc}}_{\on{co}},
\end{equation}

This equivalence is compatible with the coactions of
\begin{equation} \label{e:Chev dual Sat}
\Sph^\vee_G \overset{(\Sat^\vee_{G,\tau})^{-1}}\simeq (\Sph_\cG^{\on{spec}})^\vee.
\end{equation} 

Hence, by \secref{sss:curse}, the equivalences \eqref{e:Chev dual FLE G} and \eqref{e:dual semiinf Sat} combine to an equivalence
\begin{equation} \label{e:tensored up Sat}
\KL(G)_{\crit_G,\rho_P(\omega_X)}^{-,\on{enh}_{\on{co}}}\simeq \IndCoh^!(\Op_\cG^\mf)^{-,\on{enh}_{\on{co}}},
\end{equation} 
where
$$\KL(G)_{\crit_G,\rho_P(\omega_X)}^{-,\on{enh}_{\on{co}}}:=\KL(G)_{\crit_G}\underset{\Sph_G}\otimes 
\on{I}(G,P^-)^{\on{loc}}_{\on{co},\rho_P(\omega_X)}\simeq 
\KL(G)_{\crit_G}\overset{\Sph^\vee_G}\otimes \on{I}(G,P^-)^{\on{loc}}_{\on{co},\rho_P(\omega_X)}$$
(see Sects. \ref{sss:Jacquet KM summary} and \ref{sss:duality over Sph})
and
$$\IndCoh^!(\Op_\cG^\mf)^{-,\on{enh}_{\on{co}}}:=\IndCoh^!(\Op_\cG^{\on{mon-free}}) \overset{(\Sph_\cG^{\on{spec}})^\vee}\otimes \on{I}(\cG,\cP^-)^{\on{spec,loc}}_{\on{co}},$$
(see \secref{sss:J Op ult}). 

\sssec{}

Consider the functor
$$J^{-,\on{enh}_{\on{co}}}_{\on{KM},\rho_P(\omega_X)}:
\KL(G)_{\crit_G}\underset{\Sph_G}\otimes \on{I}(G,P^-)^{\on{loc}}_{\on{co},\rho_P(\omega_X)}\to \KL(M)_{\crit_M+\rhoch_P}.$$
(see Sects. \ref{sss:duality over Sph} and \ref{sss:tensor vs cotensor}).

\medskip

We will think of it as a functor
$$\KL(G)_{\crit_G,\rho_P(\omega_X)}^{-,\on{enh}_{\on{co}}}\to \KL(M)_{\crit_M+\rhoch_P}.$$

%

\sssec{}

Consider the functor 
$$J^{-,\on{spec},\on{enh}_{\on{co}}}_\Op:\IndCoh^!(\Op_\cG^\mf)^{-,\on{enh}_{\on{co}}} \to \IndCoh^!(\Op^\mf_{\cM,\rhoch_P}),$$
see \secref{sss:J Op ult}.

\sssec{}

The following is an enhancement of \thmref{t:local Jacquet}: 

\begin{thm} \label{t:local Jacquet enh}
The following diagram of (lax unital) factorization functors commutes
$$
\CD
\KL(M)_{\crit_M+\rhoch_P} @>{\text{\eqref{e:Chev dual FLE M}}}>> \IndCoh^!(\Op^\mf_{\cM,\rhoch_P}) \\
@A{J^{-,\on{enh}_{\on{co}}}_{\on{KM},\rho_P(\omega_X)}}AA @AA{J^{-,\on{spec},\on{enh}_{\on{co}}}_\Op}A \\
\KL(G)_{\crit_G,\rho_P(\omega_X)}^{-,\on{enh}_{\on{co}}}
@>{\text{\eqref{e:tensored up Sat}}}>> \IndCoh^!(\Op_\cG^\mf)^{-,\on{enh}_{\on{co}}},
\endCD
$$
in a way compatible with the actions of
$$\Sph_M \overset{\Sat_M}\simeq \Sph^{\on{spec}}_\cM.$$
\end{thm}

\begin{rem}
Note the formal similarity between the statement of \thmref{t:local Jacquet dual enh} and that of \thmref{t:semiinf geom Satake}.
\end{rem}

\sssec{} \label{sss:unenh Jacquet follows from enh}

Note that the statement of \thmref{t:local Jacquet} can be obtained from that of \thmref{t:local Jacquet enh} by concatenating with the commutative 
diagram
$$
\CD
\KL(G)_{\crit_G,\rho_P(\omega_X)}^{-,\on{enh}_{\on{co}}} @>{\text{\eqref{e:tensored up Sat}}}>> \IndCoh^!(\Op_\cG^\mf)^{-,\on{enh}_{\on{co}}}  \\
@A{=}AA   @AA{=}A \\
\KL(G)_{\crit_G}\overset{\Sph^\vee_G}\otimes \on{I}(G,P^-)^{\on{loc}}_{\on{co},\rho_P(\omega_X)}  
& &  \IndCoh^!(\Op_\cG^{\on{mon-free}}) \overset{(\Sph_\cG^{\on{spec}})^\vee}\otimes \on{I}(\cG,\cP^-)^{\on{spec,loc}}_{\on{co}} \\
@AAA @AAA \\
\KL(G)_{\crit_G}\otimes \on{I}(G,P^-)^{\on{loc}}_{\on{co},\rho_P(\omega_X)}   & & 
\IndCoh^!(\Op_\cG^{\on{mon-free}}) \otimes \on{I}(\cG,\cP^-)^{\on{spec,loc}}_{\on{co}}  \\
@A{\on{Id}\otimes \one_{\on{I}(G,P^-)^{\on{loc}}_{\on{co},\rho_P(\omega_X)}}}AA 
@AA{\on{Id}\otimes \one_{\on{I}(\cG,\cP^-)^{\on{spec,loc}}_{\on{co}}}}A \\
\KL(G)_{\crit_G} @>{\text{\eqref{e:Chev dual FLE G}}}>>  \IndCoh^!(\Op_\cG^{\on{mon-free}}),
\endCD
$$
where the middle vertical arrows are the left adjoints to the natural forgetful functors. 

\medskip

Indeed, after concatenation, we obtain the diagram
$$
\CD
\KL(M)_{\crit_M-\rhoch_P} @>{\text{\eqref{e:Chev dual FLE M}}}>> \IndCoh^!(\Op^\mf_{\cM,\rhoch_P}) \\
@A{J^{-,\Sph}_{\on{KM},\rho_P(\omega_X)}}AA @AA{J^{-,\on{spec}}_{\Op}}A \\
\KL(G)_{\crit_G} @>{\text{\eqref{e:Chev dual FLE G}}}>> \IndCoh^!(\Op_\cG^{\on{mon-free}}),
\endCD
$$
which identifies with diagram \eqref{e:local Jacquet} via \eqref{e:FLE duality twisted}. 

\ssec{The dual formulation}

In order to prove Theorems \ref{t:local Jacquet} and \ref{t:local Jacquet enh}, we will reformulate them in dual 
terms.  

\sssec{}

Dualizing the arrows in \eqref{e:local Jacquet}, and taking into account \eqref{e:FLE duality twisted}, we obtain that the following 
result is equivalent to \thmref{t:local Jacquet}:

\begin{thm} \label{t:local Jacquet dual}
The following diagram of (lax unital) factorization functors commutes:
\begin{equation}  \label{e:dual Jacquet diagram}
\CD
\KL(M)_{\crit_G-\rhoch_P} @>{\FLE_{M,\crit_G+\rhoch_P}\circ \tau_M}>> \IndCoh^*(\Op_{\cM,\rhoch_P}^{\on{mon-free}}) \\
@V{\Wak^{-,\Sph}_{\rho_P(\omega_X)}}VV  @VV{\on{co}\!J_\Op^{-,\on{spec}}}V \\ 
\KL(G)_{\crit_G}   @>{\on{FLE}_{G,\crit_G}\circ \tau_G}>>   \IndCoh^*(\Op_\cG^{\on{mon-free}}). 
\endCD
\end{equation} 
\end{thm} 

We now proceed to formulating the statement dual to that of \thmref{t:local Jacquet enh}. 

\sssec{}

Combining the functors $\FLE_{G,\crit_G}\circ \tau_G$ and $\Sat^{-,\semiinf}$ we obtain an equivalence 
\begin{equation} \label{e:tensored up Sat dual}
\KL(G)_{\crit_G}\underset{\Sph_G}\otimes \on{I}(G,P^-)^{\on{loc}}_{\rho_P(\omega_X)}  \simeq
\IndCoh^*(\Op_\cG^{\on{mon-free}}) \underset{\Sph_\cG^{\on{spec}}}\otimes \on{I}(\cG,\cP^-)^{\on{spec,loc}}.
\end{equation} 

Note that the functor \eqref{e:tensored up Sat} is the dual of \eqref{e:tensored up Sat dual}.

\sssec{}

Thus, we obtain that the statement dual to \thmref{t:local Jacquet enh} reads as follows:

\begin{thm} \label{t:local Jacquet dual enh}
The following diagram of (lax unital) factorization functors commutes:
\smallskip
\begin{equation}  \label{e:dual Jacquet enh diagram}
\CD
\KL(M)_{\crit_G-\rhoch_P} @>{\FLE_{M,\crit_G+\rhoch_P}\circ \tau_M}>> \IndCoh^*(\Op_{\cM,\rhoch_P}^{\on{mon-free}}) \\
@V{\Wak^{-,\on{enh}}_{\rho_P(\omega_X)}}VV @VV{\on{co}\!J^{-,\on{spec},\on{enh}}_\Op}V \\
\KL(G)_{\crit_G}\underset{\Sph_G}\otimes \on{I}(G,P^-)^{\on{loc}}_{\rho_P(\omega_X)} @>{\text{\eqref{e:tensored up Sat dual}}}>>
\IndCoh^*(\Op_\cG^{\on{mon-free}}) \underset{\Sph_\cG^{\on{spec}}}\otimes \on{I}(\cG,\cP^-)^{\on{spec,loc}},
\endCD
\end{equation}
in a way compatible with the actions of
$$\Sph_M \overset{\Sat_{M,\tau}}\simeq \Sph^{\on{spec}}_\cM.$$
\end{thm}

\sssec{}

Note that the statement of \thmref{t:local Jacquet dual} is obtained from that of \thmref{t:local Jacquet dual enh} 
by concatenating \eqref{e:dual Jacquet enh diagram} with the diagram
$$
\CD
\KL(G)_{\crit_G}\underset{\Sph_G}\otimes \on{I}(G,P^-)^{\on{loc}}_{\rho_P(\omega_X)} @>{\text{\eqref{e:tensored up Sat dual}}}>>
\IndCoh^*(\Op_\cG^{\on{mon-free}}) \underset{\Sph_\cG^{\on{spec}}}\otimes \on{I}(\cG,\cP^-)^{\on{spec,loc}} \\
@VVV @VVV \\
\KL(G)_{\crit_G}\underset{\Sph_G}\otimes \Sph_G @>{(\on{FLE}_{G,\crit_G}\circ \tau_G)\otimes \Sat_{G,\tau}}>> 
\IndCoh^*(\Op_\cG^{\on{mon-free}}) \underset{\Sph_\cG^{\on{spec}}}\otimes  \Sph_\cG^{\on{spec}} \\
@V{\sim}VV @VV{\sim}V \\
\KL(G)_{\crit_G} @>{\FLE_{G,\crit_G}\circ \tau_G}>> \IndCoh^*(\Op_\cG^{\on{mon-free}}),
\endCD
$$
where the upper vertical arrows are 
$$((-)\underset{G}\star \Delta^{-,\semiinf})^R \text{ and } ((-)\underset{\cG}\star \Delta^{-,\on{spec},\semiinf})^R,$$
respectively. 

\begin{rem} 

When working at the pointwise level, the statement of \thmref{t:local Jacquet dual enh} is close to what
is proved in \cite{FG1}, and in fact one can generalize the results of {\it loc. cit.} to deduce the full
statement of \thmref{t:local Jacquet dual enh} at the pointwise level.  

\medskip

However, the methods of  {\it loc. cit.}
do not allow to treat families over the Ran space; the reason for this is that these methods are based 
on considering the Iwahori subgroup, and the latter is a special feature of the pointwise set-up.

\end{rem}

\section{Proof of \thmref{t:local Jacquet dual enh}} \label{s:proof local Jacquet dual enh}

\ssec{Structure of the proof} 

\sssec{} \label{sss:strategy Jacquet dual enh}

The proof of \thmref{t:local Jacquet dual enh} will be roughly modeled on that of \thmref{t:semiinf geom Satake}.
It will be structured as follows: 

\medskip

In this subsection, we will define a (lax) factorization category $\bC$, 
equipped with a unital structure and (strictly unital) factorization functors
$$\KL(G)_{\crit_G}\underset{\Sph_G}\otimes \on{I}(G,P^-)^{\on{loc}}_{\rho_P(\omega_X)} 
\overset{\sF_{\on{KM}(G)}}\longrightarrow \bC \overset{\sF^{\on{spec}}_{\Op(\cG)}}\longleftarrow
\IndCoh^*(\Op_\cG^{\on{mon-free}}) \underset{\Sph_\cG^{\on{spec}}}\otimes \on{I}(\cG,\cP^-)^{\on{spec,loc}}$$
such that:

\begin{itemize}

\item Both functors $\sF_{\on{KM}(G)}$ and $\sF^{\on{spec}}_{\Op(\cG)}$ are fully faithful when restricted to the \emph{eventually coconnective subcategories},
with respect to the t-structures to be introduced; 

\medskip

\item The functors $\Wak^{-,\on{enh}}_{\rho_P(\omega_X)}$ and $\on{co}\!J^{-,\on{spec},\on{enh}}_\Op$ send compact objects to eventually coconnective objects; 

\medskip

The above two bullet points will be carried out in \secref{s:ff in plus}. 

\medskip

\item There exists a canonical identification $\sF_{\on{KM}(G)}\simeq \sF^{\on{spec}}_{\Op(\cG)}$
as functors from
$$\KL(G)_{\crit_G}\underset{\Sph_G}\otimes \on{I}(G,P^-)^{\on{loc}}_{\rho_P(\omega_X)} \overset{\text{\eqref{e:tensored up Sat dual}}}\simeq
\IndCoh^*(\Op_\cG^{\on{mon-free}}) \underset{\Sph_\cG^{\on{spec}}}\otimes \on{I}(\cG,\cP^-)^{\on{spec,loc}}$$
to $\bC$;

\medskip

This bullet point will be carried out at the end of this subsection. 

\medskip

\item There is a canonical identification between the compositions 
$$\sF_{\on{KM}(G)}\circ \Wak^{-,\on{enh}}_{\rho_P(\omega_X)} \simeq \sF^{\on{spec}}_{\Op(\cG)}\circ \on{co}\!J^{-,\on{spec},\on{enh}}_\Op$$
as factorization functors from
$$\KL(M)_{\crit_G-\rhoch_P} \overset{\FLE_{M,\crit_G+\rhoch_P}\circ \tau_M} \simeq \IndCoh^*(\Op_{\cM,\rhoch_P}^{\on{mon-free}})$$
to $\bC$.

\medskip

The proof of this property will occupy the majority of this section.

\end{itemize}

\medskip

We depict the third property above as a commutative triangle:
\begin{equation} \label{e:category C triang}
\vcenter
{\xy
(0,0)*+{\KL(G)_{\crit_G}\underset{\Sph_G}\otimes \on{I}(G,P^-)^{\on{loc}}_{\rho_P(\omega_X)} }="A";
(75,0)*+{\IndCoh^*(\Op_\cG^{\on{mon-free}}) \underset{\Sph_\cG^{\on{spec}}}\otimes \on{I}(\cG,\cP^-)^{\on{spec,loc}}}="B";
(40,-20)*+{\bC}="C";
{\ar@{->}_{\sF_{\on{KM}(G)}} "A";"C"};
{\ar@{->}^{\sF^{\on{spec}}_{\Op(\cG)}} "B";"C"};
{\ar@{->}^{\text{\eqref{e:tensored up Sat dual}}} "A";"B"};
\endxy}
\end{equation}

\medskip

We depict the fourth property above as a commutative pentagon:
\begin{equation} \label{e:category C pent}
\vcenter
{\xy
(0,40)*+{\KL(M)_{\crit_G-\rhoch_P}}="F";
(75,40)*+{\IndCoh^*(\Op_{\cM,\rhoch_P}^{\on{mon-free}})}="G";
(0,0)*+{\KL(G)_{\crit_G}\underset{\Sph_G}\otimes \on{I}(G,P^-)^{\on{loc}}_{\rho_P(\omega_X)} }="A";
(75,0)*+{\IndCoh^*(\Op_\cG^{\on{mon-free}}) \underset{\Sph_\cG^{\on{spec}}}\otimes \on{I}(\cG,\cP^-)^{\on{spec,loc}}}="B";
(40,-20)*+{\bC}="C";
{\ar@{->}_{\sF_{\on{KM}(G)}} "A";"C"};
{\ar@{->}^{\sF^{\on{spec}}_{\Op(\cG)}} "B";"C"};
{\ar@{->}_{\Wak^{-,\on{enh}}_{\rho_P(\omega_X)}} "F";"A"};
{\ar@{->}^{\on{co}\!J^{-,\on{spec},\on{enh}}_\Op} "G";"B"};
{\ar@{->}^{\FLE_{M,\crit_G+\rhoch_P}\circ \tau_M} "F";"G"};
\endxy}
\end{equation}

\medskip

The above properties clearly imply the assertion of \thmref{t:local Jacquet dual enh}. 

\sssec{}

The compatiility with the actions of 
$$\Sph_M \overset{\Sat_{M,\tau}}\simeq \Sph^{\on{spec}}_\cM$$
will be embedded into the constructions, see \secref{ss:universal action}.

\sssec{} \label{sss:defn C}

We proceed to the definition of the category $\bC$. Recall the factorization algebra
$$\Omega(R)^{\on{spec}}\in \Rep(\cG)\otimes \Rep(\cM),$$
see \secref{sss:Omega R}. 

\medskip

Let $\Omega(R)^\Op$ denote the image of $\Omega(R)^{\on{spec}}$ under the functor
$$((\fr^\reg)^*\otimes \on{Id}):\Rep(\cG)\otimes \Rep(\cM)\to \QCoh(\Op^\reg_\cG)\otimes \Rep(\cM),$$
viewed as a factorization algebra in
$$\QCoh(\Op^\reg_\cG)\otimes \Rep(\cM),$$
where $\fr^\reg$ is the projection $\Op_\cG^\reg\to \LS^\reg_\cG$, see \cite[Sect. 3.1.5]{GLC2}.

\medskip

Using the functor
$$\QCoh(\Op^\reg_\cG)\simeq \IndCoh^*(\Op^\reg_\cG)\overset{\iota^\IndCoh_*}\longrightarrow  \IndCoh^*(\Op^\mer_\cG),$$
we can regard $\Omega(R)^\Op$ as a factorization algebra in 
$$\IndCoh^*(\Op^\mer_\cG)\otimes \Rep(\cM).$$

\sssec{}

We define the (lax) factorization category $\bC$ to be
\begin{equation} \label{e:defn of C}
\Omega(R)^\Op\mod^{\on{fact}}(\IndCoh^*(\Op^\mer_\cG)\otimes \Rep(\cM)).
\end{equation} 

We now proceed to the definition of the functors 
\begin{equation} \label{e:F functors}
\sF_{\on{KM}(G)} \text{ and } \sF^{\on{spec}}_{\Op(\cG)}.
\end{equation}

\begin{rem} \label{r:OpMmf}

Informally, one should thing of $\bC$ as a \emph{partial left completion} of the \emph{would-be} category
\begin{equation} \label{e:OpMmf}
\IndCoh^*(\Op^\Mmf_{\cG,\cP^-}).
\end{equation}

The problem is, however, that we do not know how to define \eqref{e:OpMmf} as a factorization
category (cf. \secref{sss:no OpMmf}). So we have to work directly with its surrogate provided by \eqref{e:defn of C}.

\end{rem} 

\sssec{}

Consider the full subcategory
$$\IndCoh^*((\Op^\mer_\cG)^\wedge_\mf)\subset  \IndCoh^*(\Op^\mer_\cG).$$

As in \propref{p:MOp form comp}, we have a canonical identification 
\begin{equation} \label{e:Op form comp}
\IndCoh^*(\on{pt}\underset{\LS^\reg_\cG}\times \Op^\mf_\cG)\underset{\IndCoh^*(\fL_\nabla(\cG))}\otimes \Vect
\simeq
\IndCoh^*((\Op^\mer_\cG)^\wedge_\mf). 
\end{equation} 

From here we obtain an identification
\begin{equation} \label{e:Op form comp bis}
\IndCoh^*(\Op^\mf_\cG)\underset{\Sph_\cG^{\on{spec}}}\otimes \Rep(\cG)\simeq
\IndCoh^*((\Op^\mer_\cG)^\wedge_\mf). 
\end{equation} 

\sssec{} \label{sss:ident Omega R two}

Note that the lax unital factorization functor 
\begin{multline} \label{e:Rep to mer}
\Rep(\cG) \overset{\one_{\IndCoh^*(\Op^\mf_\cG)}\otimes \on{Id}}\longrightarrow \\
\to \IndCoh^*(\Op^\mf_\cG)\underset{\Sph_\cG^{\on{spec}}}\otimes \Rep(\cG)\simeq 
\IndCoh^*((\Op^\mer_\cG)^\wedge_\mf)\to \IndCoh^*(\Op^\mer_\cG)
\end{multline} 
identifies with 
$$\Rep(\cG) \overset{(\fr^\reg)^*}\to \IndCoh^*(\Op^\reg_\cG) \overset{\iota^\IndCoh_*}\longrightarrow \IndCoh^*(\Op^\mer_\cG).$$

In particular, the image of $\Omega(R)^{\on{spec}}$ along the functor 
$$\Rep(\cG) \otimes \Rep(\cM)\overset{\text{\eqref{e:Rep to mer}}\otimes \on{Id}}\longrightarrow \IndCoh^*(\Op^\mer_\cG)\otimes \Rep(\cM)$$
identifies with $\Omega(R)^\Op$.

\sssec{}

Consider the full subcategory
\begin{multline} \label{e:Omega mod form comp}
\Omega(R)^\Op\mod^{\on{fact}}\left(\IndCoh^*((\Op^\mer_\cG)^\wedge_\mf)\otimes \Rep(\cM)\right)\subset \\
\subset \Omega(R)^\Op\mod^{\on{fact}}(\IndCoh^*(\Op^\mer_\cG)\otimes \Rep(\cM)).
\end{multline}

\medskip

From the fact that $\Sph^{\on{spec}}_\cG$ is rigid and $\IndCoh^*(\Op^\mf_\cG)$ is dualizable, we obtain that the functor
\begin{multline}  \label{e:Omera R mod form comp prel}
\IndCoh^*(\Op^\mf_\cG) \underset{\Sph_\cG^{\on{spec}}}\otimes
\left(\Omega(R)^{\on{spec}}\mod^{\on{fact}}(\Rep(\cG)\otimes \Rep(\cM))\right)\to \\
\to \Omega(R)^{\on{spec}}\mod^{\on{fact}}\left((\IndCoh^*(\Op^\mf_\cG)\underset{\Sph_\cG^{\on{spec}}}\otimes \Rep(\cG))\otimes \Rep(\cM)\right)
\end{multline} 
is an equivalence, cf. \cite[Lemma B.12.9]{GLC2}.

\sssec{}

Combining \eqref{e:Omera R mod form comp prel} with \eqref{e:Op form comp bis} and \secref{sss:ident Omega R two}, we obtain an identification
\begin{multline}  \label{e:Omera R mod form comp}
\IndCoh^*(\Op^\mf_\cG) \underset{\Sph_\cG^{\on{spec}}}\otimes
\left(\Omega(R)^{\on{spec}}\mod^{\on{fact}}(\Rep(\cG)\otimes \Rep(\cM))\right) \simeq \\
\simeq \Omega(R)^\Op\mod^{\on{fact}}\left(\IndCoh^*((\Op^\mer_\cG)^\wedge_\mf)\otimes \Rep(\cM)\right).
\end{multline}

\begin{rem} \label{r:OpMmf wedge}

Recall that although we could not define \eqref{e:OpMmf} as a factorization category, we could do so for 
$\IndCoh^*((\Op^\Mmf_{\cG,\cP^-})^\wedge_\mf)$.

\medskip

In the spirit of Remark \ref{r:OpMmf}, the category 
$$\Omega(R)^\Op\mod^{\on{fact}}\left(\IndCoh^*((\Op^\mer_\cG)^\wedge_\mf)\otimes \Rep(\cM)\right)$$
can indeed be thought as a partial left completion of $\IndCoh^*((\Op^\Mmf_{\cG,\cP^-})^\wedge_\mf)$:

\end{rem} 

\sssec{}

Recall now from \secref{sss:IGP via Omega} that we have a commutative diagram of (strictly unital) $\Sph_\cG^{\on{spec}}$-linear factorization functors
\begin{equation} \label{e:Omega triang}
\vcenter
{\xy
(0,0)*+{\on{I}(G,P^-)^{\on{loc}}_{\rho_P(\omega_X)}}="A";
(75,0)*+{\on{I}(\cG,\cP^-)^{\on{spec,loc}}}="B";
(40,-20)*+{\Omega(R)^{\on{spec}}\mod^{\on{fact}}(\Rep(\cG)\otimes \Rep(\cM)).}="C";
{\ar@{->}_{\sF} "A";"C"};
{\ar@{->}^{\sF^{\on{spec}}} "B";"C"};
{\ar@{->}^{\Sat^{-,\semiinf}} "A";"B"};
\endxy}
\end{equation}


\sssec{} \label{sss:triang diag}

Tensoring the functors $\sF$ and $\sF^{\on{spec}}$ with
$$\KL(G)_{\crit_G} \overset{\FLE_{G,\crit_G}\circ \tau_G}\longrightarrow \IndCoh^*(\Op^\mf_\cG) \overset{\on{Id}}\longleftarrow \IndCoh^*(\Op^\mf_\cG),$$
respectively, we obtain functors
$$\KL(G)_{\crit_G}\underset{\Sph_G}\otimes \on{I}(G,P^-)^{\on{loc}}_{\rho_P(\omega_X)} \to
\IndCoh^*(\Op^\mf_\cG) \underset{\Sph_\cG^{\on{spec}}}\otimes
\left(\Omega(R)^{\on{spec}}\mod^{\on{fact}}(\Rep(\cG)\otimes \Rep(\cM))\right)$$
and
\begin{multline*} 
\IndCoh^*(\Op_\cG^{\on{mon-free}}) \underset{\Sph_\cG^{\on{spec}}}\otimes \on{I}(\cG,\cP^-)^{\on{spec,loc}}\to \\
\to \IndCoh^*(\Op^\mf_\cG) \underset{\Sph_\cG^{\on{spec}}}\otimes
\left(\Omega(R)^{\on{spec}}\mod^{\on{fact}}(\Rep(\cG)\otimes \Rep(\cM))\right).
\end{multline*} 

Composing with \eqref{e:Omera R mod form comp} and with \eqref{e:Omega mod form comp} 
we obtain the sought-for functors
\eqref{e:F functors}.  These functors make the diagram \eqref{e:category C triang} commute by
construction. 

\begin{rem}

Recall (see \eqref{e:Op semiinf * bis}) that we can think of 
$$\IndCoh^*(\Op_\cG^{\on{mon-free}}) \underset{\Sph_\cG^{\on{spec}}}\otimes \on{I}(\cG,\cP^-)^{\on{spec,loc}}$$
as $\IndCoh^*((\Op^\Mmf_{\cG,\cP^-})^\wedge_\mf)$.

\medskip

In terms of this identification and Remark \ref{r:OpMmf}, the functor $\sF^{\on{spec}}_{\Op(\cG)}$ should 
informally\footnote{Informally: because the target category in the next formula as not defined.} be thought of as the composition of
the inclusion
$$\IndCoh^*((\Op^\Mmf_{\cG,\cP^-})^\wedge_\mf)\hookrightarrow \IndCoh^*((\Op^\Mmf_{\cG,\cP^-}))$$
and the tautological functor
$$\IndCoh^*((\Op^\Mmf_{\cG,\cP^-}))\to \bC.$$

\medskip

More rigorously, in terms of Remark \ref{r:OpMmf wedge}, the functor $\sF^{\on{spec}}_{\Op(\cG)}$ maps 
$\IndCoh^*((\Op^\Mmf_{\cG,\cP^-})^\wedge_\mf)$ to its partial left completion.\footnote{By this we only mean that the functor
in question induces an equivalence on the eventually coconnective subcategories.}  

\end{rem}

%

\sssec{}

The rest of this section is devoted to the verification of the
commutativity of the pentagon \eqref{e:category C pent}. 

\medskip

As was mentioned earlier, the first two bullet points in \secref{sss:strategy Jacquet dual enh} will be carried out
in \secref{s:ff in plus}. 

\ssec{The composition on the spectral side} \label{ss:comp spec}

Denote 
$$\sF^{\on{spec}}_{\MOp(\cM)}:=\sF^{\on{spec}}_{\Op(\cG)}\circ \on{co}\!J^{-,\on{spec},\on{enh}}_\Op, \quad
\IndCoh^*(\Op_{\cM,\rhoch_P}^{\on{mon-free}})\to \bC.$$

As a first step towards establishing the commutativity of the pentagon \eqref{e:category C pent}, 
in this subsection, we will break down this functor into a composition of more elementary operations. 

\begin{rem} \label{r:OpMmf functor}

In terms of Remarks \ref{r:OpMmf} and \ref{r:OpMmf wedge}, the functor $\sF^{\on{spec}}_{\MOp(\cM)}$
is the composition
\begin{multline*}
\IndCoh^*(\Op_{\cM,\rhoch_P}^{\on{mon-free}})\overset{((\sfq^{\on{Miu,mon-free}})^\wedge)^!}\longrightarrow 
\IndCoh^*((\MOp^\Mmf_{\cG,\cP^-})^\wedge_\mf) \overset{\jmath^\IndCoh_*}\longrightarrow \\
\to \IndCoh^*((\Op^\Mmf_{\cG,\cP^-})^\wedge_\mf) \to \IndCoh^*(\Op^\Mmf_{\cG,\cP^-}) \to \bC.
\end{multline*} 

The idea of this section is to (informally) rewrite this functor as the composition of the 
endofunctor
\begin{multline} \label{e:restrict to formal and back}
\IndCoh^*(\Op_{\cM,\rhoch_P}^{\on{mon-free}})\overset{((\sfq^{\on{Miu,mon-free}})^\wedge)^!}\longrightarrow \\
\to \IndCoh^*((\MOp^\Mmf_{\cG,\cP^-})^\wedge_\mf) \overset{((\sfq^{\on{Miu,mon-free}})^\wedge)^\IndCoh_*}\longrightarrow
\IndCoh^*(\Op_{\cM,\rhoch_P}^{\on{mon-free}})
\end{multline} 
of $\IndCoh^*(\Op_{\cM,\rhoch_P}^{\on{mon-free}})$, followed by 
\begin{equation} \label{e:OpMmf Wak}
\IndCoh^*(\Op_{\cM,\rhoch_P}^{\on{mon-free}}) \overset{\jmath^\IndCoh_*}\longrightarrow 
\IndCoh^*(\Op^\Mmf_{\cG,\cP^-}) \to \bC.
\end{equation} 

The point is that although we cannot define the individual functors in \eqref{e:OpMmf Wak} (because we do not know how
to make sense of $\IndCoh^*(\Op^\Mmf_{\cG,\cP^-})$), we \emph{can} define their composition. The resulting functor
\begin{equation} \label{e:OpMmf Wak 1}
\IndCoh^*(\Op_{\cM,\rhoch_P}^{\on{mon-free}}) \to \bC
\end{equation} 
is the functor $\sF^{\on{spec,ext}}_{\MOp(\cM)}$ defined at the end of this section. 

\medskip

The reason for doing so is that one can express the functor \eqref{e:OpMmf Wak 1} in terms of factorization
algebras, and thus compare it more directly with its counterpart on the geometric side, thereby establishing
the commutativity of \eqref{e:category C pent}. 

\end{rem}

\sssec{}

We will first describe the functor $\on{pre-}\!\sF^{\on{spec}}_{\MOp(\cM)}$, which is the composition of $\sF^{\on{spec}}_{\MOp(\cM)}$ 
with the forgetful functor
\begin{equation} \label{e:oblv Omega R}
\oblv_{\Omega(R)^\Op}:
\Omega(R)^\Op\mod^{\on{fact}}\left(\IndCoh^*(\Op^\mer_\cG)\otimes \Rep(\cM)\right)
\to \IndCoh^*(\Op^\mer_\cG)\otimes \Rep(\cM).
\end{equation} 

Prior to that, we will describe the functor $\on{pre-pre-}\!\sF^{\on{spec}}_{\MOp(\cM)}$, equal to the composition of 
$\on{pre-}\!\sF^{\on{spec}}_{\MOp(\cM)}$ with 
\begin{equation} \label{e:inv M}
\IndCoh^*(\Op^\mer_\cG)\otimes \Rep(\cM) \overset{\on{Id}\otimes \on{inv}_\cM}\longrightarrow \IndCoh^*(\Op^\mer_\cG),
\end{equation} 
cf. \secref{sss:pre F} and Remark \ref{r:pre F}. 

\sssec{}

Namely, unwinding the definitions, we obtain that the functor $\on{pre-pre-}\!\sF^{\on{spec}}_{\MOp(\cM)}$ identifies with the composition of the
endofunctor \eqref{e:restrict to formal and back} of $\IndCoh^*(\Op_{\cM,\rhoch_P}^{\on{mon-free}})$ 
with the functor
$$\IndCoh^*(\Op_{\cM,\rhoch_P}^{\on{mon-free}})\to \IndCoh^*(\Op^\mer_\cG)$$
equal to
\begin{equation} \label{e:pre-pre-F ext}
\IndCoh^*(\Op_{\cM,\rhoch_P}^{\on{mon-free}}) \overset{\sfq^{\on{Miu}}}\simeq 
\IndCoh^*(\MOp^\Mmf_{\cG,\cP^-})\overset{(\sfp^{\on{Miu}})^\IndCoh_*}\longrightarrow 
\IndCoh^*(\Op^\mer_\cG).
\end{equation} 

\sssec{} \label{sss:pre from pre-pre}

We recover the functor $\on{pre-}\!\sF^{\on{spec}}_{\MOp(\cM)}$ from $\on{pre-pre-}\!\sF^{\on{spec}}_{\MOp(\cM)}$ using the 
self-duality of $\Rep(\cM)$.

\medskip

Namely, this is the unique $\Rep(\cM)$-linear functor, whose composition with \eqref{e:inv M} is $\on{pre-pre-}\!\sF^{\on{spec}}_{\MOp(\cM)}$.

\medskip

Explicitly, $\on{pre-}\!\sF^{\on{spec}}_{\MOp(\cM)}$ is the composition
\begin{multline*}
\IndCoh^*(\Op_{\cM,\rhoch_P}^{\on{mon-free}}) \overset{\on{Id}\otimes \on{u}_{\Rep(\cM)}}\longrightarrow
\IndCoh^*(\Op_{\cM,\rhoch_P}^{\on{mon-free}}) \otimes \Rep(\cM) \otimes \Rep(\cM) \overset{\on{act}\otimes \on{Id}}\longrightarrow \\
\to \IndCoh^*(\Op_{\cM,\rhoch_P}^{\on{mon-free}}) \otimes \Rep(\cM) \overset{\on{pre-pre-}\!\sF^{\on{spec}}_{\MOp(\cM)}\otimes \on{Id}}\longrightarrow
\IndCoh^*(\Op^\mer_\cG)\otimes \Rep(\cM),
\end{multline*}
where $\on{u}_{\Rep(\cM)}\in \Rep(\cM) \otimes \Rep(\cM)$ is the unit of the self-duality. 

\sssec{} \label{sss:recover F from A phi}

Denote by $$\CA^{\on{spec}}\in \IndCoh^*(\Op^\mer_\cG)\otimes \Rep(\cM)$$ the image of the factorization unit, i.e., 
$$\CO_{\Op^\reg_{\cM,\rhoch_P}}\in \IndCoh^*(\Op_{\cM,\rhoch_P}^{\on{mon-free}})$$
under the functor $\on{pre-}\!\sF^{\on{spec}}_{\MOp(\cM)}$.

\medskip

The lax unital structure on the functor $\sF^{\on{spec}}_{\MOp(\cM)}$ gives rise to a map 
$$\phi^{\on{spec}}:\Omega(R)^\Op\to \CA^{\on{spec}}$$
as unital factorization algebras in $\IndCoh^*(\Op^\mer_\cG)\otimes \Rep(\cM)$.

\medskip

Vice versa, given $\on{pre-}\!\sF^{\on{spec}}_{\MOp(\cM)}$, the datum of $\phi^{\on{spec}}$ recovers 
$\sF^{\on{spec}}_{\MOp(\cM)}$ as a lax unital factorization functor. 

\medskip

Namely, $\sF^{\on{spec}}_{\MOp(\cM)}$ is the composition 
of the functor
$$\IndCoh^*(\Op_{\cM,\rhoch_P}^{\on{mon-free}})\to \CA^{\on{spec}}\mod^{\on{fact}}(\IndCoh^*(\Op^\mer_\cG)\otimes \Rep(\cM))$$
induced by $\on{pre-}\!\sF^{\on{spec}}_{\MOp(\cM)}$, followed by the functor
$$\CA^{\on{spec}}\mod^{\on{fact}}(\IndCoh^*(\Op^\mer_\cG)\otimes \Rep(\cM))\to
\Omega(R)^\Op\mod^{\on{fact}}(\IndCoh^*(\Op^\mer_\cG)\otimes \Rep(\cM)),$$
given by restriction along $\phi^{\on{spec}}$. 

\sssec{}

Let us denote by $\on{pre-pre-}\!\sF^{\on{spec,ext}}_{\MOp(\cM)}$ the functor \eqref{e:pre-pre-F ext}. Note that by the same mechanism
as above, the functor $\on{pre-pre-}\!\sF^{\on{spec,ext}}_{\MOp(\cM)}$ gives rise to a (lax unital) factorization functor
$$\on{pre-}\!\sF^{\on{spec,ext}}_{\MOp(\cM)}:\IndCoh^*(\Op_{\cM,\rhoch_P}^{\on{mon-free}})\to 
\IndCoh^*(\Op^\mer_\cG)\otimes \Rep(\cM),$$
so that $\on{pre-}\!\sF^{\on{spec}}_{\MOp(\cM)}$ is the precomposition of $\on{pre-}\!\sF^{\on{spec,ext}}_{\MOp(\cM)}$ with \eqref{e:restrict to formal and back}. 

\sssec{} \label{sss:recover F from A phi ext}

Furthermore, since the functor \eqref{e:restrict to formal and back} is \emph{strictly unital}, the image of the factorization unit under
$\on{pre-}\!\sF^{\on{spec,ext}}_{\MOp(\cM)}$ is the same factorization algebra $\CA^{\on{spec}}$. 

\medskip

Hence, the homomorphism $\phi^{\on{spec}}$ gives rise to a lift of $\on{pre-}\!\sF^{\on{spec,ext}}_{\MOp(\cM)}$ to a lax unital
factorization functor 
$$\sF^{\on{spec,ext}}_{\MOp(\cM)}:\IndCoh^*(\Op_{\cM,\rhoch_P}^{\on{mon-free}})\to
\Omega(R)^\Op\mod^{\on{fact}}\left(\IndCoh^*(\Op^\mer_\cG)\otimes \Rep(\cM)\right)=:\bC,$$
so that $\sF^{\on{spec}}_{\MOp(\cM)}$ is the precomposition of $\sF^{\on{spec,ext}}_{\MOp(\cM)}$ with \eqref{e:restrict to formal and back}. 

\begin{rem}

See Remark \ref{r:OpMmf functor} for the geometric meaning of the functor $\sF^{\on{spec,ext}}_{\MOp(\cM)}$.

\end{rem}

\ssec{The composition on the geometric side} \label{ss:comp geom}

Denote 
$$\sF_{\on{KL}(M)}:=\sF_{\on{KM}(G)}\circ \Wak^{-,\on{enh}}_{\rho_P(\omega_X)}.$$

In this subsection we will rewrite this functor, along the lines of how the functor $\sF^{\on{spec}}_{\MOp(\cM)}$
was described in \secref{ss:comp spec}. 

\medskip

Namely, we will express it as a composition of an endofunctor of
$\KL(M)_{\crit_G-\rhoch_P}$ parallel to \eqref{e:restrict to formal and back}, followed by a certain functor
$$\sF^{\on{ext}}_{\on{KL}(M)}:\KL(M)_{\crit_G-\rhoch_P}\to \Omega(R)^\Op\mod^{\on{fact}}\left(\IndCoh^*(\Op^\mer_\cG)\otimes \Rep(\cM)\right)=:\bC,$$
which we will be able to access explicitly in terms of factorization algebras.

\sssec{}

First, by the same mechanism as in \secref{ss:comp spec}, we can consider the functors
$$\on{pre-pre-}\!\sF_{\on{KL}(M)}:\KL(M)_{\crit_G-\rhoch_P}\to \IndCoh^*(\Op^\mer_\cG),$$
$$\on{pre-}\!\sF_{\on{KL}(M)}:\KL(M)_{\crit_G-\rhoch_P}\to \IndCoh^*(\Op^\mer_\cG)\otimes \Rep(\cM)$$
and the factorization algebra
$$\CA\in \IndCoh^*(\Op^\mer_\cG)\otimes \Rep(\cM),$$
equipped with a homomorphism 
$$\phi:\Omega(R)^\Op\to \CA.$$

The datum of $\on{pre-pre-}\!\sF_{\on{KL}(M)}$ recovers $\on{pre-}\!\sF_{\on{KL}(M)}$
by the same mechanism as in \secref{sss:pre from pre-pre}.

\medskip

In its turn, we recover $\sF_{\on{KL}(M)}$ from $\on{pre-}\!\sF_{\on{KL}(M)}$ and the homomorphism $\phi$,
see \secref{sss:recover F from A phi}.

\sssec{}

Our current goal is to construct the counterparts on the geometric side of the functors $\on{pre-pre-}\!\sF^{\on{spec,ext}}_{\MOp(\cM)}$ 
(resp., $\on{pre-}\!\sF^{\on{spec,ext}}_{\MOp(\cM)}$, $\!\sF^{\on{spec,ext}}_{\MOp(\cM)}$), to be denoted 
$$\on{pre-pre-}\!\sF^{\on{ext}}_{\on{KL}(M)},\,\, \on{pre-}\!\sF^{\on{ext}}_{\on{KL}(M)} \text{ and }
\sF^{\on{ext}}_{\on{KL}(M)},$$
respectively. 

\sssec{}

Recall the (lax) factorization category 
$$\hg\mod_{\crit_G}^{-,\semiinf}:=(\hg\mod_{\crit_G})^{\fL(N^-_P)\cdot \fL^+(M)},$$
see \secref{sss:semiinf KM}. 

\medskip

Define the functor
$$\on{pre-pre-}\!\sF_{\on{KM}(G)}^{\on{ext}}:\hg\mod_{\crit_G}^{-,\semiinf}\to \IndCoh^*(\Op^\mer_\cG)$$
to be the composition 
\begin{multline*}
(\hg\mod_{\crit_G})^{-,\semiinf}  \to \hg\mod_{\crit_G}^{\fL^+(M)} \overset{\alpha_{\rho(\omega_X),\on{taut}}}\longrightarrow 
\hg\mod_{\crit_G,\rho(\omega_X)}^{\fL^+(M)_{\rho(\omega_X)}}  
\to \hg\mod_{\crit_G,\rho(\omega_X)} \overset{\DS^{\on{enh,rfnd}}}\longrightarrow \\
\to \IndCoh^*(\Op^\mer_\cG)
\overset{\tau_G}\to  \IndCoh^*(\Op^\mer_\cG). 
\end{multline*}
where:

\begin{itemize}

\item The first and the third arrows are the forgetful functors;

\medskip

\item The functor $\DS^{\on{enh,rfnd}}$ is as in \cite[Sect. 5.2.8]{GLC2}; 

\end{itemize} 

\medskip

From $\on{pre-pre-}\!\sF_{\on{KM}(G)}^{\on{ext}}$ we produce a $\Rep(\cM)$-linear functor 
$$\on{pre-}\!\sF_{\on{KM}(G)}^{\on{ext}}:(\hg\mod_{\crit_G})^{-,\semiinf}  \to \IndCoh^*(\Op^\mer_\cG)\otimes \Rep(\cM)$$
by the same principle as in \secref{sss:pre from pre-pre}.

\sssec{}

Let $\on{pre-}\!\sF_{\on{KM}(G)}$ and $\on{pre-pre-}\!\sF_{\on{KM}(G)}$ denote the compositions of $\sF_{\on{KM}(G)}$ with
the functors
$$\oblv_{\Omega(R)^\Op} \text{ and } (\on{Id}\otimes \on{inv}_\cM)\circ \oblv_{\Omega(R)^\Op},$$
respectively.

\medskip

Note that the functor $\on{pre-}\!\sF_{\on{KM}(G)}$ is recovered from $\on{pre-pre-}\!\sF_{\on{KM}(G)}$ by the same principle as
in \secref{sss:pre from pre-pre}.

\sssec{}

We now claim:

\begin{prop} \label{p:F KM Sph via semiinf}
The functor $\on{pre-pre-}\!\sF_{\on{KM}(G)}$ identifies canonically with the composition
\begin{multline} \label{e:F KM Sph via semiinf}
\KL(G)_{\crit_G}\underset{\Sph_G}\otimes \on{I}(G,P^-)^{\on{loc}}_{\rho_P(\omega_X)} 
\overset{\on{Id}\otimes \alpha^{-1}_{\rho_P(\omega_X),\on{taut}}}\longrightarrow \KL(G)_{\crit_G}\underset{\Sph_G}\otimes \on{I}(G,P^-)^{\on{loc}}
\overset{\Accs}\longrightarrow \\
\to (\hg\mod_{\crit_G})^{-,\semiinf} \overset{\on{pre-pre-}\!\sF_{\on{KM}(G)}^{\on{ext}}}\to \IndCoh^*(\Op^\mer_\cG),
\end{multline} 
where $\Accs$ is the functor from \eqref{e:Sph gen KM}. 
\end{prop}

\begin{proof}

Let 
$$\on{pre-pre-}\!\sF:\on{I}(G,P^-)^{\on{loc}}_{\rho_P(\omega_X)}\to \Rep(\cG)$$
be as in \secref{sss:prepre F}. We need to identify the composition \eqref{e:F KM Sph via semiinf} with the functor obtained by tensoring
$\on{pre-pre-}\!\sF$ with
$$\KL(G)_{\crit_G} \overset{\tau_G\circ \FLE_{G,\crit}}\longrightarrow \IndCoh^*(\Op^\mf_\cG)$$
over
$$\Sph_G\overset{\Sat_{G,\tau}}\simeq \Sph^{\on{spec}}_\cG.$$

\medskip

By construction, the functor \eqref{e:F KM Sph via semiinf} can be rewritten as 
\begin{multline*} 
\KL(G)_{\crit_G}\underset{\Sph_G}\otimes \on{I}(G,P^-)^{\on{loc}}_{\rho_P(\omega_X)} \to 
\KL(G)_{\crit_G}\underset{\Sph_G}\otimes \Whit_*(G) \overset{\alpha_{\rho(\omega_X),\on{taut}}\otimes \on{Id}}\longrightarrow \\
\to \KL(G)_{\crit_G,\rho(\omega_X)}\underset{\Sph_G}\otimes \Whit_*(G) 
\to \Whit_*(\hg\mod_{\crit_G,\rho(\omega_X)})
\overset{\ol\DS^{\on{enh,rfnd}}}\longrightarrow  \\
\to \IndCoh^*(\Op^\mer_\cG) \overset{\tau_G}\to  \IndCoh^*(\Op^\mer_\cG),
\end{multline*} 
where the first arrow is the tensor product of the identity functor on $\KL(G)_{\crit_G}$ with
\begin{equation} \label{e:F KM Sph via semiinf 1}
\on{I}(G,P^-)^{\on{loc}}_{\rho_P(\omega_X)} \overset{\alpha_{\rho_M(\omega_X),\on{taut}}}\longrightarrow  \on{I}(G,P^-)^{\on{loc}}_{\rho(\omega_X)} \to 
\Dmod_{\frac{1}{2}}(\Gr_{G,\rho(\omega_X)})\to \Whit_*(G).
\end{equation} 

Note, however, that the composition of \eqref{e:F KM Sph via semiinf 1} with (the inverse of)
$$\Rep(\cG)\overset{\FLE_{\cG,\infty}}\simeq \Whit_*(G)$$
is exactly the functor $\on{pre-pre-}\!\sF$.

\medskip

The required identification follows now from \cite[Theorem 6.4.5]{GLC2}.

\end{proof} 

\sssec{}

Recall (see \cite[Sect. 5.3]{GLC2}) that the category $\KL(M)_{\crit_M+\rhoch_P}$ is acted on by the monoidal category $\QCoh(\Op_{\cM,\rhoch_P}^\mf)$.
The Chevalley involution on $M$ allows us to identify 
$$\KL(M)_{\crit_M+\rhoch_P}\simeq \KL(M)_{\crit_M-\rhoch_P}.$$

Hence, we can view $\KL(M)_{\crit_M-\rhoch_P}$ also as acted on by $\QCoh(\Op_{\cM,\rhoch_P}^\mf)$.

\medskip

Let
$$(\KL(M)_{\crit_M-\rhoch_P})_{\Gmf}\subset \KL(M)_{\crit_M-\rhoch_P}$$
be the full subcategory set-theoretically supported over the image of the map $\sfq^{\on{Miu,mon-free}}$ of
\eqref{e:G-mf MOp}.

\medskip

Let
$$((\sfq^{\on{Miu,mon-free}})^\wedge)^\IndCoh_*:(\KL(M)_{\crit_M-\rhoch_P})_{\Gmf}\rightleftarrows \KL(M)_{\crit_M-\rhoch_P}:((\sfq^{\on{Miu,mon-free}})^\wedge)^!$$
denote the resulting pair of adjoint functors. 

\sssec{}

Recall now that according to \lemref{l:KM to Sph gen}, the functor
\begin{multline} \label{e:Sph gen KM rho}
\KL(G)_{\crit_G}\underset{\Sph_G}\otimes \on{I}(G,P^-)^{\on{loc}}_{\rho_P(\omega_X)}
\overset{\on{Id}\otimes \alpha^{-1}_{\rho_P(\omega_X),\on{taut}}}\simeq \\
\simeq \KL(G)_{\crit_G}\underset{\Sph_G}\otimes \on{I}(G,P^-)^{\on{loc}}
\overset{\Accs}\longrightarrow  (\hg\mod_{\crit_G})^{-,\semiinf}
\end{multline}
is fully faithful. 

\medskip

\begin{prop} \label{p:support of Wak Sph} \hfill

\smallskip

\noindent{\em(a)} The functor 
$$\Wak^{-,\semiinf}_{\rho_P(\omega_X)}:\KL(M)_{\crit_M-\rhoch_P}\to (\hg\mod_{\crit_G})^{-,\semiinf}$$
sends $(\KL(M)_{\crit_M-\rhoch_P})_{\Gmf}\subset \KL(M)_{\crit_M-\rhoch_P}$ to essential image of 
\eqref{e:Sph gen KM rho}. I.e., we have a commutative diagram
\begin{equation} \label{e:support of Wak Sph 1}
\CD
(\KL(M)_{\crit_M-\rhoch_P})_{\Gmf} @>{((\sfq^{\on{Miu,mon-free}})^\wedge)^\IndCoh_*}>> \KL(M)_{\crit_M-\rhoch_P} \\
@V{\Wak^{-,\semiinf}_{\rho_P(\omega_X)}}VV  @VV{\Wak^{-,\semiinf}_{\rho_P(\omega_X)}}V \\
\KL(G)_{\crit_G}\underset{\Sph_G}\otimes \on{I}(G,P^-)^{\on{loc}}_{\rho_P(\omega_X)} @>>{\text{\eqref{e:Sph gen KM rho}}}> (\hg\mod_{\crit_G})^{-,\semiinf}.
\endCD
\end{equation}

\smallskip

\noindent{\em(b)} 
The natural transformation in the diagram
\begin{equation} \label{e:support of Wak Sph 2}
\vcenter
{\xy
(0,0)*+{(\KL(M)_{\crit_M-\rhoch_P})_{\Gmf}}="A";
(70,0)*+{\KL(M)_{\crit_M-\rhoch_P}}="B";
(0,-20)*+{\KL(G)_{\crit_G}\underset{\Sph_G}\otimes \on{I}(G,P^-)^{\on{loc}}_{\rho_P(\omega_X)}}="C";
(70,-20)*+{(\hg\mod_{\crit_G})^{-,\semiinf}}="D";
{\ar@{->}_{((\sfq^{\on{Miu,mon-free}})^\wedge)^!} "B";"A"};
{\ar@{->} "D";"C"};
{\ar@{->}_{\Wak^{-,\semiinf}_{\rho_P(\omega_X)}} "A";"C"};
{\ar@{->}^{\Wak^{-,\semiinf}_{\rho_P(\omega_X)}} "B";"D"};
{\ar@{=>} "A";"D"};
\endxy}
\end{equation}
obtained from \eqref{e:support of Wak Sph 1} by passing to right adjoints along the horizontal arrows,
is an isomorphism.

\end{prop} 

The proof of the proposition will be given in \secref{ss:support of Wak Sph}. 

\sssec{}

Combining \propref{p:support of Wak Sph} with \corref{c:KM to Sph gen}, we obtain:

\begin{cor} \label{c:support of Wak Sph}
We have a canonical isomorphism
$$\Wak^{-,\semiinf}_{\rho_P(\omega_X)}|_{(\KL(M)_{\crit_M-\rhoch_P})_{\Gmf}} \circ ((\sfq^{\on{Miu,mon-free}})^\wedge)^!
\simeq \Wak^{-,\on{enh}}_{\rho_P(\omega_X)}$$
as functors
$$\KL(M)_{\crit_M-\rhoch_P}\to \KL(G)_{\crit_G}\underset{\Sph_G}\otimes \on{I}(G,P^-)^{\on{loc}}_{\rho_P(\omega_X)}.$$
\end{cor} 

\sssec{}

Set
$$\on{pre-pre-}\!\sF^{\on{ext}}_{\KL(M)}:= \on{pre-pre-}\!\sF_{\on{KM}(G)}^{\on{ext}}\circ \Wak^{-,\semiinf}_{\rho_P(\omega_X)}$$
as a factorization functor 
$$\KL(M)_{\crit_M-\rhoch_P}\to \IndCoh^*(\Op^\mer_\cG).$$

By combining \propref{p:F KM Sph via semiinf}, \propref{p:support of Wak Sph}(a) and \corref{c:support of Wak Sph}
we obtain an isomorphism
$$\on{pre-pre-}\!\sF_{\KL(M)} \simeq 
\on{pre-pre-}\!\sF^{\on{ext}}_{\KL(M)}\circ ((\sfq^{\on{Miu,mon-free}})^\wedge)^\IndCoh_*\circ ((\sfq^{\on{Miu,mon-free}})^\wedge)^!$$
as factorization functors
$$\KL(M)_{\crit_M-\rhoch_P}\rightrightarrows\IndCoh^*(\Op^\mer_\cG).$$

\sssec{}

Let $\on{pre-}\!\sF^{\on{ext}}_{\KL(M)}$ be the functor obtained from $\on{pre-pre-}\!\sF^{\on{ext}}_{\KL(M)}$
by the procedure in \secref{sss:pre from pre-pre}. 

\medskip

It follows formally that we have an isomorphism 
$$\on{pre-}\!\sF_{\KL(M)} \simeq 
\on{pre-}\!\sF^{\on{ext}}_{\KL(M)}\circ ((\sfq^{\on{Miu,mon-free}})^\wedge)^\IndCoh_*\circ ((\sfq^{\on{Miu,mon-free}})^\wedge)^!$$
as factorization functors
$$\KL(M)_{\crit_M-\rhoch_P}\rightrightarrows \IndCoh^*(\Op^\mer_\cG)\otimes \Rep(\cM).$$

\sssec{} \label{sss:recover F from A phi geom}

From the fact that the functors $\Accs$ and $((\sfq^{\on{Miu,mon-free}})^\wedge)^\IndCoh_*$
are strictly unital, it follows that the factorization algebra $\CA$ identifies with the image of the factorization
unit under the functor $\on{pre-}\!\sF^{\on{ext}}_{\KL(M)}$.

\medskip

Hence, the homomorphism
$$\phi:\Omega(R)^\Op\to \CA.$$
can be interpreted as a homomorphism
$$\Omega(R)^\Op\to \on{pre-}\!\sF^{\on{ext}}_{\KL(M)}(\one_{\KL(M)_{\crit_M-\rhoch_P}}).$$

Hence, by the same mechanism as in \secref{sss:recover F from A phi}, the functor $\on{pre-}\!\sF^{\on{ext}}_{\KL(M)}$ upgrades to a functor
$$\sF^{\on{ext}}_{\KL(M)}:\KL(M)_{\crit_M-\rhoch_P}\to
\Omega(R)^\Op\mod^{\on{fact}}(\IndCoh^*(\Op^\mer_\cG)\otimes \Rep(\cM))=:\bC.$$

\medskip

It follows formally that we have an isomorphism 
$$\sF_{\KL(M)} \simeq 
\sF^{\on{ext}}_{\KL(M)}\circ ((\sfq^{\on{Miu,mon-free}})^\wedge)^\IndCoh_*\circ ((\sfq^{\on{Miu,mon-free}})^\wedge)^!$$
as factorization functors
$$\KL(M)_{\crit_M-\rhoch_P}\rightrightarrows \bC.$$

\ssec{Identification of the compositions}

\sssec{} \label{ss:ident comp}

As was mentioned above, our goal is to establish the commutativity of the pentagon \eqref{e:category C pent}.
By Sects. \ref{sss:recover F from A phi ext} and \ref{sss:recover F from A phi geom}, this amounts to:

\medskip

\begin{itemize}

\item Establishing an identification
\begin{equation} \label{e:ident pre F}
\on{pre-}\!\sF^{\on{ext}}_{\KL(M)} \simeq \on{pre-}\!\sF^{\on{spec,ext}}_{\MOp(\cM)}
\end{equation} 
as factorization functors from
$$\KL(M)_{\crit_G-\rhoch_P} \overset{\FLE_{M,\crit_G+\rhoch_P}\circ \tau_M} \simeq \IndCoh^*(\Op_{\cM,\rhoch_P}^{\on{mon-free}})$$
to $\IndCoh^*(\Op^\mer_\cG)\otimes \Rep(\cM)$. 

\medskip

\item Showing that with respect to the resulting isomorphism of factorization algebras 
$$\CA=\on{pre-}\!\sF^{\on{ext}}_{\KL(M)}(\one_{\KL(M)_{\crit_G-\rhoch_P}}) \simeq 
\on{pre-}\!\sF^{\on{spec,ext}}_{\MOp(\cM)}(\one_{\IndCoh^*(\Op_{\cM,\rhoch_P}^{\on{mon-free}})})=\CA^{\on{spec}},$$
the homomorphisms $\phi$ and $\phi^{\on{spec}}$ are canonically identified.

\end{itemize} 

\sssec{}

We start by constructing the isomorphism \eqref{e:ident pre F}.

\medskip

Note, however, that according to \secref{sss:pre from pre-pre}, in order to establish \eqref{e:ident pre F}, it suffices to establish an isomorphism
\begin{equation} \label{e:ident pre pre F}
\on{pre-pre-}\!\sF^{\on{ext}}_{\KL(M)} \simeq \on{pre-pre-}\!\sF^{\on{spec,ext}}_{\MOp(\cM)}
\end{equation} 
as factorization functors with values in $\IndCoh^*(\Op^\mer_\cG)$.

\medskip

In this subsection we will construct an isomorphism between the compositions of the two functors in \eqref{e:ident pre pre F} with the forgetful functor
$$\Gamma^\IndCoh(\Op^\mer_\cG,-):\IndCoh^*(\Op^\mer_\cG)\to \Vect,$$
i.e.,
\begin{equation} \label{e:ident pre pre F pre}
\Gamma^\IndCoh(\Op^\mer_\cG,-)\circ \on{pre-pre-}\!\sF^{\on{ext}}_{\KL(M)} \simeq \Gamma^\IndCoh(\Op^\mer_\cG,-)\circ \on{pre-pre-}\!\sF^{\on{spec,ext}}_{\MOp(\cM)}.
\end{equation}

\sssec{}

By construction, the composition 
$$\Gamma^\IndCoh(\Op^\mer_\cG,-) \circ \on{pre-pre-}\!\sF_{\on{KL}(M)}^{\on{ext}}$$
identifies with
\begin{multline*} 
\KL(M)_{\crit_M-\rhoch_P}\overset{\Wak^{-,\semiinf}_{\rho_P(\omega_X)}}\longrightarrow 
(\hg\mod_{\crit_G})^{-,\semiinf}  \to \hg\mod_{\crit_G}^{\fL^+(M)} \overset{\alpha_{\rho(\omega_X),\on{taut}}}\longrightarrow \\
\to \hg\mod_{\crit_G,\rho(\omega_X)}^{\fL^+(M)_{\rho(\omega_X)}}  
\to \hg\mod_{\crit_G,\rho(\omega_X)} \overset{\DS_G}\longrightarrow \Vect.
\end{multline*}

The composition 
$$\Gamma^\IndCoh(\Op^\mer_\cG,-) \circ \on{pre-pre-}\!\sF^{\on{spec,ext}}_{\MOp(\cM)}\circ (\tau_M\circ \FLE_{M,\crit})$$
identifies with
$$\Gamma^\IndCoh(\Op^\mer_\cM,-) \circ (\tau_M\circ \FLE_{M,\crit}),$$
and hence with 
$$\KL(M)_{\crit_M-\rhoch_P} \overset{\alpha_{\rho_M(\omega_X),\on{taut}}}\longrightarrow \KL(M)_{\crit_M-\rhoch_P,\rho_M(\omega_X)}\to
\fm\mod_{\crit_M-\rhoch_P,\rho_M(\omega_X)}\overset{\DS_M}\longrightarrow \Vect.$$

Now, the two functors match by \corref{c:DS of Wak crit}.

\begin{rem}

Note that in computing the above compositions, the Chevalley involutions on either $G$ or $M$ played no role. 
They will, however, play a role for the next step in constructing \eqref{e:ident pre pre F}, see \secref{sss:center on Wak} below.

\end{rem}

\ssec{Upgrade to factorization algebras I} \label{ss:first upgrade}

\sssec{}

Consider now the functor
\begin{equation} \label{e:Gamma Op enh}
(\Gamma^\IndCoh(\Op^\mer_\cG,-))^{\on{enh}}:\IndCoh^*(\Op^\mer_\cG)\to \CO_{\Op^\reg_\cG}\mod^{\on{fact}},
\end{equation} 
see \cite[Formula (4.6)]{GLC2}.\footnote{Note that the subscript ``enh" here refers to the procedure in \cite[Sect. 4.1.6]{GLC2}, rather
than tensoring by $\on{I}(G,P^-)_{\rho_P(\omega_X)}$.} We will now refine \eqref{e:ident pre pre F pre} to an isomorphism
\begin{equation} \label{e:ident pre pre F coarse}
(\Gamma^\IndCoh(\Op^\mer_\cG,-))^{\on{enh}}
\circ \on{pre-pre-}\!\sF^{\on{ext}}_{\KL(M)} \simeq (\Gamma^\IndCoh(\Op^\mer_\cG,-))^{\on{enh}}\circ \on{pre-pre-}\!\sF^{\on{spec,ext}}_{\MOp(\cM)}.
\end{equation}

\medskip

Let $\CB$ and $\CB^{\on{spec}}$
be the factorization algebras (in $\Vect$), obtained by applying the functors in \eqref{e:ident pre pre F pre} to
$\on{Vac}(M)_{\crit_M-\rhoch_P}$, respectively. 

\sssec{}

The isomorphism \eqref{e:ident pre pre F pre} gives rise to an identification
$$\CB\simeq \CB^{\on{spec}}$$
as factorization algebras.

\medskip

Now, each side in \eqref{e:ident pre pre F coarse} gives rise to maps 
$$\psi:\CO_{\Op^\reg_\cG}\to \CB \text{ and } \psi^{\on{spec}}:\CO_{\Op^\reg_\cG}\to \CB^{\on{spec}}.$$

The datum of \eqref{e:ident pre pre F coarse} amounts to an identification 
\begin{equation} \label{e:psi vs psi spec}
\psi=\psi^{\on{spec}}.
\end{equation}

\medskip

Note that $\CB^{\on{spec}}$ (and hence $\CB$) is a \emph{classical} factorization algebra. Hence, it is enough 
to establish \eqref{e:psi vs psi spec} at the pointwise level.

\sssec{}

Denote 
\begin{equation} \label{e:twisted BW}
\BW_{\fg,\crit_G,\rho_P(\omega_X)}:=\oblv_\hg\circ \Wak_{\rho_P(\omega_X)}(\on{Vac}(M)_{\crit_M-\rhoch_P}),
\end{equation}
cf. \secref{sss:BW}. 

\medskip

Note that we have a canonical identification
$$\BW_{\fg,\crit_G,\rho_P(\omega_X)}\simeq \BW_{\fg,\crit_G}:=\oblv_\hg\circ \Wak(\on{Vac}(M)_{\crit_M-\rhoch_P})$$
(indeed, the twisted by the $Z_M$-torsor $\rho_P(\omega)$ acts as identity on vacuum objects). 

\sssec{}

Consider the corresponding homomorphism
$$\BV_{\fg,\crit_G}\to \BW_{\fg,\crit_G,\rho_P(\omega_X)}.$$

\medskip

By construction, the homomorphism $\psi$ equals
$$\CO_{\Op^\reg_\cG} \overset{\tau_G}\to \CO_{\Op^\reg_\cG} \overset{\on{FF}_G}\simeq \fz_\fg\to 
\DS_G\circ \alpha_{\rho(\omega_X),\on{taut}}(\BV_{\fg,\crit})\to \DS_G\circ \alpha_{\rho(\omega_X),\on{taut}}(\BW_{\fg,\crit_G,\rho_P(\omega_X)}).$$

\sssec{}

Note, however, that the map
\begin{equation} \label{e:center to Wak}
\fz_\fg\to \DS_G\circ \alpha_{\rho(\omega_X),\on{taut}}(\BV_{\fg,\crit})\to \DS_G\circ \alpha_{\rho(\omega_X),\on{taut}}(\BW_{\fg,\crit_G,\rho_P(\omega_X)}),
\end{equation} 
which appears in the above composition, can be rewritten as follows:

\medskip

It is known that when we view $\BW_{\fg,\crit_G,\rho_P(\omega_X)}\simeq \BW_{\fg,\crit_G}$ as a factorization module over $\fz_\fg$, 
it is \emph{commutative}\footnote{See \cite[Sect. B.10.5-6]{GLC2} what this means.} 
(see, e.g., \cite[Sect. 10]{FG2}). In particular, we have an action map
\begin{equation} \label{e:center on Wak}
\fz_\fg\otimes \BW_{\fg,\crit_G,\rho_P(\omega_X)}\to \BW_{\fg,\crit_G,\rho_P(\omega_X)}.
\end{equation} 

By transport of structure, we obtain a map 
\begin{equation} \label{e:center on BRST Wak}
\fz_\fg\otimes  \DS_G\circ \alpha_{\rho(\omega_X),\on{taut}}(\BW_{\fg,\crit_G,\rho_P(\omega_X)})\to  \DS_G\circ \alpha_{\rho(\omega_X),\on{taut}}(\BW_{\fg,\crit_G,\rho_P(\omega_X)}).
\end{equation} 

The map \eqref{e:center to Wak} is obtained from \eqref{e:center on BRST Wak} by acting on the vacuum vector.

\sssec{}

The map $\psi^{\on{spec}}$ is described as follows. We identify 
$$\CB^{\on{spec}}\simeq \DS_M\circ \alpha_{\rho_M(\omega_X),\on{taut}}(\on{Vac}(M)_{\crit_M-\rhoch_P}),$$
and in terms of this identification $\psi^{\on{spec}}$ equals
\begin{multline*}
\CO_{\Op^\reg_\cG} \overset{(\sfp^{\on{Miu}})^*}\longrightarrow \CO_{\MOp^\reg_{\cG,\cP^-}} \overset{\sfq^{\on{Miura}}}\simeq \\
\simeq \CO_{\Op^\reg_{\cM,\rhoch_P}}\overset{\on{FF}_M}\simeq \fz_{\fm,\rhoch_M(\omega_X)} \overset{\tau_M}\to 
\fz_{\fm,-\rhoch_M(\omega_X)} \simeq \DS_M\circ \alpha_{\rho_M(\omega_X),\on{taut}}(\on{Vac}(M)_{\crit_M-\rhoch_P}).
\end{multline*}

\medskip

Note the last map in the above composition, i.e., 
\begin{equation} \label{e:center to Wak M}
\fz_{\fm,-\rhoch_M(\omega_X)} \simeq \DS_M\circ \alpha_{\rho_M(\omega_X),\on{taut}}(\on{Vac}(M)_{\crit_M-\rhoch_P})
\end{equation} 
is by definition constructed as follows (see \cite[Sect. 4.8.2]{GLC2}):

\medskip

We consider $\on{Vac}(M)_{\crit_M-\rhoch_P}$ as acted on by $\fz_{\fm,-\rhoch_M(\omega_X)}$. Hence, by transport of structure, 
$\DS_M\circ \alpha_{\rho_M(\omega_X),\on{taut}}(\on{Vac}(M)_{\crit_M-\rhoch_P})$
also acquires an action of  $\fz_{\fm,-\rhoch_M(\omega_X)}$. The map \eqref{e:center to Wak M} is given by the action of $\fz_{\fm,-\rhoch_M(\omega_X)}$ on the vacuum
vector of $\DS_M\circ \alpha_{\rho_M(\omega_X),\on{taut}}(\on{Vac}(M)_{\crit_M-\rhoch_P})$. 

\sssec{} \label{sss:center on Wak}

Note that by the construction of the isomorphism \eqref{e:ident pre pre F pre}, the isomorphism $\CB\simeq \CB^{\on{spec}}$ comes from the identification
$$\DS_G\circ \alpha_{\rho(\omega_X),\on{taut}} \circ \Wak_{\rho_P(\omega_X)}(\on{Vac}(M)_{\crit_M-\rhoch_P}) \simeq
\DS_M\circ \alpha_{\rho_M(\omega_X),\on{taut}}(\on{Vac}(M)_{\crit_M-\rhoch_P}).$$

\medskip

Hence, in order to construct the identification \eqref{e:psi vs psi spec}, it suffices to prove the following: 

\medskip

The action of $\fz_\fg$ on $\BW_{\fg,\crit_G}$ in \eqref{e:center on Wak} equals the action obtained via
$$\fz_\fg\overset{\tau_G}\to
\fz_\fg\overset{\on{FF}_G}\simeq \CO_{\Op^\reg_\cG} \overset{(\sfp^{\on{Miu}})^*}\longrightarrow \CO_{\MOp^\reg_{\cG,\cP^-}} \overset{\sfq^{\on{Miura}}}
\simeq \CO_{\Op^\reg_{\cM,\rhoch_P}}\overset{\on{FF}_M}\simeq \fz_{\fm,\rhoch_M(\omega_X)} \overset{\tau_M}\to 
\fz_{\fm,-\rhoch_M(\omega_X)}$$
and the action of $\fz_{\fm,-\rhoch_M(\omega_X)}$ on 
$$\BW_{\fg,\crit_G}=\oblv_{\hg}\circ \Wak(\on{Vac}(M)_{\crit_M-\rhoch_P})$$
by transport of structure.

\medskip

However, the latter is the basic feature of the Feigin-Frenkel isomorphism\footnote{See \corref{c:Miura and BRST} for the explanation
of the appearance of the Chevalley involution in the above formula.}, see \cite[Theorem 11.3]{Fr} or \cite[Sect. 10]{FG2}.\footnote{These
references treat the case of $P=B$. The case of the general parabolic follows formally, since the vacuum module for a given Levi 
\emph{embeds} into the corresponding principal Wakimoto module.}  

\ssec{Interlude: Miura transform and the Wakimoto functor}

\sssec{}

Recall that according to \cite[Sect. 4.6]{GLC2}, we have an action of the monoidal factorization category
$$\QCoh(\Op^\mer_\cG) \overset{\Upsilon_{\Op^\mer_\cG}}\simeq \IndCoh^!(\Op^\mer_\cG)$$
on $\hg\mod_{\crit_G}$. In what follows we will twist this action by the Chevalley involution $\tau_G$. 

\medskip

Consider the induced action on the full subcategory 
$$(\hg\mod_{\crit_G})^{\fL^+(N_P)}\subset \hg\mod_{\crit_G}.$$

Consider also a similar action of
$$\QCoh(\Op^\mer_{\cM,\rhoch_P}) \overset{\Upsilon_{\Op^\mer_{\cM,\rhoch_P}}}\simeq \IndCoh^!(\Op^\mer_{\cM,\rhoch_P})$$
on $\hm\mod_{\crit_M-\rhoch_P}$, i.e., we compose the action from \cite[Sect. 4.6]{GLC2} with $\tau_M$. 

\sssec{}

Unwinding the construction in \cite[Sect. 4.6]{GLC2}, from the identification of the $\fz_\fg$-action on $\BW_{\fg,\crit_G}$ 
in \secref{sss:center on Wak} above, we obtain:

\begin{cor} \label{c:Miura and Wak}
The functor 
$$\Wak:\hm\mod_{\crit_M-\rhoch_P}\to (\hg\mod_{\crit_G})^{\fL^+(N_P)}$$
is compatible with the above actions via the monoidal functor
\begin{equation} \label{e:Miura transform IndCoh}
\IndCoh^!(\Op^\mer_\cG)\overset{(\sfp^{\on{Miu}})^!}\longrightarrow \IndCoh^!(\MOp^\mer_{\cG,\cP^-}) \overset{\sfq^{\on{Miu}}}\simeq \IndCoh^!(\Op^\mer_{\cM,\rhoch_P}).
\end{equation} 
\end{cor} 

\sssec{}

Note now that the self-duality
$$(\hg\mod_\crit)^\vee\simeq \hg\mod_\crit$$
is linear with respect to the action of $\IndCoh^!(\Op^\mer_\cG)$, \emph{up to the Chevalley involution}, see \cite[Theorem 8.2.3]{GLC2}.
A similar assertion applies to $G$ replaced by $M$.

\medskip

Hence, from \corref{c:Miura and Wak}, we obtain:

\begin{cor} \label{c:Miura and BRST}
The composite functor
$$(\hg\mod_{\crit_G})^{\fL^+(N_P)}\to \hg\mod_{\crit_G} \overset{J^{-}_{\on{KM}}}\longrightarrow \hm\mod_{\crit_M+\rhoch_P}$$
is $\IndCoh^!(\Op^\mer_\cG)$-linear via \eqref{e:Miura transform IndCoh}, where:
\begin{itemize}

\item $\IndCoh^!(\Op^\mer_\cG)$ acts on $(\hg\mod_{\crit_G})^{\fL^+(N_P)}$ as in \cite[Sect. 4.6]{GLC2} (without $\tau_G$);

\medskip

\item $\IndCoh^!(\Op^\mer_{\cM,\rhoch_P})$ acts on $\hm\mod_{\crit_M+\rhoch_P}$ by a twisted version of \cite[Sect. 4.6]{GLC2} 
(without $\tau_M$).

\end{itemize} 
\end{cor} 

\ssec{Proof of \propref{p:support of Wak Sph}} \label{ss:support of Wak Sph}

\sssec{}

First, we note that both points of the proposition are enough to check strata-wise over $\Ran$. We can further assume that we work over
a fixed point $\ul{x}\in \Ran$.

\sssec{}

Let $I_{P,\ul{x}}\subset \fL^+(G)_{\ul{x}}$ be the parahoric subgroup corresponding to $P$. Recall that the functor
$\on{Av}^{\fL^+(N_P)}_*$ defines equivalences
$$(\hg\mod_{\crit_G,\ul{x}})^{-,\semiinf}\to (\hg\mod_{\crit_G,\ul{x}})^{I_{P,\ul{x}}}$$
and
$$\Dmod_{\frac{1}{2}}(\Gr_{G,\ul{x}})^{\fL(N^-_P)_{\ul{x}}\cdot \fL^+(M)_{\ul{x}},\on{ren}} \to \Dmod_{\frac{1}{2}}(\Gr_{G,\ul{x}})^{I_{P,\ul{x}},\on{ren}}.$$

\medskip

Hence, in the statement of \propref{p:support of Wak Sph} we can replace
$$(\hg\mod_{\crit_G,\ul{x}})^{-,\semiinf} \rightsquigarrow (\hg\mod_{\crit_G})^{I_{P,\ul{x}}};$$
$$\Wak^{-,\semiinf}_{\rho_P(\omega_X)} \rightsquigarrow \Wak;$$
$$\on{I}(G,P^-)_{\rho_P(\omega_X)}  \rightsquigarrow \Dmod_{\frac{1}{2}}(\Gr_{G,\ul{x}})^{I_{P,\ul{x}},\on{ren}}.$$

\sssec{}

We regard the category $(\hg\mod_{\crit_G})^{I_{P,\ul{x}}}$ as acted on by
$$\QCoh(\Op^\mer_{\cG,\ul{x}})\overset{\Upsilon_{\Op^\mer_{\cG,\ul{x}}}}\simeq \IndCoh^!(\Op^\mer_{\cG,\ul{x}}).$$

Let 
$$(\hg\mod_{\crit_G})^{I_{P,\ul{x}}}_{\mf}\subset (\hg\mod_{\crit_G})^{I_{P,\ul{x}}}$$
be the full subcategory consisting of objects that are set-theoretically supported over $\Op^\mf_{\cG,\ul{x}}\subset \Op^\mer_{\cG,\ul{x}}$.

\medskip

Denote the above fully faithful embedding by $((\iota^\mf)^\wedge)^\IndCoh_*$, and by $((\iota^\mf)^\wedge)^!$ its right adjoint. 

\sssec{}

It follows from \corref{c:Miura and Wak} that the functor $\Wak$ sends $(\KL(M)_{\crit_M-\rhoch_P})_{\Gmf}$ to $(\hg\mod_{\crit_G})^{I_{P,\ul{x}}}_{\mf}$,
i.e., we have a commutative diagram
\begin{equation} \label{e:support of Wak usual Sph 1}
\CD
(\KL(M)_{\crit_M-\rhoch_P})_{\Gmf} @>{((\sfq^{\on{Miu,mon-free}})^\wedge)^\IndCoh_*}>> \KL(M)_{\crit_M-\rhoch_P} \\
@V{\Wak}VV  @VV{\Wak}V \\
(\hg\mod_{\crit_G})^{I_{P,\ul{x}}}_{\mf} @>{((\iota^\mf)^\wedge)^\IndCoh_*}>> (\hg\mod_{\crit_G})^{I_{P,\ul{x}}}.
\endCD
\end{equation}

Moreover, the natural transformation in the diagram
\begin{equation} \label{e:support of Wak usual Sph 2}
\vcenter
{\xy
(0,0)*+{(\KL(M)_{\crit_M-\rhoch_P})_{\Gmf}}="A";
(70,0)*+{\KL(M)_{\crit_M-\rhoch_P}}="B";
(0,-20)*+{(\hg\mod_{\crit_G})^{I_{P,\ul{x}}}_{\mf}}="C";
(70,-20)*+{(\hg\mod_{\crit_G})^{I_{P,\ul{x}}},}="D";
{\ar@{->}_{((\sfq^{\on{Miu,mon-free}})^\wedge)^!} "B";"A"};
{\ar@{->}^{((\iota^\mf)^\wedge)^!} "D";"C"};
{\ar@{->}_{\Wak} "A";"C"};
{\ar@{->}^{\Wak} "B";"D"};
\endxy}
\end{equation}
obtained by passing to right adjoints along the horizontal arrows in \eqref{e:support of Wak usual Sph 1}, is an isomorphism. 

\sssec{}

Note now that we have an inclusion of subcategories of $(\hg\mod_{\crit_G})^{I_{P,\ul{x}}}$:
\begin{equation} \label{e:mon via Sph gen}
\KL(G)_{\crit}\underset{\Sph_G}\otimes  \Dmod_{\frac{1}{2}}(\Gr_{G,\ul{x}})^{I_{P,\ul{x}},\on{ren}} \subset (\hg\mod_{\crit})^{I_{P,\ul{x}}}_{\mf}.
\end{equation}

Hence, in order to prove \propref{p:support of Wak Sph}, it suffices to show that the inclusion \eqref{e:mon via Sph gen} is an equality.
In fact, we claim:

\begin{prop} \label{p:mon via Sph gen}
The inclusion 
$$\KL(G)_{\crit}\underset{\Sph_G}\otimes  \Dmod_{\frac{1}{2}}(\Gr_{G,\ul{x}}) \subset (\hg\mod_{\crit})_{\mf}$$
is an equality.
\end{prop} 

\begin{proof}

The assertion of the proposition is equivalent to the statement that $(\hg\mod_{\crit})_{\mf}$ is generated, as a category
acted on by $\fL(G)$, by its $\fL(G)^+$-equivariant objects. 

\medskip

A version of this statement for $\overset{\circ}{I}$ (the unipotent radical of the Iwahori subgroup $I$) is \cite[Theorem 4.3]{YaD}. 
One passes from $\overset{\circ}{I}$ to $\fL(G)^+$ using the same methods as in \cite{FG3}:

\medskip

In {\it loc. cit.} it is shown that the category 
$$(\hg\mod_{\crit})^{\overset{\circ}{I}}_\reg\subset (\hg\mod_{\crit})^{\overset{\circ}{I}},$$
consisting of objects that are set-theoretically supported over $\Op^\reg_\cG\subset \Op^\mer_\cG$, 
is spherically generated. 

\medskip

The same arguments apply to all other connected components of $\Op^\mf_\cG$. 

\end{proof} 

\ssec{Upgrade to an isomorphism with values in \texorpdfstring{$\IndCoh^*(\Op^\mer_\cG)$}{IndCohOp}}

Recall that the functor \eqref{e:Gamma Op enh} is ``almost" an equivalence. More precisely, it induces an equivalence between
the eventually coconnective parts of the two sides with respect to the natural t-structures, see \cite[Lemma 4.2.3]{GLC2}. 

\medskip

Thus, the isomorphism \eqref{e:ident pre pre F coarse} constructed above implies that \eqref{e:ident pre pre F} holds
\emph{up to renormalization}. In this subsection, we will carry out the construction of \eqref{e:ident pre pre F} ``as-is". 

\sssec{}

Since the category $\KL(M)_{\crit_M-\rhoch_P}$ is compactly generated (see \cite[Sect. B.11.10]{GLC2} for what this means),
it suffices to construct an identification between the restrictions of the two functors in \eqref{e:ident pre pre F} to the
subcategory of compact objects.

\medskip

By \cite[Corollary 4.4.2(a)]{GLC2}, it suffices to show that both functors in \eqref{e:ident pre pre F} send compact objects
in $\KL(M)_{\crit_M-\rhoch_P}$ to eventually coconnective objects in $\IndCoh^*(\Op^\mer_\cG)$. 

\medskip

Since compact objects in $\KL(M)_{\crit_M-\rhoch_P}$ are eventually coconnective (over any $X^I$), it suffices to show that both functors in 
\eqref{e:ident pre pre F} have a cohomological amplitude on the left (over any $X^I$). The latter property can be checked at the pointwise on $X^I$. 

\medskip

Thus, for the duration of this subsection we will work at a fixed point $\ul{x}\in \Ran$. Moreover, by factorization,
we can assume that $\ul{x}$ is a singleton $\ul{x}=\{x\}$. 

\sssec{}

We first consider the functor $\on{pre-pre-}\!\sF^{\on{spec,ext}}_{\MOp(\cM)}$. 

\medskip

We claim that the functor in question is in fact t-exact. Indeed, $\on{pre-pre-}\!\sF^{\on{spec,ext}}_{\MOp(\cM)}$ is the functor
\eqref{e:pre-pre-F ext}, and the t-exactness follows from \cite[Corollary A.8.8]{GLC2}.

\sssec{}

We now show that the functor $\on{pre-pre-}\!\sF^{\on{ext}}_{\KL(M)}$ has a bounded cohomological amplitude. 

\medskip

As in \secref{sss:DS of usual Wak}, since the functor $\DS^{\on{enh,rfnd}}$ factors as
$$\hg\mod_{\crit_G,\rho(\omega_X)} \to \Whit_*(\hg\mod_{\crit_G,\rho(\omega_X)}) \overset{\ol\DS^{\on{enh,rfnd}}}\longrightarrow \\
\to \IndCoh^*(\Op^\mer_\cG)$$
(see \cite[Theorem 4.8.8]{GLC2}), the functor $\on{pre-pre-}\!\sF^{\on{ext}}_{\KL(M)}$ identifies with the composition
\begin{multline*} 
\KL(M)_{\crit_M-\rhoch_P} \overset{\Wak_{\rho_P(\omega_X)}}\longrightarrow 
\hg\mod_{\crit_G}^{\fL^+(P)} 
\overset{\alpha_{\rho(\omega_X),\on{taut}}}\longrightarrow 
\hg\mod_{\crit_G,\rho(\omega_X)}^{\fL^+(P)_{\rho(\omega_X)}}  
\to \hg\mod_{\crit_G,\rho(\omega_X)} \overset{\DS^{\on{enh,rfnd}}}\longrightarrow \\
\to \IndCoh^*(\Op^\mer_\cG)
\overset{\tau_G}\to  \IndCoh^*(\Op^\mer_\cG),
\end{multline*}
i.e., we can replace the semi-infinite Wakimoto functor $\Wak^{-,\semiinf}_{\rho_P(\omega_X)}$ by the ``usual" Wakimoto
functor $\Wak_{\rho_P(\omega_X)}$.

\medskip

Thus, it suffices to show that the composition 
\begin{multline*} 
\KL(M)_{\crit_M-\rhoch_P} \overset{\Wak_{\rho_P(\omega_X)}}\longrightarrow 
\hg\mod_{\crit_G}^{\fL^+(P)} 
\overset{\alpha_{\rho(\omega_X),\on{taut}}}\longrightarrow 
\hg\mod_{\crit_G,\rho(\omega_X)}^{\fL^+(P)_{\rho(\omega_X)}}  
\to \hg\mod_{\crit_G,\rho(\omega_X)} \overset{\DS^{\on{enh,rfnd}}}\longrightarrow \\
\to \IndCoh^*(\Op^\mer_\cG)
\end{multline*}
has a bounded cohomological amplitude. 

\sssec{}

It follows from the explicit shape of the functor $\Wak_{\rho_P(\omega_X)}$ that it is t-exact (at the pointwise level). In addition, from the description of 
$\Wak_{\rho_P(\omega_X)}$ in \secref{sss:usual Wak via CDO}, it follows that it factors via
$$\hg\mod_{\crit_G}^{I_P}\subset \hg\mod_{\crit_G,x}^{\fL^+(P)_x},$$
where $I_P\subset \fL^+(G)_x$ is the parahoric subgroup, corresponding to $P$.

\medskip

Hence, it suffices to show that the functor 
$$\hg\mod_{\crit_G}^{I_P} \overset{\alpha_{\rho(\omega_X),\on{taut}}}\longrightarrow 
\hg\mod_{\crit_G,\rho(\omega_X)}^{I_{P,\rho(\omega_X)}} \to \hg\mod_{\crit_G,\rho(\omega_X)} \overset{\DS^{\on{enh,rfnd}}}\longrightarrow \IndCoh^*(\Op^\mer_\cG)$$
has a bounded cohomological amplitude. 

\medskip

However, this follows from \cite[Lemma 5.3.1]{Ra2}.\footnote{More precisely, we first *-average from $\hg\mod_{\crit_G}^{I_P}$ to the ``baby Whittaker category", 
which has a bounded amplitude, then apply \cite[Lemma 5.3.1]{Ra2}, which says that the projection from the baby Whittaker category to 
$\Whit_*(\hg\mod_{\crit,\rho(\omega_X)})$ is t-exact, up to a shift.} 

\ssec{Upgrade to factorization algebras II}

So far, we have carried out the construction of the isomorphism of functors in the first bullet point 
in \secref{ss:ident comp}. In particular, we obtain an isomorphism of factorization algebras
\begin{equation} \label{e:A A spec}
\CA\overset{\nu}\simeq \CA^{\on{spec}}
\end{equation}
in $\IndCoh^*(\Op^\mer_\cG)\otimes \Rep(\cM)$.

\medskip

Thus, it remains to construct an identification of the second bullet point, i.e., an identification 
$$\phi\simeq \phi^{\on{spec}}$$
as maps from $\Omega(R)^\Op$ to the two sides of \eqref{e:A A spec}. 

\medskip

Our strategy will be as follows: we will construct \emph{another} identification
\begin{equation} \label{e:A A spec prime}
\CA\overset{\nu'}\simeq \CA^{\on{spec}}
\end{equation}
as factorization algebras in $\IndCoh^*(\Op^\mer_\cG)\otimes \Rep(\cM)$, which is \emph{a priori compatible} 
with the maps $\psi$ and $\psi^{\on{spec}}$. We will then show that $\nu'$ equals the already existing 
isomorphism $\nu$. 

\sssec{}

Note that the isomorphism \eqref{e:ident pre F} established above implies that the pentagon \eqref{e:category C pent} becomes
commutative after we apply the forgetful functor
$$\bC\to \IndCoh^*(\Op^\mer_\cG)\otimes \Rep(\cM),$$
and further the functor
\begin{equation} \label{e:Gamma Op Rep M}
\Gamma^\IndCoh(\Op^\mer_\cG,-)\otimes \on{Id}:\IndCoh^*(\Op^\mer_\cG)\otimes \Rep(\cM)\to \Rep(\cM).
\end{equation} 

In other words, the diagram \eqref{e:dual Jacquet enh diagram} becomes commutative after we compose it with
\begin{multline} \label{e:add a leg} 
\IndCoh^*(\Op_\cG^{\on{mon-free}}) \underset{\Sph_\cG^{\on{spec}}}\otimes \on{I}(\cG,\cP^-)^{\on{spec,loc}} 
\overset{\on{Id}\otimes \on{pre-}\!\sF^{\on{spec}}}\longrightarrow  \\
\to \IndCoh^*(\Op_\cG^{\on{mon-free}}) \underset{\Sph_\cG^{\on{spec}}}\otimes \Rep(\cG)\otimes \Rep(\cM)\to \\
\to \IndCoh^*(\Op^\mer_\cG) \otimes \Rep(\cM) \overset{\Gamma^\IndCoh(\Op^\mer_\cG,-)\otimes \on{Id}}\longrightarrow \Rep(\cM),
\end{multline} 
where $\on{pre-}\!\sF^{\on{spec}}$ is as in \secref{sss:pre F}. 

\sssec{}

We will take the following input from the paper \cite[Theorem 4.11]{FG1}, adapted to the parahoric setting: 

\begin{thm} \label{t:FG1} 
At the \emph{pointwise level}, the two circuits of the diagram \eqref{e:dual Jacquet enh diagram} commute when evaluated 
on the object $\on{Vac}(M)_{\crit_M-\rhoch_P}\in \KL(M)_{\crit_M-\rhoch_P}$. Moreover: 

\medskip

\noindent{\em(i)} This isomorphism is compatible with the maps to each side
from 
$$\one_{\IndCoh^*(\Op_\cG^{\on{mon-free}}) \underset{\Sph_\cG^{\on{spec}}}\otimes \on{I}(\cG,\cP^-)^{\on{spec,loc}}}
\in\IndCoh^*(\Op_\cG^{\on{mon-free}}) \underset{\Sph_\cG^{\on{spec}}}\otimes \on{I}(\cG,\cP^-)^{\on{spec,loc}};$$

\medskip

\noindent{\em(ii)}  This isomorphism is compatible with the identification obtained by applying the functor \eqref{e:add a leg}. 

\end{thm}

We will supply a proof in \secref{s:proof of Wak}. 

\sssec{}

Denote the two unital factorization algebras obtained by applying the two circuits of the diagram \eqref{e:dual Jacquet enh diagram} to $\on{Vac}(M)_{\crit_M-\rhoch_P}$ 
by $$\wt\CA \text{ and } \wt\CA^{\on{spec}},$$
respectively. 

\medskip 

We will use \thmref{t:FG1} to deduce the following:

\begin{cor} \label{c:FG1} 
The unital factorization algebras $\wt\CA$ and $\wt\CA^{\on{spec}}$ are canonically
isomorphic.
\end{cor} 

\begin{proof}

Note that for an object 
$$\CF\in \IndCoh^*((\Op^\Mmf_{\cG,\cP^-})^\wedge_\mf)\simeq \IndCoh^*(\Op_\cG^{\on{mon-free}}) \underset{\Sph_\cG^{\on{spec}}}\otimes \on{I}(\cG,\cP^-)^{\on{spec,loc}}$$
and an open subscheme 
$$U\subset \Op^\Pmf_{\cG,\cP^-},$$
it makes sense to talk about the \emph{localization} of $\CF$ on $U$; to be denoted $\CF_U$. It comes equipped with a universal map
$$\CF\to \CF_U.$$

\medskip

Take $U:=\MOp^\reg_{\cG,\cP^-}$. Tautologically, the map
$$\CO_{\Op^\Pmf_{\cG,\cP^-}}\overset{\sim}\to 
\one_{\IndCoh^*((\Op^\Mmf_{\cG,\cP^-})^\wedge_\mf)}\to \wt\CA^{\on{spec}}$$ 
identifies $\wt\CA^{\on{spec}}$ with the localization of $\CO_{\Op^\Pmf_{\cG,\cP^-}}$ on $U$. 

\medskip

In order to construct an isomorphism $\wt\CA\simeq \wt\CA^{\on{spec}}$, it suffices to show that the map
\begin{equation} \label{e:unit to tilde A}
\one_{\IndCoh^*((\Op^\Mmf_{\cG,\cP^-})^\wedge_\mf)}\to \wt\CA^{\on{spec}}
\end{equation} 
also identifies $\wt\CA$ with the localization of $\one_{\IndCoh^*((\Op^\Mmf_{\cG,\cP^-})^\wedge_\mf)}$ on $U$.

\medskip

The property of a map to be a localization on a given open can be checked strata-wise. Hence, since the map \eqref{e:unit to tilde A} is compatible
with factorization, this property can be checked at the pointwise level. 

\medskip

Hence, the required property holds follows from the isomorphism stated in \thmref{t:FG1} and the compatibility in point (i). 

\end{proof} 

\sssec{}

Applying the functor
$$\IndCoh^*(\Op_\cG^{\on{mon-free}}) \underset{\Sph_\cG^{\on{spec}}}\otimes \on{I}(\cG,\cP^-)^{\on{spec,loc}}\overset{\sF^{\on{spec}}_{\Op(\cG)}}\longrightarrow 
\Omega(R)^\Op\mod^{\on{fact}}(\IndCoh^*(\Op^\mer_\cG)\otimes \Rep(\cM))=:\bC,$$
to the isomorphism
$$\wt\CA\overset{\wt\nu'}\simeq \wt\CA^{\on{spec}}$$
of \corref{c:FG1}, we obtain the sought-for isomorphism $\nu'$ in \eqref{e:A A spec prime}. 

\medskip

To finish the proof, we need to show that $\nu'$ equals the isomorphism $\nu$ of \eqref{e:A A spec}.

\sssec{}

Note that both objects in \eqref{e:A A spec} are factorization algebras and as such they are \emph{classical}
(i.e., their values on $X_\dr$ shifted by $[-1]$ belong to the heart of the t-structure). Hence, in order to show that two 
given morphisms between them are equal, it is enough to do so at the pointwise level.

\medskip 

At the pointwise level, in order to show that two given maps between objects of 
$$(\IndCoh^*(\Op^\mer_\cG)\otimes \Rep(\cM))^\heartsuit$$
are equal, it is enough to show that this is the case after applying the functor \eqref{e:Gamma Op Rep M}. 

\medskip

The required assertion follows now the compatibility of point (ii) in \thmref{t:FG1}.

\section{Proof of \thmref{t:local Jacquet dual enh}, continuation} \label{s:ff in plus}

In the previous section, we introduced the (lax) factorization category $\bC$ and constructed
the commutative diagrams \eqref{e:category C triang} and \eqref{e:category C pent}. 

\medskip

In this section we will introduce the relevant t-structures and prove the first two bullet points
in \secref{sss:strategy Jacquet dual enh}. 

\ssec{Structure of the argument}

\sssec{} \label{sss:t-str properties}

Our proof will consist of the following steps:

\begin{itemize}

\item We will show that the functor
$$\IndCoh^*(\Op_\cG^{\on{mon-free}}) \underset{\Sph_\cG^{\on{spec}}}\otimes \on{I}(\cG,\cP^-)^{\on{spec,loc}}
\overset{\sF^{\on{spec}}_{\Op(\cG)}}\longrightarrow \bC$$
is fully-faithful ``as-is" at the pointwise level over $\Ran$ (i.e. over each $\ul{x}\in \Ran$);

\medskip

\item We will introduce a t-structure on $\IndCoh^*(\Op_\cG^{\on{mon-free}}) \underset{\Sph_\cG^{\on{spec}}}\otimes \on{I}(\cG,\cP^-)^{\on{spec,loc}}$
and show that the functor $\sF^{\on{spec}}_{\Op(\cG)}$ is fully faithful when restricted to the eventually coconnective subcategory;

\medskip

\item We will show that the functor 
$$\IndCoh^*(\Op_{\cM,\rhoch_P}^{\on{mon-free}}) \overset{\on{co}\!J^{-,\on{spec},\on{enh}}_\Op}\longrightarrow 
\IndCoh^*(\Op_\cG^{\on{mon-free}}) \underset{\Sph_\cG^{\on{spec}}}\otimes \on{I}(\cG,\cP^-)^{\on{spec,loc}}$$
has a cohomological amplitude bounded on the left.

\end{itemize}

\medskip

Let us show how the above three properties, together with the commutativity of \eqref{e:category C triang} and \eqref{e:category C pent},
which has been already established, imply the first two bullet points in \secref{sss:strategy Jacquet dual enh}, thereby completing the proof
of \thmref{t:local Jacquet dual enh}. 

\sssec{}

First, we note that the first bullet point in \secref{sss:t-str properties} combined with the commutativity of \eqref{e:category C pent}
implies that the diagram

\smallskip

\begin{equation}  \label{e:dual Jacquet enh diagram x}
\CD
\KL(M)_{\crit_G-\rhoch_P,\ul{x}} @>{\FLE_{M,\crit_G+\rhoch_P}\circ \tau_M}>> \IndCoh^*(\Op_{\cM,\rhoch_P}^{\on{mon-free}})_{\ul{x}} \\
@V{\Wak^{-,\on{enh}}_{\rho_P(\omega_X)}}VV @VV{\on{co}\!J^{-,\on{spec},\on{enh}}_\Op}V \\
\KL(G)_{\crit_G,\ul{x}}\underset{\Sph_{G,\ul{x}}}\otimes \on{I}(G,P^-)^{\on{loc}}_{\rho_P(\omega_X),\ul{x}} @>{\text{\eqref{e:tensored up Sat dual}}}>>
\IndCoh^*(\Op_{\cG,\ul{x}}^{\on{mon-free}}) \underset{\Sph_{\cG,\ul{x}}^{\on{spec}}}\otimes \on{I}(\cG,\cP^-)_{\ul{x}}^{\on{spec,loc}}.
\endCD
\end{equation}
commutes for every $\ul{x}\in \Ran$.

\medskip

I.e., we obtain that \thmref{t:local Jacquet dual enh} holds at the pointwise level. 

\sssec{}

The t-structure on 
$$\KL(G)_{\crit_G}\underset{\Sph_G}\otimes \on{I}(G,P^-)^{\on{loc}}_{\rho_P(\omega_X)}$$
is defined by transferring the structure on 
$$\IndCoh^*(\Op_\cG^{\on{mon-free}}) \underset{\Sph_\cG^{\on{spec}}}\otimes \on{I}(\cG,\cP^-)^{\on{spec,loc}}$$ 
via the equivalence \eqref{e:tensored up Sat dual}. 

\medskip

The fact that the functor $\sF_{\on{KM}(G)}$ is fully faithful when restricted to the eventually coconnective subcategory follows
tautologically from the commutativity of \eqref{e:category C triang} and the corresponding property of $\sF^{\on{spec}}_{\Op(\cG)}$.

\sssec{}

Since compact objects in $\IndCoh^*(\Op_{\cM,\rhoch_P}^{\on{mon-free}})$ are eventually coconnective, the fact that 
$\on{co}\!J^{-,\on{spec},\on{enh}}_\Op$ has a cohomological amplitude bounded on the left implies that it sends compact
objects to objects that are eventually coconnective. 

\medskip

Note that compact objects in $\KL(M)_{\crit_G-\rhoch_P}$ are also eventually coconnective. Hence, it remains to show
that the functor $\Wak^{-,\on{enh}}_{\rho_P(\omega_X)}$ has a cohomological amplitude bounded on the left. 

\medskip

Note, however, that this property can be checked pointwise on $\Ran$.  Recall that the functor $\FLE_{M,\crit_G+\rhoch_P}$
is t-exact at the pointwise level (see \cite[Corollary 6.6.9]{GLC2}). Hence, the required property holds from the commutativity
of \eqref{e:dual Jacquet enh diagram x} and the corresponding property of $\on{co}\!J^{-,\on{spec},\on{enh}}_\Op$. 

\sssec{}

Thus, the rest of this section is devoted to carrying out the steps outlined in \secref{sss:t-str properties}. 

\ssec{The t-structure on the spectral side}

We will define a t-structure on 
$$\IndCoh^*(\Op_\cG^{\on{mon-free}}) \underset{\Sph_\cG^{\on{spec}}}\otimes \on{I}(\cG,\cP^-)^{\on{spec,loc}}$$
by the same principle as in \secref{sss:t semiinf spec}.

\sssec{}

Consider the factorization space
$$\Op^\Pmf_{\cG,\cP^-}\simeq \Op^\mf_\cG\underset{\LS^\reg_\cG}\times \LS^\reg_{\cP^-},$$
and the factorization category $\IndCoh^*(\Op^\Pmf_{\cG,\cP^-})$, which we rewrite as 
\begin{equation} \label{e:Op Pmf}
\IndCoh^*(\Op^\mf_\cG)\underset{\Rep(\cG)}\otimes \Rep(\cP^-).
\end{equation}

\medskip

We endow \eqref{e:Op Pmf} with a tensor product t-structure. 

\sssec{}

Recall the equivalence 
$$\IndCoh^*(\Op_\cG^\mf)\underset{\Sph^{\on{spec}}_\cG}\otimes \on{I}(\cG,\cP^-)^{\on{spec,loc}}\simeq \IndCoh^*((\Op^\Mmf_{\cG,\cP^-})^\wedge_\mf)$$
of \eqref{e:Op semiinf * bis}, and recall the adjunction 
$$({}'\!\iota)^\IndCoh_*:
\IndCoh^*(\Op^\Pmf_{\cG,\cP^-})\rightleftarrows 
\IndCoh^*((\Op^\Mmf_{\cG,\cP^-})^\wedge_\mf):({}'\!\iota)^!$$
of \eqref{e:iota wedge *}.

\medskip

As in \propref{p:t semiinf spec} one shows that the monad $({}'\!\iota)^!\circ ({}'\!\iota)^\IndCoh_*$ acting on $\IndCoh^*(\Op^\Pmf_{\cG,\cP^-})$
is left t-exact. Hence, we obtain that
\begin{multline*} 
\IndCoh^*((\Op^\Mmf_{\cG,\cP^-})^\wedge_\mf) \overset{\eqref{e:Op semminf *}}\simeq \\
\simeq \IndCoh^*(\on{pt}\underset{\LS^\reg_\cG}\times\Op_\cG^\mf)\underset{\IndCoh^*(\fL_\nabla(\cG))}\otimes \IndCoh^*(\fL_\nabla(\cG/\cN^-_P))_{\fL^+_\nabla(\cM)}
\end{multline*} 
acquires a t-structure, uniquely characterized by the condition that the functor $({}'\!\iota)^\IndCoh_*$ is left t-exact.

\sssec{} \label{sss:iota' Op}

By construction, an object of $\IndCoh^*((\Op^\Mmf_{\cG,\cP^-})^\wedge_\mf)$ is coconnective if and only if its image under 
$({}'\!\iota)^!$ is coconnective. 

\ssec{The left t-exactness of the co-Jacquet functor}

In this subsection we show that the functor $\on{co}\!J^{-,\on{spec},\on{enh}}_\Op$ has a bounded cohomological amplitude on the left. In fact, we will show
that it is left t-exact.

\sssec{}

First, we note that since $\MOp^\Pmf_{\cG,\cP^-}$ and $(\MOp^\Mmf_{\cG,\cP^-})^\wedge_\mf)$ are (factorization) ind-placid ind-schemes, the 
categories 
$$\IndCoh^*(\MOp^\Pmf_{\cG,\cP^-}) \text{ and } \IndCoh^*((\MOp^\Mmf_{\cG,\cP^-})^\wedge_\mf)$$
have naturally defined t-structures, so that the functor $({}'\!\iota)^\IndCoh_*$ of \eqref{e:iota' MOp} is t-exact, and its right adjoint
$({}'\!\iota)^!$ is left t-exact, see \cite[Sects. A.8 and B.13.22]{GLC2}.

\sssec{}

Recall that the functor $\on{co}\!J^{-,\on{spec},\on{enh}}_\Op$ is the composition \eqref{e:defn coJ semiinf}, where the first arrow, i.e.,  
$((\sfq^{\on{Miu,mon-free}})^\wedge)^!$ is left t-exact. So it remains to show that the second arrow, i.e., 
\begin{equation} \label{e:another j}
\IndCoh^*((\MOp^\Mmf_{\cG,\cP^-})^\wedge_\mf)\overset{\jmath^\IndCoh_*}\to \IndCoh^*((\Op^\Mmf_{\cG,\cP^-})^\wedge_\mf)
\end{equation}
is also left t-exact. 

\medskip

In order to show that \eqref{e:another j} is left t-exact, by \secref{sss:iota' Op} and because the natural transformation \eqref{e:BC Miura *} is an 
isomorphism (by \lemref{l:BC Miura}), it suffices to show that the functor
$$\IndCoh^*(\MOp^\Pmf_{\cG,\cP^-}) \overset{\jmath^\IndCoh_*}\longrightarrow \IndCoh^*(\Op^\Pmf_{\cG,\cP^-})$$
is left t-exact. 

\medskip

We will show that this functor is actually t-exact.

\sssec{}

We have a commutative diagram
$$
\CD 
\IndCoh^*(\on{pt}\underset{\LS^\mer_{\cP^-}}\times \MOp^\mer_{\cG,\cP^-}) @<<< \IndCoh^*(\MOp^\Pmf_{\cG,\cP^-})  \\
@V{\jmath^\IndCoh_*}VV @VV{\jmath^\IndCoh_*}V \\
\IndCoh^*(\on{pt}\underset{\LS^\mer_{\cP^-}}\times \Op^\mer_{\cG,\cP^-}) @<<< \IndCoh^*(\Op^\Pmf_{\cG,\cP^-}),
\endCD
$$
in which the horizontal arrows are given by *-pullback. These horizontal arrows are t-exact and conservative. 

\medskip

Hence, it is enough to show that the functor
$$\IndCoh^*(\on{pt}\underset{\LS^\mer_{\cP^-}}\times \MOp^\mer_{\cG,\cP^-}) \overset{\jmath^\IndCoh_*}\longrightarrow 
\IndCoh^*(\on{pt}\underset{\LS^\mer_{\cP^-}}\times \Op^\mer_{\cG,\cP^-})$$
is t-exact.

\medskip

However, this holds by \cite[Corollary A.8.8]{GLC2}, since 
$$\jmath:\on{pt}\underset{\LS^\mer_{\cP^-}}\times \MOp^\mer_{\cG,\cP^-}\to \on{pt}\underset{\LS^\mer_{\cP^-}}\times \Op^\mer_{\cG,\cP^-}$$
is a morphism between (factorization) ind-affine ind-schemes.

\begin{rem}  \label{r:coJ  exact}
The above argument shows that the restriction of the functor $\on{co}\!J^{-,\on{spec},\on{enh}}_\Op$ to
$$\IndCoh^*((\MOp^\Mmf_{\cG,\cP^-})^\wedge_\mf)\subset \IndCoh^*(\MOp^\Mmf_{\cG,\cP^-})$$
is actually t-exact.

\medskip
 
Indeed, the connective subcategory of
$\IndCoh^*((\MOp^\Mmf_{\cG,\cP^-})^\wedge_\mf)$ is generated under colimits by the essential image of the connective subcategory of
$\IndCoh^*(\MOp^\Pmf_{\cG,\cP^-})$ along the functor
$$({}'\!\iota)^\IndCoh_*:\IndCoh^*(\MOp^\Pmf_{\cG,\cP^-})\to \IndCoh^*((\MOp^\Mmf_{\cG,\cP^-})^\wedge_\mf).$$
\end{rem} 

\ssec{Fully faithfulness on the spectral side}

In order to complete the steps in \secref{sss:t-str properties}, it remains to show that the functor $\sF^{\on{spec}}_{\Op(\cG)}$ is fully faithful at the pointwise level,
and also when restricted to the eventually coconnective subcategory in the factorization setting. 

\medskip

We will prove this in the present subsection. 

\sssec{}

Let $R_\cG$ be the regular representation of $\cG$, viewed as a commutative algebra, and hence also as a factorization algebra, in 
$\Rep(\cG)\otimes \Rep(\cG)$. 

\medskip

By a slight abuse of notation, we will denote by the same symbol the image of $R_\cG$ under the restriction functor
$$\Rep(\cG)\otimes \Rep(\cG)\to \Rep(\cG)\otimes \Rep(\cP^-).$$

\medskip

Let $R_\cG^\Op$ denote the image of $R_\cG$ under the functors
$$\Rep(\cG)\otimes \Rep(\cG)\overset{(\fr^\reg)^*}\to \QCoh(\Op^\reg_\cG)\otimes \Rep(\cG) \text{ or }
\Rep(\cG)\otimes \Rep(\cP^-)\overset{(\fr^\reg)^*}\to \QCoh(\Op^\reg_\cG)\otimes \Rep(\cP^-),$$
respectively. 

\begin{rem}

The object 
$$R_\cG^\Op\in \on{FactAlg}^{\on{untl}}(\QCoh(\Op^\reg_\cG)\otimes \Rep(\cG))$$
introduced above is closely related to the object
$$R_{\cG,\Op}\in \on{FactAlg}^{\on{untl}}(\Rep(\cG))$$
from \cite[Equation (4.8)]{GLC2}. 

\medskip

Namely,
$$(\Gamma(\Op^\reg_\cG,-)\otimes \on{Id})(R_\cG^\Op)\simeq R _{\cG,\Op}$$
as objects of $\on{FactAlg}^{\on{untl}}(\Rep(\cG))$. 

\end{rem} 

\sssec{}

Using the functor
$$\iota^\IndCoh_*:\QCoh(\Op^\reg_\cG)\to \IndCoh^*(\Op^\mer_\cG),$$
we can view $R_\cG^\Op$ as a factorization algebra in 
$$\QCoh(\Op^\mf_\cG)\otimes \Rep(\cG)  \text{ or } \QCoh(\Op^\mf_\cG)\otimes \Rep(\cP^-).$$

\medskip

Consider the resulting (lax) factorization categories
$$R_\cG^\Op\mod^{\on{fact}}\left(\IndCoh^*(\Op^\mer_\cG)\otimes \Rep(\cG)\right) \text{ and }
R_\cG^\Op\mod^{\on{fact}}\left(\IndCoh^*(\Op^\mer_\cG)\otimes \Rep(\cP^-)\right),$$
and also their full subcategories
$$R_\cG^\Op\mod^{\on{fact}}\left(\IndCoh^*((\Op^\mer_\cG)^\wedge_\mf)\otimes \Rep(\cG)\right)$$
and
$$R_\cG^\Op\mod^{\on{fact}}\left(\IndCoh^*((\Op^\mer_\cG)^\wedge_\mf)\otimes \Rep(\cP^-)\right)$$
respectively. 

\sssec{}

The operation of $(\IndCoh,*)$-pushforward along
$$\iota^\mf\times \fr:\Op^\mf_\cG\to \Op^\mer_\cG\times \LS^\reg_\cG$$
gives rise to a (lax unital) factorization functor
$$\IndCoh^*(\Op^\mf_\cG)\to \IndCoh^*(\Op^\mer_\cG)\otimes \Rep(\cG).$$

This functor sends the factorization unit, i.e., $\CO_{\Op^\reg_\cG}\in \IndCoh^*(\Op^\mf_\cG)$ to $R_\cG^\Op$. 
Hence, it upgrades to a functor
\begin{equation} \label{e:IndCoh* Op mf via R}
\IndCoh^*(\Op^\mf_\cG) \to R_\cG^\Op\mod^{\on{fact}}\left(\IndCoh^*(\Op^\mer_\cG)\otimes \Rep(\cG)\right),
\end{equation} 
see \cite[Sect. 4.1.6]{GLC2}.

\medskip

Tensoring over $\Rep(\cG)$ with $\Rep(\cP^-)$, we obtain a functor 
\begin{equation} \label{e:IndCoh* Op mf via R P}
\IndCoh^*(\Op^\Pmf_{\cG,\cP^-}) \to R_\cG^\Op\mod^{\on{fact}}\left(\IndCoh^*(\Op^\mer_\cG)\otimes \Rep(\cP^-)\right).
\end{equation} 

\sssec{}

We claim:

\begin{prop} \label{p:IndCoh* Op mf via R} \hfill

\medskip

\noindent{\em(a)} The functors \eqref{e:IndCoh* Op mf via R} and \eqref{e:IndCoh* Op mf via R P} are fully faithful 
at the pointwise level.

\medskip

\noindent{\em(b)} In the factorization setting, the functors \eqref{e:IndCoh* Op mf via R} and \eqref{e:IndCoh* Op mf via R P} induce equivalences 
between the eventually coconnective subcategories of the two sides. In particular, the restrictions of the functors \eqref{e:IndCoh* Op mf via R} and \eqref{e:IndCoh* Op mf via R P}
to the eventually coconnective subcategory of the source is fully faithful. 

\end{prop}

\begin{proof}

We first prove point (a). We will show that at the pointwise level,
the functor \eqref{e:IndCoh* Op mf via R} gives rise to an equivalence 
\begin{equation} \label{e:IndCoh* Op mf via R 1}
\IndCoh^*(\Op^\mf_\cG) \to R_\cG^\Op\mod^{\on{fact}}\left(\IndCoh^*((\Op^\mer_\cG)^\wedge_\mf)\otimes \Rep(\cG)\right).
\end{equation}

The assertion for  \eqref{e:IndCoh* Op mf via R P} would then follow by base change.

\medskip

Recall the functor 
\begin{equation} \label{e:IndCoh* Op mf via R 2}
\Sph_\cG^{\on{spec}}\to R_\cG\mod^{\on{fact}}(\Rep(\cG)\otimes \Rep(\cG))
\end{equation}
of \eqref{e:Omega R}. This functor is linear with respect to $\Sph_\cG^{\on{spec}}\otimes \Sph_\cG^{\on{spec}}$, where:

\begin{itemize}

\item $\Sph_\cG^{\on{spec}}\otimes \Sph_\cG^{\on{spec}}$ acts on $\Sph_\cG^{\on{spec}}$ naturally;

\item $\Sph_\cG^{\on{spec}}\otimes \Sph_\cG^{\on{spec}}$ acts on $R_\cG\mod^{\on{fact}}(\Rep(\cG)\otimes \Rep(\cG))$
by acting on $\Rep(\cG)\otimes \Rep(\cG)$ naturally, see \secref{sss:universal action 1}.

\end{itemize} 

\medskip

Note that the functor \eqref{e:IndCoh* Op mf via R 1} factors as
\begin{multline*}
\IndCoh^*(\Op^\mf_\cG)=
\IndCoh^*(\Op^\mf_\cG)\underset{\Sph_\cG^{\on{spec}}}\otimes \Sph_\cG^{\on{spec}}
\overset{\on{Id}\otimes \text{\eqref{e:IndCoh* Op mf via R 2}}}\longrightarrow \\
\to \IndCoh^*(\Op^\mf_\cG)\underset{\Sph_\cG^{\on{spec}}}\otimes 
\left(R_\cG\mod^{\on{fact}}(\Rep(\cG)\otimes \Rep(\cG))\right) \overset{\sim}\to \\
\overset{\sim}\to  R_\cG\mod^{\on{fact}}\left(\left(\IndCoh^*(\Op^\mf_\cG)\underset{\Sph_\cG^{\on{spec}}}\otimes \Rep(\cG)\right)\otimes \Rep(\cG)\right)\simeq \\
\simeq R_\cG^\Op\mod^{\on{fact}}\left(\IndCoh^*((\Op^\mer_\cG)^\wedge_\mf)\otimes \Rep(\cG)\right),
\end{multline*}
where:

\smallskip

\begin{itemize}

\item The third arrow is an equivalence due to the fact that $\IndCoh^*(\Op^\mf_\cG)$ is dualizable and $\Sph_\cG^{\on{spec}}$ is rigid;

\medskip

\item The last equivalence is \cite[Proposition 3.6.5]{GLC2}.
\end{itemize} 

Hence, in order to prove that \eqref{e:IndCoh* Op mf via R 1} is an equivalence at the pointwise level, it suffices to show that the functor
\begin{multline} \label{e:IndCoh* Op mf via R 3}
\IndCoh^*(\Op^\mf_\cG)\underset{\Sph_\cG^{\on{spec}}}\otimes \Sph_\cG^{\on{spec}}
\overset{\on{Id}\otimes \text{\eqref{e:IndCoh* Op mf via R 2}}}\longrightarrow \\
\to \IndCoh^*(\Op^\mf_\cG)\underset{\Sph_\cG^{\on{spec}}}\otimes 
\left(R_\cG\mod^{\on{fact}}(\Rep(\cG)\otimes \Rep(\cG))\right)
\end{multline}
has this property.

\medskip

Fix $\ul{x}\in \Ran$. Note that at the level of fibers at $\ul{x}$, the functor \eqref{e:IndCoh* Op mf via R 2}
factors as
$$\Sph_{\cG,\ul{x}}^{\on{spec}}\to \Sph^{\on{spec}}_{\cG,\on{temp},\ul{x}} \to R_\cG\mod^{\on{fact}}(\Rep(\cG)\otimes \Rep(\cG))_{\ul{x}},$$
where the second arrow is an equivalence. Indeed, both
$$\Sph^{\on{spec}}_{\cG,\on{temp},\ul{x}} \text{ and }  R_\cG\mod^{\on{fact}}(\Rep(\cG)\otimes \Rep(\cG))_{\ul{x}}$$
identify with the left completion of $\Sph_{\cG,\ul{x}}^{\on{spec}}$ with respect to the natural t-structure. 

\medskip

Hence, in order to prove that \eqref{e:IndCoh* Op mf via R 3} is a pointwise equivalence, it suffices to show that 
\begin{multline*}
\IndCoh^*(\Op^\mf_\cG)_{\ul{x}}=\IndCoh^*(\Op^\mf_\cG)_{\ul{x}}\underset{\Sph_{\cG,\ul{x}}^{\on{spec}}}\otimes \Sph_{\cG,\ul{x}}^{\on{spec}}\to \\
\to \IndCoh^*(\Op^\mf_\cG)_{\ul{x}}\underset{\Sph_{\cG,\ul{x}}^{\on{spec}}}\otimes \Sph_{\cG,\on{temp},\ul{x}}^{\on{spec}}
\end{multline*}
is an equivalence. 

\medskip

However, this is the content of \cite[Proposition 7.2.4]{GLC2}. 

\medskip

We now prove point (b). Again, by base change, it is enough to prove the assertion about \eqref{e:IndCoh* Op mf via R}. 
Let $\CO_{\Op^\reg_\cG}$ be the factorization algebra in $\Vect$ from \cite[Sect.4.4.1]{GLC2}. 

\medskip

By \cite[Corollary 4.4.2]{GLC2}, the functor
$$(\Gamma^\IndCoh(\Op^\mer_\cG,-)\otimes \on{Id}):\IndCoh^*(\Op^\mer_\cG)\otimes \Rep(\cG)\to
(\CO_{\Op^\reg_\cG}\otimes \one_{\Rep(\cG)})\mod^{\on{fact}}(\Rep(\cG))$$
induces an equivalence on eventually coconnective subcategories.

\medskip

Note that we have a tautological map
$$\CO_{\Op^\reg_\cG}\otimes \one_{\Rep(\cG)}\to (\Gamma^\IndCoh(\Op^\mer_\cG,-)\otimes \on{Id})(R^\Op_\cG)$$
as factorization algebras in $\Rep(\cG)$. Hence, the functor $\Gamma^\IndCoh(\Op^\mer_\cG,-)\otimes \on{Id}$ 
induces an equivalence between the eventually coconnective subcategory of $R_\cG^\Op\mod^{\on{fact}}\left(\IndCoh^*(\Op^\mer_\cG)\otimes \Rep(\cG)\right)$
and the eventually coconnective subcategory of 
$$(\Gamma^\IndCoh(\Op^\mer_\cG,-)\otimes \on{Id})(R^\Op_\cG)\mod^{\on{fact}}(\Rep(\cG)).$$

Note, however, that we have a tautological identification
$$(\Gamma^\IndCoh(\Op^\mer_\cG,-)\otimes \on{Id})(R^\Op_\cG)\simeq R_{\cG,\Op},$$
where $R_{\cG,\Op}$ is as in \cite[Equation (4.8)]{GLC2}. 

\medskip

With respect to this identification, the composite functor
\begin{multline*}
\IndCoh^*(\Op^\mf_\cG)  \overset{\text{\eqref{e:IndCoh* Op mf via R}}}\longrightarrow \\
\to R_\cG^\Op\mod^{\on{fact}}\left(\IndCoh^*(\Op^\mer_\cG)\otimes \Rep(\cG)\right) 
\overset{\Gamma^\IndCoh(\Op^\mer_\cG,-)\otimes \on{Id}}\longrightarrow R_{\cG,\Op}\mod^{\on{fact}}(\Rep(\cG))
\end{multline*} 
is the functor of \cite[Equation (4.10)]{GLC2}. 

\medskip

The required assertion follows now by \cite[Proposition 4.4.7]{GLC2}. 

\end{proof} 

\sssec{}

Note that the functor
$$\on{inv}_{\cN^-_\cP}:\Rep(\cP^-)\to \Rep(\cM)$$
induces a functor
$$\Rep(\cP^-)\to \Omega^{\on{spec}}\mod^{\on{fact}}(\Rep(\cM)).$$

From here we obtain a functor 
\begin{multline} \label{e:iota fact mod}
R_\cG^\Op\mod^{\on{fact}}\left(\IndCoh^*(\Op^\mer_\cG)\otimes \Rep(\cP^-)\right)\to  \\
\to \Omega(R)^\Op\mod^{\on{fact}}\left(\IndCoh^*(\Op^\mer_\cG)\otimes \Rep(\cM)\right).
\end{multline}

Unwinding, we obtain a commutative diagram

\begin{equation} \label{e:iota diag Op}
\CD
\IndCoh^*(\Op^\Pmf_{\cG,\cP^-}) @>{({}'\!\iota)^\IndCoh_*}>>  \IndCoh^*((\Op^\Mmf_{\cG,\cP^-})^\wedge_\mf) \\
@V{\text{\eqref{e:IndCoh* Op mf via R P}}}VV @VV{\sF^{\on{spec}}_{\Op(\cG)}}V \\
R_\cG^\Op\mod^{\on{fact}}\left(\IndCoh^*(\Op^\mer_\cG)\otimes \Rep(\cP^-)\right) @>{\text{\eqref{e:iota fact mod}}}>> 
\Omega(R)^\Op\mod^{\on{fact}}\left(\IndCoh^*(\Op^\mer_\cG)\otimes \Rep(\cM)\right).
\endCD
\end{equation}

\medskip

As in \propref{p:iota diag} one shows:

\begin{prop} \label{p:iota diag Op}
The functor \eqref{e:iota fact mod} admits a continuous right adjoint. Moreover, the natural transformation in the diagram
$$
\xy
(0,0)*+{\IndCoh^*(\Op^\Pmf_{\cG,\cP^-})}="A";
(75,0)*+{\IndCoh^*((\Op^\Mmf_{\cG,\cP^-})^\wedge_\mf)}="B";
(0,-30)*+{R_\cG^\Op\mod^{\on{fact}}\left(\IndCoh^*(\Op^\mer_\cG)\otimes \Rep(\cP^-)\right)}="C";
(75,-30)*+{\Omega(R)^\Op\mod^{\on{fact}}\left(\IndCoh^*(\Op^\mer_\cG)\otimes \Rep(\cM)\right),}="D";
{\ar@{->}_{\text{\eqref{e:IndCoh* Op mf via R P}}} "A";"C"};
{\ar@{->}_{({}'\!\iota)^!} "B";"A"};
{\ar@{->}^{\sF^{\on{spec}}_{\Op(\cG)}} "B";"D"};
{\ar@{->} "D";"C"};
{\ar@{=>} "A";"D"};
\endxy
$$
obtained by passing to right adjoints along the horizontal arrows in \eqref{e:iota diag Op}, is an isomorphism. 
\end{prop}

\sssec{}

We are now ready to prove the desired properties of the functor $\sF^{\on{spec}}_{\Op(\cG)}$. The adjunction \eqref{e:iota wedge *}
is monadic and induces a monadic adjunction on the eventually coconnective categories. 

\medskip

Hence, in order to prove that $\sF^{\on{spec}}_{\Op(\cG)}$
is fully faithful (``as-is" at the pointwise level, and on the eventually coconnective subcategories in the factorization context), it
suffices to prove that the monad corresponding to the functor 
$$({}'\!\iota)^\IndCoh_*:\IndCoh^*(\Op^\Pmf_{\cG,\cP^-})\to \IndCoh^*((\Op^\Mmf_{\cG,\cP^-})^\wedge_\mf)$$ maps isomorphically to the monad corresponding
to the composition
$$\sF^{\on{spec}}_{\Op(\cG)}\circ ({}'\!\iota)^\IndCoh_*.$$

However, this follows by combining \propref{p:IndCoh* Op mf via R} and \propref{p:iota diag Op}. 

\sssec{}

This concludes the proof of the properties announced in \secref{sss:strategy Jacquet dual enh}, and hence of \thmref{t:local Jacquet dual enh}. 

\newpage 

\centerline{\bf Part II: Eisenstein and constant term functors}

\bigskip

In Part II of the paper \cite{GLC2} we introduced the basic local-to-global constructions that play a role 
in the geometric Langlands theory.

\medskip 

On the geometric side these were the functors of Whittaker coefficient $$\on{coeff}_G:\Dmod_{\frac{1}{2}}(\Bun_G)\to \Whit^!(G)_\Ran,$$
and the localization functor $$\Loc_G:\KL(G)_{\crit,\Ran}\to \Dmod_{\frac{1}{2}}(\Bun_G).$$

\medskip

On the spectral side these were the spectral localization functor 
$$\Loc_\cG^{\on{spec}}:\Rep(\cG)_\Ran\to \IndCoh(\LS_\cG)$$ and the spectral Poincar\'e functor
$$\on{Poinc}^{\on{spec}}_\cG:\IndCoh(\Op^\mf_\cG)_\Ran\to \IndCoh(\LS_\cG).$$

\medskip

The main theme of the current paper is the interaction between the global (geometric and spectral) categories for $G$ with their
counterparts for Levi subgroups. We introduce the basic functors that connect these categories. 

\medskip

On the geometric side, this is the adjoint pair
$$\Eis^-_!:\Dmod_{\frac{1}{2}}(\Bun_M)\rightleftarrows \Dmod_{\frac{1}{2}}(\Bun_G):\on{CT}^-_*$$
and on the spectral side, this is the adjoint pair
$$\Eis^{-,\on{spec}}:\IndCoh(\LS_\cM) \rightleftarrows \IndCoh(\LS_\cG):\on{CT}^{-,\on{spec}}.$$

\medskip

The goal of Part II of the current paper is to study the compositions of these functors with the 
local-to-global functors mentioned above. In each case, we will express a given 
composition as a local-to-global functor for the Levi and a functor of local nature.

\medskip

The ultimate goal of these calculations is to show that the Langlands functor is compatible with
the $\Eis$ and $\on{CT}$ functors, as will be established in Part III of the present paper.

%

\section{Eisenstein and constant term functors on the geometric side} \label{s:Eis}

In this section we develop the functors that connect the automorphic category for $G$
with the automorphic categories for its Levi subgroups. 

\medskip

We start by introducing the most fundamental pair of mutually adjoint functors
$$\Eis^-_!:\Dmod_{\frac{1}{2}}(\Bun_M)\rightleftarrows \Dmod_{\frac{1}{2}}(\Bun_G):\on{CT}^-_*.$$

These functors are a direct geometric counterpart of the operators of Eisenstein series and constant
term in the theory of automorphic functions. 

\medskip

However, by contrast with the classical theory, in the geometric context, we can enhance these functors,
by (locally) tensoring up with $\on{I}(G,P^-)^{\on{loc}}$ (or its dual). This way, we obtain various versions
of enhanced Eisenstein and constant term functors. These enhanced functors will be used in the proof of 
the main theorem, \thmref{t:main}, in Part III of the paper. 

\ssec{The usual Eisenstein and constant term functors}

\sssec{}

Consider the diagram 
\begin{equation} \label{e:Eis CT diag}
\Bun_G \overset{\sfp^-}\leftarrow \Bun_{P^-} \overset{\sfq^-}\to \Bun_M. 
\end{equation} 

The functors of Eisenstein and constant term are defined by pull-push along the above diagram,
followed by various corrections explained in this subsection. 

\sssec{} \label{sss:square root global gerbe}

Consider the restrictions of the line bundles 
$$\det_{\Bun_G} \text{ and } \det_{\Bun_M}$$
along the maps in \eqref{e:Eis CT diag}. 

\medskip

Denote by $\det_{\Bun_{G,M}}$ their ratio:
$$\det_{\Bun_{G,M}}=\det_{\Bun_G}|_{\Bun_{P^-}}\otimes (\det_{\Bun_M}|_{\Bun_{P^-}})^{\otimes -1};$$
it naturally descends to $\Bun_M$. By a slight abuse of notation,
we will denote the resulting line bundle on $\Bun_M$ by the same symbol $\det_{\Bun_{G,M}}$. 

\medskip

According to \cite[Proposition 1.3.3]{GLC1}, the line bundle $\det_{\Bun_{G,M}}$ admits a canonical square root,
to be denoted $\det^{\otimes \frac{1}{2}}_{\Bun_{G,M}}$ (see \secref{sss:square root expl} for the explicit formula). 

\medskip

In particular, the pullback of the $\mu_2$-gerbes 
$$\det^{\frac{1}{2}}_{\Bun_G} \text{ and } \det^{\frac{1}{2}}_{\Bun_M}$$
along the maps \eqref{e:Eis CT diag} are canonically identified.

\sssec{}

We will denote the resulting category of twisted D-modules by
$$\Dmod_{\frac{1}{2}}(\Bun_{P^-}).$$

The morphisms \eqref{e:Eis CT diag} give rise to functors
$$(\sfp^-)^!:\Dmod_{\frac{1}{2}}(\Bun_G)\to \Dmod_{\frac{1}{2}}(\Bun_{P^-})$$
and
$$(\sfq^-)_*: \Dmod_{\frac{1}{2}}(\Bun_{P^-})\to  \Dmod_{\frac{1}{2}}(\Bun_M),$$
where the latter is continuous because the morphism $\sfq^-$ is \emph{safe}
(see \cite[Defn. 10.2.2]{DG1} for what this means). 

\sssec{The naive constant term functor}

We define the functor 
$$\on{CT}_*^{-,\on{nv}}:\Dmod_{\frac{1}{2}}(\Bun_G)\to \Dmod_{\frac{1}{2}}(\Bun_M)$$
to be
$$(\sfq^-)_*\circ (\sfp^-)^!.$$

\sssec{}

The morphism $\sfq^-$ is smooth. Hence, the functor
$$(\sfq^-)^*:\Dmod_{\frac{1}{2}}(\Bun_G)\to \Dmod_{\frac{1}{2}}(\Bun_{P^-}),$$
left adjoint to $(\sfq^-)_*$ is well-defined.

\medskip

According to \cite[Proposition 1.1.2]{DG3}, the (partially defined) functor $(\sfp^-)_!$, left adjoint to $(\sfp^-)^!$ is defined on the essential image
of $(\sfq^-)^!$. Hence, the functor $\on{CT}_*^{-,\on{nv}}$ admits a left adjoint, explicitly given by
$$\Eis^{-,\on{nv}}_!:=(\sfp^-)_!\circ (\sfq^-)^*.$$

\sssec{}

We define the corrected constant term and Eisenstein functors
$$\Eis^-_!:\Dmod_{\frac{1}{2}}(\Bun_G)\rightleftarrows \Dmod_{\frac{1}{2}}(\Bun_M):\on{CT}_*^-$$
by 
\begin{equation} \label{e:Eis ! defn}
\Eis^-_!:=\Eis^{-,\on{nv}}_![\on{shift}]
\end{equation}
and  
\begin{equation} \label{e:CT * defn}
\on{CT}_*^-:=\on{CT}_*^{-,\on{nv}}[-\on{shift}],
\end{equation}
where
\begin{equation} \label{e:global shift}
\on{shift}:=\on{dim.rel.}(\Bun_{P^-}/\Bun_M).
\end{equation}

\noindent NB: the integer $\on{dim.rel.}(\Bun_{P^-}/\Bun_M)$ depends on the connected component of $\Bun_M$. 

\sssec{The $\rho$-shift}

Let $\CP_{Z_M}$ be a $Z_M$-torsor on $X$. Translation by $\CP_{Z_M}$ is an 
automorphism
$$\on{transl}_{\CP_{Z_M}}:\Bun_M\to \Bun_M.$$

By the construction of the line bundle $\det_{\Bun_M}$, it is equivariant with respect to 
$\on{transl}_{\CP_{Z_M}}$. 

\medskip 

In particular, we obtain mutually inverse self-equivalences
$$(\on{transl}_{\CP_{Z_M}})_*:\Dmod_{\frac{1}{2}}(\Bun_M)\rightleftarrows \Dmod_{\frac{1}{2}}(\Bun_M):(\on{transl}_{\CP_{Z_M}})^*.$$

\medskip

We define the translated constant term and Eisenstein functors 
$$\Eis^-_{!,\CP_{Z_M}}:\Dmod_{\frac{1}{2}}(\Bun_G)\rightleftarrows \Dmod_{\frac{1}{2}}(\Bun_M):\on{CT}_{*,\CP_{Z_M}}^-$$
by
$$\Eis^-_{!,\CP_{Z_M}}:=\Eis^-_!\circ (\on{transl}_{\CP_{Z_M}})_* \text{ and }
\on{CT}_{*,\CP_{Z_M}}^-:=(\on{transl}_{\CP_{Z_M}})^*\circ \on{CT}_*^-.$$

\medskip

We will mostly apply the above construction when $\CP_{Z_M}=\rho_P(\omega_X)$. 

\sssec{The constant term functor for a general level} \label{sss:CT level}

Let $\kappa$ be a level for $\fg$, and recall (see \secref{sss:level for m}) that we use the same character $\kappa$
to denote the corresponding level for $\fm$. 

\medskip

We let
$$\on{CT}^-_*:\Dmod_{\crit_G+\kappa}(\Bun_G) \to \Dmod_{\crit_M+\kappa}(\Bun_M)$$ 
denote the composition
\begin{multline*} 
\Dmod_{\crit_G+\kappa}(\Bun_G)= \\
= \Dmod_{\frac{1}{2}\on{dlog}(\det(\Bun_G))+\kappa}(\Bun_G)
\overset{(\sfp^-)^!}\longrightarrow 
 \Dmod_{\frac{1}{2}\on{dlog}(\det(\Bun_G)|_{\Bun_{P^-}})+\kappa}(\Bun_{P^-}) \simeq \\
\simeq \Dmod_{\frac{1}{2}\on{dlog}(\det(\Bun_M)|_{\Bun_{P^-}})+\on{dlog}(\det^{\otimes \frac{1}{2}}_{\Bun_{G,M}})+\kappa}(\Bun_{P^-})
\overset{\otimes \det^{\otimes -\frac{1}{2}}_{\Bun_{G,M}}}\longrightarrow \\
\to \Dmod_{\frac{1}{2}\on{dlog}(\det(\Bun_M)|_{\Bun_{P^-}})+\kappa}(\Bun_{P^-})\overset{(\sfq^-)_*}\longrightarrow 
\Dmod_{\frac{1}{2}\on{dlog}(\det(\Bun_M))+\kappa}(\Bun_M)= \\
= \Dmod_{\crit_M+\kappa}(\Bun_M)\overset{[-\on{shift}]}\longrightarrow \Dmod_{\crit_M+\kappa}(\Bun_M).
\end{multline*} 

\ssec{The geometric (co-)enhanced category} \label{ss:geom global enh}

\sssec{}

Let $\CZ$ be a prestack mapping to $\Ran$. Consider the category 
$$\Dmod_{\frac{1}{2}}(\Bun_M)^{-,\on{enh}}_\CZ:=\left(\Dmod_{\frac{1}{2}}(\Bun_M)\otimes \Dmod(\CZ)\right)\underset{\Sph_{M,\CZ}}\otimes \on{I}(G,P^-)^{\on{loc}}_\CZ.$$

The assignment 
$$\CZ\rightsquigarrow \Dmod_{\frac{1}{2}}(\Bun_M)^{-,\on{enh}}_\CZ$$
is a unital crystal of categories over $\Ran$ (see \cite[Sect. 11.2.1]{GLC2} for the terminology), to be denoted $\ul\Dmod_{\frac{1}{2}}(\Bun_M)^{-,\on{enh}}$.

\sssec{}

The monadic adjunction
$$\ind_{\Sph\to \semiinf}:\Sph_M\rightleftarrows \on{I}(G,P^-)^{\on{loc}}:\oblv_{\semiinf\to \Sph}$$
gives rise to a monadic adjunction
\begin{equation} \label{e:enh adj glob}
\ind_{\on{enh}}:\Dmod_{\frac{1}{2}}(\Bun_M)\otimes \Dmod(\CZ) \rightleftarrows \Dmod_{\frac{1}{2}}(\Bun_M)^{-,\on{enh}}_\CZ:\oblv_{\on{enh}}.
\end{equation} 

\sssec{}

Here is another way to think about the category $\Dmod_{\frac{1}{2}}(\Bun_M)^{-,\on{enh}}_\CZ$.

\medskip

Let $\wt\Omega$ be as in \secref{sss:Omega tilde}. We can consider
$$\wt\Omega_\CZ\in \on{AssocAlg}(\Sph_{M,\CZ}),$$
where we consider $\Sph_{M,\CZ}$ as a monoidal category. 

\medskip

Consider $\Dmod_{\frac{1}{2}}(\Bun_M)\otimes \Dmod(\CZ)$ as a right module category over $\Sph_{M,\CZ}$. 
Then we have
$$\Dmod_{\frac{1}{2}}(\Bun_M)^{-,\on{enh}}_\CZ\simeq \wt\Omega_\CZ\mod^r(\Dmod_{\frac{1}{2}}(\Bun_M)\otimes \Dmod(\CZ)),$$
so that the adjunction \eqref{e:enh adj glob} corresponds to
$$\ind_{\wt\Omega_\CZ}:\Dmod_{\frac{1}{2}}(\Bun_M)\rightleftarrows 
\wt\Omega_\CZ\mod^r(\Dmod_{\frac{1}{2}}(\Bun_M)\otimes \Dmod(\CZ)):\oblv_{\wt\Omega_\CZ}.$$

\sssec{}

Let
\begin{equation} \label{e:Bun M enh indep}
\left(\Dmod_{\frac{1}{2}}(\Bun_M)^{-,\on{enh}}\right)_{\Ran^{\on{untl}},\on{indep}}
\end{equation} 
be the corresponding ``independent" category (see \cite[Sect. H.1]{GLC2}).

\medskip

Note that it fits into the pullback square
\begin{equation} \label{e:IGP indep}
\CD
\left(\Dmod_{\frac{1}{2}}(\Bun_M)^{-,\on{enh}}\right)_{\Ran^{\on{untl}},\on{indep}} @>>> (\Dmod_{\frac{1}{2}}(\Bun_M)^{-,\on{enh}}_\Ran \\
@V{\oblv_{\on{enh}}}VV @VV{\oblv_{\on{enh}}}V \\ 
\Dmod_{\frac{1}{2}}(\Bun_M)\otimes \Dmod(\Ran^{\on{untl}})_{\on{indep}}  @>>> \Dmod_{\frac{1}{2}}(\Bun_M) \otimes \Dmod(\Ran) \\
@V{\sim}VV \\
\Dmod_{\frac{1}{2}}(\Bun_M). 
\endCD
\end{equation}

\begin{rem} \label{r:global IGP untl} 

The contents of this remark will not be used in the sequel:

\medskip

One can show that the category $\left(\Dmod_{\frac{1}{2}}(\Bun_M)^{-,\on{enh}}\right)_{\Ran^{\on{untl}},\on{indep}}$ 
is equivalent to the category 
$\on{I}(G,P^-)^{\on{glob}}$ from \cite[Sect. 6.1]{Ga1}, defined as follows:

\medskip

For a group $H$, let $\Bun_{H,\on{gen}}$ denote the prestack of $H$-bundles defined generically on $X$.
There exists a tautological map
$$\Bun_H\to \Bun_{H,\on{gen}}.$$
 
Consider the prestack
$$\Bun_{P^-,\on{gen}}\underset{\Bun_{G,\on{gen}}}\times \Bun_G.$$

The category $\on{I}(G,P^-)^{\on{glob}}$ is by definition the fiber product
$$\Dmod_{\frac{1}{2}}(\Bun_{P^-,\on{gen}}\underset{\Bun_{G,\on{gen}}}\times \Bun_G)
\underset{\Dmod_{\frac{1}{2}}(\Bun_{P^-})}\times \Dmod_{\frac{1}{2}}(\Bun_M),$$
where:

\begin{itemize}

\item The twisting on $\Bun_{P^-,\on{gen}}\underset{\Bun_{G,\on{gen}}}\times \Bun_G$ is pulled back from $\Bun_G$;

\medskip

\item The functor
$$\Dmod_{\frac{1}{2}}(\Bun_M)\to \Dmod_{\frac{1}{2}}(\Bun_{P^-})$$
is $(\sfq^-)^*$;

\medskip

\item The functor $\Dmod_{\frac{1}{2}}(\Bun_{P^-,\on{gen}}\underset{\Bun_{G,\on{gen}}}\times \Bun_G)\to \Dmod_{\frac{1}{2}}(\Bun_{P^-})$
is !-pullback with respect to 
$$\Bun_{P^-}\to \Bun_{P^-,\on{gen}}\underset{\Bun_{G,\on{gen}}}\times \Bun_G.$$

\end{itemize} 

\medskip

In other words, $\on{I}(G,P^-)^{\on{glob}}$ is the full subcategory of 
$\Dmod_{\frac{1}{2}}(\Bun_{P^-,\on{gen}}\underset{\Bun_{G,\on{gen}}}\times \Bun_G)$,
consisting of objects, whose !-pullback to $\Bun_{P^-}$ lies in the essential image of the 
(fully faithful) functor $(\sfq^-)^*$:
\begin{equation} \label{e:IGP glob}
\CD
\on{I}(G,P^-)^{\on{glob}} @>>> \Dmod_{\frac{1}{2}}(\Bun_{P^-,\on{gen}}\underset{\Bun_{G,\on{gen}}}\times \Bun_G) \\
@VVV @VV{\text{!-pullback}}V  \\
\Dmod_{\frac{1}{2}}(\Bun_M) @>{(\sfq^-)^*}>> \Dmod_{\frac{1}{2}}(\Bun_{P^-,\on{gen}}).
\endCD
\end{equation}

%

\medskip

Under the equivalence
$$\left(\Dmod_{\frac{1}{2}}(\Bun_M)^{-,\on{enh}}\right)_{\Ran^{\on{untl}},\on{indep}} \simeq \on{I}(G,P^-)^{\on{glob}},$$
the left vertical arrow in \eqref{e:IGP indep} corresponds to the left vertical arrow in \eqref{e:IGP glob}. 

\end{rem}

\sssec{} \label{sss:geom co-enh}

We define a unital crystal of categories $\ul\Dmod_{\frac{1}{2}}(\Bun_G)^{-,\on{enh}_{\on{co}}}$ over $\Ran$ by 
$$\Dmod_{\frac{1}{2}}(\Bun_G)^{-,\on{enh}_{\on{co}}}_\CZ:=
\left(\Dmod_{\frac{1}{2}}(\Bun_G)\otimes \Dmod(\CZ)\right)\underset{\Sph_{G,\CZ}}\otimes \on{I}(G,P^-)^{\on{loc}}_{\on{co},\CZ}.$$

\ssec{The \emph{co}-enhanced constant term functor} \label{ss:co enh CT}

\sssec{}

Our first goal is to define the co-enhanced constant term functor 
\begin{equation} \label{e:CT co enh Shv}
\ul{\on{CT}}_*^{-,\on{enh}_{\on{co}}}:\ul\Dmod_{\frac{1}{2}}(\Bun_G)^{-,\on{enh}_{\on{co}}}\to
\Dmod_{\frac{1}{2}}(\Bun_M)\otimes \ul\Dmod(\Ran),
\end{equation}
which fits into the formalism of ``local-to-global" functors of \cite[Sect. 11.1]{GLC2}.

\medskip

I.e. we wish to define a family of functors 
\begin{equation} \label{e:CT co enh Z}
\on{CT}_{*,\CZ}^{-,\on{enh}_{\on{co}}}:\Dmod_{\frac{1}{2}}(\Bun_G)^{-,\on{enh}_{\on{co}}}_\CZ\to
\Dmod_{\frac{1}{2}}(\Bun_M)\otimes \Dmod(\CZ), \quad \CZ\to \Ran. 
\end{equation} 

\medskip

Following the convention of \cite[Sect. 11.1.4]{GLC2}, we will denote by $\on{CT}_*^{-,\on{enh}_{\on{co}}}$ the resulting functor
$$\Dmod_{\frac{1}{2}}(\Bun_G)^{-,\on{enh}_{\on{co}}}_\Ran\to \Dmod_{\frac{1}{2}}(\Bun_M).$$

\sssec{}

In order to simplify the notation, we will consider the case when $\CZ=\on{pt}$, and $\CZ\to \Ran$ corresponds to a point $\ul{x}\in \Ran$.
Thus, we will define a functor
\begin{equation} \label{e:CT co enh pt}
\on{CT}_{*,\ul{x}}^{-,\on{enh}_{\on{co}}}:\Dmod_{\frac{1}{2}}(\Bun_G)^{-,\on{enh}_{\on{co}}}_{\ul{x}}\to \Dmod_{\frac{1}{2}}(\Bun_M).
\end{equation} 

\sssec{}

Consider the Hecke global stack
$$\Bun_G \overset{\hl_G}\leftarrow \on{Hecke}^{\on{glob}}_{G,\ul{x}}\overset{\hr_G}\to \Bun_G.$$

Denote by $\sfs_{\ul{x}}$ the projection
$$\on{Hecke}^{\on{glob}}_{G,\ul{x}}\to \on{Hecke}^{\on{loc}}_{G,\ul{x}}:=\fL^+(G)_{\ul{x}}\backslash \Gr_{G,\ul{x}}.$$

Note that we have a canonical identification of line bundles on $\on{Hecke}^{\on{glob}}_{G,\ul{x}}$
$$\hl_G^*(\det_{\Bun_G})\simeq \hr_G^*(\det_{\Bun_G})\otimes \sfs^*(\det_{\Gr_{G,\ul{x}}}).$$

In particular, we have a canonical identification of $\mu_2$-gerbes
\begin{equation} \label{e:gerbes on Hecke}
\hl_G^*(\det^{\frac{1}{2}}_{\Bun_G})\simeq \hr_G^*(\det^{\frac{1}{2}}_{\Bun_G})\otimes \sfs^*(\det^{\frac{1}{2}}_{\Gr_{G,\ul{x}}}).
\end{equation} 

\sssec{}

Consider the fiber product 
\begin{equation} \label{e:base change Hecke}
\on{Hecke}^{\on{glob}}_{G,P^-,\ul{x}}:=\on{Hecke}^{\on{glob}}_{G,\ul{x}}\underset{\hr_G,\Bun_G}\times \Bun_{P^-}.
\end{equation} 

Denote the resulting maps by 
$$\Bun_G \overset{\hl_{G,P^-}}\leftarrow \on{Hecke}^{\on{glob}}_{G,P^-,\ul{x}} \overset{\hr_{G,P^-}}\to \Bun_{P^-}.$$

We have a natural map 
\begin{equation} \label{e:map s x}
\on{Hecke}^{\on{glob}}_{G,P^-,\ul{x}}  \to \fL^+(P^-)_{\ul{x}}\backslash \Gr_{G,\ul{x}},
\end{equation}
which we denote by the same symbol $\sfs_{\ul{x}}$. 

\medskip

By \secref{sss:square root global gerbe}, the identification \eqref{e:gerbes on Hecke} gives rise to an identification
of $\mu_2$-gerbes
\begin{equation} \label{e:ident gerbes Hecke P}
\hl_{G,P^-}^*(\det^{\frac{1}{2}}_{\Bun_G})\otimes \sfs^*(\det^{\frac{1}{2}}_{\Gr_{G,\ul{x}}}) \simeq 
(\hr_{G,P^-})^* \circ (\sfq^-)^*(\det^{\frac{1}{2}}_{\Bun_M}).
\end{equation} 

\sssec{}

Note that we can think of $\on{Hecke}^{\on{glob}}_{G,P^-,\ul{x}}$ as a twisted product
$$\Gr_{G,\ul{x}}\wt\times \Bun_{P^-}.$$

We will denote by $\sfs_{\ul{x}}^!$ the naturally defined functor 
$$\Dmod(\Gr_{G,\ul{x}})^{\fL^+(P^-)_{\ul{x}}}\to \Dmod(\on{Hecke}^{\on{glob}}_{G,P^-,\ul{x}}),\quad \CF\mapsto \CF\wt\boxtimes \omega_{\Bun_{P^-}}.$$

We will use the same symbol $\sfs_{\ul{x}}^!$ for the corresponding functor in the context of twisted D-modules. 

\sssec{} 

Thus, we have a well-defined functor 
\begin{equation} \label{e:pre pre CT enh}
(\hr_{G,P^-})_*\left((\hl_{G,P^-})^!(-)\sotimes \sfs_{\ul{x}}^!(-)\right),
\end{equation} 
which maps 
$$\Dmod_{\frac{1}{2}}(\Bun_G)\otimes \Dmod_{\frac{1}{2}}(\Gr_{G,\ul{x}})^{\fL^+(P^-)_{\ul{x}}}
\to 
\Dmod_{(\sfq^-)^*(\det^{\frac{1}{2}}_{\Bun_M})}(\Bun_{P^-}).$$

\sssec{}

Consider the functor 
\begin{equation} \label{e:pre CT enh}
\Dmod_{\frac{1}{2}}(\Bun_G)\otimes \Dmod_{\frac{1}{2}}(\Gr_{G,\ul{x}})^{\fL^+(P^-)_{\ul{x}}}\to 
\Dmod_{\frac{1}{2}}(\Bun_M)
\end{equation} 
equal to the composition of \eqref{e:pre pre CT enh} with:

\smallskip

\begin{itemize}

\item The functor 
$$\sfq^-_*:\Dmod_{(\fq^-)^*(\det^{\frac{1}{2}}_{\Bun_M})}(\Bun_{P^-})\to
\Dmod_{\det^{\frac{1}{2}}_{\Bun_M}}(\Bun_M)=:\Dmod_{\frac{1}{2}}(\Bun_M);$$

\item The functor of cohomological shift to the right by \eqref{e:global shift}
over a given connected component of $\Bun_M$. 

\end{itemize} 

\sssec{}

We claim:

\begin{lem}  \label{l:approx prel}
The functor 
\eqref{e:pre CT enh} factors via the quotient
$$\Dmod_{\frac{1}{2}}(\Gr_{G,\ul{x}})^{\fL^+(P^-)_{\ul{x}}}\simeq 
\Dmod_{\frac{1}{2}}(\Gr_{G,\ul{x}})_{\fL^+(N^-_P)_{\ul{x}}}^{\fL^+(M)_{\ul{x}}}\twoheadrightarrow 
\Dmod_{\frac{1}{2}}(\Gr_{G,\ul{x}})_{\fL(N^-_P)_{\ul{x}}}^{\fL^+(M)_{\ul{x}}}$$
along the second factor.
\end{lem} 

The proof will be given in \secref{ss:proof of approx prel}. 

\sssec{}

Composing with 
$$\on{I}(G,P^-)_{\on{co},\ul{x}}^{\on{loc}}:=\Dmod_{\frac{1}{2}}(\Gr_{G,\ul{x}})_{\fL(N^-_P)_{\ul{x}}}^{\fL^+(M)_{\ul{x}},\on{ren}}\to 
\Dmod_{\frac{1}{2}}(\Gr_{G,\ul{x}})_{\fL(N^-_P)_{\ul{x}}}^{\fL^+(M)_{\ul{x}}},$$
we obtain a functor 
\begin{equation} \label{e:pre CT enh 1}
\Dmod_{\frac{1}{2}}(\Bun_G)\otimes \on{I}(G,P^-)_{\on{co},\ul{x}}^{\on{loc}}\to 
\Dmod_{\frac{1}{2}}(\Bun_M),
\end{equation} 
which we denote by $\on{pre-CT}_{*,\ul{x}}^{-,\on{enh}_{\on{co}}}$

\sssec{}

Unwinding the construction, we obtain that the functor $\on{pre-CT}_{*,\ul{x}}^{-,\on{enh}_{\on{co}}}$ factors as 
\begin{equation} \label{e:pre CT enh 2}
\Dmod_{\frac{1}{2}}(\Bun_G)\otimes \on{I}(G,P^-)_{\on{co},\ul{x}}^{\on{loc}} 
\to  \Dmod_{\frac{1}{2}}(\Bun_G)\underset{\Sph_{G,\ul{x}}}\otimes \on{I}(G,P^-)_{\on{co},\ul{x}}^{\on{loc}}\to 
\Dmod_{\frac{1}{2}}(\Bun_M).
\end{equation} 

The second arrow in \eqref{e:pre CT enh 2} is the sought-for functor $\on{CT}_{*,\ul{x}}^{-,\on{enh}_{\on{co}}}$ of \eqref{e:CT co enh pt}. 

\begin{rem}

Note that by the same mechanism as in \lemref{l:KM to Sph gen}, the functor 
\begin{multline*}
\Dmod_{\frac{1}{2}}(\Bun_G)^{-,\on{enh}_{\on{co}}}_{\ul{x}}=\Dmod_{\frac{1}{2}}(\Bun_G)\underset{\Sph_{G,\ul{x}}}\otimes \on{I}(G,P^-)_{\on{co},\ul{x}}^{\on{loc}}\to \\
\to \Dmod_{\frac{1}{2}}(\Bun_G)\underset{\Sph_{G,\ul{x}}}\otimes \Dmod_{\frac{1}{2}}(\Gr_{G,\ul{x}})_{\fL(N^-_P)_{\ul{x}}}^{\fL^+(M)_{\ul{x}}}
\end{multline*} 
is an equivalence.

\end{rem}

\sssec{}

By construction, the composition
\begin{multline*}
\Dmod_{\frac{1}{2}}(\Bun_G)\overset{\on{Id}\otimes \one_{\on{I}(G,P^-)_{\on{co},\ul{x}}^{\on{loc}}}}\longrightarrow 
\Dmod_{\frac{1}{2}}(\Bun_G)\otimes \on{I}(G,P^-)_{\on{co},\ul{x}}^{\on{loc}} \to  \\
\to \Dmod_{\frac{1}{2}}(\Bun_G)\underset{\Sph_{G,\ul{x}}}\otimes \on{I}(G,P^-)_{\on{co},\ul{x}}^{\on{loc}}
=\Dmod_{\frac{1}{2}}(\Bun_G)^{-,\on{enh}_{\on{co}}}_{\ul{x}} \overset{\on{CT}_{*,\ul{x}}^{-,\on{enh}_{\on{co}}}}\longrightarrow  \\
\to \Dmod_{\frac{1}{2}}(\Bun_M)
\end{multline*} 
recovers the functor $\on{CT}_*^-$. 

\medskip

In families, the composition
\begin{multline*}
\Dmod_{\frac{1}{2}}(\Bun_G)\otimes \Dmod(\CZ)
\overset{\on{Id}\otimes \one_{\on{I}(G,P^-)_{\on{co},\CZ}^{\on{loc}}}}\longrightarrow 
\Dmod_{\frac{1}{2}}(\Bun_G)\otimes \on{I}(G,P^-)_{\on{co},\CZ}^{\on{loc}} \to \\
\to (\Dmod_{\frac{1}{2}}(\Bun_G)\otimes \Dmod(\CZ))\underset{\Sph_{G,\CZ}}\otimes \on{I}(G,P^-)_{\on{co},\CZ}^{\on{loc}}= 
\Dmod_{\frac{1}{2}}(\Bun_G)^{-,\on{enh}_{\on{co}}}_{\CZ} \overset{\on{CT}_{*,\CZ}^{-,\on{enh}_{\on{co}}}}\longrightarrow \\
\to \Dmod_{\frac{1}{2}}(\Bun_M)\otimes \Dmod(\CZ)
\end{multline*} 
identifies with $\on{CT}_*^-\otimes \on{Id}_{\Dmod(\CZ)}$. 

\sssec{}

By construction, the functor \eqref{e:CT co enh pt} respects the actions of $\Sph_{M,x}$ on the two sides. 

\medskip

In families, the functor \eqref{e:CT co enh Z} respects the actions of $\Sph_{M,\CZ}$ on the two sides. 

\sssec{}

Unwinding the construction, we obtain that the resulting functor \eqref{e:CT co enh Shv} has a natural
unital structure, in the sense of \cite[Sect. 11.3]{GLC2}. 

\ssec{Proof of \lemref{l:approx prel}} \label{ss:proof of approx prel}

\sssec{}

In order to simplify the exposition, we prove the assertion at the level of !-fibers at the trivial $M$-bundle $\CP^0_M\in \Bun_M$.

\medskip

Denote:
$$\on{Hecke}^{\on{glob}}_{G,N_P^-,\ul{x}}:=\on{pt}\underset{\CP^0_M,\Bun_M}\times \on{Hecke}^{\on{glob}}_{G,P^-,\ul{x}}$$
and 
$$\Bun_G \overset{\hl_{G,N_P^-}}\leftarrow \on{Hecke}^{\on{glob}}_{G,N_P^-,\ul{x}} \overset{\hr_{G,N_P^-}}\to \Bun_{N_P^-}.$$

We denote by the same symbol $\sfs_{\ul{x}}$ the resulting map
$$\on{Hecke}^{\on{glob}}_{G,N_P^-,\ul{x}}\to \fL^+(N^-_P)_{\ul{x}}\backslash \Gr_{G,\ul{x}}.$$

Consider the resulting functor
\begin{equation} \label{e:pre pre CT enh N}
(\hr_{G,N_P^-})_*\left((\hl_{G,N_P^-})^!(-)\sotimes \sfs_{\ul{x}}^!(-)\right),
\end{equation} 
which maps 
$$\Dmod_{\frac{1}{2}}(\Bun_G)\otimes \Dmod_{\frac{1}{2}}(\Gr_{G,\ul{x}})^{\fL^+(N_P^-)_{\ul{x}}}
\to \Dmod(\Bun_{N_P^-}).$$

We will show that the composite functor
\begin{equation} \label{e:pre pre CT enh N comp}
\Dmod_{\frac{1}{2}}(\Bun_G)\otimes \Dmod_{\frac{1}{2}}(\Gr_{G,\ul{x}})^{\fL^+(N_P^-)_{\ul{x}}}
\overset{\text{\eqref{e:pre pre CT enh N}}}\longrightarrow \Dmod(\Bun_{N_P^-}) \overset{\on{C}^\cdot(\Bun_{N_P^-},-)}\longrightarrow \Vect
\end{equation} 
factors via
$$\Dmod_{\frac{1}{2}}(\Gr_{G,\ul{x}})^{\fL^+(N_P^-)_{\ul{x}}}\simeq \Dmod_{\frac{1}{2}}(\Gr_{G,\ul{x}})_{\fL^+(N_P^-)_{\ul{x}}}\twoheadrightarrow
\Dmod_{\frac{1}{2}}(\Gr_{G,\ul{x}})_{\fL(N_P^-)_{\ul{x}}}$$
along the second factor.

\sssec{}

Let $\Bun_{N_P^-}^{\on{level}_{\ul{x}}}$ be the moduli stack of $N_P^-$-bundles on $X$ with a full level structure at $\ul{x}$. 
Parallel to \eqref{e:pre pre CT enh N} we have a functor
\begin{equation} \label{e:pre pre CT enh N x}
\Dmod_{\frac{1}{2}}(\Bun_G)\otimes \Dmod_{\frac{1}{2}}(\Gr_{G,\ul{x}})
\to \Dmod(\Bun_{N_P^-}^{\on{level}_{\ul{x}}}),
\end{equation}
which is $\fL(N^-_P)_{\ul{x}}$-equivariant, so that the functor \eqref{e:pre pre CT enh N} is obtained from \eqref{e:pre pre CT enh N x} by
taking $\fL^+(N^-_P)_{\ul{x}}$-invariants/coinvariants on both sides. 

\medskip

Consider the functor 
\begin{equation} \label{e:pre pre CT enh N comp x}
\Dmod_{\frac{1}{2}}(\Bun_G)\otimes \Dmod_{\frac{1}{2}}(\Gr_{G,\ul{x}})
\overset{\text{\eqref{e:pre pre CT enh N}}}\longrightarrow \Dmod(\Bun^{\on{level}_{\ul{x}}}_{N_P^-}) 
\overset{\on{C}^\cdot(\Bun^{\on{level}_{\ul{x}}}_{N_P^-},-)}\longrightarrow \Vect
\end{equation} 

It naturally factors via $\fL^+(N^-_P)_{\ul{x}}$-coinvariants along the second factor, and the resulting functor
$$\Dmod_{\frac{1}{2}}(\Bun_G)\otimes \Dmod_{\frac{1}{2}}(\Gr_{G,\ul{x}})_{\fL^+(N_P^-)_{\ul{x}}}\to \Vect$$
is the functor \eqref{e:pre pre CT enh N comp}.

\sssec{}

Now the desired assertion follows from the fact that the functor \eqref{e:pre pre CT enh N comp x} is
$\fL(N^-_P)_{\ul{x}}$-equivariant, since both arrows that comprise it have this property. 

\ssec{The enhanced Eisenstein and constant term functors} \label{e:enh CT}

In this subsection, we will use duality and adjunction manipulations to produce from $\ul{\on{CT}}_*^{-,\on{enh}_{\on{co}}}$
other enhanced functors. 

\sssec{}

By duality (see \secref{sss:duality over Sph}), the functor \eqref{e:CT co enh Z} gives rise to a $\Sph_{G,\CZ}$-linear functor
\begin{equation} \label{e:CT enh Z}
\Dmod_{\frac{1}{2}}(\Bun_G) \otimes \Dmod(\CZ) \to \Dmod_{\frac{1}{2}}(\Bun_M)^{-,\on{enh}}_\CZ.
\end{equation}

We denote the functor \eqref{e:CT enh Z} by $\on{CT}_{*,\CZ}^{-,\on{enh}}$ and refer to it as the \emph{enhanced constant term} functor.

\medskip 

We will view \eqref{e:CT enh Z} as a map between crystals of categories on $\Ran$, to be denoted $\ul{\on{CT}}_*^{-,\on{enh}}$. 

\sssec{} \label{sss:CT enh oblv}

By construction, the composition
$$\Dmod_{\frac{1}{2}}(\Bun_G) \otimes \Dmod(\CZ) \overset{\on{CT}_{*,\CZ}^{-,\on{enh}}}\longrightarrow \Dmod_{\frac{1}{2}}(\Bun_M)^{-,\on{enh}}_\CZ
\overset{\oblv_{\on{enh}}}\longrightarrow \Dmod_{\frac{1}{2}}(\Bun_M) \otimes \Dmod(\CZ)$$
identifies with $\on{CT}_*^-\otimes \on{Id}_{\Dmod(\CZ)}$. 

\sssec{}

Note that the datum of $\ul{\on{CT}}_*^{-,\on{enh}}$ (without the action of $\Sph_G$) 
is recovered from the datum of the functor 
\begin{equation} \label{e:CT enh Ran}
\on{CT}_{*,\Ran}^{-,\on{enh}}: \Dmod_{\frac{1}{2}}(\Bun_G) \to  \Dmod_{\frac{1}{2}}(\Bun_M)^{-,\on{enh}}_\Ran.
\end{equation}

\medskip

The fact that \eqref{e:CT co enh Shv} is unital implies that the essential image of \eqref{e:CT enh Ran} is contained
in the full subcategory
$$\left(\Dmod_{\frac{1}{2}}(\Bun_M)^{-,\on{enh}}\right)_{\Ran^{\on{untl}},\on{indep}}\subset \Dmod_{\frac{1}{2}}(\Bun_M)^{-,\on{enh}}_\Ran,$$
see Remark \ref{r:global IGP untl}.

\begin{rem} 

In terms of the equivalence
$$\left(\Dmod_{\frac{1}{2}}(\Bun_M)^{-,\on{enh}}\right)_{\Ran^{\on{untl}},\on{indep}}\simeq \on{I}(G,P^-)^{\on{glob}},$$
the resulting functor
$$\Dmod_{\frac{1}{2}}(\Bun_G) \to  \on{I}(G,P^-)^{\on{glob}}$$
identifies with the enhanced constant term functor of \cite[Sect. 6.3.1]{Ga1}.

\end{rem}

\sssec{}

We claim now that the functor $\on{CT}_{*,\CZ}^{-,\on{enh}}$ admits a left adjoint, to be denoted $\Eis^{-,\on{enh}}_{!,\CZ}$. 

\medskip

Indeed, it suffices to show that $\Eis^{-,\on{enh}}_{!,\CZ}$ is defined on the generators of $\Dmod_{\frac{1}{2}}(\Bun_M)^{-,\on{enh}}_\CZ$.
Hence, it is enough to show that the functor
$$\Eis^{-,\on{enh}}_{!,\CZ}\circ \ind_{\on{enh}},$$
left adjoint to $\oblv_{\on{enh}}\circ \on{CT}_{*,\CZ}^{-,\on{enh}}$ is well-defined. 

\medskip

However, by \secref{sss:CT enh oblv}, 
$$\oblv_{\on{enh}}\circ \on{CT}_{*,\CZ}^{-,\on{enh}}\simeq \on{CT}_*^-\otimes \on{Id}_{\Dmod(\CZ)},$$
and hence
$$\Eis^{-,\on{enh}}_{!,\CZ}\circ \ind_{\on{enh}}\simeq \Eis^-_!\otimes \on{Id}_{\Dmod(\CZ)}.$$

\medskip

Since $\Sph_{G,\CZ}$ is rigid (over $\CZ$), the fact that $\on{CT}_{*,\CZ}^{-,\on{enh}}$ is $\Sph_{G,\CZ}$-linear
implies that $\Eis^{-,\on{enh}}_{!,\CZ}$ inherits this structure. 

\sssec{}

The assignment 
$$\CZ\rightsquigarrow \Eis^{-,\on{enh}}_{!,\CZ}$$
is a ``local-to-global" functor in the sense of \cite[Sect. 11.1]{GLC2}. We denote it by
$$\ul\Eis^{-,\on{enh}}_!:\ul\Dmod_{\frac{1}{2}}(\Bun_M)^{-,\on{enh}}\to \Dmod_{\frac{1}{2}}(\Bun_G)\otimes \ul\Dmod(\Ran).$$

\medskip

Following the conventions of \cite[Sect. 11.1.4]{GLC2}, we will denote by $\Eis^{-,\on{enh}}_!$ the resulting functor
$$\Dmod_{\frac{1}{2}}(\Bun_M)^{-,\on{enh}}_\Ran\to \Dmod_{\frac{1}{2}}(\Bun_G).$$

\medskip

The fact that the essential image of $\on{CT}_{*,\Ran}^{-,\on{enh}}$ is contained in $$\left(\Dmod_{\frac{1}{2}}(\Bun_M)^{-,\on{enh}}\right)_{\Ran^{\on{untl}},\on{indep}}$$
implies that the functor $\ul\Eis^{-,\on{enh}}_!$ carries a unital structure. 

\sssec{} \label{sss:fact alg Eis}

Let $\CA$ be a factorization algebra in $\on{I}(G,P^-)^{\on{loc}}$. We define the functor
$$\Eis^-_{!,\CA}:\Dmod_{\frac{1}{2}}(\Bun_M)\to \Dmod_{\frac{1}{2}}(\Bun_G)$$ to be the composition 
\begin{multline*}
\Dmod_{\frac{1}{2}}(\Bun_M) \overset{\on{Id}\otimes \CA_\Ran}\to 
\Dmod_{\frac{1}{2}}(\Bun_M) \otimes \on{I}(G,P^-)^{\on{loc}}_\Ran \to \\
\to (\Dmod_{\frac{1}{2}}(\Bun_M) \otimes \Dmod(\Ran))\underset{\Sph_{M,\Ran}}\otimes \on{I}(G,P^-)^{\on{loc}}_\Ran = \\
= \Dmod_{\frac{1}{2}}(\Bun_M)^{-,\on{enh}}_\Ran \overset{\Eis^{-,\on{enh}}_!}\to \Dmod(\Bun_G).
\end{multline*} 

Taking $\CA$ to be 
$$\Delta^{-,\semiinf},\,\, \nabla^{-,\semiinf} \text{ and } \IC^{-,\semiinf},$$
respectively, we obtain the functors 
\begin{equation} \label{e:Eis enh Delta}
\Eis^-_{!,\Delta}:\Dmod_{\frac{1}{2}}(\Bun_M)\to \Dmod_{\frac{1}{2}}(\Bun_G),
\end{equation} 
\begin{equation} \label{e:Eis enh nabla}
\Eis^-_{!,\nabla}:\Dmod_{\frac{1}{2}}(\Bun_M)\to \Dmod_{\frac{1}{2}}(\Bun_G)
\end{equation} 
and
\begin{equation} \label{e:Eis enh IC}
\Eis^-_{!,\IC}:\Dmod_{\frac{1}{2}}(\Bun_M)\to \Dmod_{\frac{1}{2}}(\Bun_G).
\end{equation} 

Since $\Delta^{-,\semiinf}$ is the factorization unit in $\on{I}(G,P^-)^{\on{loc}}$, by unitality we have
$$\Eis^-_{!,\Delta}\simeq \Eis^-_!.$$

\medskip

We will shortly see that the functor \eqref{e:Eis enh IC} also identifies with an a priori
familiar functor, see \propref{p:compact Eis}. 

\medskip

A similar (but simpler) argument shows that the functor \eqref{e:Eis enh nabla} identifies with the
functor $\Eis^-_*$ from \cite[Sect. 1.1.6]{Ga3} (up to the cohomological shift by \eqref{e:global shift}). 

\sssec{} 

Similarly, for a factorization algebra $\CB$ in $\on{I}(G,P^-)^{\on{loc}}_{\on{co}}$, we obtain a functor
$$\on{CT}^-_{*,\CB}:\Dmod_{\frac{1}{2}}(\Bun_G)\to \Dmod_{\frac{1}{2}}(\Bun_M).$$

For $\CB:=\Delta^{-,\semiinf}_{\on{co}}=\one_{\on{I}(G,P^-)^{\on{loc}}_{\on{co}}}$, we recover the functor
$\on{CT}^-_*$. 

\ssec{Enhanced functors and duality}

\sssec{}

Consider now the categories 
$$\Dmod_{\frac{1}{2}}(\Bun_M)_{\on{co}} \text{ and } \Dmod_{\frac{1}{2}}(\Bun_G)_{\on{co}},$$
see \cite[Sect. 4.3.3]{DG2}.

\medskip

We introduce the following (unital) crystals of categories over the Ran space:
$$\ul\Dmod_{\frac{1}{2}}(\Bun_G)_{\on{co}}^{-,\on{enh}}:$$
$$\CZ\mapsto \Dmod_{\frac{1}{2}}(\Bun_G)^{-,\on{enh}}_{\on{co},\CZ}:=
\left(\Dmod_{\frac{1}{2}}(\Bun_G)_{\on{co}}\otimes \Dmod(\CZ)\right)\underset{\Sph_{G,\CZ}}\otimes \on{I}(G,P^-)^{\on{loc}}_\CZ,$$
and 
$$\ul\Dmod_{\frac{1}{2}}(\Bun_M)_{\on{co}}^{-,\on{enh}_{\on{co}}}:$$
$$\CZ\mapsto \Dmod_{\frac{1}{2}}(\Bun_M)^{-,\on{enh}_{\on{co}}}_{\on{co},\CZ}:=
\left(\Dmod_{\frac{1}{2}}(\Bun_M)_{\on{co}}\otimes \Dmod(\CZ)\right)\underset{\Sph_{M,\CZ}}\otimes \on{I}(G,P^-)^{\on{loc}}_{\on{co},\CZ}.$$

\sssec{}

Passing to dual functors, we obtain the functors
$$\ul\Eis^{-,\on{enh}_{\on{co}}}_{\on{co},*}:=\left(\ul{\on{CT}}_*^{-,\on{enh}}\right)^\vee, \quad
\ul\Dmod_{\frac{1}{2}}(\Bun_M)_{\on{co}}^{-,\on{enh}_{\on{co}}}\to \Dmod_{\frac{1}{2}}(\Bun_G)_{\on{co}}\otimes \ul\Dmod(\Ran),$$
$$\ul{\on{CT}}^{-,\on{enh}_{\on{co}}}_{\on{co},?}:=(\ul\Eis^{-,\on{enh}}_!)^\vee, \quad
\Dmod_{\frac{1}{2}}(\Bun_G)_{\on{co}}\otimes \ul\Dmod(\Ran)\to \ul\Dmod_{\frac{1}{2}}(\Bun_M)_{\on{co}}^{-,\on{enh}_{\on{co}}}.$$

Note that the functors $(\ul\Eis^{-,\on{enh}_{\on{co}}}_{\on{co},*},\ul{\on{CT}}^{-,\on{enh}_{\on{co}}}_{\on{co},?})$ form adjoint pair.

\medskip

Using \secref{sss:duality over Sph}, from the functor $\ul{\on{CT}}^{-,\on{enh}_{\on{co}}}_{\on{co},?}$ we produce the functor
$$\ul{\on{CT}}^{-,\on{enh}}_{\on{co},?}:\ul\Dmod_{\frac{1}{2}}(\Bun_G)_{\on{co}}^{-,\on{enh}}\to \Dmod_{\frac{1}{2}}(\Bun_M)_{\on{co}}\otimes \ul\Dmod(\Ran).$$

\sssec{}

For any $\CZ\to \Ran$, the composition
\begin{multline*}
\Dmod_{\frac{1}{2}}(\Bun_M)_{\on{co}} \overset{\on{Id}\otimes \one_{\on{I}(G,P^-)^{\on{loc}}_{\on{co},\CZ}}}\longrightarrow 
\Dmod_{\frac{1}{2}}(\Bun_M)_{\on{co}}\otimes \on{I}(G,P^-)^{\on{loc}}_{\on{co},\CZ} \to \\
\to \left(\Dmod_{\frac{1}{2}}(\Bun_M)_{\on{co}}\otimes  \Dmod(\CZ)\right)\underset{\Sph_{M,\CZ}}\otimes \on{I}(G,P^-)^{\on{loc}}_{\on{co},\CZ}=
\Dmod_{\frac{1}{2}}(\Bun_M)_{\on{co},\CZ}^{-,\on{enh}_{\on{co}}} \overset{\ul\Eis^{-,\on{enh}_{\on{co}}}_{\on{co},*}}\longrightarrow \\
\to \Dmod_{\frac{1}{2}}(\Bun_G)_{\on{co}}\otimes\Dmod(\CZ)
\end{multline*}
identifies with  
$$\Eis_{\on{co},*}\otimes \on{Id}_{\Dmod(\CZ)},$$
where $\Eis_{\on{co},*}$ is the functor from \cite[Sect. 1.4.1]{Ga3} (up to the cohomological shift by \eqref{e:global shift}). . 

\medskip

The composition
\begin{multline*}
\Dmod_{\frac{1}{2}}(\Bun_G)_{\on{co}} \overset{\on{Id}\otimes \one_{\on{I}(G,P^-)^{\on{loc}}_\CZ}}\longrightarrow 
\Dmod_{\frac{1}{2}}(\Bun_G)_{\on{co}} \otimes \on{I}(G,P^-)^{\on{loc}}_\CZ \to \\
\to \Dmod_{\frac{1}{2}}(\Bun_G)_{\on{co}} \underset{\Sph_{G,\CZ}}\otimes \on{I}(G,P^-)^{\on{loc}}_\CZ =
\Dmod_{\frac{1}{2}}(\Bun_G)_{\on{co},\CZ}^{-,\on{enh}} \overset{\on{CT}^{-,\on{enh}}_{\on{co},?}}\longrightarrow \\
\to \Dmod_{\frac{1}{2}}(\Bun_M)_{\on{co}}\otimes\Dmod(\CZ)
\end{multline*}
identifies with  
$$\on{CT}^-_{\on{co},?}\otimes \on{Id}_{\Dmod(\CZ)},$$
where 
\begin{equation} \label{e:CT ?}
\on{CT}^-_{\on{co},?}:=(\Eis^-_!)^\vee
\end{equation} 
is the functor from \cite[Sect. 1.4.5]{Ga3} (up to the cohomological shift by \eqref{e:global shift}).  

\sssec{}

The functors $\ul\Eis^{-,\on{enh}_{\on{co}}}_{\on{co},*}$ and $\ul{\on{CT}}^{-,\on{enh}}_{\on{co},?}$ are unital ``local-to-global"
functors in the sense of \cite[Sect. 11.1 and 11.3]{GLC2}. We will denote by
$$\Eis^{-,\on{enh}_{\on{co}}}_{\on{co},*} \text{ and } \on{CT}^{-,\on{enh}}_{\on{co},?}$$
the resulting functors
$$\Dmod_{\frac{1}{2}}(\Bun_M)_{\on{co},\Ran}^{-,\on{enh}_{\on{co}}}\to \Dmod_{\frac{1}{2}}(\Bun_G)_{\on{co}}$$
and
$$\Dmod_{\frac{1}{2}}(\Bun_G)_{\on{co},\Ran}^{-,\on{enh}}\to \Dmod_{\frac{1}{2}}(\Bun_M)_{\on{co}},$$
respectively.

\sssec{} \label{sss:fact alg Eis co}

By the same recipe as in \secref{sss:fact alg Eis}, for a factorization algebra $\CA$ in $\on{I}(G,P^-)^{\on{loc}}$,
we obtain a functor
$$\on{CT}^-_{\on{co},?,\CA}:\Dmod_{\frac{1}{2}}(\Bun_G)_{\on{co}}\to \Dmod_{\frac{1}{2}}(\Bun_M)_{\on{co}}.$$

Taking $\CA$ to be 
$$\Delta^{-,\semiinf},\,\, \nabla^{-,\semiinf} \text{ and } \IC^{-,\semiinf},$$
respectively, we obtain the functors 
\begin{equation} \label{e:CT co enh Delta}
\on{CT}^-_{\on{co},?,\Delta}:\Dmod_{\frac{1}{2}}(\Bun_G)_{\on{co}}\to \Dmod_{\frac{1}{2}}(\Bun_M)_{\on{co}}.
\end{equation} 
\begin{equation} \label{e:CT co enh nabla}
\on{CT}^-_{\on{co},?,\nabla}:\Dmod_{\frac{1}{2}}(\Bun_G)_{\on{co}}\to \Dmod_{\frac{1}{2}}(\Bun_M)_{\on{co}}
\end{equation} 
and
\begin{equation} \label{e:CT co enh IC}
\on{CT}^-_{\on{co},?,\IC}:\Dmod_{\frac{1}{2}}(\Bun_G)_{\on{co}}\to \Dmod_{\frac{1}{2}}(\Bun_M)_{\on{co}}.
\end{equation} 

\medskip

Since $\Delta^{-,\semiinf}$ is the factorization unit in $\on{I}(G,P^-)^{\on{loc}}$, we have
$$\on{CT}^-_{\on{co},?,\Delta}\simeq \on{CT}^-_{\on{co},?}.$$

Additionally, it is not difficult to show that
$$\on{CT}^-_{\on{co},?,\nabla}\simeq \on{CT}^-_{\on{co},*},$$
where $\on{CT}^-_{\on{co},*}\simeq (\Eis^-_*)^\vee$ is as in \cite[Sect. 1.4.3]{Ga3} (up to the cohomological shift by \eqref{e:global shift}). 

\sssec{} \label{sss:fact alg CT}

Similarly, for a factorization algebra $\CB$ in $\on{I}(G,P^-)^{\on{loc}}_{\on{co}}$, we obtain a functor
$$\Eis^-_{\on{co},*,\CB}:\Dmod_{\frac{1}{2}}(\Bun_M)_{\on{co}}\to \Dmod_{\frac{1}{2}}(\Bun_G)_{\on{co}}.$$

For $\CB=\Delta^{-,\semiinf}_{\on{co}}$ we recover the functor $\Eis^-_{\on{co},*}$.
 
\ssec{Enhanced functors and miraculous duality}

\sssec{}

Recall now the miraculous functor
$$\Mir_{\Bun_G}:\Dmod_{\frac{1}{2}}(\Bun_G)_{\on{co}}\to \Dmod_{\frac{1}{2}}(\Bun_G),$$
see \cite[Sect. 4.4.8]{DG2}.\footnote{In {\it loc. cit.}, this functor is denoted $\on{Ps-Id}_{\Bun_G,!}$.}

\begin{rem}

In \cite{Ga3}, the functor $\Mir_{\Bun_G}$ was defined in the untwisted setting, i.e., for $\Dmod(\Bun_G)$
rather than for $\Dmod_{\frac{1}{2}}(\Bun_G)$. However, the same procedure defines is also for $\Dmod_{\frac{1}{2}}(\Bun_G)$, and the resulting
functor has the same properties: the two setting are actually equivalent, see \cite[Remark 9.1.6]{GLC2}. 

\end{rem} 

\sssec{}

By construction, $\Mir_{\Bun_G}$ commutes with the natural actions of $\Sph_G$ on the two sides. 

\medskip

According to \cite[Theorem 3.1.5]{Ga3}, the functor $\Mir_{\Bun_G}$ is an equivalence. Consider also the corresponding functor for $M$. 
$$\Mir_{\Bun_M}:\Dmod_{\frac{1}{2}}(\Bun_M)_{\on{co}}\to \Dmod_{\frac{1}{2}}(\Bun_M).$$

\medskip

Denote
$$\Mir_{\Bun_G,\tau}:=\tau_G\circ \Mir_{\Bun_G} \text{ and } \Mir_{\Bun_M,\tau}:=\tau_M\circ \Mir_{\Bun_M},$$
where $\tau_G$ and $\tau_M$ are the Chevalley involutions in $G$ and $M$ respectively. 

\sssec{}

Recall now the equivalence
$$\Theta_{\on{I}(G,P^-)^{\on{loc}}}:\on{I}(G,P^-)^{\on{loc}}_{\on{co}}\to \on{I}(G,P)^{\on{loc}},$$
see \corref{c:geom intertwiner}. 

\medskip 

Composing with the Chevalley involution $\tau_G$, we obtain an equivalence that we will denote
\begin{equation} \label{e:Theta IGP glob}
\Theta_{\on{I}(G,P^-)^{\on{loc}},\tau}:\on{I}(G,P^-)^{\on{loc}}_{\on{co}}\to \on{I}(G,P^-)^{\on{loc}}.
\end{equation} 

It commutes with the actions of $\Sph_G$ and $\Sph_M$, up to the Chevalley involutions.

\medskip

Tensoring $\Theta_{\on{I}(G,P^-)^{\on{loc}},\tau}$ with $\Mir_{\Bun_M,\tau}$
we obtain an equivalence that we denote 
$$\Mir^{-,\on{enh}}_{\Bun_M,\tau}:
\ul\Dmod_{\frac{1}{2}}(\Bun_M)_{\on{co}}^{-,\on{enh}_{\on{co}}} \to \ul\Dmod_{\frac{1}{2}}(\Bun_M)^{-,\on{enh}}.$$
 
\sssec{} 
 
We now quote the following result (see \cite[Theorem 5.3.5(b)]{Ch2}):

\begin{thm} \label{t:intertwiners}
The following diagram commutes 
\begin{equation} \label{e:intertwiners}
\CD
\Dmod_{\frac{1}{2}}(\Bun_M)_{\on{co},\Ran}^{-,\on{enh}_{\on{co}}}  @>{\Mir^{-,\on{enh}}_{\Bun_M,\tau}}>> \Dmod_{\frac{1}{2}}(\Bun_M)^{-,\on{enh}}_\Ran \\
@V{\Eis^{-,\on{enh}_{\on{co}}}_{\on{co},*}[\delta_{N^-_P}]}VV @VV{\Eis^{-,\on{enh}}_![-\delta_{N^-_P}]}V \\
\Dmod_{\frac{1}{2}}(\Bun_G)_{\on{co}} @>>{\Mir_{\Bun_G,\tau}}> \Dmod_{\frac{1}{2}}(\Bun_G),
\endCD
\end{equation}
where 
$$\delta_{N^-_P}=\dim(\Bun_{N^-_P}).$$
\end{thm} 

\begin{rem}

    More precisely, \thmref{t:intertwiners} differs from \cite[Theorem 5.3.5(b)]{Ch2} in the following technical points:
    \begin{itemize}
        
        \item
            In \emph{loc.cit.}, the theorem was about the corresponding ``independent'' categories (see Remark \ref{r:global IGP untl}). Note that the functors $\Eis^{-,\on{enh}_{\on{co}}}_{\on{co},*}$ and $\Eis^{-,\on{enh}}_!$ indeed factor through these ``independent'' categories.
            
            \medskip
            
        \item
            In \emph{loc.cit.}, the equivalence $\Mir^{-,\on{enh}}_{\Bun_M,\tau}$, referred as the \emph{global inv-inv duality} (up to the Chevalley involution), was constructed by a different method, which can be shown to be equivalent to the current construction by combining \cite[Corollary 3.2.4]{Ch2} and 
            \cite[Theorem 5.1.7]{Ch2}.

    \end{itemize}

\end{rem}

\begin{rem}

\thmref{t:intertwiners} generalizes the result of \cite[Theorem 4.1.2]{Ga3}, which says that the diagram
\begin{equation} \label{e:funct eq}
\CD
\Dmod_{\frac{1}{2}}(\Bun_M)_{\on{co}} @>{\Mir_{\Bun_M,\tau}}>> \Dmod_{\frac{1}{2}}(\Bun_M) \\
@V{\Eis^-_{\on{co},*}[\delta_{N^-_P}]}VV @VV{\Eis^-_![-\delta_{N^-_P}]}V \\
\Dmod_{\frac{1}{2}}(\Bun_G)_{\on{co}} @>>{\Mir_{\Bun_G,\tau}}> \Dmod_{\frac{1}{2}}(\Bun_G),
\endCD
\end{equation}
commutes.

\end{rem}

\sssec{}

Tensoring $\Mir_{\Bun_G,\tau}$ with $(\Theta_{\on{I}(G,P^-)^{\on{loc}},\tau})^{-1}$ we obtain an equivalence that we denote 
$$\Mir^{-,\on{enh}}_{\Bun_G,\tau}:\ul\Dmod_{\frac{1}{2}}(\Bun_G)_{\on{co}}^{-,\on{enh}}\to 
\ul\Dmod_{\frac{1}{2}}(\Bun_G)_\Ran^{-,\on{enh}_{\on{co}}}.$$

\medskip

By duality, from \thmref{t:intertwiners} we obtain: 

\begin{cor} \label{c:intertwiners CT}
The following diagram commutes:
$$
\CD
\Dmod_{\frac{1}{2}}(\Bun_M)_{\on{co}} @>{\Mir_{\Bun_M,\tau}}>> \Dmod_{\frac{1}{2}}(\Bun_M) \\
@A{\on{CT}^{-,\on{enh}}_{\on{co},?}[-\delta_{N^-_P}]}AA @AA{\on{CT}_*^{-,\on{enh}_{\on{co}}}[\delta_{N^-_P}]}A \\
\Dmod_{\frac{1}{2}}(\Bun_G)_{\on{co},\Ran}^{-,\on{enh}} @>>{\Mir^{-,\on{enh}}_{\Bun_G,\tau}}> \Dmod_{\frac{1}{2}}(\Bun_G)_\Ran^{-,\on{enh}_{\on{co}}}.
\endCD
$$
\end{cor} 

In particular, we have:
\begin{cor} \label{c:intertwiners CT simple}
The following diagram commutes:
$$
\CD
\Dmod_{\frac{1}{2}}(\Bun_M)_{\on{co}} @>{\Mir_{\Bun_M,\tau}}>> \Dmod_{\frac{1}{2}}(\Bun_M) \\
@A{\on{CT}^-_{\on{co},?}[-\delta_{N^-_P}]}AA @AA{\on{CT}_*^-[\delta_{N^-_P}]}A \\
\Dmod_{\frac{1}{2}}(\Bun_G)_{\on{co}} @>>{\Mir_{\Bun_G,\tau}}> \Dmod_{\frac{1}{2}}(\Bun_G).
\endCD
$$
\end{cor} 

(Note, however, that \corref{c:intertwiners CT simple} can be directly obtained by duality from \eqref{e:funct eq}). 

\ssec{The \texorpdfstring{$\rho$}{rho}-shift} \label{ss:rho shift Eis}

\sssec{}

Let $\CP_{Z_M}$ be a $Z_M$-torsor on $X$, and consider again the auto-equivalence 
\begin{equation} \label{e:transl Bun M}
(\on{transl}_{\CP_{Z_M}})^*:\Dmod_{\frac{1}{2}}(\Bun_M) \to \Dmod_{\frac{1}{2}}(\Bun_M).
\end{equation} 

Recall now the automorphism $(\on{transl}_{\CP_{Z_M}})^*$ of 
$\Sph_M$, see \eqref{e:autom Sph}. These two automorphisms are compatible with the action of $\Sph_{M,\CZ}$ on 
$\Dmod_{\frac{1}{2}}(\Bun_M)\otimes \Dmod(\CZ)$, for $\CZ\to \Ran$. 

\sssec{}

Recall the identification 
\begin{equation} \label{e:alpha IGP again}
\on{I}(G,P^-)^{\on{loc}} \overset{\alpha_{\CP_{Z_M},\on{taut}}}\simeq \on{I}(G,P^-)^{\on{loc}}_{\CP_{Z_M}}
\end{equation}
as categories acted on by $\Sph_G$ and $\Sph_M$ (see \secref{sss:convent Sph-action on IGP} regarding our conventions
for the latter action).

\medskip

The tensor product of \eqref{e:transl Bun M} and \eqref{e:alpha IGP again} defines an identification
\begin{equation} \label{e:IGP glob twist}
\Dmod_{\frac{1}{2}}(\Bun_M)^{-,\on{enh}}_\CZ\simeq \Dmod_{\frac{1}{2}}(\Bun_M)^{-,\on{enh}}_{\CP_{Z_M},\CZ},
\end{equation} 
as categories acted on by $\Sph_{G,\CZ}$, where 
\begin{equation} \label{e:twisted enh recipient}
\Dmod_{\frac{1}{2}}(\Bun_M)^{-,\on{enh}}_{\CP_{Z_M},\CZ}:=
\left(\Dmod_{\frac{1}{2}}(\Bun_M)\otimes \Dmod(\CZ)\right)\underset{\Sph_{M,\CZ}}\otimes \on{I}(G,P^-)^{\on{loc}}_{\CP_{Z_M},\CZ}.
\end{equation}

\sssec{}

We have an adjoint pair
$$\ind_{\on{enh}}:\Dmod_{\frac{1}{2}}(\Bun_M)\otimes \Dmod(\CZ)\rightleftarrows \Dmod_{\frac{1}{2}}(\Bun_M)^{-,\on{enh}}_{\CP_{Z_M},\CZ}:\oblv_{\on{enh}}$$
so that we have a commutative diagram
$$
\CD
\Dmod_{\frac{1}{2}}(\Bun_M)^{-,\on{enh}}_{\CZ} @>{\text{\eqref{e:IGP glob twist}}}>> \Dmod_{\frac{1}{2}}(\Bun_M)^{-,\on{enh}}_{\CP_{Z_M},\CZ} \\
@V{\oblv_{\on{enh}}}VV @VV{\oblv_{\on{enh}}}V \\
\Dmod_{\frac{1}{2}}(\Bun_M)\otimes \Dmod(\CZ) @>{(\on{transl}_{\CP_{Z_M}})^*\otimes \on{Id}}>> \Dmod_{\frac{1}{2}}(\Bun_M)\otimes \Dmod(\CZ).
\endCD
$$

\sssec{}

We let 
$$\on{CT}_{*,\CP_{Z_M},\CZ}^{-,\on{enh}}:\Dmod_{\frac{1}{2}}(\Bun_G)\otimes \Dmod(\CZ)\to \Dmod_{\frac{1}{2}}(\Bun_M)^{-,\on{enh}}_{\CP_{Z_M},\CZ}$$
denote the composition
$$\Dmod_{\frac{1}{2}}(\Bun_G)\otimes \Dmod(\CZ) \overset{\on{CT}_{*,\CZ}^{-,\on{enh}}}\longrightarrow 
\Dmod_{\frac{1}{2}}(\Bun_M)^{-,\on{enh}}_{\CZ} \overset{\text{\eqref{e:IGP glob twist}}}\longrightarrow \Dmod_{\frac{1}{2}}(\Bun_M)^{-,\on{enh}}_{\CP_{Z_M},\CZ}.$$

\medskip

This functor is compatible with the actions of $\Sph_{G,\Ran}$ on the two sides. 

\medskip

By construction, the composition
$$\Dmod_{\frac{1}{2}}(\Bun_G)\otimes \Dmod(\CZ) \overset{\on{CT}_{*,\CP_{Z_M},\CZ}^{-,\on{enh}}}\longrightarrow 
\Dmod_{\frac{1}{2}}(\Bun_M)^{-,\on{enh}}_{\CP_{Z_M},\CZ}\overset{\oblv_{\on{enh}}\otimes \on{Id}}\longrightarrow 
\Dmod_{\frac{1}{2}}(\Bun_M)\otimes \Dmod(\CZ)$$
identifies with $\on{CT}_{*,\CP_{Z_M}}\otimes \on{Id}$. 

\sssec{}

We let
$$\Eis^{-,\on{enh}}_{!,\CP_{Z_M},\CZ}:\Dmod_{\frac{1}{2}}(\Bun_M)^{-,\on{enh}}_{\CP_{Z_M},\CZ}\to \Dmod_{\frac{1}{2}}(\Bun_G)\otimes \Dmod(\CZ)$$
be the left adjoint of $\on{CT}_{*,\CP_{Z_M},\CZ}^{-,\on{enh}}$.

\medskip

By construction, we have a commutative diagram
$$
\CD
\Dmod_{\frac{1}{2}}(\Bun_M)^{-,\on{enh}}_{\CZ} @>{\text{\eqref{e:IGP glob twist}}}>> \Dmod_{\frac{1}{2}}(\Bun_M)^{-,\on{enh}}_{\CP_{Z_M},\CZ} \\
@V{\Eis^{-,\on{enh}}_{!,\CZ}}VV @VV{\Eis^{-,\on{enh}}_{!,\CP_{Z_M},\CZ}}V \\
\Dmod_{\frac{1}{2}}(\Bun_G)\otimes \Dmod(\CZ) @>{\on{Id}}>> \Dmod_{\frac{1}{2}}(\Bun_G)\otimes \Dmod(\CZ).
\endCD
$$

Also, by construction, the composition
$$\Dmod_{\frac{1}{2}}(\Bun_M)\otimes \Dmod(\CZ) \overset{\ind_{\on{enh}}}\longrightarrow 
\Dmod_{\frac{1}{2}}(\Bun_M)^{-,\on{enh}}_{\CP_{Z_M},\CZ} \overset{\Eis^{-,\on{enh}}_{!,\CP_{Z_M},\CZ}}\longrightarrow 
\Dmod_{\frac{1}{2}}(\Bun_G)\otimes \Dmod(\CZ)$$
identifies with $\Eis^-_{!,\CP_{Z_M}}\otimes \on{Id}$. 

\sssec{}

Denote
$$\Dmod_{\frac{1}{2}}(\Bun_G)^{-,\on{enh}_{\on{co}}}_{\CP_{Z_M},\CZ}:=
\left(\Dmod_{\frac{1}{2}}(\Bun_G)\otimes \Dmod(\CZ)\right)\underset{\Sph_{G,\CZ}}\otimes \on{I}(G,P^-)^{\on{loc}}_{\on{co},\CP_{Z_M},\CZ}.$$

Applying \eqref{e:alpha IGP again} along the second factor we obtain an identification
\begin{equation} \label{e:IGP glob twist co}
\Dmod_{\frac{1}{2}}(\Bun_G)^{-,\on{enh}_{\on{co}}}_{\CZ}\simeq \Dmod_{\frac{1}{2}}(\Bun_G)^{-,\on{enh}_{\on{co}}}_{\CP_{Z_M},\CZ}.
\end{equation} 

Let
$$\on{CT}_{*,\CP_{Z_M},\CZ}^{-,\on{enh}_{\on{co}}}:\Dmod_{\frac{1}{2}}(\Bun_G)^{-,\on{enh}_{\on{co}}}_{\CP_{Z_M},\CZ}\to 
\Dmod_{\frac{1}{2}}(\Bun_M)\otimes \Dmod(\CZ)$$
be the functor obtained from $\on{CT}_{*,\CP_{Z_M},\CZ}^{-,\on{enh}}$ by duality. 

\medskip

By construction, this functor is compatible with the actions of $\Sph_{M,\Ran}$ on the two sides. 

\medskip

The composition
\begin{multline*}
\Dmod_{\frac{1}{2}}(\Bun_G)\otimes \Dmod(\CZ) 
\overset{\on{Id}\otimes \one_{\on{I}(G,P^-)^{\on{loc}}_{\on{co},\CP_{Z_M},\CZ}}}\longrightarrow  \\
\to \Dmod_{\frac{1}{2}}(\Bun_G)^{-,\on{enh}_{\on{co}}}_{\CP_{Z_M},\CZ}\overset{\on{CT}_{*,\CP_{Z_M},\CZ}^{-,\on{enh}_{\on{co}}}}\longrightarrow \Dmod_{\frac{1}{2}}(\Bun_M)\otimes \Dmod(\CZ)
\end{multline*}
identifies with $\on{CT}_{*,\CP_{Z_M}}\otimes \on{Id}$. 

\medskip

Unwinding, we obtain that the functor $\on{CT}_{*,\CP_{Z_M},\CZ}^{-,\on{enh}_{\on{co}}}$ identifies with the composition
\begin{multline*}
\Dmod_{\frac{1}{2}}(\Bun_G)^{-,\on{enh}_{\on{co}}}_{\CP_{Z_M},\CZ} \overset{\text{\eqref{e:IGP glob twist co}}}\simeq 
\Dmod_{\frac{1}{2}}(\Bun_G)^{-,\on{enh}_{\on{co}}}_{\CZ} \overset{\on{CT}_{*,\CZ}^{-,\on{enh}_{\on{co}}}}\longrightarrow \\
\to \Dmod_{\frac{1}{2}}(\Bun_M)\otimes \Dmod(\CZ) \overset{(\on{transl}_{\CP_{Z_M}})^*\otimes \on{Id}}\longrightarrow 
\Dmod_{\frac{1}{2}}(\Bun_M)\otimes \Dmod(\CZ).
\end{multline*}

\sssec{}

By a similar principle, we define the categories 
$$\Dmod_{\frac{1}{2}}(\Bun_G)^{-,\on{enh}}_{\on{co},\CP_{Z_M},\CZ} \text{ and } \Dmod_{\frac{1}{2}}(\Bun_M)^{-,\on{enh}_{\on{co}}}_{\on{co},\CP_{Z_M},\CZ}$$
and the functors
$$\Eis^{-,\on{enh}_{\on{co}}}_{\on{co},*,\CP_{Z_M},\CZ}:
\Dmod_{\frac{1}{2}}(\Bun_M)^{-,\on{enh}_{\on{co}}}_{\on{co},\CP_{Z_M},\CZ}\to \Dmod_{\frac{1}{2}}(\Bun_G)_{\on{co}}\otimes \Dmod(\CZ),$$
$$\on{CT}^{-,\on{enh}_{\on{co}}}_{\on{co},?,\CP_{Z_M},\CZ}:
\Dmod_{\frac{1}{2}}(\Bun_G)_{\on{co}}\otimes \Dmod(\CZ)\to \Dmod_{\frac{1}{2}}(\Bun_M)_{\on{co},\CP_{Z_M},\CZ}^{-,\on{enh}_{\on{co}}}$$
and 
$$\on{CT}^{-,\on{enh}}_{\on{co},?,\CP_{Z_M},\CZ}:
\Dmod_{\frac{1}{2}}(\Bun_G)_{\on{co},\CP_{Z_M},\CZ}^{-,\on{enh}}\to \Dmod_{\frac{1}{2}}(\Bun_M)_{\on{co}}\otimes \Dmod(\CZ).$$

\sssec{}

Let $\CA$ be a factorization algebra in $\on{I}(G,P^-)^{\on{loc}}$. Using the identification \eqref{e:alpha IGP again}, we produce from 
it a factorization algebra in $\on{I}(G,P^-)^{\on{loc}}_{\CP_{Z_M}}$, to be denoted $\CA_{\CP_{Z_M}}$. 

\medskip

The recipe in Sects. \ref{sss:fact alg Eis} and \ref{sss:fact alg CT} gives rise to functors
$$\Eis^-_{!,\CP_{Z_M},\CA}:\Dmod_{\frac{1}{2}}(\Bun_M)\to \Dmod_{\frac{1}{2}}(\Bun_G)$$
and
$$\on{CT}^-_{\on{co},?,\CP_{Z_M},\CA}:\Dmod_{\frac{1}{2}}(\Bun_G)_{\on{co}}\to \Dmod_{\frac{1}{2}}(\Bun_M)_{\on{co}}.$$

We have
$$\Eis^-_{!,\CP_{Z_M},\CA}\simeq \Eis^-_{!,\CA}\circ (\on{transl}_{\CP_{Z_M}})_* \text{ and }
\on{CT}^-_{\on{co},?,\CP_{Z_M},\CA}\simeq (\on{transl}_{\CP_{Z_M}})^*\circ \on{CT}^{-,\on{enh}}_{\on{co},?,\CA}.$$

\sssec{} \label{sss:CT with fact alg kernel}

Similarly, if $\CB$ is a factorization algebra in $\on{I}(G,P^-)^{\on{loc}}_{\on{co}}$, we construct functors
$$\Eis^-_{\on{co},*,\CP_{Z_M},\CB}:\Dmod_{\frac{1}{2}}(\Bun_M)_{\on{co}}\to \Dmod_{\frac{1}{2}}(\Bun_G)_{\on{co}}$$
and 
$$\on{CT}^-_{*,\CP_{Z_M},\CB}:\Dmod_{\frac{1}{2}}(\Bun_G)\to \Dmod_{\frac{1}{2}}(\Bun_M)$$
and we have 
$$\Eis^-_{\on{co},*,\CP_{Z_M},\CB}\simeq \Eis^{-,\on{enh}_{\on{co}}}_{\on{co},*,\CB}\circ (\on{transl}_{\CP_{Z_M}})_*,$$
$$\on{CT}^-_{*,\CP_{Z_M},\CB}\simeq (\on{transl}_{\CP_{Z_M}})^*\circ \on{CT}^-_{*,\CB}.$$

\sssec{}

Recall the equivalence $\Theta_{\on{I}(G,P^-)^{\on{loc}},\tau}$ of \eqref{e:Theta IGP glob}, and note that in the translated context it gives rise 
to an equivalence
$$\Theta_{\on{I}(G,P^-)^{\on{loc}}_{\CP_{Z_M}},\tau}:
\on{I}(G,P^-)^{\on{loc}}_{\CP_{Z_M},\on{co}}\to \on{I}(G,P^-)^{\on{loc}}_{\tau_M(\CP_{Z_M})}.$$

\medskip

Let $$\Mir^{-,\on{enh}}_{\Bun_M,\tau,\CP_{Z_M}}:
\ul\Dmod_{\frac{1}{2}}(\Bun_M)_{\on{co},\CP_{Z_M}}^{-,\on{enh}_{\on{co}}} \to \ul\Dmod_{\frac{1}{2}}(\Bun_M)^{-,\on{enh}}_{\tau_M(\CP_{Z_M})}$$
denote the tensor product of $\Mir^{-,\on{enh}}_{\Bun_M,\tau}$ and $\Theta_{\on{I}(G,P^-)^{\on{loc}}_{\CP_{Z_M}},\tau}$.

%
%

\medskip 

From \thmref{t:intertwiners} we obtain:

\begin{cor} \label{c:twisted intertwiners}
The following diagram commutes 
\begin{equation} \label{e:twisted intertwiners}
\CD
\Dmod_{\frac{1}{2}}(\Bun_M)_{\on{co},\CP_{Z_M},\Ran}^{-,\on{enh}_{\on{co}}}  @>{\Mir^{-,\on{enh}}_{\Bun_M,\tau,\CP_{Z_M}}}>> 
\Dmod_{\frac{1}{2}}(\Bun_M)^{-,\on{enh}}_{\tau_M(\CP_{Z_M}),\Ran} \\
@V{\Eis^{-,\on{enh}_{\on{co}}}_{\on{co},*,\CP_{Z_M}}[\delta_{N^-_P}]}VV @VV{\Eis^{-,\on{enh}_{\on{co}}}_{!,\tau_M(\CP_{Z_M})}[-\delta_{N^-_P}]}V \\
\Dmod_{\frac{1}{2}}(\Bun_G)_{\on{co}} @>>{\Mir_{\Bun_G,\tau}}> \Dmod_{\frac{1}{2}}(\Bun_G).
\endCD
\end{equation}
\end{cor}

\sssec{}

Let 
$$\Mir^{-,\on{enh}}_{\Bun_G,\tau,\CP_{Z_M}}:\Dmod_{\frac{1}{2}}(\Bun_G)_{\on{co},\CP_{Z_M},\Ran}^{-,\on{enh}} \to \Dmod_{\frac{1}{2}}(\Bun_G)_{\tau_M(\CP_{Z_M}),\Ran}^{-,\on{enh}_{\on{co}}}$$
be the equivalence obtained by tensoring $\Mir^{-,\on{enh}}_{\Bun_G,\tau}$ with $\Theta_{\on{I}(G,P^-)^{\on{loc}}_{\CP_{Z_M}},\tau}$. 

\medskip

By duality, from \corref{c:twisted intertwiners} we obtain: 

\begin{cor} \label{c:intertwiners twisted CT}
The following diagram commutes:
$$
\CD
\Dmod_{\frac{1}{2}}(\Bun_M)_{\on{co}} @>{\Mir_{\Bun_M,\tau}}>> \Dmod_{\frac{1}{2}}(\Bun_M) \\
@A{\on{CT}^{-,\on{enh}}_{\on{co},?,\CP_{Z_M}}[-\delta_{N^-_P}]}AA @AA{\on{CT}_{*,\tau_M(\CP_{Z_M})}^{-,\on{enh}_{\on{co}}}[\delta_{N^-_P}]}A \\
\Dmod_{\frac{1}{2}}(\Bun_G)_{\on{co},\CP_{Z_M},\Ran}^{-,\on{enh}} @>>{\Mir^{-,\on{enh}}_{\Bun_G,\tau,\CP_{Z_M}}}> \Dmod_{\frac{1}{2}}(\Bun_G)_{\tau_M(\CP_{Z_M}),\Ran}^{-,\on{enh}_{\on{co}}}.
\endCD
$$
\end{cor} 

\sssec{}

We will mostly apply the above discussion in the case when $\CP_{Z_M}=\rho_P(\omega_X)$.  Note that in this case
$\tau_M(\CP_{Z_M})=-\rho_P(\omega_X)$.

\section{Global interpretation of the enhanced functors} \label{s:Eis glob}

In this section we will give an interpretation of the functors
$$\Eis^{-,\on{enh}}_!:\Dmod_{\frac{1}{2}}(\Bun_M)^{-,\on{enh}}_\Ran\to \Dmod_{\frac{1}{2}}(\Bun_G),$$
$$\Eis^{-,\on{enh}_{\on{co}}}_{\on{co},*}: \Dmod_{\frac{1}{2}}(\Bun_M)_{\on{co},\Ran}^{-,\on{enh}_{\on{co}}}\to \Dmod_{\frac{1}{2}}(\Bun_G)_{\on{co}},$$
$$\on{CT}^{-,\on{enh}_{\on{co}}}_*:\Dmod_{\frac{1}{2}}(\Bun_G)_{\Ran}^{-,\on{enh}_{\on{co}}}\to \Dmod_{\frac{1}{2}}(\Bun_M),$$
and 
$$\on{CT}^{-,\on{enh}}_{\on{co},?}:\Dmod_{\frac{1}{2}}(\Bun_G)_{\on{co},\Ran}^{-,\on{enh}}\to \Dmod_{\frac{1}{2}}(\Bun_M)_{\on{co}}$$
via Drinfeld's compactifications.

\medskip

This interpretation will be crucial for the computation of Whittaker coefficients of enhanced Eisenstein series. 

\ssec{Drinfeld's compactifications} 

\sssec{}

Let 
$$\Bun_G \overset{\wt\sfp^-}\leftarrow \BunPmt \overset{\wt\sfq^-}\to  \Bun_M$$
be as in \cite[Sect. 1.3.6]{BG1}. Let 
$$j:\Bun_{P^-}\hookrightarrow \BunPmt$$
denote the tautological open embedding. 

\sssec{}

Let 
\begin{equation} \label{e:Dr over Ran}
\BunPmtpR\to \Ran
\end{equation} 
be the polar version of this space, where we allow the generalized $P^-$-reduction to have 
poles at the specified locus of the curve, see \cite[Sect. 6.2.1]{Ga5} (in the case $P=B$)
or \cite[Sect. 4.1.1]{BG1} (pointwise version for any $P$). 

\medskip

We denote by the same symbols the corresponding projections
$$\Bun_G \overset{\wt\sfp^-}\leftarrow  \BunPmtpR \overset{\wt\sfq^-}\to   \Bun_M.$$

\medskip

For $\CZ\to \Ran$ denote
$$\BunPmtpZ:=\BunPmtpR\underset{\Ran}\times \CZ;$$
and let
$$\Bun_G\times \CZ  \overset{\wt\sfp^-_\CZ}\leftarrow \BunPmtpZ \overset{\wt\sfq^-_\CZ}\to   \Bun_M\times \CZ$$
be the corresponding maps.

\sssec{}

We will consider three $\mu_2$-gerbes on $\BunPmtpR$:
$$\det^{\frac{1}{2}}_G:=(\wt\sfp^-)^*(\det^{\frac{1}{2}}_{\Bun_G}), \,\, \det^{\frac{1}{2}}_M:=(\wt\sfq^-)^*(\det^{\frac{1}{2}}_{\Bun_M})$$
and
$$\det^{\frac{1}{2}}_{G,M}:=\det^{\frac{1}{2}}_G\otimes \det^{-\frac{1}{2}}_M\simeq \det^{\frac{1}{2}}_G\otimes \det^{\frac{1}{2}}_M.$$

Note that according to \secref{sss:square root global gerbe}, the restriction $\det^{\frac{1}{2}}_{G,M}|_{\Bun_{P^-}}$ is canonically trivialized. 

\sssec{}  \label{sss:cats on comp begin}

Let $\CG$ be any of the above $\mu_2$-gerbes. For $\CZ\to \Ran$, we 
will consider the following categories of (twisted) D-modules on $\BunPmtpZ$
$$\Dmod_\CG(\BunPmtpZ),\,\,  \Dmod_\CG(\BunPmtpZ)_{\on{co}} \text{ and }  \Dmod_\CG(\BunPmtpZ)_{\on{co}/\Bun_M},$$
defined as follows.  

\medskip

First, we note that \eqref{e:Dr over Ran} is a relative ind-algebraic stack. For a scheme $S\to \Ran$, write 
$$\BunPmtpS\simeq \underset{i}{\on{colim}}\, \CY_{S,i},$$
where for $i\to j$, the maps 
$$f_{i,j}:\CY_{S,i}\to \CY_{S,j}$$
are closed embeddings. (Namely, the closed substacks $\CY_{S,i}\subset \BunPmtpS$ are obtained by bounding the order 
of poles of the generalized $P^-$-reduction.)

\medskip

Note that for every $i$, the map
$$\wt\sfp^-|_{\CY_{S,i}}:\CY_{S,i}\to \Bun_G\times S_i$$
is schematic.

\sssec{}

For every $S$ and $i$ as above, we let $\Dmod_\CG(\CY_{S,i})$ be the usual category of (twisted) D-modules. 

\medskip

Let 
$$\Dmod_\CG(\CY_{S,i})_{\on{co}}:=\underset{U_\CY}{\on{colim}}\, \Dmod_\CG(U_\CY)_{\on{co}},$$
where:

\begin{itemize}

\item $U_\CY$ runs over the poset of quasi-compact open substacks of $\CY_{S,i}$;

\item For an inclusion $U_\CY\overset{j}\hookrightarrow U'_\CY$, the functor $\Dmod_\CG(U_\CY)_{\on{co}}\to  \Dmod_\CG(U'_\CY)_{\on{co}}$
is $j_*$.

\end{itemize}

\sssec{}

We let 
$$\Dmod_\CG(\CY_{S,i})_{\on{co}/\Bun_M}:=\underset{U_M}{\on{lim}}\, \Dmod_\CG(\CY_{S,i}\underset{\Bun_M}\times U_M)_{\on{co}},$$
where:

\begin{itemize}

\item $U_M$ runs over the poset of quasi-compact open substacks of $\Bun_M$;

\item For an inclusion $U_M\overset{j}\hookrightarrow U'_M$, the functor 
$$\Dmod_\CG(\CY_{S,i}\underset{\Bun_M}\times U'_M)_{\on{co}}\to \Dmod_\CG(\CY_{S,i}\underset{\Bun_M}\times U_M)_{\on{co}}$$
is given by restriction.

\end{itemize}

\sssec{}

We let $\Dmod_\CG(\BunPmtpS)$ be the usual category of twisted D-modules, i.e., 
$$\Dmod_\CG(\BunPmtpS)=\underset{i}{\on{lim}}\, \Dmod_\CG(\CY_{S,i}),$$
and we let 
$$\Dmod_\CG(\BunPmtpS)_{\on{co}}:=\underset{i}{\on{lim}}\, \Dmod_\CG(\CY_{S,i})_{\on{co}},$$
$$\Dmod_\CG(\BunPmtpS)_{\on{co}/\Bun_M}:=\underset{i}{\on{lim}}\, \Dmod_\CG(\CY_{S,i})_{\on{co}/\Bun_M},$$
where in all three cases the transition functors are given by $f_{i,j}^!$. 

\medskip

Note that in all three cases we can rewrite the limits as colimits 
$$\underset{i}{\on{colim}}\, \Dmod_\CG(\CY_{S,i}),$$
$$\underset{i}{\on{colim}}\, \Dmod_\CG(\CY_{S,i})_{\on{co}}$$
and 
$$\underset{i}{\on{colim}}\, \Dmod_\CG(\CY_{S,i})_{\on{co}/\Bun_M},$$
respectively, where where the transition functors are given by $(f_{i,j})_*$. 

\sssec{} \label{sss:cats on comp end}

We let $\Dmod_\CG(\BunPmtpZ)$ be the usual category of twisted D-modules. I.e., 
$$\Dmod_\CG(\BunPmtpZ) \simeq \underset{S\to\CZ}{\on{lim}}\, \Dmod_\CG(\BunPmtpS),$$
and we let
$$\Dmod_\CG(\BunPmtpZ)_{\on{co}} \simeq \underset{S\to\CZ}{\on{lim}}\, \Dmod_\CG(\BunPmtpS)_{\on{co}},$$ 
$$\Dmod_\CG(\BunPmtpZ)_{\on{co}/\Bun_M} \simeq \underset{S\to\CZ}{\on{lim}}\, \Dmod_\CG(\BunPmtpS)_{\on{co}/\Bun_M},$$ 
where the transition functors for $g:S'\to S$ are given by $g^!$.

\medskip

If $\CZ$ is pseudo-proper (see \cite[Sect. C.4]{GLC2} for what this means), e.g., $\CZ=\Ran$ itself, we can choose $S$ to be
proper schemes. In this case the above limits can also be written as colimits 
$$\underset{S\to\CZ}{\on{colim}}\, \Dmod_\CG(\BunPmtpS),$$
$$\underset{S\to\CZ}{\on{colim}}\, \Dmod_\CG(\BunPmtpS)_{\on{co}}$$
and  
$$\underset{S\to\CZ}{\on{colim}}\, \Dmod_\CG(\BunPmtpS)_{\on{co}/\Bun_M},$$ 
where the transition functors are given by $g_*$.

\medskip

Write also
$$\Dmod_{\CG_G}(\Bun_G\times \CZ)_{\on{co}}:=
\underset{S\to\CZ}{\on{lim}}\, \Dmod_{\CG_G}(\Bun_G\times S)_{\on{co}}\simeq \Dmod_{\CG_G}(\Bun_G)_{\on{co}}\otimes \Dmod(\CZ)$$
(for a gerbe $\CG_G$ on $\Bun_G$) and 
$$\Dmod_{\CG_M}(\Bun_M\times \CZ)_{\on{co}}:=
\underset{S\to\CZ}{\on{lim}}\, \Dmod_{\CG_M}(\Bun_M\times S)_{\on{co}}\simeq \Dmod_{\CG_M}(\Bun_M)_{\on{co}}\otimes \Dmod(\CZ)$$
(for a gerbe $\CG_G$ on $\Bun_M$). 

\sssec{} \label{sss:cats on BunPt1} 

For a gerbe $\CG_G$ on $\Bun_G$, we have well-defined functors
\begin{equation} \label{e:p^! usual}
(\sfp^-_\CZ)^!:\Dmod_{\CG_G}(\Bun_G\times \CZ)\to \Dmod_{(\sfp^-)^*(\CG_G)}(\BunPmtpZ),
\end{equation}
\begin{equation} \label{e:p^! co}
(\sfp^-_\CZ)^!:\Dmod_{\CG_G}(\Bun_G\times \CZ)_{\on{co}}\to \Dmod_{(\sfp^-)^*(\CG_G)}(\BunPmtpZ)_{\on{co}},
\end{equation}
and 
\begin{equation} \label{e:p_* co}
(\sfp^-_\CZ)_*:\Dmod_{(\sfp^-)^*\CG_G}(\BunPmtpZ)_{\on{co}}\to \Dmod_{\CG_G}(\Bun_G\times \CZ)_{\on{co}}.
\end{equation}

\medskip

For a gerbe $\CG_M$ on $\Bun_M$, we have well-defined functors
\begin{equation} \label{e:q_* co}
(\sfq^-_\CZ)_*:\Dmod_{(\sfq^-)^*(\CG_M)}(\BunPmtpR)_{\on{co}}\to \Dmod_{\CG_M}(\Bun_M\times \CZ)_{\on{co}}
\end{equation}
and 
\begin{equation} \label{e:q_* rel}
(\sfq^-_\CZ)_*:\Dmod_{(\sfq^-)^*(\CG_M)}(\BunPmtpR)_{\on{co}/\Bun_M}\to \Dmod_{\CG_M}(\Bun_M\times \CZ).
\end{equation}

\sssec{}  \label{sss:cats on BunPt2} 

For a pair of gerbes $\CG_1$ and $\CG_2$ with
$$\CG_{1,2}:=\CG_1\otimes \CG_2,$$ 
the operation of $\sotimes$ product of twisted D-modules gives rise to functors
\begin{equation} \label{e:sotimes co 1}
\Dmod_{\CG_1}(\BunPmtpZ)\otimes \Dmod_{\CG_2}(\BunPmtpZ)_{\on{co}}\to \Dmod_{\CG_{1,2}}(\BunPmtpZ)_{\on{co}}
\end{equation}
and 
\begin{equation} \label{e:sotimes co 2}
\Dmod_{\CG_1}(\BunPmtpZ)\otimes \Dmod_{\CG_2}(\BunPmtpZ)_{\on{co}/\Bun_M}\to \Dmod_{\CG_{1,2}}(\BunPmtpZ)_{\on{co}/\Bun_M}.
\end{equation}

\sssec{} \label{sss:cats on BunPt3} 

In addition, for a gerbe $\CG_M$ on $\Bun_M$ we have a functor
\begin{multline} \label{e:sotimes rel}
(\sfq^-_\CZ)^!(-)\sotimes (-): \Dmod_{\CG_M}(\Bun_M\times \CZ)_{\on{co}}\otimes \Dmod_{\CG}(\BunPmtpZ)_{\on{co}/\Bun_M}\to \\
\to \Dmod_{(\sfq^-)^*(\CG_M)\otimes \CG}(\BunPmtpZ)_{\on{co}}.
\end{multline} 

\ssec{Local-to-global construction for Drinfeld's compactification, the *-version}

\sssec{}

Let 
$$\Gr_{G,\Bun_M,\Ran}:=\Gr_{G,\Ran}\wt\times \Bun_M$$
be the fibration over $\Bun_M\times \Ran$ that associates to a point
$$(\CP_M,\ul{x})\in \Bun_M\times \Ran$$
the twist $\Gr_{G,\CP_M,\ul{x}}$, 
i.e., we twist $\Gr_{G,\ul{x}}$ by the $\fL^+(G)_{\ul{x}}$-torsor induced by the \emph{canonical}\footnote{The splitting $M\to P^-$ depends on the choice of $P$,
and hence so does the map $\on{pt}/M\to \on{pt}/P^-$ as \emph{pointed} stacks. However, the map $\on{pt}/M\to \on{pt}/P^-$ as just stacks is canonically
independent of the choice of $P$.} map
\begin{equation} \label{e:induce G bundle}
\on{pt}/M\to \on{pt}/P^-\to \on{pt}/G
\end{equation} 
from the $\fL^+(M)_{\ul{x}}$-torsor $\CP_M|_{\cD_{\ul{x}}}$. 

\sssec{}

Let $\wt\sfs_\Ran$ denote the tautological projection
$$\Gr_{G,\Bun_M,\Ran}\to \fL^+(P^-)_\Ran\backslash \Gr_{G,\Ran}.$$

We will denote by $\wt\sfs^!_\Ran$ the functor 
$$\CF\mapsto \CF\wt\boxtimes \omega_{\Bun_M}, \quad\Dmod(\Gr_{G,\Ran})^{\fL^+(P^-)_\Ran}\to \Dmod(\Gr_{G,\Bun_M,\Ran}).$$

We will use the same symbol $\wt\sfs^!_\Ran$ for the corresponding functor in the context of twisted D-modules. 

\sssec{}

By the Beauville-Laszlo theorem, we can interpret $\Gr_{G,\Bun_M,\Ran}$ as the moduli space of quadruples
$$(\ul{x},\CP_G,\CP_M,\alpha),$$
where:

\begin{itemize}

\item $\ul{x}\in \Ran$;

\item $\CP_M$ is an $M$-bundle on $X$ (we will denote by $\CP^0_{G,M}$ the induced $G$-bundle by means of \eqref{e:induce G bundle});  

\item $\CP_G$ is a $G$-bundle on $X$;

\item $\alpha$ is an isomorphism between $\CP_G$ and $\CP^0_{G,M}$ over $X-\ul{x}$.

\end{itemize}

\sssec{} \label{sss:Gr to BunPt}

From the above interpretation of $\Gr_{G,\Bun_M,\Ran}$, and the fact that $\CP^0_{G,M}$ has a tautological reduction 
to $P^-$, we obtain that there exists a canonically defined map
$$\pi_\Ran:\Gr_{G,\Bun_M,\Ran}\to \BunPmtpR,$$
which commutes with the projections of both spaces to $\Ran$, $\Bun_G$ and $\Bun_M$, respectively. 

\medskip

Note also that as in \eqref{e:ident gerbes Hecke P}, we have a canonical identification of $\mu_2$-gerbes on $\Gr_{G,\Bun_M,\Ran}$
$$\pi_\Ran^*(\det^{\frac{1}{2}}_{G,M})\simeq \wt\sfs^*_\Ran(\det^{\frac{1}{2}}_{\Gr_{G,\Ran}}).$$

\sssec{}

Let $\CZ\to \Ran$. We will denote by $\pi_\CZ$ the corresponding map
$$\Gr_{G,\Bun_M,\CZ}\to \BunPmtpZ.$$

\medskip 

We claim that there is a canonically defined functor 
\begin{equation} \label{e:IGP co to tilde Ran}
\on{I}(G,P^-)^{\on{loc}}_{\on{co},\CZ} \to  \Dmod_{\det^{\frac{1}{2}}_{G,M}}(\BunPmtpZ)_{\on{co}/\Bun_M},
\end{equation} 
to be denoted $(\pi_*\circ \wt\sfs^!)^{\on{Av}}_\CZ$. 

\medskip

To simplify the notation, instead of \eqref{e:IGP co to tilde Ran}, we will describe its version when $\CZ=\on{pt}$ corresponding to $\ul{x}\in \Ran$: 
\begin{equation} \label{e:IGP co to tilde x}
(\pi_*\circ \wt\sfs^!)^{\on{Av}}_{\ul{x}}: \on{I}(G,P^-)^{\on{loc}}_{\on{co},\ul{x}} \to  \Dmod_{\det^{\frac{1}{2}}_{G,M}}(\BunPmtpx)_{\on{co}/\Bun_M}.
\end{equation} 

\noindent {\it Warning:} For the construction below, it is essential that $\ul{x}\neq \emptyset$.

\begin{rem}

The construction below could be coarsened and considered to take values in 
be given with values in $\Dmod_{\det^{\frac{1}{2}}_{G,M}}(\BunPmtpx)$, instead of 
$\Dmod_{\det^{\frac{1}{2}}_{G,M}}(\BunPmtpx)_{\on{co}/\Bun_M}$. However, we need the finer versions
(i.e., with values in $\Dmod_{\det^{\frac{1}{2}}_{G,M}}(\BunPmtpx)_{\on{co}/\Bun_M}$) in what follows. 

\medskip

Note also that we could \emph{not} further refine it to take values in $\Dmod_{\det^{\frac{1}{2}}_{G,M}}(\BunPmtpx)_{\on{co}}$:
indeed objects that we construct exhibit a $\Dmod(-)$-behavior (rather than a $\Dmod(-)_{\on{co}}$-behavior) along
$\Bun_M$.

\end{rem} 

\sssec{}

We start with the functor
$$\Dmod_{\det^{\frac{1}{2}}_{\Gr_{G,\Ran}}}(\Gr_{G,\ul{x}})^{\fL^+(P^-)}\overset{\wt\sfs^!_{\ul{x}}}\to 
\Dmod_{\wt\sfs^*(\det^{\frac{1}{2}}_{\Gr_{G,\Ran}})}(\Gr_{G,\Bun_M,\ul{x}})=
\Dmod_{\pi^*(\det^{\frac{1}{2}}_{G,M})}(\Gr_{G,\Bun_M,\ul{x}}).$$

Composing with the functor
$$(\pi_{\ul{x}})_*:\Dmod_{\wt\sfs^*(\det^{\frac{1}{2}}_{\Gr_{G,\Ran}})}(\Gr_{G,\Bun_M,\ul{x}})\to  
\Dmod_{\det^{\frac{1}{2}}_{G,M}}(\BunPmtpx)_{\on{co}/\Bun_M}$$
we obtain a functor
\begin{equation} \label{e:before averaging}
(\pi_{\ul{x}})_*\circ \wt\sfs_{\ul{x}}^!:\Dmod_{\det^{\frac{1}{2}}_{\Gr_{G,\Ran}}}(\Gr_{G,\ul{x}})^{\fL^+(P^-)}\to  
\Dmod_{\det^{\frac{1}{2}}_{G,M}}(\BunPmtpx)_{\on{co}/\Bun_M}. 
\end{equation}

\sssec{} \label{sss:push Av 1}

We now apply to \eqref{e:before averaging} the procedure analogous to one in \cite[Sect. 9.5.4]{GLC2}. Namely, let
$N^\alpha$ be sequence of group subschemes of $\fL(N^-_P)$, normalized by $\fL^+(P^-)$. 

\medskip

For a pair of indices $\alpha'\geq \alpha$ we have a natural transformation 
\begin{equation} \label{e:change index} 
(\pi_{\ul{x}})_*\circ \wt\sfs^!_{\ul{x}}\circ \on{Av}^{N^{\alpha'}}_*\to (\pi_{\ul{x}})_*\circ \wt\sfs_{\ul{x}}^!\circ \on{Av}^{N^{\alpha}}_*.
\end{equation}

As in \cite[Sect. 9.5.4]{GLC2} one shows that \eqref{e:change index} is an isomorphism for $\alpha$ large enough.\footnote{The index $\alpha$ for
which the stabilization occurs depends on the connected component of $\Bun_M$.}

\sssec{} \label{sss:push Av 2}

This implies that for any $\alpha$ that is large enough, the functor 
$$(\pi_{\ul{x}})_*\circ \wt\sfs_{\ul{x}}^!\circ \on{Av}^{N^{\alpha}}_*:\Dmod_{\det^{\frac{1}{2}}_{\Gr_{G,\Ran}}}(\Gr_{G,\ul{x}})^{\fL^+(P^-)}\to  
\Dmod_{\det^{\frac{1}{2}}_{G,M}}(\BunPmtpx)_{\on{co}/\Bun_M}$$
factors as
\begin{multline} \label{e:IGP co to tilde x pre}
\Dmod_{\det^{\frac{1}{2}}_{\Gr_{G,\Ran}}}(\Gr_{G,\ul{x}})^{\fL^+(P^-)}\simeq 
\Dmod_{\det^{\frac{1}{2}}_{\Gr_{G,\Ran}}}(\Gr_{G,\ul{x}})^{\fL^+(M)}_{\fL^+(N^-_P)}\twoheadrightarrow \\
\twoheadrightarrow \Dmod_{\det^{\frac{1}{2}}_{\Gr_{G,\Ran}}}(\Gr_{G,\ul{x}})^{\fL^+(M)}_{\fL(N^-_P)}\to \Dmod_{\det^{\frac{1}{2}}_{G,M}}(\BunPmtpx)_{\on{co}/\Bun_M}.
\end{multline}

\medskip

Precomposing the last arrow in \eqref{e:IGP co to tilde x pre} with
$$\on{I}(G,P^-)^{\on{loc}}_{\on{co},\ul{x}}\to \Dmod_{\det^{\frac{1}{2}}_{\Gr_{G,\Ran}}}(\Gr_{G,\ul{x}})^{\fL^+(M)}_{\fL(N^-_P)},$$
we obtain the sought-for functor $(\pi_*\circ \wt\sfs^!)^{\on{Av}}_{\ul{x}}$ of \eqref{e:IGP co to tilde x}. 

\sssec{} \label{sss:* BunPt unital 1}

Note that we can view the assignment 
$$\CZ\rightsquigarrow \Dmod_{\det^{\frac{1}{2}}_{G,M}}(\BunPmtpZ)_{\on{co}/\Bun_M}$$
as a unital crystal of categories over $\Ran$. 

\medskip

Namely, for $\ul{x}_1\subset \ul{x}_2$, the corresponding functor 
$$\Dmod_{\det^{\frac{1}{2}}_{G,M}}(\widetilde{\Bun}_{P^-,\ul{x}_1})_{\on{co}/\Bun_M}\to
\Dmod_{\det^{\frac{1}{2}}_{G,M}}(\widetilde{\Bun}_{P^-,\ul{x}_2})_{\on{co}/\Bun_M}$$
is given by direct image along the closed embedding
$$\widetilde{\Bun}_{P^-,\ul{x}_1}\to \widetilde{\Bun}_{P^-,\ul{x}_2}.$$

\sssec{} \label{sss:* BunPt unital 2}

By construction, the assignment
$$\CZ\rightsquigarrow (\pi_*\circ \wt\sfs^!)^{\on{Av}}_\CZ$$
has a structure of \emph{left}-lax unital structure\footnote{See \cite[Sect. I.1.2]{GLC2} what this means.}
as a functor between crystals of categories over $\Ran$. Denote this functor by
$\ul{(\pi_*\circ \wt\sfs^!)}^{\on{Av}}$. 

\medskip

However, we claim: 

\begin{prop} \label{p:* BunPt unital} 
The left-lax unital structure on $\ul{(\pi_*\circ \wt\sfs^!)}^{\on{Av}}$ is strict.
\end{prop} 

\ssec{The strong approximation argument}

In this subsection we will prove \propref{p:* BunPt unital}.

\medskip

As we shall see, the proof essentially amounts to the strong
approximation theorem applied to the \emph{unipotent} group $N^-_P$,
cf. \cite[Sect. 5.2.7 and Proposition 5.2.8]{Ga6}. 

\sssec{}

To simplify the notation, we will prove the proposition for a fixed 
pair of points $\ul{x}\subseteq \ul{x}'$ of $\Ran$. 

\medskip

Further, let us work over $1\in \Bun_M$, so that instead of $\widetilde{\Bun}_{P^-,\ul{x}}$ we are dealing 
with 
$$\widetilde{\Bun}_{N_P^-,\ul{x}}:=\widetilde{\Bun}_{P^-,\ul{x}}\underset{\Bun_M}\times \on{pt}$$
and instead of $\Gr_{G,\Bun_M,\ul{x}}$ just with $\Gr_{G,\ul{x}}$, and similarly for $\ul{x}'$. 

\sssec{}

Let $N^-_{P,\on{out},\ul{x}}$ be the group ind-scheme $\bMaps(X-\ul{x},N^-_P)$.
Let $$\fL^+(N^-_P)_{\ul{x}} \subset N^\alpha_{\ul{x}}\subset \fL(N^-_P)_{\ul{x}}$$ be a group sub\emph{scheme}
large enough so that 
$$N^-_{P,\on{out},\ul{x}}\cdot N^\alpha = \fL(N^-_P)_{\ul{x}}.$$

\medskip

By definition, the functor $(\pi_*\circ \wt\sfs^!)^{\on{Av}}_{\ul{x}}$ equals
$$(\pi_{\ul{x}})_*\circ \on{Av}_*^{N^\alpha_{\ul{x}}},$$

\medskip

Let $N^\alpha_{\ul{x}'}\subset \fL(N^-_P)_{\ul{x}'}$ be a suitable subgroup for $\ul{x}'$. 

\medskip

So,
$$(\pi_*\circ \wt\sfs^!)^{\on{Av}}_{\ul{x}'}=(\pi_{\ul{x}'})_*\circ \on{Av}_*^{N^\alpha_{\ul{x}'}}.$$

\sssec{} \label{sss:choice of subgroups}

Let 
$$N^\alpha_{\ul{x}\subseteq \ul{x}'}:=N^\alpha_{\ul{x}'}\cap \bMaps(\cD_{\ul{x}'}-\ul{x},N^-_P)\subset N^\alpha_{\ul{x}'}.$$

With no restriction of generality, we can assume that the subgroups $N^\alpha_{\ul{x}}$ and $N^\alpha_{\ul{x}'}$
are chosen so that $N^\alpha_{\ul{x}\subseteq \ul{x}'}$ equals the preimage of $N^\alpha_{\ul{x}}$ under the projection
$$\bMaps(\cD_{\ul{x}'}-\ul{x},N^-_P)\to \bMaps(\cD_{\ul{x}}-\ul{x},N^-_P)=\fL(N^-_P)_{\ul{x}}.$$

\sssec{}

The left-unital structure on $\ul{(\pi_*\circ \wt\sfs^!)}^{\on{Av}}$ is given by the system of maps
\begin{multline} \label{e:LUS}
(\pi_*\circ \wt\sfs^!)^{\on{Av}}_{\ul{x}'} \circ (\on{ins.unit}_{\ul{x}\subseteq \ul{x}'})_*= \\
= (\pi_{\ul{x}'})_*\circ \on{Av}_*^{N^\alpha_{\ul{x}'}} \circ (\on{ins.unit}_{\ul{x}\subseteq \ul{x}'})_* 
\to (\pi_{\ul{x}'})_*\circ \on{Av}_*^{N^\alpha_{\ul{x}\subseteq \ul{x}'}} \circ (\on{ins.unit}_{\ul{x}\subseteq \ul{x}'})_* 
\simeq \\
\simeq (\pi_{\ul{x}'})_*\circ  (\on{ins.unit}_{\ul{x}\subseteq \ul{x}'})_* \circ \on{Av}_*^{N^\alpha_{\ul{x}}}  \simeq
(\pi_{\ul{x}})_*\circ \on{Av}_*^{N^\alpha_{\ul{x}}}
=(\pi_*\circ \wt\sfs^!)^{\on{Av}}_{\ul{x}},
\end{multline}
where:

\begin{itemize}

\item $\on{ins.unit}_{\ul{x}\subset \ul{x}'}:\Gr_{G,\ul{x}}\to \Gr_{G,\ul{x}'}$ is the unital structure on $\Gr_G$;

\medskip

\item The next-to-last isomorphism is thanks to the choice of subgroups in \secref{sss:choice of subgroups}.

\end{itemize}

\medskip

We wish to show that the maps \eqref{e:LUS} are isomorphisms. 
I.e., we wish to show that the second arrow in \eqref{e:LUS}
is an isomorphism. 

\sssec{}

Denote 
$$N^\alpha_{\on{out},\ul{x}}:=N^-_{P,\on{out},\ul{x}}\cap N^\alpha,$$
and similarly for $\ul{x}'$.

\medskip

We factor the second arrow in \eqref{e:LUS} as
\begin{multline} \label{e:LUS 1}
(\pi_{\ul{x}'})_*\circ \on{Av}_*^{N^\alpha_{\ul{x}'}} \circ (\on{ins.unit}_{\ul{x}\subseteq \ul{x}'})_* 
\to (\pi_{\ul{x}'})_*\circ \on{Av}_*^{N^\alpha_{\on{out},\ul{x}'}} \circ \on{Av}_*^{N^\alpha_{\ul{x}\subseteq \ul{x}'}} \circ (\on{ins.unit}_{\ul{x}\subseteq \ul{x}'})_* \to \\
\to (\pi_{\ul{x}'})_* \circ \on{Av}_*^{N^\alpha_{\ul{x}\subseteq \ul{x}'}} \circ (\on{ins.unit}_{\ul{x}\subseteq \ul{x}'})_*.
\end{multline} 

We claim that both arrows in \eqref{e:LUS 1} are isomorphisms when the subgroups $N^\alpha_{\ul{x}}$ and $N^\alpha_{\ul{x}'}$
are large enough. 

\sssec{}

Up to enlarging $N^\alpha_{\ul{x}}$ and $N^\alpha_{\ul{x}'}$, we can assume that 
\begin{equation} \label{e:subgroups large}
N^\alpha_{\on{out},\ul{x}'}\cdot N^\alpha_{\ul{x}\subseteq \ul{x}'}=N^\alpha_{\ul{x}'},
\end{equation} 
which is a form of \emph{strong approximation} for $N^-_P$. 

\medskip 

The equality \eqref{e:subgroups large} implies that the first arrow in \eqref{e:LUS 1} is an isomorphism. 

\sssec{}

We now prove that the second arrow in \eqref{e:LUS 1} is an isomorphism. We rewrite it equivalently as 
$$(\pi_{\ul{x}'})_*\circ \on{Av}_*^{N^\alpha_{\on{out},\ul{x}'}} \circ  (\on{ins.unit}_{\ul{x}\subseteq \ul{x}'})_* \circ \on{Av}_*^{N^\alpha_{\ul{x}}}\to
(\pi_{\ul{x}'})_*\circ (\on{ins.unit}_{\ul{x}\subseteq \ul{x}'})_* \circ \on{Av}_*^{N^\alpha_{\ul{x}}}.$$

In fact, we claim that the map
\begin{equation} \label{e:LUS 2}
(\pi_{\ul{x}'})_*\circ \on{Av}_*^{N^\alpha_{\on{out},\ul{x}'}} \circ  (\on{ins.unit}_{\ul{x}\subseteq \ul{x}'})_*\to 
(\pi_{\ul{x}'})_*\circ (\on{ins.unit}_{\ul{x}\subseteq \ul{x}'})_* 
\end{equation} 
is an isomorphism when evaluated on $N^\alpha_{\on{out},\ul{x}}$-equivariant objects. 

\sssec{}

Indeed, we note that the map 
$$\pi_{\ul{x}}:\Gr_{G,\ul{x}}\to \widetilde{\Bun}_{N_P^-,\ul{x}}$$
factors as
\begin{equation} \label{e:factor via out}
\Gr_{G,\ul{x}}\to  N^-_{P,\on{out},\ul{x}}\backslash \Gr_{G,\ul{x}}\overset{\ol\pi_{\ul{x}}}\longrightarrow \widetilde{\Bun}_{N_P^-,\ul{x}},
\end{equation} 
and similarly for $\ul{x}'$.

\medskip

The isomorphism \eqref{e:LUS 2} follows now from the next commutative diagram:

\begin{equation} \label{e:LUS 3}
\CD
N^\alpha_{\on{out},\ul{x}'}\backslash\Biggl(\underset{N^\alpha_{\on{out},\ul{x}}}{\underline{N^\alpha_{P,\on{out},\ul{x}'}\times \Gr_{G,\ul{x}}}}\Biggr)
@>{\on{action}\circ \on{ins.unit}_{\ul{x}\subseteq \ul{x}'}}>> N^\alpha_{\on{out},\ul{x}'}\backslash\Gr_{G,\ul{x}'}\\
@V{\sim}VV    \\
N^\alpha_{\on{out},\ul{x}}\backslash \Gr_{G,\ul{x}} & &  @VV{\ol\pi_{\ul{x}'}}V  \\
@V{\ol\pi_{\ul{x}}}VV   \\
\widetilde{\Bun}_{N_P^-,\ul{x}} @>>> \widetilde{\Bun}_{N_P^-,\ul{x}'}.
\endCD
\end{equation} 

Indeed, the left- (resp. right-) hand side in \eqref{e:LUS 2} is given by pushforward along the clockwise
(resp., counterclockwise) circuit in \eqref{e:LUS 3}.  

\ssec{The global enhanced functors, the *-version}

\sssec{} \label{sss:pre CT enh}

For $\CZ\to \Ran$, define the functor 
\begin{multline*} 
\on{pre-CT}^{-,\on{enh}_{\on{co}},\on{glob}}_{*,\CZ}:\Dmod_{\frac{1}{2}}(\Bun_G) \otimes 
\on{I}(G,P^-)^{\on{loc}}_{\on{co},\CZ} \to \Dmod_{\frac{1}{2}}(\Bun_M\times \CZ)\simeq \\
\simeq \Dmod_{\frac{1}{2}}(\Bun_M)\otimes \Dmod(\CZ)
\end{multline*} 
to be 
\begin{equation} \label{e:defn of CT glob}
(\wt\sfq^-_\CZ)_*\left((\wt\sfp^-)^!(-)\sotimes (\pi_*\circ \wt\sfs^!)_\CZ^{\on{Av}}(-)\right)[\on{shift}],
\end{equation} 
where: 

\begin{itemize}

\item The functor $(\wt\sfp^-)^!$ is understood in the sense \eqref{e:p^! usual}; 

\smallskip

\item The functor $\sotimes$ is understood in the sense \eqref{e:sotimes co 2}; 

\smallskip

\item The functor $(\wt\sfq^-_\CZ)_*$ is understood in the sense \eqref{e:q_* rel}; 

\smallskip
 
\item $[\on{shift}]$ is as in \eqref{e:global shift}.

\end{itemize}

\medskip

\noindent NB: the sign of the shift is opposite to that in \eqref{e:CT * defn}.

\sssec{}

Unwinding the constructions, we obtain that the functor $\on{pre-CT}^{-,\on{enh}_{\on{co}},\on{glob}}_{*,\CZ}$ factors canonically as 
\begin{multline} \label{e:pre CT enh glob}
\Dmod_{\frac{1}{2}}(\Bun_G) \otimes \on{I}(G,P^-)^{\on{loc}}_{\on{co},\CZ} \simeq \\
\simeq \left(\Dmod_{\frac{1}{2}}(\Bun_G) \otimes \Dmod(\CZ)\right)\underset{\Dmod(\CZ)}\otimes \on{I}(G,P^-)^{\on{loc}}_{\on{co},\CZ}\to \\
\to \left(\Dmod_{\frac{1}{2}}(\Bun_G) \otimes \Dmod(\CZ)\right)\underset{\Sph_{G,\CZ}}\otimes \on{I}(G,P^-)^{\on{loc}}_{\on{co},\CZ}= \\
=\Dmod_{\frac{1}{2}}(\Bun_G)_{\CZ}^{-,\on{enh}_{\on{co}}} \to \Dmod_{\frac{1}{2}}(\Bun_M)\otimes \Dmod(\CZ).
\end{multline} 

Denote the last arrow in \eqref{e:pre CT enh glob} by 
$$\on{CT}^{-,\on{enh}_{\on{co}},\on{glob}}_{*,\CZ}: \Dmod_{\frac{1}{2}}(\Bun_G)_{\CZ}^{-,\on{enh}_{\on{co}}} \to \Dmod_{\frac{1}{2}}(\Bun_M)\otimes \Dmod(\CZ).$$

\medskip 

By construction, the functor $\on{CT}^{-,\on{enh}_{\on{co}},\on{glob}}_{*,\CZ}$ intertwines the actions of $\Sph_{M,\CZ}$ on the two sides. 

\medskip

Denote also
$$\on{CT}^{-,\on{enh}_{\on{co}},\on{glob}}_*:=(\on{Id}\otimes \on{C}^\cdot(\Ran,-))\circ \on{CT}^{-,\on{enh}_{\on{co}},\on{glob}}_{*,\Ran},\quad 
\Dmod_{\frac{1}{2}}(\Bun_G)_{\Ran}^{-,\on{enh}_{\on{co}}} \to \Dmod_{\frac{1}{2}}(\Bun_M).$$
 
\sssec{}

We claim:

\begin{prop} \label{p:loc vs global CT enh co}
There is a canonical isomorphism 
$$\on{CT}^{-,\on{enh}_{\on{co}},\on{glob}}_{*,\CZ} \simeq \on{CT}^{-,\on{enh}_{\on{co}}}_{*,\CZ},$$
as $\Sph_{M,\CZ}$-linear functors
$$\Dmod_{\frac{1}{2}}(\Bun_G)_{\CZ}^{-,\on{enh}_{\on{co}}} \rightrightarrows \Dmod_{\frac{1}{2}}(\Bun_M)\otimes \Dmod(\CZ).$$
\end{prop}

\begin{proof}

In order to simplify the notation we will establish the required isomorphism when $\CZ=\on{pt}$, i.e., 
\begin{equation} \label{e:loc vs global CT enh co}
\on{CT}^{-,\on{enh}_{\on{co}},\on{glob}}_{*,\ul{x}} \simeq \on{CT}^{-,\on{enh}_{\on{co}}}_{*,\ul{x}}
\end{equation}
as $\Sph_{M,\ul{x}}$-linear functors
$$\Dmod_{\frac{1}{2}}(\Bun_G)_{\ul{x}}^{-,\on{enh}_{\on{co}}} \rightrightarrows \Dmod_{\frac{1}{2}}(\Bun_M)$$
for $\ul{x}\in \Ran$.

\medskip

Note that there is a naturally defined map
$$\sff_{\ul{x}}:\on{Hecke}^{\on{glob}}_{G,P^-,\ul{x}}\to \BunPmtpx,$$
which respects the projections of both spaces to $\Bun_G$ and $\Bun_M$, respectively. In addition, we have
$$\sff_{\ul{x}}^*(\det^{\frac{1}{2}}_{G,M})\simeq \sfs_{\ul{x}}^*(\det^{\frac{1}{2}}_{\Gr_{G,\ul{x}}}),$$
where $\sfs_{\ul{x}}$ is the map from \eqref{e:map s x}. 

\medskip

Consider the functor
$$(\sff_{\ul{x}})_*\circ \sfs_{\ul{x}}^!:\Dmod_{\frac{1}{2}}(\Gr_{G,\ul{x}})^{\fL^+(P^-)_{\ul{x}}}\to  \Dmod_{\det^{\frac{1}{2}}_{G,M}}(\BunPmtpx)_{\on{co}/\Bun_M}.$$

\medskip

As in Sects. \ref{sss:push Av 1}-\ref{sss:push Av 2}, from $(\sff_{\ul{x}})_*\circ \sfs_{\ul{x}}^!$ we obtain a functor
$$((\sff_{\ul{x}})_*\circ \sfs_{\ul{x}}^!)^{\on{Av}}:\on{I}(G,P^-)^{\on{loc}}_{\on{co},\ul{x}}\to \Dmod_{\det^{\frac{1}{2}}_{G,M}}(\BunPmtpx)_{\on{co}/\Bun_M}.$$

We claim:

\begin{lem} \label{l:approx}
There is a canonical isomorphism 
$$((\sff_{\ul{x}})_*\circ \sfs_{\ul{x}}^!)^{\on{Av}}\simeq (\pi_*\circ \wt\sfs^!)^{\on{Av}}_{\ul{x}}[2\cdot \on{shift}].$$
\end{lem} 

The proof will be given in \secref{ss:proof of approx}. 

\medskip

The isomorphism \eqref{e:loc vs global CT enh co} follows from \lemref{l:approx} by the projection formula.

\end{proof}

\begin{cor}  \label{c:loc vs global CT enh co}
There is a canonical isomorphism 
$$\on{CT}^{-,\on{enh}_{\on{co}},\on{glob}}_* \simeq \on{CT}^{-,\on{enh}_{\on{co}}}_*,$$
as functors 
$$\Dmod_{\frac{1}{2}}(\Bun_G)_\Ran^{-,\on{enh}_{\on{co}}} \rightrightarrows \Dmod_{\frac{1}{2}}(\Bun_M).$$
\end{cor}

\sssec{}

Define the functor 
\begin{multline*}
\on{pre-Eis}^{-,\on{enh}_{\on{co}},\on{glob}}_{\on{co},*,\CZ}: \Dmod_{\frac{1}{2}}(\Bun_M)_{\on{co}} \otimes 
\on{I}(G,P^-)^{\on{loc}}_{\on{co},\CZ} \to \Dmod_{\frac{1}{2}}(\Bun_G\times \CZ)_{\on{co}} \simeq \\
\simeq \Dmod_{\frac{1}{2}}(\Bun_G)_{\on{co}}\otimes \Dmod(\CZ)
\end{multline*}
to be 
\begin{equation} \label{e:defn of Eis glob}
(\wt\sfp_\CZ^-)_*\left((\wt\sfq^-)^!(-)\sotimes (\pi_*\circ \wt\sfs^!)^{\on{Av}}_\CZ(-)\right)[\on{shift}],
\end{equation} 
where: 

\begin{itemize}

\item The functor $(\wt\sfq^-)^!(-)\sotimes (-)$ is understood in the sense of \eqref{e:sotimes rel};

\smallskip

\item The functor $(\wt\sfp_\CZ^-)_*$ is understood in the sense of \eqref{e:p_* co};

\smallskip

\item $[\on{shift}]$ is as in \eqref{e:global shift}.

\end{itemize}

\sssec{}

Unwinding the constructions, we obtain that the functor $\on{pre-Eis}^{-,\on{enh}_{\on{co}},\on{glob}}_{\on{co},*,\CZ}$ factors canonically as 
\begin{multline} \label{e:pre Eis enh glob}
\Dmod_{\frac{1}{2}}(\Bun_M)_{\on{co}} \otimes \on{I}(G,P^-)^{\on{loc}}_{\on{co},\CZ} \simeq \\
\simeq \left(\Dmod_{\frac{1}{2}}(\Bun_M)_{\on{co}} \otimes \Dmod(\CZ)\right)\underset{\Dmod(\CZ)}\otimes \on{I}(G,P^-)^{\on{loc}}_{\on{co},\CZ}\to \\
\to \left(\Dmod_{\frac{1}{2}}(\Bun_M)_{\on{co}} \otimes \Dmod(\CZ)\right)\underset{\Sph_{M,\CZ}}\otimes \on{I}(G,P^-)^{\on{loc}}_{\on{co},\CZ}= \\
=\Dmod_{\frac{1}{2}}(\Bun_M)_{\on{co},\CZ}^{-,\on{enh}_{\on{co}}} \to \Dmod_{\frac{1}{2}}(\Bun_G)_{\on{co}}\otimes \Dmod(\CZ).
\end{multline} 

Denote the last arrow in \eqref{e:pre CT enh glob} by 
$$\Eis^{-,\on{enh}_{\on{co}},\on{glob}}_{\on{co},*,\CZ}:\Dmod_{\frac{1}{2}}(\Bun_M)_{\on{co},\CZ}^{-,\on{enh}_{\on{co}}} \to 
\Dmod_{\frac{1}{2}}(\Bun_G)_{\on{co}}\otimes \Dmod(\CZ).$$

\medskip 

By construction, the functor $\Eis^{-,\on{enh}_{\on{co}}}_{\on{co},*,\CZ}$ intertwines the actions of $\Sph_{G,\CZ}$ on the two sides. 

\sssec{}

Finally, denote 
$$\Eis^{-,\on{enh}_{\on{co}},\on{glob}}_{\on{co},*}:=(\on{Id}\otimes \on{C}^\cdot(\Ran,-))\circ \Eis^{-,\on{enh}_{\on{co}},\on{glob}}_{\on{co},*,\Ran},\quad 
\Dmod_{\frac{1}{2}}(\Bun_M)_{\on{co},\Ran}^{-,\on{enh}_{\on{co}}} \to \Dmod_{\frac{1}{2}}(\Bun_G)_{\on{co}}.$$

\medskip

By duality, from \propref{p:loc vs global CT enh co} we obtain:

\begin{cor} \label{c:loc vs global Eis enh co} \hfill 

\smallskip

\noindent{\em(a)} 
There is a canonical isomorphism 
$$\Eis^{-,\on{enh}_{\on{co}},\on{glob}}_{\on{co},*,\CZ}\simeq \Eis^{-,\on{enh}_{\on{co}}}_{\on{co},*,\CZ}$$
as $\Sph_{G,\CZ}$-linear functors
$$\Dmod_{\frac{1}{2}}(\Bun_M)_{\on{co},\CZ}^{-,\on{enh}_{\on{co}}} \rightrightarrows \Dmod_{\frac{1}{2}}(\Bun_G)_{\on{co}}\otimes \Dmod(\CZ).$$

\noindent{\em(b)} 
There is a canonical isomorphism 
$$\Eis^{-,\on{enh}_{\on{co}},\on{glob}}_{\on{co},*}\simeq \Eis^{-,\on{enh}_{\on{co}}}_{\on{co},*}$$
as functors
$$\Dmod_{\frac{1}{2}}(\Bun_M)_{\on{co},\Ran}^{-,\on{enh}_{\on{co}}} \rightrightarrows \Dmod_{\frac{1}{2}}(\Bun_G)_{\on{co}}.$$

\end{cor}

\ssec{Proof of \lemref{l:approx}} \label{ss:proof of approx} 

\sssec{}

In order to simplify the exposition, we will prove the required isomorphism after restriction to the preimage of the trivial bundle
$\CP^0_M\in \Bun_M$. 

\medskip

Denote:
$$\BunNtpx:=\on{pt}\underset{\CP^0_M,\Bun_M}\times \BunPmtpx \text{ and }
\on{Hecke}^{\on{glob}}_{G,N_P^-,\ul{x}}:=\on{pt}\underset{\CP^0_M,\Bun_M}\times \on{Hecke}^{\on{glob}}_{G,P^-,\ul{x}}.$$

We keep the same symbols for the corresponding maps
$$\sff_{\ul{x}}:\on{Hecke}^{\on{glob}}_{G,N_P^-,\ul{x}}\to \BunNtpx,$$
$$\pi_{\ul{x}}:\Gr_{G,\ul{x}}\to \BunNtpx,$$
$$\sfs_{\ul{x}}:\on{Hecke}^{\on{glob}}_{G,N_P^-,\ul{x}}\to \fL^+(N_P^-)_{\ul{x}}\backslash \Gr_{G,\ul{x}},$$
and
$$\wt\sfs_{\ul{x}}: \Gr_{G,\ul{x}}\to \fL^+(N_P^-)_{\ul{x}}\backslash \Gr_{G,\ul{x}}.$$

We will prove an isomorphism 
$$((\sff_{\ul{x}})_*\circ \sfs_{\ul{x}}^!)^{\on{Av}}\simeq (\pi_*\circ \wt\sfs^!)^{\on{Av}}_{\ul{x}}[2\delta_{N^-_P}]$$
as functors
$$\Dmod_{\frac{1}{2}}(\Gr_{G,\ul{x}})^{\fL^+(N^-_P)_{\ul{x}}}\rightrightarrows \Dmod_{\frac{1}{2}}(\BunNtpx).$$

\sssec{}

Let $\Bun^{\on{level}_{\ul{x}}}_{N^-_P}$ denote the moduli space of $N^-_P$-bundles on $X$ with a full level structure 
at $\ul{x}$. We have
$$\on{Hecke}^{\on{glob}}_{G,N_P^-,\ul{x}}\simeq (\Bun^{\on{level}_{\ul{x}}}_{N^-_P}\times \Gr_{G,\ul{x}})/\fL^+(N^-_P)_{\ul{x}}.$$

For a group subscheme
$$\fL^+(N^-_P)_{\ul{x}}\subseteq N^\alpha\subset \fL(N^-_P)_{\ul{x}}$$ denote
$$\CY^\alpha:=(\Bun^{\on{level}_{\ul{x}}}_{N^-_P}\times \Gr_{G,\ul{x}})/N^\alpha.$$

We can regard $\CY^\alpha$ as a fibration with typical fiber $\Gr_{G,\ul{x}}$ over 
$$\Bun^\alpha_{N^-_P}:=\Bun^{\on{level}_{\ul{x}}}_{N^-_P}/N^\alpha,$$
i.e.,
$$\CY^\alpha\simeq \Gr_{G,\ul{x}}\wt\times\Bun^\alpha_{N^-_P}.$$

Denote by $\sfs^\alpha$ the resulting map
$$\CY^\alpha\to N^\alpha\backslash \Gr_{G,\ul{x}}.$$

Let $(\sfs^\alpha)^!$ denote the functor
$$\Dmod_{\frac{1}{2}}(\Gr_{G,\ul{x}})^{N_\alpha}\to \Dmod_{\frac{1}{2}}(\CY^\alpha), \quad 
\CF\mapsto \CF\wt\boxtimes \omega_{\Bun^\alpha_{N^-_P}}.$$

\sssec{}

Note that the map $\sff_{\ul{x}}$ factors as
$$\on{Hecke}^{\on{glob}}_{G,N_P^-,\ul{x}}\to \CY^\alpha \overset{\sff^\alpha}\to \BunNtpx,$$
where the first arrow is a fibration with contractible fibers of relative dimension 
$$d_\alpha:=\dim(N^\alpha/\fL^+(N^-_P)_{\ul{x}}).$$

\medskip

From here, we obtain that 
\begin{equation} \label{e:f alpha 1}
\sff^\alpha_*\circ (\sfs^\alpha)^! \simeq (\sff_{\ul{x}})_*\circ \sfs_{\ul{x}}^![-2d_\alpha],
\end{equation} 
as functors
$$\Dmod_{\frac{1}{2}}(\Gr_{G,\ul{x}})^{N_\alpha}\rightrightarrows \Dmod_{\frac{1}{2}}(\BunNtpx).$$

\sssec{}

Note that we have a canonically defined map
$$\sg_{\ul{x}}:\Gr_{G,\ul{x}}\to \on{Hecke}^{\on{glob}}_{G,N_P^-,\ul{x}}$$
so that
$$\sff_{\ul{x}}\circ \sg_{\ul{x}}=\pi_{\ul{x}}, \quad \Gr_{G,\ul{x}}\rightrightarrows \BunPmtpx$$
and 
$$\sfs_{\ul{x}}\circ \sg_{\ul{x}}=\wt\sfs_{\ul{x}}, \quad 
\Gr_{G,\ul{x}}\rightrightarrows \fL^+(P^-)_{\ul{x}}\backslash \Gr_{G,\ul{x}}.$$

\sssec{}

Let $\sg^\alpha$ denote the composition
$$\Gr_{G,\ul{x}} \overset{\sg_{\ul{x}}}\longrightarrow \on{Hecke}^{\on{glob}}_{G,N_P^-,\ul{x}}\to \CY_\alpha.$$

We have
\begin{equation} \label{e:f alpha 0}
\sff^\alpha\circ \sg^\alpha=\pi_{\ul{x}} \text{ and } \sfs^\alpha\circ \sg^\alpha=\wt\sfs_{\ul{x}}.
\end{equation} 

\sssec{}

The map $\sg^\alpha$ is the base change of the map
$$\on{pt}\to \Bun^\alpha_{N^-_P}.$$

When $\alpha$ is large enough, the latter map is a fibration with contractible fibers
of relative dimension $d_\alpha-\delta_{N^-_P}$. Hence, the same is true for $\sg^\alpha$. 

\medskip

Hence, from \eqref{e:f alpha 0} we obtain
\begin{equation} \label{e:f alpha 2}
(\pi_{\ul{x}})_*\circ \wt\sfs_{\ul{x}}^! \simeq (\sff^\alpha)_*\circ (\sfs^\alpha)^![2(d_\alpha-\delta_{N^-_P})]
\end{equation} 
as functors
$$\Dmod_{\frac{1}{2}}(\Gr_{G,\ul{x}})^{N_\alpha}\rightrightarrows \Dmod_{\frac{1}{2}}(\BunNtpx).$$

\sssec{}

Comparing \eqref{e:f alpha 1} and \eqref{e:f alpha 2}, we obtain 
$$(\pi_{\ul{x}})_*\circ \wt\sfs^! \simeq (\sff_{\ul{x}})_*\circ \sfs_{\ul{x}}^![-2\delta_{N^-_P}],$$
as desired.

\qed[\lemref{l:approx}]

\ssec{Local-to-global construction for Drinfeld's compactification, the !-version}

\sssec{}

Consider the (partially defined) left adjoint $(\pi_\CZ)_!$ to 
\begin{multline} \label{e:pi ! Z}
\pi_\CZ^!:\Dmod_{\det^{\frac{1}{2}}_{G,M}}(\BunPmtpZ)\to \\
\to \Dmod_{\pi^*(\det^{\frac{1}{2}}_{G,M})}(\Gr_{G,\Bun_M,\CZ})\simeq \Dmod_{\wt\sfs^*(\det^{\frac{1}{2}}_{\Gr_{G,\Ran}})}(\Gr_{G,\Bun_M,\CZ}).
\end{multline}

\sssec{}

We claim:

\begin{prop} \label{p:loc-to-glob Drinf defnd} \hfill

\smallskip

\noindent{\em(a)} 
The functor $(\pi_\CZ)_!$ is well-defined and is $(\Sph_G\otimes \Sph_M)_\CZ$-linear on the essential image of 
\begin{multline} \label{e:loc-to-glob Drinf defnd}
\on{I}(G,P^-)^{\on{loc}}_\CZ\to \Dmod_{\det^{\frac{1}{2}}_{\Gr_{G,\CZ}}}
(\Gr_{G,\CZ})^{\fL(N_P^-)_\CZ\cdot \fL^+(M)_\CZ}\hookrightarrow \\
\to \Dmod_{\det^{\frac{1}{2}}_{\Gr_{G,\CZ}}}(\Gr_{G,\ul{x}})^{\fL^+(P^-)_\CZ}\overset{\wt\sfs^!}\to \Dmod_{\wt\sfs^*(\det^{\frac{1}{2}}_{\Gr_{G,\CZ}})}(\Gr_{G,\Bun_M,\CZ}).
\end{multline}
Denote the resulting functor
$$\on{I}(G,P^-)^{\on{loc}}_\CZ\to \Dmod_{\det^{\frac{1}{2}}_{G,M}}(\BunPmtpZ)$$
by $(\pi_!\circ \wt\sfs^!)_\CZ$. 

\smallskip

\noindent{\em(b)} 
The assignment $\CZ\rightsquigarrow (\pi_!\circ \wt\sfs^!)_\CZ$ commutes with base change, i.e., for a map $f:\CZ'\to \CZ$, the
natural transformation 
$$(\pi_!\circ \wt\sfs^!)_{\CZ'}\circ f^!\to f^!\circ (\pi_!\circ \wt\sfs^!)_\CZ$$
is an isomorphism as functors
$$\on{I}(G,P^-)^{\on{loc}}_{\CZ'}\rightrightarrows \Dmod_{\det^{\frac{1}{2}}_{G,M}}(\BunPmtpZp).$$

\smallskip

\noindent{\em(c)} 
Let $\CZ=S$ be a scheme. Then the functor $(\pi_!\circ \wt\sfs^!)_S$ sends compact objects in $\on{I}(G,P^-)^{\on{loc}}_S$
to objects in $ \Dmod_{\det^{\frac{1}{2}}_{G,M}}(\BunPmtpS)$ that are ULA with respect to
$\wt\sfq^-:\BunPmtpS\to \Bun_M$. 

\end{prop}

The proof \propref{p:loc-to-glob Drinf defnd} will occupy the majority of this subsection.

\medskip

Note that to prove points (a) and (b) we can assume that $\CZ$ (resp., $\CZ$ and $\CZ'$) are schemes, to be denoted $S$ and $S'$,
respectively. 

\sssec{}

We will first show that $(\pi_!\circ \wt\sfs)_S^!$ is defined on 
$\Delta^{-,\semiinf}_S\in \on{I}(G,P^-)^{\on{loc}}_S$ and 
$$(\pi_!\circ \wt\sfs)_S^!(\Delta^{-,\semiinf}_S)\simeq (j_S)_!(\omega_{\Bun_{P^-}\times S}),$$
where:

\begin{itemize}

\item 
$j_S$ denotes the 
corresponding locally closed embedding $\Bun_{P^-}\times S\to \BunPmtpS$;

\medskip

\item The object $\omega_{\Bun_{P^-}\times S}\in \Dmod_{\det^{\frac{1}{2}}_{G,M}}(\Bun_{P^-}\times S)$
makes sense due to the trivialization of the gerbe $\det^{\frac{1}{2}}_{G,M}|_{\Bun_{P^-}}$. 

\end{itemize}

\sssec{}

Consider the unit orbit of $\fL(N^-_P)$ acting on $\Gr_G$:
$$\bi_S:(\fL(N^-_P)/\fL^+(N^-_P))_S\to \Gr_{G,S}.$$

Recall that
\begin{equation} \label{e:Delta again}
\Delta^{-,\semiinf}_S=(\bi_S)_!(\omega_{(\fL(N^-_P)/\fL^+(N^-_P))_S}).
\end{equation} 

\medskip

Consider the corresponding locally closed embedding
$$(\fL(N^-_P)/\fL^+(N^-_P))_{\Bun_M,S}\overset{\bi_{\Bun_M,S}}\hookrightarrow \Gr_{G,\Bun_M,S}.$$

From \eqref{e:Delta again}, we obtain 
$$\wt\sfs^!_S(\Delta^{-,\semiinf}_S)\simeq (\bi_{\Bun_M,S})_!(\omega_{(\fL(N^-_P)/\fL^+(N^-_P))_{\Bun_M,S}}).$$

In particular, the object $\wt\sfs^!_S(\Delta^{-,\semiinf}_S)$ is \emph{ind-holonomic}, so the partially defined functor $(\pi_S)_!$ 
is automatically defined on it. 

\sssec{}

Note that the map
$$(\fL(N^-_P)/\fL^+(N^-_P))_{\Bun_M,S}\overset{\bi_{\Bun_M,S}}\hookrightarrow \Gr_{G,\Bun_M,S} \overset{\pi_S}\to \BunPmtpS$$
factors as
$$(\fL(N^-_P)/\fL^+(N^-_P))_{\Bun_M,S}\overset{\pi^0_S}\to \Bun_{P^-}\times S \overset{j_S}\to \BunPmtpS.$$

Moreover, the above map 
$$(\fL(N^-_P)/\fL^+(N^-_P))_{\Bun_M,S}\overset{\pi^0_S}\to \Bun_{P^-}\times S$$
is a locally trivial fibration, whose fibers are ind-contractible.\footnote{It is for this assertion to hold that we need $S$ to be non-empty.}

\medskip

Hence, the naturally defined map
$$(\pi^0_S)_!(\omega_{(\fL(N^-_P)/\fL^+(N^-_P))_{\Bun_M,S}})\to \omega_{\Bun_{P^-}\times S}$$
is an isomorphism.

\medskip

Therefore.
\begin{multline*} 
(\pi_!\circ \wt\sfs)_S^!(\Delta^{-,\semiinf}_S)\simeq
(\pi_S)_!\circ (\bi_{\Bun_M,S})_!(\omega_{(\fL(N^-_P)/\fL^+(N^-_P))_{\Bun_M,S}})\simeq \\
\simeq (j_S)_!\circ (\pi^0_S)_!(\omega_{(\fL(N^-_P)/\fL^+(N^-_P))_{\Bun_M,S}})\simeq 
(j_S)_!(\omega_{\Bun_{P^-}\times S}),
\end{multline*}
as desired. 

\sssec{}

Let 
$$p_{\Gr,S}:\Gr_{G,S}\to S, \,\,  p_{\Gr,\Bun_M,S}:\Gr_{G,\Bun_M,S}\to S \text{ and } p_{\wt\Bun_{P^-},S}:\BunPmtpS\to S$$
denote the tautological projections. 

\medskip

Next we will show that the functor $(\pi_!\circ \wt\sfs^!)_S$
is defined on objects of the form 
$$\Delta^{-,\semiinf}_S\sotimes p_S^!(\CF_S), \quad \CF_S\in \Dmod(S).$$

Moreover, we will show that the naturally defined map
\begin{equation} \label{e:loc-to-glob Drinf defnd 0}
(\pi_!\circ \wt\sfs^!)_S(\Delta^{-,\semiinf}_S\sotimes p_{\Gr,S}^!(\CF_S))\to 
(\pi_!\circ \wt\sfs^!)_S(\Delta^{-,\semiinf}_S)\sotimes p_{\wt\Bun_{P^-},S}^!(\CF_S)
\end{equation}
is an isomorphism.

\sssec{}

We have
\begin{multline*} 
\wt\sfs^!(\Delta^{-,\semiinf}_S\sotimes p_{\Gr,S}^!(\CF_S))\simeq
\wt\sfs^!(\Delta^{-,\semiinf}_S)\sotimes p_{\Gr,\Bun_M,S}^!(\CF_S)\simeq \\
\simeq (\bi_{\Bun_M,S})_!(\omega_{(\fL(N^-_P)/\fL^+(N^-_P))_{\Bun_M,S}}) \sotimes p_{\Gr,\Bun_M,S}^!(\CF_S).
\end{multline*} 

Now, it follows from \cite[Proposition 1.6.3]{Ga5} that the naturally defined map
\begin{multline}  \label{e:loc-to-glob Drinf defnd 1}
(\bi_{\Bun_M,S})_!\circ (\bi_{\Bun_M,S})^!\circ p_{\Gr,\Bun_M,S}^!(\CF_S)\to \\
\to (\bi_{\Bun_M,S})_!(\omega_{(\fL(N^-_P)/\fL^+(N^-_P))_{\Bun_M,S}}) \sotimes p_{\Gr,\Bun_M,S}^!(\CF_S)
\end{multline}
is an isomorphism; in particular, the left-hand side in \eqref{e:loc-to-glob Drinf defnd 1} is defined.

\medskip

Hence, in order to prove \eqref{e:loc-to-glob Drinf defnd 0}, it suffices to show that the naturally defined map
\begin{multline}  \label{e:loc-to-glob Drinf defnd 2}
(\pi_S)_!\circ (\bi_{\Bun_M,S})_!\circ (\bi_{\Bun_M,S})^!\circ p_{\Gr,\Bun_M,S}^!(\CF_S)\to \\
\to (\pi_S)_!\circ (\bi_{\Bun_M,S})_!(\omega_{(\fL(N^-_P)/\fL^+(N^-_P))_{\Bun_M,S}})\sotimes p_{\wt\Bun_{P^-},S}^!(\CF_S)\simeq 
(j_S)_!(\omega_{\Bun_{P^-}\times S}) \sotimes p_{\wt\Bun_{P^-},S}^!(\CF_S)
\end{multline}
is an isomorphism. 

\sssec{}

We rewrite the left-hand side in \eqref{e:loc-to-glob Drinf defnd 2} as
$$(j_S)_!\circ (\pi^0_S)_!\circ (p_{\Gr,\Bun_M,S}\circ \bi_{\Bun_M,S})^!(\CF_S)\simeq
(j_S)_!\circ (\pi^0_S)_!\circ (\pi^0_S)^!\circ j_S^!\circ p_{\wt\Bun_{P^-},S}^!(\CF_S).$$

By the contractibility result mentioned above, the above expression maps isomorphically to
$$(j_S)_!\left(\omega_{\Bun_{P^-}\times S}\sotimes \left(j_S^!\circ p_{\wt\Bun_{P^-},S}^!(\CF_S)\right)\right).$$

\medskip

Thus, in order to prove \eqref{e:loc-to-glob Drinf defnd 0}, we have to show that the naturally defined map
$$(j_S)_!\left(\omega_{\Bun_{P^-}\times S}\sotimes \left(j_S^!\circ p_{\wt\Bun_{P^-},S}^!(\CF_S)\right)\right)\to
(j_S)_!(\omega_{\Bun_{P^-}\times S})\sotimes p_{\wt\Bun_{P^-},S}^!(\CF_S)$$
is an isomorphism. However, the latter is straightforward: 

\medskip

The map $j_S$ factors as a composition 
$$\Bun_{P^-}\times S\overset{j\times \on{id}}\longrightarrow  \BunPmt \times S\overset{i_S}\to \BunPmtpS,$$
where $i_S$ is a closed embedding. Hence, it is enough to show that
$$(j\times \on{id})_!(\omega_{\Bun_{P^-}}\boxtimes \CF_S)\to j_!(\omega_{\Bun_{P^-}})\boxtimes \CF_S$$
is an isomorphism, which is obvious. 

\sssec{}

Thus, we have shown that the functor $(\pi_!\circ \wt\sfs^!)_S$ is defined and $\Dmod(S)$-linear on the full subcategory of 
$$\on{I}(G,P^-)^{0,\on{loc}}_S\subset \on{I}(G,P^-)^{\on{loc}}_S$$
spanned by objects of the form
$$\Delta^{-,\semiinf}_S\sotimes p_S^!(\CF_S).$$

\sssec{}

We now prove point (a) of \propref{p:loc-to-glob Drinf defnd} on all of $\on{I}(G,P^-)^{\on{loc}}_S$. Note that the functor \eqref{e:pi ! Z}
(for $\CZ=S$) is $\Sph_{M,S}$-linear with respect to the natural actions of $\Sph_{M,S}$ on the two sides.

\medskip

Recall that $\on{I}(G,P^-)^{\on{loc}}_S$ is generated by $\on{I}(G,P^-)^{0,\on{loc}}_S$ under the action of
$\Sph_{M,S}$.

\medskip

Finally, recall that $\Sph_{M,S}$ is rigid relative to $\Dmod(S)$. Now, the assertion of \propref{p:loc-to-glob Drinf defnd}(a)
follows from the following general paradigm:

\medskip

Let $\bA^0$ and $\bA$ be a pair of monoidal categories, with $\bA^0$ symmetric monoidal and $\bA$ rigid
over $\bA^0$. Let $F:\bC\to \bD$ be a functor between $\bA$-module categories. Suppose that the (partially defined) 
left adjoint $F^L$ to $F$ is defined and is $\bA^0$-linear on a full subcategory $\bD^0\subset \bD$. 

\medskip

We claim:

\begin{lem}
Under the above circumstances, $F^L$
is defined and is $\bA$-linear on the full subcategory of $\bD$ generated by $\bD^0$ under the action of $\bA$.
\end{lem} 

\sssec{}

Point (b) of \propref{p:loc-to-glob Drinf defnd} on all of $\on{I}(G,P^-)^{\on{loc}}_S$ follows by a similar rigidity argument. 

\sssec{}

Let us now prove point (c) of \propref{p:loc-to-glob Drinf defnd}. The starting point is that the object
$$j_!(\omega_{\Bun_{P^-}}) \in \Dmod_{\det^{\frac{1}{2}}_{G,M}}(\BunPmt)$$
is ULA over $\Bun_M$. This is the assertion of \cite[Theorem 5.1.5]{BG1}.

\medskip

Now, it is easy to see that whenever $\CF\in \Dmod_{\det^{\frac{1}{2}}_{G,M}}(\BunPmt)$ is ULA over $\Bun_M$,
then so is any object of the form 
$$\CF_M\underset{M}\star (i_S)_*(\CF\boxtimes \omega_S), \quad \CF_M\in (\Sph_{M,S})^c.$$

\medskip

Now, the assertion of \propref{p:loc-to-glob Drinf defnd}(c) follows from the fact that every compact object in $\on{I}(G,P^-)^{\on{loc}}_S$
is a retract of an object object of the form
$$\CF_M\underset{M}\star \Delta^{-,\semiinf}_S, \quad \CF_M\in (\Sph_{M,S})^c$$
and that the functor $(\pi_!\circ \wt\sfs^!)_S$ is $\Sph_{M,S}$-linear. 

\qed[\propref{p:loc-to-glob Drinf defnd}]

\sssec{} \label{sss:Verdier dual on BunPt} 

Let $\CZ=S$ be a scheme, and let $U\subset \Bun_M$ be a quasi-compact open substack. 

\medskip

Note that the functor $(\pi_!\circ \wt\sfs^!)_S$ sends compact objects to objects that are \emph{locally compact} relative to $\Bun_M$,
i.e., they become compact after a pullback $U\to \Bun_M$, whenever $U$ is a scheme. 
Note that such objects are !-extended from a quasi-compact substack of $\BunPmtpS\underset{\Bun_M}\times U$, see
\cite[Proposition 2.3.7]{DG2}. 

\sssec{}

We claim: 

\begin{lem} \label{l:Verdier duality loc vs glob}
For $\CF\in (\on{I}(G,P^-)^{\on{loc}}_S)^c$ we have
\begin{equation} \label{e:Verdier duality loc vs glob}
\BD^{\on{glob}}((\pi_!\circ \wt\sfs^!)_S(\CF)) \simeq 
(\pi_*\circ \wt\sfs^!)^{\on{Av}}_\CZ(\BD^{\on{loc}}(\CF))[-2 \dim(\Bun_M)],
\end{equation}
where:

\begin{itemize}

\item $\BD^{\on{loc}}$ denote the canonical equivalence $((\on{I}(G,P^-)^{\on{loc}}_S)^c)^{\on{op}}\to (\on{I}(G,P^-)^{\on{loc}}_{\on{co},S})^c$;

\medskip

\item $\BD^{\on{glob}}$ denotes the Verdier duality functor 
$$(\Dmod_{\det^{\frac{1}{2}}_{G,M}}(\BunPmtpZ)^{\on{loc.c}/\Bun_M})^{\on{op}} \to 
\Dmod_{\det^{\frac{1}{2}}_{G,M}}(\BunPmtpZ)_{\on{co}/\Bun_M},$$
which make sense due to the !-extension property mentioned above. 

\end{itemize} 

\end{lem} 

\begin{proof}

To simplify the exposition, we will take $S=\on{pt}$, so that $S\to \Ran$ corresponds to $\ul{x}\in \Ran$.

\medskip

Any compact object $\CF\in \on{I}(G,P^-)^{\on{loc}}_{\ul{x}}$ is of the form 
$$\on{Av}_!^{\fL(N^-_P)_S}(\CF'),\quad \CF'\in \Dmod_{\frac{1}{2}}(\Gr_{G,\ul{x}})^{N^\alpha\cdot \fL^+(M)_{\ul{x}},\on{ren}},$$
where $N^\alpha\subset \fL(N^-_P)_{\ul{x}}$ is a large enough subgroup, normalized by $\fL^+(M)_{\ul{x}}$, and $\CF'$ is
compact as an object of $\Dmod_{\frac{1}{2}}(\Gr_{G,\ul{x}})^{N^\alpha}$. 

\medskip

Taking $N^\alpha$ to be as \secref{sss:push Av 2}, we have
$$(\pi_!\circ \wt\sfs^!)_{\ul{x}}(\CF)\simeq (\pi_{\ul{x}})_!(\CF'\wt\boxtimes \omega_{\Bun_M}).$$

\medskip

Note that $\BD^{\on{loc}}(\CF)$ is the image of 
$$\BD(\CF')\in \Dmod_{\frac{1}{2}}(\Gr_{G,\ul{x}})^{N^\alpha\cdot \fL^+(M)_{\ul{x}},\on{ren}}$$
under 
$$\Dmod_{\frac{1}{2}}(\Gr_{G,\ul{x}})^{N^\alpha\cdot \fL^+(M)_{\ul{x}},\on{ren}}\simeq
\Dmod_{\frac{1}{2}}(\Gr_{G,\ul{x}})_{N^\alpha}^{\cdot \fL^+(M)_{\ul{x}},\on{ren}}\twoheadrightarrow
\Dmod_{\frac{1}{2}}(\Gr_{G,\ul{x}})_{\fL(N^-_P)_{\ul{x}}}^{\fL^+(M)_{\ul{x}},\on{ren}},$$
where $\BD$ is Verdier duality on $\Gr_{G,\ul{x}}$. 

\medskip

Now, the assertion of the lemma follows from the fact that Verdier duality swaps *- and !- direct images, whenever the latter are defined,
and taking into account that
$$\BD(\CF'\wt\boxtimes \omega_{\Bun_M}) \simeq \BD(\CF')\wt\boxtimes \omega_{\Bun_M}[-2\dim(\Bun_M)].$$

\end{proof} 

\sssec{}

In particular, for $\CF$ as above and 
$$\CF_G\in  \Dmod_{\det^{\frac{1}{2}}_G}(\BunPmtpS),\,\, \CF_M\in  \Dmod_{\frac{1}{2}}(\Bun_M\times S)$$ we 
have a canonical isomorphism  
\begin{multline} \label{e:Verdier duality BunPt}
\CHom_{\Dmod_{\det^{\frac{1}{2}}_G}(\BunPmtpS)}\left((\wt\sfq^-_S)^*(\CF_M) \overset{*}\otimes (\pi_!\circ \wt\sfs^!)_S(\CF),\CF_G\right) \simeq \\
\simeq \CHom_{\Dmod_{\frac{1}{2}}(\Bun_M\times S)}\left(\CF_M,(\wt\sfq^-_S)_*\left(\CF_G\sotimes (\pi_*\circ \wt\sfs^!)^{\on{Av}}_\CZ(\BD^{\on{loc}}(\CF))\right)\right)[-2 \dim(\Bun_M)],
\end{multline}
where the operation $\overset*{\otimes}$ is well-defined thanks to \lemref{p:loc-to-glob Drinf defnd}(c). 

\sssec{} \label{sss:! BunPt unital}

As in Sects. \ref{sss:* BunPt unital 1}-\ref{sss:* BunPt unital 2}, we can view the assignment
$$\CZ\rightsquigarrow  \Dmod_{\det^{\frac{1}{2}}_{G,M}}(\BunPmtpZ)$$
as a unital crystal of categories over $\Ran$
and the assigment
$$\CZ\rightsquigarrow (\pi_!\circ \wt\sfs)_\CZ$$
as a \emph{right}-lax unital functor, to be denoted $\ul{(\pi_!\circ \wt\sfs)}$.  

\sssec{}

We claim: 

\begin{cor} \label{c:! BunPt unital} 
The right-lax unital structure on $\ul{(\pi_!\circ \wt\sfs)}$ is strict.
\end{cor} 

\begin{proof}

Follows from \propref{p:* BunPt unital} using \lemref{l:Verdier duality loc vs glob}. 

\end{proof}

\ssec{The global enhanced functors, the !-version}

\sssec{}

For $\CZ\to \Ran$, define the functor 
\begin{multline*} 
\on{pre-Eis}^{-,\on{enh},\on{glob}}_{!,\CZ}:\Dmod_{\frac{1}{2}}(\Bun_M) \otimes 
\on{I}(G,P^-)^{\on{loc}}_{\CZ} \to \Dmod_{\frac{1}{2}}(\Bun_G\times \CZ)\simeq \\
\simeq \Dmod_{\frac{1}{2}}(\Bun_G)\otimes \Dmod(\Ran)
\end{multline*} 
to be 
\begin{equation} \label{e:defn of Eis ! glob}
(\wt\sfp^-_\CZ)_*\left((\wt\sfq^-)^!(-)\sotimes (\pi_!\circ \wt\sfs^!)_\CZ(-)\right)[-\on{shift}],
\end{equation} 
where the functors  $(\wt\sfq^-)^!$,  $(\wt\sfp_\CZ^-)_*$ and $\sotimes$ are the corresponding functors
on the usual categories of (twisted) D-modules. 

\sssec{}

Recall the functor $\on{pre-CT}^{-,\on{enh}_{\on{co}},\on{glob}}_{*,\CZ}$ from \secref{sss:pre CT enh}.
By duality, it gives rise to a functor
\begin{equation} \label{e:pre CT enh dual}
\Dmod_{\frac{1}{2}}(\Bun_G\times \CZ)\to \Dmod_{\frac{1}{2}}(\Bun_M) \otimes \on{I}(G,P^-)^{\on{loc}}_{\CZ}.
\end{equation}

We claim:

\begin{lem} \label{l:Eis enh glob}
The functor $\on{pre-Eis}^{-,\on{enh},\on{glob}}_{!,\CZ}$ is the left adjoint of the functor \eqref{e:pre CT enh dual}.
\end{lem}

\begin{proof}

With no restriction of generality, we can assume that $\CZ=S$ is a scheme. In order to simplify the notation, 
we will assume that $S=\on{pt}$, corresponding to $\ul{x}\in \Ran$.

\medskip

Fix objects
$$\CF_G\in \Dmod_{\frac{1}{2}}(\Bun_G),\quad \CF_M\in \Dmod_{\frac{1}{2}}(\Bun_M),\quad 
\CF_{G,M}\in (\on{I}(G,P^-)^{\on{loc}}_{\ul{x}})^c.$$

Let $\BD^{\on{loc}}(\CF_{G,M})\in \on{I}(G,P^-)^{\on{loc}}_{\on{co},\ul{x}}$ be the ``formal dual" of $\CF_{G,M}$, see \secref{sss:Verdier dual on BunPt}. 

\medskip

We need to establish a canonical isomorphism 
\begin{multline} \label{e:Eis enh glob 1}
\CHom_{\Dmod_{\frac{1}{2}}(\Bun_G)}\left((\wt\sfp^-)_*\left((\wt\sfq^-)^!(\CF_M)\sotimes (\pi_!\circ \wt\sfs^!)_\CZ(\CF_{G,M})\right)[-\on{shift}],\CF_G\right) 
\simeq \\
\simeq \CHom_{\Dmod_{\frac{1}{2}}(\Bun_M)}\left(\CF_M,
(\wt\sfq^-_S)_*\left((\wt\sfp^-)^!(\CF_G)\sotimes (\pi_*\circ \wt\sfs^!)^{\on{Av}}_\CZ(\BD^{\on{loc}}(\CF_{G,M}))\right)[\on{shift}]\right).
\end{multline}

Since the morphism $(\wt\sfp^-_S)_*$ is ind-proper, we can rewrite the left-hand side of \eqref{e:Eis enh glob 1} as
\begin{equation} \label{e:Eis enh glob 2 LHS}
\CHom_{\Dmod_{\det_G^{\frac{1}{2}}}(\BunPmtpx)}\left((\wt\sfq^-)^!(\CF_M)\sotimes (\pi_!\circ \wt\sfs^!)_\CZ(\CF_{G,M})[-\on{shift}],
(\wt\sfp^-_S)^!(\CF_G)\right).
\end{equation} 

By \propref{p:loc-to-glob Drinf defnd}(c), we can rewrite \eqref{e:Eis enh glob 2 LHS} further as
\begin{equation} \label{e:Eis enh glob 3 LHS}
\CHom_{\Dmod_{\det_G^{\frac{1}{2}}}(\BunPmtpx)}\left((\wt\sfq^-)^*(\CF_M)\overset{*}\otimes 
(\pi_!\circ \wt\sfs^!)_\CZ(\CF_{G,M})[-\on{shift}-2 \dim(\Bun_M)],(\wt\sfp^-_S)^!(\CF_G)\right).
\end{equation}

Applying \eqref{e:Verdier duality BunPt}, we obtain that \eqref{e:Eis enh glob 3 LHS} identifies with the
right-hand side in \eqref{e:Eis enh glob 1}, as required. 

\end{proof}

\sssec{}

Unwinding the constructions, we obtain that the functor $\on{pre-Eis}^{-,\on{enh},\on{glob}}_{!,\CZ}$ factors canonically as 
\begin{multline} \label{e:pre Eis ! enh glob}
\Dmod_{\frac{1}{2}}(\Bun_M) \otimes \on{I}(G,P^-)^{\on{loc}}_{\CZ} \simeq \\
\simeq \left(\Dmod_{\frac{1}{2}}(\Bun_M) \otimes \Dmod(\CZ)\right)\underset{\Dmod(\CZ)}\otimes \on{I}(G,P^-)^{\on{loc}}_{\CZ}\to \\
\to \left(\Dmod_{\frac{1}{2}}(\Bun_M) \otimes \Dmod(\CZ)\right)\underset{\Sph_{M,\CZ}}\otimes \on{I}(G,P^-)^{\on{loc}}_{\CZ}= \\
=\Dmod_{\frac{1}{2}}(\Bun_M)_{\CZ}^{-,\on{enh}} \to \Dmod_{\frac{1}{2}}(\Bun_G)\otimes \Dmod(\CZ).
\end{multline} 

Denote the last arrow in \eqref{e:pre Eis ! enh glob} by 
$$\on{Eis}^{-,\on{enh},\on{glob}}_{!,\CZ}: \Dmod_{\frac{1}{2}}(\Bun_M)_{\CZ}^{-,\on{enh}} \to \Dmod_{\frac{1}{2}}(\Bun_G)\otimes \Dmod(\CZ).$$

\medskip 

By construction, the functor $\on{Eis}^{-,\on{enh},\on{glob}}_{!,\CZ}$ intertwines the actions of $\Sph_{G,\CZ}$ on the two sides. 

\medskip

Denote also
$$\on{Eis}^{-,\on{enh},\on{glob}}_!:=(\on{Id}\otimes \on{C}^\cdot(\Ran,-))\circ \on{Eis}^{-,\on{enh},\on{glob}}_{!,\Ran},\quad 
\Dmod_{\frac{1}{2}}(\Bun_M)_{\Ran}^{-,\on{enh}} \to \Dmod_{\frac{1}{2}}(\Bun_G).$$

\sssec{}

Note that the identification of \lemref{l:Eis enh glob} is compatible with the actions of $\Sph_{M,\CZ}$ and $\Sph_{G,\CZ}$.
Hence, combining with \propref{p:loc vs global CT enh co}, we obtain:

\begin{cor} \label{c:loc vs glob Eis ! enh} \hfill

\smallskip

\noindent{\em(a)}
There is a canonical isomorphism 
$$\on{Eis}^{-,\on{enh},\on{glob}}_{!,\CZ}\simeq \on{Eis}^{-,\on{enh}}_{!,\CZ}$$
as $\Sph_{G,\CZ}$-linear functors
$$\Dmod_{\frac{1}{2}}(\Bun_M)_{\CZ}^{-,\on{enh}} \rightrightarrows \Dmod_{\frac{1}{2}}(\Bun_G)\otimes \Dmod(\CZ).$$

\smallskip

\noindent{\em(b)}
There is a canonical isomorphism 
$$\on{Eis}^{-,\on{enh},\on{glob}}_!\simeq \on{Eis}^{-,\on{enh}}_!$$
as  functors
$$\Dmod_{\frac{1}{2}}(\Bun_M)_{\Ran}^{-,\on{enh}} \rightrightarrows \Dmod_{\frac{1}{2}}(\Bun_G).$$
\end{cor}

\sssec{}

Define the functor 
\begin{multline*} 
\on{pre-CT}^{-,\on{enh},\on{glob}}_{\on{co},?,\CZ}:\Dmod_{\frac{1}{2}}(\Bun_G)_{\on{co}} \otimes 
\on{I}(G,P^-)^{\on{loc}}_{\CZ} \to \Dmod_{\frac{1}{2}}(\Bun_M\times \CZ)_{\on{co}}\simeq \\
\simeq \Dmod_{\frac{1}{2}}(\Bun_M)_{\on{co}}\otimes \Dmod(\Ran)
\end{multline*} 
to be 
\begin{equation} \label{e:defn of CT ? glob}
(\wt\sfq^-_\CZ)_*\left((\wt\sfp^-)^!(-)\sotimes (\pi_!\circ \wt\sfs^!)_\CZ(-)\right)[-\on{shift}],
\end{equation} 
where: 

\begin{itemize}

\item The functor $(\wt\sfp^-)^!$ is understood in the sense \eqref{e:p^! co}; 

\smallskip

\item The functor $\sotimes$ is understood in the sense \eqref{e:sotimes co 1}; 

\smallskip

\item The functor $(\wt\sfq^-_\CZ)_*$ is understood in the sense \eqref{e:q_* co}; 

\smallskip
 
\item $[\on{shift}]$ is as in \eqref{e:global shift}.

\end{itemize}

\sssec{}

Unwinding the constructions, we obtain that the functor $\on{pre-CT}^{-,\on{enh},\on{glob}}_{\on{co},?,\CZ}$ factors canonically as 
\begin{multline} \label{e:pre CT enh ?  glob}
\Dmod_{\frac{1}{2}}(\Bun_G)_{\on{co}} \otimes \on{I}(G,P^-)^{\on{loc}}_{\CZ} \simeq \\
\simeq \left(\Dmod_{\frac{1}{2}}(\Bun_G)_{\on{co}} \otimes \Dmod(\CZ)\right)\underset{\Dmod(\CZ)}\otimes \on{I}(G,P^-)^{\on{loc}}_{\CZ}\to \\
\to \left(\Dmod_{\frac{1}{2}}(\Bun_G)_{\on{co}} \otimes \Dmod(\CZ)\right)\underset{\Sph_{G,\CZ}}\otimes \on{I}(G,P^-)^{\on{loc}}_{\CZ}= \\
=\Dmod_{\frac{1}{2}}(\Bun_G)_{\on{co},\CZ}^{-,\on{enh}} \to \Dmod_{\frac{1}{2}}(\Bun_M)_{\on{co}}\otimes \Dmod(\CZ).
\end{multline} 

Denote the last arrow in \eqref{e:pre CT enh ? glob} by 
$$\on{CT}^{-,\on{enh},\on{glob}}_{\on{co},?,\CZ}:\Dmod_{\frac{1}{2}}(\Bun_G)_{\on{co},\CZ}^{-,\on{enh}} \to 
\Dmod_{\frac{1}{2}}(\Bun_M)_{\on{co}}\otimes \Dmod(\CZ).$$

\medskip 

By construction, the functor $\on{CT}^{-,\on{enh},\on{glob}}_{\on{co},?,\CZ}$ intertwines the actions of $\Sph_{M,\CZ}$ on the two sides. 

\sssec{}

Denote 
$$\on{CT}^{-,\on{enh},\on{glob}}_{\on{co},?}:=(\on{Id}\otimes \on{C}^\cdot(\Ran,-))\circ \on{CT}^{-,\on{enh},\on{glob}}_{\on{co},?,\Ran},\quad 
\Dmod_{\frac{1}{2}}(\Bun_M)_{\on{co},\Ran}^{-,\on{enh}_{\on{co}}} \to \Dmod_{\frac{1}{2}}(\Bun_G)_{\on{co}}.$$

\medskip

By duality, from \propref{p:loc vs global CT enh co} we obtain:

\begin{cor} \label{c:loc vs global CT enh ? co} \hfill 

\smallskip

\noindent{\em(a)} 
There is a canonical isomorphism 
$$\on{CT}^{-,\on{enh},\on{glob}}_{\on{co},?,\CZ}\simeq \on{CT}^{-,\on{enh}}_{\on{co},?,\CZ}$$
as $\Sph_{M,\CZ}$-linear functors
$$\Dmod_{\frac{1}{2}}(\Bun_G)_{\on{co},\CZ}^{-,\on{enh}} \rightrightarrows \Dmod_{\frac{1}{2}}(\Bun_M)_{\on{co}}\otimes \Dmod(\CZ).$$

\noindent{\em(b)} 
There is a canonical isomorphism 
$$\on{CT}^{-,\on{enh},\on{glob}}_{\on{co},?}\simeq \on{CT}^{-,\on{enh}}_{\on{co},?}$$
as functors
$$\Dmod_{\frac{1}{2}}(\Bun_G)_{\on{co},\Ran}^{-,\on{enh}} \rightrightarrows \Dmod_{\frac{1}{2}}(\Bun_M)_{\on{co}}.$$

\end{cor}

\ssec{Relation to compactified Eisenstein series}

\sssec{} \label{sss:IC as twisted}

Let $\IC^{-,\semiinf,\on{glob}}$ be the (shifted) half-twisted intersection cohomology sheaf on $\BunPmt$, i.e.,
$$\IC^{-,\semiinf,\on{glob}}:=j_{!*}(\IC_{\Bun_{P^-}})[\dim(\Bun_M)]\in \Dmod_{\det^{\frac{1}{2}}_{G,M}}(\BunPmt),$$
where:

\begin{itemize}

\item $j$ denotes the embedding $\Bun_{P^-}\hookrightarrow \BunPmt$; 

\medskip

\item $\IC_{\Bun_{P^-}}[\dim(\Bun_M)]\simeq \omega_{\Bun_{P^-}}[-\on{shift}]$, where  $[\on{shift}]$ is as in \eqref{e:global shift};

\medskip

\item We regard $\omega_{\Bun_{P^-}}$ as an object of $\Dmod_{\det^{\frac{1}{2}}_{G,M}}(\Bun_{P^-})$
thanks to the trivialization of $\det^{\frac{1}{2}}_{G,M}|_{\Bun_{P^-}}.$

\end{itemize}

\sssec{}

Following \cite{BG1}, define the functor
$$\Eis^-_{!*}:\Dmod_{\frac{1}{2}}(\Bun_M)\to \Dmod_{\frac{1}{2}}(\Bun_G)$$
to be
$$(\wt\sfp^-)_*((\wt\sfq^-)^!(-)\sotimes \IC^{-,\semiinf,\on{glob}}).$$

\sssec{}

Define the functor
$$\on{CT}^-_{\on{co},!*}:=(\Eis^-_{!*})^\vee, \quad \Dmod_{\frac{1}{2}}(\Bun_G)_{\on{co}}\to \Dmod_{\frac{1}{2}}(\Bun_M)_{\on{co}}.$$

I.e.,
$$\on{CT}^-_{\on{co},!*}\simeq (\wt\sfq^-)_*((\wt\sfp^-)^!(-)\sotimes \IC^{-,\semiinf,\on{glob}}).$$

\sssec{}

Recall the functor $\Eis^-_{!,\IC}$, see \secref{sss:fact alg Eis}. We are going to prove:

\begin{prop} \label{p:compact Eis}
The functor $\Eis^-_{!*}$ identifies canonically with $\Eis^-_{!,\IC}$.
\end{prop} 

Note that as an immediate corollary we obtain:

\begin{cor} \label{c:compact CT}
The functor $\on{CT}^-_{\on{co},!*}$ identifies canonically with $\on{CT}^-_{\on{co},?,\IC}$.
\end{cor} 

\sssec{Proof of \propref{p:compact Eis}}

First, by the definition of $\IC^{-,\semiinf}$,
its support on $\Gr_G$ is contained in the closure of the unit orbit of $\fL(N^-_P)$.  This implies that the support of
$\wt\sfs^!(\IC^{-,\semiinf}_\Ran)$ on $\Gr_{G,\Bun_M,\Ran}$ is contained in $\pi^{-1}(\BunPmt\times \Ran)$. 

\medskip

In particular, the object 
\begin{equation} \label{e:loc IC semiinf}
(\pi_!\circ \wt\sfs^!)_\Ran(\IC^{-,\semiinf}_\Ran)\in  \Dmod_{\det^{\frac{1}{2}}_{G,M}}(\BunPmtpR)
\end{equation}
is supported on $\BunPmt\times \Ran\subset \BunPmtpR$.

\medskip

Hence, in order to prove the proposition, it suffices to show that if we apply
$$(\on{Id}\otimes \on{C}^\cdot(\Ran,-)):\Dmod_{\det^{\frac{1}{2}}_{G,M}}(\BunPmt\times \Ran)\to \Dmod_{\det^{\frac{1}{2}}_{G,M}}(\BunPmt)$$
to \eqref{e:loc IC semiinf}, viewed as an object of $\Dmod_{\det^{\frac{1}{2}}_{G,M}}(\BunPmt\times \Ran)$, the result is isomorphic to
$\IC^{-,\semiinf,\on{glob}}$.

\medskip

However, the latter is the assertion of \cite[Theorems 6.3.2 and 3.4.4]{Ga5} (adapted to the case of an arbitrary parabolic).

\qed[\propref{p:compact Eis}]

\section{Interaction of parabolic induction with the Whittaker model} \label{s:Eis and Whit}

In this section we will perform the computation of Whittaker coefficients of Eisenstein series. When combined 
with its spectral counterpart, this will ensure the compatibility of the Langlands functor with the Eisenstein functors. 

\ssec{Whittaker coefficients of Eisenstein series}

In this subsection we will state several assertions, all equivalent to the calculation of Whittaker coefficients of 
(usual, i.e., unenhanced) Eisenstein series. 

\sssec{}

We will prove:

\begin{thm} \label{t:CT co of Poinc}
For $\CZ\to \Ran$, the following diagram of functors commutes:
\begin{equation}  \label{e:CT co of Poinc *}
\CD
\Whit_*(G)_\CZ  @>{\on{ins.unit}_\CZ}>> \Whit_*(G)_{\CZ^{\subseteq}}  \\
@VV{\on{Poinc}_{G,*,\CZ}[2\delta_{(N_P)_{\rho_P(\omega_X)}}+\delta_{(N^-_P)_{\rho_P(\omega_X)}}]}V
@VV{J^{-,!}_{\Whit,\Theta}}V \\
\Dmod_{\frac{1}{2}}(\Bun_G)_{\on{co}}\otimes \Dmod(\CZ) & &  \Whit_*(M)_{\CZ^{\subseteq}}  \\
@VV{\on{CT}^-_{\on{co},?,\rho_P(\omega_X)}\otimes \on{Id}}V @VV{\on{Poinc}_{M,*,\CZ^{\subseteq}}}V \\
\Dmod_{\frac{1}{2}}(\Bun_M)_{\on{co}}\otimes \Dmod(\CZ) @<<{\on{Id}\otimes (\on{pr}_{\on{small},\CZ})_!}< 
\Dmod_{\frac{1}{2}}(\Bun_M)_{\on{co}}\otimes \Dmod(\CZ^{\subseteq}),
\endCD
\end{equation}

\noindent where:
\begin{itemize}

\item The notation $\on{ins.unit}_\Ran$ refers to the unital structure on $\ul{\Whit_*(G)}$, see \cite[Sect. 11.2.8]{GLC2};

\medskip

\item When forming $\on{Poinc}_{M,*,\CZ^{\subseteq}}$, we regard $\CZ^\subseteq$ as mapping to $\Ran$ by $\on{pr}_{\on{big}}$, 
see \cite[Sect. 11.2.6]{GLC2};

\medskip

\item The functor $J^{-,!}_{\Whit,\Theta}$ is one from \eqref{e:J Theta Whit};

\medskip

\item $\on{CT}^-_{\on{co},?,\rho_P(\omega_X)}:=(\Eis^-_{!,\rho_P(\omega_X)})^\vee$;

\medskip

\item $\delta_{(N_P)_{\rho_P(\omega_X)}}:=\dim(\Bun_{(N_P)_{\rho_P(\omega_X)}})$ and $\delta_{(N^-_P)_{\rho_P(\omega_X)}}:=\dim(\Bun_{(N^-_P)_{\rho_P(\omega_X)}})$. 

\end{itemize}

\end{thm} 

\sssec{}

Note that the unital property of the functor $\ul{\on{Poinc}}_{G,!}$ implies that 
that the essential image of the functor
$$\on{coeff}_G:\Dmod_{\frac{1}{2}}(\Bun_G) \to \Whit^!(G)_\Ran$$
belongs to the full subcategory
$$\Whit^!(G)_{\Ran^{\on{untl}},\on{indep}}\subset \Whit^!(G)_\Ran,$$
see \cite[Sect. H.1]{GLC2}.

\medskip

When viewed as such, we will denote it, as well as its composition with the embedding 
$$\Whit^!(G)_{\Ran^{\on{untl}},\on{indep}}\hookrightarrow \Whit^!(G)_{\Ran^{\on{untl}}},$$
by $\on{coeff}_{G,\on{untl}}$. 

\sssec{}

Recall now that the functor 
$$\on{emb.indep}_{\Whit^!(G)}:\Whit^!(G)_{\Ran^{\on{untl}},\on{indep}}\overset{\on{embed.indep}_{\Whit^!(G)}}\hookrightarrow \Whit^!(G)_{\Ran^{\on{untl}}},$$
admits a left adjoint, see \cite[Lemma H.1.10]{GLC2}.

\medskip

We claim:

\begin{cor} \label{c:coeff of Eis !}
The following diagram of functors commutes:
$$
\CD
\Whit^!(G)_{\Ran^{\on{untl}},\on{indep}} @<{\on{emb.indep}_{\Whit^!(G)}^L}<<  \Whit^!(G)_{\Ran^{\on{untl}}} 
@<{(J^{-,!}_{\Whit,\Theta})^\vee}<<  \Whit^!(M)_{\Ran^{\on{untl}}} \\
@AA{\on{coeff}_{G,\on{untl}}}[2\delta_{(N_P)_{\rho_P(\omega_X)}}+\delta_{(N^-_P)_{\rho_P(\omega_X)}}]A & & @AA{\on{coeff}_{M,\on{untl}}}A  \\
\Dmod_{\frac{1}{2}}(\Bun_G) & @<{\Eis^-_{!,\rho_P(\omega_X)}}<<  & \Dmod_{\frac{1}{2}}(\Bun_M).
\endCD
$$
\end{cor}

\begin{proof}

Passing to the dual functors in \eqref{e:CT co of Poinc *} with $\CZ=\Ran$, we obtain that it suffices to show that the functor
\begin{multline}   \label{e:coeff of Eis ! 1}
\Dmod_{\frac{1}{2}}(\Bun_M) \overset{\on{Id}\otimes \omega_{\Ran^{\subseteq}}}\longrightarrow \Dmod_{\frac{1}{2}}(\Bun_M) \otimes \Dmod(\Ran^{\subseteq})
\overset{\on{coeff}_{M,\Ran^{\subseteq}}}\longrightarrow \\
\to \Whit^!(M)_{\Ran^{\subseteq}} \overset{(J^{-,!}_{\Whit,\Theta})^\vee}\longrightarrow \Whit^!(G)_{\Ran^{\subseteq}} 
\overset{(\on{ins.unit}_\Ran)^\vee}\longrightarrow \Whit^!(G)_\Ran
\end{multline}
takes values in $\Whit^!(G)_{\Ran^{\on{untl}},\on{indep}}$, and as such identifies with
\begin{multline}   \label{e:coeff of Eis ! 2}
\Dmod_{\frac{1}{2}}(\Bun_M) \overset{\on{coeff}_{M,\on{untl}}}\longrightarrow \Whit^!(M)_{\Ran^{\on{untl}}}  \overset{(J^{-,!}_{\Whit,\Theta})^\vee}\longrightarrow \\
\to \Whit^!(G)_{\Ran^{\on{untl}}} \overset{\on{emb.indep}_{\Whit^!(G)}^L}\longrightarrow \Whit^!(G)_{\Ran^{\on{untl}},\on{indep}}. 
\end{multline} 

\medskip

Note that the composition of the first two arrows in \eqref{e:coeff of Eis ! 1} can be rewritten as
$$\Dmod_{\frac{1}{2}}(\Bun_M) \overset{\on{coeff}_{M,\on{untl}}}\longrightarrow \Whit^!(M)_{\Ran^{\on{untl}}}  \to \Whit^!(M)_{\Ran} 
\overset{(\on{pr}_{\on{big}})^!}\longrightarrow 
\Whit^!(M)_{\Ran^{\subseteq}},$$
where the second arrow is the forgetful functor.

\medskip

Hence, it suffices to establish an isomorphism between
\begin{multline*}  
\Whit^!(M)_{\Ran^{\on{untl}}}  \to \Whit^!(M)_{\Ran}  \overset{(\on{pr}_{\on{big}})^!}\longrightarrow \\
\to \Whit^!(M)_{\Ran^{\subseteq}} \overset{(J^{-,!}_{\Whit,\Theta})^\vee}\longrightarrow \Whit^!(G)_{\Ran^{\subseteq}} 
\overset{(\on{ins.unit}_\Ran)^\vee}\longrightarrow \Whit^!(G)_\Ran,
\end{multline*}
which is the same as 
\begin{multline*}  
\Whit^!(M)_{\Ran^{\on{untl}}} \overset{(J^{-,!}_{\Whit,\Theta})^\vee}\longrightarrow \Whit^!(G)_{\Ran^{\on{untl}}} \to 
\Whit^!(G)_{\Ran}  \overset{(\on{pr}_{\on{big}})^!}\longrightarrow \\
\to \Whit^!(G)_{\Ran^{\subseteq}} \overset{(\on{ins.unit}_\Ran)^\vee}\longrightarrow \Whit^!(G)_\Ran,
\end{multline*}
and 
$$\Whit^!(M)_{\Ran^{\on{untl}}}  \overset{(J^{-,!}_{\Whit,\Theta})^\vee}\longrightarrow \Whit^!(G)_{\Ran^{\on{untl}}} \overset{\on{emb.indep}_{\Whit^!(G)}^L}\longrightarrow \Whit^!(G)_{\Ran^{\on{untl}},\on{indep}}.$$

Note that the functor $\Whit^!(G)_{\Ran^{\subseteq}} \overset{(\on{ins.unit}_\Ran)^\vee}\longrightarrow \Whit^!(G)_\Ran$ is the composition
$$\Whit^!(G)_{\Ran^{\subseteq}}:=\Whit^!(G)_{\Ran^{\subseteq},\on{pr}_{\on{big}}} \to \Whit^!(G)_{\Ran^{\subseteq},\on{pr}_{\on{small}}} 
\overset{(\on{pr}_{\on{small}})_!}\longrightarrow \Whit^!(G)_\Ran,$$
where the second arrow is the dual of to the functor
\begin{equation}   \label{e:coeff of Eis ! 3}
\Whit_*(G)_{\Ran^{\subseteq},\on{pr}_{\on{small}}} \to  \Whit_*(G)_{\Ran^{\subseteq},\on{pr}_{\on{big}}}
\end{equation} 
defined by the unital structure. Note also that the dual of \eqref{e:coeff of Eis ! 3} identifies canonically with the right adjoint of
\begin{equation}   \label{e:coeff of Eis ! 4}
\Whit^!(G)_{\Ran^{\subseteq},\on{pr}_{\on{small}}} \to  \Whit^!(G)_{\Ran^{\subseteq},\on{pr}_{\on{big}}}.
\end{equation} 

Now, the required isomorphism follows from \lemref{l:indep L} below.

\end{proof} 

\begin{lem} \label{l:indep L}
In the context of \cite[Lemma H.1.10]{GLC2} assume that $\ul\bC$ is dualizable as a crystal of categories over $\Ran_{\on{untl}}$. 
Then the functor left adjoint to 
$$\on{emb.indep}:\bC^{\on{loc}}_{\Ran^{\on{untl}},\on{indep}}\hookrightarrow \bC^{\on{loc}}_{\Ran^{\on{untl}}},$$
composed with the embedding $\bC^{\on{loc}}_{\Ran^{\on{untl}},\on{indep}}\hookrightarrow \bC^{\on{loc}}_\Ran$, 
is given by
$$\bC^{\on{loc}}_{\Ran^{\on{untl}}}\to \bC^{\on{loc}}_\Ran \overset{(\on{pr}_{\on{big}})^!}\longrightarrow \bC^{\on{loc}}_{\Ran^\subseteq,\on{pr}_{\on{big}}}
\to \bC^{\on{loc}}_{\Ran^\subseteq,\on{pr}_{\on{small}}}\overset{(\on{pr}_{\on{small}})_!}\longrightarrow \bC^{\on{loc}}_\Ran,$$
where the third arrow is the functor right adjoint to the functor 
$$\bC^{\on{loc}}_{\Ran^\subseteq,\on{pr}_{\on{small}}}\to \bC^{\on{loc}}_{\Ran^\subseteq,\on{pr}_{\on{big}}},$$
defined by the unital structure on $\ul\bC^{\on{loc}}$.
\end{lem}

\begin{proof}

Follows from \cite[Proposition I.2.3]{GLC2} applied to the dual of $\bC^{\on{loc}}$. 

\end{proof} 

\sssec{}

We now deduce another consequence of \thmref{t:CT co of Poinc}, this time for the computation of the constant term of
the !-Poincar\'e functor:

\begin{cor} \label{c:CT Poinc !}
For $\CZ\to \Ran$, the following diagram of functors commutes:
$$
\CD
\Whit^!(G)_\CZ  @>{\on{ins.unit}_\CZ}>> \Whit^!(G)_{\CZ^{\subseteq}}  \\
@VV{\on{Poinc}_{G,!,\CZ}[d]}V @VV{J^{-,!}_{\Whit,\tau}}V \\
\Dmod_{\frac{1}{2}}(\Bun_G)\otimes \Dmod(\CZ) & &  \Whit^!(M)_{\CZ^{\subseteq}}  \\
@VV{\on{CT}_{*,\rho_P(\omega_X)}\otimes \on{Id}}V @VV{\on{Poinc}_{M,!,\CZ^{\subseteq}}}V \\
\Dmod_{\frac{1}{2}}(\Bun_M)\otimes \Dmod(\CZ) @<<{\on{Id}\otimes (\on{pr}_{\on{small},\CZ})_!}< 
\Dmod_{\frac{1}{2}}(\Bun_M)\otimes \Dmod(\CZ^{\subseteq}),
\endCD,
$$
where:

\begin{itemize}

\item $\on{CT}_{*,\rho_P(\omega_X)}\simeq (\on{transl}_{2\rho_P(\omega_X)})^*\circ \tau_M\circ \on{CT}^-_{*,\rho_P(\omega_X)}\circ \tau_G\simeq
\tau_M\circ (\on{transl}_{-2\rho_P(\omega_X)})^*\circ \on{CT}^-_{*,\rho_P(\omega_X)}\circ \tau_G;$

\medskip

\item $J^{-,!}_{\Whit,\tau}=\tau_M\circ J^{-,!}_{\Whit}\circ \tau_G$;

\medskip

\item $d=-2\delta_{(N_P)_{\rho_P(\omega_X)}}-\delta_{(N^-_P)_{\rho_P(\omega_X)}}+2\delta_{N^-_P}$. 

\end{itemize}

\end{cor}

\sssec{}

Before we prove \corref{c:CT Poinc !}, we record the following result from \cite[Theorem 1.1.6]{Lin}\footnote{In order to apply
the result from {\it loc. cit.} to deduce \thmref{t:Kevin}, one needs to combine it with \thmref{t:local-to-global Whit}(c) below.}:

\begin{thm} \label{t:Kevin}
The diagram 
\begin{equation} \label{e:Mir and Poinc}
\CD
\Whit^!(G)_\Ran @<{\Theta_{\Whit(G)}}<< \Whit_*(\Gr_G)_\Ran \\
@V{\Poinc_{G,!}}VV @VV{\Poinc_{G,*}[2\delta_{N_{\rho(\omega_X)}}]}V \\
\Dmod_{\frac{1}{2}}(\Bun_G) @<{\Mir_{\Bun_G}[2\delta_G]}<< \Dmod_{\frac{1}{2}}(\Bun_G) _{\on{co}},
\endCD
\end{equation} 
commutes, where 
$$\delta_G=\dim(\Bun_G).$$
\end{thm}

\begin{proof}[Proof of \corref{c:CT Poinc !}]

Obtained juxtaposing \thmref{t:CT co of Poinc}, \corref{c:intertwiners CT simple} and diagrams \eqref{e:Mir and Poinc} for $G$ and $M$. 

\end{proof}

\sssec{}

In the next several subsections we will recall the definitions of several geometric tools needed for the proof
of \thmref{t:CT co of Poinc}.

%

\ssec{The global model for the Whittaker category}

\sssec{}

Let 
$$\BunNb \overset{\ol\sfp}\to \Bun_G.$$
be as in \cite[Sect. 4.1]{Ga6}. 

\medskip

For $\CZ\to \Ran$, we let $\BunNbZ$ be a version of $\BunNb$, where we allow poles at the specified
points of the curve. When $\CZ=S$ is a scheme, the prestack $\BunNbS$ is an ind-algebraic stack. 

\medskip

We will denote by the same symbol $\ol\sfp$ the projection
$$\BunNbZ\to \Bun_G,$$
and by $\ol\sfp_\CZ$ the projection
$$\BunNbZ\to \Bun_G\times \CZ.$$

\medskip

We will denote by
$$\Dmod_{\frac{1}{2}}(\BunNbZ) \text{ and } \Dmod_{\frac{1}{2}}(\BunNbZ)_{\on{co}}$$
the corresponding categories of twisted D-modules, by the means of the $\mu_2$-gerbe $\pi^*(\det^{\frac{1}{2}}_{\Bun_G})$,
see Sects. \ref{sss:cats on comp begin}-\ref{sss:cats on comp end}. 

\sssec{}

Following \cite[Sect. 4.7]{Ga6}, one singles out a full subcategory
$$\Whit^{\on{glob}}(G)_{\CZ}\subset \Dmod_{\frac{1}{2}}(\BunNbZ).$$

\sssec{} \label{sss:Whit glob as co}

%
%
%
%

Let $\CZ = S$ be a scheme. Then as in \secref{sss:cats on comp begin}, we can write $\BunNbS$ as a union
of algebraic stacks $\CY_{S,i}$. We can write each $\CY_{S,i}$ as a union of opens $U \subset \CY_{S,i}$
stable under the loop groupoid from \cite[Sect. 4.4]{Ga6}. 

\medskip

Hence, we can make sense of
$$\Whit^{\on{glob}}(G)_{S,U}\subset \Dmod_{\frac{1}{2}}(U),$$
by imposing the Whittaker equivariance condition as in \cite[Sect. 4.7]{Ga6}. 

\medskip

Then for 
$$\Whit^{\on{glob}}(G)_i:=\Whit^{\on{glob}}(G)_{S}\cap \Dmod_{\frac{1}{2}}(CY_{S,i}),$$
we have 
$$\Whit^{\on{glob}}(G)_i\simeq \underset{i}{\on{lim}}\,  \Whit^{\on{glob}}(G)_{S,U},$$
where the limit is formed using the restriction functors. 

\medskip

We have
$$\Whit^{\on{glob}}(G)_{S}\simeq \underset{i}{\on{lim}}\, \Whit^{\on{glob}}(G)_i\simeq
 \underset{i}{\on{colim}}\, \Whit^{\on{glob}}(G)_i,$$
 where:
 
 \begin{itemize}
 
 \item The limit is formed using !-pullback functors;
 
 \item The colimit is formed using the *-pushforward functors. 
 
 \end{itemize}


\sssec{}

Define
$$\Whit^{\on{glob}}(G)_{\on{co},i} =
\underset{U}{\on{colim}} \, \Whit^{\on{glob}}(G)_{S,U},$$

\noindent where the colimit is formed using $*$-pushforward functors,
and the index runs over open substacks U as above.

\medskip

Define
$$\Whit^{\on{glob}}(G)_{\on{co},S} =
\underset{i}{\on{colim}}\,  \Whit^{\on{glob}}(G)_{\on{co},i},$$

\noindent where the colimit is formed using $*$-pushforward functors
(this time over closed embeddings). 

\medskip 

We can rewrite 
$$\Whit^{\on{glob}}(G)_{\on{co},S} \simeq
\underset{i}{\on{lim}}\,  \Whit^{\on{glob}}(G)_{\on{co},i},$$

\noindent where the limit is formed using !-pullback functors.

\medskip

The functor\footnote{See \cite[Sect. 4.4.2]{DG2}.} $\on{Ps-Id}^{\on{nv}}_{CY_{S,i}}$ gives rise to a functor
\begin{equation} \label{e:Theta i glob}
\Whit^{\on{glob}}(G)_{\on{co},i}\to \Whit^{\on{glob}}(G)_i.
\end{equation}

The functors \eqref{e:Theta i glob} combine to a functor
$$\Theta^{\on{glob}}_{\Whit(G)}:\Whit^{\on{glob}}(G)_{\on{co},S}\overset{\sim}\to \Whit^{\on{glob}}(G)_S.$$

\medskip 

We define the category $\Whit^{\on{glob}}(G)_{\on{co},\CZ}$ and the functor $\Theta^{\on{glob}}_{\Whit(G)}$ 
for a general prestack $\CZ$ mapping to $\Ran$ by right Kan extending from the case of schemes.

\sssec{}

We now note the following property of the Whittaker subcategory (see \cite[Corollary 4.7.9(a)]{Ga6}:

\medskip

For every $i$, there exists a quasi-compact open $U_0\subset \CY_{S,i}$, such for any $U\subset U_0$,
any object of $\Whit^{\on{glob}}(G)_{S,U}$ is a \emph{clean}\footnote{Clean=both $!$ and $*$} extension
of its restriction to $U_0$ along $U_0\hookrightarrow U$. 

\medskip

In particular, we obtain that the *-extension functor
$$\Whit^{\on{glob}}(G)_{S,U_0}\to \Whit^{\on{glob}}(G)_{\on{co},i}$$
is already an equivalence. 

\medskip

It follows formally that the functors \eqref{e:Theta i glob}, and hence $\Theta^{\on{glob}}_{\Whit(G)}$, 
are equivalences.

%
%
%

%
%
%
%

\sssec{}

As in \secref{sss:Gr to BunPt}, we have a naturally defined map 
$$\pi_\CZ:\Gr_{G,\rho(\omega_X),\CZ}\to \BunNbZ.$$

\medskip

Thanks to the existence of a canonical square root of the line $\fl_{G,N_{\rho(\omega_X)}}$
(see \cite[Sect. 9.1.1]{GLC1}), the pullback of the $\mu_2$-gerbe $\det^{\frac{1}{2}}_{\Bun_G}$ along $\pi_\CZ$ identifies with the
$\mu_2$-gerbe $\det^{\frac{1}{2}}_{\Gr_{G,\rho(\omega_X)}}$ on $\Gr_{G,\rho(\omega_X),\CZ}$, see
\cite[Sect. 9.3.1]{GLC2}. 

\medskip

Consider the resulting functor
$$(\pi_\CZ)_*:\Dmod_{\frac{1}{2}}(\Gr_{G,\rho(\omega_X),\CZ})\to \Dmod_{\frac{1}{2}}(\BunNbZ)_{\on{co}}.$$

As in Sects. \ref{sss:push Av 1}-\ref{sss:push Av 2}, from $(\pi_\CZ)_*$ we produce a functor
$$\Dmod_{\frac{1}{2}}(\Gr_{G,\rho(\omega_X),\CZ})_{\fL(N)_{\rho(\omega)},\chi}=:\Whit_*(G)_\CZ
\overset{((\pi_\CZ)_*)^{\on{Av},\chi}}\longrightarrow \Dmod_{\frac{1}{2}}(\BunNbZ)_{\on{co}}.$$

\medskip

Note that by the definition of the functor $\on{Poinc}_{G,*,\CZ}$, we have
$$\on{Poinc}_{G,*,\CZ}\simeq (\ol\sfp_\CZ)_*\circ ((\pi_\CZ)_*)^{\on{Av},\chi}.$$

\sssec{}

The following is \cite[Theorem 5.1.4, Proposition 5.4.2, Corollary 5.4.5]{Ga6}:

\begin{thm} \label{t:local-to-global Whit} \hfill

\smallskip

\noindent{\em(a)} The functor $((\pi_\CZ)_*)^{\on{Av},\chi}$ takes values in $\Whit^{\on{glob}}(G)_{\on{co},\CZ}$, and the resulting functor
$$\Whit_*(G)_\CZ \overset{((\pi_\CZ)_*)^{\on{Av},\chi}}\longrightarrow \Whit^{\on{glob}}(G)_{\on{co},\CZ}$$
is an equivalence.

\smallskip

\noindent{\em(b)} The restriction of the functor $\pi_\CZ^!:\Dmod_{\frac{1}{2}}(\BunNbZ) \to \Dmod_{\frac{1}{2}}(\Gr_{G,\rho(\omega_X),\CZ})$
to 
$$\Whit^{\on{glob}}(G)_\CZ\subset \Dmod_{\frac{1}{2}}(\BunNbZ)$$ takes values in 
\begin{equation} \label{e:Whit as subcat}
\Whit^!(G)_\CZ:=\Dmod_{\frac{1}{2}}(\Gr_{G,\rho(\omega_X),\CZ})^{\fL(N)_{\rho(\omega),\CZ},\chi}\subset 
\Dmod_{\frac{1}{2}}(\Gr_{G,\rho(\omega_X),\CZ}),
\end{equation} 
and the resulting functor
\begin{equation} \label{e:global-to-local Whit !}
\pi_\CZ^!:\Whit^{\on{glob}}(G)_\CZ\to \Whit^!(G)_\CZ
\end{equation} 
is an equivalence.

\medskip

\noindent{\em(b')} The inverse to the functor \eqref{e:global-to-local Whit !} is given by the restriction of the (partially defined) left adjoint
$(\pi_\CZ)_!$ of $\pi_\CZ^!$ to \eqref{e:Whit as subcat}.

\medskip

\noindent{\em(c)} The composition 
$$\pi_\CZ^!\circ \Theta^{\on{glob}}_{\Whit(G)}\circ ((\pi_\CZ)_*)^{\on{Av},\chi}:\Whit_*(G)_\CZ\to \Whit^!(G)_\CZ$$
is canonically isomorphic to $\Theta_{\Whit(G)}[-2\delta_{N_{\rho(\omega_X)}}]$.

\end{thm}

\sssec{}

As in Sects. \ref{sss:* BunPt unital 1}-\ref{sss:* BunPt unital 2} and \ref{sss:! BunPt unital} the assignments
$$\CZ\rightsquigarrow \Dmod_{\frac{1}{2}}(\BunNbZ) \text{ and }
\CZ\rightsquigarrow \Dmod_{\frac{1}{2}}(\BunNbZ)_{\on{co}}$$
have natural structures of unital crystals of categories over $\Ran$, and the assignments
$$\CZ\rightsquigarrow  ((\pi_\CZ)_*)^{\on{Av},\chi} \text{ and } \CZ\rightsquigarrow (\pi_\CZ)_!$$
are left-lax and right-lax unital functors, to be denoted,
$$\ul{\pi_*}^{\on{Av},\chi} \text{ and } \ul{\pi_!},$$
respectively.

\medskip

However, as in \propref{p:* BunPt unital} (resp., \corref{c:! BunPt unital}) one shows that these
structures are strictly unital. 

\ssec{The space of \emph{generic} reductions} 

\sssec{} \label{sss:gen red}

Let 
$$\BunNbg:=\Bun_{N,\rho(\omega_X),\on{gen}}\underset{\Bun_{G,\on{gen}}}\times \Bun_G$$ be the space of $G$-bundles 
equipped with a \emph{generic} reduction to $N_{\rho(\omega_X)}$ (see \cite{Bar} or \cite[Appendix A]{Ga5}, where the definition in families
is given). Denote by $\ol\sfp_{\on{gen}}$ the projection
$$\BunNbg\to \Bun_G.$$

\medskip

Note that we have a naturally defined map
\begin{equation} \label{e:gen vs Ran}
\sff:\BunNbR\to \BunNbg.
\end{equation} 

The following is established in \cite[Theorem A.1.10]{Ga5}:

\begin{lem} \label{l:UHA}
The map $\sff$ of \eqref{e:gen vs Ran} is pseudo-proper and \emph{universally homologically acyclic}.\footnote{See \cite[Sect. A.1.8]{Ga5} for what this means.}
\end{lem}

\sssec{}

Similarly, for $\CZ\to \Ran$ we have a map
\begin{equation} \label{e:gen vs Ran Z}
\sff_\CZ: \BunNbZsub\to  \BunNbg\times \CZ,
\end{equation} 
and parallel to \lemref{l:UHA}, one shows that the map \eqref{e:gen vs Ran Z} is also universally homologically acyclic. 

\ssec{Closures of semi-infinite orbits}

\sssec{} \label{sss:Gr +}

Let $\Gr^{\on{pos}}_G\subset \Gr_G$ be the factorization subspace equal to the closure of the unit orbit of
$\fL(N_P)$, see \cite[Sect. 6.1.5]{Ga6}. 

\medskip

We can regard $\Gr_G$ as a factorization \emph{module} space over $\Gr^{\on{pos}}_G$. For $\CZ\to \Ran$,
denote by $\Gr^{\on{pos}}_{G,\CZ^\subseteq}$ the corresponding prestack over $\CZ^\subseteq$, i.e.,
$$\Gr^{\on{pos}}_{G,\CZ^\subseteq}\underset{\CZ^\subseteq}\times \CZ=\Gr_{G,\CZ}$$
and 
$$\Gr^{\on{pos}}_{G,\CZ^\subseteq}\underset{\CZ^\subseteq}\times (\Ran\times \CZ^\subseteq)_{\on{disj}}\simeq
(\Gr^{\on{pos}}_{G,\Ran}\times \Gr^{\on{pos}}_{G,\CZ^\subseteq})\underset{\Ran\times \CZ^\subseteq}\times (\Ran\times \CZ^\subseteq)_{\on{disj}}.$$

\sssec{}  \label{sss:Gr + BunNb}

Note that we have a Cartesian square
$$
\CD
\Gr^{\on{pos}}_{G,\rho(\omega)_X,\Ran} @>>> \Gr_{G,\rho(\omega)_X,\Ran} \\
@V{\pi^{\on{pos}}_\Ran}VV @VV{\pi_\Ran}V \\
\BunNb\times \Ran @>>> \BunNbR,
\endCD
$$
and, more generally, for $\CZ\to \Ran$:
$$
\CD
\Gr^{\on{pos}}_{G,\rho(\omega)_X,\CZ^{\subseteq}} @>>> \Gr_{G,\rho(\omega)_X,\CZ^{\subseteq}} \\
@V{\pi^{\on{pos}}_{\CZ^{\subseteq}}}VV @VV{\pi_{\CZ^{\subseteq}}}V \\
\BunNbZ\underset{\CZ}\times \CZ^{\subseteq} @>>> \ol\Bun_{\rho(\omega)_X,\CZ^{\subseteq}} \\
@V{\on{id}\times \on{pr}_{\on{small},\CZ}}VV \\
\BunNbZ.
\endCD
$$

\sssec{}

We claim:

\begin{lem} \label{l:unital local and global Whit}
For $\CF\in \Whit^{\on{glob}}(G)_\CZ$, the object
$$(\pi^{\on{pos}}_{\CZ^{\subseteq}})^! \circ (\on{id}\times \on{pr}_{\on{small},\CZ})^!(\CF)\in 
\Dmod_{\frac{1}{2}}(\Gr^{\on{pos}}_{G,\rho(\omega)_X,\CZ^{\subseteq}})\subset \Dmod_{\frac{1}{2}}(\Gr_{G,\rho(\omega)_X,\CZ^{\subseteq}}),$$
viewed as an object of $\Whit^!(G)_{\CZ^{\subseteq}}$, identifies canonically with
$$\on{ins.unit}_\CZ(\pi^!_\CZ(\CF)).$$
\end{lem} 

\begin{proof}

To simplify the notation, we will assume that $\CZ=\on{pt}$, corresponding to $\ul{x}\in \Ran$, so that
$$\CZ^\subseteq=:\Ran_{\ul{x}}.$$

\medskip

The unital structure on $\Gr_G$ as a factorization space gives rise to a map
$$\on{ins.unit}_{\Gr,\ul{x}}:\Gr_{G,\rho(\omega)_X,\ul{x}}\times \Ran_{\ul{x}}\to \Gr_{G,\rho(\omega)_X,\Ran_{\ul{x}}}.$$

It follows from \cite[Theorem 6.2.5]{Ga6} that the functor of !-pullback along $\on{ins.unit}_{\Gr,\ul{x}}$ gives rise to an equivalence
between the full subcategory of $\Whit^!(G)_{\Ran_{\ul{x}}}$ consisting of objects supported on
$$\Gr^{\on{pos}}_{G,\rho(\omega)_X,\Ran_{\ul{x}}} \subset \Gr_{G,\rho(\omega)_X,\Ran_{\ul{x}}}$$
and 
$$\Whit^!(G)_{\ul{x}}\otimes \Dmod(\Ran_{\ul{x}}).$$

\medskip

Hence, it suffices to construct an isomorphism between 
$$(\on{ins.unit}_{\Gr,\ul{x}})^!\circ (\pi^{\on{pos}}_{\Ran_{\ul{x}}})^! \circ (\on{id}\times \on{pr}_{\on{small},\ul{x}})^!(\CF) \text{ and }
(\on{ins.unit}_{\Gr,\ul{x}})^!\circ \on{ins.unit}_{\ul{x}}(\pi^!_\CZ(\CF)).$$

First, we tautologically have:
$$(\on{ins.unit}_{\Gr,\ul{x}})^!\circ (\pi^{\on{pos}}_{\Ran_{\ul{x}}})^! \circ (\on{id}\times \on{pr}_{\on{small},\ul{x}})^!(\CF)\simeq
\pi^!_{\ul{x}}(\CF)\boxtimes \omega_{\Ran_{\ul{x}}}.$$

\medskip

By the definition of the unital structure on $\Whit^!(G)$, the functor
$$\on{ins.unit}_{\ul{x}}:\Whit^!(G)_{\ul{x}}\to \Whit^!(G)_{\Ran_{\ul{x}}}$$
is given by
$$\on{Av}_!^{\fL(N)_{\rho(\omega_X),\Ran_{\ul{x}}},\chi}\circ (\on{ins.unit}_{\Gr,\ul{x}})_*((-)\boxtimes \omega_{\Ran_{\ul{x}}}),$$
where $\on{Av}_!^{\fL(N)_{\rho(\omega_X),\Ran_{\ul{x}}},\chi}$ is the (partially defined) left adjoint to the embedding
$$\Whit^!(G)_{\Ran_{\ul{x}}}\hookrightarrow \Dmod_{\frac{1}{2}}(\Gr_{G,\rho(\omega)_X,\Ran_{\ul{x}}}).$$

Now, by \cite[Theorem 6.4.8]{Ga6}, the functor 
$$\on{Av}_!^{\fL(N)_{\rho(\omega_X),\Ran_{\ul{x}}},\chi}\circ (\on{ins.unit}_{\Gr,\ul{x}})_*:
\Whit^!(G)_{\ul{x}}\otimes \Dmod(\Ran_{\ul{x}})\to  \Whit^!(G)_{\Ran_{\ul{x}}}$$
is a right inverse of $(\on{ins.unit}_{\Gr,\ul{x}})^!$.

\medskip

Hence, we obtain 
$$(\on{ins.unit}_{\Gr,\ul{x}})^!(\on{ins.unit}_\CZ(\pi^!_{\ul{x}}(\CF)))\simeq \pi^!_{\ul{x}}(\CF)\boxtimes \omega_{\Ran_{\ul{x}}},$$
as desired. 

\end{proof} 

\sssec{} \label{sss:Gr neg}

Let 
$$(\Gr_G\times \Gr_M)^{\on{neg}}\subset \Gr_G\times \Gr_M$$
be the factorization space equal to the closure of the unit orbit of $\fL(P^-)$.

\medskip 

As in \secref{sss:Gr +}, for $\CZ\to \Ran$, we can form the space 
$$(\Gr_G\times \Gr_M)^{\on{neg}}_{\CZ^\subseteq}.$$

Furthermore, as in \secref{sss:Gr + BunNb}, we have a natural projection
$$\pi_\Ran:(\Gr_G\times \Gr_M)_\Ran\to \BunPmtpR,$$
so that we have Cartesian squares
$$
\CD
(\Gr_G\times \Gr_M)^{\on{neg}}_\Ran @>>> (\Gr_G\times \Gr_M)_\Ran \\
@V{\pi^{\on{neg}}_\Ran}VV @VV{\pi_\Ran}V \\
\BunPmt\times \Ran @>>> \BunPmtpR
\endCD
$$
and 
$$
\CD
(\Gr_G\times \Gr_M)^{\on{neg}}_{\CZ^\subseteq} @>>> (\Gr_G\times \Gr_M)_{\CZ^\subseteq} \\
@V{\pi^{\on{neg}}_{\CZ^\subseteq}}VV @VV{\pi_{\CZ^\subseteq}}V \\
\BunPmtpZ\underset{\CZ}\times \CZ^\subseteq @>>> \BunPmtpZsub. 
\endCD
$$

\sssec{}

Denote by $\pi^{\on{neg}}$ the composition
$$(\Gr_G\times \Gr_M)^{\on{neg}}_\Ran \overset{\pi^{\on{neg}}_\Ran}\longrightarrow \BunPmt\times \Ran\to \BunPmt.$$

\medskip

Note also that we have a canonical isomorphism of $\mu_2$-gerbes on $(\Gr_G\times \Gr_M)^{\on{neg}}_\Ran$
$$(\pi^{\on{neg}})^*(\det_{G,M}^{\frac{1}{2}})\simeq \det_{\Gr_G}^{\frac{1}{2}}\otimes \det_{\Gr_M}^{\frac{1}{2}},$$
and similarly for the map $\pi^{\on{neg}}_{\CZ^\subseteq}$. 

\sssec{} \label{sss:Delta tilde}

Let $\wt\Delta^{-,\semiinf}$ denote the image of $\Delta^{-,\semiinf}$ under 
\begin{multline*}
\on{I}(G,P^-)^{\on{loc}}\to \Dmod_{\frac{1}{2}}(\Gr_{G,\rho(\omega_X)})^{\fL(N^-_P)\cdot \fL^+(M)} \simeq \\
\simeq \Dmod_{\det_{\Gr_G}^{\frac{1}{2}}\otimes \det_{\Gr_M}^{\frac{1}{2}}}(\Gr_G\times \Gr_M)^{\fL(P^-)}
\to \Dmod_{\det_{\Gr_G}^{\frac{1}{2}}\otimes \det_{\Gr_M}^{\frac{1}{2}}}(\Gr_G\times \Gr_M).
\end{multline*}

Note that $\wt\Delta^{-,\semiinf}_\Ran$ is supported on
$$(\Gr_G\times \Gr_M)^{\on{neg}}_\Ran\subset (\Gr_G\times \Gr_M)_\Ran.$$

We have a tautologically defined map\footnote{In fact, one can show that this map is an isomorphism.}
$$(\pi^{\on{neg}})_!(\wt\Delta^{-,\semiinf}_\Ran)\to  j_!(\omega_{\Bun_{P^-}}).$$

Hence, we obtain a map
\begin{equation}  \label{e:pullback Delta}
\wt\Delta^{-,\semiinf}_\Ran\to (\pi^{\on{neg}})^!\circ j_!(\omega_{\Bun_{P^-}}).
\end{equation}

The following assertion is a variant of \cite[Corollary 3.6.6]{Ga5}:

\begin{lem} \label{l:pullback Delta}
The map \eqref{e:pullback Delta} is an isomorphism. 
\end{lem}

\ssec{Zastava spaces}

\sssec{} \label{sss:Zas}

Let $\Zas$ be the \emph{Zastava space}, defined to be the open substack of
$$\BunNb\underset{\Bun_G}\times \BunPmt,$$
corresponding to the condition that the generalized reductions of the given $G$-bundle to $N_{\rho(\omega_X)}$ and $P^-$ 
are transversal \emph{at the generic point} of $X$.

\medskip

For $\CZ\to \Ran$, we let $\Zas_\CZ$ denote the version of $\Zas$, where we allow both generalized
reductions to have poles at the specified point of the curve, i.e., $\Zas_\CZ$ is the open substack of 
$$\BunNbZ\underset{\Bun_G\times \CZ}\times \BunPmtpZ,$$
given by the generic transversality condition. 

\medskip

We will denote by
$${}'\wt\sfp^-_\CZ:\Zas_\CZ\to \BunNbZ \text{ and } {}'\ol\sfp_\CZ:\Zas_\CZ\to \BunPmtpZ$$
the corresponding maps. 

\sssec{}

We will consider the $\mu_2$-gerbes on $\Zas_\CZ$
$$\det_G^{\frac{1}{2}}, \,\, \det_M^{\frac{1}{2}} \text{ and } \det_{G,M}^{\frac{1}{2}},$$
pulled back by means of $'\ol\sfp_\CZ$ from the corresponding gerbes on $\BunPmtpZ$.

\medskip 

We will consider the corresponding categories of twisted D-modules
$$\Dmod_\CG(\Zas_\CZ),\,\, \Dmod_\CG(\Zas_\CZ)_{\on{co}},\,\, \Dmod_\CG(\Zas_\CZ)_{\on{co}/\Bun_M},$$
defined by the recipe of Sects. \ref{sss:cats on comp begin}-\ref{sss:cats on comp end}. 

\sssec{}

Note now that we have a canonically defined map\footnote{Here $N(M)$ is the maximal unipotent subgroup of the Levi $M$.}
$$\fs:\Zas\to \BunNbMg$$
(see \secref{sss:gen red}), so that there is a commutative diagram
$$
\CD
\Zas @>{\fs}>>  \BunNbMg \\
@V{'\ol\sfp}VV @VV{\ol\sfp_{M,\on{gen}}}V \\
\BunPmt @>{\wt\sfq^-}>> \Bun_M.
\endCD
$$

\medskip

Similarly, for $\CZ\to \Ran$, we have a map
$$\fs_\CZ:\Zas_\CZ\to \BunNbMg\times \CZ$$
and a commutative diagram 
\begin{equation} \label{e:Zast to config}
\CD
\Zas_\CZ @>{\fs_\CZ}>>  \BunNbMg\times \CZ \\
@V{'\ol\sfp_\CZ}VV @VV{\ol\sfp_{M,\on{gen}}\times \on{id}_\CZ}V \\
\BunPmtpZ @>{\wt\sfq^-_\CZ}>> \Bun_M\times \CZ
\endCD
\end{equation}

\sssec{}

Denote
$$\Zas^\Ran:=\Zas\underset{\BunNbMg}\times \BunNbMR$$  and 
$$\Zas^\Ran_\CZ:=\Zas_\CZ\underset{\BunNbMg\times \CZ}\times \BunNbMZsub,$$
where 
$$\BunNbMR\to \BunNbMg \text{ and } \BunNbMZsub\to \BunNbMg\times \CZ$$
are the maps \eqref{e:gen vs Ran} and \eqref{e:gen vs Ran Z}, respectively. 

\begin{rem}

The superscript $\Ran$ in the notations $\Zas^\Ran$ and $\Zas^\Ran_\CZ$ is meant to distinguish
these spaces from $\Zas_\Ran$ and $\Ran_{\CZ^\subseteq}$, respectively.

\end{rem} 

\sssec{}

Consider the fiber products
$$
\CD
\Zas^\Ran\underset{\BunNbMR}\times \Gr_{M,\rho_M(\omega_X),\Ran} @>>>  \Gr_{M,\rho(\omega_X),\Ran} \\
@VVV @VV{\pi_\Ran}V \\
\Zas^\Ran @>{'\fs}>> \BunNbMR
\endCD
$$
and
\begin{equation} \label{e:Zast base change}
\CD
\Zas^\Ran_\CZ\underset{\BunNbMZsub}\times \Gr_{M,\rho_M(\omega_X),\CZ^\subseteq} @>>>  \Gr_{M,\rho(\omega_X),\CZ^\subseteq} \\
@VVV @VV{\pi_{\CZ^\subseteq}}V \\
\Zas^\Ran_\CZ @>{'\fs_\CZ}>> \BunNbMZsub.
\endCD
\end{equation} 

\sssec{}

Recall the space $(\Gr_G\times \Gr_M)^{\on{neg}}_{\CZ^\subseteq}$, see \secref{sss:Gr neg}.
The key geometric input in the proof of \thmref{t:CT co of Poinc} is the following observation: 

\begin{lem} \label{l:basic Zastava}
There exists a canonical isomorphism
$$\Zas^\Ran_\CZ\underset{\BunNbMZsub}\times \Gr_{M,\rho_M(\omega_X),\CZ^\subseteq} \simeq 
\Gr^{\on{pos}}_{G,\rho(\omega_X),\CZ^\subseteq}\underset{\Gr_{G,\rho(\omega_X),\CZ^\subseteq}}\times 
(\Gr_{G,\rho(\omega_X)}\times \Gr_{M,\rho(\omega_X)})^{\on{neg}}_{\CZ^\subseteq},$$
so that:

\begin{itemize}

\item The map 
$$\Zas^\Ran_\CZ\underset{\BunNbMZsub}\times \Gr_{M,\rho_M(\omega_X),\CZ^\subseteq} \to
\Zas^\Ran_\CZ\to \Zas_\CZ\overset{'\wt\sfq^-_{\CZ}}\to \BunNbZ$$
identifies with
\begin{multline*}
\Gr^{\on{pos}}_{G,\rho(\omega_X),\CZ^\subseteq}\underset{\Gr_{G,\rho(\omega_X),\CZ^\subseteq}}\times 
(\Gr_{G,\rho(\omega_X)}\times \Gr_{M,\rho(\omega_X)})^{\on{neg}}_{\CZ^\subseteq}\to
\Gr^{\on{pos}}_{G,\rho(\omega_X),\CZ^\subseteq}\overset{\pi^{\on{pos}}_{\CZ^\subseteq}}\longrightarrow \\
\to \BunNbZ\underset{\CZ}\times \CZ^\subseteq\to \BunNbZ;
\end{multline*} 

\item The map 
$$\Zas^\Ran_\CZ\underset{\BunNbMZsub}\times \Gr_{M,\rho_M(\omega_X),\CZ^\subseteq} \to
\Zas^\Ran_\CZ\to \Zas_\CZ\overset{'\ol\sfp^-_{\CZ}}\to \BunPmtpZ$$
identifies with
\begin{multline*}
\Gr^{\on{pos}}_{G,\rho(\omega_X),\CZ^\subseteq}\underset{\Gr_{G,\rho(\omega_X),\CZ^\subseteq}}\times 
(\Gr_{G,\rho(\omega_X)}\times \Gr_{M,\rho(\omega_X)})^{\on{neg}}_{\CZ^\subseteq}\to \\
\to (\Gr_{G,\rho(\omega_X)}\times \Gr_{M,\rho(\omega_X)})^{\on{neg}}_{\CZ^\subseteq} \overset{\pi^{\on{neg}}_{\CZ^\subseteq}}\longrightarrow
\BunPmtpZ\underset{\CZ}\times \CZ^\subseteq\overset{\on{id}\times \on{pr}_{\on{small},\CZ}}\longrightarrow \BunPmtpZ.$$
\end{multline*} 

\end{itemize} 

\end{lem}
 
\ssec{Proof of \thmref{t:CT co of Poinc}} 

\sssec{}

Recall that $j_Z$ denotes the locally closed embedding
$$\Bun_{P^-}\times \CZ\hookrightarrow  \BunPmtpZ.$$

Note that as in \secref{sss:IC as twisted}, we can regard 
$\omega_{\Bun_{P^-}\times \CZ}$ as an object of $\Dmod_{\det^{\frac{1}{2}}_{G,M}}(\Bun_{P^-}\times \CZ)$. Hence, we can regard
$(j_\CZ)_!(\omega_{\Bun_{P^-}\times \CZ})$ as an object of $\Dmod_{\det^{\frac{1}{2}}_{G,M}}( \BunPmtpZ)$. 

\medskip

The functor
$$(\on{CT}^-_{\on{co},?,\rho_P(\omega_X)}\otimes \on{Id})\circ \on{Poinc}_{G,*,\CZ}$$
is by definition
$$(\on{transl}_{\rho_P(\omega_X)})^*\circ (\wt\sfq^-_\CZ)_*\circ \left(\left((\wt\sfp^-_\CZ)^!\circ \on{Poinc}_{G,*,\CZ}(-)\right)
\sotimes (j_\CZ)_!(\omega_{\Bun_{P^-}\times \CZ})\right)[-\on{shift}],$$
which we rewrite as
\begin{equation} \label{e:CT co of Poinc LHS 1}
(\on{transl}_{\rho_P(\omega_X)})^*\circ (\wt\sfq^-_\CZ)_*\circ 
\left(\left((\wt\sfp^-_\CZ)^!\circ (\ol\sfp_\CZ)_*\circ ((\pi_\CZ)_*)^{\on{Av},\chi}(-)\right)\sotimes (j_\CZ)_!(\omega_{\Bun_{P^-}\times \CZ})\right)[-\on{shift}].
\end{equation} 

We depict this diagrammatically as follows:

$$
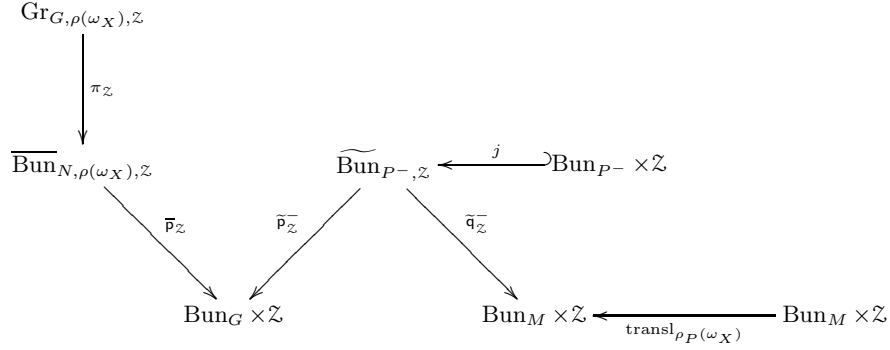
\xy
(20,0)*+{\Bun_M\times \CZ}="X";
(60,0)*+{\Bun_M\times \CZ}="X'";
(-20,0)*+{\Bun_G\times \CZ}="Y";
(0,20)*+{\BunPmtpZ}="Z";
(30,20)*+{\Bun_{P^-}\times \CZ}="W";
(-40,20)*+{\BunNbZ}="U";
(-40,40)*+{\Gr_{G,\rho(\omega_X),\CZ}}="V";
{\ar@{->}^{\ol\sfp_\CZ} "U";"Y"};
{\ar@{->}^{\wt\sfq^-_\CZ} "Z";"X"};
{\ar@{->}_{\wt\sfp^-_\CZ} "Z";"Y"};
{\ar@{_{(}->}_j "W";"Z"};  
{\ar@{->}^{\pi_\CZ} "V";"U"};
{\ar@{->}^{\on{transl}_{\rho_P(\omega_X)}} "X'";"X"};
\endxy
$$

\sssec{}

By base change for any $\CF\in \Dmod_{\frac{1}{2}}(\BunNbZ)$, we obtain a morphism
\begin{multline} \label{e:CT co of Poinc LHS 12}
(\wt\sfq^-_\CZ)_*\circ \left((\wt\sfp^-_\CZ)^!\circ (\ol\sfp_\CZ)_*(\CF)\sotimes (j_\CZ)_!(\omega_{\Bun_{P^-}\times \CZ})\right)\to \\
\to (\wt\sfq^-_\CZ)_*\circ ({}'\ol\sfp_\CZ)_*\left(({}'\wt\sfp^-_\CZ)^!(\CF)\sotimes ({}'\ol\sfp_\CZ)^!\circ (j_\CZ)_!(\omega_{\Bun_{P^-}\times \CZ})\right),
\end{multline} 
see the maps in the next diagram\footnote{Note that the long vertical arrow is an open embedding.}:

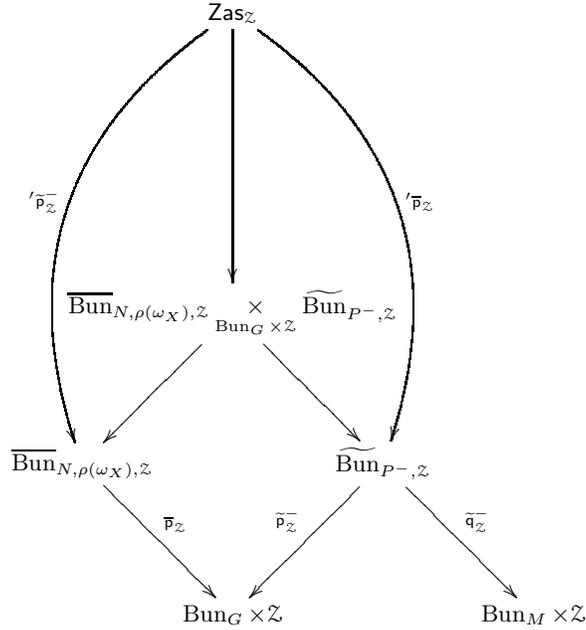
$$
\xy
(20,0)*+{\Bun_M\times \CZ}="X";
(-20,0)*+{\Bun_G\times \CZ}="Y";
(0,20)*+{\BunPmtpZ}="Z";
(-40,20)*+{\BunNbZ}="U";
(-20,40)*+{\BunNbZ\underset{\Bun_G\times \CZ}\times \BunPmtpZ}="V";
(-20,80)*+{\Zas_\CZ}="T";
{\ar@{->}^{\ol\sfp_\CZ} "U";"Y"};
{\ar@{->}^{\wt\sfq^-_\CZ} "Z";"X"};
{\ar@{->}_{\wt\sfp^-_\CZ} "Z";"Y"};
{\ar@{->} "V";"U"};
{\ar@{->} "V";"Z"};
{\ar@{->}"T";"V"};
{\ar@/_3pc/@{->}_{'\wt\sfp^-_\CZ}"T";"U"};
{\ar@/^3pc/@{->}^{'\ol\sfp_\CZ} "T";"Z"};
\endxy
$$

The following is established as in \cite[Proposition 4.1.8]{Lin}: 

\begin{lem} \label{l:open Bruhat}
The map \eqref{e:CT co of Poinc LHS 12} is an isomorphism if $\CF\in \Whit^{\on{glob}}(G)_{\on{co},\CZ}$.
\end{lem} 

\sssec{}

Applying \lemref{l:open Bruhat}, we rewrite the expression in \eqref{e:CT co of Poinc LHS 1} (without the $\rho_P$-translation and the cohomological shift) as
\begin{equation} \label{e:CT co of Poinc LHS 2}
(\wt\sfq^-_\CZ)_*\circ ({}'\ol\sfp_\CZ)_*\left(\left(({}'\wt\sfp^-_\CZ)^!\circ ((\pi_\CZ)_*)^{\on{Av},\chi}(-) \right) 
\sotimes ({}'\ol\sfp_\CZ)^!\circ (j_\CZ)_!(\omega_{\Bun_{P^-}\times \CZ})\right).
\end{equation}

We depict this by the following diagram:
\begin{equation} \label{e:CT co of Poinc diag 1}
\vcenter
{\xy
(20,0)*+{\Bun_M\times \CZ}="X";
(0,20)*+{\BunPmtpZ}="Z";
(-40,20)*+{\BunNbZ}="U";
(-20,40)*+{\Zas_\CZ}="V";
(-60,40)*+{\Gr_{G,\rho(\omega_X),\CZ}}="T";
{\ar@{->}^{\wt\sfq^-_\CZ} "Z";"X"};
{\ar@{->}^{'\ol\sfp_\CZ} "V";"Z"};
{\ar@{->}_{'\wt\sfp^-_\CZ} "V";"U"};
{\ar@{->}_{\pi_\CZ} "T";"U"};
\endxy}
\end{equation} 

\sssec{}

Using \eqref{e:Zast to config}, we rewrite \eqref{e:CT co of Poinc LHS 2} as
\begin{equation} \label{e:CT co of Poinc LHS 3}
(\ol\sfp_{M,\on{gen}}\times \on{id}_\CZ)_*\circ (\fs_\CZ)_*\left(\left(({}'\wt\sfp^-_\CZ)^!\circ ((\pi_\CZ)_*)^{\on{Av},\chi}(-) \right)
\sotimes ({}'\ol\sfp_\CZ)^!\circ (j_\CZ)_!(\omega_{\Bun_{P^-}\times \CZ})\right),
\end{equation}
where $\ol\sfp_{M,\on{gen}}$ denotes the map
$$\BunNbMg\to \Bun_M.$$

I.e., we replace \eqref{e:CT co of Poinc diag 1} by
\begin{equation} \label{e:CT co of Poinc diag 2}
\vcenter
{\xy
(20,0)*+{\Bun_M\times \CZ}="X";
(0,20)*+{\BunNbMg\times \CZ}="Z";
(-40,20)*+{\BunNbZ}="U";
(-20,40)*+{\Zas_\CZ}="V";
(-60,40)*+{\Gr_{G,\rho(\omega_X),\CZ}}="T";
{\ar@{->}^{\ol\sfp_{M,\on{gen}}\times \on{id}_\CZ} "Z";"X"};
{\ar@{->}^{\fs_\CZ} "V";"Z"};
{\ar@{->}_{'\wt\sfp^-_\CZ} "V";"U"};
{\ar@{->}_{\pi_\CZ} "T";"U"};
\endxy}
\end{equation}

\sssec{} \label{sss:CT co of Poinc diag 3}

Note that according to \lemref{l:UHA}, the map
$$(\sff_{M,\CZ})_*\circ (\sff_{M,\CZ})^!\to \on{Id}$$
is an isomorphism, where $\sff_{M,\CZ}$ denotes the map 
$$\sff_{M,\CZ}: \BunNbMZsub\to  \BunNbMg\times \CZ$$
and also its base change
$$\Zas^\Ran_\CZ\to \Zas_\CZ.$$

\medskip

Hence, we can rewrite the expression in \eqref{e:CT co of Poinc LHS 2} as
\begin{multline*}
(\on{id}\times \on{pr}_{\on{small},\CZ})_*\circ (\ol\sfp_{M,\CZ^\subseteq})_*\circ (\fs_{\CZ^{\subseteq}})_*\circ \sff_{M,\CZ}^!
\left(\left(({}'\wt\sfp^-_\CZ)^!\circ ((\pi_\CZ)_*)^{\on{Av},\chi}(-)\right) \sotimes ({}'\ol\sfp_\CZ)^!\circ (j_\CZ)_!(\omega_{\Bun_{P^-}\times \CZ})\right)\simeq \\
\simeq (\on{id}\times \on{pr}_{\on{small},\CZ})_!\circ (\ol\sfp_{M,\CZ^\subseteq})_*\circ (\fs_{\CZ^{\subseteq}})_*\circ \sff_{M,\CZ}^!
\left(\left(({}'\wt\sfp^-_\CZ)^!\circ ((\pi_\CZ)_*)^{\on{Av},\chi}(-)\right) \sotimes ({}'\ol\sfp_\CZ)^!\circ (j_\CZ)_!(\omega_{\Bun_{P^-}\times \CZ})\right),
\end{multline*}
where the isomorphism is due to the fact that $\on{pr}_{\on{small},\CZ}$ is pseudo-proper, and 
where the maps are as in the diagram below:

\begin{equation} \label{e:CT co of Poinc diag 3}
\vcenter
{\xy
(40,20)*+{\Bun_M\times \CZ^\subseteq}="X";
(60,0)*+{\Bun_M\times \CZ}="X'";
(20,40)*+{\BunNbMZsub}="Z";
(-40,20)*+{\BunNbZ}="U";
(-20,40)*+{\Zas_\CZ}="V";
(0,60)*+{\Zas^\Ran_\CZ}="S";
(-60,40)*+{\Gr_{G,\rho(\omega_X),\CZ}}="T";
(20,0)*+{\BunPmtpZ}="R";
(-10,0)*+{\Bun_{P^-}\times \CZ}="P";
{\ar@{->}^{\ol\sfp_{M,\CZ^\subseteq}} "Z";"X"};
{\ar@{->}^{\fs_{\CZ^\subseteq}} "S";"Z"};
{\ar@{->}_{'\wt\sfp^-_\CZ} "V";"U"};
{\ar@{->}_{\pi_\CZ} "T";"U"};
{\ar@{->}_{\sff_{M,\CZ}} "S";"V"};
{\ar@{^{(}->}_j "P";"R"};  
{\ar@{->}^{'\ol\sfp_\CZ} "V";"R"};  
{\ar@{->}^{\wt\sfq^-_\CZ} "R";"X'"};  
{\ar@{->}^{\on{id}\times \on{pr}_{\on{small},\CZ}} "X";"X'"};  
\endxy}
\end{equation} 

\sssec{}

Note now that since $\rho_P$ is central in $M$, we have a well-defined map 
$$\on{transl}_{\rho_P(\omega_X)}:\BunNbMZsubM\to \BunNbMZsub.$$

The next assertion is the main point in the proof of \thmref{t:CT co of Poinc}:

\begin{prop} \label{p:CT co of Poinc}
There is a canonical isomorphism
\begin{multline} \label{e:CT co of Poinc main}
(\on{transl}_{\rho_P(\omega_X)})^*\circ (\fs_{\CZ^{\subseteq}})_*\circ \sff_{M,\CZ}^!
\left(\left(({}'\wt\sfp^-_\CZ)^!\circ ((\pi_\CZ)_*)^{\on{Av},\chi}(-)\right) \sotimes ({}'\ol\sfp_\CZ)^!\circ (j_\CZ)_!(\omega_{\Bun_{P^-}\times \CZ})\right) \simeq \\
\simeq ((\pi_{M,\CZ^\subseteq})_*)^{\on{Av},\chi} \circ J^{-,!}_{\Whit,\Theta} \circ \on{ins.unit}_\CZ(-)[\on{shift}-2\delta_{(N_P)_{\rho_P(\omega_X)}}-\delta_{(N^-_P)_{\rho_P(\omega_X)}}],
\end{multline}
as functors
$$\Whit_*(G)_\CZ\rightrightarrows \Dmod_{\frac{1}{2}}(\BunNbMZsubM)_{\on{co}},$$
where $\pi_{M,\CZ^\subseteq}$ denotes the map
$$\Gr_{M,\rho_M(\omega_X),\CZ^\subseteq}\to \BunNbMZsubM.$$
\end{prop} 

\sssec{} \label{sss:proof of CT co of Poinc}

\propref{p:CT co of Poinc} combined with \secref{sss:CT co of Poinc diag 3}
implies that the expression in \eqref{e:CT co of Poinc LHS 1} identifies with
\begin{multline*}
(\on{id}\times \on{pr}_{\on{small},\CZ})_!\circ 
(\ol\sfp_{M,\CZ^\subseteq})_*\circ ((\pi_{M,\CZ^\subseteq})_*)^{\on{Av},\chi} \circ J^{-,!}_{\Whit,\Theta} \circ \on{ins.unit}_\CZ[-2\delta_{(N_P)_{\rho_P(\omega_X)}}-\delta_{(N^-_P)_{\rho_P(\omega_X)}}]\simeq \\
\simeq (\on{id}\times \on{pr}_{\on{small},\CZ})_!\circ 
\on{Poinc}_{M,*,\CZ^\subseteq}\circ J^{-,!}_{\Whit,\Theta} \circ \on{ins.unit}_\CZ[-2\delta_{(N_P)_{\rho_P(\omega_X)}}-\delta_{(N^-_P)_{\rho_P(\omega_X)}}],
\end{multline*}
establishing \thmref{t:CT co of Poinc}.

\qed[\thmref{t:CT co of Poinc}]

\begin{rem}

Note that one can regard \propref{p:CT co of Poinc} as giving an expression for $J^{-,!}_{\Whit}$ in ``finite-dimensional" terms, i.e., 
in terms of sheaves with finite-dimensional support. Namely, it says that in terms of the equivalences
$$\Whit^!(G)_\CZ\overset{\Theta_{\Whit(G)}}\simeq \Whit_*(G)_\CZ \overset{((\pi_\CZ)_*)^{\on{Av},\chi}}\simeq \Whit^{\on{glob}}(G)_{\on{co},\CZ}$$
and
$$\Whit^!(M)_{\CZ^\subseteq}\overset{\Theta_{\Whit(M)}}\simeq \Whit_*(M)_{\CZ^\subseteq} \overset{((\pi_{M,\CZ^\subseteq})_*)^{\on{Av},\chi}}\simeq \Whit^{\on{glob}}(M)_{\on{co},\CZ^\subseteq},$$
the functor
$$\Whit^!(G)_\CZ \overset{\on{ins.unit}_\CZ}\longrightarrow \Whit^!(G)_{\CZ^\subseteq} \overset{J^{-,!}_{\Whit}}\longrightarrow  \Whit^!(M)_{\CZ^\subseteq}$$
corresponds to 
$$(\on{transl}_{\rho_P(\omega_X)})^*\circ (\fs_{\CZ^{\subseteq}})_*\circ \sff_{M,\CZ}^!
\left(\left(({}'\wt\sfp^-_\CZ)^!(-)\right) \sotimes ({}'\ol\sfp_\CZ)^!\circ (j_\CZ)_!(\omega_{\Bun_{P^-}\times \CZ})\right).$$

\end{rem} 

\ssec{Proof of \propref{p:CT co of Poinc}}

\sssec{}

First, we notice that translation by $\rho_P(\omega_X)$ allows us to think about
$\Whit^!(M)_\CZ$ as a full subcategory in $\Dmod_{\frac{1}{2}}(\Gr_{M,\rho(\omega_X)})$ (rather than $\Dmod_{\frac{1}{2}}(\Gr_{M,\rho_M(\omega_X)})$.

\medskip

A similar remark applies also to $\Whit_*(M)_\CZ$, $\Whit^{\on{glob}}(M)_\CZ$ and $\Whit^{\on{glob}}(M)_{\on{co},\CZ}$.

\medskip

Thus, our goal is to establish an isomorphism
\begin{multline} \label{e:CT co of Poinc main 1}
(\fs_{\CZ^{\subseteq}})_*\circ \sff_{M,\CZ}^!
\left(\left(({}'\wt\sfp^-_\CZ)^!\circ ((\pi_\CZ)_*)^{\on{Av},\chi}(-)\right)
\sotimes ({}'\ol\sfp_\CZ)^!\circ (j_\CZ)_!(\omega_{\Bun_{P^-}\times \CZ})\right) \simeq \\
\simeq 
((\pi_{M,\CZ^\subseteq})_*)^{\on{Av},\chi} \circ J^{-,!}_{\Whit,\Theta,\rho_P} \circ \on{ins.unit}_\CZ(-)[\on{shift}-2\delta_{(N_P)_{\rho_P(\omega_X)}}-\delta_{(N^-_P)_{\rho_P(\omega_X)}}],
\end{multline}
as functors
$$\Whit_*(G)_\CZ\rightrightarrows \Dmod_{\frac{1}{2}}(\BunNbMZsub)_{\on{co}},$$
where:

\begin{itemize}

\item $\pi_{M,\CZ^\subseteq}$ denotes now the map $\Gr_{M,\rho(\omega_X),\CZ^\subseteq}\to \BunNbMZsub$; 

\medskip

\item $J^{-,!}_{\Whit,\Theta,\rho_P}:=(\on{transl}_{\rho_P(\omega_X)})_*\circ J^{-,!}_{\Whit,\Theta,}$;

\end{itemize} 

\sssec{}

First, it is straightforward to check\footnote{One inspects the actions of the corresponding groupoids, cf. \secref{sss:Jacqet to Whit}.}
that the left-hand side in \eqref{e:CT co of Poinc main 1} belongs to
$$\Whit^{\on{glob}}(M)_{\on{co},\CZ^\subseteq}\subset  \Dmod_{\frac{1}{2}}(\BunNbMZsub)_{\on{co}}.$$

Hence, by \thmref{t:local-to-global Whit}(b,c), the isomorphism \eqref{e:CT co of Poinc main 1} is equivalent to
\begin{multline} \label{e:CT co of Poinc main 2}
(\pi_{M,\CZ^\subseteq})^!\circ \Theta^{\on{glob}}_{\Whit(M)}\circ 
(\fs_{\CZ^{\subseteq}})_*\circ \sff_{M,\CZ}^!
\left(\left(({}'\wt\sfp^-_\CZ)^!\circ ((\pi_\CZ)_*)^{\on{Av},\chi}(-)\right) \sotimes ({}'\ol\sfp_\CZ)^!\circ (j_\CZ)_!(\omega_{\Bun_{P^-}\times \CZ})\right) \simeq \\
\simeq  J^{-,!}_{\Whit,\rho_P} \circ  \on{ins.unit}_\CZ\circ \Theta_{\Whit(G)}(-)[\on{shift}-2\delta_{N_{\rho(\omega_X)}}-\delta_{(N^-_P)_{\rho_P(\omega_X)}}]
\end{multline}
as functors 
$$\Whit_*(G)_\CZ\rightrightarrows \Whit_*(M)_{\CZ^\subseteq}\subset \Dmod_{\frac{1}{2}}(\Gr_{M,\rho(\omega_X),\CZ^{\subseteq}}).$$

Note that the left-hand side in \eqref{e:CT co of Poinc main 2} identifies with 
\begin{equation} \label{e:CT co of Poinc main 2.5}
(\pi_{M,\CZ^\subseteq})^!\circ  
(\fs_{\CZ^{\subseteq}})_*\circ \sff_{M,\CZ}^!
\left(\left(({}'\wt\sfp^-_\CZ)^!\circ \Theta^{\on{glob}}_{\Whit(G)} \circ ((\pi_\CZ)_*)^{\on{Av},\chi}(-)\right) 
\sotimes ({}'\ol\sfp_\CZ)^!\circ (j_\CZ)_!(\omega_{\Bun_{P^-}\times \CZ})\right).
\end{equation}

\sssec{}

Note that by \lemref{l:basic Zastava} we have a diagram in which the square is Cartesian:

\medskip

\begin{equation} \label{e:CT co of Poinc diag 4}
\vcenter
{\xy
(20,80)*+{\Gr^{\on{pos}}_{G,\rho(\omega_X),\CZ^\subseteq}\underset{\Gr_{G,\rho(\omega_X),\CZ^\subseteq}}\times 
(\Gr_{G,\rho(\omega_X)}\times \Gr_{M,\rho(\omega_X)})^{\on{neg}}_{\CZ^\subseteq}}="A";
(40,60)*+{\Gr_{M,\rho(\omega_X),\CZ^\subseteq}}="X";
(20,40)*+{\BunNbMZsub}="Z";
(-40,20)*+{\BunNbZ}="U";
(-20,40)*+{\Zas_\CZ}="V";
(0,60)*+{\Zas^\Ran_\CZ}="S";
(-60,40)*+{\Gr_{G,\rho(\omega_X),\CZ}}="T";
(20,0)*+{\BunPmtpZ}="R";
(-10,0)*+{\Bun_{P^-}\times \CZ}="P";
{\ar@{->}^{\fs_{\CZ^\subseteq}} "S";"Z"};
{\ar@{->}_{'\wt\sfp^-_\CZ} "V";"U"};
{\ar@{->}_{\pi_\CZ} "T";"U"};
{\ar@{->}_{\sff_{M,\CZ}} "S";"V"};
{\ar@{^{(}->}_j "P";"R"};  
{\ar@{->}^{'\ol\sfp_\CZ} "V";"R"};  
{\ar@{->}^{\pi_{M,\CZ^\subseteq}} "X";"Z"};  
{\ar@{->}^{'\fs_{\CZ^\subseteq}} "A";"X"};
{\ar@{->}_{'\pi_{M,\CZ^\subseteq}} "A";"S"};
\endxy}
\end{equation} 

Hence, by base change, the expression in \eqref{e:CT co of Poinc main 2.5} identifies with $({}'\fs_{\CZ^\subseteq})_*$ applied to 
\begin{equation} \label{e:CT co of Poinc LHS 4}
({}'\wt\sfp^-_\CZ\circ \sff_{M,\CZ}\circ {}'\pi_{M,\CZ^\subseteq})^!\left(\Theta^{\on{glob}}_{\Whit(G)}\circ ((\pi_\CZ)_*)^{\on{Av},\chi}(-)\right)\sotimes 
({}'\ol\sfp_\CZ\circ \sff_{M,\CZ}\circ {}'\pi_{M,\CZ^\subseteq})^!\left((j_\CZ)_!(\omega_{\Bun_{P^-}\times \CZ})\right).
\end{equation} 

\sssec{}

We claim:

\begin{prop} \label{p:CT co of Poinc next} \hfill

\smallskip

\noindent{\em(a)} The first factor in \eqref{e:CT co of Poinc LHS 4}, viewed as an object of
$$\Dmod_{\det^{\frac{1}{2}}_G}(\Gr_{G,\rho(\omega_X)}\times \Gr_{M,\rho(\omega_X)})_{\CZ^\subseteq},$$
identifies with
$$p_1^!\circ \on{ins.unit}_\CZ \circ \pi_\CZ^!\circ \Theta^{\on{glob}}_{\Whit(G)}\circ ((\pi_\CZ)_*)^{\on{Av},\chi}(-),$$
where 
$$p_1:(\Gr_{G,\rho(\omega_X)}\times \Gr_{M,\rho(\omega_X)})_{\CZ^\subseteq}\to \Gr_{G,\rho(\omega_X)}.$$

\smallskip

\noindent{\em(b)} The second factor in \eqref{e:CT co of Poinc LHS 4}, viewed as an object of
$$\Dmod_{\det^{\frac{1}{2}}_{G,M}}(\Gr_{G,\rho(\omega_X)}\times \Gr_{M,\rho(\omega_X)})_{\CZ^\subseteq},$$
identifies with the object\footnote{See \secref{sss:Delta tilde} for the notation.}
$\wt\Delta^{-,\semiinf}_{\CZ^\subseteq}$.
\end{prop}

\sssec{} 

Let us show how \propref{p:CT co of Poinc next} implies \propref{p:CT co of Poinc}. Indeed, combining 
\propref{p:CT co of Poinc next}(a) and \thmref{t:local-to-global Whit}(c), we obtain that the first factor in \eqref{e:CT co of Poinc LHS 4}
identifies with
$$p_1^!\circ \on{ins.unit}_\CZ(-)[-2\delta_{N_{\rho_P(\omega_X)}}].$$

Now, the assertion of \propref{p:CT co of Poinc} follows from the definition of the functor $J^{-,!}_{\Whit}$; 
there is a built-in cohomological shift\footnote{This shift is designed to be compatible with $[\on{shift}]$ in the definition of $\on{CT}^-_*$.}
and the additional shift by $\delta_{(N^-_P)_{\rho_P(\omega_X)}}$ comes from the translation
by $\rho_P(\omega_X)$. 

\qed[\propref{p:CT co of Poinc next}]

\sssec{Proof of \propref{p:CT co of Poinc next}(a)}

By \lemref{l:basic Zastava}, the composition
$$\Gr^{\on{pos}}_{G,\rho(\omega_X),\CZ^\subseteq}\underset{\Gr_{G,\rho(\omega_X),\CZ^\subseteq}}\times 
(\Gr_{G,\rho(\omega_X)}\times \Gr_{M,\rho(\omega_X)})^{\on{neg}}_{\CZ^\subseteq} 
\overset{'\pi_{M,\CZ^\subseteq}}\longrightarrow \Zas^\Ran_\CZ\overset{\sff_{M,\CZ}}\longrightarrow \Zas_\CZ
\overset{'\wt\sfp^-_\CZ}\longrightarrow \BunNbZ$$
is the same as 
\begin{multline*} 
\Gr^{\on{pos}}_{G,\rho(\omega_X),\CZ^\subseteq}\underset{\Gr_{G,\rho(\omega_X),\CZ^\subseteq}}\times 
(\Gr_{G,\rho(\omega_X)}\times \Gr_{M,\rho(\omega_X)})^{\on{neg}}_{\CZ^\subseteq} \overset{p_1}\to \\
\to \Gr^{\on{pos}}_{G,\rho(\omega_X),\CZ^\subseteq} \overset{\pi^{\on{pos}}_{\CZ^\subseteq}}\longrightarrow 
\BunNbZ\underset{\CZ}\times \CZ^\subseteq \overset{\on{id}\times \on{pr}_{\on{small},\CZ}}\longrightarrow \BunNbZ.
\end{multline*} 

Now, the assertion of \propref{p:CT co of Poinc next}(a) follows from \lemref{l:unital local and global Whit}. 

\sssec{Proof of \propref{p:CT co of Poinc next}(b)}

By \lemref{l:basic Zastava}, the composition
$$\Gr^{\on{pos}}_{G,\rho(\omega_X),\CZ^\subseteq}\underset{\Gr_{G,\rho(\omega_X),\CZ^\subseteq}}\times 
(\Gr_{G,\rho(\omega_X)}\times \Gr_{M,\rho(\omega_X)})^{\on{neg}}_{\CZ^\subseteq} 
\overset{'\pi_{M,\CZ^\subseteq}}\longrightarrow \Zas^\Ran_\CZ\overset{\sff_{M,\CZ}}\longrightarrow \Zas_\CZ
\overset{'\ol\sfp_\CZ}\longrightarrow \BunPmtpZ$$
is the same as 
\begin{multline*}
\Gr^{\on{pos}}_{G,\rho(\omega_X),\CZ^\subseteq}\underset{\Gr_{G,\rho(\omega_X),\CZ^\subseteq}}\times 
(\Gr_{G,\rho(\omega_X)}\times \Gr_{M,\rho(\omega_X)})^{\on{neg}}_{\CZ^\subseteq}\to \\
\to (\Gr_{G,\rho(\omega_X)}\times \Gr_{M,\rho(\omega_X)})^{\on{neg}}_{\CZ^\subseteq} \overset{\pi^{\on{neg}}_{\CZ^\subseteq}}\longrightarrow
\BunPmtpZ\underset{\CZ}\times \CZ^\subseteq\overset{\on{id}\times \on{pr}_{\on{small},\CZ}}\longrightarrow \BunPmtpZ.$$
\end{multline*} 

The assertion follows now from \lemref{l:pullback Delta}. 

\ssec{Whittaker coefficients of Eisenstein series: the enhanced version}

\sssec{}

The next result is an enhanced version of \thmref{t:CT co of Poinc}

\begin{thm} \label{t:CT co of Poinc enh}
For $\CZ\to \Ran$, the following diagram commutes
\begin{equation}  \label{e:CT co of Poinc * enh}
\CD
(\Whit_*(G)\underset{\Sph_G}\otimes \on{I}(G,P^-)^{\on{loc}}_{\rho_P(\omega_X)})_\CZ 
@>{\on{ins.unit}_\CZ}>> (\Whit_*(G)\underset{\Sph_G}\otimes \on{I}(G,P^-)^{\on{loc}}_{\rho_P(\omega_X)})_{\CZ^\subseteq}  \\
@VV{\Poinc_{G,*,\CZ}^{-,\on{enh}}[d]}V
@VV{J^{-,\on{enh}}_{\Whit,\Theta}}V  \\
\Dmod_{\frac{1}{2}}(\Bun_G)^{-,\on{enh}}_{\on{co},\rho_P(\omega_X),\CZ} & &  \Whit_*(M)_{\CZ^{\subseteq}}  \\
@VV{\on{CT}^{-,\on{enh}}_{\on{co},?,\rho_P(\omega_X),\CZ}\otimes \on{Id}}V @VV{\on{Poinc}_{M,*,\CZ^{\subseteq}}}V \\
\Dmod_{\frac{1}{2}}(\Bun_M)_{\on{co}}\otimes \Dmod(\CZ) @<<{\on{Id}\otimes (\on{pr}_{\on{small},\CZ})_!}< 
\Dmod_{\frac{1}{2}}(\Bun_M)_{\on{co}}\otimes \Dmod(\CZ^{\subseteq}),
\endCD
\end{equation}
where:

\begin{itemize}

\item $\Poinc_{G,*,\CZ}^{-,\on{enh}}$ is the functor 
\begin{multline*} 
(\Whit_*(G)\underset{\Sph_G}\otimes \on{I}(G,P^-)^{\on{loc}}_{\rho_P(\omega_X)})_\CZ
\overset{\on{Poinc}_{G,*,\CZ}\otimes \on{Id}_{\on{I}(G,P^-)_{\rho_P(\omega_X)}}}\longrightarrow \\
\to \left(\Dmod_{\frac{1}{2}}(\Bun_G)_{\on{co}}\otimes \Dmod(\CZ)\right)\underset{\Sph_{G,\CZ}}\otimes
\on{I}(G,P^-)^{\on{loc}}_{\rho_P(\omega_X),\CZ}=: 
\Dmod_{\frac{1}{2}}(\Bun_G)^{-,\on{enh}}_{\on{co},\rho_P(\omega_X),\CZ} ;
\end{multline*}

\medskip

\item $J^{-,\on{enh}}_{\Whit,\Theta}:=\Theta^{-1}_{\Whit(M)}\circ J^{-,\on{enh}}_{\Whit}\circ (\Theta_{\Whit(G)}\otimes \on{Id}_{\on{I}(G,P^-)^{\on{loc}}_{\rho_P(\omega_X)}})$;

\medskip

\item $d=2\delta_{(N_P)_{\rho_P(\omega_X)}}+\delta_{(N^-_P)_{\rho_P(\omega_X)}}$. 

\end{itemize} 

Moreover, the isomorphism of functors in \eqref{e:CT co of Poinc * enh} respects the actions of $\Sph_{M,\CZ}$. 

\end{thm}

\sssec{} \label{sss:coeff Z untl}

Note that unital property of the local-to-global functor $\ul{\on{Poinc}}_{G,!}$ implies that the functor
$$\Dmod_{\frac{1}{2}}(\Bun_G)\otimes \Dmod(\CZ) \overset{\on{Id}\otimes (\on{pr}_{\on{small},\CZ})^!}\longrightarrow
\Dmod_{\frac{1}{2}}(\Bun_G)\otimes \Dmod(\CZ^\subseteq) \overset{\on{coeff}_{G,\CZ^\subseteq}}\longrightarrow \Whit^!(G)_{\CZ^\subseteq}$$
takes values in 
$$\Whit^!(G)_{\CZ^{\subseteq,\on{untl}},\on{indep}} \subset \Whit^!(G)_{\CZ^\subseteq}.$$

We will denote the corresponding functor
$$\Dmod_{\frac{1}{2}}(\Bun_G)\otimes \Dmod(\CZ)\to \Whit^!(G)_{\CZ^{\subseteq,\on{untl}},\on{indep}},$$
as well as its composition with the embedding
$$\Whit^!(G)_{\CZ^{\subseteq,\on{untl}},\on{indep}}\hookrightarrow \Whit^!(G)_{\CZ^{\subseteq,\on{untl}}}$$
by $\on{coeff}_{G,\CZ,\on{untl}}$.

\sssec{}

Passing to dual functors as in \corref{c:coeff of Eis !}, from \thmref{t:CT co of Poinc enh} we obtain: 

\begin{cor} \label{c:coeff of Eis ! enh}
For $\CZ\to \Ran$, the following diagram commutes:
$$
\CD
\Whit^!(G)_{\CZ^{\subseteq,\on{untl}},\on{indep}} @<{\on{emb.indep}_{\Whit^!(G)}^L}<<  \Whit^!(G)_{\CZ^{\subseteq,\on{untl}}} \\
@A{\on{coeff}_{G,\CZ,\on{untl}}[d]}AA    @AA{\on{co}\!J^{-,\on{enh}}_{\Whit,\Theta}}A \\
\Dmod_{\frac{1}{2}}(\Bun_G)\otimes \Dmod(\CZ) & & (\Whit^!(M)\underset{\Sph_M}\otimes \on{I}(G,P^-)^{\on{loc}}_{\rho_P(\omega_X)})_{\CZ^{\subseteq,\on{untl}}}  \\
& & @AAA \\
@A{\Eis^{-,\on{enh}}_{!,\rho_P(\omega_X),\CZ}}AA  
\Whit^!(M)_{\CZ^{\subseteq,\on{untl}}}\underset{\Sph_{M,\CZ}}\otimes  \on{I}(G,P^-)^{\on{loc}}_{\rho_P(\omega_X),\CZ^{\subseteq,\on{untl}}} \\
& & @AA{\on{coeff}_{M,\CZ,\on{untl}}\otimes \on{ins.unit}_\CZ}A  \\
\Dmod_{\frac{1}{2}}(\Bun_M)^{-,\on{enh}}_{\rho_P(\omega_X),\CZ} @>{=}>>
(\Dmod_{\frac{1}{2}}(\Bun_M)\otimes \Dmod(\CZ))\underset{\Sph_{M,\CZ}}\otimes \on{I}(G,P^-)^{\on{loc}}_{\rho_P(\omega_X),\CZ},
\endCD
$$
where:
\begin{itemize}

\item $\on{co}\!J^{-,\on{enh}}_{\Whit,\Theta}$ is the functor obtained from $J^{-,\on{enh}}_{\Whit,\Theta}$ by duality; 

\medskip

\item $d=2\delta_{(N_P)_{\rho_P(\omega_X)}}+\delta_{(N^-_P)_{\rho_P(\omega_X)}}$. 

\end{itemize} 

Moreover, the isomorphism of functors in the above diagram respects the actions of $\Sph_{G,\CZ}$. 

\end{cor}

\ssec{Proof of \thmref{t:CT co of Poinc enh}}

\sssec{}

We will show that the two functors
$$\Whit_*(G)_\CZ\otimes \on{I}(G,P^-)^{\on{loc}}_{\rho_P(\omega_X),\CZ} \rightrightarrows \Dmod_{\frac{1}{2}}(\Bun_M)_{\on{co}}\otimes \Dmod(\CZ)$$
are canonically isomorphic. The compatibility with the actions of $\Sph_{G,\CZ}$ and $\Sph_{M,\CZ}$ will follow from the construction. 

\sssec{}

By \corref{c:loc vs global CT enh ? co}, for $\CF_W\in \Whit_*(G)_\CZ$ and $\CF_I\in \on{I}(G,P^-)^{\on{loc}}_{\rho_P(\omega_X),\CZ}$, the object 
\begin{equation} \label{e:CT co of Poinc enh LHS init}
(\on{CT}^{-,\on{enh}}_{\on{co},?,\rho_P(\omega_X),\CZ}\otimes \on{Id})\circ 
(\on{Poinc}_{G,*,\CZ}\otimes \on{Id}_{\on{I}(G,P^-)_{\rho_P(\omega_X)}})(\CF_W\otimes \CF_I)
\end{equation}
is 
\begin{equation} \label{e:CT co of Poinc enh LHS 1}
(\on{transl}_{\rho_P(\omega_X)})^*\circ 
(\wt\sfq^-_\CZ)_*\circ \left(\left((\wt\sfp^-_\CZ)^!\circ (\ol\sfp_\CZ)_*\circ ((\pi_\CZ)_*)^{\on{Av},\chi}(\CF_W)\right)
\sotimes (\pi_!\circ \wt\sfs^!)_\CZ(\CF_{I,-\rho_P})\right)[-\on{shift}],
\end{equation} 
where:

\begin{itemize}

\item The functor 
$$(\pi_!\circ \wt\sfs^!)_\CZ:\on{I}(G,P^-)^{\on{loc}}_{\rho_P(\omega_X),\CZ} \to \Dmod_{\det_{G,M}^{\frac{1}{2}}}(\BunPmtpZ)$$
is as in \propref{p:loc-to-glob Drinf defnd}(a). 

\item $\CF_{I,-\rho_P}$ denotes the image of $\CF_I$ under 
$\on{I}(G,P^-)^{\on{loc}}_{\rho_P(\omega_X),\CZ}\overset{\alpha^{-1}_{\on{taut},\rho_P(\omega_X)}}\simeq 
\on{I}(G,P^-)^{\on{loc}}_\CZ$.

\end{itemize}

\sssec{}

An analog of \lemref{l:open Bruhat} is applicable in the current situation, and we rewrite \eqref{e:CT co of Poinc enh LHS 1} as

\medskip

\begin{equation} \label{e:CT co of Poinc enh LHS 2}
(\on{transl}_{\rho_P(\omega_X)})^*\circ 
(\wt\sfq^-_\CZ)_*\circ ({}'\ol\sfp_\CZ)_*\left(\left(({}'\wt\sfp^-_\CZ)^!\circ ((\pi_\CZ)_*)^{\on{Av},\chi}(\CF_W)\right) 
\sotimes ({}'\ol\sfp_\CZ)^!\circ (\pi_!\circ \wt\sfs^!)_\CZ(\CF_{I,-\rho_P})\right)[-\on{shift}].
\end{equation}

\medskip

Further, as in the proof of \thmref{t:CT co of Poinc} given in \secref{sss:proof of CT co of Poinc}, it suffices to construct an isomorphism
\begin{multline} \label{e:CT co of Poinc enh main}
(\on{transl}_{\rho_P(\omega_X)})^*\circ (\fs_{\CZ^{\subseteq}})_*\circ \sff_{M,\CZ}^!
\left(\left(({}'\wt\sfp^-_\CZ)^!\circ ((\pi_\CZ)_*)^{\on{Av},\chi}(\CF_W)\right) \sotimes ({}'\ol\sfp_\CZ)^!\circ (\pi_!\circ \wt\sfs^!)_\CZ(\CF_{I,-\rho_P})\right) \simeq \\
\simeq ((\pi_{M,\CZ^\subseteq})_*)^{\on{Av},\chi} \circ \on{ins.unit}_\CZ\left(J^{-,\on{enh}}_{\Whit,\Theta}(\CF_W\otimes \CF_I)\right)
[\on{shift}-2\delta_{(N_P)_{\rho_P(\omega_X)}}-\delta_{(N^-_P)_{\rho_P(\omega_X)}}].
\end{multline}

\sssec{}

Finally, as in the proof of \propref{p:CT co of Poinc}, in order to prove \eqref{e:CT co of Poinc enh main}, it suffices to establish the following extension 
of \lemref{l:pullback Delta}:

\begin{prop} \label{p:pullback Delta gen}
For $\CF_I\in \on{I}(G,P^-)^{\on{loc}}_\CZ$, the object 
\begin{multline*} 
(\pi^{\on{neg}}_{\CZ^\subseteq})^!\circ (\on{id}\times \on{pr}_{\on{small},\CZ})^!\left((\pi_!\circ \wt\sfs^!)_\CZ(\CF_I)\right)\in \\
\in \Dmod_{(\pi^{\on{neg}})^*(\det_{G,M}^{\frac{1}{2}})}((\Gr_G\times \Gr_M)^{\on{neg}}_{\CZ^\subseteq})\subset
\Dmod_{\det_{\Gr_G}^{\frac{1}{2}}\otimes \det_{\Gr_M}^{\frac{1}{2}}}(\Gr_G\times \Gr_M)_{\CZ^{\subseteq}}
\end{multline*}
identifies canonically with the image of 
$$\on{ins.unit}_\CZ(\CF_I)\in \on{I}(G,P^-)^{\on{loc}}_{\CZ^\subseteq}$$
under
\begin{multline} \label{e:P-av}
\on{I}(G,P^-)^{\on{loc}}\to \Dmod_{\frac{1}{2}}(\Gr_G)^{\fL(N^-_P)\cdot \fL^+(M)} \simeq \\
\simeq \Dmod_{\det_{\Gr_G}^{\frac{1}{2}}\otimes \det_{\Gr_M}^{\frac{1}{2}}}(\Gr_G\times \Gr_M)^{\fL(P^-)}\to 
\Dmod_{\det_{\Gr_G}^{\frac{1}{2}}\otimes \det_{\Gr_M}^{\frac{1}{2}}}(\Gr_G\times \Gr_M).
\end{multline}
\end{prop}

\qed[\thmref{t:CT co of Poinc enh}]

\begin{rem}
As the proof of \propref{p:pullback Delta gen} will show, the isomorphism of functors 
$$ \on{I}(G,P^-)^{\on{loc}}_\CZ\rightrightarrows 
\Dmod_{\det_{\Gr_G}^{\frac{1}{2}}\otimes \det_{\Gr_M}^{\frac{1}{2}}}(\Gr_G\times \Gr_M)_{\CZ^{\subseteq}}$$
stated in \propref{p:pullback Delta gen} is compatible with the actions of $(\Sph_G\otimes \Sph_M)_\CZ$
on the two sides.
\end{rem}

\ssec{Proof of \propref{p:pullback Delta gen}}

\sssec{}

For $\CF_I\in \on{I}(G,P^-)^{\on{loc}}_\CZ$, let $\on{Av}_!^{\fL(M)/\fL^+(M)}(\CF_I)$ denote\footnote{We use the symbol $\on{Av}_!^{\fL(M)/\fL^+(M)}$ because
the functor $\Dmod(\Gr_G)^{\fL(N^-_P)\cdot \fL^+(M)} \to \Dmod(\Gr_G\times \Gr_M)^{\fL(P^-)}$ is given embedding along
$\Gr_G\overset{\on{id}\times 1_{\Gr_M}}\longrightarrow \Gr_G\times \Gr_M$, followed by !-averaging along $\fL(M)/\fL^+(M)$.} 
its image along the functor \eqref{e:P-av}. 

\medskip

We first construct a map in one direction:
\begin{equation} \label{e:pullback Delta gen}
\on{Av}_!^{\fL(M)/\fL^+(M)}(\on{ins.unit}_\CZ(\CF_I))\to (\pi^{\on{neg}}_{\CZ^\subseteq})^!\circ (\on{id}\times \on{pr}_{\on{small},\CZ})^!\left((\pi_!\circ \wt\sfs^!)_\CZ(\CF_I)\right).
\end{equation} 

\sssec{}

Let $\wt\pi_\CZ$ denote the map
$$(\Gr_G\times \Gr_M)_\CZ\to \BunPmtpZ.$$

Note also that we have a naturally defined map
$$\sg_\CZ:(\Gr_G\times \Gr_M)_\CZ\to \Gr_{G,\Bun_M,\CZ}$$
so that
$$\wt\pi_\CZ=\pi_\CZ\circ \sg_\CZ.$$

By construction
$$\on{Av}_!^{\fL(M)/\fL^+(M)}(\CF_I)\simeq \sg_\CZ^! \circ \wt\sfs^!_\CZ(\CF_I).$$

From here we obtain a map
$$(\wt\pi_\CZ)_!\circ \on{Av}_!^{\fL(M)/\fL^+(M)}(\CF_I)\to (\pi_!\circ \wt\sfs^!)_\CZ(\CF_I).$$

In particular, we obtain a map
\begin{multline} \label{e:pullback Delta gen 1}
(\pi^{\on{neg}}_{\CZ^\subseteq})_!\circ \on{Av}_!^{\fL(M)/\fL^+(M)}(\on{ins.unit}_\CZ(\CF_I))= \\
=(\wt\pi_{\CZ^\subseteq})_!\circ \on{Av}_!^{\fL(M)/\fL^+(M)}(\on{ins.unit}_\CZ(\CF_I))\to
(\pi_!\circ \wt\sfs^!)_{\CZ^\subseteq}(\on{ins.unit}_\CZ(\CF_I)).
\end{multline} 

\sssec{}

Using the fact that the functor $\ul{(\pi_!\circ \wt\sfs^!)}$ is strictly unital, we obtain that
the right-hand side in \eqref{e:pullback Delta gen 1} identifies with 
$$(\on{id}\times \on{pr}_{\on{small},\CZ})^!\circ (\pi_!\circ \wt\sfs^!)_\CZ(\CF_I).$$

\medskip

To summrize we obtain a map
\begin{equation} \label{e:pullback Delta gen 2}
(\pi^{\on{neg}}_{\CZ^\subseteq})_!\circ \on{Av}_!^{\fL(M)/\fL^+(M)}(\on{ins.unit}_\CZ(\CF_I))\to 
(\on{id}\times \on{pr}_{\on{small},\CZ})^!\circ (\pi_!\circ \wt\sfs^!)_\CZ(\CF_I)
\end{equation}

By adjunction, from \eqref{e:pullback Delta gen 2}, we obtain the sought-for map \eqref{e:pullback Delta gen}. 

\begin{rem}
One can show that the map 
$$(\on{id}\times \on{pr}_{\on{small},\CZ})_!\circ (\pi^{\on{neg}}_{\CZ^\subseteq})_!\circ \on{Av}_!^{\fL(M)/\fL^+(M)}(\on{ins.unit}_\CZ(\CF_I))\to 
(\pi_!\circ \wt\sfs^!)_\CZ(\CF_I),$$
obtained by adjunction from \eqref{e:pullback Delta gen 2}, is an isomorphism. 
\end{rem}

\sssec{}

We now show that the map \eqref{e:pullback Delta gen} is an isomorphism. 

\medskip

The categories 
$$\on{I}(G,P^-)^{\on{loc}}_\CZ \text{ and }\Dmod_{\det_{\Gr_G}^{\frac{1}{2}}\otimes \det_{\Gr_M}^{\frac{1}{2}}}(\Gr_G\times \Gr_M)_{\CZ^\subseteq}$$
are acted on by $\Sph_{M,\CZ}$, and the two functors in the lemma, as well as the natural transformation between them constructed above,
respect these actions. 

\medskip

Hence, in order to show that \eqref{e:pullback Delta gen} is an isomorphism, it suffices to do so for $\CF_I=\Delta^{-,\semiinf}_\CZ$, since this object
generates $\on{I}(G,P^-)^{\on{loc}}_\CZ$ under the $\Sph_{M,\CZ}$-action. 

\medskip

However, in the latter case, the map \eqref{e:pullback Delta gen} equals the map \eqref{e:pullback Delta}, and the assertion follows from \lemref{l:pullback Delta}. 

\qed[\propref{p:pullback Delta gen}]

\section{Constant term of Kac-Moody localization} \label{s:Eis and Loc}

The goal of this subsection is to express the composition of constant term and localization functors
in terms of the Jacquet functor for Kac-Moody modules. 

\medskip

Once combined with its spectral counterpart, this will be used in Part III of the paper in order to establish the 
compatibility of the Langlands functor with constant terms.  

\ssec{Linear twistings on \texorpdfstring{$\Bun_G$}{BunG}}

\sssec{}

Consider the torus $G_{\on{ab}}:=G/[G,G]$. Since $\Bun_{G_{\on{ab}}}$ is an (abelian) group-stack, it make
sense to talk about (commutative) multiplicative de Rham twistings on it. Such twistings are in bijection with
(abelian) central extensions of the Lie algebra of $\Bun_{G_{\on{ab}}}$, i.e., $\Gamma(X,\CO_X)\otimes \fg_{\on{ab}}[1]$.

\medskip

In particular, points in 
$$\left(\Gamma(X,\CO_X)\otimes \fg_{\on{ab}}\right)^*\simeq \Gamma(X,\omega_X\otimes \fg_{\on{ab}}^*)[1],$$
i.e., torsors on $X$ with respect to $\fg_{\on{ab}}^*\otimes \omega_X$, 
give rise to commutative multiplicative de Rham twistings on $\Bun_{G_{\on{ab}}}$,
(cf. \secref{sss:twist by ZcG}).

\sssec{} \label{sss:twisted Loc sign}

Given such a torsor $\CP$, we consider the de Rham twisting on $\Bun_G$ obtained
by pullback from the corresponding commutative multiplicative de Rham twisting on
$\Bun_{G_{\on{ab}}}$; we will denote it by the same symbol $\CP$. 

\medskip

For a point $\ul{x}\in \Ran$, the resulting twisting on $\Bun^{\on{level}_{\ul{x}}}_G$ 
is compatible with the \emph{inverse} of the multiplicative twisting on $\fL(G)_{\ul{x}}$ 
corresponding to the central extension of $\hg$ from  \secref{sss:twist by ZcG}. 

\medskip

Thus, for $\CP$ as above and a level $\kappa$, we can consider the corresponding category 
$\Dmod_{\kappa+\CP}(\Bun_G)$ and the localization functor
$$\ul\Loc_{G,\kappa+\CP}:\ul\KL(G)_{\kappa+\CP}\to 
\Dmod_{\kappa-\CP}(\Bun_G)\otimes \ul\Dmod(\Ran).$$

\sssec{}  \label{sss:L line bundles}

Let now $\CP_{Z^0_\cG}$ be a $Z^0_\cG$-torsor. Weil pairing with $\CP_{Z^0_\cG}$ gives rise to a line bundle on
$\Bun_{G_{\on{ab}}}$, to be denoted $\CL_{\CP_{Z^0_\cG}}$; by a slight abuse of notation, we will denote by
the same symbol $\CL_{\CP_{Z^0_\cG}}$ its pullback along
$$\Bun_G\to \Bun_{G_{\on{ab}}}.$$

\sssec{}

Note that as in \secref{sss:create g ab tors} we can regard $\on{dlog}(\CP_{Z^0_\cG})$ as a 
torsor on $X$ with respect to $\fg_{\on{ab}}^*\otimes \omega_X$.

\medskip

It is easy to see, however, that 
$$\on{dlog}(\CP_{Z^0_\cG})\simeq \on{dlog}(\CL_{\CP_{Z^0_\cG}})$$
as de Rham twistings on $\Bun_G$. 

\sssec{}

Thus, for a level $\kappa$ we have a naturally defined local-to-global functor
$$\ul\Loc_{G,\kappa+\on{dlog}(\CP_{Z^0_\cG})}:\ul\KL(G)_{\kappa+\on{dlog}(\CP_{Z^0_\cG})}\to 
\Dmod_{\kappa-\on{dlog}(\CL_{\CP_{Z^0_\cG}})}(\Bun_G)\otimes \ul\Dmod(\Ran).$$

\medskip

For $\kappa=\crit$ we will denote the resulting functor
\begin{multline*} 
\ul\KL(G)_{\crit+\on{dlog}(\CP_{Z^0_\cG})}\overset{\ul\Loc_{G,\crit+\on{dlog}(\CP_{Z^0_\cG})}}\longrightarrow 
\Dmod_{\crit-\on{dlog}(\CL_{\CP_{Z^0_\cG})}}(\Bun_G)\otimes \ul\Dmod(\Ran)\overset{(-)\otimes \CL_{\CP_{Z^0_\cG}}}\longrightarrow \\
\to \Dmod_\crit(\Bun_G)\otimes \ul\Dmod(\Ran)\simeq \Dmod_{\frac{1}{2}}(\Bun_G)\otimes \ul\Dmod(\Ran)
\end{multline*}
by $\ul\Loc_{G,\CP_{Z^0_\cG}}$.

\sssec{} \label{sss:global torsors lambda}

When $\CP_{Z^0_\cG}=\clambda(\omega_X)$ for a coweight $\clambda$ of $Z^0_\cG$ (i.e., a character of $G$),
we will use a shorthand notation
$$\Dmod_{\kappa+\clambda}(\Bun_G):=\Dmod_{\kappa+\on{dlog}(\clambda(\omega_X))}(\Bun_G)\otimes \ul\Dmod(\Ran),$$
and also 
$$\ul\Loc_{G,\kappa+\clambda}:\ul\KL(G)_{\kappa+\clambda}\to \Dmod_{\kappa-\clambda}(\Bun_G)\otimes \ul\Dmod(\Ran)$$
and 
$$\ul\Loc_{G,\clambda}:\ul\KL(G)_{\crit+\clambda}\to \Dmod_{\frac{1}{2}}(\Bun_G)\otimes \ul\Dmod(\Ran)$$
for the corresponding functors. 

\medskip

Note that as \secref{sss:lambdach twist}, the above notions are well-defined for any element $\lambda\in z_\cG:=\on{Lie}(Z^0_\cG)$. 

%
%
%
%

\sssec{}

Contents of this subsection will mostly be applied when the reductive group in question is $M$, a Levi subgroup
of the original $G$. 

\ssec{Integration along \texorpdfstring{$\Bun_{P^-}\to \Bun_M$}{int} and BRST}

\sssec{}

For a level $\kappa$ for $M$, let us denote by $\sfq^{-,\on{shift}}_*$ the functor 
\begin{multline*}
\Dmod_{\crit_G+\kappa}(\Bun_{P^-}) =  \Dmod_{\frac{1}{2}\on{dlog}(\det(\Bun_G)|_{\Bun_{P^-}})+\kappa}(\Bun_{P^-}) \simeq \\
\simeq \Dmod_{\frac{1}{2}\on{dlog}(\det(\Bun_M)|_{\Bun_{P^-}})+\on{dlog}(\det^{\otimes \frac{1}{2}}_{\Bun_{G,M}})+\kappa}(\Bun_{P^-})
\overset{\otimes \det^{\otimes -\frac{1}{2}}_{\Bun_{G,M}}}\longrightarrow \\
\to \Dmod_{\frac{1}{2}\on{dlog}(\det(\Bun_M)|_{\Bun_{P^-}})+\kappa}(\Bun_{P^-})\overset{(\sfq^-)_*}\longrightarrow 
\Dmod_{\frac{1}{2}\on{dlog}(\det(\Bun_M))+\kappa}(\Bun_M)\simeq \\
\simeq \Dmod_{\crit_M+\kappa}(\Bun_M)\overset{[-\on{shift}]}\longrightarrow \Dmod_{\crit_M+\kappa}(\Bun_M),
\end{multline*} 
where $[\on{shift}]$ is as in \eqref{e:global shift}. 

\sssec{}

Let $\fl_{N^-_P}$ be the (ungraded) line from \cite[Sect. 12.5.7]{GLC2}. I.e.,
$$\fl_{N^-_P}=\det(\fn^-_P\otimes \Gamma(X,\CO_X))\simeq \det(\Gamma(X,\CO_X))^{\otimes \dim(\fn^-_P)}\otimes \det(\fn^-_P)^{\otimes (1-g)}.$$

\medskip

Let $\CL_{\rhoch_P(\omega_X)}$ be the line bundle on $\Bun_M$ from \secref{sss:L line bundles}. 

\medskip

We will prove:

\begin{thm}  \label{t:Loc kappa}
For $\CZ\to \Ran$, the following diagram of functors canonically commutes
\begin{equation} \label{e:Loc kappa 0}
\CD
\Dmod_{\crit_G+\kappa}(\Bun_{P^-})\otimes \Dmod(\CZ) @>{\sfq^{-,\on{shift}}_*\otimes \on{Id}}>> \Dmod_{\crit_M+\kappa}(\Bun_M)\otimes  \Dmod(\CZ)  \\
@A{\Loc_{P^-,\crit_G+\kappa,\CZ}}AA  @AA{\fl_{N^-_P}\otimes ((-)\otimes \CL_{\rhoch_P(\omega_X)})\otimes (\on{pr}_{\on{small},\CZ})_!}A \\ 
\KL(P^-)_{\crit_G+\kappa,\CZ} & & \Dmod_{\crit_M+\kappa-\rhoch_P}(\Bun_M)\otimes  \Dmod(\CZ^\subseteq) \\
@V{\on{ins.vac.}_\CZ}VV  @AA{\Loc_{M,\crit_M+\kappa+\rhoch_P,\CZ^\subseteq}}A \\
\KL(P^-)_{\crit_G+\kappa,\CZ^\subseteq} @>{\BRST^-}>> \KL(M)_{\crit_M+\kappa+\rhoch_P,\CZ^\subseteq}.
\endCD
\end{equation} 
\end{thm} 

\sssec{}

Let $\on{CT}^-_*$ be as in \secref{sss:CT level}. We claim:

\begin{cor} \label{c:CT and Loc kappa}
For $\CZ\to \Ran$, the following diagram of functors canonically commutes
\begin{equation} \label{e:CT and Loc kappa}
\CD
\Dmod_{\crit_G+\kappa}(\Bun_G)\otimes \Dmod(\CZ) @>{\on{CT}^-_*\otimes \on{Id}}>> \Dmod_{\crit_M+\kappa}(\Bun_M)\otimes  \Dmod(\CZ)  \\
@A{\Loc_{G,\crit_G+\kappa,\CZ}}AA  @AA{\fl_{N^-_P}\otimes ((-)\otimes \CL_{\rhoch_P(\omega_X)})\otimes (\on{pr}_{\on{small},\CZ})_!}A \\ 
\KL(G)_{\crit_G+\kappa,\CZ} & & \Dmod_{\crit_M+\kappa-\rhoch_P}(\Bun_M)\otimes  \Dmod(\CZ) \\
@V{\on{ins.vac.}_\CZ}VV  @AA{\Loc_{M,\crit_M+\kappa+\rhoch_P,\CZ^\subseteq}}A \\
\KL(G)_{\crit_G+\kappa,\CZ^\subseteq} @>{J_{\on{KM}}^{-,\Sph}}>> \KL(M)_{\crit_M+\kappa+\rhoch_P,\CZ^\subseteq}.
\endCD
\end{equation} 
\end{cor} 

\begin{proof}

Follows from \thmref{t:Loc kappa} in the same way as 
\cite[Theorem 12.8.5]{GLC2} follows from the combination of \cite[Corollary 13.3.11 and Theorem 12.5.13]{GLC2}. 

\end{proof}

\sssec{}

In particular, specializing to the critical level, we obtain:
\begin{cor} \label{c:CT and Loc crit}
For $\CZ\to \Ran$, the following diagram of functors canonically commutes
\begin{equation} \label{e:CT and Loc crit}
\CD
\Dmod_{\frac{1}{2}}(\Bun_G)\otimes \Dmod(\CZ) @>{\on{CT}^-_*\otimes \on{Id}}>> \Dmod_{\frac{1}{2}}(\Bun_M)\otimes  \Dmod(\CZ)  \\
@A{\Loc_{G,\CZ}}AA  @AA{\fl_{N^-_P}\otimes (\on{Id}\otimes (\on{pr}_{\on{small},\CZ})_!)}A \\ 
\KL(G)_{\crit_G,\CZ} & & \Dmod_{\frac{1}{2}}(\Bun_M)\otimes  \Dmod(\CZ^\subseteq) \\
@V{\on{ins.vac.}_\CZ}VV  @AA{\Loc_{M,\rhoch_P,\CZ^\subseteq}}A \\
\KL(G)_{\crit_G,\CZ^\subseteq} @>{J_{\on{KM}}^{-,\Sph}}>> \KL(M)_{\crit_M+\rhoch_P,\CZ^\subseteq}.
\endCD
\end{equation}
\end{cor} 

\ssec{Identification of line bundles}

As a first step towards the proof of \thmref{t:Loc kappa}, we will explain the appearance of the line $\fl_{N^-_P}$. 

\sssec{}

Consider the relative cotangent bundle $T^*(\Bun_{P^-}/\Bun_M)$. It is easy to see that its determinant (viewed as an ungraded line bundle)
canonically descends to a line bundle on $\Bun_M$, to be denoted 
$$\det(T^*(\Bun_{P^-}/\Bun_M)).$$

\sssec{} \label{sss:square root expl}

Recall now (see \cite[Proposition 1.3.3]{GLC1}) that the line bundle $\det^{\otimes \frac{1}{2}}_{\Bun_{G,M}}$ on $\Bun_M$ attaches to $\CP_M\in \Bun_M$
the line
$$\det(\Gamma(X,(\fn^-_P)_{\CP_M}\otimes \omega_X^{\otimes \frac{1}{2}}))\otimes 
\det(\Gamma(X,\fn^-_P\otimes \omega_X^{\otimes \frac{1}{2}}))^{\otimes -1}.$$

Note that this is the same  
$$\det(\Gamma(X,(\fn^-_P)_{\CP_{P^-}}\otimes \omega_X^{\otimes \frac{1}{2}}))\otimes 
\det(\Gamma(X,\fn^-_P\otimes \omega_X^{\otimes \frac{1}{2}}))^{\otimes -1},$$
where $\CP_{P^-}$ is any $P^-$-bundle that projects to $\CP_M$. 

\sssec{}

We claim:

\begin{prop} \label{p:rel det}
There is a canonical isomorphism of line bundles on $\Bun_M$:
$$\det(T^*(\Bun_{P^-}/\Bun_M))\simeq \det^{\otimes \frac{1}{2}}_{\Bun_{G,M}}\otimes \CL_{\rhoch_P(\omega_X)}\otimes 
\fl_{N^-_P},$$
where
$$\rhoch_P(\omega_X):=2\rhoch_P(\omega^{\otimes \frac{1}{2}}_X).$$
\end{prop}

\begin{proof}

The tangent space to $\Bun_M$ at $\CP_M$ is
$$\Gamma(X,(\fn^-_P)_{\CP_M})[1].$$

Since 
$$\det((V[1])^*)\simeq \det(V),$$
the line bundle $\det(T^*(\Bun_{P^-}/\Bun_M))$ attaches to $\CP_M\in \Bun_M$ the line 
$$\det(\Gamma(X,(\fn^-_P)_{\CP_M})).$$

Thus, we need to establish an isomorphism
\begin{multline} \label{e:det identity}
\det(\Gamma(X,(\fn^-_P)_{\CP_M}\otimes \omega_X^{\otimes \frac{1}{2}}))\otimes 
\det(\Gamma(X,(\fn^-_P)_{\CP_M}))^{\otimes -1}\simeq \\
\simeq \det(\Gamma(X,\fn^-_P\otimes \omega_X^{\otimes \frac{1}{2}}))\otimes 
\det(\Gamma(X,\fn^-_P\otimes \CO_X))^{\otimes -1}\otimes \on{Weil}(-2\rhoch_P(\CP_M),\omega^{\otimes \frac{1}{2}}).
\end{multline}

Recall that for a vector bundle $\CE$ of rank $n$ and a line bundle $\CL$ on $X$ we have
\begin{multline} \label{e:det formula}
\det(\Gamma(X,\CE\otimes \CL))\simeq
\det(\Gamma(X,\CE))\otimes \det(\Gamma(X,\CL))^{\otimes n}\otimes \det(\Gamma(X,\CO_X))^{\otimes -n}\otimes
\on{Weil}(\det(\CE),\CL).
\end{multline} 

Hence, \eqref{e:det identity} is equivalent to
$$\on{Weil}(\on{det}((\fn^-_P)_{\CP_M}),\omega^{\otimes \frac{1}{2}})\simeq
\on{Weil}(\on{det}(\fn^-_P\otimes \CO_X),\omega^{\otimes \frac{1}{2}})\otimes \on{Weil}(-2\rhoch_P(\CP_M),\omega^{\otimes \frac{1}{2}}),$$
where $\on{Weil}(\on{det}(\fn^-_P\otimes \CO_X),\omega^{\otimes \frac{1}{2}})$ is a 
constant line.\footnote{Note that for a line $\fl$ and a line bundle $\CL$, we have $\on{Weil}(\fl\otimes \CO_X,\CL)\simeq \fl^{\deg(\CL)}$.} 

\end{proof}

\sssec{}

Thanks to \propref{p:rel det}, the operation of tensoring by $\det(T^*(\Bun_{P^-}/\Bun_M))$ gives rise to a functor
$$\Dmod_{\crit_M+\kappa}(\Bun_{P^-})\to \Dmod_{\crit_G+\kappa+\rhoch_P}(\Bun_{P^-}).$$

\medskip

Hence,  the operation of tensoring by $\det(T^*(\Bun_{P^-}))$ defines a functor
$$\Dmod_{\crit_M+\kappa}(\Bun_{P^-})\to \Dmod_{\crit_G+2\crit_M+\kappa+\rhoch_P}(\Bun_{P^-}).$$

\sssec{}

We record the following particular case of \cite[Proposition 10.5.7]{GLC2}: 

\begin{prop} \label{p:Loc and duality P}
With respect to the identification
$$\left(\KL(P^-)_{\crit_M-\kappa-\rhoch_P}\right)^\vee \simeq \KL(P^-)_{\crit_G+\kappa}$$
induced by \eqref{e:duality KL P}, and
$$(\Dmod_{\crit_M-\kappa+\rhoch_P}(\Bun_{P^-})_{\on{co}})^\vee\simeq
\Dmod_{-\crit_M+\kappa-\rhoch_P}(\Bun_{P^-}),$$
the dual of the functor\footnote{See \secref{sss:twisted Loc sign} for the sign flip.}
$$\Gamma_{P^-,\crit_M-\kappa-\rhoch_P,\CZ}:\Dmod_{\crit_M-\kappa+\rhoch_P}(\Bun_{P^-})_{\on{co}}\otimes \Dmod(\CZ)\to 
\KL(P^-)_{\crit_M-\kappa-\rhoch_P,\CZ}$$
identifies with the functor 
\begin{multline*} 
\KL(P^-)_{\crit_G+\kappa,\CZ}\overset{\Loc_{P^-,\crit_G+\kappa,\CZ}}\longrightarrow 
\Dmod_{\crit_G+\kappa}(\Bun_{P^-})\otimes \Dmod(\CZ)
\overset{((-)\otimes \det(T^*(\Bun_{P^-}))^{\otimes -1})\otimes \on{Id}}\longrightarrow \\
\to \Dmod_{-\crit_M+\kappa-\rhoch_P}(\Bun_{P^-})\otimes \Dmod(\CZ)
\overset{[-\dim(\Bun_{P^-})]}\longrightarrow \Dmod_{-\crit_M+\kappa-\rhoch_P}(\Bun_{P^-})\otimes \Dmod(\CZ),
\end{multline*}
where $\dim(\Bun_{P^-})$ refers to the dimension of a \emph{given connected component} of $\Bun_{P^-}$. 
\end{prop} 

\begin{rem}

In the statement of \propref{p:Loc and duality P}, it is not claimed that the category $$\Dmod_{\crit_M-\kappa+\rhoch_P}(\Bun_{P^-})_{\on{co}}$$
is dualizable, so the symbol $(\Dmod_{\crit_M-\kappa+\rhoch_P}(\Bun_{P^-})_{\on{co}})^\vee$ is understood as 
$$\on{Funct}_{\on{cont}}(\Dmod_{\crit_M-\kappa+\rhoch_P}(\Bun_{P^-})_{\on{co}},\Vect),$$
and similarly for the notion of the dual functor.

\end{rem} 

\ssec{Proof of \thmref{t:Loc kappa}}


\sssec{}

First, as in \cite[Sect. 12.5.10]{GLC2}, we observe that there exists a natural transformation 
\begin{equation} \label{e:Loc kappa 1}
\Res^{\hm}_{\hp^-}\circ \Gamma_{M,\crit_M-\kappa-\rhoch_P,\CZ}\to
\Gamma_{P^-,\crit_M-\kappa-\rhoch_P,\CZ}\circ ((\sfq^-)^!\otimes \on{Id})
\end{equation}
as functors
$$\Dmod_{\crit_M-\kappa+\rhoch_P}(\Bun_M)_{\on{co}}\otimes \Dmod(\CZ) \rightrightarrows \KL(P^-)_{\crit_M-\kappa-\rhoch_P,\CZ}.$$

\sssec{}

Passing to dual functors and using \propref{p:Loc and duality P} and \cite[Proposition 10.5.7]{GLC2}, from \eqref{e:Loc kappa 1}
we obtain a natural transformation

\medskip

\begin{equation} \label{e:Loc kappa P}
\vcenter
{\xy 
(0,0)*+{\Dmod_{\crit_G+\kappa}(\Bun_{P^-})\otimes \Dmod(\CZ)}="A";
(70,0)*+{\Dmod_{\crit_M+\kappa}(\Bun_M)\otimes  \Dmod(\CZ)}="B";
(0,-40)*+{\KL(P^-)_{\crit_G+\kappa,\CZ}}="C"; 
(70,-40)*+{\KL(M)_{\crit_M+\kappa+\rhoch_P,\CZ}.}="D";
{\ar@{->}^{\sfq^{-,\on{shift}}_*\otimes \on{Id}} "A";"B"};
{\ar@{->}^{\BRST^-} "C";"D"};
{\ar@{->}^{\Loc_{P^-,\crit_G+\kappa,\CZ}} "C";"A"};
{\ar@{->}_{\fl_{N^-_P}\otimes ((-)\otimes \CL_{\rhoch_P(\omega_X)})\circ \Loc_{M,\crit_M+\kappa+\rhoch_P,\CZ}} "D";"B"};
{\ar@{=>} "D";"A"};
\endxy}
\end{equation} 

Thus, we obtain a natural transformation
\begin{equation} \label{e:Loc kappa 2}
((-)\otimes \CL_{\rhoch_P(\omega_X)}\otimes \fl_{N^-_P})\circ \ul\Loc_{M,\crit_M+\kappa+\rhoch_P}\circ 
\ul{\BRST}^-\to \sfp^{-,\on{shift}}_*\circ \ul\Loc_{P^-,\crit_G+\kappa}
\end{equation} 
as local-to-global functors
$$\ul{\KL}(P^-)_{\crit_G+\kappa}\rightrightarrows \Dmod_{\crit_M+\kappa}(\Bun_M)\otimes \ul\Dmod(\Ran).$$

\sssec{}

Both sides in \eqref{e:Loc kappa 2} have natural lax unital structures and the map \eqref{e:Loc kappa 2}
is compatible with these structures. Furthermore, the lax unital structure on the right-hand side is strict.

\medskip

Hence, by \cite[Sect. 11.4.8]{GLC2}, the natural transformation \eqref{e:Loc kappa 2} extends to a natural transformation
\begin{equation} \label{e:Loc kappa 3}
((-)\otimes \CL_{\rhoch_P(\omega_X)}\otimes \fl_{N^-_P})\circ \ul\Loc_{M,\crit_M+\kappa+\rhoch_P}\circ 
(\ul{\BRST}^-)^{\int\, \on{ins.vac}}\to \sfp^{-,\on{shift}}_*\circ \ul\Loc_{P^-,\crit_G+\kappa}, 
\end{equation} 
see \cite[Sect. 11.4]{GLC2} for the notation $\int\, \on{ins.vac}$. 

\medskip

In other words, we obtain a natural transformation from the counterclockwise circuit to the 
clockwise circuit in \eqref{e:Loc kappa 0}. 

\sssec{}

We now claim that \eqref{e:Loc kappa 3} is an isomorphism. Indeed, this reduces to the statement of
\cite[Theorem 12.5.13]{GLC2} by the same mechanism as \cite[Theorem 12.8.15]{GLC2}.

\qed[\thmref{t:Loc kappa}] 

\ssec{Constant term of localization: the enhanced version}

In this section we specialize to the critical level. 

\sssec{}

Let $\Loc_{G,\CZ}^{-,\on{enh}_{\on{co}}}$ denote the functor
\begin{multline*} 
\KL(G)^{-,\on{enh}_{\on{co}}}_{\crit,\CZ}:=\KL(G)_{\crit,\CZ}\underset{\Sph_{G,\CZ}}\otimes \on{I}(G,P^-)^{\on{loc}}_{\on{co},\CZ}
\overset{\Loc_{G,Z}\otimes \on{Id}}\longrightarrow  \\
\to \left(\Dmod_{\frac{1}{2}}(\Bun_G)\otimes \Dmod(\CZ)\right)\underset{\Sph_{G,\CZ}}\otimes 
\on{I}(G,P^-)^{\on{loc}}_{\on{co},\CZ}=:
\Dmod_{\frac{1}{2}}(\Bun_G)^{-,\on{enh}_{\on{co}}}_\CZ.
\end{multline*} 

\sssec{}

We will prove the following enhanced version of \corref{c:CT and Loc crit}:

\begin{thm} \label{t:CT and Loc crit enh}
For $\CZ\to \Ran$, the following diagram of functors 
canonically commutes, in a way compatible with the $\Sph_{M,\CZ}$-actions: 
$$
\CD
\Dmod_{\frac{1}{2}}(\Bun_G)^{-,\on{enh}_{\on{co}}}_\CZ @>{\on{CT}^{-,\on{enh}_{\on{co}}}_{*,\CZ}}>> \Dmod_{\frac{1}{2}}(\Bun_M)\otimes  \Dmod(\CZ)  \\
@A{\Loc^{-,\on{enh}_{\on{co}}}_{G,\CZ}}AA  @AA{\fl_{N^-_P}\otimes (\on{Id}\otimes (\on{pr}_{\on{small},\CZ})_!)}A \\ 
\KL(G)^{-,\on{enh}_{\on{co}}}_{\crit_G,\CZ} & & \Dmod_{\frac{1}{2}}(\Bun_M)\otimes  \Dmod(\CZ^\subseteq) \\
@V{\on{ins.vac.}_\CZ}VV  @AA{\Loc_{M,\rhoch_P,\CZ^\subseteq}}A \\
\KL(G)^{-,\on{enh}_{\on{co}}}_{\crit_G,\CZ^\subseteq} @>{J_{\on{KM}}^{-,\on{enh}_{\on{co}}}}>> \KL(M)_{\crit_M+\rhoch_P,\CZ^\subseteq},
\endCD
$$
where $J_{\on{KM}}^{-,\on{enh}_{\on{co}}}$ is the functor from \secref{sss:J enh}.
\end{thm} 

\sssec{Strategy of proof} \label{sss:CT and Loc crit enh strat}

We will construct a natural transformation
\begin{equation} \label{e:CT and Loc crit enh}
\xy 
(0,0)*+{\Dmod_{\frac{1}{2}}(\Bun_G)^{-,\on{enh}_{\on{co}}}_\CZ}="A";
(70,0)*+{\Dmod_{\frac{1}{2}}(\Bun_M)\otimes  \Dmod(\CZ)}="B";
(0,-40)*+{\KL(G)^{-,\on{enh}_{\on{co}}}_{\crit_G,\CZ}}="C"; 
(70,-40)*+{\KL(M)_{\crit_M+\rhoch_P,\CZ},}="D";
{\ar@{->}^{\on{CT}^{-,\on{enh}_{\on{co}}}_{*,\CZ}} "A";"B"};
{\ar@{->}^{J_{\on{KM}}^{-,\on{enh}_{\on{co}}}} "C";"D"};
{\ar@{->}^{\Loc^{-,\on{enh}_{\on{co}}}_{G,\CZ}} "C";"A"};
{\ar@{->}_{\fl_{N^-_P}\otimes \Loc_{M,\rhoch_P,\CZ}} "D";"B"};
{\ar@{=>} "D";"A"};
\endxy
\end{equation} 
compatible with the $\Sph_{M,\CZ}$-actions, 
such that its concatenation with the commutative diagram
$$
\CD
\KL(G)^{-,\on{enh}_{\on{co}}}_{\crit_G,\CZ} @>{J_{\on{KM}}^{-,\on{enh}_{\on{co}}}}>> \KL(M)_{\crit_M+\rhoch_P,\CZ} \\
@A{\on{Id}\otimes \one_{\on{I}(G,P^-)^{\on{loc}}_{\on{co},\CZ}}}AA @AA{\on{Id}}A \\
\KL(G)_{\crit_G,\CZ}  @>{J_{\on{KM}}^{-,\Sph}}>> \KL(M)_{\crit_M+\rhoch_P,\CZ}, 
\endCD
$$
which is a natural transformation 
\begin{equation} \label{e:CT and Loc crit enh 1}
\xy
(0,0)*+{\Dmod_{\frac{1}{2}}(\Bun_G)^{-,\on{enh}_{\on{co}}}_\CZ}="A";
(70,0)*+{\Dmod_{\frac{1}{2}}(\Bun_M)\otimes  \Dmod(\CZ)}="B";
(0,-40)*+{\KL(G)_{\crit_G,\CZ}}="C"; 
(70,-40)*+{\KL(M)_{\crit_M+\rhoch_P,\CZ},}="D";
{\ar@{->}^{\Loc^{-,\on{enh}_{\on{co}}}_{G,\CZ}\circ (\on{id}\otimes \one_{\on{I}(G,P^-)^{\on{loc}}_{\on{co},\CZ}})} "C";"A"};
{\ar@{->}_{\fl_{N^-_P}\otimes \Loc_{M,\rhoch_P,\CZ}} "D";"B"};
{\ar@{->}^{\on{CT}^{-,\on{enh}_{\on{co}}}_{*,\CZ}} "A";"B"};
{\ar@{->}^{J_{\on{KM}}^{-,\Sph}} "C";"D"};
{\ar@{=>} "D";"A"};
\endxy
\end{equation} 
is obtained by concatenating the commutative diagram
\begin{equation} \label{e:CT and Loc crit enh 1.5}
\CD
\Dmod_{\frac{1}{2}}(\Bun_G)^{-,\on{enh}_{\on{co}}}_\CZ @>{\on{CT}^{-,\on{enh}_{\on{co}}}_{*,\CZ}}>>  \Dmod_{\frac{1}{2}}(\Bun_M)\otimes  \Dmod(\CZ) \\
@A{\on{Id}\otimes \one_{\on{I}(G,P^-)^{\on{loc}}_{\on{co},\CZ}}}AA @AA{\on{Id}}A \\
\Dmod_{\frac{1}{2}}(\Bun_G)\otimes  \Dmod(\CZ) @>{\on{CT}^-_{*,\CZ}}>>  \Dmod_{\frac{1}{2}}(\Bun_M)\otimes  \Dmod(\CZ) 
\endCD
\end{equation} 
with the natural transformation
\begin{equation} \label{e:CT and Loc crit enh 2}
\xy
(0,0)*+{\Dmod_{\frac{1}{2}}(\Bun_G)\otimes  \Dmod(\CZ)}="A";
(70,0)*+{\Dmod_{\frac{1}{2}}(\Bun_M)\otimes  \Dmod(\CZ)}="B";
(0,-40)*+{\KL(G)_{\crit_G,\CZ}}="C"; 
(70,-40)*+{\KL(M)_{\crit_M+\rhoch_P,\CZ},}="D";
{\ar@{->}^{\Loc_{G,\CZ}} "C";"A"};
{\ar@{->}_{\fl_{N^-_P}\otimes \Loc_{M,\rhoch_P,\CZ}} "D";"B"};
{\ar@{->}^{\on{CT}^-_{*}\otimes \on{Id}_{\Dmod(\CZ)}} "A";"B"};
{\ar@{->}^{J_{\on{KM}}^{-,\Sph}} "C";"D"};
{\ar@{=>} "D";"A"};
\endxy
\end{equation} 
obtained from \eqref{e:CT and Loc crit}.

\sssec{} \label{sss:strat CT and Loc crit enh 1}

Let us show how the existence of a diagram \eqref{e:CT and Loc crit enh} with the above properties 
implies the assertion of \thmref{t:CT and Loc crit enh}. 

\medskip

Namely, let us view the family of diagrams \eqref{e:CT and Loc crit enh} as a natural transformation
\begin{equation} \label{e:CT and Loc crit enh 1 form}
\fl_{N^-_P}\otimes \ul\Loc_{M,\rhoch_P}\circ \ul{J}_{\on{KM}}^{-,\on{enh}_{\on{co}}}\to 
\ul{\on{CT}}^{-,\on{enh}_{\on{co}}}_*\circ \ul{\Loc}^{-,\on{enh}_{\on{co}}}_G
\end{equation} 
as lax-unital local-to-global functors
\begin{equation} \label{e:CT and Loc crit enh 3}
\ul{\KL}(G)^{-,\on{enh}_{\on{co}}}_{\crit_G}\rightrightarrows 
\Dmod_{\frac{1}{2}}(\Bun_M)\otimes  \ul\Dmod(\Ran).
\end{equation} 

The right-hand side in \eqref{e:CT and Loc crit enh 3} is strictly unital, and hence \eqref{e:CT and Loc crit enh 3}
induces a natural transformation
\begin{equation} \label{e:CT and Loc crit enh 2 form}
\fl_{N^-_P}\otimes \ul\Loc_{M,\rhoch_P}\circ (\ul{J}_{\on{KM}}^{-,\on{enh}_{\on{co}}})^{\int\, \on{ins.vac}}\to 
\ul{\on{CT}}^{-,\on{enh}_{\on{co}}}_*\circ \ul{\Loc}^{-,\on{enh}_{\on{co}}}_G.
\end{equation} 

We claim that \eqref{e:CT and Loc crit enh 2 form} is an isomorphism. 

\sssec{} \label{sss:strat CT and Loc crit enh 2}

Indeed, since $\one_{\on{I}(G,P^-)^{\on{loc}}_{\on{co}}}$ generates $\on{I}(G,P^-)^{\on{loc}}_{\on{co}}$ as a $\Sph_M$-module 
category and since the natural transformation \eqref{e:CT and Loc crit enh 2 form} is compatible with the $\ul{\Sph}_M$-actions,
in order to show that \eqref{e:CT and Loc crit enh 2 form} is an isomorphism, its suffices to show that the induced natural
transformation
\begin{multline} \label{e:CT and Loc crit enh 3 form}
\fl_{N^-_P}\otimes \ul\Loc_{M,\rhoch_P}\circ \left(\ul{J}_{\on{KM}}^{-,\on{enh}_{\on{co}}}\circ 
(\on{Id}\otimes \ul\one_{\on{I}(G,P^-)^{\on{loc}}_{\on{co}}})\right)^{\int\, \on{ins.vac}}\to \\
\to \ul{\on{CT}}^{-,\on{enh}_{\on{co}}}_*\circ \ul{\Loc}^{-,\on{enh}_{\on{co}}}_G\circ (\on{Id}\otimes \ul\one_{\on{I}(G,P^-)^{\on{loc}}_{\on{co}}})
\end{multline} 
is an isomorphism as local-to-global functors 
$$\ul{\KL}(G)_{\crit_G}\rightrightarrows 
\Dmod_{\frac{1}{2}}(\Bun_M)\otimes  \ul\Dmod(\Ran)$$
is an isomorphism. 

\medskip

However, by \secref{sss:CT and Loc crit enh strat}, the natural transformation is obtained by concatenating \eqref{e:CT and Loc crit enh 1.5}
with the natural transformation
\begin{equation} \label{e:CT and Loc crit enh 4 form}
\fl_{N^-_P}\otimes \ul\Loc_{M,\rhoch_P}\circ (\ul{J}_{\on{KM}}^{-,\on{Sph}})^{\int\, \on{ins.vac}}\to
(\on{CT}^-_{*}\otimes \on{Id}_{\ul\Dmod(\Ran)})\circ \ul{\Loc}_G.
\end{equation}

Hence, it is enough to show that \eqref{e:CT and Loc crit enh 4 form} is an isomorphism. However, the latter follows by construction 
from the commutativity of \eqref{e:CT and Loc crit}. 

\qed[\thmref{t:CT and Loc crit enh}]

\ssec{Construction of the natural transformation for the enhanced functor}  \label{ss:nat transf}

To simplify the notation, we will assume that $\CZ=\on{pt}$ so that $\CZ\to \Ran$ corresponds to
a point $\ul{x}\in \Ran$. 

\medskip

Thus, our goal will be to construct a natural transformation in
\begin{equation} \label{e:CT and Loc crit enh x}
\xy 
(0,0)*+{\Dmod_{\frac{1}{2}}(\Bun_G)\underset{\Sph_{G,\ul{x}}}\otimes \on{I}(G,P^-)^{\on{loc}}_{\on{co},\ul{x}}}="A";
(70,0)*+{\Dmod_{\frac{1}{2}}(\Bun_M)}="B";
(0,-30)*+{\KL(G)_{\crit,\ul{x}}\underset{\Sph_{G,\ul{x}}}\otimes \on{I}(G,P^-)^{\on{loc}}_{\on{co},\ul{x}}}="C"; 
(70,-30)*+{\KL(M)_{\crit_M+\rhoch_P,\ul{x}},}="D";
{\ar@{->}^{\on{CT}^{-,\on{enh}_{\on{co}}}_{*,\ul{x}}} "A";"B"};
{\ar@{->}^{J_{\on{KM}}^{-,\on{enh}_{\on{co}}}} "C";"D"};
{\ar@{->}^{\Loc_{G,\ul{x}}\otimes \on{Id}} "C";"A"};
{\ar@{->}_{\fl_{N^-_P}\otimes \Loc_{M,\rhoch_P,\ul{x}}} "D";"B"};
{\ar@{=>} "D";"A"};
\endxy
\end{equation}  

\sssec{} \label{sss:nat tarnsf 1}

Let $\kappa$ be an arbitrary level. Parallel to the functor 
$$\Loc_{G,\kappa,\ul{x}}:\KL(G)_{\kappa,\ul{x}}\to \Dmod_{\kappa}(\Bun_G),$$
the construction in \cite[Sect. 10.2 and 10.3]{GLC2} naturally produces a functor 
$$\Loc_{G,\kappa,\ul{x}}:\hg\mod_{\kappa,\ul{x}}\to \Dmod_{\kappa}(\Bun^{\on{level}_{\ul{x}}}_G),$$
compatible with the actions of $\fL(G)_{\ul{x}}$. 

\medskip

In particular, for a subgroup $H\subset \fL^+(G)_{\ul{x}}$, we have a commutative diagram
\begin{equation} \label{e:localization and conv}
\CD
\Dmod_{\kappa}(\Bun_G) \otimes\Dmod_\kappa(\Gr_{G,\ul{x}})^H@>{\star}>> 
\Dmod_{\kappa}(H\backslash \Bun_G^{\on{level}_{\ul{x}}}) \\
@A{\Loc_{G,\kappa,\ul{x}}\otimes \on{Id}}AA @AA{\Loc_{G,\kappa,\ul{x}}}A \\
\KL(G)_{\kappa,\ul{x}} \otimes \Dmod_\kappa(\Gr_{G,\ul{x}})^H @>{\star}>>  \hg\mod^H_{\kappa,\ul{x}}
\endCD
\end{equation} 

Consider the particular case of \eqref{e:localization and conv} when $H=\fL^+(P^-)_{\ul{x}}$. We obtain a
commutative diagram: 
\begin{equation} \label{e:localization and P-}
\CD
\Dmod_{\kappa}(\Bun_G) \otimes \Dmod_\kappa(\Gr_{G,\ul{x}})^{\fL^+(P^-)_{\ul{x}}}
@>{\star}>> 
\Dmod_{\kappa}(\Bun_G\underset{\on{pt}/\fL^+(G)_{\ul{x}}}\times \on{pt}/\fL^+(P^-)_{\ul{x}}) \\
@A{\Loc_{G,\kappa,\ul{x}}\otimes \on{Id}}AA @AA{\Loc_{G,\kappa,\ul{x}}}A \\
\KL(G)_{\kappa,\ul{x}}  \otimes \Dmod_\kappa(\Gr_{G,\ul{x}})^{\fL^+(P^-)_{\ul{x}}} @>{\star}>> 
\hg\mod_{\kappa,\ul{x}}^{\fL^+(P^-)_{\ul{x}}}.
\endCD
\end{equation} 

\sssec{} \label{sss:nat tarnsf 2}

Note now that the map
$$\sfp^-:\Bun_{P^-}\to \Bun_G$$
naturally factors via a map
$$'\sfp^-:\Bun_{P^-}\to \Bun_G\underset{\on{pt}/\fL^+(G)_{\ul{x}}}\times \on{pt}/\fL^+(P^-)_{\ul{x}}.$$

As in \cite[Equation (12.17)]{GLC2}, we obtain a natural transformation
\begin{equation} \label{e:Bun G P level}
\xy 
(0,0)*+{\Dmod_{\kappa}(\Bun_G\underset{\on{pt}/\fL^+(G)_{\ul{x}}}\times \on{pt}/\fL^+(P^-)_{\ul{x}})}="A";
(70,0)*+{\Dmod_{\kappa}(\Bun_{P^-})}="B";
(0,-30)*+{\hg\mod_{\kappa,\ul{x}}^{\fL^+(P^-)_{\ul{x}}}}="C"; 
(70,-30)*+{\KL(P^-)_{\kappa,\ul{x}},}="D";
{\ar@{->}^{({}'\sfp^-)^!} "A";"B"};
{\ar@{->}^{\Res^{\hg}_{\hp^-}} "C";"D"};
{\ar@{->}^{\Loc_{G,\kappa,\ul{x}}} "C";"A"};
{\ar@{->}_{\Loc_{P^-,\kappa,\ul{x}}} "D";"B"};
{\ar@{=>} "D";"A"};
\endxy
\end{equation}  

\sssec{}

Recall also the diagram 
$$
\xy 
(0,0)*+{\Dmod_{\crit_G+\kappa}(\Bun_{P^-})}="A";
(70,0)*+{\Dmod_{\crit_M+\kappa}(\Bun_M)}="B";
(0,-40)*+{\KL(P^-)_{\crit_G+\kappa,\ul{x}}}="C"; 
(70,-40)*+{\KL(M)_{\crit_M+\kappa+\rhoch_P,\ul{x}},}="D";
{\ar@{->}^{\sfq^{-,\on{shift}}_*} "A";"B"};
{\ar@{->}^{\BRST^-} "C";"D"};
{\ar@{->}^{\Loc_{P^-,\crit_G+\kappa,\ul{x}}} "C";"A"};
{\ar@{->}_{\fl_{N^-_P}\otimes ((-)\otimes \CL_{\rhoch_P(\omega_X)})\circ \Loc_{M,\crit_M+\kappa+\rhoch_P,\ul{x}}} "D";"B"};
{\ar@{=>} "D";"A"};
\endxy
$$
see \eqref{e:Loc kappa P}.

\medskip

Concatenating it with diagrams \eqref{e:Bun G P level} and \eqref{e:localization and P-} we obtain a natural transformation
$$
\xy
(0,0)*+{\Dmod_{\crit_G+\kappa}(\Bun_G) \otimes \Dmod_{\crit_G+\kappa}(\Gr_{G,\ul{x}})^{\fL^+(P^-)_{\ul{x}}}}="A";
(70,0)*+{\Dmod_{\crit_M+\kappa}(\Bun_M)}="B";
(0,-40)*+{\KL(G)_{\crit_G+\kappa,\ul{x}}  \otimes \Dmod_{\crit_G+\kappa}(\Gr_{G,\ul{x}})^{\fL^+(P^-)_{\ul{x}}}}="C"; 
(70,-40)*+{\KL(M)_{\crit_M+\kappa+\rhoch_P,\ul{x}}.}="D";
{\ar@{->}"A";"B"};
{\ar@{->} "C";"D"};
{\ar@{->}^{\Loc_{G,\crit_G+\kappa,\ul{x}}\otimes \on{Id}} "C";"A"};
{\ar@{->}_{\fl_{N^-_P}\otimes ((-)\otimes \CL_{\rhoch_P(\omega_X)})\circ \Loc_{M,\crit_M+\kappa+\rhoch_P,\ul{x}}} "D";"B"};
{\ar@{=>} "D";"A"};
\endxy
$$

\sssec{}

We now specialize to the critical level, so that the latter diagram becomes
\begin{equation} \label{e:conv and CT}
\vcenter
{\xy
(0,0)*+{\Dmod_{\frac{1}{2}}(\Bun_G) \otimes \Dmod_{\frac{1}{2}}(\Gr_{G,\ul{x}})^{\fL^+(P^-)_{\ul{x}}}}="A";
(70,0)*+{\Dmod_{\frac{1}{2}}(\Bun_M)}="B";
(0,-40)*+{\KL(G)_{\crit_G,\ul{x}}  \otimes \Dmod_{\frac{1}{2}}(\Gr_{G,\ul{x}})^{\fL^+(P^-)_{\ul{x}}}}="C"; 
(70,-40)*+{\KL(M)_{\crit_M+\rhoch_P,\ul{x}}.}="D";
{\ar@{->}"A";"B"};
{\ar@{->} "C";"D"};
{\ar@{->}^{\Loc_{G,\ul{x}}\otimes \on{Id}} "C";"A"};
{\ar@{->}_{ \fl_{N^-_P}\otimes \Loc_{M,\rhoch_P,\ul{x}}} "D";"B"};
{\ar@{=>} "D";"A"};
\endxy}
\end{equation} 
 
\medskip

Note that the functor
\begin{multline*}
\Dmod_{\frac{1}{2}}(\Bun_G)  \otimes \Dmod_{\frac{1}{2}}(\Gr_{G,\ul{x}})^{\fL^+(P^-)_{\ul{x}}}  \overset{\star}\to \\
\to \Dmod_{\frac{1}{2}}(\Bun_G\underset{\on{pt}/\fL^+(G)_{\ul{x}}}\times \on{pt}/\fL^+(P^-)_{\ul{x}}) \overset{({}'\sfp^-)^!}\longrightarrow 
\Dmod_{\frac{1}{2}}(\Bun_{P^-})
\end{multline*}
identifies with the functor \eqref{e:pre pre CT enh}.

\medskip

Hence, the upper horizontal arrow in \eqref{e:conv and CT} is the functor \eqref{e:pre CT enh}. In particular, \eqref{e:conv and CT} factors via 
\begin{equation} \label{e:conv and CT coinv}
\vcenter
{\xy
(0,0)*+{\Dmod_{\frac{1}{2}}(\Bun_G) \otimes \Dmod_{\frac{1}{2}}(\Gr_{G,\ul{x}})_{\fL(N^-_P)_{\ul{x}}}^{\fL^+(M)_{\ul{x}}}}="A";
(70,0)*+{\Dmod_{\frac{1}{2}}(\Bun_M)}="B";
(0,-40)*+{\KL(G)_{\crit_G,\ul{x}} \otimes \Dmod_{\frac{1}{2}}(\Gr_{G,\ul{x}})_{\fL(N^-_P)_{\ul{x}}}^{\fL^+(M)_{\ul{x}}}}="C"; 
(70,-40)*+{\KL(M)_{\crit_M+\rhoch_P,\ul{x}}.}="D";
{\ar@{->} "A";"B"};
{\ar@{->} "C";"D"};
{\ar@{->}^{\Loc_{G,\ul{x}}\otimes \on{Id}} "C";"A"};
{\ar@{->}_{\fl_{N^-_P}\otimes \Loc_{M,\rhoch_P,\ul{x}}} "D";"B"};
{\ar@{=>} "D";"A"};
\endxy}
\end{equation} 

\sssec{}

Precomposing \eqref{e:conv and CT coinv pre} with 
$$\on{I}(G,P^-)_{\on{co},\ul{x}}^{\on{loc}}\to \Dmod_{\frac{1}{2}}(\Gr_{G,\ul{x}})_{\fL(N^-_P)_{\ul{x}}}^{\fL^+(M)_{\ul{x}}},$$
we obtain a natural transformation 
\begin{equation} \label{e:conv and CT coinv pre}
\vcenter
{\xy
(0,0)*+{\Dmod_{\frac{1}{2}}(\Bun_G) \otimes \on{I}(G,P^-)_{\on{co},\ul{x}}^{\on{loc}}}="A";
(70,0)*+{\Dmod_{\frac{1}{2}}(\Bun_M)}="B";
(0,-40)*+{\KL(G)_{\crit_G,\ul{x}}  \otimes \on{I}(G,P^-)_{\on{co},\ul{x}}^{\on{loc}}}="C"; 
(70,-40)*+{\KL(M)_{\crit_M+\rhoch_P,\ul{x}},}="D";
{\ar@{->}^{\on{pre-CT}_{*,\ul{x}}^{-,\on{enh}_{\on{co}}}} "A";"B"};
{\ar@{->} "C";"D"};
{\ar@{->}^{\Loc_{G,\ul{x}}\otimes \on{Id}} "C";"A"};
{\ar@{->}_{\fl_{N^-_P}\otimes \Loc_{M,\rhoch_P,\ul{x}}} "D";"B"};
{\ar@{=>} "D";"A"};
\endxy}
\end{equation} 

The functors and the natural transformation in \eqref{e:conv and CT coinv pre} are compatible with\footnote{I.e., the horizontal arrows 
coequalize the $\Sph_G$ actions, and the above diagram is compatible with this structure.} 
the action of $\Sph_{G,x}$. Hence, it further factors via 
\begin{equation} \label{e:conv and CT coinv 1}
\vcenter
{\xy
(0,0)*+{\Dmod_{\frac{1}{2}}(\Bun_G) \underset{\Sph_{G,\ul{x}}}\otimes \on{I}(G,P^-)_{\on{co},\ul{x}}^{\on{loc}}}="A";
(70,0)*+{\Dmod_{\frac{1}{2}}(\Bun_M)}="B";
(0,-40)*+{\KL(G)_{\crit_G,\ul{x}}  \underset{\Sph_{G,\ul{x}}}\otimes \on{I}(G,P^-)_{\on{co},\ul{x}}^{\on{loc}}}="C"; 
(70,-40)*+{\KL(M)_{\crit_M+\rhoch_P,\ul{x}},}="D";
{\ar@{->}^{\on{CT}_{*,\ul{x}}^{-,\on{enh}_{\on{co}}}} "A";"B"};
{\ar@{->}^{J_{\on{KM}}^{-,\on{enh}_{\on{co}}}} "C";"D"};
{\ar@{->}^{\Loc_{G,\ul{x}}\otimes \on{Id}} "C";"A"};
{\ar@{->}_{\fl_{N^-_P}\otimes \Loc_{M,\rhoch_P,\ul{x}}} "D";"B"};
{\ar@{=>} "D";"A"};
\endxy}
\end{equation} 
as desired. 

\medskip

The compatibility with the $\Sph_{M,x}$-action follows from the construction. 

\ssec{Addendum: Whittaker coefficients of localization, take II}

In this subsection we use the methods developed above in order to sharpen a theorem from \cite[Sect. 14]{GLC2},
which computes Whittaker coefficients of Kac-Moody localization. 

\medskip

Namely, we will show that the computation
in {\it loc. cit.} is compatible with the action of $\Sph_G$. 

\sssec{}

Let us denote by $\sP^{\on{loc,true}}_G$ the (lax unital) factorization functor
\begin{multline*}
\KL(G)_{\crit}\underset{\Sph_G}\otimes \Whit_*(G)
\overset{\alpha_{\on{taut},\rho(\omega_X)}\otimes \on{Id}}\longrightarrow \\
\to \KL(G)_{\crit,\rho(\omega_X)}\underset{\Sph_G}\otimes \Whit_*(G)\overset{\star}\to 
(\hg\mod_{\crit,\rho(\omega_X)})_{\fL(N)_{\rho(\omega_X)},\chi}\overset{\ol\DS}\longrightarrow \Vect,
\end{multline*}
where the last arrow is \cite[Equation (2.8)]{GLC2}.

\medskip

Recall that the functor $\sP^{\on{loc}}_G$ is defined to be the composition
$$\KL(G)_{\crit}\otimes \Whit_*(G)\to
\KL(G)_{\crit}\underset{\Sph_G}\otimes \Whit_*(G)\overset{\sP^{\on{loc,true}}_G}\to \Vect,$$
see \cite[Sect. 6.4.6]{GLC2}. 

\sssec{} \label{sss:Whit Loc comp}

In \cite[Theorem 14.1.4]{GLC2} it was shown that the functor
\begin{multline}  \label{e:coeff Loc LHS}
(\KL(G)_\crit\otimes \Whit_*(G))_\CZ
\overset{\Loc_{G,\CZ}\otimes \on{Id}}\longrightarrow 
(\Dmod_{\frac{1}{2}}(\Bun_G)\otimes \Dmod(\CZ))\underset{\Dmod(\CZ)} \otimes \Whit_*(G)_\CZ \to \\
\overset{\on{coeff}_{G,\CZ}\otimes \on{id}}\longrightarrow (\Whit^!(G)\otimes \Whit_*(G))_\CZ\to \Dmod(\CZ)
\end{multline}
is canonically isomorphic to
\begin{multline} \label{e:coeff Loc RHS}
(\KL(G)_\crit\otimes \Whit_*(G))_\CZ
\overset{\on{ins.unit}_\CZ}\longrightarrow (\KL(G)_\crit\otimes \Whit_*(G))_{\CZ^\subseteq} \overset{\sP^{\on{loc}}_G}\longrightarrow \\
\to \Dmod(\CZ^\subseteq)\overset{(\on{pr}_{\on{small},\CZ})_!}\longrightarrow 
\Dmod(\CZ)\overset{\fl^{\otimes -\frac{1}{2}}_{G,N_{\rho(\omega_X)}}\otimes \fl_{N_{\rho(\omega_X)}}[\delta_{N_{\rho(\omega_X)}}]}\longrightarrow \Dmod(\CZ),
\end{multline}
where:
\begin{itemize}

\item $\delta_{N_{\rho(\omega_X)}}:=\dim(\Bun_{N_{\rho(\omega_X)}})$;

\medskip

\item $\fl_{N_{\rho(\omega_X)}}$ is the line from \cite[Sect. 14.1.1]{GLC2}, i.e.,
$$\det(T^*(\Bun_{N_{\rho(\omega_X)}}))\simeq \CO_{\Bun_{N_{\rho(\omega_X)}}}\otimes \fl_{N_{\rho(\omega_X)}}[\delta_{N_{\rho(\omega_X)}}];$$

\item $\fl^{\otimes \frac{1}{2}}_{G,N_{\rho(\omega_X)}}$ is the line from \cite[Sect. 9.2.4]{GLC2}.\footnote{See \cite[Sect. 14.1.4]{GLC2} for why the twist by this
line appears.}

\end{itemize}

\sssec{}

Note now that both \eqref{e:coeff Loc LHS} and \eqref{e:coeff Loc RHS} factor via 
$$(\KL(G)_\crit\otimes \Whit_*(G))_\CZ\to (\KL(G)_\crit\underset{\Sph_G}\otimes \Whit_*(G))_\CZ\to \Dmod(\CZ).$$

Indeed, for \eqref{e:coeff Loc LHS} this factorization is provided by
\begin{multline}  \label{e:coeff Loc LHS Sph}
(\KL(G)_\crit\underset{\Sph_G}\otimes \Whit_*(G))_\CZ
\overset{\Loc_{G,\CZ}\otimes \on{Id}}\longrightarrow 
(\Dmod_{\frac{1}{2}}(\Bun_G)\otimes \Dmod(\CZ))\underset{\Sph_{G,\CZ}}\otimes \Whit_*(G)_\CZ \to \\
\overset{\on{coeff}_{G,\CZ}\otimes \on{Id}}\longrightarrow (\Whit^!(G)\underset{\Sph_{G,\CZ}}\otimes \Whit_*(G))_\CZ\to \Dmod(\CZ).
\end{multline}

For \eqref{e:coeff Loc RHS} this factorization is provided by
\begin{multline} \label{e:coeff Loc RHS Sph}
(\KL(G)_\crit\underset{\Sph_G}\otimes \Whit_*(G))_\CZ
\overset{\on{ins.unit}_\CZ}\longrightarrow (\KL(G)_\crit\underset{\Sph_G}\otimes \Whit_*(G))_{\CZ^\subseteq} 
\overset{\sP^{\on{loc,true}}_G}\longrightarrow \\
\to \Dmod(\CZ^\subseteq)\overset{(\on{pr}_{\on{small},\CZ})_!}\longrightarrow 
\Dmod(\CZ)\overset{\fl^{\otimes -\frac{1}{2}}_{G,N_{\rho(\omega_X)}}\otimes \fl_{N_{\rho(\omega_X)}}[\delta_{N_{\rho(\omega_X)}}]}\longrightarrow \Dmod(\CZ).
\end{multline}

\sssec{} \label{sss:glob true}

Denote the functor \eqref{e:coeff Loc LHS} by $\sP^{\on{glob}}_{G,\CZ}$
and the functor \eqref{e:coeff Loc LHS Sph} by $\sP^{\on{glob,true}}_{G,\CZ}$.

\medskip

We can view the assignments
$$\CZ\rightsquigarrow \sP^{\on{glob}}_{G,\CZ} \text{ and } \CZ\rightsquigarrow \sP^{\on{glob,true}}_{G,\CZ}$$
as (unital) local-to-global functors, to be denoted
$$\ul\sP^{\on{glob}}_G \text{ and } \ul\sP^{\on{glob,true}}_G,$$
respectively. 

\medskip

Note that we can think of the functor \eqref{e:coeff Loc RHS Sph}, up to tensoring by a line and a cohomological shift. 
as
$$(\ul\sP^{\on{loc,true}}_G)^{\int\, \on{ins.unit}}_\CZ.$$

\sssec{}

In this subsection we will prove the following sharpening of \cite[Theorem 14.1.4]{GLC2}:

\begin{thm} \label{t:coeff Loc}
We have a canonical isomorphism 
$$(\ul\sP^{\on{loc,true}}_G)^{\int\, \on{ins.unit}}\otimes \fl^{\otimes -\frac{1}{2}}_{G,N_{\rho(\omega_X)}}\otimes \fl_{N_{\rho(\omega_X)}}[\delta_{N_{\rho(\omega_X)}}]\simeq \ul\sP^{\on{glob,true}}_G$$
as local-to-global functors
$$\ul\KL(G)_\crit\underset{\ul\Sph_G}\otimes \ul\Whit_*(G)\rightrightarrows \Vect.$$
\end{thm} 

The rest of this subsection is devoted to the proof of this theorem. 

\sssec{}

Parallel to the functor $\on{CT}^{-,\on{enh}_{\on{co}}}_{*,\CZ}$ constructed in \secref{ss:co enh CT}, one can consider the functor
\begin{equation} \label{e:coeff enr}
(\Dmod_{\frac{1}{2}}(\Bun_G)\otimes \Dmod(\CZ))\underset{\Sph_{G,\CZ}}\otimes \Whit_*(G)_\CZ\to \Dmod(\CZ),
\end{equation} 
so that the composition
\begin{multline*}
\Dmod_{\frac{1}{2}}(\Bun_G)\otimes \Dmod(\CZ)
\overset{\on{Id}\otimes \one_{\Whit_*(G)_\CZ}}\longrightarrow \\
\to (\Dmod_{\frac{1}{2}}(\Bun_G)\otimes \Dmod(\CZ))\underset{\Sph_{G,\CZ}}\otimes \Whit_*(G)_\CZ\to \Dmod(\CZ)
\end{multline*}
identifies with
$$\on{coeff}_G^{\on{Vac,glob}}\otimes \on{Id}_{\Dmod(\CZ)},$$
where $\on{coeff}_G^{\on{Vac,glob}}$ is the functor from \cite[Sect. 9.3.6]{GLC2}. 

\sssec{}

A construction parallel to that in \secref{ss:nat transf}, combined with \cite[Sect. 12.5.16-12.5.17]{GLC2}
gives rise to a natural transformation
$$
\xy 
(0,0)*+{(\Dmod_{\frac{1}{2}}(\Bun_G)\otimes \Dmod(\CZ))\underset{\Sph_{G,\CZ}}\otimes \Whit_*(G)_\CZ}="A";
(90,0)*+{\Dmod(\CZ)}="B";
(0,-30)*+{\KL(G)_{\crit,\CZ}\underset{\Sph_{G,\CZ}}\otimes \Whit_*(G)_\CZ}="C"; 
(90,-30)*+{\Dmod(\CZ).}="D";
{\ar@{->}^{\text{\eqref{e:coeff enr}}} "A";"B"};
{\ar@{->}_{\sP^{\on{loc,true}}_G} "C";"D"};
{\ar@{->}^{\Loc_{G,\CZ}\otimes \on{Id}} "C";"A"};
{\ar@{->}_{\fl^{\otimes -\frac{1}{2}}_{G,N_{\rho(\omega_X)}}\otimes \fl_{N_{\rho(\omega_X)}}[\delta_{N_{\rho(\omega_X)}}]} "D";"B"};
{\ar@{=>} "D";"A"};
\endxy
$$

\medskip

As in Sects. \ref{sss:strat CT and Loc crit enh 1}-\ref{sss:strat CT and Loc crit enh 2}, one shows that the induced 
natural transformation 
$$
\xy 
(0,0)*+{(\Dmod_{\frac{1}{2}}(\Bun_G)\otimes \Dmod(\CZ))\underset{\Sph_{G,\CZ}}\otimes \Whit_*(G)_\CZ}="A";
(90,0)*+{\Dmod(\CZ)}="B";
(0,-30)*+{\KL(G)_{\crit,\CZ}\underset{\Sph_{G,\CZ}}\otimes \Whit_*(G)_\CZ}="C"; 
(90,-30)*+{\Dmod(\CZ)}="D";
(0,-60)*+{\KL(G)_{\crit,\CZ^\subseteq}\underset{\Sph_{G,\CZ^\subseteq}}\otimes \Whit_*(G)_{\CZ^\subseteq}}="E"; 
(90,-60)*+{\Dmod(\CZ^\subseteq)}="F";
{\ar@{->}^{\text{\eqref{e:coeff enr}}} "A";"B"};
{\ar@{->}_{\sP^{\on{loc,true}}_G} "E";"F"};
{\ar@{->}^{\Loc_{G,\CZ}\otimes \on{Id}} "C";"A"};
{\ar@{->}_{\fl^{\otimes -\frac{1}{2}}_{G,N_{\rho(\omega_X)}}\otimes \fl_{N_{\rho(\omega_X)}}[\delta_{N_{\rho(\omega_X)}}]} "D";"B"};
{\ar@{->}^{(\on{pr}_{\on{small},\CZ})_!} "F";"D"};
{\ar@{->}^{\on{ins.unit}_\CZ} "C";"E"};
{\ar@{=>} "F";"A"};
\endxy
$$
is an isomorphism. 

\sssec{}

However, now, as in \propref{p:loc vs global CT enh co} one shows that the functor \eqref{e:coeff enr} identifies 
with 
\begin{multline*}
(\Dmod_{\frac{1}{2}}(\Bun_G)\otimes \Dmod(\CZ))\underset{\Sph_{G,\CZ}}\otimes \Whit_*(G)_\CZ \to \\
\overset{\on{coeff}_{G,\CZ}\otimes \on{id}}\longrightarrow (\Whit^!(G)\underset{\Sph_{G,\CZ}}\otimes \Whit_*(G))_\CZ\to \Dmod(\CZ),
\end{multline*}
as required. 

\qed[\thmref{t:coeff Loc}]

\ssec{Localization and translation by \texorpdfstring{$Z_G$}{ZG}-torsors}

\sssec{} \label{sss:twist G-bundle global}

Let $\CP_G$ be a $G$-torsor on $X$. Parallel to \cite[Sect. 10]{GLC2}, we can consider the localization functor
$$\Loc_{G,\kappa,\CP_G,\CZ}:\KL(G)_{\kappa,\CP_G,\CZ}\to \Dmod_\kappa(\Bun_{G,\CP_G})\otimes \Dmod(\CZ),$$
which makes the diagram
\begin{equation} \label{e:G twisted loc diag}
\CD
\Dmod_\kappa(\Bun_G)\otimes \Dmod(\CZ) @>{\alpha_{\CP_G,\on{taut}}\otimes \on{Id}}>>  \Dmod_\kappa(\Bun_{G,\CP_G})\otimes \Dmod(\CZ) \\
@A{\Loc_{G,\kappa,\CZ}}AA @AA{\Loc_{G,\kappa,\CP_G,\CZ}}A \\
\KL(G)_{\kappa,\CZ} @>{\alpha_{\CP_G,\on{taut}}}>> \KL(G)_{\kappa,\CP_G,\CZ},
\endCD
\end{equation} 
commute, where in the top row, $\alpha_{\CP_G,\on{taut}}$ is the tautological identification
$$\Bun_G\overset{\alpha_{\CP_G,\on{taut}}}\longrightarrow \Bun_{G,\CP_G}.$$

\sssec{}

Assume now that $\CP_G$ is induced from a $Z_G$-torsor $\CP_{Z_G}$. In this case, we have a canonical identification
$$\Bun_{G,\CP_{Z_G}} \overset{\alpha_{\CP_{Z_G},\on{cent}}}\longrightarrow \Bun_G.$$

\medskip

Recall the notation $\kappa(\on{dlog}(\CP_{Z_G}),-)$ for the corresponding torsor over $\fg_{\on{ab}}^*\otimes \omega_X$,
see \secref{sss:kappa dlog}. We claim: 

\begin{lem} \label{l:global twistings shift}
The isomorphism $\alpha_{\CP_{Z_G},\on{cent}}$ identifies the twisting $\kappa$ on 
$\Bun_{G,\CP_{Z_G}}$ with the twisting 
$$\kappa-\kappa(\on{dlog}(\CP_{Z_G},-)$$
on $\Bun_G$. Furthermore, the diagram
\begin{equation} \label{e:ZG twisted loc diag}
\CD
\Dmod_\kappa(\Bun_{G,\CP_{Z_G}})\otimes \Dmod(\CZ) 
@>{\alpha_{\CP_{Z_G},\on{cent}}}>> \Dmod_{\kappa-\kappa(\on{dlog}(\CP_{Z_G},-)}(\Bun_G)\otimes \Dmod(\CZ)  \\
@A{\Loc_{G,\kappa,\CP_{Z_G},\CZ}}AA @AA{\Loc_{G,\kappa+\kappa(\on{dlog}(\CP_{Z_G},-),\CZ}}A \\
\KL(G)_{\kappa,\CP_{Z_G},\CZ} @>{\alpha_{\CP_{Z_G},\on{cent}}}>> \KL(G)_{\kappa+\kappa(\on{dlog}(\CP_{Z_G},-),\CZ}.
\endCD
\end{equation} 
\end{lem}

\sssec{} \label{sss:transl Bun}

Note that the composite map
$$\Bun_G\overset{\alpha_{\CP_{Z_G},\on{taut}}}\longrightarrow
\Bun_{G,\CP_{Z_G}} \overset{\alpha_{\CP_{Z_G},\on{cent}}}\longrightarrow \Bun_G$$
is the inverse of the automorphism 
$$\on{transl}_{\CP_{Z_G}}:\Bun_G\to \Bun_G,$$
given by tensoring with $\CP_{Z_G}$.

\sssec{}

From \lemref{l:global twistings shift} we obtain:

\begin{cor} \label{c:pullback transl twist}
The pullback of the twisting of the twisting $\kappa$ on $\Bun_G$ 
with respect to $\on{transl}_{\CP_{Z_G}}$
identifies canonically with the twisting $\kappa-\kappa(\on{dlog}(\CP_{Z_G},-)$. 
\end{cor} 

\sssec{}

Concatenating diagrams \eqref{e:G twisted loc diag} and \eqref{e:ZG twisted loc diag}, we obtain a commutative diagram
\begin{equation} \label{e:transl loc diag}
\CD
\Dmod_\kappa(\Bun_G)\otimes \Dmod(\CZ)  @>{(\on{transl}_{\CP_{Z_G}})^*}>> \Dmod_{\kappa-\kappa(\on{dlog}(\CP_{Z_G},-)}(\Bun_G)\otimes \Dmod(\CZ)   \\
@A{\Loc_{G,\kappa,\CZ}}AA  @AA{\Loc_{G,\kappa-\kappa(\on{dlog}(\CP_{Z_G},-),\CZ}}A \\
\KL(G)_{\kappa,\CZ} @>{(\on{transl}_{\CP_{Z_G}})^*}>>  \KL(G)_{\kappa+\kappa(\on{dlog}(\CP_{Z_G},-),\CZ},
\endCD
\end{equation} 
where the functor $(\on{transl}_{\CP_{Z_G}})^*$ in the bottom line is as in \secref{e:central un-twist of KL}. 

%
%
%

\sssec{} \label{sss:lambda twist glob}

Parallel to \secref{sss:lambda twist} for $\CP_{Z_G}$ of the form
$\lambda(\omega_X)$ for $\lambda:\BG_m\to Z_G$, we will use a shorthandnotation
$$\kappa-\kappa(\lambda,-):=\kappa-\kappa(\on{dlog}(\lambda(\omega_X),-).$$

Note that this notation is consistent with one in \secref{sss:global torsors lambda}, i.e., we can
treat $\kappa(\lambda,-)$ as a bona fide element of $z_\cg$. 

\ssec{Constant term of localization: the translated version}

\sssec{}

Let $\CP_{Z_M}$ be a $Z_M$-torsor. Let $\on{CT}^-_{*,\CP_{Z_M}}$ be the functor
$$\Dmod_{\crit_G+\kappa}(\Bun_G) \overset{\on{CT}^-_*}\longrightarrow 
\Dmod_{\crit_M+\kappa}(\Bun_M) \overset{(\on{transl}_{\CP_{Z_M}})^*}\longrightarrow
\Dmod_{\crit_M+\kappa-\kappa(\on{dlog}(\CP_{Z_M}),-)}(\Bun_M).$$

\sssec{}

Concatenating the diagrams \eqref{e:transl loc diag}, \eqref{e:CT and Loc kappa} and \eqref{e:BRST twisted CD sph}.
we obtain:

\begin{cor} \label{c:CT and Loc kappa transl}
For $\CZ\to \Ran$, the following diagram of functors canonically commutes
$$
\CD
\Dmod_{\crit_G+\kappa}(\Bun_G)\otimes \Dmod(\CZ) @>{\on{CT}^-_{*,\CP_{Z_M}}\otimes \on{Id}}>> 
\Dmod_{\crit_M+\kappa-\kappa(\on{dlog}(\CP_{Z_M}),-)}(\Bun_M)\otimes  \Dmod(\CZ)  \\
@A{\Loc_{G,\crit_G+\kappa,\CZ}}AA  
@A{\fl\otimes ((-)\otimes \CL_{\rhoch_P(\omega_X)})\otimes (\on{pr}_{\on{small},\CZ})_!}AA \\ 
\KL(G)_{\crit_G+\kappa,\CZ} & & \Dmod_{\crit_M+\kappa+\rhoch_P-\kappa(\on{dlog}(\CP_{Z_M}),-)}(\Bun_M)\otimes  \Dmod(\CZ) \\
@V{\on{ins.vac.}_\CZ}VV  @A{\Loc_{M,\crit_M+\kappa+\rhoch_P+\kappa(\on{dlog}(\CP_{Z_M}),-),\CZ^\subseteq}}AA \\
\KL(G)_{\crit_G+\kappa,\CZ^\subseteq} @>{J_{\on{KM},\CP_{Z_M}}^{-,\Sph}}>> 
\KL(M)_{\crit_M+\kappa+\rhoch_P+\kappa(\on{dlog}(\CP_{Z_M}),-),\CZ^\subseteq},
\endCD
$$
where 
$$\fl:=\fl_{N^-_P}\otimes \on{Weil}(\CP_{Z_M},\rhoch_P(\omega_X)).$$
\end{cor} 

Specializing to the critical level, we obtain:

\begin{cor} \label{c:CT and Loc crit transl}
For $\CZ\to \Ran$, the following diagram of functors canonically commutes
$$
\CD
\Dmod_{\frac{1}{2}}(\Bun_G)\otimes \Dmod(\CZ) @>{\on{CT}^-_{*,\CP_{Z_M}}\otimes \on{Id}}>> 
\Dmod_{\frac{1}{2}}(\Bun_M)\otimes  \Dmod(\CZ)  \\
@A{\Loc_{G,\CZ}}AA  @AA{\fl\otimes (\on{Id}\otimes (\on{pr}_{\on{small},\CZ})_!)}A \\ 
\KL(G)_{\crit_G,\CZ} & & \Dmod_{\frac{1}{2}}(\Bun_M)\otimes  \Dmod(\CZ^\subseteq) \\
@V{\on{ins.vac.}_\CZ}VV  @AA{\Loc_{M,\rhoch_P,\CZ^\subseteq}}A \\
\KL(G)_{\crit_G,\CZ^\subseteq} @>{J_{\on{KM},\CP_{Z_M}}^{-,\Sph}}>> 
\KL(M)_{\crit_M+\rhoch_P,\CZ^\subseteq},
\endCD
$$
where 
$$\fl:=\fl_{N^-_P}\otimes \on{Weil}(\CP_{Z_M},\rhoch_P(\omega_X)).$$
\end{cor} 

In the specific case of $\CP_{Z_M}=\rho_P(\omega_X)$ we obtain:

\begin{cor} \label{c:CT and Loc crit transl rho}
For $\CZ\to \Ran$, the following diagram of functors canonically commutes
$$
\CD
\Dmod_{\frac{1}{2}}(\Bun_G)\otimes \Dmod(\CZ) @>{\on{CT}^-_{*,\rho_P(\omega_X)}\otimes \on{Id}}>> 
\Dmod_{\frac{1}{2}}(\Bun_M)\otimes  \Dmod(\CZ)  \\
@A{\Loc_{G,\CZ}}AA  @AA{\fl\otimes (\on{Id}\otimes (\on{pr}_{\on{small},\CZ})_!)}A \\ 
\KL(G)_{\crit_G,\CZ} & & \Dmod_{\frac{1}{2}}(\Bun_M)\otimes  \Dmod(\CZ^\subseteq) \\
@V{\on{ins.vac.}_\CZ}VV  @AA{\Loc_{M,\rhoch_P,\CZ^\subseteq}}A \\
\KL(G)_{\crit_G,\CZ^\subseteq} @>{J_{\on{KM},\rho_P(\omega_X)}^{-,\Sph}}>> 
\KL(M)_{\crit_M+\rhoch_P,\CZ^\subseteq},
\endCD
$$
where 
$$\fl:=\fl_{N^-_P}\otimes \on{Weil}(\CP_{Z_M},\rhoch_P(\omega_X)).$$
\end{cor} 

\sssec{}

Similarly, from \thmref{t:CT and Loc crit enh}, we obtain: 

\begin{cor} \label{c:CT and Loc crit transl enh}
For $\CZ\to \Ran$, the following diagram of functors canonically commutes
$$
\CD
\Dmod_{\frac{1}{2}}(\Bun_G)^{-,\on{enh}_{\on{co}}}_\CZ @>{\on{CT}^{-,\on{enh}_{\on{co}}}_{*,\CP_{Z_M},\CZ}}>> 
\Dmod_{\frac{1}{2}}(\Bun_M)\otimes  \Dmod(\CZ)  \\
@A{\Loc^{-,\on{enh}_{\on{co}}}_{G,\CZ}}AA  @AA{\fl\otimes (\on{Id}\otimes (\on{pr}_{\on{small},\CZ})_!)}A \\ 
\KL(G)^{-,\on{enh}_{\on{co}}}_{\crit_G,\CZ} & & \Dmod_{\frac{1}{2}}(\Bun_M)\otimes  \Dmod(\CZ^\subseteq) \\
@V{\on{ins.vac.}_\CZ}VV  @AA{\Loc_{M,\rhoch_P,\CZ^\subseteq}}A \\
\KL(G)^{-,\on{enh}_{\on{co}}}_{\crit_G,\CZ^\subseteq}  
@>{J_{\on{KM},\CP_{Z_M}}^{-,\on{enh}_{\on{co}}}\otimes \fl_{N^-_P}\otimes \on{Weil}(\CP_{Z_M},\rhoch_P(\omega_X))}>> 
\KL(M)_{\crit_M+\rhoch_P,\CZ^\subseteq},
\endCD
$$
where 
$$\fl:=\fl_{N^-_P}\otimes \on{Weil}(\CP_{Z_M},\rhoch_P(\omega_X)).$$
\end{cor} 

In the specific case of $\CP_{Z_M}=\rho_P(\omega_X)$ we obtain: 

\begin{cor} \label{c:CT and Loc crit transl enh rho}
For $\CZ\to \Ran$, the following diagram of functors canonically commutes
$$
\CD
\Dmod_{\frac{1}{2}}(\Bun_G)^{-,\on{enh}_{\on{co}}}_\CZ @>{\on{CT}^{-,\on{enh}_{\on{co}}}_{*,\rho_P(\omega_X),\CZ}}>> 
\Dmod_{\frac{1}{2}}(\Bun_M)\otimes  \Dmod(\CZ)  \\
@A{\Loc^{-,\on{enh}_{\on{co}}}_{G,\CZ}}AA  @AA{\fl\otimes (\on{Id}\otimes (\on{pr}_{\on{small},\CZ})_!)}A \\ 
\KL(G)^{-,\on{enh}_{\on{co}}}_{\crit_G,\CZ} & & \Dmod_{\frac{1}{2}}(\Bun_M)\otimes  \Dmod(\CZ^\subseteq) \\
@V{\on{ins.vac.}_\CZ}VV  @AA{\Loc_{M,\rhoch_P,\CZ^\subseteq}}A \\
\KL(G)^{-,\on{enh}_{\on{co}}}_{\crit_G,\CZ^\subseteq}  
@>{J_{\on{KM},\rho_P(\omega_X)}^{-,\on{enh}_{\on{co}}}}>> 
\KL(M)_{\crit_M+\rhoch_P,\CZ^\subseteq},
\endCD
$$
where 
$$\fl:=\fl_{N^-_P}\otimes \on{Weil}(\rho_P(\omega_X),\rhoch_P(\omega_X)).$$
\end{cor} 

\sssec{}

For future reference, we record the following observation, which is a particular case of \propref{p:rel det}:

\begin{cor} \label{c:rel det transl}
There is a canonical isomorphism of lines 
$$
\fl_{N^-_P}\otimes \on{Weil}(\rho_P(\omega_X),\rhoch_P(\omega_X))
\simeq 
\fl_{G,N_{\rho(\omega_X)}}^{\otimes -\frac{1}{2}}\otimes \fl_{N_{\rho(\omega_X)}}\otimes 
\fl^{\otimes \frac{1}{2}}_{M,N(M)_{\rho_M(\omega_X)}}\otimes \fl^{\otimes -1}_{N(M)_{\rho_M(\omega_X)}},$$
where:

\begin{itemize}

\item $\fl_{G,N_{\rho(\omega_X)}}^{\otimes \frac{1}{2}}$ is the square root of the line $\fl_{G,N_{\rho(\omega_X)}}$,
given by \cite[Proposition 1.3.3]{GLC1}, and 
$\fl^{\otimes \frac{1}{2}}_{M,N(M)_{\rho_M(\omega_X)}}$ is the corresponding line for $M$;

\item $\fl_{N_{\rho(\omega_X)}}$ is the line from \cite[Sect. 14.1.1]{GLC2}, and $\fl_{N(M)_{\rho_M(\omega_X)}}$
is the corresponding line for $M$.

\end{itemize}

\end{cor}

\section{Spectral Eisenstein and constant term functors}  \label{s:spec CT}

In this section start dealing with the global spectral category, i.e., $\IndCoh(\LS_\cG)$,
where $\LS_\cG$ is the stack of $\cG$-local systems on $X$. 

\medskip

We develop the spectral counterparts of the constructions in Sects.
\ref{s:Eis}-\ref{s:Eis glob}, i.e., the spectral Eisenstein and constant term functors,
as well as their enhancements. 

\medskip

We will perform a spectral analog of the ``localization and constant term" calculation
from \ref{s:Eis and Loc}. The analog of the ``Whittaker vs Eisenstein" calculation\footnote{Note, however, that the Langlands functor
takes the localization functor on the geometric side to the Poincar\'e functor on the spectral side.}  
from \secref{s:Eis and Whit}
will be dealt with in \secref{s:Eis and Op}. 

\ssec{The usual spectral Eisenstein and constant term functors}

\sssec{}

Consider the diagram
$$\LS_\cG  \overset{\sfp^{-,\on{spec}}}\longleftarrow \LS_{\cP^-} \overset{\sfq^{-,\on{spec}}}\longrightarrow \LS_{\cM}.$$

The (usual) spectral Eisenstein and constant term functors are the adjoint pair
$$\Eis^{-,\on{spec}}:\IndCoh(\LS_{\cM}) \leftrightarrows \IndCoh(\LS_{\cG}):\on{CT}^{-,\on{spec}}$$
defined by
$$\Eis^{-,\on{spec}}:=(\sfp^{-,\on{spec}})^\IndCoh_*\circ (\sfq^{-,\on{spec}})^{\IndCoh,*}$$
and 
$$\on{CT}^{-,\on{spec}}:=(\sfq^{-,\on{spec}})^\IndCoh_* \circ (\sfp^{-,\on{spec}})^!,$$
where the functor $(\sfq^{-,\on{spec}})^{\IndCoh,*}$ is well-defined since the morphism $\sfq^{-,\on{spec}}$ has a finite
Tor-dimension (in fact, the morphism $\sfq^{-,\on{spec}}$ is quasi-smooth).

\sssec{}

It is shown in \cite[Proposition 13.2.6]{AG} that the functor $\Eis^{-,\on{spec}}$ sends 
$$\IndCoh_\Nilp(\LS_\cM)\to \IndCoh_\Nilp(\LS_\cG).$$

\medskip

Similarly, it follows from \cite[Theorem 7.1.3]{AG} that the functor $\on{CT}^{-,\on{spec}}$ also sends 
$$\IndCoh_\Nilp(\LS_\cG)\to \IndCoh_\Nilp(\LS_\cM).$$ 

Hence, the functors $(\Eis^{-,\on{spec}},\on{CT}^{-,\on{spec}})$ induce an adjoint pair
$$\IndCoh_\Nilp(\LS_\cM)\rightleftarrows \IndCoh_\Nilp(\LS_\cG).$$

\sssec{}

In addition to the functor $\Eis^{-,\on{spec}}$ we will consider the functor
$$\Eis^{-,\on{spec}}_\QCoh:\QCoh(\LS_{\cM})\to \QCoh(\LS_{\cG})$$
defined by
$$\Eis^{-,\on{spec}}_\QCoh:=(\sfp^{-,\on{spec}})_*\circ (\sfq^{-,\on{spec}})^*.$$

We have a commutative diagram
\begin{equation} \label{e:Eis QCoh}
\CD
\IndCoh(\LS_\cM) @>{\Psi_{\LS_\cM}}>> \QCoh(\LS_\cM) \\
@V{\Eis^{-,\on{spec}}}VV  @VV{\Eis^{-,\on{spec}}_\QCoh}V \\
\IndCoh(\LS_\cG) @>{\Psi_{\LS_\cG}}>> \QCoh(\LS_\cG). 
\endCD
\end{equation}

\sssec{}

Let 
$$\on{CT}^{-,\on{spec}}_\QCoh:\QCoh(\LS_{\cG})\to \QCoh(\LS_{\cM})$$
be the functor defined by 
$$\on{CT}^{-,\on{spec}}_\QCoh:=(\sfq^{-,\on{spec}})_*\circ (\sfp^{-,\on{spec}})^*.$$

Note that we have
\begin{equation} \label{e:CT and Eis QCoh dual}
(\Eis^{-,\on{spec}}_\QCoh)^\vee \simeq \on{CT}^{-,\on{spec}}_\QCoh
\end{equation} 
with respect to the canonical self-duality
\begin{equation} \label{e:naive self-dual QCoh}
\QCoh(\LS_{\cG})^\vee \simeq \QCoh(\LS_{\cG}),
\end{equation}
and a similar self-duality for $\cM$. 

\sssec{}

Consider the (ungraded) line bundle $\det(T^*(\LS_{\cP^-}/\LS_\cG))$ on $\LS_{\cP^-}$.
It is easy to see that it canonically descends to a line bundle on $\LS_\cM$; by a slight abuse
of notation, we will denote it by the same character $\det(T^*(\LS_{\cP^-}/\LS_\cG))$.

\medskip

We claim:

\begin{lem} \label{l:CT spec IndCoh vs QCoh}
There is a commutative diagram
$$
\CD
\QCoh(\LS_\cM) @>{\left((-)\otimes \det(T^*(\LS_{\cP^-}/\LS_\cG))\right)\circ \Xi_{\LS_\cM}[-(2g-2)\cdot \dim(\cn_P)]}>> \IndCoh(\LS_\cM)  \\
@A{\on{CT}^{-,\on{spec}}_\QCoh}AA @AA{\on{CT}^{-,\on{spec}}}A \\
\QCoh(\LS_\cG) @>{\Xi_{\LS_\cG}}>> \IndCoh(\LS_\cG). 
\endCD
$$
\end{lem} 

\begin{proof}

Since both $\LS_\cG$ and $\LS_{\cP^-}$ are quasi-smooth, we have a commutative diagram
$$
\CD
\QCoh(\LS_{\cP^-}) @>{\left((-)\otimes \det(T^*(\LS_{\cP^-}/\LS_\cG))\right)\circ \Xi_{\LS_{\cP^-}}[-(2g-2)\cdot \dim(\cn_P)]}>> \IndCoh(\LS_{\cP^-})  \\
@A{(\sfp^{-,\on{spec}})^*}AA @AA{(\sfp^{-,\on{spec}})^!}A \\
\QCoh(\LS_\cG) @>{\Xi_{\LS_\cG}}>> \IndCoh(\LS_\cG), 
\endCD
$$
where 
$$(2g-2)\cdot \dim(\cn_P)=\dim(\LS_\cG)-\dim(\LS_{\cP^-}).$$

Since $\sfq^{-,\on{spec}}$ is quasi-smooth, the diagram
$$
\CD
\QCoh(\LS_{\cP^-}) @>{\Xi_{\LS_{\cP^-}}}>> \IndCoh(\LS_{\cP^-})  \\
@V{(\sfq^{-,\on{spec}})_*}VV @VV{(\sfq^{-,\on{spec}})^\IndCoh_*}V \\
\QCoh(\LS_\cM) @>{\Xi_{\LS_\cM}}>> \IndCoh(\LS_\cM)
\endCD
$$
also commutes.

\medskip

Concatenating the two diagrams, we obtain the assertion of the lemma. 

\end{proof}

\sssec{}

Recall that Serre duality defines an identification
\begin{equation} \label{e:Serre self-dual IndCoh}
\IndCoh(\LS_{\cG})^\vee \simeq \IndCoh(\LS_{\cG}),
\end{equation}
and similarly for $\cM$.

\medskip

Since $\sfq^{-,\on{spec}}$ is quasi-smooth, with respect to these identifications, we have
\begin{equation} \label{e:CT and Eis dual}
(\on{CT}^{-,\on{spec}})^\vee \simeq \Eis^{-,\on{spec}} \circ \Bigl((-)\otimes \det(T^*(\LS_{\cP^-}/\LS_\cM))\Bigr)[(2g-2)\cdot \dim(\cn^-_P)]
\end{equation}
as functors
$$\IndCoh(\LS_\cM)\rightrightarrows \IndCoh(\LS_\cG),$$
where:

\begin{itemize}

\item By a slight abuse of notation we denote by $\det(T^*(\LS_{\cP^-}/\LS_\cM))$ the natural descent
of the same-named line bundle from $\LS_{\cP^-}$ to $\LS_\cM$;

\medskip

\item The integer $(2g-2)\cdot \dim(\cn^-_P)$ appears as $\dim(\LS_{\cP^-})-\dim(\LS_{\cM})$.

\end{itemize}

\sssec{}

For future reference we note:

\begin{lem} \label{l:rel det LS}
There is a canonical isomorphism of (ungraded) line bundles
$$\det(T^*(\LS_{\cP^-}/\LS_\cM))\simeq \CL_{2\rho_P(\omega_X)}$$
and
$$\det(T^*(\LS_{\cP^-}/\LS_\cG))\simeq \CL_{-2\rho_P(\omega_X)}$$
where by a slight abuse of notation we denote by the same symbol $\CL_{2\rho_P(\omega_X)}$ the pullback
of the corresponding line bundle along
$$\LS_\cM\to \Bun_\cM.$$
\end{lem}

\begin{proof}

Let $\CE$ be a vector bundle with a connection on $X$. First, we claim that there is as canonical isomorphism
\begin{equation} \label{e:det dR}
\det(\on{C}^\cdot(X,\CE))\simeq \on{Weil}(\det(\CE),\omega_X)^{\otimes -1}.
\end{equation}

Indeed, we calculate $\on{C}^\cdot(X,\CE)$ as 
$$\on{Fib}\left(\Gamma(X,\CE)\overset{\nabla}\to \Gamma(X,\CE\otimes \omega_X)\right),$$
and hence
\begin{equation} \label{e:det dR 1}
\det(\on{C}^\cdot(X,\CE))^{\otimes -1}\simeq \det(\Gamma(X,\CE\otimes \omega_X))\otimes \det(\Gamma(X,\CE))^{\otimes -1}.
\end{equation}

Applying \eqref{e:det formula}, we rewrite the right-hand side in \eqref{e:det dR 1} as 
\begin{equation} \label{e:det dR 2}
\det(\Gamma(X,\omega_X))^{\otimes \on{rk}(\CE)}\otimes \det(\Gamma(X,\CO_X))^{\otimes -\on{rk}(\CE)}\otimes \on{Weil}(\det(\CE),\omega_X).
\end{equation}

Applying \eqref{e:det dR 1} to $\CE=\CO_X$, we obtain 
$$\det(\Gamma(X,\omega_X))\otimes \det(\Gamma(X,\CO_X))^{\otimes -1}\simeq \det(\on{C}^\cdot(X))^{\otimes -1}.$$

Note, however, the line $\det(\on{C}^\cdot(X))$ is canonically trivialized by the symplectic structure on $\on{C}^\cdot(X)[1]$. 
Thus, the expression in \eqref{e:det dR 2} identifies with $\on{Weil}(\det(\CE),\omega_X)$, as desired. 

\medskip

We now prove 
$$\det(T^*(\LS_{\cP^-}/\LS_\cM))\simeq \CL_{2\rho_P(\omega_X)};$$
the second isomorphism is proved similarly.

\medskip 

For $\sigma=(\CP_{\cP^-},\nabla)\in \LS_{\cP^-}$, the relative tangent space $T_\sigma(\LS_{\cP^-}/\LS_\cM)$
identifies with
$$\on{C}^\cdot(X,(\cn^-_P)_\sigma)[1].$$

Hence, by \eqref{e:det dR},
\begin{equation} \label{e:rel det LS 1}
\det(T^*(\LS_{\cP^-}/\LS_\cM))_\sigma\simeq \on{Weil}(\det((\cn^-_P)_{\CP_{\cP^-}}),\omega_X)^{\otimes -1}.
\end{equation}

We have
\begin{multline*}
\on{Weil}(\det((\cn^-_P)_{\CP_{\cP^-}}),\omega_X)\simeq \on{Weil}(-2\rho_P(\CP_{\cP^-}),\omega_X)\otimes 
\on{Weil}(\det(\cn^-_P)\otimes \CO_X,\omega_X)\simeq \\
\simeq \on{Weil}(-2\rho_P(\CP_{\cP^-}),\omega_X)\otimes \det(\cn^-_P)^{\otimes (2g-2)}=
(\CL_{-2\rho_P(\omega_X)})_{\CP_{\cP^-}}\otimes \det(\cn^-_P)^{\otimes (2g-2)}.
\end{multline*} 

Finally, we claim that the line $\det(\cn^-_P)^{\otimes 2}$ is canonically trivial. Indeed, the Killing form on $\cg$ gives rise
to an identification
$$\det(\cn^-_P) \simeq \det(\cn_P)^{\otimes -1},$$
while the Chevalley involution on $\cg$ gives rise
to an identification
$$\cn^-_P\simeq \cn_P.$$

\end{proof} 

\ssec{Spectral Eisenstein functor and localization} \label{ss:Eis and Loc}

The goal of this section is to express the composition
$$\Gamma_\cG^{\on{spec},\IndCoh}\circ \Eis^{-,\on{spec}}$$
via $\Gamma_\cM^{\on{spec},\IndCoh}$ and a functor of local nature. 

\medskip

The main result is given by \corref{c:Eis spec Gamma}. 

\sssec{}

For $\CZ\to \Ran$, recall the functor
$$\Gamma^{\on{spec}}_{\cG,\Ran}:\QCoh(\LS_\cG)\otimes \Dmod(\CZ) \to \Rep(\cG)_\CZ,$$
see \cite[Sect. 17.6.2]{GLC2}. 

\medskip

Explicitly, it is given by direct image along 
$$\LS_\cG\times \CZ\overset{\on{ev}_\CZ}\longrightarrow \LS^\reg_{\cG,\CZ}.$$

\sssec{}

From the commutative diagram
\begin{equation} \label{e:Eis spec Gamma diag}
\CD
\LS^\reg_{\cG,\CZ}  @<{\sfp^{-,\on{spec}}}<< \LS^\reg_{\cP^-,\CZ}  @>{\sfq^{-,\on{spec}}}>> \LS^\reg_{\cM,\CZ} \\
@A{\on{ev}_\CZ}AA @A{\on{ev}_\CZ}AA @AA{\on{ev}_\CZ}A \\
\LS_\cG\times \CZ @<{\sfp^{-,\on{spec}}}<< \LS_{\cP^-}\times \CZ @>{\sfq^{-,\on{spec}}}>> \LS_{\cM}\times \CZ 
\endCD
\end{equation} 
we obtain a natural transformation
\begin{equation} \label{e:Eis spec Gamma nat}
\vcenter
{\xy
(0,0)*+{\QCoh(\LS_\cG)\otimes \Dmod(\CZ)}="A";
(50,0)*+{\QCoh(\LS_\cM)\otimes \Dmod(\CZ) ,}="B";
(0,20)*+{\Rep(\cG)_\CZ}="C";
(50,20)*+{\Rep(\cM)_\CZ}="D";
{\ar@{->}^{\Eis^{\on{spec}}_\QCoh} "B";"A"};
{\ar@{->}_{\on{co}\!J^{-,\on{spec},!}} "D";"C"};
{\ar@{->}^{\Gamma^{\on{spec}}_{\cG,\CZ}} "A";"C"};
{\ar@{->}_{\Gamma^{\on{spec}}_{\cM,\CZ}} "B";"D"};
{\ar@{=>} "D";"A"};
\endxy}
\end{equation} 
where 
$$\on{co}\!J^{-,\on{spec},!}: \Rep(\cM)\to \Rep(\cG)$$
is the functor dual to $J^{-,\on{spec},!}$, i.e., this is the functor 
$$\Rep(\cM)\overset{\Res^\cM_{\cP^-}}\longrightarrow \Rep(\cP^-)\overset{\on{coind}^\cG_{\cP^-}}\longrightarrow \Rep(\cG).$$

\sssec{}

Applying \eqref{e:Eis spec Gamma nat} to $\CZ=\Ran$ and taking into account that the essential image of the functor
$$\Gamma_\cG:\QCoh(\LS_\cG)\to \Rep(\cG)_\Ran$$
lands in the full subcategory
$$\Rep(\cG)_{\Ran^{\on{untl}},\on{indep}}\subset \Rep(\cG)_\Ran,$$
we obtain a natural transformation 
\begin{equation} \label{e:Eis spec Gamma nat 1}
\vcenter
{\xy
(-5,0)*+{\QCoh(\LS_\cG)}="A";
(85,0)*+{\QCoh(\LS_\cM).}="B";
(40,20)*+{\Rep(\cG)_{\Ran^{\on{untl}}}}="C";
(85,20)*+{\Rep(\cM)_{\Ran^{\on{untl}}}}="D";
(-5,20)*+{\Rep(\cG)_{\Ran^{\on{untl}},\on{indep}}}="E";
{\ar@{->}^{\Eis^{\on{spec}}_\QCoh} "B";"A"};
{\ar@{->}^{\Gamma^{\on{spec}}_\cG} "A";"E"};
{\ar@{->}_{\Gamma^{\on{spec}}_{\cM,\on{untl}}} "B";"D"};
{\ar@{->}_{\on{co}\!J^{-,\on{spec},!}} "D";"C"};
{\ar@{->}_{\on{emb.indep}_{\Rep(\cG)}^L} "C";"E"};
{\ar@{=>} "D";"A"};
\endxy}
\end{equation} 

We will prove:

\begin{thm} \label{t:Eis spec Gamma}
The natural transformation \eqref{e:Eis spec Gamma nat 1} is an isomorphism.
\end{thm}

This theorem is a spectral counterpart of \corref{c:coeff of Eis !}. We will deduce it from 
\thmref{t:CT of Loc spec QCoh}, which is a spectral analog of \thmref{t:CT co of Poinc}. 

\sssec{}

Precomposing \eqref{e:Eis spec Gamma nat} with the functors $\Psi_{\LS_\cG}$ and $\Psi_{\LS_\cM}$ and using \eqref{e:Eis QCoh}, we obtain a natural
transformation
\begin{equation} \label{e:Eis spec Gamma nat IndCoh}
\vcenter
{\xy
(0,0)*+{\IndCoh(\LS_\cG)\otimes \Dmod(\CZ)}="A";
(70,0)*+{\IndCoh(\LS_\cM)\otimes \Dmod(\CZ),}="B";
(0,30)*+{\Rep(\cG)_\CZ}="C";
(70,30)*+{\Rep(\cM)_\CZ}="D";
{\ar@{->}^{\Eis^{-,\on{spec}}} "B";"A"};
{\ar@{->}_{\on{co}\!J^{-,\on{spec},!}} "D";"C"};
{\ar@{->}^{\Gamma^{\on{spec},\IndCoh}_{\cG,\CZ}} "A";"C"};
{\ar@{->}_{\Gamma^{\on{spec},\IndCoh}_{\cM,\CZ}} "B";"D"};
{\ar@{=>} "D";"A"};
\endxy}
\end{equation} 
and further
\begin{equation} \label{e:Eis spec Gamma nat 1 IndCoh}
\vcenter
{\xy
(-5,0)*+{\IndCoh(\LS_\cG)}="A";
(90,0)*+{\IndCoh(\LS_\cM).}="B";
(45,30)*+{\Rep(\cG)_{\Ran^{\on{untl}}}}="C";
(90,30)*+{\Rep(\cM)_{\Ran^{\on{untl}}}}="D";
(-5,30)*+{\Rep(\cG)_{\Ran^{\on{untl}},\on{indep}}}="E";
{\ar@{->}^{\Eis^{-,\on{spec}}} "B";"A"};
{\ar@{->}^{\Gamma^{\on{spec},\IndCoh}_\cG} "A";"E"};
{\ar@{->}_{\Gamma^{\on{spec},\IndCoh}_{\cM,\on{untl}}} "B";"D"};
{\ar@{->}_{\on{co}\!J^{-,\on{spec},!}} "D";"C"};
{\ar@{->}_{\on{emb.indep}_{\Rep(\cG)}^L} "C";"E"};
{\ar@{=>} "D";"A"};
\endxy}
\end{equation} 

From \thmref{t:Eis spec Gamma} we obtain:
\begin{cor} \label{c:Eis spec Gamma}
The natural transformation \eqref{e:Eis spec Gamma nat 1 IndCoh} is an isomorphism.
\end{cor}

\sssec{}

In order to prove \thmref{t:Eis spec Gamma} we will reformulate it in dual terms. Recall the functor
$$\Loc^{\on{spec}}_{\cG,\CZ}:\Rep(\cG)_\CZ\to \QCoh(\LS_\cG)\otimes \Dmod(\CZ),$$
given by pullback along $\on{ev}_\CZ$, see \cite[Sect. 17.6.1]{GLC2}. 

\medskip

Diagram \eqref{e:Eis spec Gamma diag} gives rise to a natural transformation
\begin{equation} \label{e:CT spec Loc nat}
\vcenter
{\xy
(0,0)*+{\QCoh(\LS_\cG)\otimes \Dmod(\CZ)}="A";
(60,0)*+{\QCoh(\LS_\cM)\otimes \Dmod(\CZ),}="B";
(0,20)*+{\Rep(\cG)_\CZ}="C";
(60,20)*+{\Rep(\cM)_\CZ}="D";
{\ar@{<-}^{\on{CT}^{-,\on{spec}}_\QCoh} "B";"A"};
{\ar@{<-}_{J^{-,\on{spec},!}} "D";"C"};
{\ar@{<-}^{\Loc^{\on{spec}}_{\cG,\CZ}} "A";"C"};
{\ar@{<-}_{\Loc^{\on{spec}}_{\cM,\CZ}} "B";"D"};
{\ar@{=>} "D";"A"};
\endxy}
\end{equation} 
which can be viewed as obtained by passing to dual functors in \eqref{e:Eis spec Gamma nat}. 

\sssec{}

We can view the system of natural transformations \eqref{e:CT spec Loc nat} as a map
\begin{equation} \label{e:CT spec Loc nat 1}
\ul{\Loc}^{\on{spec}}_\cM \circ \ul{J}^{-,\on{spec},!} \to 
(\on{CT}^{\on{spec}}_\QCoh\otimes \on{Id})\circ \ul{\Loc}^{\on{spec}}_\cG
\end{equation} 
as lax unital local-to-global functors
$$\ul\Rep(\cG)\rightrightarrows \QCoh(\LS_\cM)\otimes \ul{\Dmod}(\Ran),$$
with the right-hand sides being strictly unital.

\medskip

By \cite[Sect. 11.4.8]{GLC2}, the natural transformatiom \eqref{e:CT spec Loc nat 1} gives rise
to a natural transformation
\begin{equation} \label{e:CT spec Loc nat 2}
\left(\ul{\Loc}^{\on{spec}}_\cM \circ \ul{J}^{-,\on{spec},!} \right)^{\int\, \on{ins.unit}}\to 
(\on{CT}^{\on{spec}}_\QCoh\otimes \on{Id})\circ \ul{\Loc}^{\on{spec}}_\cG.
\end{equation} 

In other words, we obtain a system of natural transformations

\medskip

\begin{equation} \label{e:CT spec Loc nat 3}
\vcenter
{\xy
(0,0)*+{\QCoh(\LS_\cG)\otimes \Dmod(\CZ)}="A";
(60,0)*+{\QCoh(\LS_\cM)\otimes \Dmod(\CZ)}="B";
(120,0)*+{\QCoh(\LS_\cM)\otimes \Dmod(\CZ^\subseteq).}="F";
(60,20)*+{\Rep(\cG)_{\CZ^\subseteq}}="C";
(120,20)*+{\Rep(\cM)_{\CZ^\subseteq}}="D";
(0,20)*+{\Rep(\cG)_\CZ}="E";
{\ar@{->}_{\on{CT}^{-,\on{spec}}_\QCoh\otimes \on{Id}} "A";"B"};
{\ar@{->}_{\Loc^{\on{spec}}_{\cG,\CZ}} "E";"A"};
{\ar@{->}^{\Loc^{\on{spec}}_{\cM,\CZ^\subseteq}} "D";"F"};
{\ar@{->}^{J^{-,\on{spec},!}} "C";"D"};
{\ar@{->}^{\on{ins.unit}_\CZ} "E";"C"};
{\ar@{=>} "D";"A"};
{\ar@{->}_{\on{Id} \otimes (\on{pr}_{\on{small},\CZ})_!} "F";"B"};
\endxy}
\end{equation} 

\sssec{}

We will prove:

\begin{thm} \label{t:CT of Loc spec QCoh}
The natural transformation \eqref{e:CT spec Loc nat 2} is an isomorphism.
\end{thm} 

\sssec{}

We claim that \thmref{t:CT of Loc spec QCoh} formally implies \thmref{t:Eis spec Gamma}.
Indeed, this is a duality manipulation parallel to one the proof of \corref{c:coeff of Eis !}. 

\sssec{}

Recall (see \cite[Sect. 17.6.4]{GLC2}) that by a slight abuse of notation, we use the same symbol $\Loc^{\on{spec}}_{\cG,\CZ}$
to denote the composition:
$$\Rep(\cG)_\CZ \overset{\Loc^{-,\on{spec}}_{\cG,\CZ}}\longrightarrow \QCoh(\LS_\cG)\otimes \Dmod(\CZ)
\overset{\Xi_{\LS_\cG}\otimes \on{Id}}\longrightarrow \IndCoh(\LS_\cG)\otimes \Dmod(\CZ).$$

\sssec{}

From \thmref{t:CT of Loc spec QCoh} and Lemmas \ref{l:CT spec IndCoh vs QCoh} and \ref{l:rel det LS} we obtain:

\begin{cor} \label{c:CT of Loc spec IndCoh} 
The following diagram of functors commutes: 
$$
\CD 
\Rep(\cG)_\CZ @>{\on{ins.unit}_\CZ}>> \Rep(\cG)_{\CZ^\subseteq}  @>{J^{-,\on{spec},!}}>> \Rep(\cM)_{\CZ^\subseteq} \\
@V{\Loc^{-,\on{spec}}_{\cG,\CZ}}VV & & @VV{\Loc^{-,\on{spec}}_{\cM,\CZ^\subseteq}}V \\
\IndCoh(\LS_\cG)\otimes \Dmod(\CZ) @>>> \IndCoh(\LS_\cM)\otimes \Dmod(\CZ) 
@<{\on{Id} \otimes (\on{pr}_{\on{small},\CZ})_!}<< 
\IndCoh(\LS_\cM)\otimes \Dmod(\CZ^\subseteq),
\endCD
$$
where the left bottom horizontal arrow is the tensor product of $\on{Id}_{\Dmod(\CZ)}$ with 
$$\IndCoh(\LS_\cG)\overset{\on{CT}^{-,\on{spec}}}\longrightarrow \IndCoh(\LS_\cM)
\overset{(-)\otimes \CL_{2\rho_P(\omega_X)}[(2g-2)\cdot \dim(\cn_P)]}\longrightarrow  \IndCoh(\LS_\cM).$$
\end{cor} 

\begin{rem}

\corref{c:CT of Loc spec IndCoh} is the Langlands-dual counterpart of \corref{c:CT Poinc !}. Note that 
the functor 
$$(\on{transl}_{2\rho_P(\omega_X)})^*: \Dmod_{\frac{1}{2}}(\Bun_M)\to \Dmod_{\frac{1}{2}}(\Bun_M),$$
appearing in \corref{c:CT Poinc !}, 
is the Langlands-dual of 
$$\IndCoh(\LS_\cG) \overset{(-)\otimes \CL_{2\rho_P(\omega_X)}}\longrightarrow  \IndCoh(\LS_\cM).$$

\end{rem}

\ssec{Proof of \thmref{t:CT of Loc spec QCoh}}

\sssec{}

The natural transformation \eqref{e:CT spec Loc nat} is obtained by concatenating the commutative
diagram
$$
\CD
\Rep(\cG)_\CZ @>{\Res^\cG_{\cP^-}}>> \Rep(\cP^-)_\CZ \\
@V{\Loc^{\on{spec}}_{\cG,\CZ}}VV @VV{\Loc^{\on{spec}}_{\cP^-,\CZ}}V \\
\QCoh(\LS_\cG)\otimes \Dmod(\CZ) @>{(\sfp^{-,\on{spec}})^*\otimes \on{Id}}>> \QCoh(\LS_{\cP^-})\otimes \Dmod(\CZ)
\endCD
$$
with the natural transformation
\begin{equation} \label{e:CT spec Loc nat 4}
\xy 
(0,0)*+{\QCoh(\LS_{\cP^-})\otimes \Dmod(\CZ)}="A";
(70,0)*+{\QCoh(\LS_\cM)\otimes \Dmod(\CZ).}="B";
(0,20)*+{\Rep(\cP^-)_\CZ}="C";
(70,20)*+{\Rep(\cM)_\CZ}="D";
{\ar@{<-}^{(\sfq^{-,\on{spec}})_*\otimes \on{Id}} "B";"A"};
{\ar@{<-}_{\on{inv}_{\cN^-_P}} "D";"C"};
{\ar@{<-}^{\Loc^{\on{spec}}_{\cP^-,\CZ}} "A";"C"};
{\ar@{<-}_{\Loc^{\on{spec}}_{\cM,\CZ}} "B";"D"};
{\ar@{=>} "D";"A"};
\endxy
\end{equation} 

Hence, it is enough to show that the resulting map
\begin{equation} \label{e:CT spec Loc nat 5}
\left(\ul{\Loc}^{\on{spec}}_\cM \circ \on{inv}_{\cN^-_P} \right)^{\int\, \on{ins.unit}}\to 
((\sfq^{-,\on{spec}})_* \otimes \on{Id})\circ \ul{\Loc}^{\on{spec}}_{\cP^-}
\end{equation} 
is an isomorphism.

\sssec{}

To prove that \eqref{e:CT spec Loc nat 5} is an isomorphism, it suffices to show that it becomes
such after taking *-fibers at all field-valued points. Let $\sigma$ be a $k'$-valued point of $\LS_\cM$.

\medskip

Up to base changing the situation, we can assume that $k'=k$. 

\sssec{}

Let $\cN^-_{P,\sigma}$ be the twist of the constant group-scheme with a connection with fiber $\cN^-_P$
by $\sigma$ (via the adjoint action of $\cM$ on $\cN^-_P$). Note that
$$\LS_{\cN^-_{P,\sigma}}\simeq \on{pt}\underset{\sigma,\LS_\cM}\times \LS_{\cP^-}.$$

\medskip

For $\CZ\to \Ran$, the map
$$\on{ev}_\CZ:\LS_{\cN^-_{P,\sigma}}\times \CZ\to \LS^\reg_{\cN^-_{P,\sigma},\CZ}$$
gives rise to a natural transformation
\begin{equation}  \label{e:Loc N}
\vcenter
{\xy 
(0,0)*+{\QCoh(\LS_{\cN^-_{P,\sigma}})\otimes \Dmod(\CZ)}="A";
(80,0)*+{\Dmod(\CZ).}="B";
(0,20)*+{\Rep(\cN^-_{P,\sigma})_\CZ}="C";
(80,20)*+{\Dmod(\CZ)}="D";
{\ar@{<-}^{\Gamma(\LS_{\cN^-_{P,\sigma}},-)\otimes \on{Id}} "B";"A"};
{\ar@{<-}_{\on{inv}_{\cN^-_{P,\sigma}}} "D";"C"};
{\ar@{<-}^{\Loc^{\on{spec}}_{\cN^-_{P,\sigma},\CZ}} "A";"C"};
{\ar@{<-}_{\on{Id}} "B";"D"};
{\ar@{=>} "D";"A"};
\endxy}
\end{equation} 

\medskip

The concatenation of \eqref{e:CT spec Loc nat 4} with the functor of taking the *-fiber at $\sigma$
identifies with the concatenation of the commutative diagram
$$
\CD
\Rep(\cP^-)_\CZ  @>>> \Rep(\cN^-_{P,\sigma})_\CZ \\
@V{\Loc^{\on{spec}}_{\cP^-,\CZ}}VV @VV{\Loc^{\on{spec}}_{\cN^-_{P,\sigma},\CZ}}V \\
\QCoh(\LS_{\cP^-})\otimes \Dmod(\CZ) @>>> \QCoh(\LS_{\cN^-_{P,\sigma}})\otimes \Dmod(\CZ)
\endCD
$$
(where the horizontal arrow are given by *-restriction) with the natural transformation \eqref{e:Loc N}.

\medskip

Hence, it suffices to show that the map
\begin{equation}  \label{e:Loc N 1}
(\ul{\on{inv}}_{\cN^-_{P,\sigma}})^{\int\, \on{ins.unit}}\to 
(\Gamma(\LS_{\cN^-_{P,\sigma}},-)\otimes \on{Id})\circ \ul\Loc^{\on{spec}}_{\cN^-_{P,\sigma}}
\end{equation}
as functors
$$\ul\Rep(\cN^-_{P,\sigma})\rightrightarrows \ul\Dmod(\Ran)$$ 
is an isomorphism.

\sssec{}

Since $\cN^-_{P,\sigma}$ is unipotent, the category $\Rep(\cN^-_{P,\sigma})$ is generated by the trivial representation.
Hence, in order to prove that \eqref{e:Loc N 1} is an isomorphism, it is enough to show that the map
\begin{equation}  \label{e:Loc N 2}
\on{C}^{\on{fact}}_\cdot(X,\on{inv}_{\cN^-_{P,\sigma}}(k))\to \Gamma\left(\LS_{\cN^-_{P,\sigma}},\CO_{\LS_{\cN^-_{P,\sigma}}}\right),
\end{equation}
induced by \eqref{e:Loc N 1} is an isomorphism.

\sssec{} \label{sss:Omega LS as good}

Recall the construction
$$R\mapsto ``\Spec(R)"$$
from \cite[Sect. 12.7.6]{GLC2}, and reviewed in \secref{sss:Spec quot}. 

\medskip

Since $\cN^-_P$ is unipotent, the classifying stack
$\on{pt}/\cN^-_{P,\sigma}$ identifies with 
$$``\Spec_X(\on{inv}_{\cN^-_{P,\sigma}}(k))".$$

From here it formally follows that 
$$\LS_{\cN^-_{P,\sigma}}\simeq ``\Spec(\on{C}^{\on{fact}}_\cdot(X,\on{inv}_{\cN^-_{P,\sigma}}(k)))",$$
cf. \cite[last formula in Sect. 12.7.4]{GLC2}. 

\medskip

In terms of this identification, the map \eqref{e:Loc N 2} identifies with the map \cite[Equation (12.46)]{GLC2},
which is reproduced as \eqref{e:Spec R quot map abs}.

\medskip

Hence, in order to show that \eqref{e:Loc N 2} is an isomorphism, we have to show that 
$$\on{C}^{\on{fact}}_\cdot(X,\on{inv}_{\cN^-_{P,\sigma}}(k))\in \on{ComAlg}(\Vect)$$
is \emph{as good as connective} (see \cite[Sect. 12.7.6]{GLC2} or \secref{sss:Spec quot} for what this means). 

\medskip

However, this is done by the same argument as \propref{p:as good as conn}.

\ssec{Global (co-)enhanced categories on the spectral side}

\sssec{}

Let $\CZ$ be a prestack mapping to $\Ran$. We define the category $\IndCoh(\LS_\cM)^{-,\on{enh}}_\CZ$
in a way parallel to \secref{ss:geom global enh}. Namely, 

\medskip

$$\IndCoh(\LS_\cM)^{-,\on{enh}}_\CZ:=
\left(\IndCoh(\LS_\cM)\otimes \Dmod(\CZ)\right)\underset{\Sph^{\on{spec}}_{\cM,\CZ}}\otimes 
\on{I}(\cG,\cP^-)^{\on{spec,loc}}_\CZ.$$

\medskip

The assignment
\begin{equation} \label{e:enh spec glob}
\CZ\rightsquigarrow \IndCoh(\LS_\cM)^{-,\on{enh}}_\CZ
\end{equation} 
is a crystal of categories over $\Ran$, which we will denote by 
$$\ul{\IndCoh}(\LS_\cM)^{-,\on{enh}}.$$

\sssec{}

The monadic adjunction
$$\ind^{\on{spec}}_{\Sph\to \semiinf}:\Sph^{\on{spec}}_\cM\rightleftarrows \on{I}(\cG,\cP^-)^{\on{spec,loc}}:
\oblv^{\on{spec}}_{\semiinf\to \Sph}$$
gives rise to an adjunction
\begin{equation} \label{e:enh spec glob adj}
\ind^{\on{spec}}_{\on{enh}}:\IndCoh(\LS_\cM)\otimes \Dmod(\CZ)\rightleftarrows
\IndCoh(\LS_\cM)^{-,\on{enh}}_\CZ:\oblv^{\on{spec}}_{\on{enh}}.
\end{equation} 

\sssec{} \label{sss:enh spec glob monad}

Here is another way to think about the categories \eqref{e:enh spec glob}. 

\medskip

Let $\wt\Omega^{\on{spec}}$ be as in \secref{sss:wt Omega spec}. Consider
$$\wt\Omega^{\on{spec}}_\CZ\in \on{AssocAlg}(\Sph^{\on{spec}}_{\cM,\CZ}),$$
where we consider $\Sph^{\on{spec}}_{\cM,\CZ}$ as a monoidal category. 

\medskip

Consider $\IndCoh(\LS_\cM)\otimes \Dmod(\CZ)$ as a right module category over $\Sph^{\on{spec}}_{\cM,\CZ}$. 
Then 
$$\IndCoh(\LS_\cM)^{-,\on{enh}}_\CZ\simeq \wt\Omega^{\on{spec}}_\CZ\mod^r(\IndCoh(\LS_\cM)\otimes \Dmod(\CZ)),$$
so that \eqref{e:enh spec glob adj} is the adjunction
$$\ind_{\wt\Omega_\CZ}:
\IndCoh(\LS_\cM)\otimes \Dmod(\CZ)\rightleftarrows \wt\Omega^{\on{spec}}_\CZ\mod^r(\IndCoh(\LS_\cM)\otimes \Dmod(\CZ)):\oblv_{\wt\Omega_\CZ}.$$

\sssec{}

Consider the full category 
$$\IndCoh(\LS_\cM)^{-,\on{enh}}_{\Ran^{\on{untl}},\on{indep}}\subset \IndCoh(\LS_\cM)^{-,\on{enh}}_\Ran.$$

As in \eqref{e:IGP indep}, it fits into the pullback square
\begin{equation} \label{e:IGP spec indep}
\CD
\IndCoh(\LS_\cM)^{-,\on{enh}}_{\Ran^{\on{untl}},\on{indep}} @>>> \IndCoh(\LS_\cM)^{-,\on{enh}}_\Ran \\
@V{\oblv^{\on{spec}}_{\on{enh}}}VV @VV{\oblv^{\on{spec}}_{\on{enh}}}V \\
\IndCoh(\LS_\cM)\otimes \Dmod(\Ran^{\on{untl}})_{\on{indep}} @>>> \IndCoh(\LS_\cM)\otimes \Dmod(\Ran) \\
@V{\sim}VV \\
\IndCoh(\LS_\cM).
\endCD
\end{equation} 

\sssec{} \label{sss:spec co-enh}

In a way parallel to \secref{sss:geom co-enh}, we define 
$$\IndCoh(\LS_\cG)^{-,\on{enh}_{\on{co}}}_\CZ:=
\left(\IndCoh(\LS_\cG)\otimes \Dmod(\CZ)\right)\underset{\Sph^{\on{spec}}_{\cG,\CZ}}\otimes 
\on{I}(\cG,\cP^-)^{\on{spec,loc}}_{\on{co},\CZ}.$$

\medskip

We regard the assignment 
$$\CZ\rightsquigarrow \IndCoh(\LS_\cG)^{-,\on{enh}_{\on{co}}}_\CZ$$
as a crystal of categories over $\Ran$, to be denoted
$$\ul{\IndCoh}(\LS_\cG)^{-,\on{enh}_{\on{co}}}.$$
 
\ssec{The enhanced spectral Eisenstein functor} \label{ss:Eis spec enh} 

In this subsection we define the \emph{enhanced} spectral Eisenstein functor, which is a unital local-to-global functor
$$\ul\Eis^{-,\on{spec,enh}}:\ul{\IndCoh}(\LS_\cM)^{-,\on{enh}}\to \IndCoh(\LS_\cG)\otimes \ul\Dmod(\Ran).$$

\sssec{}

Let $\LS^{\on{mer,glob}}_{\cP^-,\Ran}$ be the space of meromorphic $\cP^-$-local systems on $X$
with poles at specified points, defined as in \cite[Sect. B.7.15]{GLC2}. 

\medskip

Denote
$$\LS^{\mf_\cG,\on{glob}}_{\cP^-,\Ran}:=(\LS_\cG\times \Ran)\underset{\LS^{\on{mer,glob}}_{\cG,\Ran}}\times \LS^{\on{mer,glob}}_{\cP^-,\Ran}.$$

\medskip

We will denote by $\sfp^{-,\on{spec}}_\Ran$ the map
$$\LS^{\mf_\cG,\on{glob}}_{\cP^-,\Ran}\to \LS_\cG\times \Ran.$$

\medskip

Given $\CZ\to \Ran$, we will substitute $\Ran\rightsquigarrow \CZ$ in the subscript to denote the corresponding base-changed objects. 

\sssec{}

Let $\sfs^{\on{spec},\cG}_\Ran$ denote the naturally defined map
\begin{equation} \label{e:tricky ev}
\LS^{\mf_\cG,\on{glob}}_{\cP^-,\Ran}\to 
\LS^\reg_{\cG,\Ran}\underset{\LS^\mer_{\cG,\Ran}}\times \LS^\mer_{\cP^-,\Ran}=:\LS^{\mf_\cG}_{\cP^-,\Ran}.
\end{equation} 

Note that we have a commutative (but not Cartesian) diagram 
$$
\CD
\LS^{\mf_\cG,\on{glob}}_{\cP^-,\Ran} @>{\sfp^{-,\on{spec}}_\Ran}>> \LS_\cG\times \Ran \\
@V{\sfs^{\on{spec},\cG}_\Ran}VV @VV{\on{ev}^\cG_\Ran}V \\
\LS^{\mf_\cG}_{\cP^-,\Ran} @>>{\sfp^{-,\on{spec,loc}}_\Ran}> \LS^\reg_{\cG,\Ran},
\endCD
$$
where 
$$\sfp^{-,\on{spec,loc}}_\Ran:\LS^{\mf_\cG}_{\cP^-,\Ran}\to \LS^\reg_{\cG,\Ran}$$
is the projection on the first factor.

\sssec{} \label{sss:LS P G mer good}

We have a naturally defined map
\begin{equation} \label{e:iota G glob} 
\iota^{\on{glob},\cG}:\LS_{\cP^-}\times \Ran \to \LS^{\mf_\cG,\on{glob}}_{\cP^-,\Ran},
\end{equation} 
so that the diagram 
$$
\CD
\LS_{\cP^-}\times \Ran @>{\iota^{\on{glob},\cG}}>> \LS^{\mf_\cG,\on{glob}}_{\cP^-,\Ran} \\
@V{\on{ev}^{\cP^-}_\Ran}VV @VV{\sfs^{\on{spec},\cG}_\Ran}V \\
\LS^\reg_{\cP^-,\Ran} @>>{\iota^\cG}> \LS^{\mf_\cG}_{\cP^-,\Ran}
\endCD
$$
is Cartesian. 

\medskip

The map $\iota^{\on{glob},\cG}$ is an isomorphism at the classical level. This observation, combined with
the calculation of cotangent spaces implies that $\LS^{\mf_\cG,\on{glob}}_{\cP^-,\Ran}$
is a (derived) algebraic stack locally almost of finite type (see \cite[Chapter 1, Theorem 9.1.2]{GaRo4}). 

\medskip

In particular, we obtain 
that the category $\IndCoh(\LS^{\mf_\cG,\on{glob}}_{\cP^-,\Ran})$
is well-defined. 

\sssec{} \label{sss:global IGP gen}

We claim that there is a well-defined functor
$$\IndCoh(\LS_\cM)^{-,\on{enh}}_\CZ\to \IndCoh(\LS^{\mf_\cG,\on{glob}}_{\cP^-,\CZ}),$$
to be denoted $(\sfq^{-,\on{spec}}_\CZ)^{*,\IndCoh,\on{enh}}$. 

\medskip

We will perform this construction in the present subsection when $\CZ=\ul{x}$, and the case of a general $\CZ$ will be considered in 
Sects. \ref{sss:global IGP gen-first}-\ref{sss:global IGP gen-last}. 

\sssec{}

Consider the map
$$\LS^{\mf_\cG}_{\cP^-,\ul{x}} \to \LS^\mer_{\cP^-,\ul{x}} \to \LS^\mer_{\cM,\ul{x}},$$
and the fiber product
\begin{equation} \label{e:LS M fp}
\LS^{\on{mer,glob}}_{\cM,\ul{x}}\underset{\LS^\mer_{\cM,\ul{x}}}\times \LS^{\mf_\cG}_{\cP^-,\ul{x}}.
\end{equation}

\medskip

We claim that when $\ul{x}$ is fixed, \eqref{e:LS M fp} is an algebraic stack locally almost of finite type.
This follows from the calculation of cotangent spaces combined with the fact that the map
\begin{equation} \label{e:LS M fp 0}
\LS_\cM\underset{\LS^\reg_{\cM,\ul{x}}}\times \LS^\reg_{\cP^-,\ul{x}}\to
\LS^{\on{mer,glob}}_{\cM,\ul{x}}\underset{\LS^\mer_{\cM,\ul{x}}}\times \LS^{\mf_\cG}_{\cP^-,\ul{x}}
\end{equation}
is an isomorphism at the classical level. 

\medskip

In particular, the category $\IndCoh(\LS^{\on{mer,glob}}_{\cM,\ul{x}}\underset{\LS^\mer_{\cM,\ul{x}}}\times \LS^{\mf_\cG}_{\cP^-,\ul{x}})$
is well-defined. Moreover, as in \eqref{e:amb vs IGP spec *}, we have a canonical isomorphism
\begin{multline} \label{e:LS M fp 1}
\IndCoh(\LS_\cM)^{-,\on{enh}}_{\ul{x}}:=
\IndCoh(\LS_\cM)\underset{\Sph^{\on{spec}}_{\cM,\ul{x}}}\otimes \on{I}(\cG,\cP^-)^{\on{spec,loc}}_{\ul{x}}\simeq \\
\simeq \IndCoh(\LS^{\on{mer,glob}}_{\cM,\ul{x}}\underset{\LS^\mer_{\cM,\ul{x}}}\times \LS^{\mf_\cG}_{\cP^-,\ul{x}}).
\end{multline}

\sssec{}

Combining the maps
$$\LS^{\mf_\cG,\on{glob}}_{\cP^-,\ul{x}}\to \LS^{\on{mer,glob}}_{\cP^-,\ul{x}}\to \LS^{\on{mer,glob}}_{\cM,\ul{x}}$$
and $\sfs^{\on{spec},\cG}_{\ul{x}}$, we obtain a map
$$\LS^{\mf_\cG,\on{glob}}_{\cP^-,\ul{x}}\to \LS^{\on{mer,glob}}_{\cM,\ul{x}}\underset{\LS^\mer_{\cM,\ul{x}}}\times \LS^{\mf_\cG}_{\cP^-,\ul{x}},$$
which we denote by $\sfq^{-,\on{spec,enh}}_{\ul{x}}$. 

\medskip

We claim that $\sfq_{\ul{x}}^{-,\on{spec,enh}}$ is schematic of finite Tor-dimension. Indeed, its base change by the nil-isomorphism 
\eqref{e:LS M fp 0} is the map
$$\LS_{\cP^-}\overset{\sfq^{-,\on{spec}}_{\ul{x}}\times \on{ev}^{\cP^-}_{\ul{x}}}
\longrightarrow \LS_\cM\underset{\LS^\reg_{\cM,\ul{x}}}\times \LS^\reg_{\cP^-,\ul{x}},$$
which has the required property. 

\sssec{}

Hence, we obtain a well-defined functor
$$(\sfq^{-,\on{spec,enh}}_{\ul{x}})^{\IndCoh,*}:\IndCoh(\LS^{\on{mer,glob}}_{\cM,\ul{x}}\underset{ \LS^\mer_{\cM,\ul{x}}}\times \LS^{\mf_\cG}_{\cP^-,\ul{x}})\to
\IndCoh(\LS^{\mf_\cG,\on{glob}}_{\cP^-,\ul{x}}).$$

\medskip

Combining with \eqref{e:LS M fp 1}, we obtain the sought-for functor
$$(\sfq^{-,\on{spec}}_{\ul{x}})^{*,\IndCoh,\on{enh}}:\IndCoh(\LS_\cM)^{-,\on{enh}}_{\ul{x}}\to
\IndCoh(\LS^{\mf_\cG,\on{glob}}_{\cP^-,\ul{x}}).$$

\begin{rem} \label{r:old pullback}

Here is another way to think about the composition
\begin{multline} \label{e:LS M fp 0.5}
\IndCoh(\LS_\cM)\otimes \on{I}(\cG,\cP^-)^{\on{spec,loc}}_{\ul{x}}\to \IndCoh(\LS_\cM)^{-,\on{enh}}_{\ul{x}}\to \\
\overset{(\sfq^{-,\on{spec}}_{\ul{x}})^{*,\IndCoh,\on{enh}}}\longrightarrow \IndCoh(\LS^{\mf_\cG,\on{glob}}_{\cP^-,\ul{x}}).
\end{multline} 

\medskip

Denote:
$$\on{Hecke}^{\on{spec,glob}}_{\cG,\cP^-,\ul{x}}:=
\LS_\cG\underset{\LS^{\on{mer,glob}}_{\cG,\ul{x}}}\times \LS^{\on{mer,glob}}_{\cP^-,\ul{x}}
\underset{\LS^{\on{mer,glob}}_{\cM,\ul{x}}}\times\LS_\cM.$$

This is an algebraic stack, locally almost of finite type, equipped with maps
$$\LS_\cG \overset{\sfp^{-,\on{spec,mer}}_{\ul{x}}}\longleftarrow \on{Hecke}^{\on{spec,glob}}_{\cG,\cP^-,\ul{x}}
\overset{\sfq^{-,\on{spec,mer}}_{\ul{x}}}\longrightarrow  \LS_\cM.$$
and 
$$\sfs^{\on{spec}}_{\ul{x}}:\on{Hecke}^{\on{spec,glob}}_{\cG,\cP^-,\ul{x}}\to \on{Hecke}^{\on{spec,loc}}_{\cG,\cP^-,\ul{x}},$$
where $\on{Hecke}^{\on{spec,loc}}_{\cG,\cP^-}$ is as in \secref{sss:Hecke spec GP}. 

\medskip

We claim that the functor
\begin{equation} \label{e:old pullback}
(\sfq^{-,\on{spec,mer}}_{\ul{x}})^*(-)\overset{*}\otimes (\sfs^{\on{spec}}_{\ul{x}})^*(-),
\end{equation}
which maps 
$$\IndCoh(\LS_\cM)\otimes \on{I}(\cG,\cP^-)^{\on{spec,loc}}_{\ul{x}}=\IndCoh(\LS_\cM)\otimes \IndCoh(\on{Hecke}^{\on{spec,loc}}_{\cG,\cP^-,\ul{x}})
\to \IndCoh(\on{Hecke}^{\on{spec,glob}}_{\cG,\cP^-,\ul{x}}),$$
is well-defined. Indeed, this follows from the fact that the map
$$\on{Hecke}^{\on{spec,glob}}_{\cG,\cP^-,\ul{x}} \overset{\sfq^{-,\on{spec,mer}}_{\ul{x}}\times \sfs^{\on{spec}}_{\ul{x}}}\longrightarrow 
\LS_\cM\times \on{Hecke}^{\on{spec,loc}}_{\cG,\cP^-,\ul{x}}$$
is schematic of finite Tor-dimension. (The former follows by deformation theory using, e.g., \cite[Theorem B.2.14]{AG}; for the latter, 
we observe that the base change of the above morphism 
with respect to the nil-embedding 
$$\LS^\reg_{\cP^-,\ul{x}}\to \on{Hecke}^{\on{spec,loc}}_{\cG,\cP^-,\ul{x}}$$
is the morphism
$$\LS_{\cP^-} \overset{\sfq^{-,\on{spec}}\times \on{ev}_{\ul{x}}}\longrightarrow \LS_\cM\times \LS^\reg_{\cP^-,\ul{x}},$$
which is of finite Tor-dimension.) 

\medskip

Now, the functor \eqref{e:LS M fp 0.5} is isomorphic to the composition of \eqref{e:old pullback} with the
$(\IndCoh,*)$-pushforward functor along
$$\on{Hecke}^{\on{spec,glob}}_{\cG,\cP^-,\ul{x}}\to \LS^{\mf_\cG,\on{glob}}_{\cP^-,\ul{x}}.$$

\end{rem} 

\sssec{}

The spectral Hecke groupoid $\on{Hecke}^{\on{spec,loc}}_{\cG,\ul{x}}$  acts naturally on 
$$\LS^{\on{mer,glob}}_{\cM,\ul{x}}\underset{\LS^\mer_{\cM,\ul{x}}}\times \LS^{\mf_\cG}_{\cP^-,\ul{x}} \text{ and }
\LS^{\mf_\cG,\on{glob}}_{\cP^-,\ul{x}},$$
and the map $\sfq^{-,\on{spec,enh}}_{\ul{x}}$ is compatible with this action. 

\medskip

From here we obtain that the functor $(\sfq^{-,\on{spec}}_{\ul{x}})^{*,\IndCoh,\on{enh}}$ is compatible with the
actions of $\Sph^{\on{spec}}_{\cG,\ul{x}}$ on the two sides. The same will be true when $\ul{x}\in \Ran$
is replaced by an arbitrary $\CZ\to \Ran$. 

\sssec{}

We define the functor
$$\Eis^{-,\on{spec,enh}}_\CZ:\IndCoh(\LS_\cM)^{-,\on{enh}}_\CZ\to \IndCoh(\LS_\cG)\otimes \Dmod(\CZ)$$
as the composition
$$(\sfp^{-,\on{spec}}_\CZ)^\IndCoh_*\circ (\sfq^{-,\on{spec}}_\CZ)^{*,\IndCoh,\on{enh}}.$$

\medskip

Since both functors $(\sfq^{-,\on{spec}}_\CZ)^{*,\IndCoh,\on{enh}}$ and $(\sfp^{-,\on{spec}}_\CZ)^\IndCoh_*$ 
respect the actions of $\Sph^{\on{spec}}_{\cG,\CZ}$, we obtain that $\Eis^{-,\on{spec,enh}}_\CZ$ also 
respects the actions of  $\Sph^{\on{spec}}_{\cG,\CZ}$.

\sssec{} 

Unwinding the construction, we obtain that the composition
$$\IndCoh(\LS_\cM)\otimes \Dmod(\CZ)\overset{\ind^{\on{spec}}_{\on{enh}}}\longrightarrow 
\IndCoh(\LS_\cM)^{-,\on{enh}}_\CZ \overset{\Eis^{-,\on{spec,enh}}_\CZ}\longrightarrow \IndCoh(\LS_\cG)\otimes \Dmod(\CZ)$$
identifies canonically with $\Eis^{-,\on{spec}}\otimes \on{Id}$. 

\sssec{}

We will view the assignment
$$\CZ\rightsquigarrow \Eis^{-,\on{spec,enh}}_\CZ$$
as a local-to-global functor 
$$\ul{\IndCoh}(\LS_\cM)^{-,\on{enh}}\to \IndCoh(\LS_\cG)\otimes \ul\Dmod(\Ran),$$
to be denoted $\ul\Eis^{-,\on{spec,enh}}$.

\medskip

It has a natural lax unital structure. However, the isomorphism
$$\Eis^{-,\on{spec,enh}}_\CZ\circ \ind_{\on{enh}}\simeq \Eis^{-,\on{spec}}\otimes \on{Id}$$
implies that this lax unital structure is strict. 

\sssec{} \label{sss:fact alg enh spec}

Let $\CA$ be a factorization algebra in $\on{I}(\cG,\cP^-)^{\on{spec,loc}}$. Parallel to \secref{sss:fact alg Eis} we can consider the functor
$$\Eis^{-,\on{spec}}_\CA:\IndCoh(\LS_\cM)\to \IndCoh(\LS_\cG).$$

\ssec{Partial enhancement}

This subsection serves as a preparation to \secref{ss:comp Eis spec}. 

\sssec{}

Set
$$\IndCoh(\LS_\cM)^{-,\on{part.enh}}_\CZ:=
\left(\IndCoh(\LS_\cM)\otimes \Dmod(\CZ)\right)\underset{\Rep(\cM)_\CZ}\otimes \Rep(\cP^-)_\CZ.$$

The monadic adjunction
$$\Res^\cM_{\cP^-}:\Rep(\cM)\rightleftarrows \Rep(\cP^-):\on{inv}_{\cN^-_P}$$
gives rise to a monadic adjunction
$$\ind^{\on{spec}}_{\on{part.enh}}:\IndCoh(\LS_\cM)\otimes \Dmod(\CZ)\rightleftarrows
\IndCoh(\LS_\cM)^{-,\on{part.enh}}_\CZ:\oblv^{\on{spec}}_{\on{part.enh}}.$$

\medskip 

We will consider the assignment
$$\CZ\rightsquigarrow \IndCoh(\LS_\cM)^{-,\on{part.enh}}_\CZ$$
as a crystal of categories over $\Ran$, to be denoted $$\ul{\IndCoh}(\LS_\cM)^{-,\on{part.enh}}.$$

\sssec{}

Consider the full category 
$$\IndCoh(\LS_\cM)^{-,\on{part.enh}}_{\Ran^{\on{untl}},\on{indep}}\subset \IndCoh(\LS_\cM)^{-,\on{part.enh}}_\Ran.$$

As in \eqref{e:IGP spec indep}, it fits into the pullback square
\begin{equation} \label{e:IGP part spec indep}
\CD
\IndCoh(\LS_\cM)^{-,\on{part.enh}}_{\Ran^{\on{untl}},\on{indep}} @>>> \IndCoh(\LS_\cM)^{-,\on{part.enh}}_\Ran \\
@VVV @VV{\oblv^{\on{spec}}_{\on{part.enh}}}V \\
\IndCoh(\LS_\cM)\otimes \Dmod(\Ran^{\on{untl}})_{\on{indep}} @>>> \IndCoh(\LS_\cM)\otimes \Dmod(\Ran) \\
@V{\sim}VV \\
\IndCoh(\LS_\cM).
\endCD
\end{equation} 

\sssec{}

Consider the category 
$$\IndCoh(\LS_\cM)\underset{\QCoh(\LS_\cM)}\otimes \QCoh(\LS_{\cP^-}).$$

The operation 
$$(-)\otimes (\on{ev}^{\cP^-}_\CZ)^*(-)$$
gives rise to a functor
$$\IndCoh(\LS_\cM)^{-,\on{part.enh}}_\CZ\to \left(\IndCoh(\LS_\cM)\underset{\QCoh(\LS_\cM)}\otimes \QCoh(\LS_{\cP^-})\right)\otimes \Dmod(\CZ),$$
which we denote by $(\sfq^{-,\on{spec}}_\CZ)^{*,\on{part.enh}}$. 

\sssec{}

Note now that the functor
$$(\sfq^{-,\on{spec}})^{\IndCoh,*}:\IndCoh(\LS_\cM)\to \IndCoh(\LS_{\cP^-})$$
gives rise to a functor
\begin{equation} \label{e:IndCoh tensor up}
(\sfq^{-,\on{spec}})^{\IndCoh,*}\otimes \on{Id}:
\IndCoh(\LS_\cM)\underset{\QCoh(\LS_\cM)}\otimes \QCoh(\LS_{\cP^-})\to  \IndCoh(\LS_{\cP^-}).
\end{equation} 

\sssec{}

Consider the composition 
$$(\sfq^{-,\on{spec}}_\CZ)^{\IndCoh,*,\on{part.enh}}:=((\sfq^{-,\on{spec}})^{\IndCoh,*}\otimes \on{Id})\circ (\sfq^{-,\on{spec}}_\CZ)^{*,\on{part.enh}}.$$

\medskip

Composing further with $(\sfp^{-,\on{spec}})^\IndCoh_*$, we obtain a functor
$$(\sfp^{-,\on{spec}})^\IndCoh_*\circ (\sfq^{-,\on{spec}}_\CZ)^{\IndCoh,*,\on{part.enh}}:\IndCoh(\LS_\cM)^{-,\on{part.enh}}_\CZ\to 
\IndCoh(\LS_\cG)\otimes \Dmod(\CZ),$$
which we denote by $\Eis^{-,\on{spec,part.enh}}_\CZ$.

\medskip

Note that the composition $\Eis^{-,\on{spec,part.enh}}_\CZ\circ \ind^{\on{spec}}_{\on{part.enh}}$ identifies 
canonically with $\Eis^{-,\on{spec}}\otimes \on{Id}$. 

\sssec{}

We will view the assignment
$$\CZ\rightsquigarrow \Eis^{-,\on{spec,part.enh}}_\CZ$$
as a local-to-global functor 
$$\ul{\IndCoh}(\LS_\cM)^{-,\on{part.enh}}\to \IndCoh(\LS_\cG)\otimes \ul\Dmod(\Ran),$$
to be denoted $\ul\Eis^{-,\on{spec,part.enh}}$.

\medskip

It has a natural lax unital structure. However, the isomorphism
$$\Eis^{-,\on{spec,part.enh}}_\CZ\circ \ind_{\on{enh}}\simeq \Eis^{-,\on{spec}}\otimes \on{Id}$$
implies that this lax unital structure is strict. 

\sssec{} \label{sss:fact alg enh spec partial}

Let $\CA$ be a factorization algebra in $\Rep(\cP^-)$. Parallel to \secref{sss:fact alg Eis} we can consider the functor
$$\Eis^{-,\on{spec,enh}}_\CA:\IndCoh(\LS_\cM)\to \IndCoh(\LS_\cG).$$

\sssec{}

Recall now that we have a pair of mutually adjoint factorization functors
$$\iota_*:\Rep(\cP^-)\rightleftarrows \on{I}(\cG,\cP^-)^{\on{spec,loc}}:\iota^!.$$

It induces a pair of adjoint functors
$$(\on{Id}\otimes \iota_*):
\IndCoh(\LS_\cM)^{-,\on{part.enh}}_\CZ\leftrightarrows \IndCoh(\LS_\cM)^{-,\on{enh}}_\CZ:(\on{Id}\otimes \iota^!).$$

Unwinding the constructions, we obtain that there is a canonical isomorphism
$$(\sfq^{-,\on{spec}}_\CZ)^{\IndCoh,*,\on{part.enh}}\simeq (\sfq^{-,\on{spec}}_\CZ)^{*,\IndCoh,\on{enh}}\circ (\on{Id}\otimes \iota_*).$$

From here, we obtain an isomorphism
\begin{equation} \label{e:Eis part ou pas}
\ul\Eis^{-,\on{spec,part.enh}}\simeq \ul\Eis^{-,\on{spec,enh}}\circ (\on{Id}\otimes \iota_*).
\end{equation} 

\ssec{Partial enhancement, continued}

The contents of this subsection are optional, in that they will not be necessary in the sequel. 

\sssec{}

Consider the assignment
$$\CZ\rightsquigarrow (\sfq^{-,\on{spec}}_\CZ)^{*,\on{part.enh}}$$
as a local-to-global functor, to be denoted $(\ul\sfq^{-,\on{spec}})^{*,\on{part.enh}}$.

\medskip

It has a natural unital structure, and hence gives rise to a functor
\begin{equation} \label{e:q* indep}
\IndCoh(\LS_\cM)^{-,\on{part.enh}}_{\Ran^{\on{untl}},\on{indep}} \to 
\IndCoh(\LS_\cM)\underset{\QCoh(\LS_\cM)}\otimes \QCoh(\LS_{\cP^-}).
\end{equation}

\sssec{}

We claim:

\begin{prop} \label{p:q* indep}
The functor \eqref{e:q* indep} is an equivalence.
\end{prop} 

\begin{proof}

Denote
$$\QCoh(\LS_\cM)^{-,\on{part.enh}}_\CZ:=
\left(\QCoh(\LS_\cM)\otimes \Dmod(\CZ)\right)\underset{\Rep(\cM)_\CZ}\otimes \Rep(\cP^-)_\CZ.$$

The operation 
$$(-)\otimes (\on{ev}^{\cP^-}_\CZ)^*(-)$$
gives rise to a functor
$$\QCoh(\LS_\cM)^{-,\on{part.enh}}_\CZ\to \QCoh(\LS_{\cP^-})\otimes \Dmod(\CZ),$$
which we denote by $\sfq^{-,\on{spec}}\otimes \Loc^{\on{spec}}_{\cP^-,\CZ}$. 

\medskip

The assignment
$$\CZ\rightsquigarrow \sfq^{-,\on{spec}}\otimes \Loc^{\on{spec}}_{\cP^-,\CZ}$$
is a unital local to global functor, to be denoted
\begin{equation} \label{e:q^* loc-to-glob}
\sfq^{-,\on{spec}}\otimes \ul\Loc^{\on{spec}}_{\cP^-}:\ul\QCoh(\LS_\cM)^{-,\on{part.enh}}\to \QCoh(\LS_{\cP^-})\otimes \ul\Dmod(\Ran).
\end{equation}

\medskip

It suffices to show that the resulting functor
\begin{equation} \label{e:q* indep 1}
\QCoh(\LS_\cM)^{-,\on{part.enh}}_{\Ran^{\on{untl}},\on{indep}} \to \QCoh(\LS_{\cP^-})
\end{equation}
is an equivalence. 

\medskip

To simplify the exposition, we will show that the functor \eqref{e:q* indep 1} becomes 
an equivalence after applying the operation
$$\Vect\underset{\QCoh(\LS_\cM)}\otimes -,$$
where $\QCoh(\LS_\cM)\to \Vect$ corresponds to the trivial local system.

\medskip

The resulting functor is the functor 
\begin{equation} \label{e:q* indep 2}
\Rep(\cN^-_P)_{\Ran^{\on{untl}},\on{indep}}\to \QCoh(\LS_{\cN^-_P}),
\end{equation}
induced by
$$\ul\Loc^{\on{spec}}_{\cN^-_P}:\ul\Rep(\cN^-_P)\to \QCoh(\LS_{\cN^-_P})\otimes \ul\Dmod(\Ran).$$

We rewrite 
$$\Rep(\cN^-_P)\simeq \CA\mod^{\on{com}},$$
where $\CA=\on{inv}_{\cN^-_P}(k)$. 

\medskip

By \secref{sss:Omega LS as good}, we have
$$\QCoh(\LS_{\cN^-_P})\simeq \on{C}^{\on{fact}}(X,\CA)\mod.$$

\medskip

With respect to these identifications, the functor \eqref{e:q* indep 2} is the functor
\begin{equation} \label{e:q* indep 3}
(\CA\mod^{\on{com}})_{\Ran^{\on{untl}},\on{indep}}\to \on{C}^{\on{fact}}(X,\CA)\mod
\end{equation}
of \cite[Equation (15.32)]{GLC2}.

\medskip

The fact that \eqref{e:q* indep 3} is an equivalence is the content of \cite[Remark 15.6.16]{GLC2}.

\end{proof}

\sssec{}

Next we observe:

\begin{lem} \label{l:IndCoh tensor up}
The functor \eqref{e:IndCoh tensor up} defines an equivalence
$$\IndCoh(\LS_\cM)\underset{\QCoh(\LS_\cM)}\otimes \QCoh(\LS_{\cP^-})\simeq  \IndCoh_{\cM}(\LS_{\cP^-}),$$
where $$\IndCoh_{\cM}(\LS_{\cP^-})\subset \IndCoh(\LS_{\cP^-})$$ is the full subcategory, consisting of objects
with singular support contained in 
$$\Sing(\LS_\cM)\underset{\LS_\cM}\times \LS_{\cP^-}\subset \Sing(\LS_{\cP^-}).$$
\end{lem}

\begin{proof}

This is a particular case of \cite[Corollary 7.6.2]{AG}. 

\end{proof}

\sssec{}

Combining \propref{p:q* indep} and \lemref{l:IndCoh tensor up}, we obtain:

\begin{cor} \label{c:IndCoh tensor up}
The functor
$$\IndCoh(\LS_\cM)^{-,\on{part.enh}}_{\Ran^{\on{untl}},\on{indep}} \to \IndCoh(\LS_{\cP^-}),$$
defined by the (unital) local-to-global functor 
$$\CZ\rightsquigarrow (\sfq^{-,\on{spec}}_\CZ)^{\IndCoh,*,\on{part.enh}},$$
is an equivalence onto $\IndCoh_{\cM}(\LS_{\cP^-})$.
\end{cor}  

\sssec{}

To summarize, we obtain that the datum of the functor $\ul\Eis^{-,\on{spec,part.enh}}$ is equivalent to the datum of
the functor
$$\IndCoh_{\cM}(\LS_{\cP^-})\hookrightarrow \IndCoh(\LS_{\cP^-})\overset{(\sfp^{-,\on{spec}})^\IndCoh_*}\longrightarrow 
\IndCoh(\LS_\cG).$$

Next, we will describe a similar interpretation of the functor $\ul\Eis^{-,\on{spec,enh}}$. 

\sssec{}

Consider the fiber product 
$$(\LS_{\cP^-})_\dr\underset{(\LS_{\cG})_\dr}\times \LS_{\cG}.$$ 

It comes equipped with a natural map
$$\LS_{\cP^-}\to (\LS_{\cP^-})_\dr\underset{(\LS_{\cG})_\dr}\times \LS_{\cG}.$$ 

\medskip

Let 
$$\IndCoh_\cM\left((\LS_{\cP^-})_\dr\underset{(\LS_{\cG})_\dr}\times \LS_{\cG}\right)\subset 
\IndCoh\left((\LS_{\cP^-})_\dr\underset{(\LS_{\cG})_\dr}\times \LS_{\cG}\right)$$
be the full subcategory consisting of objects, whose image under the !-pullback functor
\begin{equation} \label{e:Nick}
\IndCoh\left((\LS_{\cP^-})_\dr\underset{(\LS_{\cG})_\dr}\times \LS_{\cG}\right)\to \IndCoh(\LS_{\cP^-})
\end{equation} 
is contained in 
$$\IndCoh_{\cM}(\LS_{\cP^-})\subset \IndCoh(\LS_{\cP^-}).$$

\sssec{}

The following theorem can be deduced from the results of \cite{Ro}: 

\begin{thm} \label{t:Nick}
The category $\IndCoh(\LS_\cM)^{-,\on{enh}}_{\Ran^{\on{untl}},\on{indep}}$
identifies canonically with 
$$\on{I}(\cG,\cP^-)^{\on{spec,glob}}:=\IndCoh_\cM\left((\LS_{\cP^-})_\dr\underset{(\LS_{\cG})_\dr}\times \LS_{\cG}\right).$$

Under this identification, the forgetful functor
$$(\on{Id}\otimes \iota^!):\IndCoh(\LS_\cM)^{-,\on{enh}}_{\Ran^{\on{untl}},\on{indep}}\to \IndCoh(\LS_\cM)^{-,\on{part.enh}}_{\Ran^{\on{untl}},\on{indep}}$$
corresponds to the functor
$$\IndCoh_\cM\left((\LS_{\cP^-})_\dr\underset{(\LS_{\cG})_\dr}\times \LS_{\cG}\right)\to \IndCoh_{\cM}(\LS_{\cP^-}),$$
induced by \eqref{e:Nick}, and the functor 
$$\IndCoh(\LS_\cM)^{-,\on{enh}}_{\Ran^{\on{untl}},\on{indep}}\to \IndCoh(\LS_\cG),$$
induced by $\ul\Eis^{-,\on{spec,enh}}$ identifies with 
$$\IndCoh_\cM\left((\LS_{\cP^-})_\dr\underset{(\LS_{\cG})_\dr}\times \LS_{\cG}\right)\hookrightarrow 
\IndCoh\left((\LS_{\cP^-})_\dr\underset{(\LS_{\cG})_\dr}\times \LS_{\cG}\right)\to  \IndCoh(\LS_\cG),$$
where the second arrow is the $(\IndCoh,*)$-pushforward functor. 

\end{thm}

\ssec{The ``compactified" spectral Eisenstein functor} \label{ss:comp Eis spec}

In this subsection we will perform a calculation (formulated in \propref{p:Eis !* spec}) that will be used
in the proof of the main theorem of this paper, \thmref{t:main}. 

\sssec{}

Consider the construction from \secref{sss:fact alg enh spec} applied to 
$$\CA:=\IC^{-,\on{spec},\semiinf}.$$

Denote the resulting functor $\IndCoh(\LS_\cM)\to \IndCoh(\LS_\cG)$ by $\Eis^{-,\on{spec}}_{\IC}$. 

\sssec{} \label{sss:defn Eis spec IC}

Consider now the functor 
$$\Eis^{-,\on{spec}}_{!*}:\IndCoh(\LS_\cM)\to \IndCoh(\LS_\cG),$$
equal to the composition 
\begin{multline*}
\IndCoh(\LS_\cM)=\IndCoh(\LS_\cM)\underset{\QCoh(\LS_\cM)}\otimes \QCoh(\LS_\cM)\to 
\IndCoh(\LS_\cM)\underset{\QCoh(\LS_\cM)}\otimes \QCoh(\LS_{\cP^-}) \to \\
\overset{(\sfq^{-,\on{spec}})^{\IndCoh,*}\otimes \on{Id}}\longrightarrow 
\IndCoh(\LS_{\cP^-})\overset{(\sfp^{-,\on{spec}})^\IndCoh_*}\longrightarrow \IndCoh(\LS_\cG),
\end{multline*} 
where the second arrow is the tensor product of $\on{Id}_{\IndCoh(\LS_\cM)}$ with the functor
\begin{equation} \label{e:LS M to P}
\QCoh(\LS_\cM)\to \QCoh(\LS_{\cP^-}),
\end{equation} 
given by \emph{direct} image along the canonical map $\LS_\cM\to \LS_{\cP^-}$. 

\sssec{}

We are going to prove:

\begin{prop} \label{p:Eis !* spec}
The functors $\Eis^{-,\on{spec}}_{\IC}$ and $\Eis^{-,\on{spec}}_{!*}$ are canonically isomorphic. 
\end{prop}

The rest of this subsection is devoted to the proof of \propref{p:Eis !* spec}. 

\sssec{} \label{sss:Eis !* spec 1}

Note that using the functor $\sfq^{-,\on{spec}}\otimes \ul\Loc^{\on{spec}}_{\cP^-}$ of \eqref{e:q^* loc-to-glob}, we can associate to
$$\CB\in \on{FactAlg}(\Rep(\cP^-))$$ a $\QCoh(\LS_\cM)$-linear functor
$$(\sfq^{-,\on{spec}}\otimes \Loc^{\on{spec}}_{\cP^-})_\CB:\QCoh(\LS_\cM)\to \QCoh(\LS_{\cP^-}),$$
where $\QCoh(\LS_{\cP^-})$ is regarded as a $\QCoh(\LS_\cM)$-module category via $(\sfq^{-,\on{spec}})^*$. 

\medskip

Take $\CB$ to be the direct image of $\one_{\Rep(\cM)}$ under the factorization functor
$$\Rep(\cM)=\QCoh(\on{pt}/\cM)\to \QCoh(\on{pt}/\cP^-)=\Rep(\cP^-),$$
where the arrow is the functor of direct image along 
$$\on{pt}/\cM\to \on{pt}/\cP^-.$$

\medskip

By \eqref{e:Eis part ou pas}, it suffices to show that the resulting functor $(\sfq^{-,\on{spec}}\otimes \Loc^{\on{spec}}_{\cP^-})_\CB$
identifies with \eqref{e:LS M to P}. 

\sssec{}

To simplify the exposition, we will perform the corresponding calculation after applying the operation
$$\Vect\underset{\QCoh(\LS_\cM)}\otimes -,$$
where $\QCoh(\LS_\cM)\to \Vect$ corresponds to the trivial local system.

\sssec{}

Thus, we need to show that the functor $\Vect\to \QCoh(\LS_{\cN^-_P})$ corresponding to the local-to-global functor
$$\ul\Loc^{\on{spec}}_{\cN^-_P}:\ul\Rep(\cN^-_P)\to \QCoh(\LS_{\cN^-_P})\otimes \ul\Dmod(\Ran)$$
and 
$$\CB\in \on{FactAlg}(\Rep(\cN^-_P))$$
equal to the regular representation $R_{\N^-_P}$, identifies canonically with the functor of 
direct image along
\begin{equation} \label{e:triv loc sys N}
\on{pt}\to \LS_{\cN^-_P},
\end{equation}
corresponding to the trivial $\cN^-_P$-local system. 

\sssec{}

The latter is a standard calculation, valid for any algebraic group $\sH$, cf. \cite[Theorem 12.6.3]{AGKRRV}.

\medskip

Let us, however, give a different argument, specific to the unipotent case.

\sssec{} \label{sss:Eis !* spec end} 

Recall that according to \secref{sss:Omega LS as good}, we have
$$\QCoh(\LS_{\cN^-_P})=\on{C}^{\on{fact}}(X,\CA)\mod,$$
where $\CA:=\on{inv}_{\cN^-_P}(k)$. 

\medskip

Under the equivalence
$$\Rep(\cN^-_P)\simeq \CA\mod^{\on{fact}},$$
the object 
$$\CB\in \on{ComAlg}(\on{FactAlg}(\Rep(\cN^-_P)))$$
corresponds to the augmentation
\begin{equation} \label{e:augm loc}
\CA\to k.
\end{equation} 

\medskip

Therefore, the functor
$$\Loc^{\on{spec}}_{\cN^-_P,\CB}:\Vect\to \QCoh(\LS_{\cN^-_P})$$
corresponds to the homomorphism
\begin{equation} \label{e:augm glob}
\on{C}^{\on{fact}}(X,\CA)\overset{\text{\eqref{e:augm loc}}}\longrightarrow \on{C}^{\on{fact}}(X,k)\simeq k,
\end{equation} 
i.e., to the natural augmentation on $\on{C}^{\on{fact}}(X,\CA)$. 

\medskip

However, the homomorphism \eqref{e:augm glob} corresponds under
$$\LS_{\cN^-_P}\simeq ``\Spec(\on{C}^{\on{fact}}(X,\CA))"$$
to the map \eqref{e:triv loc sys N}. 

\qed[\propref{p:Eis !* spec}]

\ssec{Spectral Eisenstein functor and localization: the enhanced version}

\sssec{}

In this section we will provide an enhancement of \corref{c:Eis spec Gamma}. Namely, we will express 
the composition
$$\Gamma_\cG^{\on{spec},\IndCoh}\circ \Eis^{-,\on{spec,enh}},$$
via $\Gamma_\cM^{\on{spec},\IndCoh}$ and a functor of local nature.

\medskip

This calculation is parallel to one in \secref{s:Eis and Loc}. 

\sssec{}

Recall the factorization functor 
$$J^{-,\on{spec,pre-enh}}:\Rep(\cG)\otimes \on{I}(\cG,\cP^-)^{\on{spec,loc}}\to \Rep(\cM)$$
of \eqref{e:gen Jacquet functor spec}. 

\medskip

Let us denote by
$$\on{co}\!J^{-,\on{spec,pre-enh}}:\Rep(\cM)\otimes \on{I}(\cG,\cP^-)^{\on{spec,loc}}\to \Rep(\cG)$$
the functor obtained from $J^{-,\on{spec,pre-enh}}$ via the canonical self-dualities on $\Rep(\cG)$ and $\Rep(\cM)$, respectively. 

\medskip

Explicitly, $\on{co}\!J^{-,\on{spec,pre-enh}}$ is given by
$$\CF_M,\CF_{G,P^-}\mapsto (p_1)^\IndCoh_*\circ (p_2\times \on{Id})^{\IndCoh,*}(\CF_M\otimes \CF_{G,P^-}),$$
where
$$\CF_M\in \Rep(\cM)\simeq \QCoh(\LS^\reg_\cM)=:\IndCoh^*(\LS^\reg_\cM) \text{ and }
\CF_{G,P^-}\in \on{I}(\cG,\cP^-)^{\on{spec,loc}},$$
and $p_1$ and $p_2$ are the two maps in \eqref{e:gen Jacquet functor spec diagram}. 

\medskip

The functor $\on{co}\!J^{-,\on{spec,pre-enh}}$ factors naturally as
\begin{equation} \label{e:co J spec another}
\on{co}\!J^{-,\on{spec,pre-enh}}:\Rep(\cM)\otimes \on{I}(\cG,\cP^-)^{\on{spec,loc}} 
\to \Rep(\cM)\underset{\Sph^{\on{spec}}_\cM}\otimes \on{I}(\cG,\cP^-)^{\on{spec,loc}} \to \Rep(\cG).
\end{equation}

We denote the second arrow in \eqref{e:co J spec another} by $\on{co}\!J^{-,\on{spec,enh}}$. 

\sssec{} \label{sss:LS P G mf}

The functor $\on{co}\!J^{-,\on{spec,enh}}$ can also be interpreted as follows. Recall that according to \secref{sss:LS P G mf app}, 
we have
\begin{equation} \label{e:LS P G mf}
\Rep(\cM)\underset{\Sph^{\on{spec}}_{\cM}}\otimes \on{I}(\cG,\cP^-)^{\on{spec,loc}}
\simeq \IndCoh^*(\LS^{\mf_\cG}_{\cP^-}).
\end{equation} 

In terms of this identification, the functor $\on{co}\!J^{-,\on{spec,enh}}$ corresponds to the functor
$$\IndCoh^*(\LS^{\mf_\cG}_{\cP^-})\overset{(\sfp^{-,\on{spec,loc}}_\Ran)^\IndCoh_*}\longrightarrow 
\IndCoh^*(\LS^\reg_\cG)\simeq \QCoh(\LS^\reg_\cG)\simeq \Rep(\cG).$$

\sssec{}

From the commutative (but not Cartesian!) diagram 
$$
\CD
\LS^{\mf_\cG,\on{glob}}_{\cP^-,\CZ} @>{\sfq^{-,\on{spec,enh}}_\CZ}>> 
\LS^{\on{mer,glob}}_{\cM,\CZ}\underset{\LS^\mer_{\cM,\CZ}}\times \LS^{\mf_\cG}_{\cP^-,\CZ} \\
@V{\sfp^{-,\on{spec}}_\CZ}VV @VV{\on{ev}^\cM_\CZ\times \on{id}}V \\
\LS_\cG & &  \LS^\mer_{\cM,\CZ}\underset{\LS^\mer_{\cM,\CZ}}\times \LS^{\mf_\cG}_{\cP^-,\CZ}\simeq \LS^{\mf_\cG}_{\cP^-,\CZ} \\
@V{\on{ev}^\cG_\CZ}VV @VV{\sfp^{-,\on{spec,loc}}_\CZ}V \\
\LS^\reg_{\cG,\CZ} @>{\on{id}}>> \LS^\reg_{\cG,\CZ} 
\endCD
$$
we obtain a natural transformation 
\begin{equation} \label{e:Gamma Eis enh 1}
(\sfp^{-,\on{spec,loc}}_\CZ)^\IndCoh_*\circ (\on{ev}^\cM_\CZ\times \on{id})^\IndCoh_*\to 
(\on{ev}^\cG_\CZ)^\IndCoh_*\circ (\sfp^{-,\on{spec}}_\CZ)^\IndCoh_*\circ (\sfq^{-,\on{spec,enh}}_\CZ)^{\IndCoh,*}
\end{equation}
as functors
$$\IndCoh^*(\LS^{\on{mer,glob}}_{\cM,\CZ}\underset{\LS^\mer_{\cM,\CZ}}\times \LS^{\mf_\cG}_{\cP^-,\CZ})\rightrightarrows 
\IndCoh^*(\LS^\reg_\cG)_\CZ.$$

\sssec{}

Using the identification \eqref{e:LS P G mf}, we interpret \eqref{e:Gamma Eis enh 1} as a natural transformation
\begin{equation} \label{e:Gamma Eis enh 2}
\on{co}\!J^{-,\on{spec,enh}}\circ (\Gamma_\cM^{\on{spec},\IndCoh} \otimes \on{Id})\to \Gamma^{\on{spec},\IndCoh}_{\cG,\CZ}\circ \Eis^{-,\on{spec,enh}}_\CZ
\end{equation}
as functors
$$\IndCoh(\LS_\cM)^{-,\on{enh}}_\CZ\rightrightarrows \Rep(\cG)_\CZ,$$
i.e.,
\begin{equation} \label{e:Eis enh spec Gamma nat IndCoh}
\vcenter
{\xy
(0,0)*+{\IndCoh(\LS_\cG)\otimes \Dmod(\CZ)}="A";
(75,0)*+{\IndCoh(\LS_\cM)^{-,\on{enh}}_\CZ,}="B";
(0,30)*+{\Rep(\cG)_\CZ}="C";
(75,30)*+{(\Rep(\cM)\underset{\Sph^{\on{spec}}_\cM}\otimes \on{I}(\cG,\cP^-)^{\on{spec,loc}})_\CZ}="D";
{\ar@{->}^{\Eis^{-,\on{spec,enh}}} "B";"A"};
{\ar@{->}_{\on{co}\!J^{-,\on{spec,enh}}} "D";"C"};
{\ar@{->}^{\Gamma^{\on{spec},\IndCoh}_{\cG,\CZ}} "A";"C"};
{\ar@{->}_{\Gamma^{\on{spec},\IndCoh}_\cM\otimes \on{Id}} "B";"D"};
{\ar@{=>} "D";"A"};
\endxy}
\end{equation} 

\sssec{} \label{sss:nat trans enh spec}

Note that when we precompose both sides of \eqref{e:Gamma Eis enh 2} with
$$\ind^{\on{spec}}_{\on{enh}}:\IndCoh(\LS_\cM)\otimes \Dmod(\CZ)\to \IndCoh(\LS_\cM)^{-,\on{enh}}_\CZ,$$
we recover the natural transformation \eqref{e:Eis spec Gamma nat IndCoh}. 

\sssec{}

From \eqref{e:Eis enh spec Gamma nat IndCoh} we obtain a natural transformation from the counterclockwise to the
clockwise composition in the diagram

\medskip

\begin{equation} \label{e:Eis enh spec Gamma nat IndCoh 1}
\vcenter
{\xy
(0,0)*+{\Rep(\cG)_{\CZ^{\subseteq,\on{untl}},\on{indep}}}="A";
(90,0)*+{\Rep(\cG)_{\CZ^{\subseteq,\on{untl}}}}="B";
(0,-30)*+{\IndCoh(\LS_\cG)\otimes \Dmod(\CZ)}="C";
(90,-30)*+{(\Rep(\cM)\underset{\Sph_\cM^{\on{spec}}}\otimes \on{I}(\cG,\cP^-)^{\on{spec,loc}})_{\CZ^{\subseteq,\on{untl}}}}="D";
(90,-50)*+{\Rep(\cM)_{\CZ^{\subseteq,\on{untl}}}\underset{\Sph_{\cM,\CZ}^{\on{spec}}}\otimes 
\on{I}(\cG,\cP^-)^{\on{spec,loc}}_{\CZ^{\subseteq,\on{untl}}}}="E";
(0,-70)*+{\IndCoh(\LS_\cM)^{-,\on{enh}}_\CZ}="F";
(90,-70)*+{(\IndCoh(\LS_\cM)\otimes \Dmod(\CZ))\underset{\Sph_{\cM,\CZ}^{\on{spec}}}\otimes \on{I}(\cG,\cP^-)^{\on{spec,loc}}_\CZ}="G";
{\ar@{->}^{\on{emb.indep}_{\Rep(\cG)}^L} "B";"A"};
{\ar@{->}^{\Gamma^{\on{spec},\IndCoh}_{\cG,\CZ,\on{untl}}} "C";"A"};
{\ar@{->}^{\Eis^{-,\on{spec,enh}}} "F";"C"}; 
{\ar@{->}^{=} "F";"G"}; 
{\ar@{->}_{\on{co}\!J^{-,\on{spec,enh}}} "D";"B"}; 
{\ar@{->} "E";"D"}; 
{\ar@{->}_{\Gamma^{\on{spec},\IndCoh}_{\cM,\CZ,\on{untl}}\otimes \on{ins.unit}_\CZ} "G";"E"}; 
\endxy}
\end{equation} 
where the notations $\Gamma^{\on{spec},\IndCoh}_{\cG,\CZ,\on{untl}}$ and $\Gamma^{\on{spec},\IndCoh}_{\cM,\CZ,\on{untl}}$
have the same meaning as in \secref{sss:coeff Z untl}.

\sssec{}

We claim: 

\begin{thm} \label{t:Eis spec enh Gamma}
The natural transformation in \eqref{e:Eis enh spec Gamma nat IndCoh 1} is an isomorphism.
\end{thm} 

\begin{proof}

It is enough to show that the natural transformation in \eqref{e:Eis enh spec Gamma nat IndCoh 1} becomes an isomorphism
after precomposing with 
$$\ind^{\on{spec}}_{\on{enh}}:\IndCoh(\LS_\cM)\otimes \Dmod(\CZ)\to \IndCoh(\LS_\cM)^{-,\on{enh}}_\CZ.$$

However, by \secref{sss:nat trans enh spec}, up to tensoring with $\Dmod(\CZ)$, the resulting natural 
transformation equals one in \eqref{e:Eis spec Gamma nat 1 IndCoh}.
Now, the latter natural transformation is an isomorphism by \corref{c:Eis spec Gamma}.  

\end{proof} 

\ssec{The (co-)enhanced spectral constant term functor} \label{ss:CT spec enh} 

The goal of this section is to define the co-enhanced spectral constant term functor
$$\ul{\on{CT}}^{-,\on{spec},\on{enh}_{\on{co}}}:\ul\IndCoh(\LS_\cG)^{-,\on{enh}_{\on{co}}}\to \IndCoh(\LS_\cM)\otimes \ul\Dmod(\Ran),$$
where the left-hand side is the category defined in \secref{sss:spec co-enh}. 

\sssec{}

Denote
$$\LS^{\mf_\cM,\on{glob}}_{\cP^-,\Ran}:=(\LS_\cM\times \Ran)\underset{\LS^{\on{mer,glob}}_{\cM,\Ran}}\times \LS^{\on{mer,glob}}_{\cP^-,\Ran}.$$

\medskip

We will denote by $\sfq^{-,\on{spec}}_\Ran$ the map
$$\LS^{\mf_\cM,\on{glob}}_{\cP^-,\Ran}\to \LS_\cM\times \Ran.$$

\sssec{}

Let $\sfs^{\on{spec},\cM}_\Ran$ denote the naturally defined map
\begin{equation} \label{e:tricky ev M}
\LS^{\mf_\cM,\on{glob}}_{\cP^-,\Ran}\to 
\LS^\reg_{\cM,\Ran}\underset{\LS^\mer_{\cM,\Ran}}\times \LS^\mer_{\cP^-,\Ran}=:\LS^{\mf_\cM}_{\cP^-,\Ran}.
\end{equation} 

Note that we have a commutative (but not Cartesian) diagram 
$$
\CD
\LS^{\mf_\cM,\on{glob}}_{\cP^-,\Ran} @>{\sfq^{-,\on{spec}}_\Ran}>> \LS_\cM\times \Ran \\
@V{\sfs^{\on{spec},\cM}_\Ran}VV @VV{\on{ev}^\cM_\Ran}V \\
\LS^{\mf_\cM}_{\cP^-,\Ran} @>>{\sfq^{-,\on{spec,loc}}_\Ran}> \LS^\reg_{\cM,\Ran},
\endCD
$$
where 
$$\sfq^{-,\on{spec,loc}}_\Ran:\LS^{\mf_\cM}_{\cP^-,\Ran}\to \LS^\reg_{\cM,\Ran}$$
is the projection on the first factor.

\sssec{}

We have a naturally defined map
\begin{equation} \label{e:iota M glob}
\iota^{\on{glob},\cM}:\LS_{\cP^-}\times \Ran \to \LS^{\mf_\cM,\on{glob}}_{\cP^-,\Ran},
\end{equation} 
so that the diagram 
$$
\CD
\LS_{\cP^-}\times \Ran @>{\iota^{\on{glob,\cM}}}>> \LS^{\mf_\cM,\on{glob}}_{\cP^-,\Ran} \\
@V{\on{ev}^{\cP^-}_\Ran}VV @VV{\sfs^{\on{spec},\cM}_\Ran}V \\
\LS^\reg_{\cP^-,\Ran} @>>{\iota^\cM}> \LS^{\mf_\cM}_{\cP^-,\Ran}
\endCD
$$
is Cartesian.

\medskip

Note, however, that unlike the case of \eqref{e:iota G glob}, the map \eqref{e:iota M glob} is
\emph{not} a nil-isomorphism. 

\sssec{}

We let 
\begin{equation} \label{e:LS P G mf compl}
\left(\LS^{\mf_\cM,\on{glob}}_{\cP^-,\Ran}\right){}^\wedge_\mf
\end{equation} 
denote the formal completion of $\LS^{\mf_\cM,\on{glob}}_{\cP^-,\Ran}$ along the map \eqref{e:iota M glob}.

\medskip

By the same mechanism as in \secref{sss:LS P G mer good}, we obtain that \eqref{e:LS P G mf compl} is a formal algebraic stack, 
locally almost of finite type. In particular, the category 
$$\IndCoh\left(\left(\LS^{\mf_\cM,\on{glob}}_{\cP^-,\Ran}\right){}^\wedge_\mf\right)$$
is well-defined. 

\begin{rem} 

For a fixed $\ul{x}$, the prestack $\LS^{\mf_\cM,\on{glob}}_{\cP^-,\ul{x}}$ is actually an algebraic stack,
locally almost of finite type, so we could work directly with the $\IndCoh$ category on it. 

\medskip

However, we do not know how to define $\IndCoh$ of $\LS^{\mf_\cM,\on{glob}}_{\cP^-,\CZ}$ for 
general $\CZ\to \Ran$, cf. Remark \ref{r:no LS N mer}.  

\end{rem} 

\sssec{}

We claim that that there is a well-defined functor
$$\IndCoh(\LS_\cG)^{-,\on{enh}_{\on{co}}}_\CZ\to \IndCoh\left(\left(\LS^{\mf_\cM,\on{glob}}_{\cP^-,\CZ}\right){}^\wedge_\mf\right),$$
to be denoted $(\sfp^{-,\on{spec}}_\CZ)^{!,\on{enh}}$. 

\medskip

We will perform this construction in this subsection when $\CZ=\ul{x}$, and the case of a general $\CZ$ will be considered in 
\secref{sss:global IGP gen M}. 

\sssec{}

Let $\left(\LS^{\on{mer,glob}}_{\cG,\ul{x}}\right){}^\wedge_\mf$ be the formal completion of $\LS^{\on{mer,glob}}_{\cG,\ul{x}}$ along
$$\LS_\cG\to \LS^{\on{mer,glob}}_{\cG,\ul{x}}.$$

\medskip

Consider the map
$$\LS^{\mf_\cM}_{\cP^-,\ul{x}} \to \LS^\mer_{\cP^-,\ul{x}} \to \LS^\mer_{\cG,\ul{x}},$$
and the fiber product
\begin{equation} \label{e:LS G fp}
\left(\LS^{\on{mer,glob}}_{\cG,\ul{x}}\right){}^\wedge_\mf\underset{\LS^\mer_{\cG,\ul{x}}}\times \LS^{\mf_\cM}_{\cP^-,\ul{x}}.
\end{equation}

\medskip

We claim that when $\ul{x}$ is fixed, \eqref{e:LS G fp} is an algebraic stack locally almost of finite type.
This follows from deformation theory combined with the fact that the map
\begin{equation} \label{e:LS G fp 0}
\LS_\cG\underset{\LS^\reg_{\cG,\ul{x}}}\times \LS^\reg_{\cP^-,\ul{x}}\to
\left(\LS^{\on{mer,glob}}_{\cG,\ul{x}}\right){}^\wedge_\mf\underset{\LS^\mer_{\cG,\ul{x}}}\times \LS^{\mf_\cM}_{\cP^-,\ul{x}}
\end{equation}
is a nil-isomorphism. 

\medskip

In particular, the category 
$$\IndCoh\left(\left(\LS^{\on{mer,glob}}_{\cG,\ul{x}}\right){}^\wedge_\mf\underset{\LS^\mer_{\cG,\ul{x}}}\times \LS^{\mf_\cM}_{\cP^-,\ul{x}}\right)$$
is well-defined. Moreover, as in \eqref{e:amb vs IGP spec !}, we have a canonical isomorphism
\begin{multline} \label{e:LS G fp 1}
\IndCoh(\LS_\cG)^{-,\on{enh}_{\on{co}}}_{\ul{x}}:=
\IndCoh(\LS_\cG)\overset{\Sph^{\on{spec}}_{\cG,\ul{x}}}\otimes \on{I}(\cG,\cP^-)^{\on{spec,loc}}_{\on{co},\ul{x}}\simeq \\
\simeq \IndCoh\left(\left(\LS^{\on{mer,glob}}_{\cG,\ul{x}}\right){}^\wedge_\mf\underset{\LS^\mer_{\cG,\ul{x}}}\times \LS^{\mf_\cM}_{\cP^-,\ul{x}}\right).
\end{multline}

\sssec{}

Combining the maps
$$\LS^{\mf_\cM,\on{glob}}_{\cP^-,\ul{x}}\to \LS^{\on{mer,glob}}_{\cP^-,\ul{x}}\to \LS^{\on{mer,glob}}_{\cG,\ul{x}}$$
and $\sfs^{\on{spec},\cM}_{\ul{x}}$, we obtain a map
\begin{equation} \label{e:LS G fp 1.5}
\left(\LS^{\mf_\cM,\on{glob}}_{\cP^-,\ul{x}}\right){}^\wedge_\mf
\to \left(\LS^{\on{mer,glob}}_{\cG,\ul{x}}\right){}^\wedge_\mf\underset{\LS^\mer_{\cG,\ul{x}}}\times \LS^{\mf_\cM}_{\cP^-,\ul{x}}
\end{equation}
which we denote by $\sfp^{-,\on{spec,enh}}_{\ul{x}}$. 

\medskip

Thus, we can consider the functor
\begin{multline*}
(\sfp^{-,\on{spec,enh}}_{\ul{x}})^!:
\IndCoh\left(\left(\LS^{\on{mer,glob}}_{\cG,\ul{x}}\right){}^\wedge_\mf\underset{\LS^\mer_{\cG,\ul{x}}}\times \LS^{\mf_\cM}_{\cP^-,\ul{x}}\right)\to  \\
\to \IndCoh\left(\left(\LS^{\mf_\cM,\on{glob}}_{\cP^-,\ul{x}}\right){}^\wedge_\mf\right).
\end{multline*} 

\medskip

Combining with \eqref{e:LS G fp 1}, we obtain the sought-for functor
$$(\sfp^{-,\on{spec}}_{\ul{x}})^{!,\on{enh}}:\IndCoh(\LS_\cG)^{-,\on{enh}_{\on{co}}}_{\ul{x}}\to
\IndCoh\left(\left(\LS^{\mf_\cM,\on{glob}}_{\cP^-,\ul{x}}\right){}^\wedge_\mf\right).$$

\begin{rem}

Parallel to Remark \ref{r:old pullback}, the composition\footnote{In the formula below, the first arrow is the left adjoint
of the forgetful functor.}
\begin{multline*}
\IndCoh(\LS_\cG)\otimes \on{I}(\cG,\cP^-)^{\on{spec,loc}}_{\on{co}}\to \IndCoh(\LS_\cG)^{-,\on{enh}_{\on{co}}}_{\ul{x}}
\overset{(\sfp^{-,\on{spec}}_{\ul{x}})^{!,\on{enh}}}\longrightarrow \\
\to \IndCoh\left(\left(\LS^{\mf_\cM,\on{glob}}_{\cP^-,\ul{x}}\right){}^\wedge_\mf\right)
\end{multline*} 
can be described using the stack $\on{Hecke}^{\on{spec,glob}}_{\cG,\cP^-,\ul{x}}$. 

\end{rem}

\sssec{}

The spectral Hecke groupoid $\on{Hecke}^{\on{spec,loc}}_{\cM,\ul{x}}$  acts naturally on both sides of \eqref{e:LS G fp 1.5}
and the map $\sfp^{-,\on{spec,enh}}_{\ul{x}}$ is compatible with this action. 

\medskip

From here we obtain that the functor $(\sfp^{-,\on{spec}}_{\ul{x}})^{!,\on{enh}}$ is compatible with the
actions of $\Sph^{\on{spec}}_{\cM,\ul{x}}$ on the two sides. The same will be true when $\ul{x}\in \Ran$
is replaced by an arbitrary $\CZ\to \Ran$. 

\sssec{}

We define the functor
$$\on{CT}^{-,\on{spec},\on{enh}_{\on{co}}}_\CZ:\IndCoh(\LS_\cG)^{-,\on{enh}_{\on{co}}}_\CZ\to \IndCoh(\LS_\cM)\otimes \Dmod(\CZ)$$
as the composition
$$(\sfq^{-,\on{spec}}_\CZ)^\IndCoh_*\circ (\sfp^{-,\on{spec}}_\CZ)^{!,\on{enh}}.$$

\medskip

Since both functors $(\sfp^{-,\on{spec}}_\CZ)^{!,\on{enh}}$ and $(\sfq^{-,\on{spec}}_\CZ)^\IndCoh_*$ 
respect the actions of $\Sph^{\on{spec}}_{\cM,\ul{x}}$, we obtain that $\on{CT}^{-,\on{spec},\on{enh}_{\on{co}}}_\CZ$ also 
respects the actions of  $\Sph^{\on{spec}}_{\cM,\CZ}$.

\sssec{} 

Unwinding the construction, we obtain that the composition
\begin{equation} \label{e:ind enh co}
\IndCoh(\LS_\cG)\otimes \Dmod(\CZ)\to 
\IndCoh(\LS_\cG)^{-,\on{enh}_{\on{co}}}_\CZ \overset{\on{CT}^{-,\on{spec},\on{enh}_{\on{co}}}_\CZ}\longrightarrow 
\IndCoh(\LS_\cM)\otimes \Dmod(\CZ)
\end{equation} 
identifies canonically with $\on{CT}^{-,\on{spec}}\otimes \on{Id}$. 

\sssec{}

We will view the assignment
$$\CZ\rightsquigarrow \on{CT}^{-,\on{spec},\on{enh}_{\on{co}}}_\CZ$$
as a local-to-global functor 
$$\ul{\IndCoh}(\LS_\cG)^{-,\on{enh}_{\on{co}}}\to \IndCoh(\LS_\cM)\otimes \ul\Dmod(\Ran),$$
to be denoted $\ul{\on{CT}}^{-,\on{spec},\on{enh}_{\on{co}}}$.

\medskip

It has a natural lax unital structure. However, the isomorphism between \eqref{e:ind enh co} and $\on{CT}^{-,\on{spec}}\otimes \on{Id}$
implies that this lax unital structure is strict. 

\sssec{} \label{sss:fact alg enh spec co}

Let $\CA$ be a factorization algebra in $\on{I}(\cG,\cP^-)_{\on{co}}^{\on{spec,loc}}$. Parallel to \secref{sss:fact alg Eis} we can consider the functor
$$\on{CT}^{-,\on{spec}}_\CA:\IndCoh(\LS_\cG)\to \IndCoh(\LS_\cM).$$

\ssec{Partial enhancement for the constant term functor}

\sssec{}

Set
$$\IndCoh(\LS_\cG)^{-,\on{part.enh}_{\on{co}}}_\CZ:=
\left(\IndCoh(\LS_\cG)\otimes \Dmod(\CZ)\right)\underset{\Rep(\cG)_\CZ}\otimes \Rep(\cP^-)_\CZ.$$


\medskip 

We will consider the assignment
$$\CZ\rightsquigarrow \IndCoh(\LS_\cG)^{-,\on{part.enh}_{\on{co}}}_\CZ$$
as a crystal of categories over $\Ran$, to be denoted $$\ul{\IndCoh}(\LS_\cG)^{-,\on{part.enh}_{\on{co}}}.$$

\sssec{}

Consider the functor
$$(\sfp^{-,\on{spec}}_\CZ)^{!,\on{part.enh}_{\on{co}}}:\IndCoh(\LS_\cG)^{-,\on{part.enh}_{\on{co}}}_\CZ\to
\IndCoh(\LS_{\cP^-})\otimes \Dmod(\CZ)$$
defined as 
$$(\sfp^{-,\on{spec}})^!(-)\otimes (\on{ev}^{\cP^-}_\CZ)^*(-).$$

\medskip

Set
$$\on{CT}^{-,\on{spec},\on{enh}_{\on{co}}}_\CZ:=((\sfq^{-,\on{spec}})^\IndCoh_*\otimes \on{Id})\circ (\sfp^{-,\on{spec}}_\CZ)^{!,\on{part.enh}_{\on{co}}};$$
this is a functor 
$$\IndCoh(\LS_\cG)^{-,\on{part.enh}_{\on{co}}}_\CZ\to \IndCoh(\LS_\cM)\otimes \Dmod(\CZ).$$

\medskip

The composition
\begin{equation} \label{e:CT part enh comp}
\IndCoh(\LS_\cG)\otimes \Dmod(\CZ)\to \IndCoh(\LS_\cG)^{-,\on{part.enh}_{\on{co}}}_\CZ 
\overset{\on{CT}^{-,\on{spec},\on{enh}_{\on{co}}}_\CZ }\longrightarrow \IndCoh(\LS_\cM)\otimes \Dmod(\CZ)
\end{equation}
identifies canonically with $\on{CT}^{-,\on{spec}}\otimes \on{Id}$. 

\sssec{}

We will view the assignment
$$\CZ\rightsquigarrow \on{CT}^{-,\on{spec,part.enh}_{\on{co}}}_\CZ$$
as a local-to-global functor 
$$\ul{\IndCoh}(\LS_\cG)^{-,\on{part.enh}_{\on{co}}}\to \IndCoh(\LS_\cM)\otimes \ul\Dmod(\Ran),$$
to be denoted $\ul{\on{CT}}^{-,\on{spec,part.enh}_{\on{co}}}$.

\medskip

It has a natural lax unital structure. However, the isomorphism between \eqref{e:CT part enh comp} and
$\on{CT}^{-,\on{spec}}\otimes \on{Id}$
implies that this lax unital structure is strict. 

\sssec{} \label{sss:fact alg enh spec partial co}

Let $\CB$ be a factorization algebra in $\Rep(\cP^-)$. Parallel to \secref{sss:fact alg Eis} we can consider the functor
$$\on{CT}^{-,\on{spec}}_\CB:\IndCoh(\LS_\cG)\to \IndCoh(\LS_\cM).$$

\sssec{}

Recall the factorization functor
$$\iota^\IndCoh_*:\IndCoh^!(\LS^\reg_{\cP^-})\to \on{I}(\cG,\cP^-)^{\on{spec,loc}}_{\on{co}}$$
introduced in \secref{sss:Delta spec co to Delta spec}. 

\medskip

Since $\Rep(\cG)$ is rigid, we can interpret $\IndCoh(\LS_\cG)^{-,\on{part.enh}_{\on{co}}}_\CZ$ as
$$\left(\IndCoh(\LS_\cG)\otimes \Dmod(\CZ)\right)\overset{\Rep(\cG)_\CZ}\otimes \IndCoh^!(\LS^\reg_{\cP^-})_\CZ.$$

\medskip

Tensoring $\iota_*$ with the identity functor on $\IndCoh(\LS_\cG)$, we obtain a functor
$$(\on{Id}\otimes \iota_*):\IndCoh(\LS_\cG)^{-,\on{part.enh}_{\on{co}}}_\CZ\to
\IndCoh(\LS_\cG)^{-,\on{part.enh}}_\CZ.$$

\sssec{}

Unwinding the constructions, we obtain a canonical isomorphism
\begin{equation} \label{e:enh vs part spectral CT pre}
(\sfp^{-,\on{spec}}_{\ul{x}})^{!,\on{enh}} \circ (\on{Id}\otimes \iota_*)\simeq 
(\iota^{\on{glob},\cM})^\IndCoh_*\circ (\sfp^{-,\on{spec}}_\CZ)^{!,\on{part.enh}_{\on{co}}},
\end{equation} 
as functors
$$\IndCoh(\LS_\cG)^{-,\on{part.enh}_{\on{co}}}_\CZ\rightrightarrows \IndCoh\left(\left(\LS^{\mf_\cM,\on{glob}}_{\cP^-,\CZ}\right){}^\wedge_\mf\right),$$
where we view $(\iota^{\on{glob},\cM})^\IndCoh_*$ as a functor
$$\IndCoh(\LS_{\cP^-})\otimes \Dmod(\CZ)\to \IndCoh\left(\left(\LS^{\mf_\cM,\on{glob}}_{\cP^-,\CZ}\right){}^\wedge_\mf\right).$$

From here we obtain an isomorphism
\begin{equation} \label{e:enh vs part spectral CT}
\ul{\on{CT}}^{-,\on{spec,enh}_{\on{co}}}\circ (\on{Id}\otimes \iota_*)\simeq \ul{\on{CT}}^{-,\on{spec,part.enh}_{\on{co}}}
\end{equation} 
as local-to-global functors 
$$\ul{\IndCoh}(\LS_\cG)^{-,\on{part.enh}_{\on{co}}}\rightrightarrows \IndCoh(\LS_\cM)\otimes \ul\Dmod(\Ran).$$

\sssec{}

Note that, thanks to \eqref{e:enh vs part spectral CT}, for a factorization algebra $\CB$ in $\Rep(\cB^-)$, we have
$$\on{CT}^{-,\on{spec}}_\CB\simeq \on{CT}^{-,\on{spec}}_{\iota^\IndCoh_*(\CB)},$$
where:

\begin{itemize}

\item The left-hand side is understood in the sense of \secref{sss:fact alg enh spec partial co};

\item The right-hand side  is understood in the sense of \secref{sss:fact alg enh spec}. 

\end{itemize} 

\ssec{The ``compactified" spectral constant term functor} \label{ss:comp spec CT}

\sssec{}

In the framework of \secref{sss:fact alg enh spec co}, take $\CA$ to be the object
$$\IC^{-,\on{spec},\semiinf}_{\on{co}}\in \on{I}(\cG,\cP^-)^{\on{spec,loc}}_{\on{co}},$$
introduced in \secref{sss:semiinf IC spec co}. 

\medskip

Denote the resulting functor $\on{CT}^{-,\on{spec}}_\CA$ by
$$\on{CT}^{-,\on{spec}}_{\IC}:\IndCoh(\LS_\cG)\to \IndCoh(\LS_\cM).$$

\sssec{} \label{sss:comp spec CT}

Consider now the functor
$$\on{CT}^{-,\on{spec}}_{!*}:\IndCoh(\LS_\cG)\to \IndCoh(\LS_\cM),$$
defined as follows.

\medskip

Note that since the functor
$$(\sfq^{-,\on{spec}})^\IndCoh_*:\IndCoh(\LS_{\cP^-})\to \IndCoh(\LS_\cM)$$
is $\QCoh(\LS_\cM)$-linear, it upgrades to a functor
$$(\sfq^{-,\on{spec}})^{\IndCoh,\on{enh}}_*:\IndCoh(\LS_{\cP^-})\to 
\IndCoh(\LS_\cM)\underset{\QCoh(\LS_\cM)}\otimes \QCoh(\LS_{\cP^-}),$$
so that the composition
\begin{multline*} 
\IndCoh(\LS_{\cP^-})\overset{(\sfq^{-,\on{spec}})^{\IndCoh,\on{enh}}_*}\longrightarrow 
\IndCoh(\LS_\cM)\underset{\QCoh(\LS_\cM)}\otimes \QCoh(\LS_{\cP^-}) \overset{\on{Id}\otimes \sfq^{-,\on{spec}}_*}\longrightarrow \\
\to \IndCoh(\LS_\cM)\underset{\QCoh(\LS_\cM)}\otimes \QCoh(\LS_\cM)=\IndCoh(\LS_\cM)
\end{multline*} 
is the original functor $(\sfq^{-,\on{spec}})^\IndCoh_*$. 

\medskip

We let $\on{CT}^{-,\on{spec}}_{!*}$ be the composition
\begin{multline} \label{e:CT spec IC defn}
\IndCoh(\LS_\cG) \overset{(\sfp^{-,\on{spec}})^!}\longrightarrow \IndCoh(\LS_{\cP^-}) 
\overset{(\sfq^{-,\on{spec}})^{\IndCoh,\on{enh}}_*}\longrightarrow \\
\to \IndCoh(\LS_\cM)\underset{\QCoh(\LS_\cM)}\otimes \QCoh(\LS_{\cP^-}) \to \\
\to \IndCoh(\LS_\cM)\underset{\QCoh(\LS_\cM)}\otimes \QCoh(\LS_\cM)=\IndCoh(\LS_\cM),
\end{multline}
where the third arrow is the tensor product of $\on{Id}_{\IndCoh(\LS_\cM)}$ with the functor
$$\QCoh(\LS_{\cP^-}) \to \QCoh(\LS_\cM),$$
given by *-pullback along 
$$\LS_\cM\to \LS_{\cP^-},$$
cf. \secref{sss:defn Eis spec IC}. 

\sssec{}

From \lemref{l:rel det LS}, we obtain:

\begin{lem} \label{l:spec CT and Eis IC}
With respect to the Serre duality identifications
$$\IndCoh(\LS_\cG)^\vee \simeq \IndCoh(\LS_\cG) \text{ and } \IndCoh(\LS_\cM)^\vee \simeq \IndCoh(\LS_\cM),$$
there is a canonical isomorphism
$$(\on{CT}^{-,\on{spec}}_{!*})^\vee \simeq \Eis^{-,\on{spec}}_{!*}\circ 
\Bigl((-)\otimes \CL_{2\rho_P(\omega_X)}\Bigr)[(2g-2)\cdot \dim(\cn^-_P)].$$
\end{lem}

\sssec{}

We are going to prove:

\begin{prop} \label{p:CT spec IC}
The functors $\on{CT}^{-,\on{spec}}_{\IC}$ and $\on{CT}^{-,\on{spec}}_{!*}$
are canonically isomorphic. 
\end{prop} 

The rest of this subsection is devoted to the proof of this proposition. 

\sssec{}

Let $\CB$ be a factorization algebra in $\Rep(\cP^-)$. By \eqref{e:enh vs part spectral CT},
it suffices to establish an isomorphism
$$\on{CT}^{-,\on{spec}}_{!*}\simeq \on{CT}^{-,\on{spec}}_{\CB}$$
for $\on{CT}^{-,\on{spec}}_{\CB}$ as in \secref{sss:fact alg enh spec partial co} when $\CB$ is 
as in \secref{sss:Eis !* spec 1}. 

\sssec{}

By construction, for any $\CB\in \on{FactAlg}(\Rep(\cP^-))$, we have 
$$\on{CT}^{-,\on{spec}}_{\CB}\simeq 
(\sfq^{-,\on{spec}})^\IndCoh_*\left((\sfp^{-,\on{spec}})^!(-)\otimes \Loc^{\on{spec}}_{\cP^-,\CB}(k)\right).$$

\medskip

In Sects. \ref{sss:Eis !* spec 1}-\ref{sss:Eis !* spec end} it was shown that for $\CB$ being the image of $\one_{\Rep(\cM)}$ along
$$\Rep(\cM) \overset{\on{coind}^{\cP^-}_\cM}\longrightarrow \Rep(\cP^-),$$
the object 
$$\Loc^{\on{spec}}_{\cP^-,\CB}(k)\in \QCoh(\LS_{\cP^-})$$
identifies with the direct image of $\CO_{\LS_\cM}$ along
$$\LS_\cM\to \LS_{\cP^-};$$
by a slight abuse of notation, we will denote this object by the same symbol $\CO_{\LS_\cM}$. 

\sssec{}

Thus, it remains to show that the functor
$$\IndCoh(\LS_{\cP^-}) \overset{(-)\otimes \CO_{\LS_\cM}}\longrightarrow \IndCoh(\LS_{\cP^-}) \overset{(\sfq^{-,\on{spec}})^\IndCoh_*}\longrightarrow 
\IndCoh(\LS_\cM)$$ identifies with the composition of the last three arrows in \eqref{e:CT spec IC defn}. 

\medskip

However, this follows from the fact that both functors in question can be identified with 
\begin{multline*}
\IndCoh(\LS_{\cP^-}) \overset{(\sfq^{-,\on{spec}})^{\IndCoh,\on{enh}}_*}\longrightarrow 
\IndCoh(\LS_\cM)\underset{\QCoh(\LS_\cM)}\otimes \QCoh(\LS_{\cP^-}) \overset{\on{Id}\otimes ((-)\otimes \CO_{\LS_\cM})}\longrightarrow \\
\to \IndCoh(\LS_\cM)\underset{\QCoh(\LS_\cM)}\otimes \QCoh(\LS_{\cP^-}) 
\overset{\on{Id}\otimes (\sfq^{-,\on{spec}})_*}\longrightarrow \\
\to \IndCoh(\LS_\cM)\underset{\QCoh(\LS_\cM)}\otimes \QCoh(\LS_\cM)=\IndCoh(\LS_\cM).
\end{multline*}

\qed[\propref{p:CT spec IC}]

\ssec{Enhancement and duality on the spectral side}

The contents of this subsection will not be necessary for the sequel.

\sssec{}

Unwinding the construction of the functor $\on{CT}_\CZ^{-,\on{spec},\on{enh}_{\on{co}}}$, we obtain that it can be obtained
from the functor $\Eis_\CZ^{-,\on{spec},\on{enh}}$ as follows:

\medskip

The functor $\Eis_\CZ^{-,\on{spec},\on{enh}}$ preserves compactness; hence admits a right adjoint, which is a functor
\begin{multline} \label{e:Eis spec enh Ra}
\IndCoh(\LS_\cG)\otimes \Dmod(\CZ) \to \\
\to (\IndCoh(\LS_\cM)\otimes \Dmod(\CZ))\underset{\Sph^{\on{spec}}_{\cM,\CZ}}\otimes 
\on{I}(\cG,\cP^-)^{\on{spec,loc}}_\CZ=:\IndCoh(\LS_\cM)^{-,\on{enh}}_\CZ.
\end{multline}

Denote the above functor \eqref{e:Eis spec enh Ra} by $\on{CT}_\CZ^{-,\on{spec},\on{enh}}$. 

\medskip 

Using the relative rigidity of $\Sph^{\on{spec}}_{\cM,\CZ}$ and $\Sph^{\on{spec}}_{\cG,\CZ}$ over $\CZ$, from $\on{CT}_\CZ^{-,\on{spec},\on{enh}}$,
we obtain a functor
$$(\IndCoh(\LS_\cG)\otimes \Dmod(\CZ))\overset{\Sph^{\on{spec}}_{\cM,\CZ}}\otimes \on{I}(\cG,\cP^-)^{\on{spec,loc}}_{\on{co},\CZ}\to 
\IndCoh(\LS_\cM)\otimes \Dmod(\CZ),$$
which is our functor $\on{CT}_\CZ^{-,\on{spec},\on{enh}_{\on{co}}}$.

\sssec{}

Note now that there is another variant for enhancement: we can take the \emph{dual} of $\on{CT}^{-,\on{enh}}_{*,\CZ}$
and thus obtain a functor
\begin{multline*} 
\IndCoh(\LS_\cM)^{-,\on{enh}_{\on{co}}}:=
(\IndCoh(\LS_\cM)\otimes \Dmod(\CZ))\overset{\Sph^{\on{spec}}_{\cM,\CZ}}\otimes 
\on{I}(\cG,\cP^-)^{\on{spec,loc}}_{\on{co},\CZ}\to \\
\to \IndCoh(\LS_\cG)\otimes \Dmod(\CZ),
\end{multline*}
to be denoted $\Eis_\CZ^{-,\on{spec},\on{enh}_{\on{co}}}$. 

\medskip

The following functors are each others counterparts on the spectral and geometric sides, respectively:
$$\on{CT}_\CZ^{-,\on{spec},\on{enh}} \leftrightarrow \on{CT}^{-,\on{enh}}_{*,\CZ}$$
$$\on{CT}_\CZ^{-,\on{spec},\on{enh}_{\on{co}}}  \leftrightarrow \on{CT}^{-,\on{enh}_{\on{co}}}_{*,\CZ}$$
$$\Eis_\CZ^{-,\on{spec},\on{enh}} \leftrightarrow \Eis_{!,\CZ}^{-,\on{enh}}$$
and 
$$\Eis_\CZ^{-,\on{spec},\on{enh}_{\on{co}}} \leftrightarrow \Eis^{-,\on{enh}_{\on{co}}}_{\on{co},*,\CZ}.$$

\medskip

Recall now that there was another kind of relationship between the functors $\Eis_!^{-,\on{enh}}$ and $\Eis^{-,\on{enh}_{\on{co}}}_{\on{co},*}$,
namely, one given by \thmref{t:intertwiners}. 

\medskip

Hence, one expects a similar relationship between the functors $\Eis_\CZ^{-,\on{spec},\on{enh}}$
and $\Eis_\CZ^{-,\on{spec},\on{enh}_{\on{co}}}$.

\sssec{}

Recall the equivalence
$$\on{I}(\cG,\cP^-)^{\on{spec,loc}}_{\on{co},\CZ} \simeq \on{I}(\cG,\cP^-)^{\on{spec,loc}}_{\CZ}$$
of \eqref{e:self-duality I G P spec}; denote it by $\Theta_{\on{I}(\cG,\cP^-)^{\on{spec,loc}}}$. 

\medskip

Recall also that the monoidal category $\Sph^{\on{spec}}_{\cM,\CZ}$ is rigid relative to $\CZ$,
so we can identify $$\overset{\Sph^{\on{spec}}_{\cM,\CZ}}\otimes \text{ and } \underset{\Sph^{\on{spec}}_{\cM,\CZ}}\otimes.$$

\medskip

Let $\Theta_{\on{I}(\cG,\cP^-)^{\on{spec,glob}}}$ denote the equivalence
$$\IndCoh(\LS_\cM)^{-,\on{enh}_{\on{co}}}_\CZ\to \IndCoh(\LS_\cM)^{-,\on{enh}}_\CZ$$
obtained by tensoring $\Theta_{\on{I}(\cG,\cP^-)^{\on{spec,loc}}}$ with the functor
$$((-)\otimes \CL_{2\rho_P(\omega_X)})[(2g-2)\cdot \dim(\cn^-_P)]:
\IndCoh(\LS_\cM)\to \IndCoh(\LS_\cM)$$
over\footnote{Note that according to Remark \ref{r:IGP spec duality}
it is the latter functor that makes the actions of $\Sph^{\on{spec}}_{\cM,\CZ}$ on the two sides compatible.}
$\Sph^{\on{spec}}_{\cM,\CZ}$. 

\sssec{}

We propose:

\begin{conj} \label{c:Theta spec glob}
The following diagram of functors commutes:
$$
\CD
\IndCoh(\LS_\cM)^{-,\on{enh}_{\on{co}}}_\CZ @>{\Theta_{\on{I}(\cG,\cP^-)^{\on{spec,glob}}}}>> \IndCoh(\LS_\cM)^{-,\on{enh}}_\CZ \\
@V{\Eis_\CZ^{-,\on{spec},\on{enh}_{\on{co}}}}VV @VV{\Eis_\CZ^{-,\on{spec},\on{enh}}}V \\
\IndCoh(\LS_\cG)\otimes \Dmod(\CZ) @>>{\on{id}}> \IndCoh(\LS_\cG)\otimes \Dmod(\CZ), 
\endCD
$$
in a way compatible with the $\Sph^{\on{spec}}_{\cG,\CZ}$-actions.
\end{conj}

\begin{rem}
As a reality check, note, that according to \eqref{e:CT and Eis dual}, 
the diagram in \conjref{c:Theta spec glob} does become commutative when juxtaposed with the
commutative square
$$
\CD
\IndCoh(\LS_\cM)\otimes \Dmod(\CZ) @>{((-)\otimes \CL_{2\rho_P(\omega_X)})[(2g-2)\cdot \dim(\cn^-_P)]}>> \IndCoh(\LS_\cM)\otimes \Dmod(\CZ) \\
@V{\ind_{\on{enh}_{\on{co}}}}VV @VV{\ind_{\on{enh}}}V \\
\IndCoh(\LS_\cM)^{-,\on{enh}_{\on{co}}}_\CZ @>{\Theta_{\on{I}(\cG,\cP^-)^{\on{spec,glob}}}}>> \IndCoh(\LS_\cM)^{-,\on{enh}}_\CZ.
\endCD
$$
\end{rem}

\sssec{}

Note that the functor $\Theta_{\on{I}(\cG,\cP^-)^{\on{spec,glob}}}$ induces a self-duality on $\IndCoh(\LS_\cM)^{-,\on{enh}}_\Ran$,
which in turn induces a self-duality on $\IndCoh(\LS_\cM)^{-,\on{enh}}_{\Ran^{\on{untl}},\on{indep}}$. 

\medskip

Recall now that according to \thmref{t:Nick}, we have a canonical equivalence
$$\IndCoh(\LS_\cM)^{-,\on{enh}}_{\Ran^{\on{untl}},\on{indep}}\simeq \IndCoh_\cM\left((\LS_{\cP^-})_\dr\underset{(\LS_{\cG})_\dr}\times \LS_{\cG}\right).$$

Thus, $\Theta_{\on{I}(\cG,\cP^-)^{\on{spec,glob}}}$ induces a self-duality on 
$\IndCoh_\cM\left((\LS_{\cP^-})_\dr\underset{(\LS_{\cG})_\dr}\times \LS_{\cG}\right)$.

\sssec{}

We propose:

\begin{conj} \label{c:Nick dual}
The self-duality on $\IndCoh_\cM\left((\LS_{\cP^-})_\dr\underset{(\LS_{\cG})_\dr}\times \LS_{\cG}\right)$ induced by 
the functor $\Theta_{\on{I}(\cG,\cP^-)^{\on{spec,glob}}}$
equals the canonical (Serre) self-duality defined on $\IndCoh$ of a prestack locally almost of finite type.
\end{conj}

\begin{rem}

Note that \conjref{c:Nick dual} implies \conjref{c:Theta spec glob}.

\end{rem} 

\section{Spectral constant term functor and opers} \label{s:Eis and Op}

The goal of this section is to express the composition 
$$\on{CT}^{-,\on{spec}}\circ \on{Poinc}^{\on{spec}}_{\cG,!}$$
via the (local) functor $J_\Op^{-,\on{spec}}$ followed by $\on{Poinc}^{\on{spec}}_{\cM,!}$,
and similarly for the enhanced versions of the functors involved.

\medskip

This will be \emph{parallel} to the ``Whittaker vs Eisenstein" calculation in \secref{s:Eis and Whit}. That said, this calculation 
is the \emph{Langlands dual} of the ``localization vs constant term" calculation from \secref{s:Eis and Loc}.

\ssec{Statement of the result--the unenhanced version}

In this subsection we will state a theorem that relates $\on{CT}^{-,\on{spec}}\circ \on{Poinc}^{\on{spec}}_{\cG,!}$
and $\on{Poinc}^{\on{spec}}_{\cM,!}\circ J_\Op^{-,\on{spec}}$.

\sssec{}

Recall the (unital) local-to-global functor 
$$\ul{\on{Poinc}}^{\on{spec}}_{\cG,!}:\ul\IndCoh^!(\Op^\mf_\cG)\to \IndCoh(\LS_\cG)\otimes \ul\Dmod(\Ran),$$
see \cite[Sect. 17.4.1]{GLC2} or \secref{sss:Poinc spec} below. 

\medskip

For a $Z^0_\cG$-torsor $\CP_{Z^0_\cG}$, recall the space of translated opers, see \secref{sss:transl oper}.
By a similar token, we obtain a functor
$$\ul{\on{Poinc}}^{\on{spec}}_{\cG,!,\CP_{Z^0_\cG}}:\ul\IndCoh^!(\Op^\mf_{\cG,\CP_{Z^0_\cG}})\to \IndCoh(\LS_\cG)\otimes \ul\Dmod(\Ran).$$

In particular, we can consider the functor
$$\ul{\on{Poinc}}^{\on{spec}}_{\cM,!,\rhoch_P}:\ul\IndCoh^!(\Op^\mf_{\cM,\rhoch_P})\to \IndCoh(\LS_\cM)\otimes \ul\Dmod(\Ran).$$

\sssec{}

On the one hand, we consider the (unital) local-to-global functor
$$\on{CT}^{-,\on{spec}}\circ \ul{\on{Poinc}}^{\on{spec}}_{\cG,!}:\ul\IndCoh^!(\Op^\mf_\cG) \to \IndCoh(\LS_\cM)\otimes \ul\Dmod(\Ran).$$

On the other hand, we consider the \emph{lax} unital local-to-global functor
$$\ul{\on{Poinc}}^{\on{spec}}_{\cM,!,\rhoch_P}\circ J_\Op^{-,\on{spec}}:\ul\IndCoh^!(\Op^\mf_\cG) \to \IndCoh(\LS_\cM)\otimes \ul\Dmod(\Ran),$$
where 
$$J_\Op^{-,\on{spec}}:\IndCoh^!(\Op^\mf_\cG)\to \IndCoh^!(\Op^\mf_{\cM,\rhoch_P})$$
is the factorization functor from \secref{sss:J opers !}.

\medskip

Our first goal will be to construct a natural transformation
\begin{equation} \label{e:Op nat trans}
\ul{\on{Poinc}}^{\on{spec}}_{\cM,!,\rhoch_P}\circ J_\Op^{-,\on{spec}}\to \on{CT}^{-,\on{spec}}\circ \ul{\on{Poinc}}^{\on{spec}}_{\cG,!}.
\end{equation} 

\medskip

Fix $\CZ\to \Ran$. Thus, we wish to construct a natural transformation
\begin{equation} \label{e:Op nat trans x}
\on{Poinc}^{\on{spec}}_{\cM,!,\rhoch_P,\CZ}\circ J_\Op^{-,\on{spec}}\to \on{CT}^{-,\on{spec}}\circ \on{Poinc}^{\on{spec}}_{\cG,!,\CZ}
\end{equation} 
as functors
$$\IndCoh^!(\Op^\mf_{\cG,\CZ})\rightrightarrows \IndCoh(\LS_\cM)\otimes \Dmod(\CZ).$$

\sssec{} \label{sss:Poinc spec} 

We begin by rewriting the right-hand side in \eqref{e:Op nat trans x}. Recall the space 
$$\Op^{\mf,\on{glob}}_{\cG,\CZ}:=\Op^{\mer,\on{glob}}_{\cG,\CZ}\underset{\LS^{\mer,\on{glob}}_{\cG,\CZ}}\times (\LS_\cG\times \CZ)$$
(see \cite[Sect. 15.1.1]{GLC2} for the notation).

\medskip

Recall (see \cite[Sect. 17.4.1]{GLC2})
that the functor $\on{Poinc}^{\on{spec}}_{\cG,!,\CZ}$
is given by !-pull followed by $(\IndCoh,*)$-push along the diagram
$$
\CD
& & \LS_\cG\times \CZ \\
& & @AA{\fr^{\on{glob}}_\cG}A \\
\Op^\mf_{\cG,\CZ} @<{\on{ev}^{\Op_\cG}_{\CZ}}<< \Op^{\mf,\on{glob}}_{\cG,\CZ}.
\endCD
$$

\medskip

Then the right-hand side in \eqref{e:Op nat trans x} is given by the diagram
\begin{equation} \label{e:Op nat trans 1}
\CD
& & & & \LS_\cM \times \CZ \\
& & & & @AA{\sfq^{-,\on{spec}}}A \\
& & \LS_\cG\times \CZ @<<{\sfp^{-,\on{spec}}}< \LS_{\cP^-} \times \CZ\\
& & @AA{\fr^{\on{glob}}_\cG}A \\
\Op^\mf_{\cG,\CZ} @<{\on{ev}^{\Op_\cG}_{\CZ}}<< \Op^{\mf,\on{glob}}_{\cG,\CZ},
\endCD
\end{equation}
in which we apply !-pullback along the horizontal arrows and $(\IndCoh,*)$-pushforward along the vertical
ones. 

\medskip

By base change, we can replace \eqref{e:Op nat trans 1} by the diagram
\begin{equation} \label{e:Op nat trans 2}
\CD
& & & & \LS_\cM\times \CZ \\
& & & & @AA{\sfq^{-,\on{spec}}}A \\
& & & & \LS_{\cP^-} \times \CZ\\
& & & & @AA{'\fr^{\on{glob}}_\cG}A \\
\Op^\mf_{\cG,\CZ} @<<{\on{ev}^{\Op_\cG}_{\CZ}}< \Op^{\mf,\on{glob}}_{\cG,\CZ} @<<{'\sfp^{-,\on{spec}}}< 
\Op^{\mf,\on{glob}}_{\cG,\CZ}\underset{\LS_\cG\times \CZ}\times (\LS_{\cP^-}\times \CZ)
\endCD
\end{equation}

\sssec{}

Recall the affine D-scheme $\MOp_{\cG,\cP^-}$ (see \secref{sss:Miura}). Consider the corresponding ind-scheme
$$\MOp^{\mer,\on{glob}}_{\cG,\cP^-,\CZ}$$
and the fiber product 
$$\MOp^{\mer,\mf_{\cP^-},\on{glob}}_{\cG,\cP^-,\CZ}:=
\MOp^{\mer,\on{glob}}_{\cG,\cP^-,\CZ}\underset{\LS^{\mer,\on{glob}}_{\cP^-,\CZ}}\times (\LS_{\cP^-}\times \CZ).$$

We have a naturally defined map
\begin{equation} \label{e:j glob}
\jmath^{\on{glob}}:
\MOp^{\mer,\mf_{\cP^-},\on{glob}}_{\cG,\cP^-,\CZ}\to
\Op^{\mf,\on{glob}}_{\cG,\CZ}\underset{\LS_\cG\times \CZ}\times (\LS_{\cP^-}\times \CZ).
\end{equation} 

\sssec{}

The following observation is key for our construction:

\begin{prop} \label{p:j glob}
The map $\jmath^{\on{glob}}$ of \eqref{e:j glob} is ind-proper.\footnote{In fact, as the proof will show, this map is an ind-closed embedding.}
\end{prop}

In a sense, this proposition explains the mechanism of the appearance of Miura opers. 

\medskip

The proof of \propref{p:j glob} will be given in \secref{sss:j glob}. We assume it temporarily and proceed with the
construction of the natural transformation \eqref{e:Op nat trans x}. 

\sssec{} \label{sss:crucial nat trans}

From \propref{p:j glob} we obtain an adjoint pair $((\jmath^{\on{glob}})^\IndCoh_*,(\jmath^{\on{glob}})^!)$, and hence 
a natural transformation
$$(\jmath^{\on{glob}})^\IndCoh_*\circ (\jmath^{\on{glob}})^!\to \on{Id}.$$

\medskip

Hence, the right-hand side in \eqref{e:Op nat trans x} receives a map from the functor
\begin{equation} \label{e:Op nat trans 3}
(\sfq^{-,\on{spec}})^\IndCoh_*\circ  
({}'\fr^{\on{glob}}_\cG)^\IndCoh_*\circ (\jmath^{\on{glob}})^\IndCoh_*\circ (\jmath^{\on{glob}})^!\circ ({}'\sfp^{-,\on{spec}})^!\circ (\on{ev}^{\Op_\cG}_{\CZ})^!.
\end{equation} 

We will now show that the functor \eqref{e:Op nat trans 3} identifies with the left-hand side in \eqref{e:Op nat trans x}. 

\sssec{} \label{sss:CT Op first}

Let $\sfp^{\on{Miura,glob}}$ denote the composition
$$\MOp^{\mer,\mf_{\cP^-},\on{glob}}_{\cG,\cP^-,\CZ}\overset{\jmath^{\on{glob}}}\longrightarrow 
\Op^{\mf,\on{glob}}_{\cG,\CZ}\underset{\LS_\cG\times \CZ}\times (\LS_{\cP^-}\times \CZ)\overset{'\sfp^{-,\on{spec}}}\longrightarrow 
\Op^{\mf,\on{glob}}_{\cG,\CZ}.$$

Note also that we have a naturally defined map 
$$\sfq^{\on{Miura,glob}}:
\MOp^{\mer,\mf_{\cP^-},\on{glob}}_{\cG,\cP^-,\CZ}\to 
\Op^{\mf,\on{glob}}_{\cM,\rhoch_P,\CZ}$$
so that the composition
$$\MOp^{\mer,\mf_{\cP^-},\on{glob}}_{\cG,\cP^-,\CZ}\overset{\jmath^{\on{glob}}}\longrightarrow 
\Op^{\mf,\on{glob}}_{\cG,\CZ}\underset{\LS_\cG\times \CZ}\times (\LS_{\cP^-}\times \CZ) \overset{'\fr^{\on{glob}}_\cG}\longrightarrow
\LS_{\cP^-}\times \CZ\overset{\sfp^{-,\on{spec}}\times \on{id}}\longrightarrow \LS_\cM\times \CZ$$
equals
$$\MOp^{\mer,\mf_{\cP^-},\on{glob}}_{\cG,\cP^-,\CZ}
\overset{\sfq^{\on{Miura,glob}}}\longrightarrow \Op^{\mf,\on{glob}}_{\cM,\rhoch_P,\CZ} \overset{\fr^{\on{glob}}_\cM}\to \LS_\cM\times \CZ.$$

\medskip

With these notations, the functor \eqref{e:Op nat trans 3} is given by pull-push along the diagram
\begin{equation} \label{e:Op nat trans 4}
\CD
& & & & \LS_\cM\times \CZ \\ 
& & & & @AA{\fr^{\on{glob}}_\cM}A \\
& & & & \Op^{\mf,\on{glob}}_{\cM,\rhoch_P,\CZ} \\
& & & & @AA{\sfq^{\on{Miura,glob}}}A \\
\Op^\mf_{\cG,\CZ} @<<{\on{ev}^{\Op_\cG}_{\CZ}}< \Op^{\mf,\on{glob}}_{\cG,\CZ} @<<{\sfp^{\on{Miura,glob}}}< 
\MOp^{\mer,\mf_{\cP^-},\on{glob}}_{\cG,\cP^-,\CZ}. 
\endCD
\end{equation} 

\sssec{}

Note that the evaluation map
$$\on{ev}^{\MOp}_{\CZ}:\MOp^{\mer,\on{glob}}_{\cG,\cP^-,\CZ}\to \MOp^\mer_{\cG,\cP^-,\CZ}$$
gives rise to a map
$$\MOp^{\mer,\mf_{\cP^-},\on{glob}}_{\cG,\cP^-,\CZ}\to \MOp^{\mer,\mf_{\cP^-}}_{\cG,\cP^-,\CZ},$$
where 
$$\MOp^{\mer,\mf_{\cP^-}}_{\cG,\cP^-,\CZ}:=\MOp^\mer_{\cG,\cP^-,\CZ}\underset{\LS^\mer_{\cP^-,\CZ}}\times \LS^\reg_{\cP^-,\CZ},$$
see \secref{sss:Miura mf}. Denote resulting map by the same symbol $\on{ev}^{\MOp}_{\CZ}$.

\medskip

Then the horizontal composition in \eqref{e:Op nat trans 4} equals
$$\MOp^{\mer,\mf_{\cP^-},\on{glob}}_{\cG,\cP^-,\CZ}\overset{\on{ev}^{\MOp}_{\CZ}}
\longrightarrow 
\MOp^{\mer,\mf_{\cP^-}}_{\cG,\cP^-,\CZ}\overset{\sfp^{\on{Miura},\mf}}\longrightarrow
\Op^\mf_{\cG,\CZ},$$
where $\sfp^{\on{Miura},\mf}$ is as in \secref{sss:Miura map mf}. 

\medskip

Hence, we can replace \eqref{e:Op nat trans 4} by 
\begin{equation} \label{e:Op nat trans 5}
\CD
& & & & \LS_\cM \times \CZ \\ 
& & & & @AA{\fr^{\on{glob}}_\cM}A \\
& & & & \Op^{\mf,\on{glob}}_{\cM,\rhoch_P,\CZ} \\
& & & & @AA{\sfq^{\on{Miura,glob}}}A \\
\Op^\mf_{\cG,\CZ} @<<{\sfp^{\on{Miura},\mf}}<  \MOp^{\mer,\mf_{\cP^-}}_{\cG,\cP^-,\CZ} @<<{\on{ev}^{\MOp}_{\CZ}}< 
\MOp^{\mer,\mf_{\cP^-},\on{glob}}_{\cG,\cP^-,\CZ}. 
\endCD
\end{equation} 

\sssec{} \label{sss:CT Op last}

Note now that we have a \emph{Cartesian} square
$$
\CD
\Op^{\mf}_{\cM,\rhoch_P,\CZ} @<{\on{ev}^{\Op_\cM}_{\CZ}}<<   \Op^{\mf,\on{glob}}_{\cM,\rhoch_P,\CZ}  \\
@A{\sfq^{\on{Miura},\mf}}AA @AA{\sfq^{\on{Miura,glob}}}A \\
\MOp^{\mer,\mf_{\cP^-}}_{\cG,\cP^-,\CZ} @<<{\on{ev}^{\MOp}_{\CZ}}< \MOp^{\mer,\mf_{\cP^-},\on{glob}}_{\cG,\cP^-,\CZ}.
\endCD
$$

Hence, we can replace \eqref{e:Op nat trans 5} by 
\begin{equation} \label{e:Op nat trans 6}
\CD
& & & & \LS_\cM \times \CZ \\ 
& & & & @AA{\fr^{\on{glob}}_\cM}A \\
& & \Op^{\mf}_{\cM,\rhoch_P,\CZ}  @<<{\on{ev}^{\Op_\cM}_{\CZ}}<   \Op^{\mf,\on{glob}}_{\cM,\rhoch_P,\CZ} \\
& & @A{\sfq^{\on{Miura},\mf}}AA \\
\Op^\mf_{\cG,\CZ} @<<{\sfp^{\on{Miura},\mf}}<  \MOp^{\mer,\mf_{\cP^-}}_{\cG,\cP^-,\CZ}. 
\endCD
\end{equation} 

Finally, the functor defined by \eqref{e:Op nat trans 6} is by definition the left-hand side in \eqref{e:Op nat trans x}. 

\sssec{}

Thus, we have completed the construction of the natural transformation \eqref{e:Op nat trans}.

\medskip

By \cite[Sect. 11.4.8]{GLC2}, the natural transformation \eqref{e:Op nat trans} gives rise to a natural transformation
\begin{equation} \label{e:Op nat trans int}
\left(\ul{\on{Poinc}}^{\on{spec}}_{\cM,!,\rhoch_P}\circ J_\Op^{-,\on{spec}}\right)^{\int, \on{ins.unit}}
\to \on{CT}^{-,\on{spec}}\circ \ul{\on{Poinc}}^{\on{spec}}_{\cG,!}.
\end{equation} 

\sssec{}

The main result of this section reads:

\begin{thm} \label{t:CT opers}
The natural transformation \eqref{e:Op nat trans int} is an isomorphism.
\end{thm} 

\ssec{Proof of \thmref{t:CT opers}}

\sssec{}

Fix $\CZ\to \Ran$. We wish to show that the natural transformation from the counterclockwise to
the clockwise circuit in the diagram
\begin{equation} \label{e:Op nat trans int 1}
\vcenter
{\xy
(0,0)*+{\IndCoh(\LS_\cG)\otimes \Dmod(\CZ)}="A";
(80,0)*+{\IndCoh(\LS_\cM)\otimes \Dmod(\CZ)}="B";
(0,-20)*+{\IndCoh^!(\Op^\mf_\cG)_\CZ}="C";
(80,-20)*+{\IndCoh(\LS_\cM)\otimes \Dmod(\CZ^\subseteq)}="D";
(80,-40)*+{\IndCoh^!(\Op^\mf_{\cM,\rhoch_P})_{\CZ^\subseteq}}="E";
(0,-40)*+{\IndCoh^!(\Op^\mf_\cG)_{\CZ^\subseteq}}="F";
{\ar@{->}^{\on{CT}^{-,\on{spec}}\otimes \on{Id}} "A";"B"};
{\ar@{->}^{\on{Poinc}^{\on{spec}}_{\cG,!,\CZ}} "C";"A"};
{\ar@{->}^{\on{ins.unit}_\CZ} "C";"F"};
{\ar@{->}^{J_\Op^{-,\on{spec}}} "F";"E"};
{\ar@{->}_{\on{Poinc}^{\on{spec}}_{\cM,!,\rhoch_P,\CZ^\subseteq}} "E";"D"};
{\ar@{->}_{\on{Id}\otimes (\on{pr}_{\on{small},\CZ})_!} "D";"B"};
\endxy} 
\end{equation}
arising from \eqref{e:Op nat trans} is an isomorphism.


\sssec{}

Let $\Op_{\cG,\CZ^\subseteq}^{\mf\rightsquigarrow\reg}$ be as in \cite[Sect. 15.2.5]{GLC2}. I.e., it is the space that associates
to a point of $\CZ^\subseteq$ mapping to $(\ul{x}\subseteq \ul{x}')\in \Ran^{\subseteq}$ the space of opers on $\cD_{\ul{x}'}-\ul{x}$
that are monodromy-free around $\ul{x}$, i.e., 
$$\Op_\cG(\cD_{\ul{x}'}-\ul{x})\underset{\LS_\cG(\cD_{\ul{x}'}-\ul{x})}\times \LS_\cG(\cD_{\ul{x}'}).$$

\medskip

Consider the corresponding maps
$$\Op^\mf_{\cG,\CZ}\overset{\on{pr}^{\Op_\cG}_{\on{small},\CZ}}\longleftarrow 
\Op_{\cG,\CZ^\subseteq}^{\mf\rightsquigarrow\reg} \overset{\on{pr}^{\Op_\cG}_{\on{big},\CZ}}\longrightarrow 
\Op^\mf_{\cG,\CZ^\subseteq}.$$

Note that the map $\on{pr}^{\Op_\cG}_{\on{big},\CZ}$ is an ind-closed embedding. 

\sssec{}

Recall that the functor 
$$\on{ins.unit}_\CZ:\IndCoh^!(\Op^\mf_\cG)_\CZ\to \IndCoh^!(\Op^\mf_\cG)_{\CZ^\subseteq}$$
is given by pull-push along the diagram
$$
\CD
& & \Op^\mf_{\cG,\CZ^\subseteq} \\
& & @AA{\on{pr}^{\Op_\cG}_{\on{big},\CZ}}A \\
\Op^\mf_{\cG,\CZ} @<<{\on{pr}^{\Op_\cG}_{\on{small},\CZ}}< \Op_{\cG,\CZ^\subseteq}^{\mf\rightsquigarrow\reg}.
\endCD
$$

\sssec{}

By Sects. \ref{sss:CT Op first}-\ref{sss:CT Op last}, the counterclockwise 
composition in \eqref{e:Op nat trans int 1} is given by pull-push along the following diagram

\medskip

\begin{equation} \label{e:Op nat trans int 2}
\CD
& & & & & & \LS_\cM\times \CZ \\
& & & & & & @AA{\on{Id}\otimes (\on{pr}_{\on{small},\CZ})_!}A \\
& & & & & & \LS_\cM\times \CZ^\subseteq \\
& & & & & & @AAA \\
& & \Op^\mf_{\cG,\CZ^\subseteq} @<<{\sfp^{\on{Miura},\mf}}<  
\MOp^{\mer,\mf_{\cP^-}}_{\cG,\cP^-,\CZ^\subseteq} @<<{\on{ev}^{\MOp}_{\CZ}}<  \MOp^{\mer,\mf_{\cP^-},\on{glob}}_{\cG,\cP^-,\CZ^\subseteq} \\
& & @AA{\on{pr}^{\Op_\cG}_{\on{big},\CZ}}A \\
\Op^\mf_{\cG,\CZ} @<{\on{pr}^{\Op_\cG}_{\on{small},\CZ}}<< \Op_{\cG,\CZ^\subseteq}^{\mf\rightsquigarrow\reg}. 
\endCD
\end{equation}

\medskip

Note that the fiber product
$$\Op_{\cG,\CZ^\subseteq}^{\mf\rightsquigarrow\reg}\underset{ \Op^\mf_{\cG,\CZ^\subseteq}}\times 
\MOp^{\mer,\mf_{\cP^-},\on{glob}}_{\cG,\cP^-,\CZ^\subseteq}$$
identifies with
\begin{multline} \label{e:MOp tilde}
\wt\MOp^{\mer,\mf_{\cP^-},\on{glob}}_{\cG,\cP^-,\CZ^\subseteq}
:= \\
=\left(\Op^{\mf,\on{glob}}_{\cG,\CZ}\underset{\CZ,\on{pr}_{\on{small},\CZ}}\times \CZ^\subseteq\right) 
\underset{\Op^{\mf,\on{glob}}_{\cG,\CZ^\subseteq}}\times \MOp^{\mer,\mf_{\cP^-},\on{glob}}_{\cG,\cP^-,\CZ^\subseteq},
\end{multline}
where 
$$\Op^{\mf,\on{glob}}_{\cG,\CZ}\underset{\CZ,\on{pr}_{\on{small},\CZ}}\times \CZ^\subseteq\to 
\Op^{\mf,\on{glob}}_{\cG,\CZ^\subseteq}$$
is the natural restriction map.

\medskip

Hence, we can replace the diagram \eqref{e:Op nat trans int 2} by 

\begin{equation} \label{e:Op nat trans int 3}
\CD
& & & & & & \LS_\cM\times \CZ \\
& & & & & & @AA{\on{Id}\otimes (\on{pr}_{\on{small},\CZ})_!}A \\
& & & & & & \LS_\cM\times \CZ^\subseteq \\
& & & & & & @AAA \\
\Op^\mf_{\cG,\CZ} @<{\on{ev}^{\Op_\cG}_\CZ}<< \Op^{\mf,\on{glob}}_{\cG,\CZ} 
@<<< \Op^{\mf,\on{glob}}_{\cG,\CZ}\underset{\CZ,\on{pr}_{\on{small},\CZ}}\times \CZ^\subseteq @<<< 
\wt\MOp^{\mer,\mf_{\cP^-},\on{glob}}_{\cG,\cP^-,\CZ^\subseteq}.
\endCD
\end{equation}

\sssec{}

Note that we have a naturally defined map 
$$\wt\jmath^{\on{glob}}:\wt\MOp^{\mer,\mf_{\cP^-},\on{glob}}_{\cG,\cP^-,\CZ^\subseteq}\to
\Op^{\mf,\on{glob}}_{\cG,\CZ}\underset{\LS_\cG\times \CZ}\times (\LS_{\cP^-}\times \CZ),$$
such that:

\begin{itemize}

\item The horizontal composition in \eqref{e:Op nat trans int 3} equals the composition of
$\wt\jmath^{\on{glob}}$ with the horizontal composition in \eqref{e:Op nat trans 2};

\medskip

\item The vertical composition in \eqref{e:Op nat trans int 3} equals the composition of
$\wt\jmath^{\on{glob}}$ with the vertical composition in \eqref{e:Op nat trans 2}.

\end{itemize}

It follows formally from \propref{p:j glob} that the map $\wt\jmath^{\on{glob}}$ is pseudo-proper.

\medskip

Unwinding, we obtain that the natural transformation in \eqref{e:Op nat trans int 1} equals
\begin{multline*}
(\sfq^{-,\on{spec}})^\IndCoh_*\circ  
({}'\fr^{\on{glob}}_\cG)^\IndCoh_*\circ (\wt\jmath^{\on{glob}})^\IndCoh_*\circ 
(\wt\jmath^{\on{glob}})^!\circ ({}'\sfp^{-,\on{spec}})^!\circ (\on{ev}^{\Op_\cG}_{\CZ})^!\to \\
\to (\sfq^{-,\on{spec}})^\IndCoh_*\circ  
({}'\fr^{\on{glob}}_\cG)^\IndCoh_*\circ ({}'\sfp^{-,\on{spec}})^!\circ (\on{ev}^{\Op_\cG}_{\CZ})^!,
\end{multline*} 
induced by the counit of the adjunction
\begin{equation} \label{e:p contr}
(\wt\jmath^{\on{glob}})^\IndCoh_*\circ (\wt\jmath^{\on{glob}})^!\to \on{Id}. 
\end{equation}

\sssec{}

We will prove:

\begin{prop} \label{p:j contr}
The natural transformation \eqref{e:p contr} is an isomorphism.
\end{prop}

The proof will be given in \secref{sss:j contr}.

\sssec{}

From \propref{p:j contr}, we obtain that 
the natural transformation in \eqref{e:Op nat trans int 1} is an isomorphism as well, as desired.

\qed[\thmref{t:CT opers}]

\ssec{How do Miura opers appear?} \label{ss:j glob}

\sssec{}

In \secref{sss:Zas} we introduced the Zastava space, $\Zas$. Here we will use a variant of this construction, defined as follows. 

\medskip

For $\CZ\to \Ran$, let $\Bun_{\cG,\CZ}^{\mer,\on{glob}}$ (resp., $\Bun_{\cB,\CZ}^{\mer,\on{glob}}$)
be the space of $\cG$- (resp., $\cB$-) bundles on $X$ defined away from the marked point. 

\medskip

Consider the fiber product
$$\Bun_{\cB,\CZ}^{\mer,\on{glob}}\underset{\Bun_{\cG,\CZ}^{\mer,\on{glob}}}\times (\Bun_{\cP^-}\times \CZ).$$

\medskip

Let $\on{Zast}_\CZ$ be its open subfunctor corresponding to the condition that the 
$\cB$-reduction and the $\cP^-$-reduction are transversal at the
generic point of the curve.

\sssec{}

We have a tautological forgetful map
\begin{equation} \label{e:Op to Zast}
\Op^{\mf,\on{glob}}_{\cG,\CZ}\underset{\LS_\cG\times \CZ}\times (\LS_{\cP^-}\times \CZ)\to 
\Bun_{\cB,\CZ}^{\mer,\on{glob}}\underset{\Bun_{\cG,\CZ}^{\mer,\on{glob}}}\times (\Bun_{\cP^-}\times \CZ).
\end{equation}

The following observation is the reason Miura opers make an appearance: 

\begin{prop} \label{p:Op to Zast}
The image of \eqref{e:Op to Zast} lies in $\on{Zast}_\CZ$.
\end{prop}

The proof will be given in \secref{ss:Op to Zast}. 

\sssec{}

Let 
$$\on{Zast}^\circ_{\CZ^\subseteq}\subset \on{Zast}_\CZ\underset{\CZ,\on{pr}_{\on{small},\CZ}}\times \CZ^\subseteq$$ 
be the subfunctor given by the following condition: for $z\in \CZ^\subseteq$ mapping
to $(\ul{x}\subseteq \ul{x}')\in \Ran^\subseteq$, we require that the $\cB$-reduction and the
$\cP^-$-reduction be transversal on $X-\ul{x}'$. 

\medskip

We have a forgetful map
$$\wt\jmath_{\on{Zast}}:\on{Zast}^\circ_{\CZ^\subseteq}\to \on{Zast}_\CZ,$$
which fits into the \emph{Cartesian} square
\begin{equation} \label{e:j Cart}
\CD
\wt\MOp^{\mer,\mf_{\cP^-},\on{glob}}_{\cG,\cP^-,\CZ^\subseteq} @>>>  \on{Zast}^\circ_{\CZ^\subseteq} \\
@V{\wt\jmath^{\on{glob}}}VV @VV{\wt\jmath_{\on{Zast}}}V \\
\Op^{\mf,\on{glob}}_{\cG,\CZ}\underset{\LS_\cG\times \CZ}\times (\LS_{\cP^-}\times \CZ) @>>> \on{Zast}_\CZ,
\endCD
\end{equation} 
where $\wt\MOp^{\mer,\mf_{\cP^-},\on{glob}}_{\cG,\cP^-,\CZ^\subseteq}$ is as in \eqref{e:MOp tilde}.

\sssec{}

We will prove:

\begin{lem} \label{l:J Zast} \hfill 

\smallskip

\noindent{\em(a)} The map $\wt\jmath_{\on{Zast}}$ is pseudo-proper;

\smallskip

\noindent{\em(b)} The fibers of $\wt\jmath_{\on{Zast}}$ over geometric points of $\on{Zast}_\CZ$
are $\CO$-contractible.

\end{lem}

Let us show how \lemref{l:J Zast} implies Propositions \ref{p:j glob} and \ref{p:j contr}.

\sssec{Proof of \propref{p:j glob}} \label{sss:j glob}

The map $\jmath^{\on{glob}}$ is ind-schematic (and injective). Therefore, in order to show that
it is ind-proper (in fact, an ind-closed embedding), it suffices to show that it is pseudo-proper.

\medskip

By \lemref{l:J Zast}(a) and \eqref{e:j Cart}, the map $\wt\jmath^{\on{glob}}$ is pseudo-proper.
Now, the map $\jmath^{\on{glob}}$ is the composition
\begin{multline*}
\MOp^{\mer,\mf_{\cP^-},\on{glob}}_{\cG,\cP^-,\CZ}\simeq
\CZ\underset{\CZ^\subseteq}\times \wt\MOp^{\mer,\mf_{\cP^-},\on{glob}}_{\cG,\cP^-,\CZ^\subseteq} \to \\
\to \wt\MOp^{\mer,\mf_{\cP^-},\on{glob}}_{\cG,\cP^-,\CZ^\subseteq} 
\overset{\wt\jmath^{\on{glob}}}\longrightarrow
\Op^{\mf,\on{glob}}_{\cG,\CZ}\underset{\LS_\cG}\times (\LS_{\cP^-}\times \CZ),
\end{multline*}
where the second arrow is the base change
of 
$$\CZ\overset{\on{diag}_\CZ}\to \CZ^\subseteq$$
and thus is a closed embedding. 

\qed[\propref{p:j glob}]

\sssec{Proof of \propref{p:j contr}} \label{sss:j contr}

As was remarked above, the map $\wt\jmath^{\on{glob}}$ is pseudo-proper; hence, 
the functor $(\wt\jmath^{\on{glob}})^\IndCoh_*$ satisfies base change. Hence, in order to show
that \eqref{e:p contr} is an isomorphism, it is enough to do so at the level of fibers over field-valued 
points. 

\medskip

By \eqref{e:j Cart}, the assertion follows now from \lemref{l:J Zast}(b).

\qed[\propref{p:j contr}]

\sssec{Proof of \lemref{l:J Zast}(a)}

It suffices to show that the inclusion
$$\on{Zast}^\circ_{\CZ^\subseteq}\hookrightarrow \on{Zast}_\CZ\underset{\CZ,\on{pr}_{\on{small},\CZ}}\times \CZ^\subseteq$$ 
is an ind-closed embedding.

\medskip

Let us be  given $\ul{x}\in X$ and

\smallskip

\begin{itemize}

\item A $\cB$-bundle $\CP_\cB$ on $X-\ul{x}$;

\medskip

\item A $\cP^-$-bundle $\CP_{\cP^-}$ on $X-\ul{x}$;

\item An isomorphism $\cG\overset{\cB}\times \CP_\cB\overset{\beta}\simeq \cG\overset{\cP^-}\times \CP_{\cP^-}$
over $X-\ul{x}$ that is \emph{generically transversal}.

\end{itemize}

\medskip

We have to show that the locus of $\ul{x}'\supseteq \ul{x}$ such that $\beta$ be transversal over $X-\ul{x}'$
is an ind-closed subfunctor of $\Ran_{\ul{x}}$. 

\medskip

This reduces to the following situation: let $\gamma:\CL_1\to \CL_2$ be a \emph{non-zero} map between line bundles
on $X$ with poles at $\ul{x}$. Then the locus of $\ul{x}'\supseteq \ul{x}$ such that $\gamma$ is an isomorphism on
$X-\ul{x}'$ is an ind-closed subfunctor of $\Ran_{\ul{x}}$. 

\medskip

Up to modifying $\CL_2$ at $\ul{x}$, we can assume that $\gamma$ is a regular map. In this case the above locus in
$\Ran_{\ul{x}}$ is the support of $\on{coker}(\gamma)$. 

\qed[\lemref{l:J Zast}(a)]

\sssec{Proof of \lemref{l:J Zast}(b)}

Fix a geometric point $(\ul{x},\CP_\cB,\CP_{\cP^-},\beta)\in \on{Zast}_\CZ$ as above. Let $U\subset X-\ul{x}$
be the locus over which $\beta$ is transversal. Write $U=X-\ul{y}$. 

\medskip

Then the fiber of $\wt\jmath^{\on{Zast}}$ over the above point is $\Ran_{\ul{y}}$. Now the assertion follows
from the homological contractibility of $\Ran_{\ul{y}}$, since
$$\Ran_{\ul{y}}\simeq (\Ran_{\ul{y}})_\dr.$$

\qed[\lemref{l:J Zast}(b)]

\ssec{Opers force transversality} \label{ss:Op to Zast}

In this subsection we prove \propref{p:Op to Zast}. 

\sssec{}

The assertion of \propref{p:Op to Zast} can be checked at the level of geometric points. Let 
$$(\CP_\cG,\nabla,\CP_\cB,\CP_{\cP^-})$$ 
be an oper equipped with a reduction to $\cP^-$ of the underlying local system over an open curve $U\subset X$.

\medskip

Assume by contradiction that $\CP_\cB$ and $\CP_{\cP^-}$ are not generically transversal. Then, up to shrinking $U$,
we can assume that the relative position of $\CP_\cB$ and $\CP_{\cP^-}$ corresponds to a non-unit element 
$w\in W/W_M$, where $W$ is the Weyl group. 

\sssec{}

Consider the fiber square of vector bundles on $U$:
$$
\CD
\cg_{\CP_\cG}/\cb_{\CP_\cB}\cap \cp^-_{\CP_{\cP^-}} @>>> \cg_{\CP_\cG}/\cb_{\CP_\cB} \\
@VVV @VVV \\
\cg_{\CP_\cG}/\cp^-_{\CP_{\cP^-}} @>>> \cg_{\CP_\cG}/\cb_{\CP_\cB}+ \cp^-_{\CP_{\cP^-}}. 
\endCD
$$

The data of the connection $\nabla$ and that of reductions of $\CP_\cG$ to $\CP_\cB$ and $\CP_{\cP^-}$
gives rise to a section
\begin{equation} \label{e:oper mod P}
\nabla\,\on{mod}\,\cb_{\CP_\cB}\cap \cp^-_{\CP_{\cP^-}}
\end{equation} 
of $(\cg_{\CP_\cG}/\cb_{\CP_\cB}\cap \cp^-_{\CP_{\cP^-}})\otimes \omega_X$.

\medskip

By assumption, the image of \eqref{e:oper mod P} in $(\cg_{\CP_\cG}/\cp^-_{\CP_{\cP^-}})\otimes \omega_X$ vanishes, and hence, so does its
further image on $(\cg_{\CP_\cG}/\cb_{\CP_\cB}+ \cp^-_{\CP_{\cP^-}})\otimes \omega_X$.

\sssec{}

The image of \eqref{e:oper mod P} in $(\cg_{\CP_\cG}/\cb_{\CP_\cB})\otimes \omega_X$ is given by
the oper condition. 

\medskip

However, it is easy to see that for $w\neq 1$ the map 
$$(\cg_{\CP_\cG}/\cb_{\CP_\cB})\otimes \omega_X\to (\cg_{\CP_\cG}/\cb_{\CP_\cB}+ \cp^-_{\CP_{\cP^-}})\otimes \omega_X$$
is non-zero on oper elements. 

\qed[\propref{p:Op to Zast}]

\ssec{The enhanced version of spectral Poincar\'e vs constant term compatibility}

\sssec{} \label{sss:Poinc spec enh}

Since the functor 
$$\on{Poinc}^{\on{spec}}_{\cG,!,\CZ}:\IndCoh^!(\Op^\mf_{\cG,\CZ})\to \IndCoh(\LS_\cG)\otimes \Dmod(\CZ)$$
respects the actions of $\Sph^{\on{spec}}_{\cG,\CZ}$, it induces a functor
$$\on{Poinc}^{-,\on{spec,enh}_{\on{co}}}_{\cG,!,\CZ}:\IndCoh^!(\Op_\cG^{\on{mon-free}})^{-,\on{enh}_{\on{co}}}_\CZ\to
\IndCoh(\LS_\cG)^{-,\on{enh}_{\on{co}}}_\CZ,$$
where:

\begin{itemize}

\item $\IndCoh^!(\Op_\cG^{\on{mon-free}})^{-,\on{enh}_{\on{co}}}$ is as in \secref{sss:J Op ult}. 

\medskip

\item $\IndCoh(\LS_\cG)^{-,\on{enh}_{\on{co}}}_\CZ$ is as in \secref{sss:spec co-enh}. 

\end{itemize} 

Denote by 
$$\ul{\on{Poinc}}^{-,\on{spec,enh}_{\on{co}}}_{\cG,!}:\ul\IndCoh^!(\Op_\cG^{\on{mon-free}})^{-,\on{enh}_{\on{co}}}\to
\ul\IndCoh(\LS_\cG)^{-,\on{enh}_{\on{co}}}$$
the resulting functor between crystals of categories.

\medskip

The unital structure on $\ul{\on{Poinc}}^{-,\on{spec}}_{\cG,!}$ induces one on $\ul{\on{Poinc}}^{-,\on{spec,enh}_{\on{co}}}_{\cG,!}$. 

\sssec{}

Let
$$\ul{\on{CT}}^{-,\on{spec,enh}_{\on{co}}}:\ul\IndCoh(\LS_\cG)^{-,\on{enh}_{\on{co}}}\to
\IndCoh(\LS_\cM)\otimes \ul\Dmod(\Ran)$$
be as in \secref{ss:CT spec enh}.

\medskip

Let 
$$J^{-,\on{spec},\on{enh}_{\on{co}}}_\Op:\IndCoh^!(\Op_\cG^{\on{mon-free}})^{-,\on{enh}_{\on{co}}}\to
\IndCoh^!(\Op^\mf_{\cM,\rhoch_P})$$
be as in \secref{sss:J Op ult}. 

\medskip

Our current goal is to construct a natural transformation
\begin{equation} \label{e:Op nat trans enh}
\ul{\on{Poinc}}^{\on{spec}}_{\cM,!,\rhoch_P}\circ \ul{J}_\Op^{-,\on{spec,enh}_{\on{co}}}\to 
\ul{\on{CT}}^{-,\on{spec,enh}_{\on{co}}}\circ \ul{\on{Poinc}}^{-,\on{spec,enh}_{\on{co}}}_{\cG,!}
\end{equation} 
as local-to-global functors 
$$\ul\IndCoh^!(\Op_\cG^{\on{mon-free}})^{-,\on{enh}_{\on{co}}}\rightrightarrows \IndCoh(\LS_\cM)\otimes \ul\Dmod(\Ran),$$
compatible with the lax unital structure. 

\sssec{}

Given $\CZ$, we wish to construct a natural transformation
\begin{equation} \label{e:Op nat trans enh Z}
\on{Poinc}^{\on{spec}}_{\cM,!,\rhoch_P,\CZ}\circ J_\Op^{-,\on{spec,enh}_{\on{co}}}\to 
\on{CT}^{-,\on{spec,enh}_{\on{co}}}_\CZ\circ \on{Poinc}^{-,\on{spec,enh}_{\on{co}}}_{\cG,!,\CZ}
\end{equation} 
as functors
$$\IndCoh^!(\Op_\cG^{\on{mon-free}})^{-,\on{enh}_{\on{co}}}_\CZ\rightrightarrows \IndCoh(\LS_\cM)\otimes \Dmod(\CZ).$$

\medskip

By construction, the right-hand side in \eqref{e:Op nat trans enh Z} is given by pull-push along the diagram

\medskip

$$
\CD
& & & & \LS_\cM\times \CZ \\
& & & & @AA{\sfq^{-,\on{spec}}_\CZ}A \\
& & \left(\LS^{\on{mer,glob}}_{\cG,\CZ}\right){}^\wedge_\mf\underset{\LS^\mer_{\cG,\CZ}}\times \LS^{\mf_\cM}_{\cP^-,\CZ} @<<{\sfp^{-,\on{spec,enh}}_\CZ}<
\left(\LS^{\mf_\cM,\on{glob}}_{\cP^-,\CZ}\right){}^\wedge_\mf \\
& & @AAA \\
(\Op^\Mmf_{\cG,\cP^-,\CZ})^\wedge_\mf @<<< (\Op^{\Mmf,\on{glob}}_{\cG,\cP^-,\CZ})^\wedge_\mf,
\endCD
$$
where:

\begin{itemize}

\item $(\Op^\Mmf_{\cG,\cP^-,\CZ})^\wedge_\mf$ is as in \secref{ss:form coml Op GP}, i.e., is the formal completion of
$$\Op^\mer_{\cG,\CZ}\underset{\LS^\mer_{\cG,\CZ}}\times \LS^{\mf_\cM}_{\cP^-,\CZ}$$ along the natural map from 
$\Op^\mer_{\cG,\CZ}\underset{\LS^\mer_{\cG,\CZ}}\times \LS^\reg_{\cP^-,\CZ}$;

\item $(\Op^{\Mmf,\on{glob}}_{\cG,\cP^-,\CZ})^\wedge_\mf$ denotes the formal completion of
$$\Op^{\mer,\on{glob}}_{\cG,\CZ}\underset{\LS^\mer_{\cG,\CZ}}\times \LS^{\mf_\cM}_{\cP^-,\CZ}$$ along the natural map from
$\Op^{\mer,\on{glob}}_{\cG,\CZ}\underset{\LS^\mer_{\cG,\CZ}}\times \LS^\reg_{\cP^-,\CZ}$.

\item The left horizontal arrow is induced by the restriction map $\Op^{\mer,\on{glob}}_{\cG,\CZ}\to \Op^{\mer}_{\cG,\CZ}$;

\medskip

\item The middle vertical arrow is induced by the map $\Op^{\mer,\on{glob}}_{\cG,\CZ}\to \LS^{\mer,\on{glob}}_{\cG,\CZ}$.

\end{itemize}

\sssec{}

As in \secref{sss:crucial nat trans}, the above functor receives a natural transformation from the functor given by pull-push along the diagram
\begin{equation} \label{e:Op nat trans enh 1}
\CD
& & & & \LS_\cM\times \CZ \\
& & & & @AA{\sfq^{-,\on{spec}}_\CZ}A \\
& & & & \left(\LS^{\mf_\cM,\on{glob}}_{\cP^-,\ul{x}}\right){}^\wedge_\mf \\
& & & & @AAA \\
(\Op^\Mmf_{\cG,\cP^-,\CZ})^\wedge_\mf @<<< (\Op^{\Mmf,\on{glob}}_{\cG,\cP^-,\CZ})^\wedge_\mf @<<< 
(\MOp^{\Mmf,\on{glob}}_{\cG,\cP^-,\CZ})^\wedge_\mf,
\endCD
\end{equation} 

\medskip

\noindent where $(\MOp^{\Mmf,\on{glob}}_{\cG,\cP^-,\CZ})^\wedge_\mf$ is the completion of
$$\MOp^{\mer,\on{glob}}_{\cG,\cP^-,\CZ}\underset{\LS^\mer_{\cM,\CZ}}\times \LS^\reg_{\cM,\CZ}$$
along the natural map from 
$$\MOp^{\mer,\mf_{\cP^-},\on{glob}}_{\cG,\cP^-,\CZ}:=
\MOp^{\mer,\on{glob}}_{\cG,\cP^-,\CZ}\underset{\LS^{\mer,\on{glob}}_{\cP^-,\CZ}}\times (\LS_{\cP^-}\times \CZ)\simeq
\MOp^{\mer,\on{glob}}_{\cG,\cP^-,\CZ}\underset{\LS^\mer_{\cP^-,\CZ}}\times \LS^\reg_{\cP^-,\CZ}.$$

\sssec{}

We now show that the latter functor equals the left-hand side in \eqref{e:Op nat trans enh Z}. Indeed, 
we can replace the diagram \eqref{e:Op nat trans enh 1} by 
\begin{equation} \label{e:Op nat trans enh 2}
\CD
& & & & \LS_\cM\times \CZ \\
& & & & @AAA \\
& & & & \Op^{\mf,\on{glob}}_{\cM,\rhoch_P,\CZ} \\ 
& & & & @AAA \\
(\Op^\Mmf_{\cG,\cP^-,\CZ})^\wedge_\mf @<<< (\MOp^\Mmf_{\cG,\cP^-,\CZ})^\wedge_\mf @<<< (\MOp^{\Mmf,\on{glob}}_{\cG,\cP^-,\CZ})^\wedge_\mf,
\endCD
\end{equation} 

\medskip

\noindent where $(\MOp^\Mmf_{\cG,\cP^-,\CZ})^\wedge_\mf$ is as \secref{sss:form compl Miura}, i.e., is the completion of
$$\MOp^{\mer}_{\cG,\cP^-,\CZ}\underset{\LS^\mer_{\cM,\CZ}}\times \LS^\reg_{\cM,\CZ}$$
along the natural map from 
$$\MOp^{\mer}_{\cG,\cP^-,\CZ}\underset{\LS^\mer_{\cP^-,\CZ}}\times \LS^\reg_{\cP^-,\CZ}.$$

\sssec{}

We note now that we have a Cartesian square 
$$
\CD
\Op^\mf_{\cM,\rhoch_P,\CZ}  @<<< \Op^{\mf,\on{glob}}_{\cM,\rhoch_P,\CZ}  \\
@AAA @AAA \\
(\MOp^\Mmf_{\cG,\cP^-,\CZ})^\wedge_\mf @<<< (\MOp^{\Mmf,\on{glob}}_{\cG,\cP^-,\CZ})^\wedge_\mf.
\endCD
$$

This diagram satisfies base change for !-pullbacks and *-pushforwards. Hence, we can further replace \eqref{e:Op nat trans enh 2}
by 
$$
\CD
& & & & \LS_\cM\times \CZ \\
& & & & @AAA \\
& & \Op^\mf_{\cM,\rhoch_P,\CZ}  @<<< \Op^{\mf,\on{glob}}_{\cM,\rhoch_P,\CZ} \\ 
& & @AAA \\
(\MOp^\Mmf_{\cG,\cP^-,\CZ})^\wedge_\mf @<<< (\MOp^{\Mmf,\on{glob}}_{\cG,\cP^-,\CZ})^\wedge_\mf.
\endCD
$$

The functor define by the latter diagram is by definition the left-hand side in \eqref{e:Op nat trans enh Z}.

\sssec{} \label{sss:ind enh op}

Thus, we have completed the construction of the natural transformation \eqref{e:Op nat trans enh}.

\medskip

By construction, if we precompose \eqref{e:Op nat trans enh} on both sides with the functor
\begin{multline}  \label{e:ind enh op}
\IndCoh^!(\Op^\mf_\cG)  \overset{\on{Id}\otimes \one_{\on{I}(\cG,\cP^-)^{\on{spec,loc}}_{\on{co}}}}\longrightarrow 
\IndCoh^!(\Op^\mf_\cG) \otimes \on{I}(\cG,\cP^-)^{\on{spec,loc}}_{\on{co}}\to \\
\to \IndCoh^!(\Op_\cG^\mf)\overset{(\Sph^{\on{spec}}_\cG)^\vee}\otimes \on{I}(\cG,\cP^-)^{\on{spec,loc}}_{\on{co}} =
\IndCoh^!(\Op_\cG^{\on{mon-free}})^{-,\on{enh}_{\on{co}}},
\end{multline}
we recover the natural transformation \eqref{e:Op nat trans}. 

\sssec{}

Since the right-hand side in \eqref{e:Op nat trans enh} is a (strictly) unital local-to-global functor, by \cite[Sect. 11.4.8]{GLC2}, 
from \eqref{e:Op nat trans enh} we obtain a natural transformation
\begin{equation} \label{e:Op nat trans enh int}
\left(\ul{\on{Poinc}}^{\on{spec}}_{\cM,!,\rhoch_P}\circ \ul{J}_\Op^{-,\on{spec,enh}_{\on{co}}}\right)^{\int\, \on{ins.unit}}\to 
\ul{\on{CT}}^{-,\on{spec,enh}_{\on{co}}}\circ \ul{\on{Poinc}}^{-,\on{spec,enh}_{\on{co}}}_{\cG,!}.
\end{equation} 

We claim:

\begin{thm} \label{t:CT opers enh}
The natural transformation \eqref{e:Op nat trans enh int} is an isomorphism.
\end{thm} 
 
\begin{proof} 

By construction, the categories involved carry an action of $\ul\Sph_\cM^{\on{spec}}$, and the two functors,
as well as the natural transformation \eqref{e:Op nat trans enh int} are compatible with the action.

\medskip

Hence, since $\IndCoh(\Op_\cG^{\on{mon-free}})^{-,\on{enh}_{\on{co}}}$ is generated under the action of $\Sph_\cM^{\on{spec}}$
by the image of \eqref{e:ind enh op}, it suffices to show that \eqref{e:Op nat trans enh int} becomes an isomorphism
after precomposing with \eqref{e:ind enh op}.

\medskip

However, by \secref{sss:ind enh op}, the resulting natural transformation identifies with \eqref{e:Op nat trans}. Hence,
the required assertion follows from \thmref{t:CT opers}.

\end{proof}

\ssec{Addendum: spectral global sections and the spectral Poincar\'e functor}

In this subsection we will prove a spectral counterpart of \thmref{t:coeff Loc}. 

\sssec{}

Let us denote by $\sP^{\on{spec,loc,true}}_\cG$ the (lax unital) factorization functor
\begin{multline*}
\IndCoh^*(\Op^\mf_\cG)\underset{\Sph^{\on{spec}}_\cG}\otimes \Rep(\cG) \simeq
\IndCoh^*(\Op^\mer_\cG)_\mf\hookrightarrow \\
\hookrightarrow \IndCoh^*(\Op^\mer_\cG) \overset{\Gamma^\IndCoh(\Op^\mer_\cG,-)}\longrightarrow \Vect,
\end{multline*}
where the first arrow is \cite[Equation (3.17)]{GLC2}. 

\medskip

Recall that the functor $\sP^{\on{spec,loc}}_\cG$ is defined to be the composition
$$\IndCoh^*(\Op^\mf_\cG)\otimes \Rep(\cG) \to 
\IndCoh^*(\Op^\mf_\cG)\underset{\Sph^{\on{spec}}_\cG}\otimes \Rep(\cG) \overset{\sP^{\on{spec,loc,true}}_\cG}
\longrightarrow \Vect,$$
see \cite[Sect. 6.4.7]{GLC2}. 

\sssec{} \label{sss:Poinc Gamma comp}

In \cite[Theorem 17.7.2]{GLC2} it was shown that the functor
\begin{multline}  \label{e:Poinc Gamma LHS}
(\IndCoh^*(\Op^\mf_\cG)\otimes \Rep(\cG))_\CZ
\overset{\on{Poinc}_{\cG,*,\CZ}\otimes \on{Id}}\longrightarrow 
(\IndCoh(\LS_\cG)\otimes \Dmod(\CZ))\underset{\Dmod(\CZ)} \otimes \Rep(\cG)_\CZ \to \\
\overset{\Gamma^{\on{spec},\IndCoh}_\CZ\otimes \on{Id}}\longrightarrow (\Rep(\cG)\otimes \Rep(\cG))_\CZ\to \Dmod(\CZ)
\end{multline}
is canonically isomorphic to
\begin{multline} \label{e:Poinc Gamma RHS}
(\IndCoh^*(\Op^\mf_\cG)\otimes \Rep(\cG))_\CZ
\overset{\on{ins.unit}_\CZ}\longrightarrow (\IndCoh^*(\Op^\mf_\cG)\otimes \Rep(\cG))_{\CZ^\subseteq}  \overset{\sP^{\on{spec,loc}}_G}\longrightarrow \\
\to \Dmod(\CZ^\subseteq)\overset{(\on{pr}_{\on{small},\CZ})_!}\longrightarrow \Dmod(\CZ).
\end{multline}

\sssec{}

Note now that both \eqref{e:Poinc Gamma LHS} and \eqref{e:Poinc Gamma RHS} factor via 
$$(\IndCoh^*(\Op^\mf_\cG)\otimes \Rep(\cG))_\CZ\to (\IndCoh^*(\Op^\mf_\cG)\underset{\Sph^{\on{spec}}_\cG}\otimes \Rep(\cG))_\CZ \to \Dmod(\CZ).$$

Indeed, for \eqref{e:Poinc Gamma LHS} this factorization is provided by
\begin{multline}  \label{e:Poinc Gamma LHS Sph}
(\IndCoh^*(\Op^\mf_\cG)\underset{\Sph^{\on{spec}}_\cG}\otimes \Rep(\cG))_\CZ\overset{\on{Poinc}_{\cG,*,\CZ}\otimes \on{Id}}\longrightarrow \\
\to (\IndCoh(\LS_\cG)\otimes \Dmod(\CZ))\underset{\Sph^{\on{spec}}_{\cG,\CZ}}\otimes \Rep(\cG)_\CZ 
\overset{\Gamma^{\on{spec},\IndCoh}_\CZ\otimes \on{Id}}\longrightarrow \\
\to (\Rep(\cG)\underset{\Sph^{\on{spec}}_\cG}\otimes \Rep(\cG))_\CZ\to \Dmod(\CZ). 
\end{multline} 

For \eqref{e:Poinc Gamma RHS} this factorization is provided by
\begin{multline} \label{e:Poinc Gamma RHS Sph}
(\IndCoh^*(\Op^\mf_\cG)\underset{\Sph^{\on{spec}}_\cG}\otimes \Rep(\cG))_\CZ
\overset{\on{ins.unit}_\CZ}\longrightarrow \\
\to (\IndCoh^*(\Op^\mf_\cG)\underset{\Sph^{\on{spec}}_\cG}\otimes \Rep(\cG))_{\CZ^\subseteq} 
\overset{\sP^{\on{spec,loc,true}}_\cG}\longrightarrow 
 \Dmod(\CZ^\subseteq)\overset{(\on{pr}_{\on{small},\CZ})_!}\longrightarrow \Dmod(\CZ).
\end{multline}

\sssec{} \label{sss:glob true spec}

Denote the functor \eqref{e:Poinc Gamma LHS} by $\sP^{\on{spec,glob}}_{\cG,\CZ}$
and the functor \eqref{e:Poinc Gamma LHS Sph} by $\sP^{\on{spec,glob,true}}_{\cG,\CZ}$.

\medskip

Note that we can think of the functor \eqref{e:Poinc Gamma RHS Sph} as
$$(\ul\sP^{\on{spec,loc,true}}_\cG)^{\int\, \on{ins.unit}}_\CZ.$$

\sssec{}

In this subsection we will prove the following sharpening of \cite[Theorem 14.1.4]{GLC2}:

\begin{thm} \label{t:Poinc Gamma}
There is a canonical isomorphism 
$$\sP^{\on{spec,glob,true}}_{\cG,\CZ}\simeq (\ul\sP^{\on{spec,loc,true}}_\cG)^{\int\, \on{ins.unit}}_\CZ$$
as local-to-global functors
$$\ul\IndCoh^*(\Op^\mf_\cG)\underset{\ul\Sph^{\on{spec}}_\cG}\otimes \ul\Rep(\cG)\to \ul\Dmod(\Ran).$$
\end{thm} 

The rest of this subsection is devoted to the proof of this theorem. 

\sssec{}

Fix $\CZ\to \Ran$. Unwinding the definitions, we can rewrite the functor $\sP^{\on{spec,glob,true}}_{\cG,\CZ}$
as *-pull followed by *-push along the following diagram
$$
\CD
& & \CZ \\
& & @AAA \\
& & (\LS^{\mer,\on{glob}}_{\cG,\CZ})^\wedge_\mf \\
& & @AAA \\
(\Op^\mer_{\cG,\CZ})^\wedge_\mf @<<< (\Op^{\mer,\on{glob}}_{\cG,\CZ})^\wedge_\mf,
\endCD
$$
while $\sP^{\on{spec,loc,true}}_{\cG,\CZ}$ is just the *-pushforward along
$$(\Op^\mer_{\cG,\CZ})^\wedge_\mf\to \CZ.$$

\sssec{}

This gives rise to a natural transformation
\begin{equation} \label{e:Poinc Gamma 1}
\ul\sP^{\on{spec,loc,true}}_\cG\to \ul\sP^{\on{spec,glob,true}}_\cG,
\end{equation}
as lax unital local-to-global functors
$$\ul\IndCoh^*(\Op^\mf_\cG)\underset{\ul\Sph^{\on{spec}}_\cG}\otimes \ul\Rep(\cG)\rightrightarrows \Vect,$$
where the right-hand side in \eqref{e:Poinc Gamma 1} is strictly unital.

\medskip

We claim that the induced natural transformation
\begin{equation} \label{e:Poinc Gamma 2}
(\ul\sP^{\on{spec,loc,true}}_\cG)^{\int\, \on{ins.unut}}\to \ul\sP^{\on{spec,glob,true}}_\cG
\end{equation}
is an isomorphism.

\sssec{}

Indeed, when we precompose the functors on the two sides of \eqref{e:Poinc Gamma 1} with 
$$\ul\IndCoh^*(\Op^\mf_\cG)\otimes \ul\Rep(\cG)\to 
\ul\IndCoh^*(\Op^\mf_\cG)\underset{\ul\Sph^{\on{spec}}_\cG}\otimes \ul\Rep(\cG),$$
we obtain a natural transformation 
$$\ul\sP^{\on{spec,loc}}_\cG\to \ul\sP^{\on{spec,glob}}_\cG,$$
(between lax unital functors, with the right-hand side striclty unital), 
such that the induced natural transformation
$$(\ul\sP^{\on{spec,loc}}_\cG)^{\int\, \on{ins.unit}}\to \ul\sP^{\on{spec,glob}}_\cG$$
is the isomorphism of \cite[Theorem 17.7.2]{GLC2}.

\medskip 

This formally implies that \eqref{e:Poinc Gamma 2} is also an isomorphism. 

\qed[\thmref{t:Poinc Gamma}]

\newpage

\centerline{\bf Part III: Applications to the Langlands functor}

\bigskip

In this Part we use the material developed in Part II of the paper to establish properties
of the Langlands functor
$$\BL_G::\Dmod_{\frac{1}{2}}(\Bun_G)\to \IndCoh_\Nilp(\LS_\cG)$$
vis-a-vis functors that connect $G$ to Levi quotients of its parabolic subgroups. 

\medskip

We prove the following:

\begin{itemize}

\item It intertwines the $\Eis_!$ functor on the geometric side with the $\Eis^{\on{spec}}$
functor on the spectral side (up to some twists and shifts), see \thmref{t:L and Eis},
and similarly for the enhanced versions; 

\smallskip

\item It intertwines the $\on{CT}_*$ functor on the geometric side with the $\on{CT}^{\on{spec}}$
functor on the spectral side (up to some twists and shifts), see \thmref{t:CT compat}, 
and similarly for the enhanced versions; 

\smallskip

\item The Langlands functor admits a left adjoint, see \thmref{t:left adjoint};

\medskip

\item The left adjoint of the Langlands functor identifies with its dual (up to some twists and shifts),
see \thmref{t:left adjoint as dual}.

\end{itemize}

\bigskip

The main result of this Part and of the entire paper, \thmref{t:main}, says that the Langlands functor
$\BL_G$ and its left adjoint $\BL_G^L$ define mutually inverse equivalences between the following 
full subcategories on the geometric and spectral sides, respectively:

\begin{itemize}

\item On the geometric side, this is the subcategory $$\Dmod_{\frac{1}{2}}(\Bun_G)_{\Eis}\subset \Dmod_{\frac{1}{2}}(\Bun_G),$$
generated by the essential images of the Eisenstein functors $\Eis_!$ for all \emph{proper} parabolics;

\medskip

\item On the spectral side, this is the full subcategory $$\IndCoh_\Nilp(\LS_\cG)_{\on{red}}\subset \IndCoh_\Nilp(\LS_\cG),$$
consisting of objects, set-theoretically supported on the locus $\LS^{\on{red}}_\cG\subset \LS_\cG$, consisting of
\emph{reducible} local systems.

\end{itemize}

\medskip 

We will reduce the assertion of \thmref{t:main} to \propref{p:AP 2}, which says that 
the functors
$$\on{CT}^{-,\on{spec}}_{!*} \text{ and } \on{CT}^{-,\on{spec}}_{!*}\circ \BL_G\circ \BL_G^L,$$
are canonically isomorphic, where $\on{CT}^{-,\on{spec}}_{!*}$ is the ``compactified" 
spectral constant term functor from \secref{sss:comp spec CT}.

\bigskip

\section{Langlands functor and Eisenstein series} \label{s:L and Eis}

In this section we establish the compatibility of the Langlands functor with the functors 
$\Eis_!$ and $\Eis^{\on{spec}}$, as well as their enhancements.

\ssec{The Langlands functor: recollections}

In this subsection we recall some material from \cite[Sect. 18]{GLC2}.

\sssec{}

The Langlands functor
$$\BL_G:\Dmod_{\frac{1}{2}}(\Bun_G)\to \IndCoh_\Nilp(\LS_\cG)$$
was introduced in \cite[Sect. 1.6.6]{GLC1}. 

\sssec{Whittaker compatibility}

Combining \cite[Corollary 1.6.5]{GLC1} with \cite[Proposition 18.1.6]{GLC2}
(and the fully faithfulness of the functor $\Gamma_\cG^{\on{spec}}$), 
we obtain that 
the functor $\BL_G$ is characterized uniquely by the following two properties:

\begin{itemize}

\item It makes the diagram
\begin{equation} \label{e:Whit compat}
\CD
\Whit^!(G)_\Ran @>{\on{CS}_G}>>  \Rep(\cG)_\Ran \\
@A{\on{coeff}_G[2\delta_{N_{\rho(\omega_X)}}]}AA @AA{\Gamma^{\on{spec},\IndCoh}_\cG}A  \\
\Dmod_{\frac{1}{2}}(\Bun_G) @>{\BL_G}>> \IndCoh_\Nilp(\LS_\cG)
\endCD
\end{equation} 
commute;

\medskip

\item It sends compact objects in $\Dmod_{\frac{1}{2}}(\Bun_G)$ to $\IndCoh_\Nilp(\LS_\cG)^{>-\infty}$.

\end{itemize}

\sssec{Compatibility with the Hecke action }

According to \cite[Corollary 18.3.10]{GLC2}, the functor $\BL_G$ is compatible with the spherical Hecke action in the sense
that the functor
$$\Dmod_{\frac{1}{2}}(\Bun_G)\otimes \Dmod(\Ran) \overset{\BL_G\otimes \on{Id}}\longrightarrow 
\IndCoh_\Nilp(\LS_\cG)\otimes \Dmod(\Ran)$$
is compatible with the actions of
\begin{equation} \label{e:Sat}
\Sph_{G,\Ran} \overset{\Sat_G}\simeq \Sph^{\on{spec}}_{\cG,\Ran}
\end{equation} 
on each side, respectively.

\sssec{Compatibility with Kac-Moody localization}

According to \cite[Theorem 18.5.2]{GLC2}, we have the following commutative diagram
\begin{equation} \label{e:Loc compat}
\CD
\Dmod_{\frac{1}{2}}(\Bun_G) @>{\BL_G}>> \IndCoh_\Nilp(\LS_\cG) \\
@A{\Loc_G\otimes \fl[-d]}AA @AA{\on{Poinc}^{\on{spec}}_{G,*}}A \\
\KL(G)_{\crit,\Ran} @>{\FLE_{G,\crit}}>> \IndCoh^*(\Op^\mf_\cG)_\Ran,
\endCD
\end{equation}
where
$$\fl:=\fl^{\otimes \frac{1}{2}}_{G,N_{\rho(\omega_X)}}\otimes \fl^{\otimes -1}_{N_{\rho(\omega_X)}} \text{ and }
d=\delta_{N_{\rho(\omega_X)}}$$
for the lines $\fl^{\otimes \frac{1}{2}}_{G,N_{\rho(\omega_X)}}$ and $\fl_{N_{\rho(\omega_X)}}$ from \secref{sss:Whit Loc comp}. 

\sssec{}

Let us rewrite \eqref{e:Loc compat} equivalently as 
\begin{equation} \label{e:Loc compat Ran}
\CD
\Dmod_{\frac{1}{2}}(\Bun_G)\otimes \Dmod(\Ran) @>{\BL_G\otimes \on{Id}}>> \IndCoh_\Nilp(\LS_\cG)\otimes \Dmod(\Ran) \\
@A{\Loc_{G,\Ran}\otimes \fl[-d]}AA @AA{\on{Poinc}^{\on{spec}}_{G,*,\Ran}}A \\
\KL(G)_{\crit,\Ran} @>{\FLE_{G,\crit}}>> \IndCoh^*(\Op^\mf_\cG)_\Ran,
\endCD
\end{equation}
and note that left column is a functor between $\Sph_{G,\Ran}$-module categories, and the right column is a functor
between $\Sph^{\on{spec}}_{\cG,\Ran}$-module categories, while the horizontal arrows intertwine these actions via
\eqref{e:Sat}.

\medskip

We claim now that the statement of \cite[Theorem 18.5.2]{GLC2} can be complemented as follows:

\begin{thm} \label{t:Langlands and Loc Sat}
The datum of commutativity in \eqref{e:Loc compat Ran} is compatible with the Hecke actions via \eqref{e:Sat}.
\end{thm} 

The proof will be given in \secref{ss:Langlands and Loc Sat}. 

\sssec{The twisted version}

Let $\CP_{Z^0_\cG}$ be a $Z^0_\cG$-torsor on $X$. Recall the twisted FLE:
$$\on{FLE}_{G,\crit+\on{dlog}(\CP_{Z^0_\cG})}:\KL(G)_{\crit+\on{dlog}(\CP_{Z^0_\cG})}\to 
\IndCoh^*(\Op_{\cG,\CP_{Z^0_\cG}}^{\on{mon-free}}),$$
see \secref{e:twisted FLE}.

\medskip

Using \eqref{e:twisted FLE compat} (instead of the untwisted version, i.e., \cite[Theorem 6.4.5]{GLC2}), one proves the
following variant of 
\cite[Theorem 18.5.2]{GLC2} :

\begin{thm} \label{t:Langlands and Loc Sat tw}
We have a commutative diagram
$$
\CD
\Dmod_{\frac{1}{2}}(\Bun_G)@>{\BL_G}>> \IndCoh_\Nilp(\LS_\cG) \\
@A{\Loc_{G,\CP_{Z^0_\cG}}\otimes \fl[-d]}AA 
@AA{\on{Poinc}^{\on{spec}}_{G,*,\CP_{Z^0_\cG}}}A \\
\KL(G)_{\crit+\on{dlog}(\CP_{Z^0_\cG}),\Ran} @>{\FLE_{G,\crit+\on{dlog}(\CP_{Z^0_\cG})}}>> \IndCoh^*(\Op^\mf_{\cG,\CP_{Z^0_\cG}})_\Ran.
\endCD
$$
\end{thm}

\sssec{}

We will apply \thmref{t:Langlands and Loc Sat tw} when the group in question is the Levi subgroup $M$ of the original $G$
and $\CP_{Z^0_\cG}=\rhoch_P(\omega_X)$. In this case, the commutative diagram in \thmref{t:Langlands and Loc Sat tw} reads as
\begin{equation}  \label{e:Langlands and Loc Levi}
\hskip1cm
\CD
\Dmod_{\frac{1}{2}}(\Bun_M) @>{\BL_M}>> \IndCoh_\Nilp(\LS_\cM) \\
@A{\Loc_{M,\rhoch_P}\otimes \fl[-d]}AA
@AA{\on{Poinc}^{\on{spec}}_{G,*,\rhoch_P}}A \\
\KL(M)_{\crit_M+\rhoch_P,\Ran} @>{\FLE_{M,\crit+\rhoch_P}}>> \IndCoh^*(\Op^\mf_{\cM,\rhoch_P})_\Ran,
\endCD
\end{equation} 
where
$$\fl:=\fl^{\otimes \frac{1}{2}}_{M,N(M)_{\rho_M(\omega_X)}}\otimes \fl^{\otimes -1}_{N(M)_{\rho_M(\omega_X)}}$$
and 
$$d:=\delta_{N(M)_{\rho_M(\omega_X)}}.$$

\ssec{Compatibility of the Langlands functor with Eisenstein series} 

In this subsection we establish the compatibility of the Langlands functor with (usual) Eisenstein series.

\sssec{}

We are going to prove:

\begin{thm} \label{t:L and Eis}
The following diagram commutes
\begin{equation} \label{e:L and Eis}
\CD
\Dmod_{\frac{1}{2}}(\Bun_G) @>{\BL_G}>> \IndCoh_\Nilp(\LS_\cG) \\
@A{\Eis^-_{!,\rho_P(\omega_X)}[\delta_{(N^-_P)_{\rho_P(\omega_X)}}]}AA @AA{\Eis^{-,\on{spec}}}A \\
\Dmod_{\frac{1}{2}}(\Bun_M) @>{\BL_M}>> \IndCoh_\Nilp(\LS_\cM).
\endCD
\end{equation}
\end{thm} 

The rest of this subsection is devoted to the proof of this theorem.

\sssec{}

We claim that in order to prove the commutativity of \eqref{e:L and Eis}, it sufficient to prove the
commutativity of
\begin{equation} \label{e:L and Eis coarse}
\CD
\Dmod_{\frac{1}{2}}(\Bun_G) @>{\Psi_{\LS_\cG}\circ \BL_G}>> \QCoh(\LS_\cG) \\
@A{\Eis^-_{!,\rho_P(\omega_X)}[\delta_{(N^-_P)_{\rho_P(\omega_X)}}]}AA @AA{\Psi_{\LS_\cG}\circ \Eis^{-,\on{spec}}}A \\
\Dmod_{\frac{1}{2}}(\Bun_M) @>{\BL_M}>> \IndCoh_\Nilp(\LS_\cM).
\endCD
\end{equation}

Indeed, it is sufficient to show that \eqref{e:L and Eis} commutes, when evaluated on compact objects
of $\Dmod_{\frac{1}{2}}(\Bun_M)$. 

\medskip

We claim that both circuits in \eqref{e:L and Eis} send $\Dmod_{\frac{1}{2}}(\Bun_M)^c$ to
$\IndCoh_\Nilp(\LS_\cG)^{>-\infty}$. This would prove the required assertion, since the functor $\Psi_{\LS_\cG}$
is fully faithful when restricted to $\IndCoh_\Nilp(\LS_\cG)^{>-\infty}$.

\medskip

For the clockwise circuit, the claim follows from the fact that the functor $\Eis^-_{!,\rho_P(\omega_X)}$
preserves compactness (it admits a right adjoint, namely, $\on{CT}^-_{*,\rho_P(\omega_X)}$) and the fact that 
$\BL_G$ sends compact objects to $\IndCoh_\Nilp(\LS_\cG)^{>-\infty}$.

\medskip

For the counterclockwise circuit, this follows from the fact that the functor $\BL_M$ sends compact objects to 
$\IndCoh_\Nilp(\LS_\cM)^{>-\infty}$ and the fact that the functor $\Eis^{-,\on{spec}}$ is left t-exact
(up to a cohomological shift). 

\sssec{}

Next, since the functor 
$$\Gamma^{\on{spec}}_\cG:\QCoh(\LS_\cG)\to \Rep(\cG)_\Ran$$
is fully faithful, in order to prove the commutatiion of \eqref{e:L and Eis coarse}, it suffices to prove 
the commutativity of 
$$
\CD
\Dmod_{\frac{1}{2}}(\Bun_G) @>{\Gamma^{\on{spec}}_\cG\circ \Psi_{\LS_\cG}\circ \BL_G}>>  \Rep(\cG)_\Ran \\
@A{\Eis^-_{!,\rho_P(\omega_X)}[\delta_{(N^-_P)_{\rho_P(\omega_X)}}]}AA @AA{\Gamma^{\on{spec}}_\cG\circ \Psi_{\LS_\cG}\circ \Eis^{-,\on{spec}}}A \\
\Dmod_{\frac{1}{2}}(\Bun_M) @>{\BL_M}>> \IndCoh_\Nilp(\LS_\cM),
\endCD
$$
which is the same as 
\begin{equation} \label{e:L and Eis Whit}
\CD
\Dmod_{\frac{1}{2}}(\Bun_G) @>{\Gamma^{\on{spec},\IndCoh}_\cG\circ \BL_G}>> \Rep(\cG)_\Ran \\
@A{\Eis^-_{!,\rho_P(\omega_X)}[\delta_{(N^-_P)_{\rho_P(\omega_X)}}]}AA @AA{\Gamma^{\on{spec},\IndCoh}_\cG \circ \Eis^{-,\on{spec}}}A \\
\Dmod_{\frac{1}{2}}(\Bun_M) @>{\BL_M}>> \IndCoh_\Nilp(\LS_\cM),
\endCD
\end{equation}

\medskip

Hence, given the commutation of \eqref{e:Whit compat}, we obtain that it suffices to prove the commutation of
\begin{equation} \label{e:L and Eis 1}
\CD
\Whit^!(G)_\Ran @>{\on{CS}_G}>>  \Rep(\cG)_\Ran \\
@A{\on{coeff}_G[2\delta_{N_{\rho(\omega_X)}}]}AA @AA{\Gamma^{\on{spec},\IndCoh}_\cG}A  \\
\Dmod_{\frac{1}{2}}(\Bun_G) & & \IndCoh_\Nilp(\LS_\cG) \\
@A{\Eis^-_{!,\rho_P(\omega_X)}[\delta_{(N^-_P)_{\rho_P(\omega_X)}}]}AA @AA{\Eis^{-,\on{spec}}}A \\
\Dmod_{\frac{1}{2}}(\Bun_M) @>{\BL_M}>> \IndCoh_\Nilp(\LS_\cM).
\endCD
\end{equation} 

\sssec{} \label{sss:use Part II}

We now use the material from Part II of the paper. Namely, by \corref{c:coeff of Eis !}, 
the left vertical composition in \eqref{e:L and Eis 1} identifies with
\begin{multline} \label{e:L and Eis geom}
\Dmod_{\frac{1}{2}}(\Bun_M) \overset{\on{coeff}_{M,\on{untl}}}\longrightarrow  \Whit^!(M)_{\Ran^{\on{untl}}} 
\overset{(J^{-,!}_{\Whit,\Theta})^\vee}\longrightarrow \\
\to \Whit^!(G)_{\Ran^{\on{untl}}} \overset{\on{emb.indep}_{\Whit^!(G)}^L}\longrightarrow \Whit^!(G)_{\Ran^{\on{untl}},\on{indep}}\to \Whit^!(G)_\Ran.
\end{multline}

Similarly, by \corref{c:Eis spec Gamma}, the right vertical composition in \eqref{e:L and Eis 1} identifies with
\begin{multline} \label{e:L and Eis spec}
\IndCoh_\Nilp(\LS_\cM)\overset{\Gamma^{\on{spec},\IndCoh}_{\cM,\on{untl}}}\longrightarrow \Rep(\cM)_{\Ran^{\on{untl}}}
\overset{\on{co}\!J^{-,\on{spec},!}}\longrightarrow \\
\to \Rep(\cG)_{\Ran^{\on{untl}}} \overset{\on{emb.indep}_{\Rep(\cG)}^L}\longrightarrow \Rep(\cG)_{\Ran^{\on{untl}},\on{indep}}\to \Rep(\cG)_\Ran.
\end{multline}

Thus, it suffices to establish the commutation of each of the three squares in
\begin{equation} \label{e:L and Eis 2}
\CD
\Whit^!(G)_{\Ran^{\on{untl}},\on{indep}} @>{\on{CS}_G}>> \Rep(\cG)_{\Ran^{\on{untl}},\on{indep}} \\
@A{\on{emb.indep}_{\Whit^!(G)}^L}AA @AA{\on{emb.indep}_{\Rep(\cG)}^L}A \\
\Whit^!(G)_{\Ran^{\on{untl}}}  @>{\on{CS}_G}>> \Rep(\cG)_{\Ran^{\on{untl}}} \\
@A{(J^{-,!}_{\Whit,\Theta})^\vee}AA  @AA{\on{co}\!J^{-,\on{spec},!}}A \\
\Whit^!(M)_{\Ran^{\on{untl}}}  @>{\on{CS}_M}>>  \Rep(\cM)_{\Ran^{\on{untl}}} \\
@A{\on{coeff}_{M,\on{untl}}}AA @AA{\Gamma^{\on{spec},\IndCoh}_{\cM,\on{untl}}}A \\
\Dmod_{\frac{1}{2}}(\Bun_M) @>{\BL_G}>> \IndCoh_\Nilp(\LS_\cM). 
\endCD
\end{equation} 

\sssec{}

The upper square in \eqref{e:L and Eis 2} commutes tautologically (since $\on{CS}_G$ is an equivalence of
unital factorization categories). 

\medskip

The lower square is \eqref{e:Whit compat} for $M$.

\sssec{}

Finally, we use the input from Part I of the paper. Namely, the commutative of the middle square in \eqref{e:L and Eis 2}
is obtained by duality from \eqref{e:!-Jacquet on Whit Theta}.

\qed[\thmref{t:L and Eis}]

\begin{rem}

Note that one can depict the proof of \thmref{t:L and Eis} given above as a cube: 

\begin{equation} \label{e:Eis cube}
\xy
(0,0)*+{\Whit^!(M)_{\Ran^{\on{untl}}}}="X";
(40,20)*+{\Whit^!(G)_{\Ran^{\on{untl}},\on{indep}}}="Y";
(70,0)*+{\Rep(\cM)_{\Ran^{\on{untl}}}}="X'";
(110,20)*+{\Rep(\cG)_{\Ran^{\on{untl}},\on{indep}}}="Y'";
{\ar@{->} "X";"Y"};
{\ar@{->} "X'";"Y'"};
{\ar@{->}^{\on{CS}_M} "X";"X'"};
{\ar@{->}^{\on{CS}_G} "Y";"Y'"};
(0,-50)*+{\Dmod_{\frac{1}{2}}(\Bun_M)}="Z";
(40,-30)*+{\Dmod_{\frac{1}{2}}(\Bun_G)}="W";
(70,-50)*+{\IndCoh_\Nilp(\LS_\cM(X)),}="Z'";
(110,-30)*+{\IndCoh_\Nilp(\LS_\cG(X))}="W'";
{\ar@{->}_{\Eis^-{!,\rho_P(\omega_X)}} "Z";"W"};
{\ar@{->}_{\Eis^{-,\on{spec}}} "Z'";"W'"};
{\ar@{->}_{\BL_M} "Z";"Z'"};
{\ar@{->}^{\BL_G} "W";"W'"};
{\ar@{->}^{\on{coeff}_{M,\on{untl}}} "Z";"X"};
{\ar@{->}_{\on{coeff}_{G,\on{untl}}} "W";"Y"};
{\ar@{->}^{\Gamma^{\on{spec},\IndCoh}_{\cM,\on{untl}}} "Z'";"X'"};
{\ar@{->}^{\Gamma^{\on{spec},\IndCoh}_{\cG,\on{untl}}} "W'";"Y'"};
\endxy
\end{equation}
where the upper slanted arrows are
$$\Whit^!(M)_{\Ran^{\on{untl}}} \overset{(J^{-,!}_{\Whit,\Theta})^\vee}\longrightarrow 
\Whit^!(G)_{\Ran^{\on{untl}}} \overset{\on{emb.indep}_{\Whit^!(G)}^L}\longrightarrow \Whit^!(G)_{\Ran^{\on{untl}},\on{indep}}$$
and 
$$\Rep(\cM)_{\Ran^{\on{untl}}}\overset{\on{co}\!J^{-,\on{spec},!}}\longrightarrow \Rep(\cG)_{\Ran^{\on{untl}}}
\overset{\on{emb.indep}_{\Rep(\cG)}^L}\longrightarrow \Rep(\cG)_{\Ran^{\on{untl}},\on{indep}},$$
respectively. 

\medskip

\thmref{t:L and Eis} asserts that the bottom lid of \eqref{e:Eis cube} commutes. The proof consists of showing that all the other
faces commute: 

\begin{itemize}

\item The front and back face encode the compatibility of the Langlands functor with Whittaker coefficients;

\item The side faces encode the compatibility of the Eisenstein functor with the corresponding local-to-global
functor on each side of Langlands correspondence;

\item The top lid encodes the compatibility of the Casselman-Shalika equivalence with Jacquet functors. 

\end{itemize}

\end{rem} 

\ssec{Compatibility of the Langlands functor with Eisenstein series--the enhanced version} 

In this subsection we will prove a generalization of \thmref{t:L and Eis}, which establishes the compatibility 
of the Langlands functor with \emph{enhanced} Eisenstein series.

\sssec{} \label{sss:L M enh}

Fix $\CZ\to \Ran$. Since the Langlands functor $\BL_M$ is compatible with the Hecke actions, we can tensor it with $\Sat^\semiinf$ and thus obtain 
a functor $$\BL_{M,\CZ}^{-,\on{enh}}:\Dmod_{\frac{1}{2}}(\Bun_M)_\CZ^{-,\on{enh}}\to \IndCoh_\Nilp(\LS_\cM)^{-,\on{enh}}_\CZ.$$

\medskip

The fact that the functor $\Eis^{-,\on{spec,enh}}$ sends $\IndCoh_\Nilp(\LS_\cM)$ to $\IndCoh_\Nilp(\LS_\cG)$
implies that the functor $\Eis^{-,\on{spec,enh}}_\CZ$ sends $\IndCoh_\Nilp(\LS_\cM)^{-,\on{enh}}_\CZ$ to
$\IndCoh_\Nilp(\LS_\cG) \otimes \Dmod(\CZ)$ (indeed, the essential image if $\ind^{\on{spec}}_{\on{enh}}$ generates
the target). 

\sssec{}

We are going to prove the following:

\begin{thm} \label{t:L and Eis enh}
The following diagram commutes
\begin{equation} \label{e:L and Eis enh}
\hskip1cm
\CD
\Dmod_{\frac{1}{2}}(\Bun_G)\otimes \Dmod(\CZ) @>{\BL_G\otimes \on{Id}_\CZ}>> \IndCoh_\Nilp(\LS_\cG) \otimes \Dmod(\CZ) \\
@A{\Eis^{-,\on{enh}}_{!,\rho_P(\omega_X),\CZ}[\delta_{(N^-_P)_{\rho_P(\omega_X)}}]}AA @AA{\Eis^{-,\on{spec,enh}}_\CZ}A \\
\Dmod_{\frac{1}{2}}(\Bun_M)^{-,\on{enh}}_\CZ @>{\BL^{-,\on{enh}}_{M,\CZ}}>> \IndCoh_\Nilp(\LS_\cM)^{-,\on{enh}}_\CZ,
\endCD
\end{equation}
in a way compatible with the Hecke actions via \eqref{e:Sat}. 
\end{thm} 

\begin{proof}

The proof follows closely that of \thmref{t:L and Eis}. The fact that both circuits in \eqref{e:L and Eis}
send compact objects in $\Dmod_{\frac{1}{2}}(\Bun_M)$ to eventually coconnective objects in $\IndCoh_\Nilp(\LS_\cG)$
formally implies that the same is true for \eqref{e:L and Eis enh},

\medskip

Hence, we obtain that it is sufficient to show that the diagram
\begin{equation} \label{e:L and Eis enh 1}
\CD
\Whit^!(G)_{\CZ^{\subseteq,\on{untl}},\on{indep}} @>{\on{CS}_G}>>  \Rep(\cG)_{\CZ^{\subseteq,\on{untl}},\on{indep}} \\
@A{\on{coeff}_{G,\CZ,\on{untl}}[2\delta_{N_{\rho(\omega_X)}}]}AA @AA{\Gamma^{\on{spec},\IndCoh}_{\cG,\on{untl}}}A  \\
\Dmod_{\frac{1}{2}}(\Bun_G)\otimes \Dmod(\CZ) & & \IndCoh_\Nilp(\LS_\cG) \otimes \Dmod(\CZ)\\
@A{\Eis^{-,\on{enh}}_{!,\rho_P(\omega_X),\CZ}[\delta_{(N^-_P)_{\rho_P(\omega_X)}}]}AA @AA{\Eis^{-,\on{spec,enh}}_\CZ}A \\
\Dmod_{\frac{1}{2}}(\Bun_M)^{-,\on{enh}}_\CZ @>{\BL^{-,\on{enh}}_{M,\CZ}}>> \IndCoh_\Nilp(\LS_\cM)^{-,\on{enh}}_\CZ
\endCD
\end{equation} 
commutes, in a way compatible with the Hecke actions via \eqref{e:Sat}. 

\medskip

By \corref{c:coeff of Eis ! enh} and \thmref{t:Eis spec enh Gamma}, we replace \eqref{e:L and Eis enh 1} by

\medskip

\begin{equation} \label{e:L and Eis enh 2}
\CD
\Whit^!(G)_{\CZ^{\subseteq,\on{untl}},\on{indep}} @>{\on{CS}_G}>>  \Rep(\cG)_{\CZ^{\subseteq,\on{untl}},\on{indep}} \\
@A{\on{emb.indep}_{\Whit^!(G)}^L}AA  @AA{\on{emb.indep}_{\Rep(\cG)}^L} A \\
\Whit^!(G)_{\CZ^{\subseteq,\on{untl}}} @>{\on{CS}_G}>> \Rep(\cG)_{\CZ^{\subseteq,\on{untl}}} \\
@A{\on{co}\!J^{-,\on{enh}}_{\Whit,\Theta}}AAA @AA{\on{co}\!J^{-,\on{spec,enh}}}A \\
(\Whit^!(M)\underset{\Sph_M}\otimes \on{I}(G,P^-)^{\on{loc}}_{\rho_P(\omega_X)})_{\CZ^{\subseteq,\on{untl}}} @>{\on{CS}_M\otimes \Sat^\semiinf}>> 
(\Rep(\cM)\underset{\Sph_\cM^{\on{spec}}}\otimes \on{I}(\cG,\cP^-)^{\on{spec,loc}})_{\CZ^{\subseteq,\on{untl}}} \\
@AAA @AAA \\
\Whit^!(M)_{\CZ^{\subseteq,\on{untl}}}\underset{\Sph_{M,\CZ}}\otimes  \on{I}(G,P^-)^{\on{loc}}_{\rho_P(\omega_X),\CZ^{\subseteq,\on{untl}}} 
@>{\on{CS}_M\otimes \Sat^\semiinf}>> 
\Rep(\cM)_{\CZ^{\subseteq,\on{untl}}}\underset{\Sph_{\cM,\CZ}^{\on{spec}}}\otimes \on{I}(\cG,\cP^-)^{\on{spec,loc}}_{\CZ^{\subseteq,\on{untl}}} \\
@A{\on{coeff}_{M,\CZ,\on{untl}}\otimes \on{ins.unit}_\CZ}AA @AA{\Gamma^{\on{spec},\IndCoh}_{\cM,\CZ,\on{untl}}\otimes \on{ins.unit}_\CZ}A \\
(\Dmod_{\frac{1}{2}}(\Bun_M)\otimes \Dmod(\CZ))\underset{\Sph_{M,\CZ}}\otimes \on{I}(G,P^-)^{\on{loc}}_{\rho_P(\omega_X),\CZ} 
@>{\BL_M\otimes \Sat^\semiinf}>>
(\IndCoh(\LS_\cM)\otimes \Dmod(\CZ))\underset{\Sph_{\cM,\CZ}^{\on{spec}}}\otimes \on{I}(\cG,\cP^-)^{\on{spec,loc}}_\CZ \\
@A{=}AA @AA{=}A \\
\Dmod_{\frac{1}{2}}(\Bun_M)^{-,\on{enh}}_{\rho_P(\omega_X),\CZ}  @>{\BL^{-,\on{enh}}_{M,\CZ}}>> \IndCoh_\Nilp(\LS_\cM)^{-,\on{enh}}_\CZ. 
\endCD
\end{equation}

We claim that all the inner squares in \eqref{e:L and Eis enh 2} commute. This is tautological for the top, middle and bottom squares.
For the second square from the bottom, the commutation follows from the fact that diagram \eqref{e:Whit compat} for $M$ commutes,
compatibly with Hecke actions. 

\medskip

For the second square from the top, the commutation follows from the commutation of 
$$
\CD 
\Whit^!(G) @>{\on{CS}_G}>> \Rep(\cG) \\
@AA{\on{co}\!J^{-,\on{enh}}_{\Whit,\Theta}}A  @AA{\on{co}\!J^{-,\on{spec,enh}}}A \\
\Whit^!(M)\underset{\Sph_M}\otimes \on{I}(G,P^-)^{\on{loc}}_{\rho_P(\omega_X)} @>{\on{CS}_M\otimes \Sat^\semiinf}>> 
\Rep(\cM)\underset{\Sph_\cM^{\on{spec}}}\otimes \on{I}(\cG,\cP^-)^{\on{spec,loc}},
\endCD
$$
which is obtained by duality from \eqref{e:semiinf CS diag}. 

\end{proof}

\sssec{}

Combining \thmref{t:L and Eis enh} with Propositions \ref{p:semiinf Sat IC}, \ref{p:compact Eis} and \ref{p:Eis !* spec},
we obtain: 

\begin{cor} \label{c:L and Eis !*}
The following diagram commutes
\begin{equation} \label{e:L and Eis !*}
\CD
\Dmod_{\frac{1}{2}}(\Bun_G) @>{\BL_G}>> \IndCoh_\Nilp(\LS_\cG) \\
@A{\Eis^-_{!*,\rho_P(\omega_X)}[\delta_{(N^-_P)_{\rho_P(\omega_X)}}]}AA @AA{\Eis^{-,\on{spec}}_{!*}}A \\
\Dmod_{\frac{1}{2}}(\Bun_M) @>{\BL_M}>> \IndCoh_\Nilp(\LS_\cM).
\endCD
\end{equation}
\end{cor} 

\ssec{Proof of \thmref{t:Langlands and Loc Sat}} \label{ss:Langlands and Loc Sat}

\sssec{}

We need to show that for $\CZ\to \Ran$, the diagram
$$ 
\CD
\Dmod_{\frac{1}{2}}(\Bun_G)\otimes \Dmod(\CZ) @>{\BL_G\otimes \on{Id}}>> \IndCoh_\Nilp(\LS_\cG)\otimes \Dmod(\CZ) \\
@A{\Loc_{G,\CZ}\otimes \fl[-d]}AA @AA{\on{Poinc}^{\on{spec}}_{G,*,\CZ}}A \\
\KL(G)_{\crit,\CZ} @>{\FLE_{G,\crit}}>> \IndCoh^*(\Op^\mf_\cG)_\CZ,
\endCD
$$
commutes, in a way compatible with actions of
\begin{equation} \label{e:Sat Z}
\Sph_{G,\CZ} \overset{\Sat}\simeq \Sph^{\on{spec}}_{\cG,\CZ}.
\end{equation} 

\medskip

By the same logic as in the proof of \cite[Theorem 18.5.2]{GLC2}, it is enough to show that the diagram
\begin{equation} \label{e:Langlands and Loc Sat 1}
\CD 
\Whit^!(G)_{\CZ^\subseteq} @>>> \Rep(\cG)_{\CZ^\subseteq} \\
@A{\on{coeff}_{G,\CZ,\on{untl}}}AA @AA{\Gamma^{\on{spec},\IndCoh}_{\cG,\CZ,\on{untl}}}A \\
\Dmod_{\frac{1}{2}}(\Bun_G)\otimes \Dmod(\CZ) & & \IndCoh_\Nilp(\LS_\cG)\otimes \Dmod(\CZ) \\
@A{\Loc_{G,\CZ}\otimes \fl[d]}AA @AA{\on{Poinc}^{\on{spec}}_{G,*,\CZ}}A \\
\KL(G)_{\crit,\CZ} @>{\FLE_{G,\crit}}>> \IndCoh^*(\Op^\mf_\cG)_\CZ
\endCD
\end{equation}
commutes in a way compatible with the actions of \eqref{e:Sat Z}, where $\on{coeff}_{G,\CZ,\on{untl}}$ is as in 
\secref{sss:coeff Z untl}, and similarly for $\Gamma^{\on{spec},\IndCoh}_{\cG,\CZ,\on{untl}}$. 

\sssec{}

By duality, the datum of the left column in \eqref{e:Langlands and Loc Sat 1} as a functor equipped with an action of $\Sph_{G,\CZ}$ is
equivalent to the datum of the $\Dmod(\CZ)$-linear functor
\begin{multline}  \label{e:Langlands and Loc Sat geom}
\KL(G)_{\crit,\CZ} \underset{\Sph_{G,\CZ}}\otimes \Whit_*(G)_{\CZ^\subseteq} \overset{(\Loc_{G,\CZ}\otimes \fl[d])\otimes \on{Id}}\longrightarrow \\
\to (\Dmod_{\frac{1}{2}}(\Bun_G)\otimes \Dmod(\CZ)) \underset{\Sph_{G,\CZ}}\otimes \Whit_*(G)_{\CZ^\subseteq} \to \\
\overset{\on{coeff}_{G,\CZ,\on{untl}}\otimes \on{Id}}\longrightarrow 
\Whit^!(G)_{\CZ^\subseteq}\underset{\Sph_{G,\CZ}}\otimes \Whit_*(G)_{\CZ^\subseteq} 
\to \Whit^!(G)_{\CZ^\subseteq}\underset{\Sph_{G,\CZ^\subseteq}}\otimes \Whit_*(G)_{\CZ^\subseteq} \to \\
\to \Dmod(\CZ^\subseteq)
\overset{(\on{pr}_{\on{small},\CZ})_!}\longrightarrow \Dmod(\CZ),
\end{multline} 
where the second-to-last functor is the canonical pairing.

\medskip

Similarly, the datum of the right column in  \eqref{e:Langlands and Loc Sat 1} as a functor equipped with an action of $\Sph^{\on{spec}}_{\cG,\CZ}$ is
equivalent to the datum of the $\Dmod(\CZ)$-linear functor
\begin{multline}  \label{e:Langlands and Loc Sat spec}
\IndCoh^*(\Op^\mf_\cG)_\CZ\underset{\Sph^{\on{spec}}_{\cG,\CZ}}\otimes \Rep(\cG)_{\CZ^\subseteq}
\overset{\on{Poinc}^{\on{spec}}_{G,*,\CZ}\otimes \on{Id}}\longrightarrow \\
\to (\IndCoh_\Nilp(\LS_\cG)\otimes \Dmod(\CZ)) \underset{\Sph^{\on{spec}}_{\cG,\CZ}}\otimes \Rep(\cG)_{\CZ^\subseteq}\to \\
\overset{\Gamma^{\on{spec},\IndCoh}_{\cG,\CZ,\on{untl}}\otimes \on{Id}}\longrightarrow 
\Rep(\cG)_{\CZ^\subseteq}\underset{\Sph^{\on{spec}}_{\cG,\CZ}}\otimes \Rep(\cG)_{\CZ^\subseteq}\to 
\Rep(\cG)_{\CZ^\subseteq}\underset{\Sph^{\on{spec}}_{\cG,\CZ^\subseteq}}\otimes \Rep(\cG)_{\CZ^\subseteq}\to \\
\to \Dmod(\CZ^\subseteq)
\overset{(\on{pr}_{\on{small},\CZ})_!}\longrightarrow \Dmod(\CZ),
\end{multline}
where the second-to-last functor is the canonical pairing, corresponding to the self-duality of $\Rep(\cG)$. 

\medskip

Hence, in order to establish the commutativity of \eqref{e:Langlands and Loc Sat 1}, it suffices to show that the functors 
\eqref{e:Langlands and Loc Sat geom} and \eqref{e:Langlands and Loc Sat spec} match up under
\begin{equation} \label{e:Langlands and Loc Sat 2}
\KL(G)_{\crit,\CZ} \underset{\Sph_{G,\CZ}}\otimes \Whit_*(G)_{\CZ^\subseteq} 
\overset{\FLE_{G,\crit}\otimes \FLE^{-1}_{\cG,\infty}}\simeq 
\IndCoh^*(\Op^\mf_\cG)_\CZ\underset{\Sph^{\on{spec}}_{\cG,\CZ}}\otimes \Rep(\cG)_{\CZ^\subseteq}.
\end{equation} 

\sssec{} 

We can rewrite \eqref{e:Langlands and Loc Sat geom} as 
\begin{multline}  \label{e:Langlands and Loc Sat geom 1}
\KL(G)_{\crit,\CZ} \underset{\Sph_{G,\CZ}}\otimes \Whit_*(G)_{\CZ^\subseteq} \to 
\KL(G)_{\crit,\CZ^\subseteq} \underset{\Sph_{G,\CZ^\subseteq}}\otimes \Whit_*(G)_{\CZ^\subseteq} \to \\
\overset{\sP^{\on{glob,true}}_{G,\CZ^\subseteq}\otimes \fl[d]}\longrightarrow 
\Dmod(\CZ^\subseteq)\overset{(\on{pr}_{\on{small},\CZ})_!}\longrightarrow \Dmod(\CZ),
\end{multline} 
where for $\CZ'\to \Ran$ we denote by $\sP^{\on{glob,true}}_{G,\CZ'}$ the functor
\begin{multline*}
\KL(G)_{\crit,\CZ'} \underset{\Sph_{G,\CZ'}}\otimes \Whit_*(G)_{\CZ'}  
\overset{\Loc_{G,\CZ'}\otimes \on{Id}}\longrightarrow 
(\Dmod_{\frac{1}{2}}(\Bun_G)\otimes \Dmod(\CZ')) \underset{\Sph_{G,\CZ'}}\otimes \Whit_*(G)_{\CZ'}\to \\
\overset{\on{coeff}_{G,\CZ'}\otimes \on{Id}}\longrightarrow 
\Rep(\cG)_{\CZ'}\underset{\Sph^{\on{spec}}_{\cG,\CZ'}}\otimes \Rep(\cG)_{\CZ'}
\to \Dmod(\CZ'),
\end{multline*} 
see \secref{sss:glob true}.

\medskip

Similarly, we can rewrite 
\begin{multline}  \label{e:Langlands and Loc Sat spec 1}
\IndCoh^*(\Op^\mf_\cG)_\CZ\underset{\Sph^{\on{spec}}_{\cG,\CZ}}\otimes \Rep(\cG)_{\CZ^\subseteq}\to
\IndCoh^*(\Op^\mf_\cG)_{\CZ^\subseteq}\underset{\Sph^{\on{spec}}_{\cG,\CZ^\subseteq}}\otimes \Rep(\cG)_{\CZ^\subseteq} \to \\
\overset{\sP^{\on{spec,glob,true}}_{\cG,\CZ^\subseteq}}\longrightarrow \Dmod(\CZ^\subseteq)\overset{(\on{pr}_{\on{small},\CZ})_!}\longrightarrow \Dmod(\CZ),
\end{multline} 
where for $\CZ'\to \Ran$ we denote by  $\sP^{\on{spec,glob}}_{\cG,\CZ'}$ is the functor
\begin{multline*}
\IndCoh^*(\Op^\mf_\cG)_{\CZ'}\underset{\Sph^{\on{spec}}_{\cG,\CZ'}}\otimes \Rep(\cG)_{\CZ'} \to \\
\overset{\on{Poinc}^{\on{spec}}_{G,*,\CZ'}\otimes \on{Id}}\longrightarrow 
(\IndCoh_\Nilp(\LS_\cG)\otimes \Dmod(\CZ')) \underset{\Sph^{\on{spec}}_{\cG,\CZ'}}\otimes \Rep(\cG)_{\CZ'} 
\overset{\Gamma^{\on{spec},\IndCoh}_{\cG,\CZ'}\otimes \on{Id}}\longrightarrow \\
\to \Rep(\cG)_{\CZ'}\underset{\Sph^{\on{spec}}_{\cG,\CZ'}}\otimes \Rep(\cG)_{\CZ'} \to 
\Dmod(\CZ'),
\end{multline*} 
see \secref{sss:glob true spec}.

\sssec{}

Hence, in order to establish the commutativity of \eqref{e:Langlands and Loc Sat 1}, it suffices to show that the functors 
$$\sP^{\on{glob,true}}_{G,\CZ'}\otimes \fl[d] \text{ and } \sP^{\on{spec,glob,true}}_{\cG,\CZ'}$$
match up under
\begin{equation} \label{e:Langlands and Loc Sat 3}
\KL(G)_{\crit,\CZ'} \underset{\Sph_{G,\CZ'}}\otimes \Whit_*(G)_{\CZ'} 
\overset{\FLE_{G,\crit}\otimes \FLE^{-1}_{\cG,\infty}}\simeq 
\IndCoh^*(\Op^\mf_\cG)_{\CZ'}\underset{\Sph^{\on{spec}}_{\cG,\CZ'}}\otimes \Rep(\cG)_{\CZ'}.
\end{equation} 
for any $\CZ'\to \Ran$. 

\sssec{}

We now use Theorems \ref{t:coeff Loc} and \ref{t:Poinc Gamma} to identify 
$$\sP^{\on{glob,true}}_{G,\CZ'}\otimes \fl[d]\simeq (\ul\sP^{\on{loc,true}})^{\int\,\on{ins.unit}}_{G,\CZ'}$$
and 
$$\sP^{\on{spec,glob,true}}_{\cG,\CZ'}\simeq (\ul\sP^{\on{spec,loc,true}})^{\int\,\on{ins.unit}}_{\cG,\CZ'},$$
respectively. 

\medskip

Hence, it suffices to show that under \eqref{e:Langlands and Loc Sat 3}, the functors
$$\sP^{\on{loc,true}}_\cG \text{ and } \sP^{\on{spec,loc,true}}_\cG$$
match up.

\medskip

However, this follows from \cite[Theorem 6.4.5(b)]{GLC2} by applying the functor $\Gamma^\IndCoh(\Op^\mer_\cG,-)$.

\qed[\thmref{t:Langlands and Loc Sat}]

\section{Langlands functor and constant terms} \label{s:L and CT}

Starting from this section, we will assume the validity of the geometric Langlands conjecture
(as stated in \cite[Conjecture 1.6.7]{GLC1}) for the Levi quotients of all proper parabolic subgroups of $G$. 

\medskip

We will establish a result pertaining to compatibility between the Langlands functor and
the operation of constant term. This result will be central to our proof that the Langlands
functor induces an equivalence on Eisenstein-generated categories. 

\ssec{Statement of the compatibility result}

\sssec{}

Passing to right adjoints along the vertical arrows in the commutative diagram \eqref{e:L and Eis}, we obtain a natural transformation
\begin{equation} \label{e:CT Nat trans}
\vcenter
{\xy
(0,0)*+{\Dmod_{\frac{1}{2}}(\Bun_G)}="X";
(50,0)*+{\IndCoh_\Nilp(\LS_\cG)}="Y";
(0,-30)*+{\Dmod_{\frac{1}{2}}(\Bun_M)}="Z";
(50,-30)*+{\IndCoh_\Nilp(\LS_\cM)}="W";
{\ar@{->}_{\BL_M} "Z";"W"};
{\ar@{->}^{\BL_G} "X";"Y"};
{\ar@{->}_{\on{CT}^-_{*,\rho_P(\omega_X)}[-\delta_{(N^-_P)_{\rho_P(\omega_X)}}]} "X";"Z"};
{\ar@{->}^{\on{CT}^{-,\on{spec}}} "Y";"W"};
{\ar@{=>} "Z";"Y"};
\endxy}
\end{equation}

At the end of this paper (see \secref{sss:L and CT}) we will prove:

\begin{thm} \label{t:L and CT}
The natural transformation \eqref{e:CT Nat trans} is an isomorphism.
\end{thm}

\sssec{}

However, the proof of \thmref{t:L and CT} that we will give will not be direct: we will deduce \thmref{t:L and CT} 
once we prove that the functor $\BL_G$ is an equivalence on Eisenstein-generated categories 
(\thmref{t:main}), and this will rely on the inductive assumption that GLC holds for proper Levi subgroups of $G$. 

\medskip

In this section we will prove that the two circuits in \eqref{e:CT Nat trans} are \emph{abstractly} isomorphic, see 
\thmref{t:CT compat}. Although we will give a canonical construction of the isomorphism, its relation to the natural transformation 
\eqref{e:CT Nat trans} is not evident. We will use the existence of an isomorphism (and its enhancement, \thmref{t:CT compat enh})
in the proof of \thmref{t:main}. 

\begin{rem} \label{r:the scalar}

We repeat that we do not know whether the isomorphisms of Theorems \ref{t:CT compat} and \ref{t:L and CT}
are the same. However, one can show that they differ at most by a (non-zero) scalar.

\medskip

The potential discrepancy between (analogues of) isomorphisms of functors in 
Theorems \ref{t:L and CT} and \ref{t:CT compat} is present also in the quantum geometric Langlands program.

\end{rem} 

\sssec{} 

We form the following cube whose commutativity we analyze in what follows: 

\begin{equation} \label{e:CT cube}
\xy
(0,0)*+{\Dmod_{\frac{1}{2}}(\Bun_G)}="X";
(40,20)*+{\Dmod_{\frac{1}{2}}(\Bun_M)}="Y";
(70,0)*+{ \IndCoh_\Nilp(\LS_\cG)}="X'";
(110,20)*+{\IndCoh_\Nilp(\LS_\cM)}="Y'";
{\ar@{->}^{\on{CT}^-_{*,\rho_P(\omega_X)}[-\delta_{(N_P^-)_{\rho_P(\omega_X)}}]} "X";"Y"};
{\ar@{->}_{\on{CT}^{-,\on{spec}}} "X'";"Y'"};
{\ar@{->}^{\BL_G} "X";"X'"};
{\ar@{->}^{\BL_M} "Y";"Y'"};
(0,-50)*+{\KL(G)_{\crit_G,\Ran}}="Z";
(40,-30)*+{\KL(M)_{\crit_M+\rhoch_P,\Ran^{\subseteq}}}="W";
(70,-50)*+{\IndCoh^*(\Op^\mf_\cG)_\Ran.}="Z'";
(110,-30)*+{\IndCoh^*(\Op^\mf_{\cM,\rhoch_P})_{\Ran^{\subseteq}}}="W'";
{\ar@{->}_{J^{-,\Sph}_{\on{KM},\rho_P(\omega_X)}\circ \on{ins.vac}} "Z";"W"};
{\ar@{->}_{J^{-,\on{spec}}_{\Op,\Theta}\circ \on{ins.unit}} "Z'";"W'"};
{\ar@{->}_{\FLE_{G,\crit}} "Z";"Z'"};
{\ar@{->}_{\FLE_{M,{\crit+\rhoch_P}}} "W";"W'"};
{\ar@{->}^{\Loc_G} "Z";"X"};
{\ar@{->}_{'\!\Loc_M} "W";"Y"};
{\ar@{->}_{\on{Poinc}^{\on{spec}}_{\cG,*}} "Z'";"X'"};
{\ar@{->}_{'\!\on{Poinc}^{\on{spec}}_{\cM,*,\rhoch_P}} "W'";"Y'"};
\endxy
\end{equation} 
where:

\begin{itemize}

\item $'\!\Loc_M:=(\on{Id}\otimes \on{C}^\cdot(\Ran^\subseteq,-))\circ \Loc_{M,\Ran^\subseteq}$;

\item $'\!\on{Poinc}^{\on{spec}}_{\cM,*,\rhoch_P}:=(\on{Id}\otimes \on{C}^\cdot(\Ran^\subseteq,-))\circ \on{Poinc}^{\on{spec}}_{\cM,*,\rhoch_P,\Ran^\subseteq}$. 

\end{itemize}

\medskip

The commutation of its upper lid in \eqref{e:CT cube} is what we eventually want to establish. 
We will see that various results from Part II this paper amount to commutativity of the other five faces up to various shifts and lines. 
We will use these five faces to establish the commutation of its upper lid, see \thmref{t:CT compat} for the precise statement. 

\sssec{} \label{sss:specify shifts and lines}

We now analyze its 2-skeleton:

\begin{itemize}

\item By \thmref{t:local Jacquet}, the bottom lid commutes;

\medskip

\item By \corref{c:CT and Loc crit transl rho}, the left face commutes after tensoring the clockwise circuit with
$$\fl^{\otimes -1}_{N^-_P}\otimes \on{Weil}(\rho_P(\omega_X),\rhoch_P(\omega_X))^{\otimes -1}[\delta_{(N_P^-)_{\rho_P(\omega_X)}}];$$

\item By \eqref{e:Loc compat}, the front face commutes after tensoring the clockwise circuit with
$$\fl^{\otimes \frac{1}{2}}_{G,N_{\rho(\omega_X)}}\otimes \fl^{\otimes -1}_{N_{\rho(\omega_X)}}[-\delta_{N_{\rho(\omega_X)}}];$$

\item By \eqref{e:Langlands and Loc Levi}, the back face commutes after tensoring the clockwise circuit with
$$\fl^{\otimes \frac{1}{2}}_{M,N(M)_{\rho_M(\omega_X)}}\otimes \fl^{\otimes -1}_{N(M)_{\rho_M(\omega_X)}}[-\delta_{N(M)_{\rho_M(\omega_X)}}];$$

\item By \thmref{t:CT opers} and \cite[Theorem 17.4.7]{GLC2}, 
the right face commutes after tensoring the clockwise circuit with $$\fl^{\otimes -1}_{\on{Kost}(\cG)}[\delta_G]$$
and the counterclockwise cirucit with $$\fl^{\otimes -1}_{\on{Kost}(\cM)}[\delta_M],$$ where the line $\fl_{\on{Kost}(\cG)}$
is as in \cite[Sect.17.2.2]{GLC2}, and similarly for $\cM$.

\end{itemize}

\sssec{}

From \secref{sss:specify shifts and lines} above, we obtain that the two circuits in the upper lid in \eqref{e:CT cube} commute
\emph{after} we precompose them with the functor $\Loc_G$, \emph{up to} a cohomological shifts and twist by a constant line.

\medskip

We will now show that this cohomological shift and line are actually trivial. 

\sssec{}

The cancelation of cohomological shifts this amounts to the equality
$$-\delta_{(N_P^-)_{\rho_P(\omega_X)}}-\delta_{N_{\rho(\omega_X)}}+\delta_G=
-\delta_{N(M)_{\rho_M(\omega_X)}}+\delta_M$$
which is equivalent to the valid equality
\begin{multline*} 
\delta_G-\delta_M= \\
= -\dim(\Gamma(X,\fg_{\rho_P(\omega_X)}))+\dim(\Gamma(X,\fm_{\rho_P(\omega_X)}))=
-\dim(\Gamma(X,(\fn_P)_{\rho_P(\omega_X)}))-\dim(\Gamma(X,(\fn^-_P)_{\rho_P(\omega_X)}))= \\
=\delta_{(N_P)_{\rho_P(\omega_X)}}+\delta_{(N_P^-)_{\rho_P(\omega_X)}}=
\delta_{(N^-_P)_{\rho_P(\omega_X)}}+\delta_{N_{\rho(\omega_X)}}-\delta_{N(M)_{\rho_M(\omega_X)}}.
\end{multline*} 

\sssec{}

Let us now show that twists by the constant lines specified in \secref{sss:specify shifts and lines} cancel out. This follows
by combining \corref{c:rel det transl} with the following assertion (and its counterpart for $M$):\footnote{Note that the Killing form on $\fg$
restricts to a non-degenerate form on $z_\fm/z_\fg$, thereby trivializing the line 
$\det(z_\fm/z_\fg)^{\otimes 2}\simeq \det(z_\fm)^{\otimes 2}\otimes \det(z_\fg)^{\otimes -2}$.} 

\begin{prop} \label{p:two lines}
There exists a canonical isomorphism of lines
$$\fl_{\on{Kost}(\cG)}\simeq \fl_{G,N_{\rho(\omega_X)}}\otimes \fl^{-\otimes 2}_{N_{\rho(\omega_X)}}
\otimes \det(\Gamma(X,\CO_X)\otimes \fg) \otimes ((z_\fg)^{\otimes 2})^{\otimes (1-g)}.$$
\end{prop}

\begin{proof}

First, using the fact that $z_\fg$ is the dual of $z_\cg$ and 
$$\fl_{\on{Kost}(\cG)}\simeq \fl_{\on{Kost}(\cG_{\on{ad}})}\otimes \det(\Gamma(X,\CO_X)\otimes z_\cg)$$
and 
$$\det(\Gamma(X,\CO_X)\otimes \fg)\simeq \det(\Gamma(X,\CO_X)\otimes \fg_{\on{ad}})\otimes 
\det(\Gamma(X,\CO_X)\otimes z_\fg),$$
in order to prove the proposition, we can replace $G$ by $G_{\on{ad}}$. Thus, we can assume that $G$ is semi-simple. 
In this case, using the Killing form, we identify $\fa(\cg)$ with $\fa(\fg)$.

\medskip

We have:
$$\fl_{G,N_{\rho(\omega_X)}}\otimes \det(\Gamma(X,\CO_X)\otimes \fg)=
\det(\Gamma(X,\fg_{\rho(\omega_X)}))$$
and
$$\fl_{N_{\rho(\omega_X)}}=\det(\Gamma(X,\fn_{\rho(\omega_X)})).$$

Thus, we need to establish an isomorphism
\begin{equation} \label{e:two lines 1}
\det(\Gamma(X,\fa(\fg)_{\omega_X})) \otimes \det(\Gamma(X,\fn_{\rho(\omega_X)}))
\simeq \det(\Gamma(X,(\fg/\fn)_{\rho(\omega_X)})).
\end{equation}

\medskip

Decompose 
$$\fg\simeq \underset{e\in \BZ^{\geq 0}}\oplus\, V^e.$$
with respect to the action of the principal $SL_2$.\footnote{We emphasize that in the above formula, $e$ is a non-negative inteter,
and not the ``$e$" element in the principal $SL_2$.}
Decompose each $V^e$ into its weight spaces
$$V^e=\underset{n}\oplus\, V^e(n).$$

\medskip

Note that we can identify
$$\fa(\fg)\simeq \underset{e}\oplus\, V^e(e),$$
Note, however, that canonical $\BG_m$-action on $\fa(\fg)$ is shifted by $1$ relative to\footnote{See \cite[Sect. 3.1.9]{BD1}.}
the action of $\BG_m\hookrightarrow SL_2$ on $V^e(e)$, so that 
$$\det(\Gamma(X,\fa(\fg)_{\omega_X})) \simeq \underset{e}\otimes\, 
\det(\Gamma(X,V^e(e)_{\rho(\omega_X)}\otimes \omega_X)).$$

\medskip

To prove \eqref{e:two lines 1}, it suffices to show that for each $e$, we have:
\begin{equation} \label{e:Morozov}
\det(\Gamma(X,V^e(e)_{\rho(\omega_X)}\otimes \omega_X)) \otimes \det(\Gamma(X,\underset{n>0}\oplus V^e(n)_{\rho(\omega_X)}))\simeq
\det(\Gamma(X,\underset{n\leq 0}\oplus V^e(n)_{\rho(\omega_X)})),
\end{equation} 
where we note that 
$$V^e(n)_{\rho(\omega_X)}\simeq V^e(n)\otimes \omega_X^{\otimes n}.$$

\medskip

We will show that
\begin{equation} \label{e:Morozov 1}
\det(\Gamma(X,V^e(e)_{\rho(\omega_X)}\otimes \omega_X)) \simeq 
\det(\Gamma(X,V^e(-e)_{\rho(\omega_X)}))
\end{equation}
and for every $n>0$
\begin{equation} \label{e:Morozov 2}
\det(\Gamma(X,V^e(n)_{\rho(\omega_X)}))\simeq 
\det(\Gamma(X,V^e(-n+2)_{\rho(\omega_X)})).
\end{equation} 

The Killing form identifies $V^e(e)$ with the dual of $V^e(-e)$, so that by Serre
duality
$$\Gamma(X,V^e(e)_{\rho(\omega_X)}\otimes \omega_X)^\vee \simeq \Gamma(X,V^e(-e)_{\rho(\omega_X)})[1].$$
This implies \eqref{e:Morozov 1}.

\medskip

Similarly, for $n>0$, the Killing form and the action of the positive generator of ${\mathfrak{sl}_2}$ identifies
$V^e(n)$ with the dual vector space of $V^e(-n+2)$, and hence by Serre duality
$$\Gamma(X,V^e(n)_{\rho(\omega_X)}))^\vee\simeq \Gamma(X,V^e(-n+2)_{\rho(\omega_X)})[1].$$
This implies \eqref{e:Morozov 2}.

\end{proof}

\begin{rem}
The isomorphism constructed above was canonical. However, it relied on certain conventions. 
For example, we used weight lowering operators $f^i$ at one point, but sometimes one finds weight lowering operators 
$\frac{f^i}{i!}$ are better choices. These sorts of conventions do not affect our main results, but they do affect the scalar discussed in 
Remark \ref{r:the scalar}, and are one reason we are not confident that that scalar is one.
\end{rem} 

\sssec{}

We are now ready to state the main result of this section:

\begin{thm} \label{t:CT compat}
Assuming the validity of GLC for $M$, there exists a unique datum of commutation for the upper lid in \eqref{e:CT cube}, i.e., the square 
\begin{equation} \label{e:L and CT}
\CD
\Dmod_{\frac{1}{2}}(\Bun_G) @>{\BL_G}>> \IndCoh_\Nilp(\LS_\cG) \\
@V{\on{CT}^-_{*,\rho_P(\omega_X)}[-\delta_{(N^-_P)_{\rho_P(\omega_X)}}]}VV @VV{\on{CT}^{-,\on{spec}}}V \\
\Dmod_{\frac{1}{2}}(\Bun_M) @>{\BL_M}>> \IndCoh_\Nilp(\LS_\cM), 
\endCD
\end{equation}
that makes \eqref{e:CT cube} commute (with the specified cohomological shifts and twists by constant lines along the edges).
\end{thm} 

\sssec{}

We now commence the proof of \thmref{t:CT compat}. 

\medskip

The commutation of the five of the faces of the cube \eqref{e:CT cube} established above implies that
the \emph{outer} diagram in 

\medskip

\begin{equation} \label{e:top lid precomposed}
\xy
(0,0)*+{\Dmod_{\frac{1}{2}}(\Bun_G)}="X";
(0,30)*+{\Dmod_{\frac{1}{2}}(\Bun_M)}="Y";
(50,0)*+{ \IndCoh_\Nilp(\LS_\cG)}="X'";
(50,30)*+{\IndCoh_\Nilp(\LS_\cM)}="Y'";
(-50,-22)*+{\KL(G)_{\crit,\Ran}}="W";
{\ar@{->}_{\on{CT}^-_{*,\rho_P(\omega_X)}[-\delta_{(N^-_P)_{\rho_P(\omega_X)}}]} "X";"Y"};
{\ar@{->}_{\on{CT}^{-,\on{spec}}} "X'";"Y'"};
{\ar@{->}^{\BL_G} "X";"X'"};
{\ar@{->}^{\BL_M} "Y";"Y'"};
{\ar@{->}_{\Loc_G} "W";"X"}
{\ar@{->}^{\on{CT}^-_{*,\rho_P(\omega_X)}\circ  \Loc_G[-\delta_{N(P^-)_{\rho_P(\omega_X)}}]} "W";"Y"}
{\ar@{->}_{\BL_G\circ \Loc_G} "W";"X'"}
\endxy
\end{equation} 
is endowed with a commutation datum. 

\medskip

The statement of \thmref{t:CT compat} is equivalent to the fact that this datum comes from 
a uniquely defined commutation datum of the inner square in \eqref{e:top lid precomposed}.

\medskip

For expositional purposes, we first consider the case when $P=B$. 

\ssec{Proof of \thmref{t:CT compat} for \texorpdfstring{$P=B$}{PB}}  \label{ss:CT B}

\sssec{}  \label{sss:CT B}

The category $\Dmod_{\frac{1}{2}}(\Bun_T)\simeq \Dmod(\Bun_T)$ splits as a direct sum according to connected components
of $\Bun_T$, which are indexed by the coweight lattice of $T$. For each coweight $\mu$,
let $\on{CT}^{-,\mu}_{*,\rho_P(\omega_X)}$ denote the corresponding direct summand of $\on{CT}^-_{*,\rho_P(\omega_X)}$. 

\medskip

Let $\on{CT}^{-,\on{spec},\mu}$ denote the corresponding direct summand of $\on{CT}^{-,\on{spec}}$. 
It corresponds to the direct summand $\QCoh(\LS_\cT)^\mu$ of
$$\QCoh(\LS_\cT)=\IndCoh_{\{0\}}(\LS_\cT)=\IndCoh_\Nilp(\LS_\cT)$$
consisting of objects, on which the action of $\cT$ by 1-automorphisms\footnote{I.e., Here and below ``1-automorphisms" of a given object means 
the group of automorphisms of the identity map on this object.} of $\LS_\cT$ has 
character $-\mu$ (here we regard $\mu$ as a weight of $\cT$), see \cite[Sect. 1.5.3]{GLC1}

\medskip

Thus, in proving \thmref{t:CT compat}, instead of the diagram \eqref{e:top lid precomposed}, we can consider

\medskip

\begin{equation} \label{e:top lid precomposed mu}
\xy
(0,0)*+{\Dmod_{\frac{1}{2}}(\Bun_G)}="X";
(0,30)*+{\Dmod_{\frac{1}{2}}(\Bun^\mu_M)}="Y";
(50,0)*+{ \IndCoh_\Nilp(\LS_\cG(X))}="X'";
(50,30)*+{\IndCoh_\Nilp(\LS_\cM(X))^\mu}="Y'";
(-60,-22)*+{\KL(G)_{\crit,\Ran}}="W";
{\ar@{->}_{\on{CT}^{-,\mu}_{*,\rho_P(\omega_X)}[-\delta_{N(P^-)_{\rho_P(\omega_X)}}]} "X";"Y"};
{\ar@{->}_{\on{CT}^{-,\on{spec},\mu}} "X'";"Y'"};
{\ar@{->}^{\BL_G} "X";"X'"};
{\ar@{->}^{\BL^\mu_M} "Y";"Y'"};
{\ar@{->}_{\Loc_G} "W";"X"}
{\ar@{->}^{\on{CT}^{-,\mu}_{*,\rho_P(\omega_X)}\circ  \Loc_G[-\delta_{N(P^-)_{\rho_P(\omega_X)}}]} "W";"Y"}
{\ar@{->}_{\BL_G\circ \Loc_G} "W";"X'"}
\endxy
\end{equation} 
for a fixed $\mu$. 

%
%
%
%

\sssec{}

For a fixed $\mu$, let $\lambda\in \Lambda_G^{+,\BQ}$ be large enough so that the image of the map
$$\Bun^\mu_B\to \Bun_G$$
is contained in the open union of Harder-Narasimhan strata $\Bun_G^{(<\lambda)}$ (see \cite[Sect. 7.4]{DG2}
for our conventions regarding the parameterization of the Harder-Narasimhan strata). 

\medskip

By construction, we have:

\begin{lem} \label{l:Eis doesnt matter}
The functor $\on{CT}^{-,\mu}_{*,\rho_P(\omega_X)}$ factors as
$$\Dmod_{\frac{1}{2}}(\Bun_G)\twoheadrightarrow \Dmod_{\frac{1}{2}}(\Bun_G^{(<\lambda)}) 
\overset{(\on{CT}^{-,\mu}_{*,\rho_P(\omega_X)})^{(<\lambda)}}\longrightarrow 
\Dmod_{\frac{1}{2}}(\Bun_T),$$
where the first arrow is the restriction functor.
\end{lem} 

\sssec{}

Let $P'$ be a standard parabolic in $G$ with Levi quotient $M'$. Recall that $\Lambda_{G,P'}$
denotes the quotient of $\Lambda$ by the root lattice of $M'$, i.e.,
$$\Lambda_{G,P'}\simeq \pi_{1,\on{alg}}(M')\simeq \pi_0(\Bun_{M'}).$$

\medskip

Recall (see \cite[Sects. 7.1.3-7.1.5]{DG2}) that we can view $\Lambda_{G,P'}$ as a \emph{subset} of $\Lambda_G^{\BQ}$.
Denote $$\Lambda^+_{G,P'}:=\Lambda_{G,P'}\cap \Lambda_G^{+,\BQ}.$$

\sssec{}

Let
$$\Bun_{M'}^{\underset{G}{\not<}\lambda}\subset \Bun_{M'}$$
be the union of connected components, indexed by coweights $\lambda'\in \Lambda^+_{G,P'}$
with $$\lambda'\underset{G}{\not<} \lambda.$$

\medskip

Note that the functor 
$$\Eis_!:\Dmod_{\frac{1}{2}}(\Bun_{M'})\to \Dmod_{\frac{1}{2}}(\Bun_G)$$
induces a functor
\begin{equation} \label{e:Eis < lambda}
\Dmod_{\frac{1}{2}}(\Bun_{M'}^{\underset{G}{\not<}\lambda})\to \Dmod_{\frac{1}{2}}(\Bun_G^{(\not<\lambda)}),
\end{equation} 
where 
$$\Bun_G^{(\not<\lambda)}\subset \Bun_G$$
is the closed substack equal to the complement of $\Bun_G^{(<\lambda)}$.

\sssec{}

We claim:

\begin{lem} \label{l:HS gen}
For $\lambda$ large enough, the essential image of \eqref{e:Eis < lambda} generates the target.
\end{lem} 

\begin{proof}

Let $\lambda$ satisfy
$$\langle \lambda,\check\alpha_i\rangle \geq \on{max}(0,2g-2)$$
for all simple roots $\check\alpha_i$.

\medskip

Then according to \cite[Sect. 3]{DG2}, we can partition the set $\Lambda_G^{(\not<\lambda)}$ into locally closed subsets\footnote{Here ``locally closed"
refers to the order relation $\underset{G}\leq$ on $\Lambda_G$, see \cite[Appendix A]{DG2}.} $S$, such that for 
$$\Bun_G^{(S)}:=\underset{\mu\in S}\bigcup\, \Bun_G^{(\mu)}$$
the following holds: 

\medskip

For every $S$ there exists a parabolic $P'$ such that for
$$\Bun_{P'}^{(S)}:=\Bun_{P'}\underset{\Bun_{M'}}\times \Bun^{(S)}_{M'},$$

\begin{itemize}

\item The map $\sfp:\Bun_{P'}^{(S)}\to \Bun_G^{(S)}$ is an isomorphism;

\item The map $\sfq:\Bun_{P'}^{(S)}\to \Bun_M^{(S)}$ is a unipotent gerbe.\footnote{See \cite[Definition 10.3.5]{DG1} for what this means.}

\end{itemize}

\medskip

It is clear that for $S$ as above, the functor $\Eis_!$ induces an equivalence
$$\Dmod_{\frac{1}{2}}(\Bun_{M'}^{(S)})\to \Dmod_{\frac{1}{2}}(\Bun_G^{(S)}).$$

\end{proof} 

\sssec{}

Let $\lambda$ be sufficiently large. Applying \lemref{l:HS gen}, 
we can identify $\Dmod_{\frac{1}{2}}(\Bun_G^{(<\lambda)})$ with the quotient of $\Dmod_{\frac{1}{2}}(\Bun_G)$ by the full subcategory
generated by the essential images of $\Dmod_{\frac{1}{2}}(\Bun_{M'}^{\underset{G}{\not<}\lambda})$ along the Eisenstein functors
$$\Eis_!:\Dmod_{\frac{1}{2}}(\Bun_{M'})\to \Dmod_{\frac{1}{2}}(\Bun_G)$$
for all standard proper parabolics $P'$. 

\sssec{} \label{sss:quot of LS}

For $\lambda'\in  \Lambda_{G,P'}$, let $\IndCoh_\Nilp(\LS_{\cM'})^{\lambda'}$ be the direct summand of 
$\IndCoh_\Nilp(\LS_{\cM'})$ consisting of objects on which the action of $Z_{\cM'}$ by 1-automorphisms\footnote{See footnote
in \secref{sss:CT B} for the meaning of ``1-automorphisms".} of $\LS_{\cM'}$ has 
character $-\lambda'$. 

\medskip

Let $\IndCoh_\Nilp(\LS_\cG)^{<(\lambda)}$ denote the quotient of $\IndCoh_\Nilp(\LS_\cG)$ by the full
subcategory generated 
by the essential images of\footnote{The shift by $2(g-1)\cdot \rho_{P'}$ in the formula below is due to the fact that 
on the geometric side in \thmref{t:L and Eis}, we are dealing with the functor $\Eis^-_{!,\rho_P(\omega_X)}$ rather than just $\Eis^-_!$.}  
$$\IndCoh_\Nilp(\LS_{M'})^{\lambda'}, \quad \lambda'\not<\lambda+2(g-1)\cdot \rho_{P'}.$$

\medskip

We will prove:

\begin{prop} \label{p:Eis spec doesnt matter}
For a fixed $\mu$, and $\lambda$ large enough, for every standard parabolic $P'$ and 
$\lambda'\in \Lambda_{G,P'}$ satisfying $\lambda'\not<\lambda$, the functor
$$\IndCoh_\Nilp(\LS_{\cM'})^{\lambda'}\overset{\Eis^{\on{spec}}}\longrightarrow 
\IndCoh_\Nilp(\LS_\cG)\overset{\on{CT}^{-,\on{spec},\mu}}\longrightarrow \QCoh(\LS_\cT)^\mu$$
vanishes.
\end{prop} 

\sssec{}

Assuming \propref{p:Eis spec doesnt matter} temporarily, we obtain that for $\lambda$ sufficiently large, 
the functor $\on{CT}^{-,\on{spec},\mu}$ factors as
$$\IndCoh_\Nilp(\LS_\cG)\twoheadrightarrow \IndCoh_\Nilp(\LS_\cG)^{(<\lambda)}
\overset{(\on{CT}^{-,\on{spec},\mu})^{(<\lambda)}}\longrightarrow \QCoh(\LS_\cT)^\mu.$$

\sssec{}

The compatibility of the Langlands functor with the Eisenstein functors given by \thmref{t:L and Eis}
implies that the functor $\BL_G$
descends to a well-defined functor
$$\BL_G^{(<\lambda)}:\Dmod_{\frac{1}{2}}(\Bun_G^{(<\lambda)})\to \IndCoh_\Nilp(\LS_\cG)^{(<\lambda)}.$$

We obtain that the commutativity datum for \eqref{e:top lid precomposed mu} is equivalent to that of
the commutativity datum for the inner square in 

\medskip

\begin{equation} \label{e:top lid precomposed mu mod lambda}
\xy
(0,0)*+{\Dmod_{\frac{1}{2}}(\Bun^{(<\lambda)}_G)}="X";
(0,30)*+{\Dmod_{\frac{1}{2}}(\Bun^\mu_M)}="Y";
(50,0)*+{ \IndCoh_\Nilp(\LS_\cG)^{(<\lambda)}}="X'";
(50,30)*+{\IndCoh_\Nilp(\LS_\cM)^\mu,}="Y'";
(-60,-40)*+{\KL(G)_{\crit,\Ran}}="W";
(-30,-20)*+{\Dmod_{\frac{1}{2}}(\Bun_G)}="Z";
{\ar@{->}_{(\on{CT}^{-,\mu}_{*,\rho_P(\omega_X)})^{(<\lambda)}[-\delta_{N(P^-)_{\rho_P(\omega_X)}}]} "X";"Y"};
{\ar@{->}_{(\on{CT}^{-,\on{spec},\mu})^{(<\lambda)}} "X'";"Y'"};
{\ar@{->}^{\BL^{(<\lambda)}_G} "X";"X'"};
{\ar@{->}^{\BL^\mu_M} "Y";"Y'"};
{\ar@{->}_{\Loc_G} "W";"Z"}
{\ar@{->}"W";"Y"}
{\ar@{->} "W";"X'"}
{\ar@{->} "Z";"X"}
\endxy
\end{equation} 
compatible with the existing commutativity datum for the outer diagram. 

%
%
%
%
%
%
%

\medskip

Furthermore, we can replace $\Dmod_{\frac{1}{2}}(\Bun^{(<\lambda)}_G)$ by its connected component 
$\Dmod_{\frac{1}{2}}(\Bun^{(<\lambda),\ol\mu}_G)$
corresponding to the
image $\ol\mu$ of $\mu$ under
$$\Lambda\to \Lambda/\Lambda_{\on{s.c.}}\simeq \pi_{1,\on{alg}}(G)\simeq \pi_0(\Bun_G).$$

\medskip

Now, the required assertion follows now from the next observation:

\begin{lem} \label{l:qc quot}
The functor
$$\KL(G)_{\crit,\Ran}\overset{\Loc_G}\to \Dmod_{\frac{1}{2}}(\Bun_G)\twoheadrightarrow \Dmod_{\frac{1}{2}}(\Bun^{(<\lambda),\ol\mu}_G)$$
is a Verdier quotient. 
\end{lem} 

\begin{proof}

This is a particular case of \cite[Theorem 13.4.2]{GLC2}. 

\end{proof}

\ssec{Proof of \propref{p:Eis spec doesnt matter}}

For the duration of the proof, we will change the notation from $P'$ to $P$. 

\sssec{}

By base change, the functor  
$$\IndCoh(\LS_\cM)\overset{\Eis^{\on{spec}}}\longrightarrow 
\IndCoh(\LS_\cG)\overset{\on{CT}^{-,\on{spec}}}\longrightarrow \IndCoh(\LS_\cT)$$
can be rewritten as the composition of:

\smallskip

\begin{itemize}

\item *-pullback along $\LS_\cP\to \LS_\cM$;

\smallskip

\item !-pullback along $\LS_\cP\underset{\LS_\cG}\times \LS_{\cB^-}\to \LS_\cP$;

\smallskip

\item *-pushforward along $\LS_\cP\underset{\LS_\cG}\times \LS_{\cB^-}\to \LS_{\cB^-}\to \LS_\cT$.

\end{itemize}

However, since the morphism $\LS_\cP\to \LS_\cM$ is quasi-smooth, up to shifting the degree, we can replace
the *-pullback by the !-pullback. So the functor in question becomes !-pullback followed by *-pushforward along the diagram
$$
\xy
(0,0)*+{\LS_\cM}="X";
(40,0)*+{\LS_\cT}="Y";
(20,20)*+{\LS_\cP\underset{\LS_\cG}\times \LS_{\cB^-}.}="Z";
{\ar@{->} "Z";"X"};
{\ar@{->} "Z";"Y"};
\endxy
$$

\sssec{}

We decompose $\LS_\cP\underset{\LS_\cG}\times \LS_{\cB^-}$ according to relative positions of the two
reductions, which are indexed by the elements of $W/W_M$. 

\medskip

For each $w\in W/W_M$, let
$$(\LS_\cP\underset{\LS_\cG}\times \LS_{\cB^-})_w\subset \LS_\cP\underset{\LS_\cG}\times \LS_{\cB^-}$$
denote the corresponding locally closed substack, and let 
$$(\LS_\cP\underset{\LS_\cG}\times \LS_{\cB^-})_w^\wedge$$
denote its formal completion inside $\LS_\cP\underset{\LS_\cG}\times \LS_{\cB^-}$. 

\medskip

We will show that for every $w$ and $\lambda$ large enough, the pull-push functor along

\begin{equation} \label{e:geom lemma compl}
\xy
(0,0)*+{\LS_\cM}="X";
(40,0)*+{\LS_\cT}="Y";
(20,20)*+{(\LS_\cP\underset{\LS_\cG}\times \LS_{\cB^-})_w^\wedge}="Z";
{\ar@{->} "Z";"X"};
{\ar@{->} "Z";"Y"};
\endxy
\end{equation} 
has the property that its $(\lambda',\mu)$ component vanishes for 
\begin{equation} \label{e:lambda ineq}
\lambda'\not< \lambda, \quad \lambda'\in \Lambda^+_{G,P}.
\end{equation} 

\medskip

We will first establish the corresponding fact for the diagram 
\begin{equation} \label{e:geom lemma w}
\xy
(0,0)*+{\LS_\cM}="X";
(40,0)*+{\LS_\cT.}="Y";
(20,20)*+{(\LS_\cP\underset{\LS_\cG}\times \LS_{\cB^-})_w}="Z";
{\ar@{->} "Z";"X"};
{\ar@{->} "Z";"Y"};
\endxy
\end{equation} 

\sssec{}

Let $w\in W$ the shortest representative of the given coset. Consider the subgroup
$$\sH_w:=\cP\cap w^{-1}(\cB^-)\overset{w}\simeq w(\cP)\cap \cB^-.$$

\medskip

We can identify
$$(\LS_\cP\underset{\LS_\cG}\times \LS_{\cB^-})_w\simeq \LS_{\sH_w},$$
so that \eqref{e:geom lemma w} identifies with 
\begin{equation} \label{e:geom lemma w 1}
\xy
(0,0)*+{\LS_\cM}="X";
(40,0)*+{\LS_\cT,}="Y";
(20,20)*+{\LS_{\sH_w}}="Z";
{\ar@{->} "Z";"X"};
{\ar@{->} "Z";"Y"};
\endxy
\end{equation} 
where the two projections correspond to
\begin{equation} \label{e:H to P}
\sH_w=\cP\cap w^{-1}(\cB^-)\to \cP\to \cM
\end{equation}
\begin{equation} \label{e:H to B-}
\sH_w=\cP\cap w^{-1}(\cB^-)\overset{w}\simeq w(\cP)\cap \cB^-\to \cB^-\to \cT,
\end{equation}
respectively.

\sssec{} \label{sss:Borel case start}

Note that diagram \eqref{e:geom lemma w 1} naturally factors as
\begin{equation} \label{e:geom lemma simple}
\xy
(0,0)*+{\LS_\cM}="X";
(40,0)*+{\LS_\cT.}="Y";
(20,20)*+{\LS_{\cB^-(M)}}="W";
(20,40)*+{\LS_{\sH_w}}="Z";
{\ar@{->}_{\sfp^{-,\on{spec}}_M} "W";"X"};
{\ar@{->}^{w\circ \sfq^{-,\on{spec}}_M} "W";"Y"};
{\ar@{->}^{h} "Z";"W"};
\endxy
\end{equation} 

\medskip

Pull-push along \eqref{e:geom lemma w 1} identifies with
\begin{equation} \label{e:pull push interm}
(w\circ \sfq^{-,\on{spec}}_M)^\IndCoh_*\left((\sfp^{-,\on{spec}}_M)^!(-)\sotimes h^\IndCoh_*(\omega_{\LS_{\sH_w}})\right).
\end{equation} 

\sssec{}

Note that the stack $\LS_{\sH_w}$ and the map $h$ are quasi-smooth 
(the latter because $H_w\to \cB^-(M)$ is surjective\footnote{Indeed, $H_w\cap \cM=\cB^-(M)$.}). 
Hence, up to shift by a fixed weight, we can replace 
\eqref{e:pull push interm} by
\begin{equation} \label{e:pull push interm O}
(w\circ \sfq^{-,\on{spec}}_M)^\IndCoh_*\left((\sfp^{-,\on{spec}}_M)^!(-)\otimes h_*(\CO_{\LS_{\sH_w}})\right),
\end{equation} 
where:

\begin{itemize}

\item $h_*(\CO_{\LS_{\sH_w}})$ is viewed as an object of $\QCoh(\LS_{\cB^-(M)})$;

\smallskip

\item $\otimes$ denotes the tensor product functor $\IndCoh(-)\otimes \QCoh(-)\to \IndCoh(-)$.

\end{itemize} 

\sssec{} \label{sss:positivity}

Consider the action of $Z_\cM$ by 1-automorphisms on the stacks in 
$$
\xy
(0,0)*+{\LS_\cM}="X";
(40,0)*+{\LS_\cT.}="Y";
(20,20)*+{\LS_{\cB^-(M)}}="Z";
{\ar@{->}_{\sfp^{-,\on{spec}}_M}  "Z";"X"};
{\ar@{->}^{w\circ \sfq^{-,\on{spec}}_M}  "Z";"Y"};
\endxy
$$

The characters of $Z_\cM$ that appear in $h_*(\CO_{\LS_{\sH_w}})$ 
are sums of characters of the form
\begin{equation} \label{e:which alpha}
-\alpha,\quad \alpha \text{ is a root in } \cn_P\cap w^{-1}(\cn^-).
\end{equation}

\medskip

Hence, for $\CF\in \IndCoh_\Nilp(\LS_\cM)^{\lambda'}$, the characters of $Z_\cM$
on the pull-push of $Z_\cM$ along \eqref{e:geom lemma w} are of the form
$$-w(\lambda' + \underset{\alpha}\Sigma\, n_\alpha\cdot \alpha),$$
for $\alpha$ in \eqref{e:which alpha} and $n_\alpha\in \BZ^{\geq 0}$. 

\medskip

Hence, for the $\mu$-direct summand of the above pull-push to be non-zero, we must have
$$-w(\lambda' + \underset{\alpha}\Sigma\, n_\alpha\cdot \alpha)=\mu|_{Z_\cM}.$$

I.e., 
\begin{equation} \label{e:lambda ineq bis}
\lambda' + \underset{\alpha}\Sigma\, n_\alpha\cdot \alpha=-w^{-1}(\mu)|_{Z_\cM}.
\end{equation}

\sssec{} \label{sss:positivity bis}

However, \eqref{e:lambda ineq bis} is impossible once \eqref{e:lambda ineq} is satisfied with $\lambda$ large enough. 

\medskip

Indeed, choose $\lambda$ so that
\begin{equation} \label{e:lambda how large}
\langle \lambda,\check\omega_i\rangle > \langle -w^{-1}(\mu),\check\omega_i\rangle
\end{equation}
for all fundamental coweights $\check\omega_i$.

\medskip

If \eqref{e:lambda ineq} is satisfied, there exists a fundamental coweight $\check\omega_i$ central in $\cM$, 
such that
$$\langle \lambda',\check\omega_i\rangle \geq \langle \lambda,\check\omega_i\rangle.$$

Pairing the left hand side of \eqref{e:lambda ineq bis} with this $\check\omega_i$, we obtain
$$\langle \lambda',\check\omega_i\rangle +  \underset{\alpha}\Sigma\, n_\alpha\cdot \langle \alpha,\check\omega_i\rangle\geq 
\langle \lambda',\check\omega_i\rangle \geq 
\langle \lambda,\check\omega_i\rangle>\langle -w^{-1}(\mu),\check\omega_i\rangle,$$
contradicting \eqref{e:lambda ineq bis}.

\sssec{}  \label{sss:Borel case end}

We now prove the assertion for the pull-push along \eqref{e:geom lemma compl}. This functor admits
a filtration with subquotients of the form
$$(w\circ \sfq^{-,\on{spec}}_M)^\IndCoh_*\left((\sfp^{-,\on{spec}}_M)^!(-)\sotimes 
h^\IndCoh_*\left(\omega_{\LS_{\sH_w}}\otimes \Sym(\on{Norm}_w)\right)\right),$$
where $\on{Norm}_w$ is the normal bundle to 
$\LS_{\sH_w}$ inside $\LS_\cP\underset{\LS_\cG}\times \LS_{\cB^-}$. 

\medskip

The fiber of this normal at a given $\sigma\in \LS_{\sH_w}$ is
$$\on{Cofib}\left(\on{C}^\cdot(X,(\sh_w)_\sigma)[1]\to \on{C}^\cdot(X,(\cp\underset{\cg}\times w^{-1}(\cb^-))_\sigma)[1]\right)
\simeq \on{C}^\cdot\biggl(X,\biggl(\cg/(\cp+w^{-1}(\cb^-))\biggr)_\sigma\biggr)[1].$$

Hence, the characters of $Z_\cM$ on 
$$h^\IndCoh_*\left(\omega_{\LS_{\sH_w}}\otimes \Sym(\on{Norm}_w)\right)$$
are sums of characters of the form \eqref{e:which alpha} and also those of the form 
\begin{equation} \label{e:which beta}
\beta,\quad \beta \text{ is a root in } \cg/(\cp+w^{-1}(\cb^-)).
\end{equation}

\medskip

The proof is concluded by the same 
argument as in Sects. \ref{sss:positivity}-\ref{sss:positivity bis}. 

\qed[\propref{p:Eis spec doesnt matter}]

\ssec{Proof of \thmref{t:CT compat} for a general Levi}

\sssec{}

Recall that we are assuming the validity of the geometric Langlands conjecture for M.

\begin{rem}

What we really need to assume for the proof to go through 
is a certain property of the category of $\IndCoh_\Nilp(\LS_{\cM})$,
see \secref{sss:IndCoh as lim}. This property takes place purely on the spectral side,
and it follows from GLC. 

\medskip

Let us formulate this property for $\cG$. For $\lambda\in \Lambda_G^{+,\BQ}$, consider the Verdier quotient category
$$\IndCoh_\Nilp(\LS_\cG)\twoheadrightarrow \IndCoh_\Nilp(\LS_\cG)^{(<\lambda)}, \quad \lambda\in \Lambda_G^{+,\BQ},$$
where we kill the subcategory generated by the essential images of the functors
$$\Eis^{\on{spec}}:\IndCoh_\Nilp(\LS_{\cM'})^{\lambda'}\to \IndCoh_\Nilp(\LS_\cG), \quad 
\lambda'\in \Lambda^+_{G,P'},\,\, \lambda'\underset{G}{\not<} \lambda-2(g-1)\cdot \rho_{P'}.$$
for all standard parabolics $P'$ of $G$.

\medskip

What we need is that the functor
$$\IndCoh_\Nilp(\LS_\cG) \to \underset{\lambda, \leq_G}{\on{lim}}\, \IndCoh_\Nilp(\LS_\cG)^{(<\lambda)}$$
is an equivalence. 

\medskip

With this property at hand, we will imitate the argument from \secref{ss:CT B} essentially word for word. 

\end{rem} 

\begin{rem}

For $\lambda\in \Lambda_G^{+,\BQ}$, denote by  $\QCoh(\LS_\cG)^{(<\lambda)}$ the corresponding quotient of $\QCoh(\LS_\cG)$
so that we have a commutative diagram.
\begin{equation} \label{e:QCoh lambda}
\CD
\QCoh(\LS_\cG) @<{\Xi_{0,\Nilp}^R}<<  \IndCoh_\Nilp(\LS_\cG) \\
@VVV @VVV \\
\QCoh(\LS_\cG)^{(<\lambda)} @<{\Xi_{0,\Nilp}^R}<< \IndCoh_\Nilp(\LS_\cG)^{(<\lambda)},
\endCD
\end{equation} 

Assuming GLC for $G$, it follows from the localization argument given below that the category
$\IndCoh_\Nilp(\LS_\cG)^{(<\lambda)}$ is generated by the essential image of 
$$\QCoh(\LS_\cG)\overset{\Xi_{0,\Nilp}}\to \IndCoh_\Nilp(\LS_\cG)\to \IndCoh_\Nilp(\LS_\cG)^{(<\lambda)}.$$

This implies that the bottom horizontal arrow in \eqref{e:QCoh lambda} is actually an equivalence. 

\medskip

In particular, we obtain that the category $\IndCoh_\Nilp(\LS_\cG)$ can be recovered 
from the usual $\QCoh(\LS_\cG)$ 
also as
$$\underset{\lambda, \leq_G}{\on{lim}}\, \QCoh(\LS_\cG)^{(<\lambda)}.$$

\end{rem} 

\sssec{}

Fix $\mu\in \Lambda_M^+$, and let 
$$\Bun_M^{(<\mu)}\subset \Bun_M$$
be the quasi-compact open equal to the union of Harder-Narasimhan strata $\Bun_M^{(\mu')}$ with
$$\mu'\underset{M}{<} \mu.$$

\medskip

We consider $\Dmod_{\frac{1}{2}}(\Bun_M^{(<\mu)})$ as a quotient of $\Dmod_{\frac{1}{2}}(\Bun_M)$. 

\medskip

Assuming that $\mu$ is 
dominant enough (as a coweight of $M$), by \lemref{l:HS gen}, the kernel of the projection
$$\Dmod_{\frac{1}{2}}(\Bun_M)\to \Dmod_{\frac{1}{2}}(\Bun_M^{(<\mu)})$$
is generated by the essential images of 
$$\Dmod_{\frac{1}{2}}(\Bun^{\mu'}_{M'}), \quad \mu'\in \Lambda^+_{G,P'}, 
\quad \mu'\underset{M}{\not<} \mu+2(g-1)\cdot\rho_{P'}$$
along the functors
$$\Eis_{!,-\rho_{P'}(\omega_X)}:\Dmod_{\frac{1}{2}}(\Bun_{M'})\to \Dmod_{\frac{1}{2}}(\Bun_M),$$
where $M'$ is the Levi of a standard parabolic $P'$ of $M$.

\sssec{} \label{sss:trunc LS M}

Let $\IndCoh_\Nilp(\LS_\cM)^{(<\mu)}$ denote the quotient of $\IndCoh_\Nilp(\LS_\cM)$ by the full subcategory 
generated by the essential images of 
$$\IndCoh_\Nilp(\LS_{\cM'})^{\mu'}, \quad \mu'\in \Lambda^+_{G,P'}, 
\quad \mu'\underset{M}{\not<} \mu+2(g-1)\cdot\rho_{P'}$$
along the functors
$$\Eis^{\on{spec}}:\IndCoh_\Nilp(\LS_{\cM'})\to \IndCoh_\Nilp(\LS_\cM).$$

\sssec{}

The compatibility of the Langlands functor $\BL_M$ with Eisenstein series implies that there exists a commutative diagram
$$
\CD
\Dmod_{\frac{1}{2}}(\Bun_M^{(<\mu)}) @>{\BL_M^{(<\mu)}}>{\sim}> \IndCoh_\Nilp(\LS_\cM)^{(<\mu)} \\
@AAA @AAA \\
\Dmod_{\frac{1}{2}}(\Bun_M) @>{\sim}>{\BL_M}> \IndCoh_\Nilp(\LS_\cM). 
\endCD
$$

\sssec{} \label{sss:IndCoh as lim}

Since
$$\Dmod_{\frac{1}{2}}(\Bun_M)\to \underset{\mu}{\on{lim}}\, \Dmod_{\frac{1}{2}}(\Bun_M^{(<\mu)})$$
is an equivalence, we obtain that 
$$\IndCoh_\Nilp(\LS_\cM)\to \underset{\mu}{\on{lim}}\, \IndCoh_\Nilp(\LS_\cM)^{(<\mu)}$$
is also an equivalence. 

\medskip

Hence, in order to construct a datum of commutativity for \eqref{e:top lid precomposed},
it enough to construct a compatible data of commutativity for the diagrams (for varying $\mu$)
\begin{equation} \label{e:top lid precomposed U}
\xy
(0,0)*+{\Dmod_{\frac{1}{2}}(\Bun_G)}="X";
(0,30)*+{\Dmod_{\frac{1}{2}}(\Bun_M^{(<\mu)})}="Y";
(0,15)*+{\Dmod_{\frac{1}{2}}(\Bun_M)}="Y1";
(50,0)*+{ \IndCoh_\Nilp(\LS_\cG)}="X'";
(50,30)*+{\IndCoh_\Nilp(\LS_\cM)^{(<\mu)}}="Y'";
(50,15)*+{\IndCoh_\Nilp(\LS_\cM)}="Y1'";
(-50,-30)*+{\KL(G)_{\crit,\Ran},}="W";
{\ar@{->}_{\on{CT}^-_{*,\rho_P(\omega_X)}[-\delta_{N(P^-)_{\rho_P(\omega_X)}}]} "X";"Y1"};
{\ar@{->}_{\on{CT}^{-,\on{spec}}} "X'";"Y1'"};
{\ar@{->}^{\BL_G} "X";"X'"};
{\ar@{->}^{\BL^{(<\mu)}_M} "Y";"Y'"};
{\ar@{->}_{\Loc_G} "W";"X"}
{\ar@{->}"W";"Y"}
{\ar@{->} "W";"X'"}
{\ar@{->} "Y1";"Y"}
{\ar@{->} "Y1'";"Y'"}
\endxy
\end{equation} 
compatible with the given data of commutativity for the outer diagram. 

\sssec{}

Let $\Bun^{(<\lambda)}_G \subset \Bun_G$ be a quasi-compact open union of Harder-Narasimhan strata, such that the functor
$$\Dmod_{\frac{1}{2}}(\Bun_G)\overset{\on{CT}^-_{*,\rho_P(\omega_X)}}\longrightarrow 
\Dmod_{\frac{1}{2}}(\Bun_M)\to \Dmod_{\frac{1}{2}}(\Bun_M^{(<\mu)})$$
factors via the quotient
$$\Dmod_{\frac{1}{2}}(\Bun_G)\twoheadrightarrow \Dmod_{\frac{1}{2}}(\Bun^{(<\lambda)}_G).$$

Denote the resulting functor 
$$\Dmod_{\frac{1}{2}}(\Bun^{(<\lambda)}_G)\to \Dmod_{\frac{1}{2}}(\Bun_M^{(<\mu)})$$
by $(\on{CT}^-_{*,\rho_P(\omega_X)})^{(<\mu),(<\lambda)}$. 

\sssec{}

Let 
$$ \IndCoh_\Nilp(\LS_\cG)\twoheadrightarrow  \IndCoh_\Nilp(\LS_\cG)^{(<\lambda)}$$
denote the corresponding quotient, see \secref{sss:quot of LS}.

\bigskip

The following is a generalization of \propref{p:Eis spec doesnt matter}: 

\begin{prop} \label{p:kill Eis again}
For a fixed $\mu$ and $\lambda$ large enough, for every standard parabolic $P'$ and 
$\lambda'\in \Lambda_{G,P'}$ satisfying $\lambda'\underset{G}{\not<} \lambda$, the functor
$$\IndCoh_\Nilp(\LS_{\cM'})^{\lambda'}\overset{\Eis^{\on{spec}}}\longrightarrow  
\IndCoh_\Nilp(\LS_\cG)\overset{\on{CT}^{-,\on{spec}}}\longrightarrow \IndCoh_\Nilp(\LS_\cM)\to \IndCoh_\Nilp(\LS_\cM)^{(<\mu)}$$
vanishes.
\end{prop}

Let us assume this proposition for a moment and finish the proof of \thmref{t:CT compat}. 

\sssec{}

From \propref{p:kill Eis again}, combined with \thmref{t:L and Eis} and \lemref{l:HS gen}, we obtain: 

\begin{cor} \label{c:kill Eis again}
For $\lambda$ large enough, 
the composite functor
$$\Dmod_{\frac{1}{2}}(\Bun_G)\overset{\BL_G}\longrightarrow  \IndCoh_\Nilp(\LS_\cG)
\overset{\on{CT}^{-,\on{spec}}}\longrightarrow \IndCoh_\Nilp(\LS_\cM)\to \IndCoh_\Nilp(\LS_\cM)^{(<\mu)}$$
also factors via the quotient
$$\Dmod_{\frac{1}{2}}(\Bun_G)\twoheadrightarrow \Dmod_{\frac{1}{2}}(\Bun^{(<\lambda)}_G).$$
\end{cor}

\sssec{}

Denote the resulting functor
$$\Dmod_{\frac{1}{2}}(\Bun^{(<\lambda)}_G)\to \IndCoh_\Nilp(\LS_\cM)^{(<\mu)}$$ 
by
$$(\on{CT}^{-,\on{spec}}\circ \BL_G)^{(<\mu),(<\lambda)}.$$

We obtain that a datum of commutativity for \eqref{e:top lid precomposed U}
is equivalent to that for the diagram

\medskip

\begin{equation} \label{e:curved diag1}
\xy
(0,0)*+{\Dmod_{\frac{1}{2}}(\Bun^{(<\lambda)}_G)}="X";
(30,30)*+{\IndCoh_\Nilp(\LS_\cM)^{(<\mu)}}="Y'";
(-15,-15)*+{\Dmod_{\frac{1}{2}}(\Bun_G)}="Z";
(-30,-30)*+{\KL(G)_{\crit,\Ran}}="W";
{\ar@{->}@/^3pc/^{\BL^{(<\mu)}_M\circ (\on{CT}^-_{*,\rho_P(\omega_X)})^{(<\mu),(<\lambda)}[-\delta_{N(P^-)_{\rho_P(\omega_X)}}]} "X";"Y'"};
{\ar@{->}@/_3pc/_{(\on{CT}^{-,\on{spec}}\circ \BL_G)^{(<\mu),(<\lambda)}} "X";"Y'"};
{\ar@{->}^{\Loc_G}  "W";"Z"};
{\ar@{->} "Z";"X"}
\endxy
\end{equation}

However, this follows again from \lemref{l:qc quot}.

\ssec{Proof of \propref{p:kill Eis again}}

The proof proceeds along the same lines as that of \propref{p:Eis spec doesnt matter}, using a generalization
of diagram \eqref{e:geom lemma simple}, explained below.

\sssec{}

Let $P_1$ and $P_2$
be a pair of standard parabolics of $G$ with Levi quotients $M_1$ and $M_2$, respectively. For an element
$$w\in W_1\backslash W/W_2,$$
let 
$$(\LS_{\cP_1}\underset{\LS_{\cG}}\times \LS_{\cP_2})_w\subset \LS_{\cP_1}\underset{\LS_{\cG}}\times \LS_{\cP_2}$$
be the corresponding locally closed substack.

\medskip

Let $w$ be the shortest representative of the given double coset. Consider the subgroup
$$\sH_w:=\cP_1\cap w^{-1}(\cP_2) \overset{w}\simeq w(\cP_1)\cap \cP_2,$$
so that we can identify
$$(\LS_{\cP_1}\underset{\LS_{\cG}}\times \LS_{\cP_2})_w\simeq \LS_{\sH_w}.$$

\sssec{}

Set
$$\cP'_1:=\cP_1\cap w^{-1}(\cP_2)/(\cN_{P_1}\cap w^{-1}(\cP_2)) \text{ and }
\cP'_2:=w(\cP_1)\cap \cP_2/(w(\cP_1)\cap \cN_{P_2}).$$

The projections
$$\cP'_1\to \cP_1/\cN_{P_1}=\cM_1 \text{ and } \cP'_2\to \cP_2/\cN_{P_2}=\cM_2$$
realize $\cP'_1$ and $\cP'_2$ as standard parabolics in $\cM_1$ and $\cM_2$, respectively.

\medskip

The element $w$ gives rise to an identification of the Levi quotients
$$\cM'_1\overset{w}\simeq \cM'_2$$
of $\cP'_1$ and $\cP'_2$, respectively.

\medskip

Finally, note also that we have a Cartesian diagram
$$
\xy
(0,0)*+{\LS_{\cP_1\cap w^{-1}(\cP_2)/(N(\cP_1)\cap w^{-1}(N(\cP_2)))}\overset{w}\sim \LS_{w(\cP_1)\cap \cP_2/(w(N(\cP_1))\cap N(\cP_2))}}="W";
(-20,-20)*+{\LS_{\cP'_1}}="Z_1";
(20,-20)*+{\LS_{\cP'_2}}="Z_2";
(0,-40)*+{\LS_{\cM'_1}\overset{w}\sim \LS_{\cM_2'}.}="R";
{\ar@{->} "W";"Z_1"};
{\ar@{->} "W";"Z_2"};
{\ar@{->} "Z_1";"R"};
{\ar@{->} "Z_2";"R"};
\endxy
$$

\sssec{}

Then the diagram 
$$
\xy
(0,0)*+{(\LS_{\cP_1}\underset{\LS_{\cG}}\times \LS_{\cP_2})_w}="X";
(-20,-20)*+{\LS_{\cM_1}}="Y_1";
(20,-20)*+{\LS_{\cM_2}}="Y_2";
{\ar@{->} "X";"Y_1"};
{\ar@{->} "X";"Y_2"}
\endxy
$$
can be factored as
$$
\xy
(0,40)*+{(\LS_{\cP_1}\underset{\LS_{\cG}}\times \LS_{\cP_2})_w}="X";
(0,20)*+{\LS_{\sH_w}}="Q";
(0,0)*+{\LS_{\cP_1\cap w^{-1}(\cP_2)/(N(\cP_1)\cap w^{-1}(N(\cP_2)))}\overset{w}\sim \LS_{w(\cP_1)\cap \cP_2/(w(N(\cP_1))\cap N(\cP_2))}}="W";
(-30,-30)*+{\LS_{\cP'_1}}="Z_1";
(30,-30)*+{\LS_{\cP'_2}}="Z_2";
(-60,-60)*+{\LS_{\cM_1}}="Y_1";
(60,-60)*+{\LS_{\cM_2}.}="Y_2";
(0,-60)*+{\LS_{\cM'_1}\overset{w}\sim \LS_{\cM_2'}}="R";
{\ar@{->}^{\sim } "X";"Q"};
{\ar@{->}^{h} "Q";"W"};
{\ar@{->} "W";"Z_1"};
{\ar@{->} "W";"Z_2"};
{\ar@{->} "Z_1";"Y_1"};
{\ar@{->} "Z_2";"Y_2"};
{\ar@{->} "Z_1";"R"};
{\ar@{->} "Z_2";"R"};
\endxy
$$

\sssec{}

From this point, the argument proceeds as in Sects. \ref{sss:Borel case start}-\ref{sss:Borel case end}. 

\ssec{Constant term compatibility: the enhanced version}

\sssec{}

Given $\CZ\to \Ran$, consider the 1-skeleton of the following cube:

\begin{equation} \label{e:CT cube enh}
\xy
(0,0)*+{\Dmod_{\frac{1}{2}}(\Bun_G)^{-,\on{enh}_{\on{co}}}_\CZ}="X";
(40,20)*+{\Dmod_{\frac{1}{2}}(\Bun_M)\otimes \Dmod(\CZ)}="Y";
(80,0)*+{\IndCoh_\Nilp(\LS_\cG)^{-,\on{enh}_{\on{co}}}_\CZ}="X'";
(120,20)*+{\IndCoh_\Nilp(\LS_\cM)\otimes \Dmod(\CZ)}="Y'";
{\ar@{->}^{A_{G,M}^{\on{glob}}} "X";"Y"};
{\ar@{->}_{A_{\cG,\cM}^{\on{spec,glob}}} "X'";"Y'"};
{\ar@{->}^{C_G^{\on{glob}}} "X";"X'"};
{\ar@{->}^{C_M^{\on{glob}}} "Y";"Y'"};
(0,-50)*+{\KL(G)_{\crit_G,\CZ^\subseteq}\overset{\Sph_{G,\CZ}}\otimes \on{I}(G,P^-)^{\on{loc}}_{\on{co},\rho_P(\omega_X),\CZ}}="Z";
(40,-30)*+{\KL(M)_{\crit_M+\rhoch_P,(\CZ^{\subseteq})^\subseteq}}="W";
(80,-50)*+{\IndCoh^*(\Op^\mf_\cG)_{\CZ^\subseteq}
\overset{\Sph^{\on{spec}}_{\cG,\CZ}}\otimes \on{I}(\cG,\cP^-)_{\on{co},\CZ}^{\on{spec,loc}},}="Z'";
(120,-30)*+{\IndCoh^*(\Op^\mf_{\cM,\rhoch_P})_{(\CZ^\subseteq)^\subseteq}}="W'";
{\ar@{->}^{A_{G,M}^{\on{loc}}} "Z";"W"};
{\ar@{->}^{A_{\cG,\cM}^{\on{spec,loc}}} "Z'";"W'"};
{\ar@{->}^{C^{\on{loc}}_G} "Z";"Z'"};
{\ar@{->}_{C^{\on{loc}}_M} "W";"W'"};
{\ar@{->}^{B_G} "Z";"X"};
{\ar@{->}_{B_M} "W";"Y"};
{\ar@{->}_{B^{\on{spec}}_\cG} "Z'";"X'"};
{\ar@{->}_{B^{\on{spec}}_\cM} "W'";"Y'"};
\endxy
\end{equation} 
where:

\begin{itemize}


\medskip

\item $A_{G,M}^{\on{glob}}:={\on{CT}^{-,\on{enh}_{\on{co}}}_{*,\rho_P(\omega_X),\CZ}[-\delta_{(N_P^-)_{\rho_P(\omega_X)}}]}$

\item $A_{\cG,\cM}^{\on{spec,glob}}:=\on{CT}^{-,\on{spec},\on{enh}_{\on{co}}}_{\CZ}$;

\smallskip

\item $A_{G,M}^{\on{loc}}$ is the composition
\begin{multline*}
\KL(G)_{\crit_G,\CZ^\subseteq}\overset{\Sph_{G,\CZ}}\otimes \on{I}(G,P^-)^{\on{loc}}_{\on{co},\rho_P(\omega_X),\CZ}
\overset{\on{Id}\otimes \on{ins.unit}}\longrightarrow 
\KL(G)_{\crit_G,\CZ^\subseteq}\overset{\Sph_{G,\CZ^\subseteq}}\otimes \on{I}(G,P^-)^{\on{loc}}_{\on{co},\rho_P(\omega_X),\CZ^\subseteq}= \\
=(\KL(G)_{\crit_G}\overset{\Sph_G}\otimes \on{I}(G,P^-)^{\on{loc}}_{\on{co},\rho_P(\omega_X)})_{\CZ^\subseteq} 
\overset{\on{ins.unit}}\longrightarrow 
(\KL(G)_{\crit_G}\overset{\Sph_G}\otimes \on{I}(G,P^-)^{\on{loc}}_{\on{co},\rho_P(\omega_X)})_{(\CZ^\subseteq)^\subseteq}= \\
=\KL(G)_{\crit_G,(\CZ^{\subseteq})^\subseteq}^{-,\on{enh}_{\on{co}}} \overset{J^{-,\on{enh}_{\on{co}}}_{\on{KM},\rho_P(\omega_X)}}
\longrightarrow  \KL(M)_{\crit_M+\rhoch_P,(\CZ^{\subseteq})^\subseteq};
\end{multline*}

\item $A_{\cG,\cM}^{\on{spec,loc}}$ is the composition
\begin{multline*}
\IndCoh^*(\Op^\mf_\cG)_{\CZ\subseteq}\overset{\Sph^{\on{spec}}_{\cG,\CZ}}\otimes \on{I}(\cG,\cP^-)_\CZ^{\on{spec,loc}}
\overset{\on{Id}\otimes \on{ins.unit}}\longrightarrow \\
\to \IndCoh^*(\Op^\mf_\cG)_{\CZ^\subseteq}\overset{\Sph^{\on{spec}}_{\cG,\CZ^\subseteq}}
\otimes \on{I}(\cG,\cP^-)_{\CZ^\subseteq}^{\on{spec,loc}}=\\
=(\IndCoh^*(\Op^\mf_\cG)\overset{\Sph^{\on{spec}}_{\cG}}
\otimes \on{I}(\cG,\cP^-)^{\on{spec,loc}})_{\CZ^\subseteq}  \overset{\on{ins.unit}}\longrightarrow  \\
\to  (\IndCoh^*(\Op^\mf_\cG)\overset{\Sph^{\on{spec}}_{\cG}}
\otimes \on{I}(\cG,\cP^-)^{\on{spec,loc}})_{(\CZ^\subseteq)^\subseteq} = \\
=:\IndCoh^*(\Op^\mf_\cG)^{-,\on{enh}_{\on{co}}}_{(\CZ^\subseteq)^\subseteq}
\overset{J^{-,\on{spec},\on{enh}_{\on{co}}}_{\Op,\Theta}}\longrightarrow \IndCoh^*(\Op^\mf_{\cM,\rhoch_P})_{(\CZ^\subseteq)^\subseteq},
\end{multline*} 
where $J^{-,\on{spec},\on{enh}_{\on{co}}}_{\Op,\Theta}$ is the composition
\begin{multline*} 
\IndCoh^*(\Op^\mf_\cG)^{-,\on{enh}_{\on{co}}}:=\IndCoh^*(\Op^\mf_\cG)\overset{\Sph^{\on{spec}}_{\cG}} \otimes \on{I}(\cG,\cP^-)^{\on{spec,loc}} 
\overset{\Theta_{\Op^\mf_\cG}\otimes \on{Id}}\simeq \\
\simeq\IndCoh^!(\Op^\mf_\cG)^{-,\on{enh}_{\on{co}}}\overset{\Sph^{\on{spec}}_{\cG}} \otimes \on{I}(\cG,\cP^-)^{\on{spec,loc}} :=
\IndCoh^!(\Op^\mf_\cG)^{-,\on{enh}_{\on{co}}}
\overset{J^{-,\on{spec},\on{enh}_{\on{co}}}_{\Op}}\longrightarrow \\
\to \IndCoh^!(\Op^\mf_{\cM,\rhoch_P})\overset{\Theta^{-1}_{\Op^\mf_\cM}}\simeq \IndCoh^*(\Op^\mf_{\cM,\rhoch_P});
\end{multline*} 

\item $B_G$ is the composition
\begin{multline*}
\KL(G)_{\crit_G,\CZ^\subseteq}\overset{\Sph_{G,\CZ}}\otimes \on{I}(G,P^-)^{\on{loc}}_{\on{co},\rho_P(\omega_X),\CZ}
\overset{\Loc_{G,\CZ^\subseteq}\otimes \on{Id}}\longrightarrow \\
\to (\Dmod_{\frac{1}{2}}(\Bun_G)\otimes \Dmod(\CZ^{\subseteq}))\overset{\Sph_{G,\CZ}}\otimes \on{I}(G,P^-)^{\on{loc}}_{\on{co},\rho_P(\omega_X),\CZ}\to \\
\overset{\on{Id}\otimes (\on{pr}_{\on{small},\CZ})_!\otimes \on{Id}}\longrightarrow 
(\Dmod_{\frac{1}{2}}(\Bun_G)\otimes \Dmod(\CZ))\overset{\Sph_{G,\CZ}}\otimes \on{I}(G,P^-)^{\on{loc}}_{\on{co},\rho_P(\omega_X),\CZ}=
\Dmod_{\frac{1}{2}}(\Bun_G)^{-,\on{enh}_{\on{co}}}_\CZ;
\end{multline*}

\item $B_M$ is the composition
\begin{multline*}
\KL(M)_{\crit_M+\rhoch_P,(\CZ^{\subseteq})^\subseteq}\overset{\Loc_{M,(\CZ^\subseteq)^\subseteq}}\longrightarrow \\
\to \Dmod_{\frac{1}{2}}(\Bun_M)\otimes \Dmod((\CZ^\subseteq)^\subseteq)
\overset{\on{Id}\otimes (\on{pr}_{\on{small}^2,\CZ})_!}\longrightarrow \Dmod_{\frac{1}{2}}(\Bun_M)\otimes \Dmod(\CZ);
\end{multline*}
where $\on{pr}_{\on{small}^2,\CZ}$ is the projection $(\CZ^\subseteq)^\subseteq\to \CZ$;

\medskip

\item $B^{\on{spec}}_\cG$ is the composition
\begin{multline*}
\IndCoh^*(\Op^\mf_\cG)_{\CZ^\subseteq}
\overset{\Sph^{\on{spec}}_{\cG,\CZ}}\otimes \on{I}(\cG,\cP^-)_{\on{co},\CZ}^{\on{spec,loc}}
\overset{\on{Poinc}_{\cG,*,\CZ^\subseteq}\otimes \on{Id}}\longrightarrow \\
\to (\IndCoh_\Nilp(\LS_\cG)\otimes \Dmod(\CZ^{\subseteq}))\overset{\Sph^{\on{spec}}_{\cG,\CZ}}\otimes 
\on{I}(\cG,\cP^-)^{\on{spec,loc}}_{\on{co},\CZ}\to \\
\overset{\on{Id}\otimes (\on{pr}_{\on{small},\CZ})_!\otimes \on{Id}}\longrightarrow 
(\IndCoh_\Nilp(\LS_\cG)\otimes \Dmod(\CZ))\overset{\Sph^{\on{spec}}_{\cG,\CZ}}\otimes 
\on{I}(\cG,\cP^-)^{\on{spec,loc}}_{\on{co},\CZ}=
\IndCoh_\Nilp(\LS_\cG)^{-,\on{enh}_{\on{co}}}_\CZ
\end{multline*} 

\item $B^{\on{spec}}_\cM$ is the composition
\begin{multline*}
\IndCoh^*(\Op^\mf_{\cM,\rhoch_P})_{(\CZ^\subseteq)^\subseteq}
\overset{\on{Poinc}_{\cM,*,\rhoch_P,(\CZ^\subseteq)^\subseteq}}\longrightarrow \\
\to \IndCoh_\Nilp(\LS_\cM)\otimes \Dmod((\CZ^\subseteq)^\subseteq) \overset{\on{Id}\otimes (\on{pr}_{\on{small}^2,\CZ})_!}\longrightarrow 
\IndCoh_\Nilp(\LS_\cM)\otimes \Dmod(\CZ);
\end{multline*} 

\item $C^{\on{loc}}_G$ is the tensor product of $\FLE^{-,\on{enh}_{\on{co}}}_{G,\crit}$ and $((\Sat^\semiinf)^\vee)^{-1}$,
cf. \eqref{e:tensored up Sat};

\medskip

\item $C^{\on{loc}}_M:=\FLE_{M,\crit+\rhoch_P}$;

\medskip

\item $C_G^{\on{glob}}$ is the functor $\BL^{-,\on{enh}_{\on{co}}}_{G,\CZ}$, defined as the 
tensor product of $\BL_G$ and $((\Sat^\semiinf)^\vee)^{-1}$, cf. \secref{sss:L M enh}; 

\medskip

\item $C_M^{\on{glob}}:=\BL_M\otimes \on{Id}$;

\end{itemize}

\sssec{}

We claim that the five faces of this diagram commute (up to cohomological shifts and tensor products by constant lines)
in the same way as in diagram \eqref{e:CT cube}.

\medskip

The only points of difference are:

\begin{itemize} 

\item For the commutation of the bottom lid, we use \thmref{t:local Jacquet enh} instead of \thmref{t:local Jacquet};

\item For the commutation of the left face, we use \corref{c:CT and Loc crit transl enh rho} instead of
\corref{c:CT and Loc crit transl rho} ; 

\item For the commutation of the right face we use \thmref{t:CT opers enh} instead of \thmref{t:CT opers};

\item The front face commutes thanks to \thmref{t:Langlands and Loc Sat}.

\end{itemize} 

\sssec{}

By construction, all the above data are compatible with the actions of
\begin{equation} \label{e:Sat M again}
\Sph_{M,\CZ} \overset{\Sat_M}\simeq \Sph^{\on{spec}}_{\cM,\CZ}
\end{equation} 

\sssec{}

We will prove the following extension of \thmref{t:CT compat}: 

\begin{thm} \label{t:CT compat enh}
Assuming the validity of GLC for $M$, there exists a unique datum of commutativity for the upper lid in \eqref{e:CT cube enh}, i.e., the square 
\begin{equation} \label{e:L and CT enh}
\hskip1cm
\CD
\Dmod_{\frac{1}{2}}(\Bun_G)_\CZ^{-,\on{enh}_{\on{co}}} @>{\BL^{-,\on{enh}_{\on{co}}}_{G,\CZ}}>> \IndCoh_\Nilp(\LS_\cG)^{-,\on{enh}_{\on{co}}}_\CZ \\
@V{\on{CT}^{-,\on{enh}_{\on{co}}}_{*,\rho_P(\omega_X),\CZ}[-\delta_{(N^-_P)_{\rho_P(\omega_X)}}]}VV 
@VV{\on{CT}^{-,\on{spec,enh}_{\on{co}}}_\CZ\otimes \on{Id}}V \\
\Dmod_{\frac{1}{2}}(\Bun_M)\otimes \Dmod(\CZ) @>{\BL_M\otimes \on{Id}}>> \IndCoh_\Nilp(\LS_\cM)\otimes \Dmod(\CZ), 
\endCD
\end{equation}
that makes \eqref{e:CT cube enh} commute (with the specified cohomological shifts and twists by constant lines along the edges), 
in a way compatible with the actions of \eqref{e:Sat M again}.
\end{thm} 

\medskip

The rest of this subsection will be devoted to the proof of \thmref{t:CT compat enh}. 

\sssec{}

Consider the diagram 

\medskip

\begin{equation} \label{e:top lid precomposed enh}
\xy
(0,0)*+{\Dmod_{\frac{1}{2}}(\Bun_G)^{-,\on{enh}_{\on{co}}}_\CZ}="X";
(0,30)*+{\Dmod_{\frac{1}{2}}(\Bun_M)\otimes \Dmod(\CZ)}="Y";
(70,0)*+{ \IndCoh_\Nilp(\LS_\cG)^{-,\on{enh}_{\on{co}}}_\CZ}="X'";
(70,30)*+{\IndCoh_\Nilp(\LS_\cM)\otimes \Dmod(\CZ)}="Y'";
(-50,-22)*+{\KL(G)_{\crit_G,\CZ^\subseteq}\overset{\Sph_{G,\CZ}}\otimes \on{I}(G,P^-)^{\on{loc}}_{\on{co},\rho_P(\omega_X),\CZ},}="W";
{\ar@{->}_{\on{CT}^{-,\on{enh}_{\on{co}}}_{*,\rho_P(\omega_X),\CZ}[-\delta_{(N^-_P)_{\rho_P(\omega_X)}}]} "X";"Y"};
{\ar@{->}_{\on{CT}^{-,\on{spec,enh}_{\on{co}}}_\CZ} "X'";"Y'"};
{\ar@{->}^{\BL^{-,\on{enh}_{\on{co}}}_{G,\CZ}} "X";"X'"};
{\ar@{->}^{\BL_M\otimes \on{Id}} "Y";"Y'"};
{\ar@{->}_{B_G} "W";"X"}
{\ar@{->} "W";"Y"}
{\ar@{->} "W";"X'"}
\endxy
\end{equation} 
where:

\begin{itemize}

\item The arrow marked as $B_G$ is as in \eqref{e:CT cube enh};

\item The two triangles with slanted arrows commute.

\end{itemize}

The commutation of the fives faces in \eqref{e:CT cube enh} implies that the outer diagram in \eqref{e:top lid precomposed enh}
commutes. The statement of \thmref{t:CT compat enh} is that this datum comes from 
a uniquely defined commutativity datum of the inner square in \eqref{e:top lid precomposed}.

\medskip

We consider \eqref{e:top lid precomposed enh} as a diagram of categories acted on by 
$$\Sph_M\overset{\Sat_M}\simeq \Sph_\cM^{\on{spec}}.$$

\sssec{}

By duality, the required property of \eqref{e:top lid precomposed enh} is equivalent to a parallel property of the following diagram
of categories acted on by
$$\Sph_G\overset{\Sat_G}\simeq \Sph_\cG^{\on{spec}}:$$
$$
\xy
(0,0)*+{\Dmod_{\frac{1}{2}}(\Bun_G)\otimes \Dmod(\CZ)}="X";
(0,30)*+{\Dmod_{\frac{1}{2}}(\Bun_M)^{-,\on{enh}}_\CZ}="Y";
(70,0)*+{\IndCoh_\Nilp(\LS_\cG)\otimes \Dmod(\CZ)}="X'";
(70,30)*+{\IndCoh_\Nilp(\LS_\cM)^{-,\on{enh}}_\CZ,}="Y'";
(-50,-22)*+{\KL(G)_{\crit,\CZ^\subseteq}}="W";
{\ar@{->}_{\on{CT}^{-,\on{enh},\CZ}_{*,\rho_P(\omega_X)}[-\delta_{(N^-_P)_{\rho_P(\omega_X)}}]} "X";"Y"};
{\ar@{->}_{\on{CT}^{-,\on{spec,enh}}_\CZ} "X'";"Y'"};
{\ar@{->}^{\BL_G\otimes \on{Id}} "X";"X'"};
{\ar@{->}^{\BL^{-,\on{enh}}_{M,\CZ}} "Y";"Y'"};
{\ar@{->}^{\Loc'_{G,\CZ^\subseteq}} "W";"X"}
{\ar@{->} "W";"Y"}
{\ar@{->} "W";"X'"}
\endxy
$$
where:

\begin{itemize}

\item $\Loc'_{G,\CZ^\subseteq}:=(\on{pr}_{\on{small},\CZ})_!\circ \Loc_{G,\CZ^\subseteq}$;

\smallskip

\item The functor $\on{CT}^{-,\on{spec,enh}}_\CZ$ is obtained from $\on{CT}^{-,\on{spec,enh}_{\on{co}}}_\CZ$
by duality.

\end{itemize}

\sssec{}

In other words, we obtain a diagram
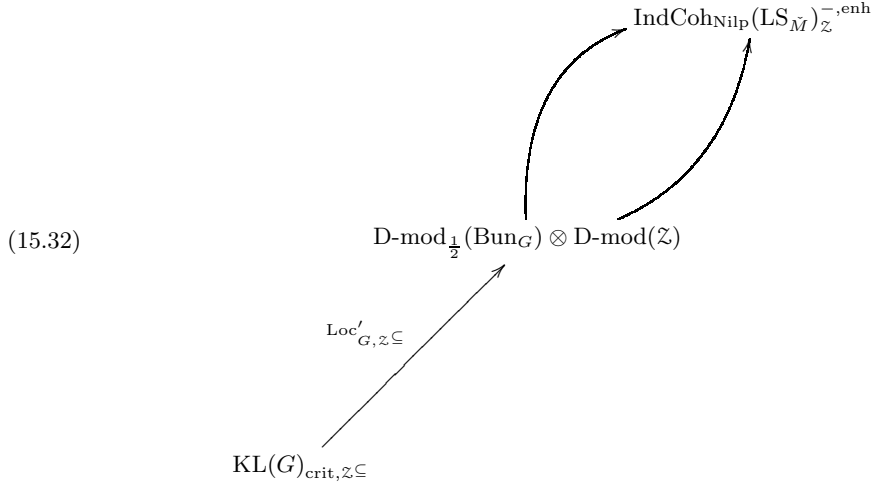
\begin{equation} \label{e:curved diag2}
\xy
(0,0)*+{\Dmod_{\frac{1}{2}}(\Bun_G)\otimes \Dmod(\CZ)}="X";
(30,30)*+{\IndCoh_\Nilp(\LS_\cM)^{-,\on{enh}}_\CZ}="Y'";
(-30,-30)*+{\KL(G)_{\crit,\CZ^\subseteq}}="Z";
{\ar@{->}@/^3pc/ "X";"Y'"};
{\ar@{->}@/_2pc/ "X";"Y'"};
{\ar@{->}^{\Loc'_{G,\CZ^\subseteq}}  "Z";"X"};
\endxy
\end{equation} 
in which the two curved arrows  
become isomorphic after precomposing with the functor $\Loc'_{G,\CZ^\subseteq}$. We wish to show that this isomorphism
comes from a uniquely defined isomorphism between the curved
arrows themselves.

\sssec{}

Recall that the forgetful functor
$$\IndCoh_\Nilp(\LS_\cM)^{-,\on{enh}}_\CZ\to \IndCoh_\Nilp(\LS_\cM)\otimes \Dmod(\CZ)$$
is monadic, and that the corresponding monad is given by the action of the associative algebra object
$$\wt\Omega^{\on{spec}}_\CZ\in \Sph_{\cM,\CZ}^{\on{spec}},$$
see \secref{sss:enh spec glob monad}. 

\medskip

The endofunctor $\sF$ of $\IndCoh_\Nilp(\LS_\cM)\otimes \Dmod(\CZ)$
underlying the above monad is naturally filtered\footnote{In fact, this filtration comes from a grading.} 
by the poset $\Lambda^{\on{pos}}_{G,P}$
(see \secref{sss:Omega tilde spec augm}, where the notation is introduced). 
Write 
$$\sF \simeq \underset{\lambda\in \Lambda^{\on{pos}}_{G,P}}{\on{colim}}\, \sF_\lambda.$$

\sssec{}

Let $\sF_{\ul\lambda}$ be a composite of a finite collection of the functors $\sF_\lambda$. Let
$$\sF_{\ul\lambda}\mod$$ denote the category of 
$$\{x\in \IndCoh_\Nilp(\LS_\cM)\otimes \Dmod(\CZ),\alpha:\sF_{\ul\lambda}(x)\to x\}.$$

\medskip

The category $\IndCoh_\Nilp(\LS_\cM)^{-,\on{enh}}_\CZ$ can be identified with a limit of categories of the form
$\sF_{\ul\lambda}\mod$. 

\sssec{}

We will show that for each $\ul\lambda$ \emph{separately}, there exists a \emph{unique} isomorphism between the 
composition of the two curved arrows in \eqref{e:curved diag2} with the projection
$$\IndCoh_\Nilp(\LS_\cM)^{-,\on{enh}}_\CZ\to \sF_{\ul\lambda}\mod,$$
so that after pre-composing with $\Loc'_{G,\CZ^\subseteq}$ we obtain the already existing isomorphism.

\medskip

The uniqueness assertion implies that these isomorphisms give rise to a uniquely defined isomorphism for
the diagram \eqref{e:curved diag2} itself. 

\sssec{}

For a sufficiently large element $\mu\in \Lambda^+_M$, let $\IndCoh_\Nilp(\LS_\cM)^{(<\mu)}$
be corresponding the quotient of $\IndCoh_\Nilp(\LS_\cM)$ (see \secref{sss:trunc LS M}), so that 
$$\IndCoh_\Nilp(\LS_\cM)\to \underset{\mu}{\on{lim}}\, \IndCoh_\Nilp(\LS_\cM)^{(<\mu)}$$
is an equivalence, see \secref{sss:IndCoh as lim}.

\medskip

The functors $\sF_\lambda$ (and, hence, their compositions $\sF_{\ul\lambda}$) have the following property: 

\medskip

For every $\mu$, the composition
\begin{multline*}
\IndCoh_\Nilp(\LS_\cM)\otimes \Dmod(\CZ) \overset{\sF_{\ul\lambda}}\to  \IndCoh_\Nilp(\LS_\cM) \otimes \Dmod(\CZ)\to \\
\to \IndCoh_\Nilp(\LS_\cM)^{(<\mu)}\otimes \Dmod(\CZ)
\end{multline*}
factors as
\begin{multline*}
\IndCoh_\Nilp(\LS_\cM) \otimes \Dmod(\CZ)\twoheadrightarrow  \IndCoh_\Nilp(\LS_\cM)^{(<\mu')}\otimes \Dmod(\CZ)
\overset{\sF_{\ul\lambda}^{\mu',\mu}}\longrightarrow  \\
\to \IndCoh_\Nilp(\LS_\cM)^{(<\mu)}\otimes \Dmod(\CZ)
\end{multline*}
for some sufficiently large $\mu'$.  

\sssec{}

For $\ul\lambda$ and $\mu\leq \mu'$ as above, let 
$$\sF_{\ul\lambda}^{\mu',\mu}\mod$$
denote the category of 
$$\{x_{\mu'}\in \IndCoh_\Nilp(\LS_\cM)^{(<\mu')}\otimes \Dmod(\CZ),\alpha:\sF^{\mu',\mu}_{\ul\lambda}(x_{\mu'})\to x_{\mu}\},$$
where $x_{\mu}$ denotes the image of $x_{\mu'}$ along the projection
$$\IndCoh_\Nilp(\LS_\cM)^{(<\mu')}\otimes \Dmod(\CZ)\to \IndCoh_\Nilp(\LS_\cM)^{(<\mu)}\otimes \Dmod(\CZ).$$

\medskip

We have:
$$\sF_{\ul\lambda}\mod \simeq \underset{\mu\leq \mu'}{\on{lim}}\, \sF_{\ul\lambda}^{\mu',\mu}\mod.$$

Hence, it is enough to show that for each $\mu\leq \mu'$, there exists a \emph{unique} isomorphism between the 
composition of the two curved arrows in \eqref{e:curved diag2} with the projection
$$\IndCoh_\Nilp(\LS_\cM)^{-,\on{enh}}_\CZ\to \sF^{\mu',\mu}_{\ul\lambda}\mod,$$
so that after pre-composing with $\Loc'_{G,\CZ^{\subseteq}}$ we obtain the already existing isomorphism.

\sssec{}

We now note that the for a fixed $\mu\leq \mu'$, the resulting two arrows 
$$\Dmod_{\frac{1}{2}}(\Bun_G)\rightrightarrows \sF^{\mu',\mu}_{\ul\lambda}\mod$$
both factor as
$$\Dmod_{\frac{1}{2}}(\Bun_G)\to \Dmod_{\frac{1}{2}}(\Bun_G)^{(<\nu),\ol\mu'}
\rightrightarrows \sF^{\mu',\mu}_{\ul\lambda}\mod$$
for some $\nu\in \Lambda_G^+$. 

\medskip

Indeed, since the functor
$$\sF^{\mu',\mu}_{\ul\lambda}\mod\to \IndCoh_\Nilp(\LS_\cM)^{(<\mu')}\otimes \Dmod(\CZ)$$
is conservative, this follows from the corresponding property of the two curved arrows
$$\Dmod_{\frac{1}{2}}(\Bun_G)\rightrightarrows  \IndCoh_\Nilp(\LS_\cM)^{(<\mu')}$$
in \eqref{e:curved diag1} established in the course of the proof of \thmref{t:CT compat}. 

\sssec{}

Now, the required assertion follows the fact that the functor
$$\KL(G)_{\crit,\CZ^\subseteq} \overset{\Loc'_{G,\CZ^\subseteq}}
\longrightarrow \Dmod_{\frac{1}{2}}(\Bun_G)\otimes \Dmod(\CZ)\to \Dmod_{\frac{1}{2}}(\Bun_G)^{(<\nu),\ol\mu'}\otimes \Dmod(\CZ)$$
is a Verdier quotient, which follows from \cite[Theorem 13.4.2]{GLC2}. 

\qed[\thmref{t:CT compat enh}]

\section{The left adjoint of the Langlands functor} \label{s:left adj}

In this section we start reaping the benefits from the work done until this point.

\begin{itemize}

\item We show that the functor $\BL_G$ admits a left adjoint (to be denoted $\BL_G^L$), which is
also compatible with the geometric and spectral Eisenstein series functors;

\smallskip

\item We show that, up to a cohomological shift, the functor $\BL_G^L$ identifies with the composition
\begin{multline} \label{e:LL as dual}
\IndCoh_\Nilp(\LS_\cG)\overset{\on{Serre}}\simeq 
\IndCoh_\Nilp(\LS_\cG)^\vee \overset{\BL_G^\vee}\longrightarrow \Dmod_{\frac{1}{2}}(\Bun_G)^\vee
\overset{\on{Verdier}}\simeq \\
\simeq \Dmod_{\frac{1}{2}}(\Bun_G)_{\on{co}}\overset{\Mir_{\Bun_G}}\longrightarrow 
\Dmod_{\frac{1}{2}}(\Bun_G) \overset{\tau_G}\to  \Dmod_{\frac{1}{2}}(\Bun_G),
\end{multline}
where $\Mir_{\Bun_G}$ is the \emph{Miraculous functor}, and $\tau_G$ is the Chevalley involution;

\smallskip

\item We show that the composition $\BL_G\circ \BL^L_G$, which is an endofunctor of $\IndCoh_\Nilp(\LS_\cG)$
is given by tensor product by an associative algebra object $\CA_G\in \QCoh(\LS_\cG)$. 

\end{itemize}

\ssec{The existence of the left adjoint}

\sssec{}

The goal of this subsection is to prove the following statement:

\begin{thm} \label{t:left adjoint}
The functor $\BL_G$ admits a left adjoint. Moreover, for every standard parabolic $P$, we have a
commutative diagram
\begin{equation} \label{e:Eis and Phi}
\CD
\Dmod_{\frac{1}{2}}(\Bun_M) @<{\BL_M^L}<< \IndCoh_\Nilp(\LS_\cM)  \\
@V{\Eis^-_{!,\rho_P(\omega_X)}[\delta_{N(P^-)_{\rho_P(\omega_X)}}]}VV @VV{\Eis^{-,\on{spec}}}V \\
\Dmod_{\frac{1}{2}}(\Bun_G) @<{\BL_G^L}<<  \IndCoh_\Nilp(\LS_\cG).
\endCD
\end{equation}
\end{thm} 

The rest of the subsection is devoted to the proof of this result. 

\sssec{}

As was shown in \cite[Theorem 13.3.6]{AG}, the category $ \IndCoh_\Nilp(\LS_\cG)$ is generated by 
the essential 
images of $\QCoh(\LS_\cM)\subset \IndCoh_\Nilp(\LS_\cM)$ along the functors 
$$\Eis^{-,\on{spec}}:\IndCoh_\Nilp(\LS_\cM) \to  \IndCoh_\Nilp(\LS_\cG).$$

\medskip

In \thmref{t:CT compat}, we have constructed a commutative square
$$
\CD
\Dmod_{\frac{1}{2}}(\Bun_M) @>{\BL_M}>> \IndCoh_\Nilp(\LS_\cM)  \\
@A{\on{CT}^-_{*,\rho_P(\omega_X)}[-\delta_{N(P^-)_{\rho_P(\omega_X)}}]}AA @AA{\on{CT}^{-,\on{spec}}}A \\
\Dmod_{\frac{1}{2}}(\Bun_G) @>{\BL_G}>>  \IndCoh_\Nilp(\LS_\cG),
\endCD
$$
in which the vertical arrows are the right adjoints to the ones in \eqref{e:Eis and Phi}. 

\medskip

It follows formally that in order to prove that $\BL_G^L$ exists and makes \eqref{e:Eis and Phi}
commute, it suffices to show that for every $M$, the (a priori partially defined) left adjoint $\BL_M^L$,
is actually defined on $\QCoh(\LS_\cM)\subset \IndCoh_\Nilp(\LS_\cM)$.

\sssec{}

Up to changing the notation, we can assume that $M=G$. However, then the existence of 
$\BL_G^L|_{\QCoh(\LS_\cG)}$ was built into the construction of $\BL_G$: this is the functor
$\BL^L_{G,\on{temp}}$ of \cite[Formula (1.10)]{GLC1}. 

\medskip

Namely, $\BL^L_{G,\on{temp}}$ 
is obtained by the spectral action of $\QCoh(\LS_\cG)$ on the object
$$\on{Poinc}_{G,!}^{\on{Vac,glob}}\in \Dmod_{\frac{1}{2}}(\Bun_G).$$

\qed[\thmref{t:left adjoint}]

\begin{rem}

Note that by passing to left adjoints along the horizontal arrows in \eqref{e:L and Eis}, we obtain a natural transformation
\begin{equation} \label{e:LL Eis nat trans}
\vcenter
{\xy
(0,0)*+{\Dmod_{\frac{1}{2}}(\Bun_G)}="X";
(50,0)*+{\IndCoh_\Nilp(\LS_\cG)}="Y";
(0,-30)*+{\Dmod_{\frac{1}{2}}(\Bun_M)}="Z";
(50,-30)*+{\IndCoh_\Nilp(\LS_\cM)}="W";
{\ar@{<-}_{\BL^L_M} "Z";"W"};
{\ar@{<-}^{\BL^L_G} "X";"Y"};
{\ar@{<-}_{\Eis^-_{!,\rho_P(\omega_X)}[\delta_{(N^-_P)_{\rho_P(\omega_X)}}]} "X";"Z"};
{\ar@{<-}^{\Eis^{-,\on{spec}}} "Y";"W"};
{\ar@{<=} "Z";"Y"};
\endxy}
\end{equation}

Alternatively, \eqref{e:LL Eis nat trans} can be viewed as obtained by passing to left adjoints along all the arrows
in \eqref{e:CT Nat trans}.

\medskip

Once we prove \thmref{t:L and CT}, we will have also proved that the natural transformation in \eqref{e:LL Eis nat trans}
is an isomorphism.

\medskip

However, we do not know whether this isomorphism equals one constructed in \thmref{e:Eis and Phi}.

\end{rem} 

\ssec{The left adjoint as a dual}

\sssec{}

Recall that Verdier duality allows us to identify
$$\Dmod_{\frac{1}{2}}(\Bun_G)^\vee \simeq \Dmod_{\frac{1}{2}}(\Bun_G)_{\on{co}}.$$

Consider the functor
$$\Dmod_{\frac{1}{2}}(\Bun_G)_{\on{co}}\overset{\BL_G^\vee}\longleftarrow  \IndCoh_\Nilp(\LS_\cG(X)),$$
dual to $\BL_G^\vee$, where we identify $ \IndCoh_\Nilp(\LS_\cG(X))$ with its own dual via Serre duality.

\medskip

Define the functor
$$\Phi_G:  \IndCoh_\Nilp(\LS_\cG(X))\longrightarrow \Dmod(\Bun_G)$$
to be 
$$\Phi_G:= \tau_G\circ \Mir_{\Bun_G}\circ \BL_G^\vee\otimes \fl_{z_\cg}[4\delta_G-2\delta_{N_{\rho(\omega_X)}}],$$
where 
$$\fl_{z_\cg}:=\det(T^*_1(\LS_{Z^0_\cG}))\simeq\det(z_\cg)^{\otimes (2-2g)}.$$

\sssec{}

We are going to prove:

\begin{thm} \label{t:left adjoint as dual}
There exists a canonical isomorphism 
$$\Phi_G\simeq \BL_G^L.$$
\end{thm}

The rest of this subsection is devoted to the proof of \thmref{t:left adjoint as dual}. 

\sssec{}

Both functors appearing in \thmref{t:left adjoint as dual} are (automatically) compatible with 
the Hecke actions via
$$\Sph_G\overset{\Sat_G}\simeq \Sph_\cG^{\on{spec}},$$
see \cite[Sect. 1.8]{GLC2}.

\medskip

Hence, they both give rise to functors
$$\Dmod_{\frac{1}{2}}(\Bun_G) _{\on{temp}}\overset{\BL^L_{G,\on{temp}}}{\underset{\Phi_{G,\on{temp}}}\leftleftarrows}
\QCoh(\LS_\cG)$$
that make both diagrams
$$
\begin{tikzcd}[row sep = large, column sep = scriptsize]
\Dmod_{\frac{1}{2}}(\Bun_G)
&&
\IndCoh_\Nilp(\LS_\cG)
\arrow[ll,shift left,shorten <= 1ex, shorten >= 1ex,"\Phi_G"]
\arrow[ll,shift right,shorten <= 1ex, shorten >= 1ex,swap,"\BL_G^L"]
\\
\Dmod_{\frac{1}{2}}(\Bun_G)_{\on{temp}}
\arrow[u,"\bu"]
&&
\QCoh(\LS_\cG)
\arrow[ll,shift left,shorten <= 1ex, shorten >= 1ex,"\Phi_{G,\on{temp}}"]
\arrow[ll,shift right,swap,shorten <= 1ex, shorten >= 1ex,"\BL_{G,\on{temp}}^L"]
\arrow[u,swap,"\Xi_{0,\Nilp}"]
\end{tikzcd}
$$
and
\begin{equation} \label{e:Phi project}
\begin{tikzcd}[row sep = large, column sep = scriptsize]
\Dmod_{\frac{1}{2}}(\Bun_G)
\arrow[d,"\bu^R"]
&&
\IndCoh_\Nilp(\LS_\cG)
\arrow[ll,shift left,shorten <= 1ex, shorten >= 1ex,"\Phi_G"]
\arrow[ll,shift right,swap,shorten <= 1ex, shorten >= 1ex,"\BL_G^L"]
\arrow[d,swap,"(\Xi_{0,\Nilp})^R"]
\\
\Dmod_{\frac{1}{2}}(\Bun_G)_{\on{temp}}
&&
\QCoh(\LS_\cG)
\arrow[ll,shift left,shorten <= 1ex, shorten >= 1ex, "\Phi_{G,\on{temp}}"]
\arrow[ll,shift right,swap,shorten <= 1ex, shorten >= 1ex,"\BL_{G,\on{temp}}^L"]
\end{tikzcd}
\end{equation}
commute, where $(\bu,\bu^R)$ are as in \cite[Sect. 18.4.1]{GLC2}.


\sssec{}

We will first show that the functors $\BL^L_{G,\on{temp}}$ and $\Phi_{G,\on{temp}}$
are (canonically) isomorphic. 

\medskip

Since the functor
$$\Loc_\cG^{\on{spec}}:\Rep(\cG)_\Ran\to \QCoh(\LS_\cG)$$ is a Verdier quotient, it suffices to establish 
an equivalence 
\begin{equation} \label{e:LL and Phi Loc}
\BL^L_{G,\on{temp}}\circ \Loc_\cG^{\on{spec}} \simeq \Phi_{G,\on{temp}} \circ \Loc_\cG^{\on{spec}}.
\end{equation} 

By construction, the functor $\BL^L_{G,\on{temp}}$ makes the diagram 
$$
\CD
\Whit^!(G)_\Ran @<{(\on{CS}_G)^{-1}}<< \Rep(\cG)_{\Ran} \\
@V{\Poinc_{G,!}[-2\delta_{N_{\rho(\omega_X)}}]}VV @VV{\Loc_\cG^{\on{spec}}}V \\
\Dmod_{\frac{1}{2}}(\Bun_G) _{\on{temp}} @<{\BL^L_{G,\on{temp}}}<< \QCoh(\LS_\cG)
\endCD
$$
commute, see \cite[Sect. 1.4.2]{GLC1}. 

\medskip 

Hence, it suffices to check that the diagram
$$
\CD
\Whit^!(G)_\Ran @<{(\on{CS}_G)^{-1}}<< \Rep(\cG)_{\Ran} \\
@V{\Poinc_{G,!}[-2\delta_{N_{\rho(\omega_X)}}]}VV @VV{\Loc_\cG^{\on{spec}}}V \\
\Dmod_{\frac{1}{2}}(\Bun_G) _{\on{temp}} 
@<{\tau_G\circ \Mir_{\Bun_G}\circ \BL_G^\vee|_{\QCoh(\LS_\cG)}\otimes \fl_{z_\cg}[4\delta_G-2\delta_{N_{\rho(\omega_X)}}]}<< \QCoh(\LS_\cG)
\endCD
$$
commutes as well. 


\sssec{}

Recall (see \cite[Lemma 1.4.12]{GLC2}) that 
$$\tau_G\circ (\on{CS}_G)^{-1}\simeq \Theta_{\Whit(G)}\circ \FLE_{\cG,\infty}.$$ Thus, we need to establish the commutativity of 
\begin{equation} \label{e:Phi com diag temp}
\CD
\Whit^!(G)_\Ran @<{\Theta_{\Whit(G)}\circ \FLE_{\cG,\infty}}<< \Rep(\cG)_{\Ran} \\
@V{\Poinc_{G,!}}VV @VV{\Loc_\cG^{\on{spec}}}V \\
\Dmod_{\frac{1}{2}}(\Bun_G) _{\on{temp}} @<{\Mir_{\Bun_G}\circ \BL_G^\vee|_{\QCoh(\LS_\cG)}\otimes \fl_{z_\cg}[4\delta_G]}<< \QCoh(\LS_\cG).
\endCD
\end{equation} 

\sssec{}

Note that $\LS_\cG$ is quasi-smooth and \emph{almost}\footnote{The symplectic structure depends on the choice of a
non-degenerate symmetric bilinear form on $z_\fg$.}
canonically derived-symplectic, so that
\begin{equation} \label{e:omega LS}
\omega_{\LS_\cG}\simeq \CO_{\LS_\cG}\otimes \fl_{z_\cg}[\dim(\LS_\cG)].
\end{equation} 

Indeed, consider the short exact sequence
$$1\to Z^0_\cG\to \cG\to \cG_1\to 1,$$
so that 
$$\LS_\cG/\LS_{Z^0_\cG}\simeq \LS_{\cG_1}.$$

Note that $T^*(\LS_{Z^0_\cG})$ is a constant bundle with fiber $T^*_1(\LS_{Z^0_\cG})$, while 
$\LS_{\cG_1}$ is \emph{canonically} derived-symplectic: for $\sigma\in \LS_{\cG_1}$
we have
$$T_\sigma(\LS_{\cG_1})\simeq \on{C}^\cdot(X,(\cg_1)_\sigma)[1],$$
so
$$(T_\sigma(\LS_{\cG_1}))^*\simeq (\on{C}^\cdot(X,(\cg_1)_\sigma)[1])^*
\overset{\on{Verdier}}\simeq \on{C}^\cdot(X,(\cg_1)^*_\sigma)[1]\overset{\on{Killing\,Form}}\simeq 
\on{C}^\cdot(X,(\cg_1)_\sigma)[1]\simeq T_\sigma(\LS_{\cG_1}).$$

\sssec{}

Recall that the functors
$$\Loc_\cG^{\on{spec}}: \Rep(\cG)_{\Ran}\leftrightarrows \QCoh(\LS_\cG):\Gamma^{\on{spec}}_\cG$$
are mutually dual, when we identify $\QCoh(\LS_\cG)$ with its own dual via the usual duality 
(i.e., usual dualization on perfect complexes).

\medskip

Hence, when we use the identification
$$\QCoh(\LS_\cG)^\vee\simeq \QCoh(\LS_\cG),$$
induced by the Serre duality on $\IndCoh(\LS_\cG)$, the dual of the functor $\Gamma^{\on{spec}}_\cG$
becomes identified with 
$$\Loc_\cG^{\on{spec}}\otimes \fl_{z_\cg}[\dim(\LS_\cG)],$$
where $\dim(\LS_\cG)=2\dim(\Bun_G)=2\delta_G$, and $\fl_{z_\cg}$ appears via \eqref{e:omega LS}. 

\sssec{}

Hence, by taking the duals in the commutative diagram
$$
\CD
\Whit^!(G)_\Ran @>{\on{CS}_G=(\FLE_{\cG,\infty})^\vee}>>  \Rep(\cG)_{\Ran} \\
@A{\on{coeff}_G[2\delta_{N_{\rho(\omega_X)}}]}AA @AA{\Gamma^{\on{spec}}_\cG}A \\
\Dmod_{\frac{1}{2}}(\Bun_G) _{\on{temp}} @>{\BL_{G,\on{temp}}}>> \QCoh(\LS_\cG),
\endCD
$$
we obtain a commutative diagram
\begin{equation} \label{e:dual of L}
\CD
\Whit_*(\Gr_G)_\Ran @<{\FLE_{\cG,\infty}}<<  \Rep(\cG)_{\Ran} \\
@V{\Poinc_{G,*}[2\delta_{N_{\rho(\omega_X)}}]}VV @VV{\Loc_\cG^{\on{spec}}\otimes \fl_{z_\cg}[2\delta_G]}V \\
\Dmod_{\frac{1}{2}}(\Bun_G) _{\on{co,temp}} @<{\BL_G^\vee|_{\QCoh(\LS_\cG)}}<< \QCoh(\LS_\cG). 
\endCD
\end{equation} 

\sssec{}

We now apply \thmref{t:Kevin}. Concatenating the diagrams 
\eqref{e:dual of L} and \eqref{e:Mir and Poinc} horizontally, we obtain the desired commutative diagram \eqref{e:Phi com diag temp}.  

\medskip

Thus, we have established the isomorphism 
\begin{equation} \label{e:left adjoint as dual temp}
\BL^L_{G,\on{temp}}\simeq \Phi_{G,\on{temp}}.
\end{equation}

\sssec{}

We will now deduce
$$\BL^L_{G,\on{temp}}\simeq \Phi_{G,\on{temp}} \,\Rightarrow\, \BL_G^L\simeq \Phi_G.$$

The functor $\BL_G^L$ being a left adjoint, preserves compactness. We will prove: 

\begin{lem} \label{l:Phi pres comp}
The functor $\Phi_G$ preserves compactness.
\end{lem}

\sssec{}

Recall now that according to \cite[Proposition 5.2.3]{GLC1}, the functor $\bu^R$ is fully faithful on compacts. 

\medskip

Hence, the assertion of \thmref{t:left adjoint as dual} follows from \lemref{l:Phi pres comp} and \eqref{e:left adjoint as dual temp}
via the commutative diagram \eqref{e:Phi project}. 

\ssec{Proof of \lemref{l:Phi pres comp}}

\sssec{}

It is enough to show that $\Phi_G$ sends objects of the form $\Eis^{-,\on{spec}}(\CF)$, for $\CF\in \QCoh(\LS_\cM)^c$
to compacts. 

%

\sssec{}

Recall that the functors 
$$\Eis^{-,\on{spec}}:\IndCoh_\Nilp(\LS_\cM)\leftrightarrows \IndCoh_\Nilp(\LS_\cG):\on{CT}^{-,\on{spec}}$$
are mutually dual, up to tensoring by a line bundle on $\LS_\cM$, see \eqref{e:CT and Eis dual}. 

\medskip 

Hence 
$$\BL_G^\vee\circ \Eis^{-,\on{spec}}$$
is isomorphic to 
$$(\on{CT}^{-,\on{spec}}\circ \BL_G)^\vee,$$
up to tensoring by a line bundle.

\medskip

Hence, combining with \thmref{t:CT compat}, we obtain that $\BL_G^\vee\circ \Eis^{-,\on{spec}}$ is isomorphic to 
$$(\BL_M\circ \on{CT}^-_*)^\vee\simeq \Eis^-_{\on{co},*}\circ \BL_M^\vee,$$
again up to tensoring by a line bundle.

\sssec{}

Now, by \eqref{e:funct eq},
$$\tau_G\circ \Mir_{\Bun_G}\circ \Eis^-_{\on{co},*} \simeq  \Eis^-_!\circ \Mir_{\Bun_M}\circ \tau_M.$$

Combining, we obtain that
$$\Phi_G\circ \Eis^{-,\on{spec}} \simeq \Eis^-_!\circ \Phi_M,$$
again up to tensoring by a line bundle.

\medskip

Now, the assertion follows from the fact that the functors $\Phi_M$ send compacts in $\QCoh(\LS_\cM)$ to compacts in
$\Dmod_{\frac{1}{2}}(\Bun_M)$, which follows from 
the isomorphism $\BL^L_{M,\on{temp}}\simeq \Phi_{M,\on{temp}}$. 

\qed[\lemref{l:Phi pres comp}]

\ssec{The composition \texorpdfstring{$\BL_G\circ \BL_G^L$}{comp}}

Recall that the geometric Langlands conjecture says that the functor $\BL_G$ is an equivalence. We can now reformulate 
this as saying that the unit of the adjunction
$$\on{Id}\to \BL_G\circ \BL_G^L$$
is an isomorphism as endofunctors of $\IndCoh_\Nilp(\LS_\cG)$, combined with the fact that 
the functor $\BL_G$ is conservative.  

\medskip

In this subsection we commence the study of the composition $\BL_G\circ \BL_G^L$. 

\sssec{}

Note that $\BL_G\circ \BL_G^L$, viewed as an endofunctor of $ \IndCoh_\Nilp(\LS_\cG)$, is $\QCoh(\LS_\cG)$-linear. Hence, it
it is a priori given by an object in 
$$\CA_G\in  \IndCoh_\Nilp(\LS_\cG)\underset{\QCoh(\LS_\cG)}\otimes  \IndCoh_\Nilp(\LS_\cG).$$

The goal of this subsection is to prove the following assertion:

\begin{thm} \label{t:AG}
The object $\CA_G$ belongs to
\begin{multline*} 
\QCoh(\LS_\cG)\simeq 
 \QCoh(\LS_\cG)\underset{\QCoh(\LS_\cG)}\otimes \QCoh(\LS_\cG)\subset \\
\subset  \IndCoh_\Nilp(\LS_\cG)\underset{\QCoh(\LS_\cG)}\otimes  \IndCoh_\Nilp(\LS_\cG).
\end{multline*}
 \end{thm}
 
In other words, this theorem implies that we have an isomorphism
$$\BL_G\circ \BL_G^L\simeq \CA_G\underset{\CO_{\LS_\cG}}\otimes -, \quad \CA_G\in  \QCoh(\LS_\cG).$$

\begin{rem} \label{r:overall strategy}

Once \thmref{t:AG} is proved, and given the fact that the functor $\BL_G$ is conservative, 
we will have interpreted GLC (i.e., \cite[Conjecture 1.6.7]{GLC1}) 
as the statement that the unit of the adjunction
\begin{equation} \label{e:unit AG}
\CO_{\LS_\cG}\to \CA_G
\end{equation}
is an isomorphism in $\QCoh(\LS_\cG)$.

\medskip

In the next section, we will show that the map \eqref{e:unit AG} becomes an isomorphism when restricted
to the \emph{reducible} locus of $\LS_\cG$. This will reduce GLC to the study of the restriction
of $\CA_G$ to the \emph{irreducible} locus $\LS^{\on{irred}}_\cG$.

\medskip

In the subsequent paper of this series (\cite{GLC4}), we will show that 
$$\CA_{G,\on{irred}}:=\CA_G|_{\LS^{\on{irred}}_\cG}$$
is a \emph{classical} vector bundle, equipped with a flat connection. 

\medskip

In the last paper of the series (\cite{GLC5}), we will deduce from this that \eqref{e:unit AG} is an isomorphism also over
$\LS^{\on{irred}}_\cG$, thereby concluding the proof of GLC.

\end{rem} 

\sssec{}

We now begin the proof of \thmref{t:AG}.

\medskip

Let $j:\LS^{\on{irred}}_\cG\subset \LS_\cG$ denote the embedding of the locus of irreducible
local systems. We denote by $i:\LS^{\on{red}}_\cG\to \LS_\cG$ be the embedding of its 
complement (with any scheme-theoretic structure). 

\medskip

Let 
\begin{equation} \label{e:red subcat}
\IndCoh_\Nilp(\LS_\cG)_{\on{red}}\subset  \IndCoh_\Nilp(\LS_\cG)
\end{equation}
be the full subcategory of objects set-theoretically supported on $\LS^{\on{red}}_\cG$.

\medskip

In other words,
$$\IndCoh_\Nilp(\LS_\cG)_{\on{red}}=\on{ker}(j^*: \IndCoh_\Nilp(\LS_\cG)\to \IndCoh_\Nilp(\LS^{\on{irred}}_\cG)).$$

Denote by $\wh{i}_!$ the tautological embedding,
$$\IndCoh_\Nilp(\LS_\cG)_{\on{red}}\hookrightarrow \IndCoh_\Nilp(\LS_\cG)$$
and by $\wh{i}^!$ its right adjoint. 
$$\wh{i}_!: \IndCoh_\Nilp(\LS_\cG)_{\on{red}}\rightleftarrows  \IndCoh_\Nilp(\LS_\cG):\wh{i}^!.$$

\sssec{}

In order to prove that $\CA_G$ belongs to $\QCoh(\LS_\cG)$, it suffices to show that 
$$(j^*\otimes j^*)(\CA_G)\in \QCoh(\LS^{\on{irred}}_\cG)$$
and 
\begin{equation} \label{e:AG on red}
(\wh{i}^!\otimes \wh{i}^!)(\CA_G)\in \QCoh(\LS_\cG)_{\on{red}}.
\end{equation}

\medskip

The former is automatic, since over $\LS^{\on{irred}}_\cG$, the nilpotent cone vanishes, so the embedding  
$$\QCoh(\LS^{\on{irred}}_\cG) \subseteq \IndCoh_\Nilp(\LS^{\on{irred}}_\cG)$$
is an equality.

\begin{rem}

Alternatively, it is easy to reduce the assertion of \thmref{t:AG} to the case when $G$ is semi-simple. In this case
$\LS^{\on{irred}}_\cG$ is smooth, so the inclusion 
$$\QCoh(\LS^{\on{irred}}_\cG) \subseteq \IndCoh(\LS^{\on{irred}}_\cG)$$
is also an equality. 

\end{rem}

\sssec{}

We now proceed to proving \eqref{e:AG on red}. 

\medskip

Since both functors $\BL_G$ and $\BL_G^L$ are compatible with the inclusions and projections
$$\IndCoh_\Nilp(\LS_\cG)\leftrightarrows \QCoh(\LS_G) \text{ and } 
\Dmod_{\frac{1}{2}}(\Bun_G) \leftrightarrows \Dmod_{\frac{1}{2}}(\Bun_G)_{\on{temp}},$$
so is the composition $\BL_G\circ \BL_G^L$. 

\medskip

Denote the induced endofunctor of $\QCoh(\LS_G)$ by $\CA_{G,\on{temp}}$. In other words,
$$\BL_{G,\on{temp}}\circ \BL^L_{G,\on{temp}}=:\CA_{G,\on{temp}}\underset{\CO_{\LS_\cG}}\otimes -, \quad \CA_{G,\on{temp}}\in  \QCoh(\LS_\cG).$$

\medskip

In order to prove \eqref{e:AG on red}, it suffices to show that 
$$(\BL_G\circ \BL_G^L)|_{ \IndCoh_\Nilp(\LS_\cG)_{\on{red}}}\simeq \CA_{G,\on{temp}}\underset{\CO_{\LS_\cG}}\otimes -$$
as endofunctors of $ \IndCoh_\Nilp(\LS_\cG)_{\on{red}}$. 

\medskip

To that end, it suffices to establish a functorial isomorphism
\begin{equation} \label{e:ten prod AG}
\BL_G\circ \BL_G^L(\CM) \simeq \CA_{G,\on{temp}}\underset{\CO_{\LS_\cG}}\otimes \CM, \quad \CM\in  \IndCoh_\Nilp(\LS_\cG)_{\on{red}}^c.
\end{equation}

\sssec{}

For $\CM$ as above, let 
$$\CM_{\on{temp}} \text{ and } (\BL_G\circ \BL_G^L(\CM))_{\on{temp}}$$
denote the projections of $\CM$ and $\BL_G\circ \BL_G^L(\CM)$, respectively, along
\begin{equation} \label{e:proj to QCoh red}
 \IndCoh_\Nilp(\LS_\cG)_{\on{red}}\to \QCoh(\LS_\cG)_{\on{red}},
\end{equation}
where
\begin{equation} \label{e:qcoh red}
\QCoh(\LS_\cG)_{\on{red}}:=\on{ker}\left(j^*:\QCoh(\LS_\cG)\to \QCoh(\LS^{\on{irred}}_\cG)\right).
\end{equation}

\medskip

By construction
\begin{equation} \label{e:temp ten prod}
(\BL_G\circ \BL_G^L(\CM))_{\on{temp}} \simeq \CA_{G,\on{temp}}\underset{\CO_{\LS_\cG}}\otimes \CM_{\on{temp}}.
\end{equation}

Thus, we wish to show that \eqref{e:temp ten prod} implies \eqref{e:ten prod AG}. 

\sssec{}

We claim:

\begin{prop} \label{p:endo red pres comp}
The functor $(\BL_G\circ \BL^L_G)|_{ \IndCoh_\Nilp(\LS_\cG)_{\on{red}}}$ 
preserves compactness. 
\end{prop}

Let us assume this proposition temporarily and finish the proof of
\thmref{t:AG}. 

\sssec{}

We claim that \propref{p:endo red pres comp} implies the following:

\begin{cor}  \label{c:endo red pres comp}
The restriction of $\CA_{G,\on{temp}}$ to the formal completion $(\LS_\cG)^\wedge_{\LS^{\on{red}}_\cG}$
is perfect as an object of $\QCoh((\LS_\cG)^\wedge_{\LS^{\on{red}}_\cG})$.
\end{cor} 

\begin{proof}[Proof of \corref{c:endo red pres comp}]

The assertion of \corref{c:endo red pres comp} can be reformulated as saying that
$$\CA_{G,\on{temp}}|_{(\LS_\cG)^\wedge_{\LS^{\on{red}}_\cG})}$$
is dualizable as an object of $\QCoh((\LS_\cG)^\wedge_{\LS^{\on{red}}_\cG})$, and equivalently, as saying that 
the functor
$$\CA_{G,\on{temp}}\underset{\CO_{\LS_\cG}}\otimes -: \QCoh((\LS_\cG)^\wedge_{\LS^{\on{red}}_\cG})\to
\QCoh((\LS_\cG)^\wedge_{\LS^{\on{red}}_\cG})$$
admits a continuous right adjoint. 

\medskip

Note that restriction along
$$(\LS_\cG)^\wedge_{\LS^{\on{red}}_\cG}\hookrightarrow \LS_\cG$$
defines an equivalence
$$\QCoh(\LS_\cG)_{\on{red}}\to \QCoh((\LS_\cG)^\wedge_{\LS^{\on{red}}_\cG}),$$
where 
$$\QCoh(\LS_\cG)_{\on{red}}\subset \QCoh(\LS_\cG)$$ is as in \eqref{e:qcoh red}. 

\medskip

Hence, it suffices to show that the functor
\begin{equation} \label{e:tensor AG}
\CA_{G,\on{temp}}\underset{\CO_{\LS_\cG}}\otimes -: \QCoh(\LS_\cG)_{\on{red}}\to \QCoh(\LS_\cG)_{\on{red}}
\end{equation}
admits a continuous right adjoint.

\medskip

It is easy to see, however, that \eqref{e:tensor AG} admits a continuous right adjoint if and only if
its composition with 
$$\Xi_{\LS_\cG}|_{\QCoh(\LS_\cG)_{\on{red}}}:\QCoh(\LS_\cG)_{\on{red}}\to \IndCoh_\Nilp(\LS_\cG)_{\on{red}}$$
does. 

\medskip

However, the latter follows from the commutative diagram
$$
\CD
\IndCoh_\Nilp(\LS_\cG)_{\on{red}} @>{\BL_G\circ \BL_G^L}>>  \IndCoh_\Nilp(\LS_\cG)_{\on{red}} \\
@A{\Xi_{\LS_\cG}|_{\QCoh(\LS_\cG)_{\on{red}}}}AA @AA{\Xi_{\LS_\cG}|_{\QCoh(\LS_\cG)_{\on{red}}}}A  \\
\QCoh(\LS_\cG)_{\on{red}} @>{\CA_{G,\on{temp}}\underset{\CO_{\LS_\cG}}\otimes -}>> \QCoh(\LS_\cG)_{\on{red}}.
\endCD
$$

Indeed, the top horizontal arrow admits a continuous right adjoint by \propref{p:endo red pres comp}, since the 
category $\IndCoh_\Nilp(\LS_\cG)_{\on{red}}$ is compactly generated.  And the right adjoint of the left vertical arrow
is the (continuous) functor $\Psi_{\LS_\cG}|_{\IndCoh_\Nilp(\LS_\cG)_{\on{red}}}$. 

\end{proof} 

\sssec{}

Let $\CM$ be as in \eqref{e:ten prod AG}. By \propref{p:endo red pres comp}, 
$$\BL_G\circ \BL_G^L(\CM)\in \Coh(\LS_\cG)_{\on{red}}:=\Coh(\LS_\cG)\cap \IndCoh(\LS_\cG)_{\on{red}},$$
and by \corref{c:endo red pres comp}, the object 
$\CA_{G,\on{temp}}\underset{\CO_{\LS_\cG}}\otimes \CM$ also belongs to this subcategory.

\medskip

The projections of these objects along \eqref{e:proj to QCoh red} are identified by \eqref{e:temp ten prod}. 
The isomorphism \eqref{e:ten prod AG} now follows, since the restriction of \eqref{e:proj to QCoh red} 
to
$$\Coh(\LS_\cG)_{\on{red}}\subset \IndCoh(\LS_\cG)_{\on{red}}$$
is fully faithful. 

\qed[\thmref{t:AG}]

\ssec{Proof of \propref{p:endo red pres comp}} 

\sssec{} \label{sss:Eis part}

Let 
$$\Dmod_{\frac{1}{2}}(\Bun_G)_{\Eis}\subset \Dmod_{\frac{1}{2}}(\Bun_G)$$
denote the full subcategory generated by the essential images of the functors
$$\Eis_!:\Dmod_{\frac{1}{2}}(\Bun_M)\to \Dmod_{\frac{1}{2}}(\Bun_G)$$
for Levi quotients $M$ of \emph{proper} parabolics $P\subset G$.

\sssec{}

We claim that the functors
$$\BL_G:\Dmod_{\frac{1}{2}}(\Bun_G)\leftrightarrows  \IndCoh_\Nilp(\LS_\cG):\BL_G^L$$
send the full subcategories
$$\Dmod_{\frac{1}{2}}(\Bun_G)_{\Eis}\subset \Dmod_{\frac{1}{2}}(\Bun_G) \text{ and }  \IndCoh_\Nilp(\LS_\cG)_{\on{red}}\subset  \IndCoh_\Nilp(\LS_\cG)$$
to one another, and the resulting functors
$$\BL_G:\Dmod_{\frac{1}{2}}(\Bun_G)_{\Eis}\leftrightarrows  \IndCoh_\Nilp(\LS_\cG)_{\on{red}}:\BL_G^L$$
preserve compactness. 

\medskip

This would imply that $(\BL_G\circ \BL^L_G)|_{ \IndCoh_\Nilp(\LS_\cG)_{\on{red}}}$ also preserves compactness. 

\sssec{}

For the functor $\BL_G$, this follows from its compatibility with the Eisenstein procedure, expressed by \thmref{t:L and Eis},
and the assumption that GLC holds for $M$. 

\sssec{}

To prove the assertion for $\BL_G^L$, we note that the subcategory \eqref{e:red subcat} is generated by the essential 
images of the functors 
$$\Eis^{\on{spec}}:\IndCoh_\Nilp(\LS_\cM)\to  \IndCoh_\Nilp(\LS_\cG)$$
for proper parabolics (indeed, the collection of their right adjoints (i.e., the functors $\on{CT}^{\on{spec}}$) is
conservative on $ \IndCoh_\Nilp(\LS_\cG)_{\on{red}}$).

\medskip

Now, the assertion concerning $\BL_G^L$ follows from the commutative diagram \eqref{e:Eis and Phi}.

\qed[\propref{p:endo red pres comp}]

\section{The Langlands functor is an equivalence on the Eisenstein part(s)} \label{s:Eis equiv}

In this section we state and prove the main result of this paper, \thmref{t:main}, which says that
the Langlands functor $\BL_G$ induces an equivalence between the Eisenstein-generated categories 
on the two sides. 

\ssec{Statement of the main result}

\sssec{}

In \secref{sss:Eis part}, we have considered the pair of (mutually adjoint) functors
\begin{equation} \label{e:Lang Eis}
\BL_G:\Dmod_{\frac{1}{2}}(\Bun_G)_{\Eis}\leftrightarrows  \IndCoh_\Nilp(\LS_\cG)_{\on{red}}:\BL_G^L.
\end{equation}

We will prove: 

\begin{mainthm} \label{t:main}
Let us assume that GLC is valid for Levi subgroups of all proper parabolics of $G$.
Then the adjoint functors in \eqref{e:Lang Eis} are mutually inverse equivalences.
\end{mainthm} 

\sssec{Proof of \thmref{t:L and CT}} \label{sss:L and CT}

Let us assume \thmref{t:main} for a moment and deduce from it \thmref{t:L and CT}. 

\medskip

By passing to the right adjoints along the vertical arrows in
$$
\CD
\Dmod_{\frac{1}{2}}(\Bun_M) @>{\BL_M}>{\sim}> \IndCoh_\Nilp(\LS_\cM) \\
@V{\Eis^-_{!,\rho_P(\omega_X)}[\delta_{N(P^-)_{\rho_P(\omega_X)}}]}VV @VV{\Eis^{-,\on{spec}}}V \\
\Dmod_{\frac{1}{2}}(\Bun_G)_{\Eis} @>{\sim}>{\BL_G}>  \IndCoh_\Nilp(\LS_\cG)_{\on{red}},
\endCD
$$
we obtain that the natural transformation in question is an equivalence, once both
circuits in the diagram \eqref{e:CT Nat trans} are restricted to
$$\Dmod_{\frac{1}{2}}(\Bun_G)_{\Eis}\subset \Dmod_{\frac{1}{2}}(\Bun_G).$$

\medskip

Hence, in order to prove the corollary, it suffices to show that both circuits vanish when restricted
to 
$$\Dmod_{\frac{1}{2}}(\Bun_G)_{\on{cusp}}:=(\Dmod_{\frac{1}{2}}(\Bun_G)_{\Eis})^\perp.$$

The vanishing is tautological for the counterclockwise circuit. For the vanishing of the clockwise circuit, it suffices to show 
that there exists \emph{some} isomorphism
$$\on{CT}^{-,\on{spec}}\circ \BL_G \simeq \BL_M\circ \on{CT}^-_{*,\rho_P(\omega_X)}[-\delta_{N(P^-)_{\rho_P(\omega_X)}}].$$

However, the latter is given by \thmref{t:CT compat}.

\qed[\thmref{t:L and CT}]

\begin{rem}

As was already mentioned, we do not know whether the isomorphism
$$\on{CT}^{-,\on{spec}}\circ \BL_G \simeq \BL_M\circ \on{CT}^-_{*,\rho_P(\omega_X)}[-\delta_{N(P^-)_{\rho_P(\omega_X)}}],$$
constructed in \thmref{t:CT compat} above equals the one given by 
\thmref{t:L and CT}.

\end{rem}

\begin{rem} \label{r:red to irred}

Note that \thmref{t:main} implies that GLC is equivalent to the statement that the functors
$(\BL_G^L,\BL_G)$ define mutually inverse equivalences
\begin{equation} \label{e:L on irred}
\BL_{G,\on{irred}}:
\Dmod_{\frac{1}{2}}(\Bun_G)_{\on{cusp}} \leftrightarrows \IndCoh_\Nilp(\LS^{\on{irred}}_\cG)\simeq \QCoh(\LS^{\on{irred}}_\cG):\BL^L_{G,\on{irred}}.
\end{equation} 

As was mentioned in Remark \ref{r:overall strategy}, this is equivalent to showing that the map
\begin{equation} \label{e:unit AG irred}
\CO_{\LS^{\on{irred}}_\cG}\to \CA_{G,\on{irred}},
\end{equation}
induced by \eqref{e:unit AG}, is an isomorphism. 

\end{rem}

\ssec{Constant terms and the left adjoint of the Langlands functor}

\sssec{}

As a preparation to the proof of \thmref{t:main}, we will need the following:

\begin{prop} \label{p:* forth}
The following diagram of functors commutes
$$
\CD
\Dmod_{\frac{1}{2}}(\Bun_M) @<{\BL^L_M}<< \IndCoh_\Nilp(\LS_\cM) \\
@A{\on{CT}^-_{*,\rho_P(\omega_X)}[-\delta_{(N^-_P)_{\rho_P(\omega_X)}}]}AA @AA{\on{CT}^{-,\on{spec}}}A \\
\Dmod_{\frac{1}{2}}(\Bun_G) @<<{\BL^L_G}< \IndCoh_\Nilp(\LS_\cG). 
\endCD
$$
\end{prop} 

The rest of this subsection is devoted to the proof of this proposition.

\sssec{}

Applying \thmref{t:left adjoint as dual} for $G$ and $M$, 
we need to show that the following diagram commutes
$$
\CD
\Dmod_{\frac{1}{2}}(\Bun_M) @<{\Mir_{M,\tau}}<< \Dmod_{\frac{1}{2}}(\Bun_M)_{\on{co}}  @<{\BL^\vee_M}<< \IndCoh_\Nilp(\LS_\cM) \\
@A{\on{CT}^-_{*,\rho_P(\omega_X)}}AA & & @AA{\on{CT}^{-,\on{spec}}}A \\
\Dmod_{\frac{1}{2}}(\Bun_G) @<<{\Mir_{G,\tau}}< \Dmod_{\frac{1}{2}}(\Bun_G)_{\on{co}}  @<<{\BL^\vee_G}< \IndCoh_\Nilp(\LS_\cG),
\endCD
$$
up to a cohomological shift by 
$$4\delta_G-2\delta_{N_{\rho(\omega_X)}}-4\delta_M+2\delta_{N(M)_{\rho_M(\omega_X)}}-\delta_{(N^-_P)_{\rho_P(\omega_X)}}=
3(\delta_G-\delta_M)-\delta_{(N_P)_{\rho_P(\omega)_X}}.$$

\medskip

\noindent NB: the (constant) lines $\fl_{z_\cg}$ and $\fl_{z_\cm}$ cancel out, since the Killing form of $\cg$
is non-degenerate on $z_\cm/z_\cg$. 

\medskip

Passing to dual functors, we obtain that it suffices to show that the following diagram commutes
$$
\CD
\Dmod_{\frac{1}{2}}(\Bun_M)_{\on{co}} @>{\Mir_{M,\tau}}>> \Dmod_{\frac{1}{2}}(\Bun_M)  @>{\BL_M}>> \IndCoh_\Nilp(\LS_\cM) \\
@V{(\on{CT}^-_{*,\rho_P(\omega_X)})^\vee}VV & & @VV{(\on{CT}^{-,\on{spec}})^\vee}V \\
\Dmod_{\frac{1}{2}}(\Bun_G)_{\on{co}} @>>{\Mir_{G,\tau}}> \Dmod_{\frac{1}{2}}(\Bun_G)  @>>{\BL_G}> \IndCoh_\Nilp(\LS_\cG),
\endCD
$$
up to the same cohomological shift. 

\sssec{}

Recall that according to \eqref{e:CT and Eis dual} and \lemref{l:rel det LS}, with respect to the Serre duality identifications
$$\IndCoh_\Nilp(\LS_\cG)^\vee\simeq \IndCoh_\Nilp(\LS_\cG) \text{ and } \IndCoh_\Nilp(\LS_\cM)^\vee\simeq \IndCoh_\Nilp(\LS_\cM),$$
we have a canonical identification
$$(\on{CT}^{-,\on{spec}})^\vee \simeq \Eis^{-,\on{spec}}\circ (\CL_{2\rho_P(\omega_X)}\otimes (-)),$$
up to a cohomological shift by 
$$\dim(\cn^-_P)\cdot (2g-2).$$

\medskip

Also, by construction, we have
$$(\on{CT}^-_{*,\rho_P(\omega_X)})^\vee\simeq \Eis^-_{\on{co},*,\rho_P(\omega_X)}.$$

\medskip

Hence, we need to show that the following diagram commutes
\begin{equation} \label{e:* forth}
\CD
\Dmod_{\frac{1}{2}}(\Bun_M)_{\on{co}} @>{\Mir_{M,\tau}}>> \Dmod_{\frac{1}{2}}(\Bun_M)  @>{\BL_M}>> \IndCoh_\Nilp(\LS_\cM) \\
@V{\Eis^-_{\on{co},*,\rho_P(\omega_X)}}VV & & @VV{\Eis^{-,\on{spec}}\circ (\CL_{2\rho_P(\omega_X)}\otimes (-))}V \\
\Dmod_{\frac{1}{2}}(\Bun_G)_{\on{co}} @>>{\Mir_{G,\tau}}> \Dmod_{\frac{1}{2}}(\Bun_G)  @>>{\BL_G}> \IndCoh_\Nilp(\LS_\cG),
\endCD
\end{equation} 
up to a cohomological shift by 
$$3(\delta_G-\delta_M)-\delta_{(N_P)_{\rho_P(\omega)_X}}-\dim(\cn^-_P)\cdot (2g-2)=
2(\delta_G-\delta_M)-\delta_{(N_P)_{\rho_P(\omega)_X}}.$$

\sssec{}

Note that by \eqref{e:funct eq}, we have the following diagram 
\begin{equation} \label{e:* forth 1}
\CD
\Dmod_{\frac{1}{2}}(\Bun_M)_{\on{co}} @>{\Mir_{M,\tau}}>> \Dmod_{\frac{1}{2}}(\Bun_M) \\
@V{\Eis^-_{\on{co},*,\rho_P(\omega_X)}}VV @VV{\Eis^-_{!,-\rho_P(\omega_X)}}V \\
\Dmod_{\frac{1}{2}}(\Bun_G)_{\on{co}} @>>{\Mir_{G,\tau}}> \Dmod_{\frac{1}{2}}(\Bun_G), 
\endCD
\end{equation}
which commutes up to a cohomological shift by $2\delta_{N^-_P}$. 

\medskip

Hence, we obtain that it suffices to show that the following diagram commutes
\begin{equation} \label{e:* forth 2}
\CD
\Dmod_{\frac{1}{2}}(\Bun_M) @>{\BL_M}>> \IndCoh_\Nilp(\LS_\cM)  \\
@V{\Eis^-_{!,-\rho_P(\omega_X)}}VV @VV{\Eis^{-,\on{spec}}\circ (\CL_{2\rho_P(\omega_X)}\otimes (-))}V \\
\Dmod_{\frac{1}{2}}(\Bun_G)  @>>{\BL_G}> \IndCoh_\Nilp(\LS_\cG),
\endCD
\end{equation} 
up to a cohomological shift by 
$$2(\delta_G-\delta_M)-\delta_{(N_P)_{\rho_P(\omega)_X}}-2\delta_{N^-_P}=
(\delta_G-\delta_M)-\delta_{(N_P)_{\rho_P(\omega)_X}}=\delta_{(N^-_P)_{\rho_P(\omega)_X}}.$$

\sssec{} 

Note that the following diagram commutes
\begin{equation} \label{e:L and transl}
\CD
\Dmod_{\frac{1}{2}}(\Bun_M) @>{\BL_M}>> \IndCoh_\Nilp(\LS_\cM)  \\
@V{(\on{transl}_{-2\rho_P(\omega_X)})_*}VV @VV{\CL_{2\rho_P(\omega_X)}\otimes (-)}V \\
\Dmod_{\frac{1}{2}}(\Bun_M) @>{\BL_M}>> \IndCoh_\Nilp(\LS_\cM), 
\endCD
\end{equation} 
see \cite[Sect. 1.5.3]{GLC1}. 

\medskip

We also note that 
$$\Eis^-_{!,-\rho_P(\omega_X)}\simeq \Eis^-_{!,\rho_P(\omega_X)}\circ (\on{transl}_{-2\rho_P(\omega_X)})_*.$$

\medskip

Hence, the commutation of \eqref{e:!* forth 2} is equivalent to the commutation of 
\begin{equation} \label{e:* forth 3}
\CD
\Dmod_{\frac{1}{2}}(\Bun_M) @>{\BL_M}>> \IndCoh_\Nilp(\LS_\cM)  \\
@V{\Eis^-_{!,\rho_P(\omega_X)}[\delta_{(N^-_P)_{\rho_P(\omega_X)}}]}VV @VV{\Eis^{-,\on{spec}}}V \\
\Dmod_{\frac{1}{2}}(\Bun_G)  @>>{\BL_G}> \IndCoh_\Nilp(\LS_\cG).
\endCD
\end{equation} 

However, the latter was established in \thmref{t:L and Eis}. 

\qed[\propref{p:* forth}]

\sssec{}

Combining with \thmref{t:CT compat} we obtain:

\begin{cor} \label{c:* back and forth}
We have a canonical isomorphism 
$$\on{CT}^{-,\on{spec}}\simeq \on{CT}^{-,\on{spec}}\circ \BL_G\circ \BL_G^L$$
as functors 
$$\IndCoh_\Nilp(\LS_\cG)\rightrightarrows \IndCoh_\Nilp(\LS_\cM).$$
\end{cor} 

In particular, we have:

\begin{cor} \label{c:* back and forth AG}
We have an isomorphism in $\QCoh(\LS_\cM)\subset \IndCoh_\Nilp(\LS_\cM)$
$$\on{CT}^{-,\on{spec}}(\CA_G)\simeq \on{CT}^{-,\on{spec}}(\CO_{\LS_\cG}).$$
\end{cor}

\ssec{Proof of \thmref{t:main}}

\sssec{}

In the course of the proof of \propref{p:endo red pres comp}, we have seen that the essential image of
each of the functors in \eqref{e:Lang Eis} generates the target category.

\medskip

In particular, we obtain that the functor $\BL_G|_{\Dmod_{\frac{1}{2}}(\Bun_G)_{\Eis}}$ is conservative.

\medskip

Hence, in order to prove \thmref{t:main}, it suffices to show that the functor $\BL_G^L|_{\IndCoh_\Nilp(\LS_\cG)_{\on{red}}}$ is
fully faithful.

\medskip

I.e., we need to show that the unit of the adjunction
\begin{equation} \label{e:unit AG red}
\on{Id}_{ \IndCoh_\Nilp(\LS_\cG)_{\on{red}}}\to (\BL_G\circ \BL_G^L)|_{\IndCoh_\Nilp(\LS_\cG)_{\on{red}}}
\end{equation}
is an isomorphism.

\sssec{}

By \thmref{t:AG}, the endofunctor $\BL_G\circ \BL_G^L$ is given by tensor product with an object
$$\CA_G\in \QCoh(\LS_\cG).$$

\medskip

Moreover, the structure of monad on $\BL_G\circ \BL_G^L$ corresponds to a structure of associative
algebra on $\CA_G$ as an object of the (symmetric) monoidal category $\QCoh(\LS_\cG)$.

\medskip

Under this identification, the unit of the adjunction
$$\on{Id}_{ \IndCoh_\Nilp(\LS_\cG)}\to \BL_G\circ \BL_G^L$$
corresponds to the map of associative algebras
\begin{equation} \label{e:unit AG again}
\CO_{\LS_\cG}\to \CA_G.
\end{equation}

\sssec{}

In order to show that \eqref{e:unit AG red} is an isomorphism, it suffices to show that 
the *-restriction of the map \eqref{e:unit AG again} to the formal completion
$$(\LS_\cG)^\wedge_{\LS^{\on{red}}_\cG} \subset \LS_\cG$$
is an isomorphism.

\sssec{}

Note that the composition
$$\QCoh(\LS_\cG)_{\on{red}}\hookrightarrow \QCoh(\LS_\cG)\overset{\wh{i}^*}\to \QCoh((\LS_\cG)^\wedge_{\LS^{\on{red}}_\cG})$$
is an equivalence. 

\medskip

Note also that the collection of functors
$$(\sfp^{-,\on{spec}})^*: \QCoh(\LS_\cG)\to \QCoh(\LS_{\cP^-}),$$
for proper standard parabolics of $G$, is conservative on $\QCoh(\LS_\cG)_{\on{red}}$. 

\medskip

Hence, it suffices to show that for each proper parabolic as above, the induced map
\begin{equation} \label{e:AP}
\CO_{\LS_{\cP^-}}\to (\sfp^{-,\on{spec}})^*(\CA_G)
\end{equation} 
is an isomorphism.

\sssec{}

We will deduce this from the following statement:

\begin{prop}   \label{p:AP 1}
For every proper parabolic, the object  
$$(\sfp^{-,\on{spec}})^*(\CA_G)\in \QCoh(\LS_{\cP^-})$$
is a line bundle. 
\end{prop} 

Indeed, the fact that \propref{p:AP 1} implies that \eqref{e:AP} is an isomorphism
follows from the next observation:

\begin{lem} \label{l:line bundle is enough}
Let $\CY$ be an algebraic stack, and let $\CA_\CY$ be a unital associative algebra object in $\QCoh(\CY)$.
Assume that $\CA_\CY$ is a line bundle. Then the unit map $\CO_\CY\to \CA_\CY$ is an isomorphism.
\end{lem} 

\begin{proof}[Proof 1]

It suffices to show that the unit map is an isomorphism at the level of fibers at field-valued points.
This reduces the assertion to the case when $\CY$ is the spectrum of a field, in which case
it is obvious.

\end{proof} 

\begin{proof}[Proof 2]

Let $\CA$ be a unital associative algebra object in a symmetric monoidal stable $\infty$-category $\bA$.
Suppose that, as a plain object, $\CA$ is invertible. We claim that in this case, the unit map
$$u:\one_\bA\to \CA$$
is an isomorphism.

\medskip

Indeed, tensoring
$$m:\CA\otimes \CA\to \CA$$
on the left with $\CA^{\otimes -1}$, we obtain a map 
$$v:\CA\to \one_\bA.$$

Let us show that the maps $u$ and $v$ are each other's inverses. The fact that $v\circ u=\on{id}$ is obvious. 
Hence, we obtain that $\one_\bA$ is a retract of $\CA$.
Write
$$\CA=\one_\bA\oplus \CB.$$

We wish to show that $\CB=0$. It suffices to show that the tautological map $\alpha:\CB\to \CA$ is
null-homotopic. Consider the composition: 
$$\CB \overset{u \otimes \on{id}}{\longrightarrow} \CA \otimes \CB
\overset{\on{id} \otimes \alpha}{\longrightarrow} 
\CA \otimes \CA \overset{m}{\longrightarrow} \CA$$

\medskip

On the one hand, it equals $\alpha$. On the other hand, since $m = \on{id} \otimes v$ and $v \circ \alpha = 0$, we obtained that 
the above composition vanishes. 

\end{proof} 

\sssec{}

Consider now the map
\begin{equation} \label{e:LS M to P again}
\LS_\cM\to \LS_{\cP^-}
\end{equation} 
induced by the canonical map
\begin{equation} \label{e:B M to P again}
\on{pt}/\cM\to \on{pt}/\cP^-.
\end{equation}

\medskip

We will deduce \propref{p:AP 1} from the following assertion:

\begin{prop}   \label{p:AP 2}
The *-pullback of $(\sfp^{-,\on{spec}})^*(\CA_G)$ along \eqref{e:LS M to P again} is a line
bundle on $\LS_\cM$. 
\end{prop}

We will prove \propref{p:AP 2} in \secref{ss:L and comp CT}. The implication 
\propref{p:AP 2} $\Rightarrow$ \propref{p:AP 1} will be carried out in \secref{ss:line bundles}
below. 

\sssec{}

Once Propositions \ref{p:AP 1} and \ref{p:AP 2} are proved, the proof of \thmref{t:main} would be complete. 

\ssec{Proof of \propref{p:AP 1}} \label{ss:line bundles} 

%
%
%
%
%
%
%
%
%
%
%
%

\sssec{} \label{sss:contracting}

Choose a coweight
$$\check\theta:\BG_m\to Z_\cM,$$
which is \emph{anti-dominant and regular}, i.e., its action on any $\alpha\in \cn^-_P$ 
is a strictly positive integer.

\medskip

The adjoint action of $\BG_m$ on $\cP^-$ gives rise to an action of $\BG_m$ on $\LS_{\cP^-}$. This action has the following properties
(cf. \cite[Sect. 11.2]{DG2}): 

\begin{itemize}

\item(a) It respects the projection $\sfq^{-,\on{spec}}:\LS_{\cP^-}\to \LS_\cM$ (because the adjoint action of $\check\theta$ on
$\cM$ is trivial);

\item(b) It is \emph{contracting}, i.e., the action of $\BG_m$ extends to an action of the monoid $\BA^1$, with the fixed-point locus 
being the section of $\sfq^{-,\on{spec}}$ given by \eqref{e:LS M to P again};

\item(c) It is isomorphic to the trivial action (because the $\BG_m$-action on $\cP^-$ is inner).

\end{itemize}

\medskip


\sssec{} \label{sss:discr}

Note that the discrepancy between the trivialization of the $\BG_m$-action on $\LS_{\cP^-}$ and the data of $\BG_m$-equivariance 
on the projection $\sfq^{-,\on{spec}}$ is given by the action of $\BG_m$ by 1-automorphisms on $\LS_\cM$.

\medskip

In particular, for $\CF\in \QCoh(\LS_{\cP^-})$, the $\BG_m$-action on $(\sfq^{-,\on{spec}})_*(\CF)$ coming from: 

\begin{itemize}

\item the $\BG_m$-equivariant structure on $\CF$ given by point (c) above, and 

\item the structure of $\BG_m$-invariance on $\sfq^{-,\on{spec}}$ from point (a)

\end{itemize}

\noindent equals the $\BG_m$-action on $(\sfq^{-,\on{spec}})_*(\CF)$ as an object of $\QCoh(\LS_\cM)$,
coming from the action of $\BG_m$ by 1-automorphisms on $\LS_\cM$.

\sssec{}

Consider
$$(\sfq^{-,\on{spec}})_*(\CO_{\LS_{\cP^-}})$$
as a graded connective commutative algebra in $\QCoh(\LS_\cM)$. 

\medskip

Note also that $(\sfq^{-,\on{spec}})_*(\CO_{\LS_{\cP^-}})$
is \emph{as good as connective} (see \cite[Sect. 12.7.6]{GLC2} or \secref{sss:Spec quot} for what this means). In
particular, the functor
$$(\sfq^{-,\on{spec}})_*:\QCoh(\LS_{\cP^-})\to \QCoh(\LS_\cM)$$
upgrades to an equivalence
\begin{equation} \label{e:LS P as good}
\QCoh(\LS_{\cP^-})\to (\sfq^{-,\on{spec}})_*(\CO_{\LS_{\cP^-}})\mod(\QCoh(\LS_\cM)).
\end{equation} 

\medskip

By property (b) in \secref{sss:contracting}, the grading on $(\sfq^{-,\on{spec}})_*(\CO_{\LS_{\cP^-}})$
is non-negative, with the degree $0$ part isomorphic to $\CO_{\LS_\cM}$, which is a retract of 
$(\sfq^{-,\on{spec}})_*(\CO_{\LS_{\cP^-}})$ via \eqref{e:LS M to P again}.
 
\medskip
 
In particular, the functor of pullback along \eqref{e:LS M to P again} corresponds under the equivalence \eqref{e:LS P as good}
to the functor
$$\CO_{\LS_\cM}\underset{(\sfq^{-,\on{spec}})_*(\CO_{\LS_{\cP^-}})}\otimes (-).$$

%

\sssec{}

We now claim:

\begin{lem} \label{l:grading pos}
The grading on $(\sfq^{-,\on{spec}})_*\circ (\sfp^{-,\on{spec}})^*(\CA_G)$ is non-negative.
\end{lem}

Assuming the lemma for a moment, the fact that $(\sfp^{-,\on{spec}})^*(\CA_G)$ is a line bundle
follows from \propref{p:AP 2} using the next version of the graded Nakayama lemma:

\begin{lem} \label{l:gr Nakayama}
Let $A$ be a non-negatively graded commutative algebra. Denote by $A_0$ its degree $0$ component.
Let $M$ be a non-negatively graded $A$-module. Suppose that the $A_0$-module 
$$M_0:=A_0\underset{A}\otimes M$$
is invertible. Then $M$ is itself invertible. 
\end{lem} 

\qed[\propref{p:AP 1}]

\sssec{Proof of \lemref{l:grading pos}} 

%
%
%
%
%
%



Recall that, according to \secref{sss:discr}, this grading equals the grading on $(\sfq^{-,\on{spec}})_*\circ (\sfp^{-,\on{spec}})^*(\CA_G)$
as an object of $\QCoh(\LS_\cM)$ coming from the action of $\BG_m$ by 1-automorphisms on $\LS_\cM$.

\medskip

Note that
$$(\sfq^{-,\on{spec}})_*\circ (\sfp^{-,\on{spec}})^*(\CA_G)\simeq
((\sfq^{-,\on{spec}})_*\circ (\sfp^{-,\on{spec}})^!(\CA_G))\otimes \CL$$
for a particular (cohomologically shifted) line bundle on $\LS_\cM$. 

\medskip

We also have
$$(\sfq^{-,\on{spec}})_*\circ (\sfp^{-,\on{spec}})^*(\CO_{\LS_\cG})\simeq
((\sfq^{-,\on{spec}})_*\circ (\sfp^{-,\on{spec}})^!(\CO_{\LS_\cG}))\otimes \CL$$
for the same line bundle $\CL$.

\medskip

Now, according to \corref{c:* back and forth AG}, we have an isomorphism
$$(\sfq^{-,\on{spec}})_*\circ (\sfp^{-,\on{spec}})^!(\CA_G)\simeq (\sfq^{-,\on{spec}})_*\circ (\sfp^{-,\on{spec}})^!(\CO_{\LS_\cG}).$$

Hence, we also have
$$(\sfq^{-,\on{spec}})_*\circ (\sfp^{-,\on{spec}})^*(\CA_G)\simeq (\sfq^{-,\on{spec}})_*\circ (\sfp^{-,\on{spec}})^*(\CO_{\LS_\cG}).$$

\medskip

Hence, the required assertion follows from the fact that the 
grading on 
$$(\sfq^{-,\on{spec}})_*\circ (\sfp^{-,\on{spec}})^*(\CO_{\LS_\cG})\simeq (\sfq^{-,\on{spec}})_*(\CO_{\LS_{\cP^-}})$$
is non-negative.

\qed[\lemref{l:grading pos}]

%
%
%
%
%
%

\ssec{Langlands functor and the ``compactified" constant term functor(s)} \label{ss:L and comp CT}

The goal of this subsection is to prove \propref{p:AP 2}. We will reduce it to a statement
about the compatibility of the Langlands functor and its left adjoint with the ``compactified" constant
term functors. 

\sssec{}

Since both $\LS_\cG$ and $\LS_{\cP^-}$ are quasi-smooth, and 
$$\CA_G\in \QCoh(\LS_\cG)\subset \IndCoh(\LS_\cG),$$
when proving \propref{p:AP 2}, instead of $(\sfp^{-,\on{spec}})^*(\CA_G)$, we can consider 
$(\sfp^{-,\on{spec}})^!(\CA_G)$ (the two differ by tensoring by a line bundle). 

\medskip

Note that the resulting object of 
$$\QCoh(\LS_\cM)\subset \IndCoh(\LS_\cM)$$
is by definition
$$\on{CT}^{-,\on{spec}}_{!*}(\CA_G),$$
see \secref{sss:comp spec CT}.

\medskip

Note also that by definition,
$$\CA_G\simeq \BL_G\circ \BL_G^L(\CO_{\LS_\cG}).$$

\medskip

Hence, in order to prove \propref{p:AP 2}, it suffices to show the following: 

\begin{prop} \label{p:!* back and forth}
There is a canonical isomorphism 
$$\on{CT}^{-,\on{spec}}_{!*}\circ \BL_G\circ \BL_G^L\simeq \on{CT}^{-,\on{spec}}_{!*}$$
as functors
$$\IndCoh_\Nilp(\LS_\cG)\rightrightarrows \IndCoh_\Nilp(\LS_\cM).$$
\end{prop}

\sssec{}

Recall the equivalence 
$$\Theta_{\on{I}(G,P^-)^{\on{loc}},\tau}: \on{I}(G,P^-)^{\on{loc}}_{\on{co}}\simeq \on{I}(G,P^-)^{\on{loc}},$$
see \eqref{e:Theta IGP glob}. Let 
$$\IC^{-,\semiinf}_{\on{co}}\in \on{I}(G,P^-)^{\on{loc}}_{\on{co}}$$
be the factorization algebra object equal to the image of $\IC^{-,\semiinf}$ under $\Theta_{\on{I}(G,P^-)^{\on{loc}},\tau}$. 

\medskip

Let
$$\on{CT}^-_{!*,\rho_P(\omega_X)}:\Dmod_{\frac{1}{2}}(\Bun_G)\to \Dmod_{\frac{1}{2}}(\Bun_M)$$
be the functor as in \secref{sss:CT with fact alg kernel}
with $\CB:=\IC^-_{\on{co}}$. 

\sssec{}

\propref{p:!* back and forth} is obtained by combining the following two assertions:

\begin{prop} \label{p:!* back}
The following diagram of functors commutes
$$
\CD
\Dmod_{\frac{1}{2}}(\Bun_M) @>{\BL_M}>> \IndCoh_\Nilp(\LS_\cM) \\
@A{\on{CT}^-_{!*,\rho_P(\omega_X)}[-\delta_{(N^-_P)_{\rho_P(\omega_X)}}]}AA @AA{\on{CT}^{-,\on{spec}}_{!*}}A \\
\Dmod_{\frac{1}{2}}(\Bun_G) @>>{\BL_G}> \IndCoh_\Nilp(\LS_\cG). 
\endCD
$$
\end{prop} 

\begin{prop} \label{p:!* forth}
The following diagram of functors commutes
$$
\CD
\Dmod_{\frac{1}{2}}(\Bun_M) @<{\BL^L_M}<< \IndCoh_\Nilp(\LS_\cM) \\
@A{\on{CT}^-_{!*,\rho_P(\omega_X)}[-\delta_{(N^-_P)_{\rho_P(\omega_X)}}]}AA @AA{\on{CT}^{-,\on{spec}}_{!*}}A \\
\Dmod_{\frac{1}{2}}(\Bun_G) @<<{\BL^L_G}< \IndCoh_\Nilp(\LS_\cG). 
\endCD
$$
\end{prop} 

\sssec{Proof of \propref{p:!* back}}

By \propref{p:CT spec IC}, we rewrite $\on{CT}^{-,\on{spec}}_{!*}$ as $\on{CT}^{-,\on{spec}}_{\IC}$.

\medskip

Note now that the object $\IC^{-,\semiinf}_{\on{co}}$ can be also described as the image of
$$\IC^{-,\on{spec},\semiinf}_{\on{co}}\in \on{I}(\cG,\cP^-)^{\on{spec,loc}}_{\on{co}}$$ 
under the equivalence $((\Sat^\semiinf)^\vee)^{-1}$.

\medskip

The assertion of the proposition follows now from \thmref{t:CT compat enh}. 

\qed[\propref{p:!* back}]

\sssec{} \label{sss:proof of forth start}

The rest of this subsection is devoted to the proof of \propref{p:!* forth}. In fact, it will follow line-by-line the proof 
of \propref{p:* forth}. 

\medskip

\noindent NB: in the argument below we will not specify the exact amounts of cohomological shifts, because
they are the same as in the proof of \propref{p:* forth}.  

\medskip

Applying \thmref{t:left adjoint as dual} for $G$ and $M$, 
we need to show that the following diagram commutes
$$
\CD
\Dmod_{\frac{1}{2}}(\Bun_M) @<{\Mir_{M,\tau}}<< \Dmod_{\frac{1}{2}}(\Bun_M)_{\on{co}}  @<{\BL^\vee_M}<< \IndCoh_\Nilp(\LS_\cM) \\
@A{\on{CT}^-_{!*,\rho_P(\omega_X)}}AA & & @AA{\on{CT}^{-,\on{spec}}_{!*}}A \\
\Dmod_{\frac{1}{2}}(\Bun_G) @<<{\Mir_{G,\tau}}< \Dmod_{\frac{1}{2}}(\Bun_G)_{\on{co}}  @<<{\BL^\vee_G}< \IndCoh_\Nilp(\LS_\cG),
\endCD
$$
up to a cohomological shift. 

\medskip

Passing to dual functors, we obtain that it suffices to show that the following diagram commutes
$$
\CD
\Dmod_{\frac{1}{2}}(\Bun_M)_{\on{co}} @>{\Mir_{M,\tau}}>> \Dmod_{\frac{1}{2}}(\Bun_M)  @>{\BL_M}>> \IndCoh_\Nilp(\LS_\cM) \\
@V{(\on{CT}^-_{!*,\rho_P(\omega_X)})^\vee}VV & & @VV{(\on{CT}^{-,\on{spec}}_{!*})^\vee}V \\
\Dmod_{\frac{1}{2}}(\Bun_G)_{\on{co}} @>>{\Mir_{G,\tau}}> \Dmod_{\frac{1}{2}}(\Bun_G)  @>>{\BL_G}> \IndCoh_\Nilp(\LS_\cG),
\endCD
$$
up to the same cohomological shift. 

\sssec{}

Recall that according to \lemref{l:spec CT and Eis IC}, with respect to the Serre duality identifications
$$\IndCoh_\Nilp(\LS_\cG)^\vee\simeq \IndCoh_\Nilp(\LS_\cG) \text{ and } \IndCoh_\Nilp(\LS_\cM)^\vee\simeq \IndCoh_\Nilp(\LS_\cM),$$
we have a canonical identification
$$(\on{CT}^{-,\on{spec}}_{!*})^\vee \simeq \Eis^{-,\on{spec}}_{!*}\circ (\CL_{2\rho_P(\omega_X)}\otimes (-)),$$
up to a cohomological shift. 

\medskip

Also, by construction, we have
$$(\on{CT}^-_{!*,\rho_P(\omega_X)})^\vee\simeq \Eis^-_{\on{co},*,\rho_P(\omega_X),\IC},$$
where $\Eis^-_{\on{co},*,\rho_P(\omega_X),\IC}$ is the functor from \secref{sss:CT with fact alg kernel}
with $\CB:=\IC^-_{\on{co}}$. 

\medskip

Hence, we need to show that the following diagram commutes
\begin{equation} \label{e:!* forth}
\CD
\Dmod_{\frac{1}{2}}(\Bun_M)_{\on{co}} @>{\Mir_{M,\tau}}>> \Dmod_{\frac{1}{2}}(\Bun_M)  @>{\BL_M}>> \IndCoh_\Nilp(\LS_\cM) \\
@V{\Eis^-_{\on{co},*,\rho_P(\omega_X),\IC}}VV & & @VV{\Eis^{-,\on{spec}}_{!*}\circ (\CL_{2\rho_P(\omega_X)}\otimes (-))}V \\
\Dmod_{\frac{1}{2}}(\Bun_G)_{\on{co}} @>>{\Mir_{G,\tau}}> \Dmod_{\frac{1}{2}}(\Bun_G)  @>>{\BL_G}> \IndCoh_\Nilp(\LS_\cG),
\endCD
\end{equation} 
up to a cohomological shift. 

\sssec{}

Note that by \corref{c:twisted intertwiners}, we have the following diagram 
\begin{equation} \label{e:!* forth 1}
\CD
\Dmod_{\frac{1}{2}}(\Bun_M)_{\on{co}} @>{\Mir_{M,\tau}}>> \Dmod_{\frac{1}{2}}(\Bun_M) \\
@V{\Eis^-_{\on{co},*,\rho_P(\omega_X),\IC}}VV @VV{\Eis^-_{!,-\rho_P(\omega_X),\IC}}V \\
\Dmod_{\frac{1}{2}}(\Bun_G)_{\on{co}} @>>{\Mir_{G,\tau}}> \Dmod_{\frac{1}{2}}(\Bun_G), 
\endCD
\end{equation}
which commutes up to a cohomological shift. 

\medskip

Hence, we obtain that it suffices to show that the following diagram commutes
\begin{equation} \label{e:!* forth 2}
\CD
\Dmod_{\frac{1}{2}}(\Bun_M) @>{\BL_M}>> \IndCoh_\Nilp(\LS_\cM)  \\
@V{\Eis^-_{!,-\rho_P(\omega_X),\IC}}VV @VV{\Eis^{-,\on{spec}}_{!*}\circ (\CL_{2\rho_P(\omega_X)}\otimes (-))}V \\
\Dmod_{\frac{1}{2}}(\Bun_G)  @>>{\BL_G}> \IndCoh_\Nilp(\LS_\cG),
\endCD
\end{equation} 
up to a cohomological shift. 

\sssec{} \label{sss:proof of forth end}

Recall the commutative diagram \eqref{e:L and transl}. We also note that 
$$\Eis^-_{!,-\rho_P(\omega_X),\IC}\simeq \Eis^-_{!,\rho_P(\omega_X),\IC}\circ (\on{transl}_{-2\rho_P(\omega_X)})_*.$$

\medskip

Hence, the commutation of \eqref{e:!* forth 2} is equivalent to the commutation of 
\begin{equation} \label{e:!* forth 3}
\CD
\Dmod_{\frac{1}{2}}(\Bun_M) @>{\BL_M}>> \IndCoh_\Nilp(\LS_\cM)  \\
@V{\Eis^-_{!,\rho_P(\omega_X),\IC}[\delta_{(N^-_P)_{\rho_P(\omega_X)}}]}VV @VV{\Eis^{-,\on{spec}}_{!*}}V \\
\Dmod_{\frac{1}{2}}(\Bun_G)  @>>{\BL_G}> \IndCoh_\Nilp(\LS_\cG).
\endCD
\end{equation} 

However, the latter was established in \corref{c:L and Eis !*}. 

\qed[\propref{p:!* forth}]

%
%
%
%
%
%
%

\newpage 

\appendix

\centerline{\bf Appendices}

\bigskip

\centerline{By J.~Campbell, L.~Chen, J.~Faergeman, D.~Gaitsgory, K.~Lin, S.~Raskin and N.~Rozenblyum}

\section{The spectral semi-infinite category/ies over the Ran space} \label{s:semiinf spec}

We will construct the category $$\IndCoh^*(\LS^\reg_\cG\underset{\LS^\mer_\cG}\times \LS^\mer_{\cP^-}\underset{\LS^\mer_\cM}\times \LS^\reg_\cM)$$
by emulating the recipe of \cite[Sect. E.2]{GLC2} for the construction of 
$$\IndCoh^*(\LS^\reg_\sH\underset{\LS^\mer_\sH}\times \LS^\reg_\sH)$$
for an algebraic group $\sH$. 

\medskip

Namely, we will first define the factorization category
\begin{equation} \label{e:defn of L G/N}
\IndCoh^*(\fL_\nabla(\cG/\cN^-_P)),
\end{equation} 
equipped with an action of the factorization group-scheme $\fL^+_\nabla(\cG)\times \fL^+_\nabla(\cM)$, 
and set 
\begin{equation} \label{e:IGP spec defn}
\IndCoh^*(\LS^\reg_\cG\underset{\LS^\mer_\cG}\times \LS^\mer_{\cP^-}\underset{\LS^\mer_\cM}\times \LS^\reg_\cM):=
\IndCoh^*(\fL_\nabla(\cG/\cN^-_P))_{\fL^+_\nabla(\cG)\times \fL^+_\nabla(\cM)}
\end{equation} 

Note, however, that the definition of \eqref{e:defn of L G/N} (even over an individual $X^I$)
is not immediate, as it does not fit into the paradigm developed in \cite[Sect. A.5]{GLC2}, since $\cG/\cN^-_P$ is not 
affine. 

\ssec{Horizontal loops into the (parabolic) base affine space}

\sssec{}

Let $Y$ be a quasi-affine scheme almost of finite type. Even though it is not affine, we can apply the definition from \cite[Sect. B.4.6]{GLC2}
and consider the factorization space $\fL_\nabla(Y)$.

\medskip

We claim:

\begin{prop} \label{p:loops indsch}
The prestack $\fL_\nabla(Y)_\Ran\to \Ran$ is a relative ind-scheme over $\Ran$. 
\end{prop}

\begin{proof}

The map
\begin{equation} \label{e:hor jets into loops}
\fL^+_\nabla(Y)_\Ran\to \fL_\nabla(Y)_\Ran
\end{equation} 
induces an isomorphism of the underlying classical prestacks,\footnote{Indeed, this is true for $Y=\BA^1$, hence for any affine $Y$, and hence for
any quasi-affine $Y$ as well.} while 
$$\fL^+_\nabla(Y)_\Ran\to \Ran$$
is a relative scheme. 

\medskip

Hence, in order to show that $\fL_\nabla(Y)$ is an ind-scheme, it suffices to show that it admits
a connective deformation theory (see \cite[Chapter 2, Theorem 1.3.12]{GaRo4}). However, the latter is a straightforward
verification.

\end{proof} 

\sssec{} \label{sss:open in hor loops}

Note that since \eqref{e:hor jets into loops} is an isomorphism at the classical level, 
$\fL_\nabla(Y)_\Ran\to \Ran$ is a relative \emph{formal} scheme\footnote{In the sense of \cite[Chapter 2, Sect. 0.2.3]{GaRo4}.}
over $\Ran$.

\medskip

Additionally, we obtain that if $U\subset Y$ is an open subscheme (resp., a Zariski cover of $Y$), then 
$$\fL_\nabla(U)_\Ran\to \fL_\nabla(Y)_\Ran$$
is an open embedding (resp., Zariski cover).

\sssec{}

We apply the above discussion to $Y$ being the (parabolic) base affine space $\cG/\cN^-_P$. Thus,
we obtain a (factorization) ind-scheme $\fL_\nabla(\cG/\cN^-_P)$.

\begin{rem}

Here is another way to see that $\fL_\nabla(\cG/\cN^-_P)$ is a relative ind-scheme. Let $\ol{\cG/\cN^-_P}$
denote the affine closure of $\cG/\cN^-_P$.

\medskip

Then, by \secref{sss:open in hor loops}, $\fL_\nabla(\cG/\cN^-_P)$ is the open subscheme of $\fL_\nabla(\ol{\cG/\cN^-_P})$, whose underlying
classical scheme equals that of
$$\fL^+_\nabla(\cG/\cN^-_P)\subset \fL^+_\nabla(\ol{\cG/\cN^-_P}).$$

\end{rem}

\sssec{}

The relevance of $\fL_\nabla(\cG/\cN^-_P)$ for us is explained by the following assertion:

\begin{lem} \label{l:loops into base affine as fiber product}
There is a canonical isomorphism
$$\fL_\nabla(\cG/\cN^-_P)\simeq \on{pt}\underset{\LS^\mer_\cG}\times \LS^\mer_{\cP^-}\underset{\LS^\mer_\cM}\times \on{pt}.$$
\end{lem} 

\begin{proof}

Parallel to that of \cite[Lemma E.2.2]{GLC2}.

\end{proof}

\begin{cor} \label{c:loops into base affine as fiber product}
We have a canonical isomorphism
$$\fL^+_\nabla(\cG)\backslash \fL_\nabla(\cG/\cN^-_P)/\fL^+_\nabla(\cM)\simeq
\LS^\reg_\cG\underset{\LS^\mer_\cG}\times \LS^\mer_{\cP^-}\underset{\LS^\mer_\cM}\times \LS^\reg_\cM.$$
\end{cor} 

\ssec{Ind-coherent sheaves on horizontal loop and arc spaces} \label{ss:IndCoh on loops and arcs}

\sssec{}

Recall the construction from \cite[Sect. A.5]{GLC2} that attaches to an ind-affine ind-scheme $\CY$ the category $\IndCoh^*(\CY)$. 

\sssec{} \label{sss:IndCoh * open}

Note that if $\CY^0\overset{j}\hookrightarrow \CY$ is an open embedding, then the functor
$$j^\IndCoh_*: \IndCoh^*(\CY^0)\to  \IndCoh^*(\CY)$$
admits a left adjoint, denoted $j^{\IndCoh,*}$ (see \cite[Sect. A.7.6]{GLC2}). Moreover, it is easy to see that $j^\IndCoh_*$ is fully faithful, making $j^{\IndCoh,*}$ 
a localization.\footnote{The fact that the counit of the adjunction $\on{Id}\to j^{\IndCoh,*}\circ j_{\IndCoh,*}$ is an isomorphism follows by base change.}

\sssec{}

We now record the following observation to the effect that the construction
$$\CY\rightsquigarrow \IndCoh^*(\CY)$$ satisfies Zariski descent.

\medskip

Let $\CU\to \CY$ is a Zariski cover (i.e., $\CU$ is a finite disjoint union of open subfunctors of $\CY$,
each of which is itself an ind-affine ind-scheme, that cover it set-theoretically). 

\begin{lem} \label{l:descent for IndCoh^*}  \hfill

\smallskip

\noindent{\em(a)}
The functor
$$\IndCoh^*(\CY)\to \on{Tot}(\IndCoh^*(\CU^\bullet)),$$
given by $(\IndCoh,*)$-restriction, is an equivalence.

\smallskip

\noindent{\em(b)}
The functor
$$|\IndCoh^*(\CU^\bullet)|\to \IndCoh^*(\CY),$$
given by $(\IndCoh,*)$-pushforward, is also an equivalence.

\end{lem} 

\sssec{} \label{sss:IndCoh* indsch}

Note that \lemref{l:descent for IndCoh^*} allows us to extend the assignment 
$$\CY\rightsquigarrow \IndCoh^*(\CY)$$
to (not necessarily ind-affine) ind-schemes that admit finite open covers by ind-affine ones.

\sssec{} \label{sss:IndCoh on arcs}

Let $I$ be a finite set and consider the scheme $X^I$. Let $Y$ be a quasi-affine
scheme of finite type. By \lemref{l:descent for IndCoh^*}, we have a well-defined
category
$$\IndCoh^*(\fL^+_\nabla(Y)_{X^I}).$$

\medskip

However, we do not know whether for a general $Y$, the assignment
$$I\rightsquigarrow \IndCoh^*(\fL^+_\nabla(Y))_{X^I}$$
can be extended to a datum of factorization category. Namely, we do not know how to establish the analogs 
of \cite[Lemmas B.13.18 and B.13.20]{GLC2}.

\sssec{}

We claim, however:

\begin{lem} \label{l:IndCoh vs QCoh arcs}
Assume that $Y$ is smooth. 
The functor 
$$\Psi_{\fL^+_\nabla(Y)}:\IndCoh^*(\fL^+_\nabla(Y))\to \QCoh(\fL^+_\nabla(Y))$$
is an equivalence. 
\end{lem}

\begin{proof}
Follows from \cite[Proposition E.12.1]{GLC2}.
\end{proof}

The above lemma allows us to define $\IndCoh^*(\fL^+_\nabla(Y))$ as a (unital) factorization category for a smooth $Y$.
Namely, \cite[Lemmas B.13.18 and B.13.20]{GLC2} hold automatically for $\QCoh(-)$. 

\medskip

Note that by \lemref{l:IndCoh vs QCoh arcs}
$$\IndCoh^*(\fL^+_\nabla(Y))\simeq \QCoh(\fL^+_\nabla(Y)),$$
as factorization categories. 

\sssec{} \label{sss:IndCoh* indsch fact good}

Let $Y$ be a smooth quasi-affine scheme of finite type. We now wish to define the (factorization) category $\IndCoh^*(\fL_\nabla(Y))$. 

\medskip

Combining \lemref{l:descent for IndCoh^*} and \propref{p:loops indsch}, 
we obtain that for a fixed finite set $I$, we have a well-defined category
$$\IndCoh^*(\fL_\nabla(Y)_{X^I}).$$

\medskip

However, we do not know whether for a general $Y$, the assignment
$$I\rightsquigarrow \IndCoh^*(\fL_\nabla(Y))_{X^I}$$
can be extended to a datum of factorization category (for the same reason as in \secref{sss:IndCoh on arcs} above). 

\medskip

That said, we \emph{can} establish the required assertions in the following cases:

\medskip

\noindent{(i)} When $Y$ is an algebraic group $\sH$.

\medskip

\noindent{(ii)} When $Y$ is an open subset of $\BA^n$. 

\medskip

Indeed, case (i) was considered in \cite[Sect. E.2]{GLC2}. For case (ii), if $Y=\BA^n$, this is a particular case of (i),
and in general the required properties are inherited under open embeddings, see \secref{sss:IndCoh * open}. 

\sssec{}  \label{sss:IndCoh* indsch fact}

In particular, we obtain that if $Y$ is a quasi-affine scheme almost of finite type \emph{that admits a cover by open subsets isomorphic
to open subsets of $\BA^n$}, then we have a well-defined unital factorization category $\IndCoh^*(\fL_\nabla(Y))$. 

\ssec{Ind-coherent sheaves on \texorpdfstring{$\fL_\nabla(\cG/\cN^-_P)$}{loopsGNbar}} \label{ss:IndCoh loopsGNbar}

In this subsection we define $$\IndCoh^*(\fL_\nabla(\cG/\cN^-_P))$$ as a factorization category. 
We will give two different constructions, which will ultimately be equivalent. 

\sssec{}

For the first definition of $\IndCoh^*(\fL_\nabla(\cG/\cN^-_P))$, we start with the unital factorization category
$$\IndCoh^*(\fL_\nabla(\cG)).$$

It is equipped with an action of the monoidal unital factorization category 
$$\IndCoh^*(\fL_\nabla(\cG\times \cG))\simeq \IndCoh^*(\fL_\nabla(\cG))\otimes \IndCoh^*(\fL_\nabla(\cG)),$$ 
and the action morphism has a naturally defined lax unital structure. 

\medskip

In particular, we can restrict this action to $\IndCoh^*(\fL_\nabla(\cN^-_P))$, which we consider as mapping to the second factor. 

\sssec{}

We set
\begin{equation} \label{e:defn of L G/N 1}
'\IndCoh^*(\fL_\nabla(\cG/\cN^-_P)):=\IndCoh^*(\fL_\nabla(\cG))\underset{\IndCoh^*(\fL_\nabla(\cN^-_P))}\otimes \Vect.
\end{equation}

\medskip

Note that by construction, $'\IndCoh^*(\fL_\nabla(\cG/\cN^-_P))$ carries an action of
\begin{equation} \label{e:loops G M}
\IndCoh^*(\fL_\nabla(\cG))\otimes \IndCoh^*(\fL_\nabla(\cM)).
\end{equation} 

\sssec{}  \label{sss:2nd defn IndCoh G/NP}

For the second definition, we apply the paradigm of \secref{sss:IndCoh* indsch fact} directly to $Y=\cG/\cN^-_P$, 
where the open cover of $\cG/\cN^-_P$ is provided by translates of the big Schubert cell.

\medskip

Thus, we obtain a well-defined unital factorization category 
$$\IndCoh^*(\fL_\nabla(\cG/\cN^-_P)).$$

Again, by construction, this category carries a monoidal action of \eqref{e:loops G M}. 

\sssec{}

It follows that $\IndCoh^*(\fL_\nabla(\cG/\cN^-_P))$ as constructed in \secref{sss:2nd defn IndCoh G/NP} 
comes equipped with a (strictly unital) factorization functor 
\begin{equation} \label{e:arcs to loops G/NP}
\iota^\IndCoh_*:\IndCoh^*(\fL^+_\nabla(\cG/\cN^-_P))\to \IndCoh^*(\fL_\nabla(\cG/\cN^-_P)),
\end{equation}
given by $(\IndCoh,*)$-pushforward.

\medskip

In particular, the factorization unit in $\IndCoh^*(\fL_\nabla(\cG/\cN^-_P))$ is the image of 
$$\CO_{\fL^+_\nabla(\cG/\cN^-_P)}\in \IndCoh^*(\fL^+_\nabla(\cG/\cN^-_P))$$
under the above functor.

\sssec{}

Note that the operation of $(\IndCoh,*)$-pushforward gives rise to a functor
\begin{equation} \label{e:IndCoh'}
'\IndCoh^*(\fL_\nabla(\cG/\cN^-_P))\to \IndCoh^*(\fL_\nabla(\cG/\cN^-_P)).
\end{equation} 

\begin{prop} \label{p:IndCoh'}
The functor \eqref{e:IndCoh'} is an equivalence.
\end{prop}

\begin{proof}

Let $U$ be the cover of $\cG/\cN^-_P$ by translates of the big Schubert cell. 
Set 
$$\wt{U}:=U\underset{\cG/\cN^-_P}\times \cG.$$

By construction, 
$$\IndCoh^*(\fL_\nabla(\cG/\cN^-_P))\simeq |\IndCoh^*(\fL_\nabla(U^\bullet))|.$$

\medskip

By \lemref{l:descent for IndCoh^*}(b), we have:
$$\IndCoh^*(\fL_\nabla(\cG))\simeq  |\IndCoh^*(\fL_\nabla(\wt{U}^\bullet))|.$$

Hence, 
$$'\IndCoh^*(\fL_\nabla(\cG/\cN^-_P))\simeq |\IndCoh^*(\fL_\nabla(\wt{U}^\bullet))\underset{\IndCoh^*(\fL_\nabla(\cN^-_P))}\otimes \Vect|.$$

Therefore, it suffices to show that for every connected component $U'$ of $U^n$ (for every $n$), the functor
$$\IndCoh^*(\fL_\nabla(\wt{U}'))\underset{\IndCoh^*(\fL_\nabla(\cN^-_P))}\otimes \Vect\to 
\IndCoh^*(\fL_\nabla(U'))$$
is an equivalence. 

\medskip

Note that for every $U'$ as above, the $\cN^-_P$-torsor $\wt{U}'\to U'$ splits as a product
$$U'\times \cN^-_P.$$

Hence, it suffices to show that
$$\IndCoh^*(\fL_\nabla(U')\times \fL_\nabla(\cN^-_P)) \underset{\IndCoh^*(\fL_\nabla(\cN^-_P))}\otimes \Vect\to \IndCoh^*(\fL_\nabla(U'))$$
is an equivalence. 

\medskip

However, from the fact that $\IndCoh^*(\fL_\nabla(\cN^-_P))$ is dualizable (see \secref{sss:dual gen loops} below), we obtain that the functor
$$\IndCoh^*(\fL_\nabla(U'))\otimes \IndCoh^*(\fL_\nabla(\cN^-_P))\to \IndCoh^*(\fL_\nabla(U')\times \fL_\nabla(\cN^-_P))$$
is an equivalence,\footnote{The dualizability is used in commuting limits involved in the definition of $\IndCoh^*(-)$ with colimits and tensor products,  
cf. \cite[Chapter 3, Proposition 3.1.7.]{GaRo3}} and the assertion follows.

\end{proof}

\ssec{Compact generation and duality}

\sssec{}

Recall the notion of compact generation for a factorization category, see \cite[Sect. B.11.10]{GLC2}. 

\medskip

Let $Y$ be a smooth scheme. Consider the factorization category $\IndCoh^*(\fL^+_\nabla(Y))$. 
It follows from \lemref{l:IndCoh vs QCoh arcs} that it is compactly generated. 

\sssec{} 

We now claim:

\begin{prop} \label{p:arcs to loops G/NP}
The functor \eqref{e:arcs to loops G/NP} admits a (conservative) right adjoint as a factorization functor,
to be denoted $\iota^!$. 
\end{prop} 

\begin{proof}

It suffices to show that if $U$ is an open affine in $\BA^n$, then the functor
\begin{equation} \label{e:iota U}
\iota^\IndCoh_*:\IndCoh^*(\fL^+_\nabla(U))\to \IndCoh^*(\fL_\nabla(U))
\end{equation}
admits a conservative right adjoint as a factorization functor and that the 
the formation of these right adjoints is compatible with $(\IndCoh,*)$-pushforwards 
along further open embeddings $U'\hookrightarrow U$). 

\medskip 

We first consider the case when $U=\BA^n$. 
In the latter case, taking $\fL^+_\nabla(\BA^n)$-invariants on both sides (cf. \cite[Lemma E.2.10]{GLC2}),
we reduce the assertion to that of the functor
\begin{equation} \label{e:arcs to loops G/NP bis}
\Vect\simeq \IndCoh(\on{pt})\to \fL_\nabla(\BA^n)/\fL^+_\nabla(\BA^n).
\end{equation}

However, the functor \eqref{e:arcs to loops G/NP bis} is the $(\IndCoh,*)$-pushforward in 
the laft category, and the assertion follows.

\medskip

We now claim that if the functor \eqref{e:iota U} admits a right adjoint (as a factorization functor),
then the same will be true for any open $U'\subset U$, and that the formation of these right adjoints
is compatible with $(\IndCoh,*)$-pushforwards along open embeddings $j:U'\hookrightarrow U$). 
Indeed, this follows from \secref{sss:IndCoh * open} by passing to right adjoints in the commutative diagram
$$
\CD
\IndCoh^*(\fL^+_\nabla(U)) @>{j^{\IndCoh,*}}>> \IndCoh^*(\fL^+_\nabla(U')) \\
@V{\iota^\IndCoh_*}VV @VV{\iota^\IndCoh_*}V \\
\IndCoh^*(\fL_\nabla(U)) @>>{j^{\IndCoh,*}}> \IndCoh^*(\fL_\nabla(U')). 
 \endCD
 $$

\end{proof} 

\sssec{} \label{sss:comp gen loops} 

Combining \propref{p:arcs to loops G/NP} with \lemref{l:IndCoh vs QCoh arcs}, we obtain:

\begin{cor} \label{c:IndCoh loops comp gen}
The factorization category $\IndCoh^*(\fL_\nabla(\cG/\cN^-_P))$ is compactly generated. 
\end{cor}

\sssec{} \label{sss:dual gen loops} 

Being compactly generated, the category $\IndCoh^*(\fL_\nabla(\cG/\cN^-_P))$ is dualizable as a
factorization category (see \cite[Sect. B.11.9]{GLC2} for what this means). 

\medskip

Moreover, it is easy to see that it is dualizable as a \emph{unital} factorization category
(see \cite[C.11.5]{GLC2} for what this means).

\medskip

We will now describe explicitly the dual of $\IndCoh^*(\fL_\nabla(\cG/\cN^-_P))$.

\sssec{}

Recall (see \cite[Sect. A.4]{GLC2}) that to any prestack $\CY$ we can attach a category $\IndCoh^!(\CY)$.

\medskip

We take 
$$\CY:=\fL_\nabla(Y),$$
viewed as a factorization space, where $Y$ is a scheme of finite type. We would like to extend the assignment
$$I\in \on{fSet}\, \rightsquigarrow \, \IndCoh^!(\fL_\nabla(Y))_I$$
to a datum of factorization category. 

\medskip

In order to do so, we need to prove assertions parallel to \cite[Lemmas B.13.11 and B.13.14]{GLC2}. We do not
know how to do this in general. 

\medskip

However, these assertions do hold in the cases considered in \secref{sss:IndCoh* indsch fact good}, see \cite[Sect. 4]{GLC2}.

\sssec{}

An assertion parallel to \lemref{l:descent for IndCoh^*} holds also for $\IndCoh^!(-)$. Hence, by the same mechanism as in 
\secref{sss:2nd defn IndCoh G/NP}, we obtain a well-defined factorization category 
$$\IndCoh^!(\fL_\nabla(\cG/\cN^-_P)).$$

\medskip

We now claim:

\begin{prop} \label{p:IndCoh^! dual of IndCoh^*}
The factorization categories $$\IndCoh^*(\fL_\nabla(\cG/\cN^-_P)) \text{ and } \IndCoh^!(\fL_\nabla(\cG/\cN^-_P))$$ are naturally dual.
\end{prop} 

\begin{proof} 

By construction, we need to establish a duality between 
$$\IndCoh^*(\fL_\nabla(\cG/\cN^-_P))_{X^I} \text{ and } \IndCoh^!(\fL_\nabla(\cG/\cN^-_P))_{X^I}$$
for finite sets $I$, compatibly with maps between finite sets. 

\medskip

Note that, up to Zariski sheafification, $\fL_\nabla(\cG/\cN^-_P)_{X^I}$ can be identified with the quotient of
$\fL_\nabla(\cG)_{X^I}$ by $\fL_\nabla(\cN^-_P)_{X^I}$. Hence, by the $\IndCoh^!(-)$-version of \lemref{l:descent for IndCoh^*},
we can identify
$$\IndCoh^!(\fL_\nabla(\cG/\cN^-_P)_{X^I})\simeq \on{Funct}_{\IndCoh^!(\fL_\nabla(\cN^-_P))_{X^I}}(\QCoh(X^I),\IndCoh^!(\fL_\nabla(\cG))_{X^I}),$$
where $\IndCoh^!(\fL_\nabla(\cN^-_P))_{X^I}$ is viewed as a comonoidal category via the group-structure on $\fL_\nabla(\cN^-_P)_{X^I}$.

\medskip

Now the assertion of the proposition follows from the identifications
$$\IndCoh^!(\fL_\nabla(\cN^-_P))_{X^I}\simeq \left(\IndCoh^*(\fL_\nabla(\cN^-_P))_{X^I}\right)^\vee$$
as comonoidal categories 
and 
$$\IndCoh^!(\fL_\nabla(\cG))_{X^I}\simeq \left(\IndCoh^*(\fL_\nabla(\cG))_{X^I}\right)^\vee$$
as modules over them, see \cite[Sect. E.4.2]{GLC2}.

\end{proof}

\ssec{Construction of the spherical spectral semi-infinite category} \label{ss:semiinf spec}

\sssec{}

We regard $\IndCoh^*(\fL_\nabla(\cG/\cN^-_P))$ as a module category for the monoidal factorization category
$$\IndCoh^*(\fL^+_\nabla(\cG))\otimes \IndCoh^*(\fL^+_\nabla(\cM)).$$ 

Moreover, the action functor has a natural lax unital structure. 

\medskip

Recall now that for an algebraic group $\sH$, the functor
$$\Psi_{\fL^+_\nabla(\sH)}: \IndCoh^*(\fL^+_\nabla(\sH))\to \QCoh(\fL^+_\nabla(\sH))$$
is an equivalence (see \lemref{l:IndCoh vs QCoh arcs}). Hence, we can regard $\IndCoh^*(\fL_\nabla(\cG/\cN^-_P))$ 
as acted on by the group-scheme $\fL^+_\nabla(\cG)\times \fL^+_\nabla(\cM)$.

%

\sssec{}

Given \corref{c:loops into base affine as fiber product}, we define the category 
\begin{equation} \label{e:semiinf spec App}
\IndCoh^*(\LS^\reg_\cG\underset{\LS^\mer_\cG}\times \LS^\mer_{\cP^-}\underset{\LS^\mer_\cM}\times \LS^\reg_\cM)
\end{equation}
by formula \eqref{e:IGP spec defn}. 

\medskip

We now proceed to establishing its properties.

\sssec{} \label{sss:semiinf spec oblv}

Note that we can identify
\begin{multline*} 
\IndCoh^*(\fL^+_\nabla(\cG/\cN^-_P))_{\fL^+_\nabla(\cG)\times \fL^+_\nabla(\cM)}\simeq  \\
\simeq \QCoh(\fL^+_\nabla(\cG/\cN^-_P))_{\fL^+_\nabla(\cG)\times \fL^+_\nabla(\cM)}\simeq \QCoh(\on{pt}/\fL^+_\nabla(\cP^-))=
\QCoh(\LS^\reg_{\cP^-})=:\IndCoh^*(\LS^\reg_{\cP^-}).
\end{multline*} 

Hence, the (monadic) adjunction 
$$\iota^\IndCoh_*:\IndCoh^*(\fL^+_\nabla(\cG/\cN^-_P))\rightleftarrows \IndCoh^*(\fL_\nabla(\cG/\cN^-_P)):\iota^!$$
gives rise to a monadic adjunction
\begin{equation} \label{e:arcs to loops semiinf}
\iota^\IndCoh_*:\IndCoh^*(\LS^\reg_{\cP^-})\rightleftarrows 
\IndCoh^*(\LS^\reg_\cG\underset{\LS^\mer_\cG}\times \LS^\mer_{\cP^-}\underset{\LS^\mer_\cM}\times \LS^\reg_\cM):\iota^!.
\end{equation} 

\sssec{}

A straightforward verification shows that the functor $\iota^\IndCoh_*$ in \eqref{e:arcs to loops semiinf} is strictly unital. 

\medskip

In particular, the factorization unit in $\IndCoh^*(\LS^\reg_\cG\underset{\LS^\mer_\cG}\times \LS^\mer_{\cP^-}\underset{\LS^\mer_\cM}\times \LS^\reg_\cM)$
is given by
$$\iota^\IndCoh_*(\CO_{\LS^\reg_{\cP^-}}).$$

\sssec{}

The fact that the factorization category $\QCoh(\LS^\reg_{\cP^-})$ is compactly generated, combined with the monadic adjunction
\eqref{e:arcs to loops semiinf}, implies that the category \eqref{e:semiinf spec App} is compactly generated. 

\medskip

Namely, compact generators of \eqref{e:semiinf spec App} are provided as the images under $\iota^\IndCoh_*$ of compact 
objects in $\QCoh(\LS^\reg_{\cP^-})\simeq \Rep(\cP^-)$. 

\sssec{} \label{sss:dual semiinf spec}

Note that the operation of coinvariants with respect to $\fL^+_\nabla(\cG)\times \fL^+_\nabla(\cM)$ preserves dualizability
(see \cite[Proposition E.1.4]{GLC2}). Hence, we obtain that the category \eqref{e:semiinf spec App} is dualizable.

\medskip

Furthermore, it follows formally from \propref{p:IndCoh^! dual of IndCoh^*} that its dual identifies with
$$\IndCoh^!(\fL_\nabla(\cG/\cN^-_P))^{\fL^+_\nabla(\cG)\times \fL^+_\nabla(\cM)}\simeq
\IndCoh^!(\LS^\reg_\cG\underset{\LS^\mer_\cG}\times \LS^\mer_{\cP^-}\underset{\LS^\mer_\cM}\times \LS^\reg_\cM).$$

\sssec{} \label{sss:gen Jacquet} 

Note that the functor
$$\Gamma^\IndCoh(\fL_\nabla(\cG),-):\IndCoh^*(\fL_\nabla(\cG))\to \Vect$$
is $\fL_\nabla(\cG)\times \fL_\nabla(\cG)$-invariant, and hence gives rise to a functor 
$$\Gamma^\IndCoh(\fL_\nabla(\cG/\cN^-_P),-):\IndCoh^*(\fL_\nabla(\cG/\cN^-_P))\to \Vect,$$
which is in turn $\fL_\nabla(\cG)\times \fL_\nabla(\cM)$-invariant.

\medskip

In particular, the functor $\Gamma^\IndCoh(\fL_\nabla(\cG/\cN^-_P),-)$ gives rise to a functor
\begin{multline*} 
\IndCoh^*(\LS^\reg_\cG\underset{\LS^\mer_\cG}\times \LS^\mer_{\cP^-}\underset{\LS^\mer_\cM}\times \LS^\reg_\cM):=\\
:=\IndCoh^*(\fL_\nabla(\cG/\cN^-_P))_{\fL^+_\nabla(\cG)\times \fL^+_\nabla(\cM)}\to \Vect_{\fL^+_\nabla(\cM)}=\Rep(\cM).
\end{multline*} 

This is the sought-for functor $(p_2)^\IndCoh_*$ in formula \eqref{e:gen Jacquet functor spec}. 

\sssec{}  \label{sss:Hecke prelim} 

Recall the construction from \cite[Sect. E.6]{GLC2}. It can be reformulated as follows: given a (factorization)
category $\bC$ equipped with an action of the monoidal (factorization) category $\IndCoh^*(\fL_\nabla(\sH))$, then the (factorization)
category 
$$\bC_{\fL^+_\nabla(\sH)}\simeq \bC\underset{\IndCoh^*(\fL^+_\nabla(\sH))}\otimes \Vect$$ acquires an action of
$$\IndCoh^*(\fL_\nabla(\sH))_{\fL^+_\nabla(\sH)\times \fL^+_\nabla(\sH)}=:\Sph_\sH^{\on{spec}}.$$

\medskip

Moreover, if $\bC$ was unital and the action of $\IndCoh^*(\fL_\nabla(\sH))$ on $\bC$ had a lax unital structure, then the resulting
action of $\Sph_\sH^{\on{spec}}$ on $\bC_{\fL^+_\nabla(\sH)}$ is strictly unital. 

\sssec{}  \label{sss:semiinf Hecke} 

In particular, the action of 
$$\IndCoh^*(\fL_\nabla(\cG))\otimes \IndCoh^*(\fL_\nabla(\cM))$$ 
on $\IndCoh^*(\fL_\nabla(\cG/\cN^-_P))$ gives rise to an action of
$$\IndCoh^*(\fL_\nabla(\cG))_{\fL^+_\nabla(\cG)\times \fL^+_\nabla(\cG)}\otimes \IndCoh^*(\fL_\nabla(\cM))_{\fL^+_\nabla(\cM)\times \fL^+_\nabla(\cM)}=
\Sph_\cG^{\on{spec}}\otimes \Sph_\cM^{\on{spec}}$$
on 
$$\IndCoh^*(\fL_\nabla(\cG/\cN^-_P))_{\fL^+_\nabla(\cG)\times \fL^+_\nabla(\cM)}=:
\IndCoh^*(\LS^\reg_\cG\underset{\LS^\mer_\cG}\times \LS^\mer_{\cP^-}\underset{\LS^\mer_\cM}\times \LS^\reg_\cM).$$

\medskip

This provides the action of the spherical categories promised in \secref{sss:spec Hecke action on semiinf}. 

\sssec{} \label{sss:gen Jacquet Sph} 

Note that the functor
\begin{equation} \label{e:gen Jacquet App}
\IndCoh^*(\fL_\nabla(\cG/\cN^-_P))_{\fL^+_\nabla(\cM)}\to \Vect_{\fL^+_\nabla(\cM)}
\end{equation}
considered in \secref{sss:gen Jacquet} is compatible with the actions of $\IndCoh^*(\fL_\nabla(\cM))_{\fL^+_\nabla(\cM)\times \fL^+_\nabla(\cM)}$
on the two sides, and it is invariant with respect to the action of $\IndCoh^*(\fL_\nabla(\cG))$ on the left-hand side. 

\medskip

This implies that the functor
\begin{multline*}
\Vect_{\fL^+_\nabla(\cG)}\otimes \IndCoh^*(\fL_\nabla(\cG/\cN^-_P))_{\fL^+_\nabla(\cG)\times \fL^+_\nabla(\cM)}
\overset{\otimes}\to \\
\to \IndCoh^*(\fL_\nabla(\cG/\cN^-_P))_{\fL^+_\nabla(\cG)\times \fL^+_\nabla(\cM)}\to \Vect_{\fL^+_\nabla(\cM)}
\end{multline*}
factors as
\begin{multline*}
\Vect_{\fL^+_\nabla(\cG)}\otimes \IndCoh^*(\fL_\nabla(\cG/\cN^-_P))_{\fL^+_\nabla(\cG)\times \fL^+_\nabla(\cM)}\to \\
\to \Vect_{\fL^+_\nabla(\cG)} \underset{\IndCoh^*(\fL_\nabla(\cG))_{\fL^+_\nabla(\cG)\times \fL^+_\nabla(\cG)}}\otimes 
\IndCoh^*(\fL_\nabla(\cG/\cN^-_P))_{\fL^+_\nabla(\cG)\times \fL^+_\nabla(\cM)}\to  \Vect_{\fL^+_\nabla(\cM)},
\end{multline*}
where the second arrow is still compatible with the actions of $\IndCoh^*(\fL_\nabla(\cM))_{\fL^+_\nabla(\cM)\times \fL^+_\nabla(\cM)}$. 

\medskip

In other words, we obtain the desired factorization of the functor 
\begin{equation} \label{e:J spec pre}
J^{-,\on{spec,pre-enh}}:\Rep(\cG)\otimes \on{I}(\cG,\cP^-)^{\on{spec,loc}}\to \Rep(\cM)
\end{equation} 
via a functor 
$$J^{-,\on{spec,enh}}: \Rep(\cG)\underset{\Sph^{\on{spec}}_\cG}\otimes \on{I}(\cG,\cP^-)^{\on{spec,loc}}\to \Rep(\cM),$$
which is $\Sph^{\on{spec}}_\cM$-linear. 


\sssec{}

In the rest of this subsection we will construct the pair of adjoint functors from \secref{sss:Omega tilde spec Ran}. 
Consider the fiber product $\LS^\reg_{\cP^-}\underset{\LS^\mer_\cM}\times \LS^\reg_\cM$. 

\medskip

As in \cite[Sect. E.2]{GLC2}, we \emph{define} 
$$\IndCoh^*(\LS^\reg_{\cP^-}\underset{\LS^\mer_\cM}\times \LS^\reg_\cM):=\IndCoh^*(\fL_\nabla(\cM))_{\fL^+_\nabla(\cP^-)\times \fL^+_\nabla(\cM)}.$$

\medskip

Let us first construct the $\Sph_\cM^{\on{spec}}$-linear monadic adjunction
$$\IndCoh^*(\on{Hecke}^{\on{spec,loc}}_{\cM}):=
\IndCoh^*(\fL_\nabla(\cM))_{\fL^+_\nabla(\cM)\times \fL^+_\nabla(\cM)}
\rightleftarrows \IndCoh^*(\LS^\reg_{\cP^-}\underset{\LS^\mer_\cM}\times \LS^\reg_\cM).$$

\sssec{}

We have the tautological functor
\begin{equation} \label{e:from LS P base change M prel}
\IndCoh^*(\fL_\nabla(\cM))_{\fL^+_\nabla(\cP^-)}\to
\IndCoh^*(\fL_\nabla(\cM))_{\fL^+_\nabla(\cM)},
\end{equation} 
where the coinvariants are taken with respect to the action by left translations.

\medskip

It is compatible with the actions of $\IndCoh^*(\fL_\nabla(\cM))$ by right translations, and hence induces a functor
\begin{equation} \label{e:from LS P base change M}
\IndCoh^*(\fL_\nabla(\cM))_{\fL^+_\nabla(\cP^-)\times \fL^+_\nabla(\cM)}\to
\IndCoh^*(\fL_\nabla(\cM))_{\fL^+_\nabla(\cM)\times \fL^+_\nabla(\cM)},
\end{equation} 
compatible with the actions of $\IndCoh^*(\fL_\nabla(\cM))_{\fL^+_\nabla(\cM)\times \fL^+_\nabla(\cM)}=:\Sph_\cM^{\on{spec}}$. 

\medskip

It is easy to see that the functor \eqref{e:from LS P base change M} is conservative.\footnote{E.g., conservativity can be checked pointwise.}

\sssec{}

We claim that the functor \eqref{e:from LS P base change M} admits a left adjoint. To prove this, it suffices to show that
the functor \eqref{e:from LS P base change M prel} admits a left adjoint, and, equivalently, that the functor
\begin{equation} \label{e:from LS P base change M prel N}
\IndCoh^*(\fL_\nabla(\cM))_{\fL^+_\nabla(\cN^-_P)}\to \IndCoh^*(\fL_\nabla(\cM))
\end{equation}
admits a left adjoint. 

\medskip

However, the latter is obvious from the identification
$$\IndCoh^*(\fL_\nabla(\cM))_{\fL^+_\nabla(\cN^-_P)}\simeq \IndCoh^*(\fL_\nabla(\cM))\otimes \Rep(\cN^-_P),$$
so that the functor \eqref{e:from LS P base change M prel N} is induced by the functor of $\cN^-_P$-invariants
$$\Rep(\cN^-_P)\to \Vect,$$
whose left adjoint is the functor of the trivial representation. 

\sssec{}

We now construct the $\Sph_\cM^{\on{spec}}$-linear monadic adjunction
$$\IndCoh^*(\LS^\reg_{\cP^-}\underset{\LS^\mer_\cM}\times \LS^\reg_\cM) \rightleftarrows \IndCoh^*(\on{Hecke}^{\on{spec,loc}}_{\cG,\cP^-}).$$

\medskip

The map $\cM\to \cG/\cN^-_P$ gives rise to the $(\IndCoh,*)$-pushforward functor 
$$\IndCoh^*(\fL_\nabla(\cM))\to \IndCoh^*(\fL_\nabla(\cG/\cN^-_P)),$$
which commutes with the actions of $\IndCoh^*(\fL_\nabla(\cM))$ on the right, and further to a functor
$$\IndCoh^*(\fL_\nabla(\cM))_{\fL^+_\nabla(\cP^-)}\to \IndCoh^*(\fL_\nabla(\cG/\cN^-_P))_{\fL^+_\nabla(\cG)},$$
which is still compatible with the actions of $\IndCoh^*(\fL_\nabla(\cM))$ on the right.

\medskip

Passing to $\fL^+_\nabla(\cM)$-coinvariants, we obtain a functor
\begin{equation} \label{e:from LS P base change G}
\IndCoh^*(\fL_\nabla(\cM))_{\fL^+_\nabla(\cP^-)\times \fL^+_\nabla(\cM)}\to \IndCoh^*(\fL_\nabla(\cG/\cN^-_P))_{\fL^+_\nabla(\cG)\times \fL^+_\nabla(\cM)}
\end{equation} 
compatible with the actions of $\IndCoh^*(\fL_\nabla(\cM))_{\fL^+_\nabla(\cM)\times \fL^+_\nabla(\cM)}$. 

\sssec{}

We claim that the functor \eqref{e:from LS P base change G} admits a (conservative) right adjoint as a factorization functor.
To prove this, it suffices to show that the functor \eqref{e:from LS P base change G} preserves compactness. 

\medskip

For that we note that the category 
$$\IndCoh^*(\fL_\nabla(\cM))_{\fL^+_\nabla(\cP^-)\times \fL^+_\nabla(\cM)}$$
is compactly generated by the essential image of the functor
\begin{equation} \label{e:arcs M to loops}
\QCoh(\LS^\reg_{\cP^-})\simeq \IndCoh^*(\fL^+_\nabla(\cM))_{\fL^+_\nabla(\cP^-)\times \fL^+_\nabla(\cM)} \overset{\iota^\IndCoh_*}\to
\IndCoh^*(\fL_\nabla(\cM))_{\fL^+_\nabla(\cP^-)\times \fL^+_\nabla(\cM)}.
\end{equation}

Now, the composition of \eqref{e:arcs M to loops} and \eqref{e:from LS P base change G} is the left adjoint in
\eqref{e:arcs to loops semiinf}, which preserves compactness.

\ssec{Construction of the category  \texorpdfstring{$\IndCoh^*((\LS^\Mmf_{\cP^-})^\wedge_\mf)$}{amb}} \label{ss:ambient spec semiinf}

In this subsection we define another semi-infinite spectral category, namely, $\IndCoh^*((\LS^\Mmf_{\cP^-})^\wedge_\mf)$,
as was promised in \secref{ss:semiinf CS}. 

\sssec{} \label{sss:ambient spec semiinf}

Let $(\LS^\Mmf_{\cP^-})^\wedge_\mf$ be as in \secref{sss:ambient spec semiinf init}. 

\medskip

Our current goal is to define 
$$\IndCoh^*\left((\LS^\Mmf_{\cP^-})^\wedge_\mf\right) 
\text{ and } \IndCoh^!\left((\LS^\Mmf_{\cP^-})^\wedge_\mf\right)$$
as unital factorization categories.

\medskip

The idea of the definition is as follows: we have
$$\LS^\reg_{\cP^-}\simeq \on{pt}/\fL^+_\nabla(\cP^-)=\on{pt}/\fL^+_\nabla(\cN^-_P)\cdot \fL^+_\nabla(\cM).$$

Hence, we can identify $(\LS^\Mmf_{\cP^-})^\wedge_\mf$ with the \'etale quotient 
$$\on{pt}/\fL_\nabla(\cN^-_P)\cdot \fL^+_\nabla(\cM),$$
and we will use this presentation to define $\IndCoh^*(-)$ and $\IndCoh^!(-)$ of $(\LS^\Mmf_{\cP^-})^\wedge_\mf$. 

\sssec{}

We set
$$\IndCoh^*\left((\LS^\Mmf_{\cP^-})^\wedge_\mf\right) :=
\Vect\underset{\IndCoh^*(\fL_\nabla(\cN^-_P)\cdot \fL^+_\nabla(\cM))}\otimes \Vect$$ 
and 
$$\IndCoh^!\left((\LS^\Mmf_{\cP^-})^\wedge_\mf\right) :=
\on{Funct}_{\IndCoh^!(\fL_\nabla(\cN^-_P)\cdot \fL^+_\nabla(\cM))}(\Vect,\Vect).$$ 

\sssec{} \label{sss:ambient via IGP}

For future reference, we note that we have
\begin{multline*}
\Vect\underset{\IndCoh^*(\fL_\nabla(\cN^-_P)\cdot \fL^+_\nabla(\cM))}\otimes \Vect \simeq \\
\simeq 
\IndCoh^*(\fL_\nabla(\cG))\underset{\IndCoh^*(\fL_\nabla(\cG)\times \fL_\nabla(\cN^-_P)\cdot \fL^+_\nabla(\cM))}\otimes \Vect \simeq \\
\simeq \IndCoh^*(\fL_\nabla(\cG/\cN^-_P))\underset{\IndCoh^*(\fL_\nabla(\cG)\times \fL^+_\nabla(\cM))}\otimes \Vect,
\end{multline*}
and so  
$$\IndCoh^*\left((\LS^\Mmf_{\cP^-})^\wedge_\mf\right) \simeq
\IndCoh^*(\fL_\nabla(\cG/\cN^-_P))\underset{\IndCoh^*(\fL_\nabla(\cG)\times \fL^+_\nabla(\cM))}\otimes \Vect.$$

Similarly,
$$\IndCoh^!\left((\LS^\Mmf_{\cP^-})^\wedge_\mf\right) \simeq
\on{Funct}_{\IndCoh^!(\fL_\nabla(\cG)\times \fL^+_\nabla(\cM))}(\Vect,\IndCoh^!(\fL_\nabla(\cG/\cN^-_P))).$$

In particular, we have
\begin{equation} \label{e:amb vs IGP spec *}
\IndCoh^*\left((\LS^\Mmf_{\cP^-})^\wedge_\mf\right) \simeq
\Rep(\cG)\underset{\Sph_\cG^{\on{spec}}}\otimes \on{I}(\cG,\cP^-)^{\on{spec,loc}}
\end{equation} 
and
\begin{equation} \label{e:amb vs IGP spec !}
\IndCoh^!\left((\LS^\Mmf_{\cP^-})^\wedge_\mf\right) \simeq
\Rep(\cG)\overset{(\Sph_\cG^{\on{spec}})^\vee}\otimes \on{I}(\cG,\cP^-)^{\on{spec,loc}}_{\on{co}},
\end{equation} 
where:

\begin{itemize}

\item $(\Sph_\cG^{\on{spec}})^\vee$ is the dual of $\Sph_\cG^{\on{spec}}$, viewed as a comonoidal category;

\item The notation $\bC_1\overset{\bB}\otimes \bC_2$ is an in \secref{sss:tensor vs cotensor}. 

\end{itemize}

\sssec{} \label{sss:P to amb semiinf}

Let $'\iota$ denote the map
$$\LS^\reg_{\cP^-}\to (\LS^\Mmf_{\cP^-})^\wedge_\mf.$$

We define the (strictly unital) factorization functor
$$({}'\!\iota)^\IndCoh_*:\IndCoh^*(\LS^\reg_{\cP^-})\to \IndCoh^*\left((\LS^\Mmf_{\cP^-})^\wedge_\mf\right)$$
to be the composition 
\begin{multline} \label{e:iota into ambient form}
\IndCoh^*(\LS^\reg_{\cP^-})\simeq \Vect\underset{\IndCoh^*(\fL^+_\nabla(\cP^-))}\otimes \Vect=
\Vect\underset{\IndCoh^*(\fL^+_\nabla(\cN^-_P)\cdot \fL^+_\nabla(\cM))}\otimes \Vect\to \\
\to \Vect\underset{\IndCoh^*(\fL_\nabla(\cN^-_P)\cdot \fL^+_\nabla(\cM))}\otimes \Vect=
\IndCoh^*\left((\LS^\Mmf_{\cP^-})^\wedge_\mf\right).
\end{multline}

\medskip

We claim that the functor \eqref{e:iota into ambient form}, admits a right adjoint as a factorization functor, to be denoted 
$({}'\!\iota)^!$. 

\sssec{} \label{sss:right adj coinv}

More generally, we claim that if $\IndCoh^*(\fL_\nabla(\sH))$ acts on a sheaf of categories $\ul\bC$ over $\Ran$,
then then functor 
$$\ul\bC\underset{\IndCoh^*(\fL^+_\nabla(\sH))}\otimes \Vect\to \ul\bC\underset{\IndCoh^*(\fL_\nabla(\sH))}\otimes \Vect$$
admits a right adjoint as a factorization functor. 

\medskip

Indeed, we write the above functor as 
\begin{multline} 
\ul\bC\underset{\IndCoh^*(\fL^+_\nabla(\sH))}\otimes \Vect\simeq 
\ul\bC\underset{\IndCoh^*(\fL^+_\nabla(\sH)\times \fL_\nabla(\sH))}\otimes \IndCoh^*( \fL_\nabla(\sH))\simeq \\
\simeq \ul\bC\underset{\IndCoh^*(\fL_\nabla(\sH))}\otimes \IndCoh^*(\fL_\nabla(\sH)/\fL^+_\nabla(\sH))
\overset{\on{Id}\otimes \Gamma^\IndCoh(\fL_\nabla(\sH)/\fL^+_\nabla(\sH),-)}\longrightarrow
\ul\bC\underset{\IndCoh^*(\fL_\nabla(\sH))}\otimes \Vect,
\end{multline} 
and the assertion follows from the fact that the functor 
$$\Gamma^\IndCoh(\fL_\nabla(\sH)/\fL^+_\nabla(\sH),-):\IndCoh^*(\fL_\nabla(\sH)/\fL^+_\nabla(\sH))\to \Vect$$
admits a $\fL_\nabla(\sH)$-equivariant right adjoint.

\sssec{}

Thus, we obtain that the (factorization) category 
$$\IndCoh^*\left((\LS^\Mmf_{\cP^-})^\wedge_\mf\right)$$
is compactly generated, and in particular, is dualizable.

\medskip

By the same mechanism as in \secref{sss:dual semiinf spec}, we obtain an identification
$$\IndCoh^*\left((\LS^\Mmf_{\cP^-})^\wedge_\mf\right)^\vee \simeq
\IndCoh^!\left((\LS^\Mmf_{\cP^-})^\wedge_\mf\right).$$

\ssec{t-structures on the spectral semi-infinite categories}

\sssec{} \label{sss:t semiinf spec}

We claim that the factorization category $\IndCoh^*\left((\LS^\Mmf_{\cP^-})^\wedge_\mf\right)$ carries a (unique) t-structure, so that the 
functor 
$$\Rep(\cP^-) \simeq \IndCoh^*(\LS^\reg_{\cP^-})\overset{({}'\!\iota)^\IndCoh_*}\longrightarrow \IndCoh^*\left((\LS^\Mmf_{\cP^-})^\wedge_\mf\right)$$
is t-exact. 

\medskip

Since the adjunction 
$$({}'\!\iota)^\IndCoh_*:\IndCoh^*(\LS^\reg_{\cP^-})\rightleftarrows \IndCoh^*\left((\LS^\Mmf_{\cP^-})^\wedge_\mf\right):({}'\!\iota)^!$$
monadic, this assertion is equivalent to the following:

\begin{prop} \label{p:t semiinf spec}
The endofunctor of $\Rep(\cP^-)\simeq \IndCoh^*(\LS^\reg_{\cP^-})$ underlying the monad
$$({}'\!\iota)^!\circ ({}'\!\iota)^\IndCoh_*$$
is left t-exact.
\end{prop} 

\begin{proof} 

We have a commutative diagram

\smallskip

\begin{equation} \label{e:t semiinf spec t-exact} 
\CD
\IndCoh(\fL_\nabla(\cP^-)/\fL^+_\nabla(\cP^-)) @>{(\sfq^{-,\on{spec}})^\IndCoh_*}>> \IndCoh(\fL_\nabla(\cM)/\fL^+_\nabla(\cM)) \\
@VVV @VVV \\
\Vect\underset{\IndCoh^*(\fL_\nabla(\cP^-))}\otimes \IndCoh(\fL_\nabla(\cP^-)/\fL^+_\nabla(\cP^-)) @>>> 
\Vect\underset{\IndCoh^*(\fL_\nabla(\cM))}\otimes  \IndCoh(\fL_\nabla(\cP^-)/\fL^+_\nabla(\cM)) \\
@V{\sim}VV @VV{\sim}V \\
\IndCoh^*(\LS^\reg_{\cP^-}) @>{({}'\!\iota)^\IndCoh_*}>> \IndCoh^*\left((\LS^\Mmf_{\cP^-})^\wedge_\mf\right),
\endCD
\end{equation} 
such that the natural transformation in diagram
$$
\xy
(0,0)*+{\IndCoh(\fL_\nabla(\cP^-)/\fL^+_\nabla(\cP^-))}="A";
(80,0)*+{\IndCoh(\fL_\nabla(\cM)/\fL^+_\nabla(\cM))}="B";
(0,-30)*+{\IndCoh^*(\LS^\reg_{\cP^-})}="C";
(80,-30)*+{\IndCoh^*\left((\LS^\Mmf_{\cP^-})^\wedge_\mf\right),}="D";
{\ar@{->} "A";"C"};
{\ar@{->}_{(\sfq^{-,\on{spec}})^!} "B";"A"};
{\ar@{->} "B";"D"};
{\ar@{->}^{({}'\!\iota)^!} "D";"C"};
{\ar@{=>} "A";"D"};
\endxy
$$
obtained from \eqref{e:t semiinf spec t-exact} by passing to right adjoints along the horizontal arrows, is an isomorphism.

\medskip

We claim that the left vertical arrow in diagram \eqref{e:t semiinf spec t-exact} sends the coconnective subcategory in $\IndCoh(\fL_\nabla(\cP^-)/\fL^+_\nabla(\cP^-))$
surjectively onto the coconnective subcategory in 
$$\IndCoh^*(\LS^\reg_{\cP^-})\simeq \QCoh(\LS^\reg_{\cP^-})\simeq \Rep(\cP^-).$$ 
Indeed, this functor admits a right inverse given by direct image along
$$\on{pt}/\fL^+_\nabla(\cP^-)\to \fL_\nabla(\cP^-)/\fL^+_\nabla(\cP^-).$$

Hence, in order to prove that $({}'\!\iota)^!\circ ({}'\!\iota)^\IndCoh_*$ is left t-exact, it suffices to show that the endofunctor of $\IndCoh(\fL_\nabla(\cP^-)/\fL^+_\nabla(\cP^-))$
underlying the monad 
$$(\sfq^{-,\on{spec}})^!\circ (\sfq^{-,\on{spec}})^\IndCoh_*$$
is left t-exact.

\medskip

However, this follows from the fact that the functor
$$(\sfq^{-,\on{spec}})^\IndCoh_*:\IndCoh(\fL_\nabla(\cP^-)/\fL^+_\nabla(\cP^-))\to \IndCoh(\fL_\nabla(\cM)/\fL^+_\nabla(\cM))$$
is left t-exact.

\end{proof} 

\sssec{} \label{sss:t IGP spec}

Similarly, we claim that the category $\on{I}(\cG,\cP^-)^{\on{spec,loc}}$
carries a (unique) t-structure, so that the functor 
$$\iota^\IndCoh_*:\Rep(\cP^-)\to \on{I}(\cG,\cP^-)^{\on{spec,loc}}$$
is t-exact, where $\iota$ denotes here the map
$$\LS^\reg_{\cP^-}\to \LS^\reg_\cG\underset{\LS^\mer_\cG}\times \LS^\mer_{\cP^-}\underset{\LS^\mer_\cM}\times \LS^\reg_\cM.$$

\medskip

Since the adjunction in \eqref{e:arcs to loops semiinf} is monadic, the statement is equivalent to the fact that the endofunctor of
$\Rep(\cP^-)$ underlying the monad $\iota^!\circ \iota^\IndCoh_*$ is left t-exact. This is proved in a way similar to \propref{p:t semiinf spec}. 

\ssec{The global counterpart} \label{ss:global IGP gen}

The goal of this subsection is to construct the functors from Sects. \ref{ss:Eis spec enh} and \ref{ss:CT spec enh} in families over $\Ran$. 

\sssec{} \label{sss:LS P G mf app}

First, let us explain how to use the techniques
developed earlier in this section to define the factorization category $\IndCoh^*(\LS^{\mf_\cG}_{\cP^-})$. 

\medskip

Namely,
$$\IndCoh^*(\LS^{\mf_\cG}_{\cP^-}):=\IndCoh^*(\fL_\nabla(\cG))\underset{\IndCoh^*(\fL_\nabla^+(\cG))\otimes \IndCoh^*(\fL_\nabla(\cP^-))}\otimes \Vect.$$

Note that by the same principle as in \eqref{e:amb vs IGP spec *}, we have:
$$\IndCoh^*(\LS^{\mf_\cG}_{\cP^-})\simeq \Rep(\cM)\underset{\Sph_\cM^{\on{spec}}}\otimes \on{I}(\cG,\cP^-)^{\on{spec,loc}}.$$

\sssec{} \label{sss:global IGP gen-first}

Next, we proceed to the construction of the functor
$$(\sfq^{-,\on{spec}})^{\IndCoh,*,\on{enh}}_\CZ:\IndCoh(\LS_\cM)^{-,\on{enh}}_\CZ\to \IndCoh(\LS^{\mf_\cG,\on{glob}}_{\cP^-,\CZ})$$
from \secref{sss:global IGP gen}. 

\medskip

Consider the fiber product 
$$\LS^{\on{mer,glob}}_{\cM,\CZ}\underset{\LS^\mer_{\cM,\CZ}}\times \LS^{\mf_\cG}_{\cP^-,\CZ}.$$

Since 
$$(\LS_\cM\times \CZ)\underset{\LS^\reg_{\cM,\CZ}}\times \LS^\reg_{\cP^-,\CZ}\to 
\LS^{\on{mer,glob}}_{\cM,\CZ}\underset{\LS^\mer_{\cM,\CZ}}\times \LS^{\mf_\cG}_{\cP^-,\CZ}$$
is an isomorphism at the classical level, as in \secref{sss:ambient spec semiinf}, we can interpret it as
$$\left((\LS_\cM\times \CZ)^{\on{level}}\underset{\CZ}\times (\fL_\nabla^+(\cG)_\CZ\backslash \fL_\nabla(\cG/\cN^-_P)_\CZ)\right)/\fL_\nabla(\cM)_\CZ,$$
where $(\LS_\cM\times \CZ)^{\on{level}}$ is as in \cite[Sect. E.7.3]{GLC2}.

\sssec{}

Recall that 
$$\CZ\rightsquigarrow \IndCoh^*((\LS_\cM\times \CZ)^{\on{level}})$$
is well-defined as a crystal of categories over $\Ran$, equipped with action of 
$$\CZ\rightsquigarrow \IndCoh^*(\fL_\nabla(\cM))_\CZ.$$

Explicitly
$$\IndCoh^*((\LS_\cM\times \CZ)^{\on{level}})\simeq \left(\IndCoh(\LS_\cM)\otimes \Dmod(\CZ)\right)\underset{\IndCoh^*(\fL_\nabla^+(\cM))_\CZ}
\otimes \Dmod(\CZ).$$

\sssec{}

Set
\begin{multline}  \label{e:global IGP gen defn}
\IndCoh^*\left(\LS^{\on{mer,glob}}_{\cM,\CZ}\underset{\LS^\mer_{\cM,\CZ}}\times \LS^{\mf_\cG}_{\cP^-,\CZ}\right):=\\
=\left(\IndCoh^*((\LS_\cM\times \CZ)^{\on{level}}) \underset{\Dmod(\CZ)}\otimes \IndCoh^*(\fL_\nabla(\cG/\cN^-_P))_\CZ\right)
\underset{\IndCoh^*(\fL_\nabla^+(G)\times \fL_\nabla(\cM))_\CZ}\otimes \Dmod(\CZ).
\end{multline} 

We can rewrite the above expression also as
\begin{equation}  \label{e:global IGP gen defn 1}
\left(\IndCoh^*((\LS_\cM\times \CZ)^{\on{level}}) \underset{\Dmod(\CZ)}\otimes \IndCoh^*(\fL_\nabla(\cG))_\CZ\right)
\underset{\IndCoh^*(\fL_\nabla^+(G)\times \fL_\nabla(\cP^-))_\CZ}\otimes \Dmod(\CZ).
\end{equation} 

\sssec{}

Note, however, that by the same principle as in \eqref{e:amb vs IGP spec *}, we have
\begin{multline*} 
\IndCoh^*\left(\LS^{\on{mer,glob}}_{\cM,\CZ}\underset{\LS^\mer_{\cM,\CZ}}\times \LS^{\mf_\cG}_{\cP^-,\CZ}\right)\simeq  \\
\simeq (\IndCoh(\LS_\cM)\otimes \Dmod(\CZ))\underset{\Sph^{\on{spec}}_{\cM,\CZ}}\otimes 
\IndCoh^*(\fL_\nabla^+(G)\backslash \fL_\nabla(\cG/\cN^-_P)/\fL_\nabla^+(\cM))_\CZ = \\
=:\IndCoh(\LS_\cM)^{-,\on{enh}}_\CZ.
\end{multline*} 

\sssec{}

By the same principle as in \secref{sss:ambient spec semiinf}, we also have 
$$\LS^{\mf_\cG,\on{glob}}_{\cP^-,\CZ}\simeq 
\left((\LS_{\cP^-}\times \CZ)^{\on{level}} \underset{\CZ}\times (\fL_\nabla^+(\cG)\backslash \fL_\nabla(\cG))_\CZ\right)/\fL_\nabla(P)_\CZ.$$

Hence, as in \propref{p:IndCoh'}, we have
\begin{multline}  \label{e:LS G P mer}
\IndCoh(\LS^{\mf_\cG,\on{glob}}_{\cP^-,\CZ})\simeq \\
\simeq \left(\IndCoh^*((\LS_{\cP^-}\times \CZ)^{\on{level}}) \underset{\Dmod(\CZ)}\otimes \IndCoh^*(\fL_\nabla(\cG))_\CZ\right)
\underset{\IndCoh^*(\fL_\nabla^+(G)\times \fL_\nabla(\cP^-))_\CZ}\otimes \Dmod(\CZ).
\end{multline}

\sssec{} \label{sss:global IGP gen-last}

Hence, we can can interpret the sought-for functor $(\sfq^{-,\on{spec}})^{\IndCoh,*,\on{enh}}_\CZ$ as a functor
\begin{multline} \label{e:sfq - enh}
\left(\IndCoh^*((\LS_\cM\times \CZ)^{\on{level}}) \underset{\Dmod(\CZ)}\otimes \IndCoh^*(\fL_\nabla(\cG))_\CZ\right)
\underset{\IndCoh^*(\fL_\nabla^+(G)\times \fL_\nabla(\cP^-))_\CZ}\otimes \Dmod(\CZ)\to \\
\to \left(\IndCoh^*((\LS_{\cP^-}\times \CZ)^{\on{level}}) \underset{\Dmod(\CZ)}\otimes \IndCoh^*(\fL_\nabla(\cG))_\CZ\right)
\underset{\IndCoh^*(\fL_\nabla^+(G)\times \fL_\nabla(\cP^-))_\CZ}\otimes \Dmod(\CZ).
\end{multline}

The latter functor is induced by the functor
$$(\sfq^{-,\on{spec,level}}_\CZ)^{\IndCoh,*}: 
\IndCoh^*((\LS_\cM\times \CZ)^{\on{level}}) \to \IndCoh^*((\LS_{\cP^-}\times \CZ)^{\on{level}}))$$
left adjoint to the 
$$(\sfq^{-,\on{spec,level}}_\CZ)^\IndCoh_*:\IndCoh^*((\LS_{\cP^-}\times \CZ)^{\on{level}})\to
\IndCoh^*((\LS_\cM\times \CZ)^{\on{level}}).$$

Explicitly, the functor $(\sfq^{-,\on{spec,level}}_\CZ)^{\IndCoh,*}$ equals 
\begin{multline*}
\IndCoh^*((\LS_\cM\times \CZ)^{\on{level}}) \simeq
(\IndCoh^*(\LS_\cM)\otimes \Dmod(\CZ))\underset{\IndCoh^*(\fL_\nabla^+(\cM))}\otimes \Dmod(\CZ)
\overset{(\sfq^{-,\on{spec}})^*\otimes \on{Id}}\longrightarrow  \\
\to (\IndCoh^*(\LS_{\cP^-})\otimes \Dmod(\CZ))\underset{\IndCoh^*(\fL_\nabla^+(\cM))}\otimes \Dmod(\CZ)\to \\
\to (\IndCoh^*(\LS_{\cP^-})\otimes \Dmod(\CZ))\underset{\IndCoh^*(\fL_\nabla^+(\cP^-))}\otimes \Dmod(\CZ)\simeq 
\IndCoh^*((\LS_{\cP^-}\times \CZ)^{\on{level}}).
\end{multline*} 

\sssec{} \label{sss:global IGP gen M}

We now proceed to the definition of the functor
$$(\sfp^{-,\on{spec}}_\CZ)^{!,\on{enh}}:\IndCoh(\LS_\cG)^{-,\on{enh}_{\on{co}}}_\CZ\to \IndCoh(\LS^{\mf_\cM,\on{glob}}_{\cP^-,\CZ}).$$

First, as above, we identify 
\begin{multline*}
\IndCoh(\LS_\cG)^{-,\on{enh}_{\on{co}}}_\CZ
\simeq \IndCoh^!\left(\left(\LS^{\on{mer,glob}}_{\cG,\CZ}\right){}^\wedge_\mf\underset{\LS^\mer_{\cG,\CZ}}\times \LS^{\mf_\cM}_{\cP^-,\CZ}\right):=\\
=\on{Funct}_{\IndCoh^!(\fL_\nabla(\cN^-_P)_\CZ\cdot \fL_\nabla^+(\cM)_\CZ)}\left(\Dmod(\CZ),\IndCoh^!((\LS_\cG\times \CZ)^{\on{level}})\right). 
\end{multline*}

Furthermore, we identify
$$\IndCoh(\LS^{\mf_\cM,\on{glob}}_{\cP^-,\CZ})\simeq 
\on{Funct}_{\IndCoh^!(\fL_\nabla(\cN^-_P)_\CZ\cdot \fL_\nabla^+(\cM)_\CZ)}\left(\Dmod(\CZ),\IndCoh^!((\LS_{\cP^-}\times \CZ)^{\on{level}})\right).$$

In terms of these identifications, the sought-for functor $(\sfp^{-,\on{spec}}_\CZ)^{!,\on{enh}}$ is induced by the !-pullback functor
$$\IndCoh^!((\LS_{\cP^-}\times \CZ)^{\on{level}})\to \IndCoh^!((\LS_\cG\times \CZ)^{\on{level}}).$$

\section{Proof of semi-infinite Casselman-Shalika} \label{s:Sith}

The goal of this section is to prove \thmref{t:semiinf CS}, which serves as a preparation for the proof of 
\thmref{t:semiinf geom Satake}. 

\medskip

Note that \thmref{t:semiinf CS} ``crosses the Langlands bridge", i.e., we we need to find a way to connect
a category on the geometric side to a category on the spectral side. This is done by interpreting\footnote{Here ``interpreting"=``mapping
\emph{almost} fully faithfully."} both sides as the category of factorization modules over the \emph{same} factorization algebra.

\ssec{Orientation session}

The material in this section is closely related what has already appeared in the existing literature, 
To orient the reader, we highlight these connections and the points of departure. 

\sssec{}

First, the origin of these ideas is \cite{AB}, which proves a version of these results in the pointwise/Iwahori setting.

\sssec{} \label{sss:orient}

Second, \cite{Ra4} considers a close variant of the main result here. There are two key differences. 

\medskip

One is that \cite{Ra4} considers $\fL(N^-)\cdot \fL^+(T)$-\emph{coinvariants}, where here we consider 
\emph{invariants} (and a general parabolic). 

\medskip

The other is that  \cite{Ra4} considers the factorization algebra 
$\Upsilon$ (cf. \emph{loc. cit}.), where we consider the factorization algebra $\Omega$. We encourage the reader to think that 
$\Omega$ is fundamentally connected to invariants while $\Upsilon$ is fundamentally connected to coinvariants; 
we hope to study this idea further in future works.

\sssec{}

Finally, \cite{YaR} proves a pointwise (and quantum!) version of our results, i.e., recovers this work recovers the main 
result of \cite{AB} but using factorization algebras in the main constructions.

\sssec{}

Compared to the previous works, the main novelty of the present text is a ULA-ness result for $\Omega$, 
\thmref{t:Omega comp} below, which builds on finiteness results from \cite{CR}. 
A parallel result for $\Upsilon$ was given in \cite{Ra3} and plays a key technical role in \cite{Ra4}. 
Using this result, we are able to reduce \thmref{t:semiinf CS} to the pointwise theorem of \cite{YaR}.

\ssec{The semi-infinite Whittaker category}

In this section we will interpret the geometric side of \thmref{t:semiinf CS} as 
the category of factorization modules over a factorization algebra. 

\medskip

While doing so, we will provide a more conceptual point of view on the category in 
the geometric side of \thmref{t:semiinf CS}. 

\sssec{}

Consider the factorization category
$$\Dmod_{\frac{1}{2}}(\fL(G)).$$

Imposing the condition of $\fL(N^-_P)\cdot \fL^+(M)$-equivariance with respect to the action by right translations, we obtain a unital \emph{lax}\footnote{It is possible
that the lax factorization structure is strict, but we do not know how to prove this.} 
factorization category
$$\Dmod_{\frac{1}{2}}(\fL(G))^{-,\semiinf}:=\Dmod_{\frac{1}{2}}(\fL(G))^{\fL(N^-_P)\cdot \fL^+(M)}.$$

\sssec{}

The category $\Dmod_{\frac{1}{2}}(\fL(G))^{-,\semiinf}$ is a equipped with a (strictly unital) factorization functor
\begin{multline} \label{e:access}
\Dmod_{\frac{1}{2}}(\Gr_G)\underset{\Sph_G}\otimes \on{I}(G,P^-) \to
\Dmod_{\frac{1}{2}}(\Gr_G)\underset{\Sph_G}\otimes \Dmod_{\frac{1}{2}}(\Gr_G)^{\fL(N^-_P)\cdot \fL^+(M)}\to \\
\to \Dmod_{\frac{1}{2}}(\fL(G))^{\fL(N^-_P)\cdot \fL^+(M)}=:\Dmod_{\frac{1}{2}}(\fL(G))^{-,\semiinf},
\end{multline}
to be denoted $\Accs$. 

\medskip

By the same principle as in \lemref{l:KM to Sph gen}, one shows that the functor $\Accs$ is fully faithful.  We will
denote its essential image by 
$$\Dmod_{\frac{1}{2}}(\fL(G))^{-,\semiinfAccs}\subset \Dmod_{\frac{1}{2}}(\fL(G))^{-,\semiinf};$$
and refer to it as the \emph{spherically accessible} subcategory. 

\medskip

The lax factorization structure on $\Dmod_{\frac{1}{2}}(\fL(G))^{-,\semiinf}$ induces a genuine factorization structure on 
$\Dmod_{\frac{1}{2}}(\fL(G))^{-,\semiinfAccs}$. 

\sssec{}

We have a naturally defined functor 
$$J^{-,\semiinf}:\Dmod_{\frac{1}{2}}(\fL(G))^{-,\semiinf}\to \Dmod_{\frac{1}{2}}(\Gr_M),$$
given by
\begin{multline*} 
\Dmod_{\frac{1}{2}}(\fL(G))^{-,\semiinf}\overset{\sim}\to \left(\Dmod_{\frac{1}{2}}(\fL(G))\otimes \Dmod_{\frac{1}{2}}(\Gr_M)\right)^{\fL(P^-)} \to \\
\to \Dmod_{\frac{1}{2}}(\fL(G))\otimes \Dmod_{\frac{1}{2}}(\Gr_M) \to 
\Dmod_{\frac{1}{2}}(\Gr_M)\overset{\on{shift}}\longrightarrow \Dmod_{\frac{1}{2}}(\Gr_M),
\end{multline*}
where:

\begin{itemize}

\item The second arrow is the functor of forgetting the $\fL(P^-)$-equivariance;

\item The third arrow is given by !-restriction along $1_{\fL(G)}\to \fL(G)$;

\item The last arrow is the cohomological shift by $\langle \lambda,2\rhoch_P\rangle$ on 
$\Gr_M^\lambda$.

\end{itemize}

We observe:

\begin{lem}  \label{l:Gr G to M}
The precomposition of the functor $J^{-,\semiinf}$ with the functor \eqref{e:access} is the functor
$J_\Gr^{-,\on{enh}}$ of \eqref{e:J Gr untwisted}. 
\end{lem}

\begin{proof} 

Unwinding the definitions, the assertion boils down to the fact that composition
$$\Dmod_{\frac{1}{2}}(\Gr_G)\underset{\Sph_G}\otimes \Dmod_{\frac{1}{2}}(\Gr_G)\overset{\star}\to 
\Dmod_{\frac{1}{2}}(\fL(G))\overset{!\text{-fiber\,at\,}1_{\fL(G)}}\longrightarrow \Vect$$
identifies with the Verdier duality pairing.

\end{proof}

\sssec{}

By a similar token, we can consider the twisted versions
$$\Dmod_{\frac{1}{2}}(\fL(G)_{\rho(\omega_X)}) \text{ and } \Dmod_{\frac{1}{2}}(\fL(G)_{\rho(\omega_X)})^{-,\semiinf},$$
respectively, and the corresponding functors
\begin{equation} \label{e:access twisted}
\Accs:\Dmod_{\frac{1}{2}}(\Gr_{G,\rho(\omega_X)})\underset{\Sph_G}\otimes \on{I}(G,P^-)_{\rho_P(\omega_X)} \to \Dmod_{\frac{1}{2}}(\fL(G)_{\rho(\omega_X)}) 
\end{equation}
and 
\begin{equation} \label{e:Jacquet on LG}
J^{-,\semiinf}: \Dmod_{\frac{1}{2}}(\fL(G)_{\rho(\omega_X)})^{-,\semiinf}\to \Dmod_{\frac{1}{2}}(\Gr_{M,\rho_M(\omega_X)}),
\end{equation} 
see \secref{sss:tricky twisted enh Jacq}. 

\sssec{}

Set
$$\Whit^!(\fL(G))^{-,\semiinf}:=\Whit^!(\Dmod_{\frac{1}{2}}(\fL(G)_{\rho(\omega_X)})^{-,\semiinf}),$$
see \cite[Sect. Sect. 1.3.3]{GLC2} for the notation $\Whit^!(-)$. 

\begin{rem}

One could equivalently define $\Whit^!(\fL(G))^{-,\semiinf}$ as
$$\left(\Whit^!(\Dmod(\fL(G)_{\rho(\omega_X)}))\right)^{-,\semiinf}:=
\left(\Whit^!(\Dmod(\fL(G)_{\rho(\omega_X)}))\right)^{\fL(N^-_P)_{\rho(\omega_X)}\cdot \fL^+(M)_{\rho(\omega_X)}}.$$

\end{rem} 

\sssec{}

The functor \eqref{e:access twisted} induces a fully faithful (strictly unital) factorization functor
\begin{equation} \label{e:access Whit}
\Accs:\Whit^!(G) \underset{\Sph_G}\otimes \on{I}(G,P^-)_{\rho_P(\omega_X)} \to \Whit^!(\fL(G))^{-,\semiinf}. 
\end{equation} 

We will denote its essential image by  
$$\Whit^!(\fL(G))^{-,\semiinfAccs}\subset \Whit^!(\fL(G))^{-,\semiinf},$$
and refer to it as the \emph{spherically accessible} subcategory.

\medskip

The category $\Whit^!(\fL(G))^{-,\semiinf}$ carries a lax factorization structure, and it induces a genuine
factorization structure on $\Whit^!(\fL(G))^{-,\semiinfAccs}$. 

\medskip

The functor $J^{-,\semiinf}$ induces a functor
$$\Whit^!(\fL(G))^{-,\semiinf}\to \Whit^!(M),$$
to be denoted $J_{\Whit}^{-,\semiinf}$.

\medskip

By \lemref{l:Gr G to M}, the composition of the functors $\Accs$ and $J_{\Whit}^{-,\semiinf}$
is the functor $J_{\Whit}^{-,\on{enh}}$ of \eqref{e:Whit semiinf ult}. 

\sssec{}

Let $\Omega$ be the (unital) factorization algebra in $\Whit^!(M)$ equal to the image of the factorization unit 
under the functor $J_{\Whit}^{-,\semiinf}$.

\medskip

Note that this is the same as the image of the factorization unit in $\Whit^!(\fL(G))^{-,\semiinfAccs}$ along the functor 
$J_{\Whit}^{-,\on{enh}}$.

\medskip

Still equivalently, this is the same as the image of the unit in $\Whit^!(G)$ along the functor\footnote{Because  
because the averaging functor $\Whit^!(G) \to \Whit^!(\fL(G))^{-,\semiinfAccs}$
is strictly unital.} $J_{\Whit}^{-,!}$.

\sssec{}

Thus, the functor $J_{\Whit}^{-,\semiinf}$ gives rise to a (strictly unital) factorization functor
\begin{equation} \label{e:J Whit semiinf enh}
\Whit^!(\fL(G))^{-,\semiinf}\to \Omega\mod^{\on{fact}}(\Whit^!(M)),
\end{equation} 
to be denoted $J_{\Whit,\Omega}^{-,\semiinf}$. 



\sssec{}

We will prove: 

\begin{thm} \label{t:Whit vs Omega} \hfill

\smallskip

\noindent{\em(a)} The functor $J_{\Whit,\Omega}^{-,\semiinf}$ is a \emph{pointwise} equivalence, i.e.,
it induces an equivalence 
$$\Whit^!(\fL(G))^{-,\semiinf}_{\ul{x}}\to \Omega\mod^{\on{fact}}(\Whit^!(M))_{\ul{x}}$$
for every $\ul{x}\in \Ran$.

\smallskip

\noindent{\em(b)} The restriction of $J_{\Whit,\Omega}^{-,\semiinf}$ to $\Whit^!(\fL(G))^{-,\semiinfAccs}$ 
is fully faithful, i.e., for every affine scheme $S\to \Ran$,
the corresponding functor
$$\Whit^!(\fL(G))^{-,\semiinfAccs}_S\to \Omega\mod^{\on{fact}}(\Whit^!(M))_S$$
is fully faithful.

\end{thm} 

Point (a) of \thmref{t:Whit vs Omega} will be proved in \secref{ss:Ruotao}; as we shall explain in {\it loc. cit.}, 
\thmref{t:Whit vs Omega}(a) is essentially a paraphrase of the (parabolic version of the) \emph{classical}\footnote{Trivial=no twising}
case of the main result of \cite{YaR}, namely Theorem 6.4.8 in {\it loc. cit.}.

\medskip

Point (b) \thmref{t:Whit vs Omega} will be proved in \secref{sss:proof Whit vs Omega}. 

\begin{rem}

In \thmref{t:Whit vs Omega bis} we will provide a more precise version of \thmref{t:Whit vs Omega}(b). Namely, we will
describe explicitly the essential image of the restriction of $J_{\Whit,\Omega}^{-,\semiinf}$ to $\Whit^!(\fL(G))^{-,\semiinfAccs}$.

\end{rem} 

\ssec{Digression: monoidal actions of factorization categories} \label{ss:universal action}

\sssec{} \label{sss:universal action 0}

Let $\bA$ be a factorization category, and let $\bC$ be a monoidal factorization category 
that acts on $\bA$. 

\medskip

We will assume that both $\bC$ and $\bA$ are unital as factorization categories. Moreover, we will
assume that the monoidal structure on $\bC$ and the action of $\bC$ on $\bA$ are strictly unital.

\medskip

In particular, the monoidal unit in $\bC$, thought of as a factorization functor $\Vect\to \bC$, 
is unital, i.e., the monoidal unit in $\bC$ equals the factorization unit. 

\sssec{} \label{sss:universal action 1}

Let $\CA$ be a factorization algebra in $\bA$. Note that in this case the action of $\bC$ canonically extends to
an action of $\bC$ on the (lax factorization) category 
$$\CA\mod^{\on{fact}}(\bA),$$
in a way compatible with the forgetful functor.

\medskip

Indeed, the corresponding binary operation is given by
$$\bC\otimes \CA\mod^{\on{fact}}(\bA)\simeq \one_\bC\mod^{\on{fact}}(\bC)\otimes \CA\mod^{\on{fact}}(\bA)\to \CA\mod^{\on{fact}}(\bA),$$
where the second arrow is induced by action $\bC\otimes \bA\to \bA$, and it sends 
$$\one_\bC\otimes \CA\mapsto \CA,$$
by assumption. 

\sssec{} \label{sss:universal action 2}

Let now $\bA_1$ and $\bA_2$ be a pair of factorization categories, both acted on by $\bC$. 
Let $$\Phi:\bA_1\to \bA_2$$
be a factorization functor, compatible with the action of $\bC$. 

\medskip

Let $\CA$ be a factorization algebra in $\bA_1$, so that $\Phi$ extends to a functor 
$$\Phi_\CA:\CA\mod^{\on{fact}}(\bA_1)\to \Phi(\CA)\mod^{\on{fact}}(\bA_2).$$ 

Then the functor $\Phi_\CA$ is compatible with the actions of $\bC$ on the two sides.

\sssec{} \label{sss:action of Sph on Omega mod 1}

We apply the paradigm of \secref{sss:universal action 1} to $\bA:=\Whit^!(M)$ and $\bC:=\Sph_M$. From \secref{sss:universal action 1}
we obtain that the (lax factorization) category 
$\Omega\mod^{\on{fact}}(\Whit^!(M))$ carries a canonical $\Sph_M$-action.

\sssec{} \label{sss:action of Sph on Omega mod 2}

We apply the paradigm of \secref{sss:universal action 2} to $\Phi:=J_{\Whit}^{-,\semiinf}$ and $\CA$ being the factorization unit. 
We obtain that the functor $J_{\Whit,\Omega}^{-,\semiinf}$ is compatible with the actions of $\Sph_M$
on the two sides. 

\ssec{The spectral counterpart}

In this subsection we will interpret the spectral side in \thmref{t:semiinf CS} as the category of modules
over a factorization algebra. 

\sssec{}

Recall the (unital) factorization category $\IndCoh^*((\LS^\Mmf_{\cP^-})^\wedge_\mf)$, which is equipped with a
(strictly unital) factorization functor 
\begin{equation} \label{e:i for LS N}
\Rep(\cP^-) \simeq \QCoh(\LS_{\cP^-}^\reg)
\simeq \IndCoh^*(\LS_{\cP^-}^\reg) \overset{({}'\!\iota)^\IndCoh_*}\longrightarrow \IndCoh^*((\LS^\Mmf_{\cP^-})^\wedge_\mf).
\end{equation} 

\medskip

Recall the (lax unital) factorization functor
$$\IndCoh^*((\LS^\Mmf_{\cP^-})^\wedge_\mf)\overset{((\sfq^{-,\on{spec}})^\wedge_\mf)^\IndCoh_*}\longrightarrow \IndCoh^*(\LS_\cM^\reg)\simeq \QCoh(\LS_\cM^\reg) \simeq \Rep(\cM),$$
see \secref{sss:functor to Rep M}. 

\medskip

By construction, the composition $((\sfq^{-,\on{spec}})^\wedge_\mf)^\IndCoh_*\circ ({}'\!\iota)^\IndCoh_*$ is the functor of direct image 
\begin{equation} \label{e:sfp - spec}
(\sfq^{-,\on{spec}})^\IndCoh_*:\IndCoh^*(\LS_{\cP^-}^\reg)\to \IndCoh^*(\LS_\cM^\reg),
\end{equation} 
i.e., this is the functor
$$\on{inv}_{\cN^-_P}:\Rep(\cP^-)\to \Rep(\cM).$$

\sssec{}

In particular, we obtain that the functor 
$((\sfq^{-,\on{spec}})^\wedge_\mf)^\IndCoh_*$ sends the factorization unit to 
the (commutative) factorization algebra $\Omega^{\on{spec}}\in\Rep(\cM)$, see \secref{sss:Omega spec}.

\medskip

Thus, the functor $((\sfq^{-,\on{spec}})^\wedge_\mf)^\IndCoh_*$  induces a (strictly unital) factorization functor
$$\IndCoh^*((\LS^\Mmf_{\cP^-})^\wedge_\mf)\to \Omega^{\on{spec}}\mod^{\on{fact}}(\Rep(\cM)),$$
to be denoted by $\sQ^\wedge_\mf$.

\sssec{} \label{sss:action of Sph on Omega spec mod}

Recall the action of $\Sph_\cM^{\on{spec}}$ on $\Rep(\cM)$ (see \cite[Sect. E.6.5]{GLC2}). 
As in \secref{sss:action of Sph on Omega mod 1}, we obtain that the category $\Omega^{\on{spec}}\mod^{\on{fact}}(\Rep(\cM))$ 
carries a canonical action of $\Sph_\cM^{\on{spec}}$.

\medskip

As in \secref{sss:action of Sph on Omega mod 2}, we obtain that the functor $\sQ^\wedge_\mf$ respects the actions of 
$\Sph_\cM^{\on{spec}}$ on the two sides. 

\sssec{}

We will prove:

\begin{thm} \label{t:LS N vs Omega}
The functor $\sQ^\wedge_\mf$ is fully faithful.
\end{thm} 

The proof will be given in \secref{sss:proof LS N vs Omega bis}. 

\begin{rem}

One should view \thmref{t:LS N vs Omega} as a spectral counterpart of point (b) of \thmref{t:Whit vs Omega}.
We will now proceed to formulating a spectral counterpart of \thmref{t:Whit vs Omega}(a). 

\medskip

In addition, in \thmref{t:LS N vs Omega bis} we will describe explicitly the essential image of the functor
\thmref{t:LS N vs Omega}, making this a spectral counterpart of \thmref{t:Whit vs Omega bis}. 

\end{rem}

\sssec{}

Fix a point $\ul{x}\in \Ran$, and note that we can regard 
$$\IndCoh(\LS^\Mmf_{\cP^-,\ul{x}})\simeq \QCoh(\LS^\Mmf_{\cP^-,\ul{x}})$$ as a factorization module category 
over
$$\IndCoh^*(\LS^\reg_{\cP^-})\simeq \QCoh(\LS^\reg_{\cP^-})\simeq \Rep(\cP^-)$$
at $\ul{x}$.

\medskip

The map
$$\sfq^{-,\on{spec}}:\LS^\Mmf_{\cP^-,\ul{x}}\to \LS^\reg_{\cM,\ul{x}}$$
gives rise to a functor
\begin{equation} \label{e:sfp - spec x}
\IndCoh(\LS^\Mmf_{\cP^-,\ul{x}})\overset{(\sfq^{-,\on{spec}})^\IndCoh_*}\longrightarrow \IndCoh^*(\LS^\reg_{\cM,\ul{x}})\simeq 
\QCoh(\LS^\reg_{\cM,\ul{x}})\simeq \Rep(\cM)_{\ul{x}}.
\end{equation} 

Moreover, this functor is compatible with the factorization module structure on the two sides with respect to
$\IndCoh^*(\LS^\reg_{\cP^-})$ and $\IndCoh^*(\LS^\reg_\cM)$, respectively, via the functor \eqref{e:sfp - spec},
see \cite[Sect. B.12.10]{GLC2} for what this means.

\medskip

In particular, by \cite[Sect. C.14.12]{GLC2}, we obtain that the functor \eqref{e:sfp - spec x} gives rise to a functor
\begin{equation} \label{e:Jacquet spec full}
\IndCoh(\LS^\Mmf_{\cP^-,\ul{x}})\to \Omega^{\on{spec}}\mod^{\on{fact}}(\Rep(\cM))_{\ul{x}},
\end{equation} 
to be denoted $\sQ$. 

\sssec{}

The spectral counterpart of \thmref{t:Whit vs Omega}(a) reads:

\begin{thm} \label{t:LS N vs Omega pt}
The functor $\sQ$ of \eqref{e:Jacquet spec full} is an equivalence.
\end{thm} 

The proof will be given in \secref{s:LS N vs Omega pt}. 

\ssec{Identification of factorization algebras}

\sssec{}

The computational core for the proof of \thmref{t:semiinf CS} is provided by the following assertion:

\begin{thm} \label{t:ident Omega}
The equivalence
$$\on{CS}_{\cM,\tau}:\Whit^!(M)\to \Rep(\cM)$$
maps the factorization algebras $\Omega$ and $\Omega^{\on{spec}}$
to one another.
\end{thm}

This theorem is a particular case of \thmref{t:ident Omega bis} just below. 

\sssec{}

We will now discuss an extension of \thmref{t:ident Omega}, which will be used for the 
proof of \thmref{t:semiinf geom Satake}:

\begin{thm} \label{t:ident Omega bis}
The following diagram of (lax unital) factorization functors commutes:
\begin{equation} \label{e:ident Omega bis}
\CD
\Whit^!(G) @>{\on{CS}_{G,\tau}}>>  \Rep(\cG) \\
@V{J^{-,!}_{\Whit}}VV @VV{\on{inv}_{\cN^-_P}}V \\
\Whit^!(M) @>{\on{CS}_{M,\tau}}>>  \Rep(\cM).
\endCD
\end{equation} 
\end{thm}

Note that what we have just stated as \thmref{t:ident Omega bis} is identical to
\corref{c:!-Jacquet on Whit}, see, Remark \ref{r:! Jaq}.

\begin{rem}
Note that the assertion of \thmref{t:ident Omega} is a particular case of that of \thmref{t:ident Omega bis}. Namely, the two sides in 
\thmref{t:ident Omega} are obtained by evaluating the two circuits in \eqref{e:ident Omega bis} on the factorization unit.
\end{rem}

\begin{rem}

A detailed proof of \thmref{t:ident Omega bis} will appear in \cite{FH}. Here is the main idea:

\medskip

In the setting of coinvariants and $\Upsilon$ (see \secref{sss:orient}), as analogue of \thmref{t:ident Omega bis}
was proved in \cite[Theorem 4.15.1]{Ra4} by reducing to a Zastava space calculation, performed \cite[Theorem 7.9.1]{Ra3}. 
This was done in the $P=B$ case, and the case of a general parabolic follows the same ideas. 

\medskip

The invariants/$\Omega$-version follows formally from the coinvariants/$\Upsilon$ case by Verdier duality, once 
we combine with the \emph{cleanness} result \cite[Lemma 4.1.10]{Lin}.\footnote{It is referred to as ``purity" in {\it loc.cit.}}

\end{rem}

\ssec{Finiteness properties of  \texorpdfstring{$\Omega^{\on{spec}}$}{Omega spec fin}}

In this subsection we will state a crucial result, \thmref{t:Omega comp}, pertaining to finiteness properties of 
$\Omega^{\on{spec}}$, viewed as a factorization module over itself. 

\sssec{}

In \secref{s:Omega comp} we will prove: 

\begin{thm} \label{t:Omega comp}
The factorization algebra $\Omega^{\on{spec}}$ is ULA as\footnote{See \cite[Appendix B]{Ra3} or \cite[Sect. 2]{CR} for a detailed discussion of ULA objects in factorization categories.} 
 an object of $\Omega^{\on{spec}}\mod^{\on{fact}}(\Rep(\cM))$, i.e.,
for every affine scheme $S\to \Ran$, the object
$$\Omega^{\on{spec}}_S\in \Omega^{\on{spec}}\mod^{\on{fact}}(\Rep(\cM))_S$$
is compact.
\end{thm} 

\sssec{} \label{sss:access fact}

Recall the action of $\Sph_\cM^{\on{spec}}$ on $\Omega^{\on{spec}}\mod^{\on{fact}}(\Rep(\cM))$ (see \secref{sss:action of Sph on Omega spec mod});
denote the corresponding binary operation by $-\underset{\cM}\star -$. 

\medskip

For every $S\to \Ran$, denote by
$$\Omega^{\on{spec}}\mod^{\on{fact}}(\Rep(\cM))_S^{\Accs}\subset \Omega^{\on{spec}}\mod^{\on{fact}}(\Rep(\cM))_S$$
the full subcategory generated by objects of the form 
\begin{equation} \label{e:gens of accs}
\CF\underset{\cM}\star \Omega^{\on{spec}}_S \in \Omega^{\on{spec}}\mod^{\on{fact}}(\Rep(\cM))_S, \quad \CF\in (\Sph_\cM^{\on{spec}})^c_S.
\end{equation} 

Equivalently, one can take $\CF$ to lie in the essential image of the functor $\on{nv}:\Rep(\cM)^c\to \Sph_\cM^{\on{spec}}$. 

\sssec{}

Note that \thmref{t:Omega comp} (combined with the fact that $(\Sph_\cM^{\on{spec}})_S$ is rigid)
implies that the objects \eqref{e:gens of accs} are compact. 

\medskip

The assignment
$$S\rightsquigarrow \Omega^{\on{spec}}\mod^{\on{fact}}(\Rep(\cM))_S^{\Accs}$$
extends to a (lax) factorization category, which we will denote by
$$\Omega^{\on{spec}}\mod^{\on{fact}}(\Rep(\cM))^{\Accs}.$$

\sssec{}

Recall that the functor $\sQ^\wedge_\mf$ is compatible with the actions of
$\Sph_\cM^{\on{spec}}$ on the two sides (see \secref{sss:action of Sph on Omega spec mod}). 

\medskip

Now, it is easy to see that the category $\IndCoh^*((\LS^\Mmf_{\cP^-})^\wedge_\mf)$ is generated by objects of the form
$$\CF\underset{\cM}\star ({}'\!\iota)^\IndCoh_*(\one_{\LS^\reg_{\cP^-}}), \quad \CF\in \Sph_\cM^{\on{spec}},$$
where $({}'\!\iota)^\IndCoh_*$ is as in \eqref{e:i for LS N}. 

\medskip

Hence, we obtain that the essential image of the functor $\sQ^\wedge_\mf$ is contained in the subcategory
$\Omega^{\on{spec}}\mod^{\on{fact}}(\Rep(\cM))^{\Accs}$.

\sssec{}

We have the following more precise version of \thmref{t:LS N vs Omega}:

\begin{thm} \label{t:LS N vs Omega bis}
The functor $\sQ^\wedge_\mf$ induces an equivalence
$$\IndCoh^*((\LS^\Mmf_{\cP^-})^\wedge_\mf)\to \Omega^{\on{spec}}\mod^{\on{fact}}(\Rep(\cM))^{\Accs}.$$
\end{thm}


\sssec{Proof of \thmref{t:LS N vs Omega bis}} \label{sss:proof LS N vs Omega bis} 

Before we launch into the proof, note that \thmref{t:LS N vs Omega bis} tautologically implies \thmref{t:LS N vs Omega}. 

\medskip

We will prove \thmref{t:LS N vs Omega bis} assuming \thmref{t:LS N vs Omega pt}, which will be proved independently
(see \secref{s:LS N vs Omega pt}). 

\medskip

We need to show that for every affine scheme $S\to \Ran$, the corresponding functor
\begin{equation} \label{e:sfq wedge accs S}
\IndCoh^*((\LS^\Mmf_{\cP^-})^\wedge_\mf)_S\to \Omega^{\on{spec}}\mod^{\on{fact}}(\Rep(\cM))_S^{\Accs}
\end{equation}
is an equivalence. With no restriction of generality we can assume that $S$ is eventually coconnective
(in fact, it is enough to take $S$ to be smooth). 

\medskip

By construction, the essential image of \eqref{e:sfq wedge accs S} generates the target. Hence, it is enough to
show that \eqref{e:sfq wedge accs S} is fully faithful.

\medskip

Since the objects \eqref{e:gens of accs} are compact, we obtain that the functor \eqref{e:sfq wedge accs S} preserves compactness.
Hence, it admits a continuous right adjoint, which is automatically $\QCoh(S)$-linear. 

\medskip

We need to show that the unit of the adjunction is an isomorphism. In order to do so, since $S$ was assumed eventually coconnective,
it is enough to do so at the level of fibers at geometric points of $S$. I.e., we obtain that it is enough to show that for every $\ul{x}\in \Ran$,
the corresponding functor
$$\sQ^\wedge_\mf:\IndCoh^*((\LS^\Mmf_{\cP^-})^\wedge_\mf)_{\ul{x}}\to \Omega^{\on{spec}}\mod^{\on{fact}}(\Rep(\cM))^{\Accs}_{\ul{x}}$$
is fully faithful.

\medskip

Note, however, that we have a commutative diagram
$$
\CD
\IndCoh^*((\LS^\Mmf_{\cP^-})^\wedge_\mf)_{\ul{x}} @>{\sQ^\wedge_\mf}>>
\Omega^{\on{spec}}\mod^{\on{fact}}(\Rep(\cM))^{\Accs}_{\ul{x}} \\
@VVV @VVV \\
\IndCoh(\LS^\Mmf_{\cP^-,\ul{x}}) @>{\sQ}>>  \Omega\mod^{\on{fact}}(\Rep(\cM))_{\ul{x}},
\endCD
$$
where the vertical arrows are the tautological fully faithful embeddings. 

\medskip

Hence, we obtain that the fully faithfulness of $\sQ^\wedge_\mf$ at the pointwise level
follows from \thmref{t:LS N vs Omega pt}. 

\qed[\thmref{t:LS N vs Omega bis}]

\ssec{Finiteness properties of  \texorpdfstring{$\Omega$}{Omega fin}}

In this subsection we will establish results parallel to those in the previous subsection, but now
on the geometric side.

\sssec{}

First, combining Theorems \ref{t:ident Omega} and \ref{t:Omega comp} we obtain: 

\begin{cor} \label{c:Omega comp}
The factorization algebra $\Omega$ is ULA as an object of $\Omega\mod^{\on{fact}}(\Whit^!(M))$.
\end{cor}

\sssec{} \label{sss:Sph acts on factmod}

We define the (lax) factorization category 
$$\Omega\mod^{\on{fact}}(\Whit^!(M))^{\Accs}\subset \Omega\mod^{\on{fact}}(\Whit^!(M))$$
by the procedure parallel to that of \secref{sss:access fact}. 

\sssec{}

Since the essential image of the functor
$$\ind_{\Sph_M\to\semiinf}:\Sph_M\to \on{I}(G,P^-)_{\rho_P(\omega_X)}$$
generates the target, we obtain that the functor $J_{\Whit,\Omega}^{-,\semiinf}$ maps 
$$\Whit^!(\fL(G))^{-,\semiinfAccs}\subset \Whit^!(\fL(G))^{-,\semiinf}$$ to the full subcategory 
$$\Omega\mod^{\on{fact}}(\Whit^!(M))^{\Accs}\subset \Omega\mod^{\on{fact}}(\Whit^!(M)).$$

\sssec{}

We claim:

\begin{thm} \label{t:Whit vs Omega bis}
The resulting functor
$$J_{\Whit,\Omega}^{-,\semiinf}:\Whit^!(\fL(G))^{-,\semiinfAccs}\to \Omega\mod^{\on{fact}}(\Whit^!(M))^{\Accs}$$
is an equivalence.
\end{thm} 

\begin{proof} 

Parallel to the proof of \thmref{t:LS N vs Omega bis} given in \secref{sss:proof LS N vs Omega bis}, using 
\thmref{t:Whit vs Omega}(a), which will be proved independently.

\end{proof}

\sssec{} \label{sss:proof Whit vs Omega}

Note that the assertion of \thmref{t:Whit vs Omega bis} tautologically implies that of \thmref{t:Whit vs Omega}(b). 

\ssec{Proof of \thmref{t:semiinf CS}}

We are now ready to prove \thmref{t:semiinf CS}. 

\medskip

In fact, this is an immediate consequence of the combination 
of Theorems \ref{t:ident Omega}, \ref{t:Whit vs Omega bis} and \ref{t:LS N vs Omega bis}.

\medskip

Namely, $\on{CS}^{-,\semiinf}_\tau$ is the unique functor that makes the following diagram commute
$$
\CD
\Whit^!(\fL(G))^{-,\semiinfAccs} @>{\on{CS}^{-,\semiinf}_\tau}>> \IndCoh^*((\LS^\Mmf_{\cP^-})^\wedge_\mf) \\
@V{J_{\Whit,\Omega}^{-,\semiinf}}V{\sim}V @V{\sim}V{\sQ^\wedge_\mf}V \\
\Omega\mod^{\on{fact}}(\Whit(M))^{\Accs} @>{\sim}>{\on{CS}_{M,\tau}}> \Omega^{\on{spec}}\mod^{\on{fact}}(\Rep(\cM))^{\Accs}. 
\endCD
$$

The functor $\on{CS}^{-,\semiinf}_\tau$ is compatible with the actions of 
$$\Sph_M\overset{\Sat_{M,\tau}}\simeq \Sph_\cM^{\on{spec}},$$
because the functor $\on{CS}_{M,\tau}$ is. 

\qed[\thmref{t:semiinf CS}]

\ssec{Proof of \thmref{t:Whit vs Omega}(a)} \label{ss:Ruotao}

\sssec{}

As in the proof of \lemref{l:KM to Sph gen}, we can replace the category $\Whit^!(\fL(G))^{-,\semiinf}$ by
$$\Whit^!(\on{Fl}^{\on{aff},I^-_P}_G):=\Whit^!\left(\Dmod_{\frac{1}{2}}((\on{Fl}^{\on{aff},I^-_P}_G)_{\rho(\omega_X)})\right),$$
where $\on{Fl}^{\on{aff},I^-_P}_G$ is the partial affine flag scheme
$$\on{Fl}^{\on{aff},I^-_P}_G:=\fL(G)/I^-_P$$
at $\ul{x}$. 

\medskip

Under this identification, the functor $J_{\Whit,\Omega}^{-,\semiinf}$ corresponds to (the parahoric version of) the functor
\begin{equation} \label{e:YaR}
\Whit^!(\on{Fl}^{\on{aff},I^-_P}_G)\to \Omega\mod^{\on{fact}}(\Whit^!(M))_{\ul{x}}
\end{equation} 
of \cite[Equation (6.4.7)]{YaR}. 

\medskip

Hence, the assertion follows from (the parahoric version of) \cite[Theorem 6.4.8]{YaR}. 

\qed[\thmref{t:Whit vs Omega}(a)]

\begin{rem}

When $P=B$, the existence of \emph{an} equivalence between $\Whit^!(\on{Fl}^{\on{aff},I^-_P}_G)$
and $\Omega\mod^{\on{fact}}(\Whit^!(M))_{\ul{x}}$, interpreted as $\IndCoh(\LS^\Mmf_{\cP^-,\ul{x}})$, was established in \cite{AB}. 

\medskip

One can show, however, that their equivalence coincides with one provided by \eqref{e:YaR}. 

\end{rem}

\section{Sheaves on the (pointwise) \texorpdfstring{$\LS_{\sH}^\mer$}{LSunip} for a unipotent group} \label{s:LS N vs Omega pt}

The goal of this section is to prove \thmref{t:LS N vs Omega pt}. 

\medskip

While doing so, we will prove a general
statement concerning the category of ind/quasi-coherent sheaves on the stack $\LS_{\sH,\ol{x}}^\mer$
for a unipotent group $\sH$.

\ssec{The statement for unipotent groups}

\sssec{}

Applying base change with respect to $\on{pt}\to \LS^\reg_\cM$, we reduce \thmref{t:LS N vs Omega pt} to the following assertion,
see \thmref{t:LS H} below, valid for any unipotent group, to be denoted $\sH$, where in our case
$$\sH:=\cN^-_P.$$

\sssec{}

Consider $\IndCoh(\LS^\mer_{\sH,\ul{x}})$ as a factorization module category with respect to
$$\IndCoh^*(\LS^\reg_\sH)\simeq \QCoh(\LS^\reg_\sH)\simeq \Rep(\sH).$$

The functor 
$$\Gamma(\LS^\mer_{\sH,\ul{x}},-)^\IndCoh:\IndCoh(\LS^\mer_{\sH,\ul{x}})\to \Vect$$
is compatible with factorization module structure with respect to the functor
\begin{equation} \label{e:Gamma unip}
\Gamma(\LS_\sH^\reg,-)^\IndCoh:\IndCoh^*(\LS^\reg_\sH)\to \Vect.
\end{equation} 

\sssec{}

The functor \eqref{e:Gamma unip} sends the factorization unit to the (commutative) unital factorization algebra $\Omega^{\on{spec}}_\sH$ (in $\Vect$)
corresponding to the commutative algebra
$$\on{inv}_\sH(k)\in \on{ComAlg}(\Vect),$$
where the latter is the same as the cohomological Chevalley complex 
$\on{C}^\cdot_{\on{chev}}(\sh)$ of $\sh$ with trivial coefficients. 

\medskip

Hence, the functor $\Gamma(\LS^\mer_{\sH,\ul{x}},-)^\IndCoh$ gives rise to a functor
$$\IndCoh(\LS^\mer_{\sH,\ul{x}})\to
\Omega^{\on{spec}}_\sH\mod^{\on{fact}}(\Vect)_{\ul{x}}=:
\Omega^{\on{spec}}_\sH\mod^{\on{fact}}_{\ul{x}},$$
to be denoted $\sQ$.

\sssec{}

We will prove: 

\begin{thm} \label{t:LS H}
The functor $\sQ$ is an equivalence.
\end{thm} 

The proof occupies the rest of this section.

\ssec{Method of proof}

\sssec{}

In order to simplify the notation, we will assume that $\ul{x}$ is a singleton $\{x\}$. However,
the entire discussion is equally applicable to a general $\ul{x}\in \Ran$. 

%
%

\medskip

Recall that $\LS^\mer_{\sH,x}$ is isomorphic to the algebraic stack $\sh/\on{Ad}(\sH)$. In particular, it is smooth,
so the functor
$$\Psi_{\LS^\mer_{\sH,x}}:\IndCoh(\LS^\mer_{\sH,\ul{x}})\to \QCoh(\LS^\mer_{\sH,\ul{x}})$$
is an equivalence. 

\medskip

Since $\sH$ is unipotent, the functor
$$\Gamma(\LS^\mer_{\sH,\ul{x}},-):\QCoh(\LS^\mer_{\sH,\ul{x}})\to \Vect$$
is conservative.

\medskip

Hence, $\CO_{\LS^\mer_{\sH,\ul{x}}}$ is a compact generator of $\IndCoh(\LS^\mer_{\sH,\ul{x}})$.

\sssec{}

Let 
$$\sP_{\Omega^{\on{spec}}_\sH,x}\in \on{Pro}(\Omega^{\on{spec}}_\sH\mod^{\on{fact}}_{\ul{x}})$$
be as in \cite[Sect. D.2.5]{GLC2}. I.e., this is the (pro)-object that (pro)-corepresents the forgetful functor
$$\oblv_{\Omega^{\on{spec}}_\sH}:\Omega^{\on{spec}}_\sH\mod^{\on{fact}}_{\ul{x}}\to \Vect,$$
see \cite[Proposition D.2.7]{GLC2}.

\medskip

We have a tautologically defined map:

\begin{equation} \label{e:P to Omega}
\sP_{\Omega^{\on{spec}}_\sH,x}\to \sQ(\CO_{\LS^\mer_{\sH,\ul{x}}}).
\end{equation}

\sssec{}

We will prove:

\begin{prop} \label{p:P to Omega}
The map \eqref{e:P to Omega} is an isomorphism.
\end{prop} 

Clearly, \propref{p:P to Omega} implies \thmref{t:LS H}. 

\begin{rem}
Note that \propref{p:P to Omega} in particular says that $\sP_{\Omega^{\on{spec}}_\sH,x}$ is an \emph{object} of 
$\Omega^{\on{spec}}_\sH\mod^{\on{fact}}_{\ul{x}}$ (rather than a pro-object). 

\medskip

But we will establish this fact \emph{a priori}.

\end{rem} 

\ssec{Digression: factorization modules over holonomic factorization algebras}

\sssec{} \label{sss:constr ext}

Let $\CA$ be a factorization algebra on $X$. Recall the category $\on{Modif}(\CA)$, see \cite[Sect. D.2.4]{GLC2}.

\medskip

Assume now that $\CA_X$ is \emph{holonomic}. We claim that in this case, the category $\on{Modif}(\CA)$
has an initial object, to be denoted $j_{!,\on{fact}}(\CA|_{X\!-\!x})$. 

\medskip

Namely, consider $\CA_\Ran\in \Dmod(\Ran)$, 
and 
\begin{equation} \label{e:constr ext}
(j_\Ran)_!\circ (j_\Ran)^!(\CA_\Ran)\in \Dmod(\Ran),
\end{equation}
where:

\begin{itemize}

\item $j$ denotes the open embedding $X\!-\!x\hookrightarrow X$;

\item $j_\Ran$ denotes the (open embedding) of the Ran space of $X-x$ into the Ran space of $X$.

\end{itemize} 

It is easy to see that the factorization structure on $\CA_\Ran$ uniquely extends to one on \eqref{e:constr ext},
thereby making \eqref{e:constr ext} into an object of $\on{Modif}(\CA)$. This is our sought-for $j_{!,\on{fact}}(\CA|_{X\!-\!x})$.

\medskip

The fact that it is the initial object in 
this category also follows from the construction. 

\medskip 

If $\CA$ is unital, then $j_{!,\on{fact}}(\CA|_{X\!-\!x})$ is automatically unital: this follows from the fact that the map
$$\omega_\Ran\to (j_\Ran)_!\circ j_\Ran^!(\omega_\Ran)$$
is an isomorphism\footnote{This can be seen by reduction to the case of $X=\BA^1$ and the contraction principle.} 
and \cite[Proposition C.7.13]{GLC2}. 

\sssec{}

Let $j_{!,\on{fact}}(\CA|_{X\!-\!x})_x$ denote the !-fiber of $j_{!,\on{fact}}(\CA|_{X\!-\!x})$, viewed as an object of $\CA\mod^{\on{fact}}_x$. 
By \cite[Proposition D.2.7]{GLC2}, we obtain that $j_{!,\on{fact}}(\CA|_{X\!-\!x})_x$ is the object that corepresents the forgetful functor
$$\oblv_\CA:\CA\mod^{\on{fact}}_x\to \Vect.$$

I.e.,
$$j_{!,\on{fact}}(\CA|_{X\!-\!x})_x\simeq \sP_{\CA,x},$$
where $\sP_{\CA,x}$ is as in \cite[Sect. D.2.5]{GLC2}.

\sssec{} \label{sss:constr ext com}

Assume now that $\CA$ is commutative, i.e., of the form $\on{Fact}(A)$, where
$$A\in \on{ComAlg}(\Dmod(X)),$$
see \cite[Sect. B.10.1]{GLC2}.

\medskip

By assumption, $A$ is holonomic as a D-module on $X$. In this case, the construction of 
$j_{!,\on{fact}}(\CA|_{X\!-\!x})$ given by \eqref{e:constr ext} implies that it has a natural structure of
commutative algebra object in $\on{FactAlg}^{\on{untl}}(X)$, and hence by \cite[Proposition C.8.6]{GLC2}, is of the form $\on{Fact}(j_{!,\on{ComAlg}}\circ j^!(A))$,
where $j_{!,\on{ComAlg}}\circ j^!(A)$ is the initial object in the category
$$\on{Modif}(A):=\{A'\in \on{ComAlg}(\Dmod(X)),\,\, A'|_{X\!-\!x}\simeq A|_{X\!-\!x}\}.$$

\sssec{}

Note that
\begin{equation} \label{e:com fiber}
j_{!,\on{fact}}(\CA|_{X\!-\!x})_x\simeq (j_{!,\on{ComAlg}}\circ j^!(A))_x
\end{equation} 
as objects of
$$\CA\mod^{\on{fact}}_x \simeq j_{!,\on{fact}}(\CA|_{X\!-\!x})\mod^{\on{fact}}_x,$$
where we view $(j_{!,\on{ComAlg}}\circ j^!(A))_x$ as an object of $j_{!,\on{fact}}(\CA|_{X\!-\!x})\mod^{\on{fact}}_x$ via
$$(j_{!,\on{ComAlg}}\circ j^!(A))_x\mod =j_{!,\on{fact}}(\CA|_{X\!-\!x})\mod^{\on{com}}_x\to 
j_{!,\on{fact}}(\CA|_{X\!-\!x})\mod^{\on{fact}}_x,$$
see \cite[Sect. B.10.5-6]{GLC2}. 

\ssec{Proof of \propref{p:P to Omega}}

\sssec{}

We apply the discussion in \secref{sss:constr ext com} to $\CA:=\Omega^{\on{spec}}_\sH$, so that $A=\on{inv}_\sH(k)\otimes \CO_{X_\dr}$. 
Consider the corresponding objects
$$j_{!,\on{fact}}(\Omega^{\on{spec}}_\sH|_{X\!-\!x})$$
and 
$$j_{!,\on{ComAlg}}\circ j^!(\on{inv}_\sH(k)\otimes \CO_{X_\dr})\in \on{ComAlg}(\Dmod(X)).$$

We have:
\begin{equation} \label{e:P as fiber}
j_{!,\on{fact}}(\Omega^{\on{spec}}_\sH|_{X\!-\!x})_x\simeq \sP_{\Omega^{\on{spec}}_\sH,x}
\end{equation} 
as objects of $\Omega^{\on{spec}}_\sH\mod^{\on{fact}}_x$. 

\medskip

In particular, we obtain that $\sP_{\Omega^{\on{spec}}_\sH,x}$ is an \emph{object} (rather than a pro-object) of
$\Omega^{\on{spec}}_\sH\mod^{\on{fact}}_x$. 

\sssec{}

In order to prove \propref{p:P to Omega}, it suffices to show that the map \eqref{e:P to Omega} becomes an isomorphism
after applying the forgetful functor $\oblv_{\Omega^{\on{spec}}_\sH}$. 

\medskip

Thus, using the identifications \eqref{e:com fiber} and \eqref{e:P as fiber}, our goal is to show that the map
\begin{equation} \label{e:P to Omega oblv}
(j_{!,\on{ComAlg}}\circ j^!(\on{inv}_\sH(k)\otimes \CO_{X_\dr}))_x\to \Gamma(\LS_{H,x}^\mer,\CO_{\LS_{H,x}^\mer})
\end{equation} 
is an isomorphism. 

\sssec{} \label{sss:Spec quot}

Recall the construction
$$(R\in \on{ComAlg}(\Vect)) \rightsquigarrow (``\Spec(R)"\in \on{PreStk}),$$
see \cite[Sect. 12.7.6]{GLC2}. Namely, for an affine test-scheme $\Spec(R')$, 
$$\Maps_{\on{PreStk}}(\Spec(R'),``\Spec(R)"):=\Maps_{\on{ComAlg}(\Vect)}(R,R').$$

\medskip

We have a tautologically defined map
\begin{equation} \label{e:Spec R quot map abs}
R\to \Gamma(``\Spec(R)",\CO_{``\Spec(R)"}),
\end{equation}
see \cite[Equation (12.46)]{GLC2}.

\medskip

Recall also (see \cite[Sect. 12.7.6]{GLC2}) that $R$ is said to be ``as good as connective" if:

\begin{itemize}

\item The map \eqref{e:Spec R quot map abs} is an isomorphism;

\medskip

\item The induced functor $\QCoh(``\Spec(R)")\to R\mod$ is an equivalence.

\end{itemize} 

\medskip

This construction makes sense in the relative situation: for a prestack $\CY$ and $R\in \on{ComAlg}(\QCoh(\CY))$
we can consider
$$``\Spec_\CY(R)"\in \on{PreStk}/\CY,$$
so that for $(S\overset{y}\to \CY)\in \affSch_{/\CY}$, we have 
$$S\underset{\CY}\times ``\Spec_\CY(R)"=``\Spec(\Gamma(S,y^*(R)))".$$


\sssec{}

We will deduce the fact that \eqref{e:P to Omega oblv} from 

\begin{prop} \label{p:LS H as quot Spec}
There is a canonical isomorphism
$$\LS_{H,x}^\mer \simeq ``\Spec\left((j_{!,\on{ComAlg}}\circ j^!(\on{inv}_\sH(k)\otimes \CO_{X_\dr}))_x\right)",$$
so that the map \eqref{e:P to Omega oblv} identifies with the map \eqref{e:Spec R quot map abs}.
\end{prop}

\begin{prop} \label{p:as good as conn}
The commutative algebra $$(j_{!,\on{ComAlg}}\circ j^!(\on{inv}_\sH(k)\otimes \omega_X))_x$$
is \emph{as good as connective}.
\end{prop}

\ssec{Proof of \propref{p:LS H as quot Spec}}

\sssec{}

Consider $\on{pt}/\sH\times X_\dr$ as a constant prestack over $X_\dr$. We can identify
$$\LS_{H,x}^\mer\simeq \fL_x(\on{pt}/\sH).$$

Explicitly, for $R'\in \on{ComAlg}(\Vect^{\leq 0})$, we have
\begin{equation} \label{e:LS H as loops}
\Maps(\Spec(R'),\LS_{H,x}^\mer)\simeq
\Maps_{\on{PreStk}/X_\dr}(\Spec_{X_\dr}(R'\ppart),\on{pt}/\sH\times X_\dr),
\end{equation} 
where: 

\begin{itemize}

\item $t$ is a local coordinate at $x$;

\medskip

\item We regard $R\ppart$ as an object of $\on{ComAlg}(\QCoh(X_\dr)^{\leq 0})$.

\end{itemize}

\sssec{}

Since $\sH$ is unipotent, we have a canonical identification
$$\on{pt}/\sH \simeq ``\Spec(\on{inv}_\sH(k))".$$

Hence, we can rewrite the right-hand side in \eqref{e:LS H as loops} as
$$\Maps_{\on{ComAlg}(\QCoh(X_\dr))}(\on{inv}_\sH(k)\otimes \CO_{X_\dr},R'\ppart),$$
which identifies by adjunction with 
$$\Maps_{\on{ComAlg}(\QCoh(X_\dr))}(j_{!,\on{ComAlg}}\circ j^!(\on{inv}_\sH(k)\otimes \CO_{X_\dr}),R'\qqart),$$
and which we rewrite further as
$$\Maps_{\on{ComAlg}(\Vect)}((j_{!,\on{ComAlg}}\circ j^!(\on{inv}_\sH(k)\otimes \CO_{X_\dr}))_x,R'),$$
which is by definition the same as
$$\Maps_{\on{PreStk}}\left(\Spec(R'),``\Spec\left((j_{!,\on{ComAlg}}\circ j^!(\on{inv}_\sH(k)\otimes \CO_{X_\dr}))_x\right)"\right).$$

The identification of the maps \eqref{e:P to Omega oblv} and \eqref{e:Spec R quot map abs} follows from the construction. 

\qed[\propref{p:LS H as quot Spec}]

\ssec{Proof of \propref{p:as good as conn}}

\sssec{}

We identify tautologically
$$(j_{!,\on{ComAlg}}\circ j^!(\on{inv}_\sH(k)\otimes \omega_X))_x\simeq ((j_{\Ran})_!\circ (j_{\Ran})^!(\Omega^{\on{spec}}_\sH))_{\{x\}},$$
where $\{x\}$ is viewed as a point of $\Ran$. 

\medskip

The assertion of the proposition \'etale-local, so we can assume that $$(x\in X)=(0\in \BA^1).$$ In this case, by the contraction principle, we can replace the operation
of taking the !-fiber at the point $\{0\}\in \Ran(\BA^1)$ by $\on{C}^\cdot_c(\Ran(\BA^1),-)$.

\medskip

Further,
$$\on{C}^\cdot_c(\Ran(\BA^1),(j_{\Ran})_!\circ (j_{\Ran})^!(\Omega^{\on{spec}}_\sH))\simeq \on{C}^\cdot_c(\Ran(\BA^1\!-\!0),(j_\Ran)^!(\Omega^{\on{spec}}_\sH)).$$

\medskip

The rest of the proof repeats the proof of the assertion in \cite[Sect. 12.7.7]{GLC2} given in {\it loc. cit.} Sects. 12.7.8-12.7.14. 

\sssec{}

Namely, using \cite[Lemma 12.7.10]{GLC2} we reduce to the case when $\sh$ is abelian. In the latter case
$$\on{inv}_\sH(k)\simeq \Sym(\sh^*[-1]),$$
and 
$$\on{C}^\cdot_c(\Ran(\BA^1\!-\!0),j_\Ran^!(\Omega^{\on{spec}}_\sH)))\simeq \Sym(\sh^*[-1])\otimes \Sym(\sh^*).$$

\noindent (The last isomorphism is due to the fact that for any connected $Y$, the "constant D-algebra" functor 
$\on{ComAlg} \to \on{ComAlg}(\Dmod(Y))$ admits a (partially defined) left adjoint calculated as $\on{C}^\cdot_c(\Ran(Y),\on{Fact}(-))$,
while this left adjoint sends $\Sym(\CM)$, $\CM\in \Dmod(Y)$ to $\Sym(\on{C}^\cdot_c(Y,\CM))$.) 

\qed[\propref{p:as good as conn}]

\begin{rem}
Note that the factorization algebra that we have denoted by $\Omega^{\on{spec}}_\sH$ is \emph{not} the same as
the factorization algebra $\Omega(\sh)$ in \cite[Sect. 12.5.4]{GLC2}.

\medskip

Namely, the latter is $\on{Fact}(\on{C}^\cdot_{\on{chev}}(\sh\otimes \on{D}_X))$,
while the former is $\on{Fact}(\on{C}^\cdot_{\on{chev}}(\sh\otimes k_X))$, where $k_X\in \Dmod(X)$ is ``the constant sheaf"
and $\on{C}^\cdot_{\on{chev}}$ denotes the cohomological Chevalley complex of a Lie-* algebra, see \cite[Sect. 1.4.10]{BD2}.
\end{rem} 

\section{Compactness of \texorpdfstring{$\Omega^{\on{spec}}$}{Omcomp}} \label{s:Omega comp} 

The goal of this section is to prove \thmref{t:Omega comp}. 

\ssec{The shearing trick} \label{ss:shearing trick}

\sssec{}

Consider the cocharacter $2\rhoch_P$ of $Z_\cM$. Note that the symmetric monoidal category $\Rep(\cM)$ admits an automorphism, 
given by sending an object $V$ with central character $\lambda$ to 
$$V[2\langle 2\rhoch_P,\lambda\rangle].$$

Denote this automorphism by $\on{Shear}_\cM$. We refer it to as the \emph{shearing} automorphism of $\Rep(\cM)$. 

\sssec{}

Denote
$$\on{inv}_{\cN^-_P}(k)^{\Rightarrow}:=\on{Shear}_\cM(\on{inv}_{\cN^-_P}(k))\in \on{ComAlg}(\Rep(\cM))$$ and  
$$\Omega^{\on{spec},\Rightarrow}:=\on{Shear}_\cM(\Omega^{\on{spec}})\in \on{FactAlg}^{\on{untl}}(X,\Rep(\cM)),$$
so that 
$$\Omega^{\on{spec},\Rightarrow}\simeq \on{Fact}(\on{inv}_{\cN^-_P}(k)^{\Rightarrow}).$$

\medskip

Clearly, \thmref{t:Omega comp} is equivalent to the assertion that $\Omega^{\on{spec},\Rightarrow}$ is $ULA$
as an object of 
$$\Omega^{\on{spec},\Rightarrow}\mod^{\on{fact}}(\Rep(\cM)).$$

\sssec{}

We observe:

\begin{lem} \label{l:sheared Omega is conn}
The object 
$$\on{inv}_{\cN^-_P}(k)^{\Rightarrow}\in \Rep(\cM)$$
is connective.
\end{lem} 

\begin{proof}

Note that $H^i(\on{inv}_{\cN^-_P}(k))$ is a subquotient of $\Lambda^i((\cn^-_P)^*)$. The result follows now from the fact that
the weights of $2\rhoch_P$ on $(\cn^-_P)^*$ are strictly positive.

\end{proof}

By \cite[Proposition B.9.18]{GLC2}, from \lemref{l:sheared Omega is conn} we obtain that the lax factorization category
$\Omega^{\on{spec},\Rightarrow}\mod^{\on{fact}}(\Rep(\cM))$ carries a t-structure, such that the forgetful functor
$$\oblv_{\Omega^{\on{spec},\Rightarrow}}:\Omega^{\on{spec},\Rightarrow}\mod^{\on{fact}}(\Rep(\cM))\to \Rep(\cM)$$
is t-exact. 

\medskip

Moreover, $\Omega^{\on{spec},\Rightarrow}\mod^{\on{fact}}(\Rep(\cM))$ is left-complete in this t-structure. 

\sssec{} \label{sss:almost compact}

Let $\bD$ be a DG category equipped with a t-structure, compatible with filtered colimits (i.e., the subcategory $\bD^{\geq 0}$ is preserved by
filtered colimts). Recall that an object $\bd\in \bD$ is said to be \emph{almost compact} if the functor
$$\CHom_\bD(\bd,-):\bD\to \Vect$$
is preserved by filtered colimits \emph{when restricted} to $\bD^{\geq -n}$ \emph{for every individual} $n\in \BN$.

\sssec{}

Let now $\bA$ be a lax factorization category equipped with a t-structure\footnote{See \cite[Sect. B.11.11]{GLC2} for what this means}, 
and let let $\CA$ be a factorization algebra object in it. 

\medskip

We shall say that $\CA$
is \emph{almost ULA} in $\CA\mod^{\on{fact}}(\bA)$ if for every affine scheme $S\to \Ran$, the object
$\CA_S\in \CA\mod^{\on{fact}}(\bA)_S$ is almost compact. 

\sssec{}

We claim:

\begin{prop} \label{p:Omega almost ULA}
The factorization algebra  $\Omega^{\on{spec},\Rightarrow}$ is almost ULA in 
$\Omega^{\on{spec},\Rightarrow}\mod^{\on{fact}}(\Rep(\cM))$.
\end{prop}

\begin{proof}

Consider the forgeful functor
$$\oblv_\cM:\Rep(\cM)\to \Vect.$$

It is easy to see that it suffices to show that the factorization algebra $\oblv_\cM(\Omega^{\on{spec},\Rightarrow})$ is almost ULA
in $\oblv_\cM(\Omega^{\on{spec},\Rightarrow})\mod^{\on{fact}}$.

\medskip

Since $\oblv_\cM(\on{inv}_{\cN^-_P}(k)^{\Rightarrow})$ is finite-dimensional, it is almost of finite type. Now
the required assertion follows from \cite[Theorem 4.6.1]{CR}.

\end{proof} 

\ssec{From \emph{almost} compactness to compactness}

In the previous subsection we proved that for every affine scheme $S\to \Ran$, the 
object $$\Omega^{\on{spec},\Rightarrow}_S\in \Omega^{\on{spec},\Rightarrow}\mod^{\on{fact}}(\Rep(\cM))_S$$
is almost compact.

\medskip

In this subsection, we will show that it is actually compact, implying \thmref{t:Omega comp}.

\sssec{}

Let $\bD$ be as in \secref{sss:almost compact}. We shall say that an object $\bd\in \bD$ has a \emph{finite cohomological dimension}
if there exist an integer $k$, such that
$$\Hom_\bD(\bd,\bd'[k'])=0 \text{ for } \bd'\in \bD^\heartsuit \text{ and } k'> k,$$
where
$$\Hom_\bD(-,-):=H^0(\CHom_\bD(-,-)).$$

Such an integer $k$ is called an upper bound on the cohomological dimension. 

\medskip

The following lemma is well-known: 

\begin{lem} \label{l:compact vs almost compact}
Assume that $\bD$ is left-complete in its t-structure. Let $\bd\in \bD$ be almost compact.
Then $\bd$ is compact if and only if it has a finite cohomological dimension of $\bd$. 
\end{lem} 

We include the proof for completeness. 

\begin{proof} 

Suppose that $\bd$ has a finite cohomological dimension; let $k$ be an upper bound on its cohomological dimension.
Note that by the left-completeness assumption, we have
\begin{equation} \label{e:Hom into low}
\Hom_\bD(\bd,\bd')=0 \text{ for any } \bd'\in \bD^{<-(k+1)}.
\end{equation} 

Indeed, 
$$\Hom_\bD(\bd,\bd') = H^0(\CHom_\bD(\bd,\bd'))\simeq 
H^0(\underset{n}{\on{lim}}\, \CHom_\bD(\bd,\tau^{\geq -n}(\bd')),$$
while each $\CHom_\bD(\bd,\tau^{\geq -n}(\bd'))$ has a vanishing $H^0$ and $H^{-1}$
(the latter is needed to kill the term arising from $R^1lim$). 

\medskip

Let 
$$\underset{i}{\on{colim}}\, \bd'_i\to \bd'$$
be a filtered colimit diagram in $\bD$. We need to show that
$$\underset{i}{\on{colim}}\, \Hom_\bD(\bd,\bd'_i)\to \Hom_\bD(\bd,\bd')$$
is an isomorphism, where the colimit in the left-hand side is taken in $\Vect$.

\medskip

Consider the map between the long exact cohomology sequences corresponding to the fiber sequences
$$i\mapsto \, \tau^{< -(k+1)}(\bd'_i) \to \bd'_i\to \tau^{\geq -(k+1)}(\bd'_i)$$ 
and
$$\tau^{< -(k+1)}(\bd')\to \bd'\to \tau^{\geq -(k+1)}(\bd').$$ 

We have 
$$\underset{i}{\on{colim}}\, H^n(\CHom_\bD(\bd,\tau^{< -(k+1)}(\bd'_i))=0 \text{ and } H^n(\CHom_\bD(\bd,\tau^{< -(k+1)}(\bd'))=0
\text{ for } n\geq 0,$$
by \eqref{e:Hom into low}. 

\medskip

The map 
$$\underset{i}{\on{colim}}\, H^n(\CHom_\bD(\bd,\tau^{\geq -(k+1)}(\bd'_i))\to H^n(\CHom_\bD(\bd,\tau^{\geq -(k+1)}(\bd')) \text{ for all } n$$
is an isomorphism, by the assumption of almost compactness.

\medskip

Hence, the map
$$\underset{i}{\on{colim}}\, H^0(\CHom_\bD(\bd,\bd'_i))\to H^0(\CHom_\bD(\bd,\bd'))$$
is an isomorphism, as required. 

\medskip

Assume now that the cohomological dimension of $\bd$ is infinite.
Let $\bd'_k\in \bD^\heartsuit$ be a sequence of objects and $n_k$ an increasing sequence of natural numbers such that
$$\Hom_\bD(\bd,\bd'_k[n_k])\neq 0.$$

Let $\alpha_k\in \Hom_\bD(\bd,\bd'_k[n_k])$
be non-zero elements. 

\medskip

Set 
$$\bd':=\underset{k\in \BN}\oplus\, \bd'_k[n_k] \simeq \underset{k}{\on{colim}}\, \left(\underset{1,...,k}\oplus\, \bd'_k[n_k]\right).$$

By the left completeness assumption, the map
$$\bd'\to \underset{k}{\on{lim}}\, \tau^{\geq -n_k}(\bd') \simeq \underset{k}{\on{lim}}\, \left(\underset{1,...,k}\oplus\, \bd'_k[n_k]\right)\simeq 
\underset{k}{\on{lim}}\, \left(\underset{1,...,k}\Pi\, \bd'_k[n_k]\right)=\underset{k}\Pi\, \bd'_k[n_k]$$ 
is an isomorphism. Hence, the elements $\{\alpha_k\}$ give rise to an element 
$\alpha\in \Hom_\bD(\bd,\bd')$, 
which does not lie in the image of
$$\underset{k}{\on{colim}}\, \Hom_\bD\left(\bd,\left(\underset{1,...,k}\oplus\, \bd'_k[n_k]\right)\right).$$
Hence $\bd$ is not compact. 

\end{proof} 

\sssec{}

Thus, we obtain that in order to prove that $\Omega^{\on{spec},\Rightarrow}_S$ is compact in 
$\Omega^{\on{spec},\Rightarrow}\mod^{\on{fact}}(\Rep(\cM))_S$, it suffices to show that it has a finite cohomological
dimension.

\medskip

Since $\Ran$ is the colimit of prestacks of the form $X^I_\dr$, $I\in \fSet$, 
with no restriction of generality we can assume that $S$ is smooth. 

\medskip

Consider the following general paradigm. Let $\bD$ be an $\QCoh(S)$-linear category (i.e., $\CO_S$ maps to
the Hochschild homology\footnote{A.k.a., derived Bernstein center.} of $\bD$). 

\medskip

Let $\bD$ be equipped with a t-structure, compatible with filtered colimits. For every $f:S'\to S$, denote
$$\bD':=\QCoh(S')\underset{\QCoh(S)}\otimes \bD.$$
We equip $\bD'$ with the tensor product t-structure.

\sssec{}

We claim:

\begin{prop} \label{p:base change fin c.d.}
Assume that $\bd$ is almost compact and that $S$ is smooth. Then $\bd$ is compact if and only if 
there exists an integer $n$ such that for every geometric point $s:\Spec(k')\to S$, 
the object 
$$s^*(\bd)\in \bD_s:=\Vect_{k'}\underset{\QCoh(S)}\otimes \bD$$
has cohomological dimension bounded by $n$. 
\end{prop}

\sssec{}

Before we prove \propref{p:base change fin c.d.}, let us finish the proof of \thmref{t:Omega comp}.

\medskip

By \propref{p:base change fin c.d.}, it is sufficient to show that for a fixed $\ul{x}\in \Ran$, the 
object
$$\Omega^{\on{spec},\Rightarrow}_{\ul{x}}\in \Omega^{\on{spec},\Rightarrow}\mod^{\on{fact}}(\Rep(\cM))_{\ul{x}}$$
has cohomological dimension uniformly bounded by some integer $n$.

\medskip

By \thmref{t:LS N vs Omega pt} (which has been proved independently), 
$$\Omega^{\on{spec}}\mod^{\on{fact}}(\Rep(\cM))_{\ul{x}}\simeq \IndCoh(\LS^\Mmf_{\cP^-,\ul{x}}),$$
and under this equivalence, the object 
$$\Omega^{\on{spec}}_{\ul{x}}\in \Omega^{\on{spec}}\mod^{\on{fact}}(\Rep(\cM))_{\ul{x}}$$
corresponds to
$$i^\IndCoh_*(\CO_{\LS^\reg_{\cP^-,\ul{x}}})\in \IndCoh(\LS^\Mmf_{\cP^-,\ul{x}}).$$
The object $i^\IndCoh_*(\CO_{\LS^\reg_{\cP^-,\ul{x}}})$ is compact. 

\medskip

Hence, $\Omega^{\on{spec},\Rightarrow}_{\ul{x}}$ is compact as an object of $\Omega^{\on{spec},\Rightarrow}\mod^{\on{fact}}(\Rep(\cM))_{\ul{x}}$.
By \lemref{l:compact vs almost compact}, it has a finite cohomological dimension. Let $n_1$ be an integer that bounds its 
cohomological dimension when $\ul{x}$ is a singleton $x$ (it is easy to see that this category does not depend on $x$, up to 
extending scalars). 

\medskip

Let $S\to \Ran$ factor through $X^I$. We take
$$n:=|I| \cdot n_1.$$

By the above, this choice of $n$ has the required property. 

\qed[\thmref{t:Omega comp}]

\ssec{Proof of \propref{p:base change fin c.d.}}

%
%

\sssec{}

For any $\bd,\bd'\in \bD$, let $\ul\CHom_\bD(\bd,\bd')\in \QCoh(S)$ be their inner Hom object, i.e.,
$$\CHom_{\QCoh(S)}(\CF,\ul\CHom_\bD(\bd,\bd')):=\CHom_\bD(\CF\otimes \bd,\bd').$$

\medskip

For $\CF\in \QCoh(S)$, we have a naturally defined map
\begin{equation} \label{e:tensor inner Hom}
\CF\otimes \ul\CHom_\bD(\bd,\bd')\to \ul\CHom_\bD(\bd,\CF\otimes \bd')
\end{equation} 

\sssec{}

We claim:

\begin{lem} \label{l:inner Hom ten}
Assume that $\bd$ is almost compact and $\bd'\in \bD^{>-\infty}$ and $\CF$ is of finite Tor-dimension.
Then the map \eqref{e:tensor inner Hom} is an isomorphism.
\end{lem}

\begin{proof}

The assumption that $\CF$ is of finite Tor-dimension implies that it can be written as a finite colimit
of flat objects. So, we can assume that $\CF$ is flat.

\medskip

In this case, by Lazard's lemma, we can write $\CF$ as a filtered colimit of objects $\CF_i$, where each $\CF_i$ is finitely generated
and free.  It suffices to show that the map
$$\underset{i}{\on{colim}}\, \ul\CHom_\bD(\bd,\CF_i\otimes \bd')\to \ul\CHom_\bD(\bd,\CF\otimes \bd')$$
is an isomorphism, which is equivalent to
$$\underset{i}{\on{colim}}\, \CHom_\bD(\bd,\CF_i\otimes \bd')\to \CHom_\bD(\bd,\CF\otimes \bd')$$
being an isomorphism. 

\medskip

However, the latter holds by the assumption that $\bd$ is almost compact.

\end{proof} 

\begin{cor} \label{c:inner Hom ten}
Let $f:S'\to S$ be a morphism of finite Tor-dimension between affine schemes. Then for $\bd$ almost compact and $\bd'\in \bD^{>-\infty}$, 
the natural map
$$f^*(\ul\CHom_{\bD}(\bd,\bd'))\to \ul\CHom_{\bD'}(f^*(\bd),f^*(\bd'))$$
is an isomorphism.
\end{cor}

\begin{proof}

It is enough to show that the induced map
\begin{equation} \label{e:dir im inner Hom}
f_*\circ f^*(\ul\CHom_{\bD}(\bd,\bd')) \to f_*\left(\ul\CHom_{\bD'}(f^*(\bd),f^*(\bd'))\right)
\end{equation}
is an isomorphism.

\medskip

By the projection formula
$$f_*\circ f^*(\ul\CHom_{\bD}(\bd,\bd'))\simeq f_*(\CO_{S'})\otimes \ul\CHom_{\bD}(\bd,\bd').$$

Furthermore, for any $\wt\bd'\in \bD'$, we have tautologically
$$f_*\left(\ul\CHom_{\bD'}(f^*(\bd),\wt\bd')\right)\simeq \ul\CHom_{\bD}(\bd,f_*(\wt\bd')).$$

In particular,
$$f_*\left(\ul\CHom_{\bD'}(f^*(\bd),f^*(\bd'))\right)\simeq \ul\CHom_{\bD}(\bd,f_*\circ f^*(\bd))\simeq
\ul\CHom_{\bD}(\bd,f_*(\CO_{S'})\otimes \bd).$$

\medskip

By construction, the following diagram commutes
$$
\CD
f_*\circ f^*(\ul\CHom_{\bD}(\bd,\bd'))  @>{\text{\eqref{e:dir im inner Hom}}}>>  f_*\left(\ul\CHom_{\bD'}(f^*(\bd),f^*(\bd'))\right) \\
@V{\sim}VV @VV{\sim}V \\
f_*(\CO_{S'})\otimes \ul\CHom_{\bD}(\bd,\bd') @>{\text{\eqref{e:tensor inner Hom}}}>> \ul\CHom_{\bD}(\bd,f_*(\CO_{S'})\otimes \bd).
\endCD
$$

Hence, the required assertion follows from \lemref{l:inner Hom ten}.

\end{proof}

\sssec{}

We are now ready to prove \propref{p:base change fin c.d.}. The ``only if" direction is evident. We now prove the ``if"
direction.

\medskip

Set $m:=\dim(S)$. Let $n$ be as in the statement of \propref{p:base change fin c.d.}.
We will show that the cohomological dimension of $\bd$ is bounded by $n+2m$. 

\medskip

Fix an object $\bd'\in \bD^\heartsuit$. We wish to show that $\CHom_\bD(\bd,\bd')\in \Vect$ lives in degrees $\leq (k+2 n)$. Consider the object
$$\ul\CHom_\bD(\bd,\bd')\in \QCoh(S).$$

It suffices to show that $\ul\CHom_\bD(\bd,\bd')$ belongs to $\QCoh(S)^{\leq (n+2m)}$. By \cite[Proposition 2.3.2]{GLC1}, it suffices to show that 
for every geometric point $s:\Spec(k')\to S$, 
$$s^*(\ul\CHom_\bD(\bd,\bd'))\in \Vect_{k'}^{\leq n+m}.$$

By \corref{c:inner Hom ten}, we have
$$s^*(\ul\CHom_\bD(\bd,\bd'))\simeq \Hom_{\bD_s}(s^*(\bd),s^*(\bd')).$$

Hence, it suffices to show that the functor 
$$s^*:\bD\to \bD_{s}$$
has cohomological amplitude bounded on the left by $m$.

\medskip

We factor the morphism $s$ as
$$\Spec(k')\overset{s'}\to  S':=\Spec(k')\underset{\Spec(k)}\times S\to S.$$

The pullback functor
$$\bD\to \bD':=\Vect_{k'}\underset{\Vect}\otimes \bD\simeq \QCoh(S')\underset{\QCoh(S)}\otimes \bD$$
is t-exact. Hence, it suffices to show that the functor
$$(s')^*:\bD'\to \bD_s$$
has cohomological amplitude bounded on the left by $m$.

\medskip

The morphism $s'$ is a regular closed embedding of codimension $m$. Hence, we have
$$(s')^*\simeq (s')^![m]$$
(up to tensoring by a line)
where $(s')^!$ is the \emph{right} adjoint to the functor
$$s'_*:\bD'\to \bD_s.$$

Since $s'_*$ is t-exact, the functor $(s')^!$ is left t-exact, implying the desired assertion. 
 
\qed[\propref{p:base change fin c.d.}]

\section{Proof of semi-infinite geometric Satake}  \label{s:semiinf Sat}

The goal of this section is to prove \thmref{t:semiinf geom Satake}. 

\medskip

Like \thmref{t:semiinf CS}, in order to prove \thmref{t:semiinf geom Satake} we need to cross the Langlands bridge,
and the overall strategy of the proof will be similar to that of \thmref{t:semiinf CS}: we will interpret both sides as
the category of factorization modules over the same factorization algebra, but in the factorization category
$\Rep(\cG)\otimes \Rep(\cM)$, rather than just $\Rep(\cM)$. 

\ssec{Another orientation session}

As with \secref{s:Sith}, the material here has a close relationship to existing literature.

\medskip

The first such theorem is \cite{ABG}, which proves a pointwise assertion in the Iwahori context. 
A factorization theorem was recently obtained in \cite{CR}, which used coinvariants and $\Upsilon$ instead of invariants and $\Omega$. 

\medskip

The argument given here is similar to the one in \cite{CR}, making the necessary modifications and reducing to \thmref{t:semiinf CS}
instead of \cite{Ra4}.

\medskip

There is one additional result that was not previously considered: the identification of the semi-infinite IC sheaf 
\propref{p:semiinf Sat IC}, which is proved in \secref{ss:semiinf Sat IC}. 

\ssec{The factorization algebra \texorpdfstring{$\Omega(R)^{\on{spec}}$}{OmegaR}}

In this subsection we will interpret the spectral side in \thmref{t:semiinf geom Satake} as the category
of factorization modules 
over a factorization algebra in $\Rep(\cG)\otimes \Rep(\cM)$.

\sssec{}

Consider the symmetric monoidal categories 
$$\Rep(\cG)\otimes \Rep(\cG),\,\,  \Rep(\cG)\otimes \Rep(\cP^-) \text{ and } \Rep(\cG)\otimes \Rep(\cM),$$
which we turn into factorization categories by the recipe of \cite[Sect. B.11.8]{GLC2}.

\medskip

Let $R_\cG$ be the regular representation of $\cG$, viewed as a (commutative) factorization algebra in 
$\Rep(\cG)\otimes \Rep(\cG)$. 

\medskip

By a slight abuse of notation, we will denote by the same character $R_\cG$
its image under the restriction functor
$$\Rep(\cG)\otimes \Rep(\cG)\to  \Rep(\cG)\otimes \Rep(\cP^-).$$

\sssec{} \label{sss:Omega R}

Consider the (lax unital) factorization functor
$$\on{inv}_{\cN^-_P}:\Rep(\cP^-) \to \Rep(\cM).$$

Denote
$$\Omega(R)^{\on{spec}}:=(\on{Id}\otimes \on{inv}_{\cN^-_P})(R_\cG)\in  \on{FactAlg}^{\on{untl}}(\Rep(\cG)\otimes \Rep(\cM)).$$

\sssec{} 

Our current goal is to construct a (strictly unital) factorization functor 
$$\sF^{\on{spec}}:\on{I}(\cG,\cP^-)^{\on{spec,loc}}\to \Omega(R)^{\on{spec}}\mod^{\on{fact}}(\Rep(\cG)\otimes \Rep(\cM))$$

We first construct a (lax unital) factorization functor 
$$\on{pre-}\!\sF^{\on{spec}}:\on{I}(\cG,\cP^-)^{\on{spec,loc}}\to \Rep(\cG)\otimes \Rep(\cM).$$

\sssec{} \label{sss:pre F}

We define the functor $\on{pre-}\!\sF^{\on{spec}}$ as the $(\IndCoh,*)$-pushforward functor along the morphism
$$\LS^\reg_\cG\underset{\LS^\mer_\cG}\times \LS^\mer_{\cP^-}\underset{\LS^\mer_\cM}\times \LS^\reg_\cM\to 
\LS^\reg_\cG\times \LS^\reg_\cM.$$

Explicitly, in terms of \eqref{e:IGP spec defn}, the functor $\on{pre-}\!\sF^{\on{spec}}$ is given by
$$\IndCoh^*(\fL_\nabla(\cG/\cN^-_P))_{\fL^+_\nabla(\cG)\times \fL^+_\nabla(\cM)} \overset{\Gamma(\fL_\nabla(\cG/\cN^-_P),-)^\IndCoh}\longrightarrow
\Vect_{\fL^+_\nabla(\cG)\times \fL^+_\nabla(\cM)}\simeq
\Rep(\cG)\otimes \Rep(\cM).$$

\medskip

By construction, the functor $\on{pre-}\!\sF^{\on{spec}}$ is linear with respect to actions of
$$\Sph_\cG^{\on{spec}}\otimes \Sph_\cM^{\on{spec}}$$
on the two sides. 

\sssec{} \label{sss:F spec and contract}

Note that we can view the functor $\on{pre-}\!\sF^{\on{spec}}$ also as follows: it corresponds to the functor $J^{-,\on{spec,pre-enh}}$
of \eqref{e:J spec pre} via the self-duality of $\Rep(\cG)$. 

\medskip 

Consider the composition
\begin{equation} \label{e:spec and contract}
\Rep(\cG) \underset{\Sph^{\on{spec}}_\cG}\otimes \on{I}(\cG,\cP^-)^{\on{spec,loc}}
\overset{\on{Id}\otimes \on{pre-}\!\sF^{\on{spec}}}\longrightarrow \Rep(\cG) \underset{\Sph^{\on{spec}}_\cG}\otimes (\Rep(\cG)\otimes \Rep(\cM))\to \Rep(\cM),
\end{equation}
where the second arrow, i.e.,
\begin{equation} \label{e:spec and contract 1}
\Rep(\cG) \underset{\Sph^{\on{spec}}_\cG}\otimes (\Rep(\cG)\otimes \Rep(\cM))\to \Rep(\cM),
\end{equation} 
is induced by the duality pairing
$$\Rep(\cG)\underset{\Sph^{\on{spec}}_\cG}\otimes \Rep(\cG) \to \Vect.$$

By construction, the functor \eqref{e:spec and contract} identifies with the functor 
$J^{-,\on{spec,enh}}$. 

\begin{rem}  \label{r:pre F}

Let us view $\on{I}(\cG,\cP^-)^{\on{spec,loc}}$ as acted on by $\Rep(\cM)$ via pullback along
$$\LS^\reg_\cG\underset{\LS^\mer_\cG}\times \LS^\mer_{\cP^-}\underset{\LS^\mer_\cM}\times \LS^\reg_\cM\to \LS^\reg_\cM.$$

When we view $\on{pre-}\!\sF^{\on{spec}}$ as a $\Sph^{\on{spec}}_\cG\otimes \Rep(\cM)$-linear functor, it can be uniquely recovered
from the $\Sph^{\on{spec}}_\cG$-linear functor 
$$\on{pre-pre-}\!\sF^{\on{spec}}:\on{I}(\cG,\cP^-)^{\on{spec,loc}}\to \Rep(\cG),$$
equal to the composition of $\on{pre-}\!\sF^{\on{spec}}$ and 
$$(\on{Id}\otimes \on{inv}_\cM):\Rep(\cG)\otimes \Rep(\cM)\to \Rep(\cG).$$

\end{rem}

\sssec{}

Recall that the functor
$$\iota^\IndCoh_*:\Rep(\cP^-)\to \on{I}(\cG,\cP^-)^{\on{spec,loc}}$$
is strictly unital.

\medskip

The composite functor
$$\on{pre-}\!\sF^{\on{spec}}\circ \iota^\IndCoh_*:\Rep(\cP^-)\to \Rep(\cG)\otimes \Rep(\cM)$$
is the functor of direct image along
$$\LS^\reg_{\cP^-}\to \LS^\reg_{\cG}\times \LS^\reg_{\cM},$$
i.e., the functor of induction along
$$\cP^-\to \cG\times \cM.$$

In particular, we obtain that the functor $\on{pre-}\!\sF^{\on{spec}}$ sends the factorization unit to
$\Omega(R)^{\on{spec}}$. Hence, the functor $\on{pre-}\!\sF^{\on{spec}}$ upgrades to a functor 
$$\on{I}(\cG,\cP^-)^{\on{spec,loc}}\to \Omega(R)^{\on{spec}}\mod^{\on{fact}}(\Rep(\cG)\otimes \Rep(\cM)),$$
which is the sought-for functor $\sF^{\on{spec}}$.

\sssec{}

As in \secref{sss:action of Sph on Omega mod 1}, the action of $\Sph_\cG^{\on{spec}}\otimes \Sph_\cM^{\on{spec}}$ on 
$\Rep(\cG)\otimes \Rep(\cM)$ automatically extends to an action on 
$$\Omega(R)^{\on{spec}}\mod^{\on{fact}}(\Rep(\cG)\otimes \Rep(\cM)).$$

\medskip

Furthermore, as in \secref{sss:action of Sph on Omega mod 2}, 
the $\Sph_\cG^{\on{spec}}\otimes \Sph_\cM^{\on{spec}}$-linear structure on $\on{pre-}\!\sF^{\on{spec}}$
automatically extends to a $\Sph_\cG^{\on{spec}}\otimes \Sph_\cM^{\on{spec}}$-linear structure on the functor $\sF^{\on{spec}}$.

\sssec{} \label{sss:F spec and contract enh}

Note that the image of
$$\one_{\Rep(\cG)}\otimes \Omega(R)^{\on{spec}}\in \on{FactAlg}^{\on{untl}}\left(\Rep(\cG) \underset{\Sph^{\on{spec}}_\cG}\otimes (\Rep(\cG)\otimes \Rep(\cM))\right)$$
along the functor \eqref{e:spec and contract 1} identifies with
$$\Omega^{\on{spec}}\in \on{FactAlg}^{\on{untl}}(\Rep(\cM)).$$

Hence, the functor \eqref{e:spec and contract} upgrades to a functor
\begin{equation} \label{e:F spec and contract enh}
\Rep(\cG) \underset{\Sph^{\on{spec}}_\cG}\otimes \on{I}(\cG,\cP^-)^{\on{spec,loc}}\to \Omega^{\on{spec}}\mod^{\on{fact}}(\Rep(\cM)).
\end{equation}

Unwinding the definitions, we obtain that the functor \eqref{e:F spec and contract enh} identifies with the composition 
\begin{multline*} 
\Rep(\cG) \underset{\Sph^{\on{spec}}_\cG}\otimes \on{I}(\cG,\cP^-)^{\on{spec,loc}} \overset{\text{\eqref{e:amb vs IGP spec *}}}\longrightarrow \\
\to \IndCoh^*((\LS^\Mmf_{\cP^-})^\wedge_\mf) \overset{\sQ^\wedge_\mf}\longrightarrow \Omega^{\on{spec}}\mod^{\on{fact}}(\Rep(\cM)).
\end{multline*}

\ssec{Fully faithfulness on the eventually coconnective part}

\sssec{}

Recall that the factorization category $\on{I}(\cG,\cP^-)^{\on{spec,loc}}$ carries a naturally defined t-structure,
\secref{sss:t IGP spec}

\medskip

The goal of this subsection is to prove the following key assertion:

\begin{prop} \label{p:F spec ff on pos}
The functor $\sF^{\on{spec}}$ is fully faithful when restricted to the eventually coconnective subcategory of 
$\on{I}(\cG,\cP^-)^{\on{spec,loc}}$, i.e., for every affine scheme $S$, the corresponding functor
$$(\on{I}(\cG,\cP^-)^{\on{spec,loc}}_S)^{>-\infty}\to (\Omega(R)^{\on{spec}}\mod^{\on{fact}}(\Rep(\cG)\otimes \Rep(\cM))_S)^{>-\infty}$$
is fully faithful. 
\end{prop}

The rest of this subsection is devoted to the proof of this proposition.

\sssec{}

Recall the functor
\begin{equation} \label{e:Omega R}
\Sph_\cG^{\on{spec}}\to R_\cG\mod^{\on{fact}}(\Rep(\cG)\otimes \Rep(\cG))
\end{equation}
of \cite[Equation (E.30)]{GLC2}, and recall that it induces an equivalence between the eventually coconnective subcategories
of the two sides (by \cite[Proposition E.9.4]{GLC2}). 

\medskip

Applying to \eqref{e:Omega R} the tensoring-up procedure $-\underset{\Rep(\cG)}\otimes \Rep(\cP^-)$, we obtain a functor
\begin{equation} \label{e:Omega R P}
\Sph_\cG^{\on{spec}}\underset{\Rep(\cG)}\otimes \Rep(\cP^-) \to R_\cG\mod^{\on{fact}}(\Rep(\cG)\otimes \Rep(\cP^-)),
\end{equation}
which also induces an equivalence between the eventually coconnective subcategories
of the two sides, which follows formally from the corresponding property of \eqref{e:Omega R}
(or, alternatively, by the same argument as in \cite[Proposition E.9.4]{GLC2}). 

\sssec{}

We will think of $\Rep(\cG)\underset{\Sph_\cG^{\on{spec}}}\otimes \Rep(\cP^-)$ as
$$\IndCoh^*(\LS^\reg_\cG\underset{\LS^\mer_\cG}\times \LS^\reg_{\cP^-}).$$

Denote by $\iota$ the map
$$\LS^\reg_{\cP^-}\to \LS^\mer_{\cP^-}\underset{\LS^\mer_\cM}\times \LS^\reg_\cM,$$
and also the base-changed map
$$\LS^\reg_\cG\underset{\LS^\mer_\cG}\times \LS^\reg_{\cP^-}\to
\LS^\reg_\cG\underset{\LS^\mer_\cG}\times \LS^\mer_{\cP^-}\underset{\LS^\mer_\cM}\times \LS^\reg_\cM.$$

Consider the resulting functor
\begin{equation} \label{e:H I GP spec}
\iota^\IndCoh_*:\IndCoh^*(\LS^\reg_\cG\underset{\LS^\mer_\cG}\times \LS^\reg_{\cP^-}) \to
\IndCoh^*(\LS^\reg_\cG\underset{\LS^\mer_\cG}\times \LS^\mer_{\cP^-}\underset{\LS^\mer_\cM}\times \LS^\reg_\cM).
\end{equation} 

\medskip

Similar to \secref{sss:semiinf spec oblv}, one shows that the functor $\iota^\IndCoh_*$ of \eqref{e:H I GP spec}
admits a continuous right adjoint, giving rise to a monadic adjunction $(\iota^\IndCoh_*,\iota^!)$.

\sssec{}

Note that by construction, the following diagram commutes: 
\begin{equation} \label{e:iota diag}
\CD
\Sph_\cG^{\on{spec}}\underset{\Rep(\cG)}\otimes \Rep(\cP^-) @>{\iota^\IndCoh_*}>> \on{I}(\cG,\cP^-)^{\on{spec,loc}} \\
@V{\text{\eqref{e:Omega R P}}}VV @VV{\sF^{\on{spec}}}V \\
R_\cG\mod^{\on{fact}}(\Rep(\cG)\otimes \Rep(\cP)) @>{\on{Id}\otimes \on{inv}_{\cN^-_P}}>> 
\Omega(R)^{\on{spec}}\mod^{\on{fact}}(\Rep(\cG)\otimes \Rep(\cM)).
\endCD
\end{equation} 

We will prove:

\begin{prop}   \label{p:iota diag} \hfill 

\smallskip

\noindent{\em(a)}
The bottom horizontal arrow in \eqref{e:iota diag} admits a continuous right adjoint, so that the resulting monad on
$R_\cG\mod^{\on{fact}}(\Rep(\cG)\otimes \Rep(\cP))$ is left t-exact. 

\smallskip

\noindent{\em(b)}
The natural transformation in the diagram
\begin{equation} \label{e:iota diag adj}
\vcenter
{\xy
(0,0)*+{\Sph_\cG^{\on{spec}}\underset{\Rep(\cG)}\otimes \Rep(\cP^-) }="A";
(75,0)*+{\on{I}(\cG,\cP^-)^{\on{spec,loc}}}="B";
(0,-30)*+{R_\cG\mod^{\on{fact}}(\Rep(\cG)\otimes \Rep(\cP)) }="C";
(75,-30)*+{\Omega(R)^{\on{spec}}\mod^{\on{fact}}(\Rep(\cG)\otimes \Rep(\cM)),}="D";
{\ar@{->}_{\text{\eqref{e:Omega R P}}} "A";"C"};
{\ar@{->}_{\iota^!} "B";"A"};
{\ar@{->}^{\sF^{\on{spec}}} "B";"D"};
{\ar@{->} "D";"C"};
{\ar@{=>} "A";"D"};
\endxy}
\end{equation}
obtained from \eqref{e:iota diag} by passing to right adjoints along the horizontal arrows, is an isomorphism. 
\end{prop}

\sssec{}

Before we prove \propref{p:iota diag}, let us show how it implies the statement of \propref{p:F spec ff on pos}. 

\medskip

Point (a) of \propref{p:iota diag} implies that the category $\Omega(R)^{\on{spec}}\mod^{\on{fact}}(\Rep(\cG)\otimes \Rep(\cM))$
carries a t-structure, uniquely characterized by the fact that the functor
$$\on{Id}\otimes \on{inv}_{\cN^-_P}:R_\cG\mod^{\on{fact}}(\Rep(\cG)\otimes \Rep(\cP))\to 
\Omega(R)^{\on{spec}}\mod^{\on{fact}}(\Rep(\cG)\otimes \Rep(\cM))$$
is t-exact.

\medskip

We will show that the functor $\sF^{\on{spec}}$ induces an equivalence between the eventually coconnective subcategories 
of the two sides. 

\medskip

The functor $\iota^\IndCoh_*$ is t-exact. Hence, it identifies $(\on{I}(\cG,\cP^-)^{\on{spec,loc}})^{>-\infty}$ with
$$(\iota^!\circ \iota^\IndCoh_*)\mod\left(\left(\Sph_\cG^{\on{spec}}\underset{\Rep(\cG)}\otimes \Rep(\cP^-)\right)^{>-\infty}\right).$$ 

Similarly, the functor $\on{Id}\otimes \on{inv}_{\cN^-_P}$ identifies $\left(\Omega(R)^{\on{spec}}\mod^{\on{fact}}(\Rep(\cG)\otimes \Rep(\cM))\right)^{>-\infty}$
with the category of modules over the corresponding monad acting on $\left(R_\cG\mod^{\on{fact}}(\Rep(\cG)\otimes \Rep(\cP))\right)^{>-\infty}$.

\medskip

Now, point (b) of \propref{p:iota diag} implies that the functor \eqref{e:Omega R P} intertwines the two monads, i.e., the equivalence 
of \eqref{e:Omega R P} extends to an equivalence between the corresponding categories of modules.

\qed[\propref{p:F spec ff on pos}]

\ssec{Proof of \propref{p:iota diag}}

\sssec{}  \label{sss:Phi adj}

Let $\Phi:\bC\to \bD$ be a (strictly unital) functor between unital lax factorization categories.
Assume that $\Phi$ admits a right adjoint.\footnote{Of course, this right adjoint is supposed to be continuous when evaluated on every $S\to \Ran$.} 
Then this right adjoint, to be denoted $\Phi^R$,
carries a natural (lax unital) factorization structure (see \cite[Sect. C.11.21]{GLC2}).

\medskip

Let $\CA$ be a (unital) factorization algebra in $\bC$. Then functor $\Phi$ induces a (strictly unital)
factorization functor between lax factorization categories
\begin{equation} \label{e:Phi ind funct}
\Phi_\CA:\CA\mod^{\on{fact}}(\bC)\to \Phi(\CA)\mod^{\on{fact}}(\bD),
\end{equation}
which makes the diagram
\begin{equation} \label{e:Phi ind funct oblv}
\CD
\CA\mod^{\on{fact}}(\bC) @>{\Phi_\CA}>> \Phi(\CA)\mod^{\on{fact}}(\bD) \\
@V{\oblv_\CA}VV @VV{\oblv_{\Phi(\CA)}}V \\
\bC @>>{\Phi}> \bD
\endCD
\end{equation}
commute.

\medskip

We claim that the functor $\Phi_\CA$ of \eqref{e:Phi ind funct} admits a continuous right adjoint, and that the natural transformation 
in the diagram
$$
\xy
(0,0)*+{\CA\mod^{\on{fact}}(\bC)}="A";
(40,0)*+{\Phi(\CA)\mod^{\on{fact}}(\bD)}="B";
(0,-15)*+{\bC}="C";
(40,-15)*+{\bD,}="D";
{\ar@{->}_{\oblv_\CA} "A";"C"};
{\ar@{->}_{(\Phi_\CA)^R} "B";"A"};
{\ar@{->}^{\oblv_{\Phi(\CA)}} "B";"D"};
{\ar@{->}^{\Phi^R} "D";"C"};
{\ar@{=>} "A";"D"};
\endxy
$$
obtained from \eqref{e:Phi ind funct oblv} by passing to right adjoints along the horizontal arrows, is an isomorphism. 

\medskip

Indeed, it is easy to see that the right adjoint $(\Phi_\CA)^R$ of $\Phi_\CA$ is the composition of the functor 
$$(\Phi^R)_{\Phi(\CA)}:\Phi(\CA)\mod^{\on{fact}}(\bD)\to \Phi^R(\Phi(\CA))\mod^{\on{fact}}(\bC),$$
followed by the functor of restriction
$$\Phi^R(\Phi(\CA))\mod^{\on{fact}}(\bC)\to \CA\mod^{\on{fact}}(\bC)$$
along the homomorphism of factorization algebras
$$\CA\to \Phi^R(\Phi(\CA)).$$

\sssec{} \label{sss:alg adj}

We will apply the above paradigm as follows. We will start with a (lax unital) factorization functor
between unital factorization categories:
$$\Phi':\bC\to \bD',$$
and we will set $\bD:=\Phi'(\one_\bC)\mod^{\on{fact}}(\bD')$. Let $\Phi$ denote the resulting functor $\bC\to \bD$. 

\medskip

Note that for any $\CA\in \on{FactAlg}^{\on{untl}}(\bC)$, we have 
$$\Phi(\CA)\mod^{\on{fact}}(\bD)\simeq \Phi'(\CA)\mod^{\on{fact}}(\bD').$$

\sssec{}

We take
$$\bC:=\Rep(\cG)\otimes \Rep(\cP^-),\,\, \bD':=\Rep(\cG)\otimes \Rep(\cM),\,\, \Phi':=\on{Id}\otimes \on{inv}_{\cN^-_P} \text{ and } \CA:=R_\cG,$$
so that 
$$\Phi'(\one_\bC)=\one_{\Rep(\cG)}\otimes \Omega^{\on{spec}}.$$

It follows from \thmref{t:Omega comp} that the functor $\on{Id}\otimes \on{inv}_{\cN^-_P}$ preserves compactness, and hence
admits a right adjoint. 

\medskip

This implies the existence of a continuous right adjoint in point (a) of \propref{p:iota diag}. 

\sssec{}

We now show that the monad in point (a) of \propref{p:iota diag} is left t-exact. By \secref{sss:Phi adj}, it suffices to show that the monad
corresponding to the functor
$$\on{Id}\otimes \on{inv}_{\cN^-_P}:\Rep(\cG)\otimes \Rep(\cP^-)\to (\one_{\Rep(\cG)}\otimes \Omega^{\on{spec}})\mod^{\on{fact}}(\Rep(\cG)\otimes \Rep(\cM))$$
is left t-exact. 

\medskip

It suffices to show that the monad corresponding to the functor
$$\on{inv}_{\cN^-_P}:\Rep(\cP^-)\to \Omega^{\on{spec}}\mod^{\on{fact}}(\Rep(\cM))$$
is left t-exact. 

\medskip

However, by \thmref{t:LS N vs Omega bis}, the latter adjunction identifies with
$$({}'\!\iota)^\IndCoh_*:\IndCoh^*(\LS^\reg_{\cP^-})\rightleftarrows \IndCoh^*((\LS^\Mmf_{\cP^-})^\wedge_\mf):({}'\!\iota)^!,$$
and $({}'\!\iota)^!\circ ({}'\!\iota)^\IndCoh_*$ is left t-exact by \propref{p:t semiinf spec}. 

\sssec{}

We now prove point (b) of \propref{p:iota diag}. By \secref{sss:Phi adj}, it suffices to show that the natural transformation in the diagram
\begin{equation} \label{e:iota diag 1 adj}
\xy
(0,0)*+{\Sph_\cG^{\on{spec}}\underset{\Rep(\cG)}\otimes \Rep(\cP^-) }="A";
(65,0)*+{\on{I}(\cG,\cP^-)^{\on{spec,loc}}}="B";
(0,-20)*+{\Rep(\cG)\otimes \Rep(\cP^-)}="C";
(65,-20)*+{(\one_{\Rep(\cG)}\otimes \Omega^{\on{spec}})\mod^{\on{fact}}(\Rep(\cG)\otimes \Rep(\cM)),}="D";
{\ar@{->} "A";"C"};
{\ar@{->}_{\iota^!} "B";"A"};
{\ar@{->}^{\on{pre-}\!\sF^{\on{spec}}} "B";"D"};
{\ar@{->} "D";"C"};
{\ar@{=>} "A";"D"};
\endxy
\end{equation} 
obtained by passing to right adjoints along the horizontal arrows in the commutative diagram 
\begin{equation} \label{e:iota diag 1}
\xy
(0,0)*+{\Sph_\cG^{\on{spec}}\underset{\Rep(\cG)}\otimes \Rep(\cP^-) }="A";
(65,0)*+{\on{I}(\cG,\cP^-)^{\on{spec,loc}}}="B";
(0,-20)*+{\Rep(\cG)\otimes \Rep(\cP^-)}="C";
(65,-20)*+{(\one_{\Rep(\cG)}\otimes \Omega^{\on{spec}})\mod^{\on{fact}}(\Rep(\cG)\otimes \Rep(\cM)),}="D";
{\ar@{->} "A";"C"};
{\ar@{->}^{\iota^\IndCoh_*} "A";"B"};
{\ar@{->}^{\on{pre-}\!\sF^{\on{spec}}} "B";"D"};
{\ar@{->} "C";"D"};
\endxy
\end{equation} 
is an isomorphism. 

\medskip

Still equivalently, it suffices to show that the natural transformation in the diagram
\begin{equation} \label{e:iota diag 2 adj}
\xy
(0,0)*+{\Sph_\cG^{\on{spec}}\underset{\Rep(\cG)}\otimes \Rep(\cP^-) }="A";
(65,0)*+{\on{I}(\cG,\cP^-)^{\on{spec,loc}}}="B";
(0,-20)*+{\Rep(\cP^-)}="C";
(65,-20)*+{\Omega^{\on{spec}}\mod^{\on{fact}}(\Rep(\cM)),}="D";
{\ar@{->} "A";"C"};
{\ar@{->}_{\iota^!} "B";"A"};
{\ar@{->} "B";"D"};
{\ar@{->} "D";"C"};
{\ar@{=>} "A";"D"};
\endxy
\end{equation} 
obtained by passing to right adjoints along the horizontal arrows in the commutative diagram 
\begin{equation} \label{e:iota diag 2}
\xy
(0,0)*+{\Sph_\cG^{\on{spec}}\underset{\Rep(\cG)}\otimes \Rep(\cP^-) }="A";
(65,0)*+{\on{I}(\cG,\cP^-)^{\on{spec,loc}}}="B";
(0,-20)*+{\Rep(\cP^-)}="C";
(65,-20)*+{\Omega^{\on{spec}}\mod^{\on{fact}}(\Rep(\cM)),}="D";
{\ar@{->} "A";"C"};
{\ar@{->}^{\iota^\IndCoh_*} "A";"B"};
{\ar@{->}"B";"D"};
{\ar@{->} "C";"D"};
\endxy
\end{equation} 
is an isomorphism. 

\sssec{}

Using \thmref{t:LS N vs Omega bis}, we rewrite diagram \eqref{e:iota diag 2} as 
$$
\xy
(0,0)*+{\IndCoh^*(\LS^\reg_\cG\underset{\LS^\mer_\cG}\times \LS^\reg_{\cP^-})}="A";
(75,0)*+{\IndCoh^*(\LS^\reg_\cG\underset{\LS^\mer_\cG}\times \LS^\mer_{\cP^-}\underset{\LS^\mer_\cM}\times \LS^\reg_\cM)}="B";
(0,-20)*+{\IndCoh^*(\LS^\reg_{\cP^-})}="C";
(75,-20)*+{\IndCoh^*((\LS^\Mmf_{\cP^-})^\wedge_\mf),}="D";
{\ar@{->} "A";"C"};
{\ar@{->}^{\iota^\IndCoh_*} "A";"B"};
{\ar@{->}"B";"D"};
{\ar@{->}^{({}'\!\iota)^\IndCoh_*} "C";"D"};
\endxy
$$
in which the vertical arrows are the $(\IndCoh,*)$-pushforward functors. The required assertion follows now by base change. 

\ssec{The functor \texorpdfstring{$\sF^{\on{geom}}$}{Fgeom}}

In this subsection we will map the geometric side in \thmref{t:semiinf geom Satake} to the category of modules
over a factorization algebra in $\Whit_*(G)\otimes \Whit^!(M)$.

\sssec{} 

Recall the functors
$$J_{\Whit}^{-,\on{pre-enh}}:\Whit^!(G)\otimes \on{I}(G,P^-)^{\on{loc}}_{\rho_P(\omega_X)}\to \Whit^!(M)$$ and 
$$J_{\Whit}^{-,\on{enh}}:\Whit^!(G)\underset{\Sph_G}\otimes \on{I}(G,P^-)^{\on{loc}}_{\rho_P(\omega_X)}\to \Whit^!(M).$$
from \eqref{e:Whit semiinf pair}.

\medskip

By duality, the functor $J_{\Whit}^{-,\on{pre-enh}}$ gives rise to a functor
\begin{equation} \label{e:F geom}
\on{I}(G,P^-)^{\on{loc}}_{\rho_P(\omega_X)}\to \Whit_*(G)\otimes \Whit^!(M),
\end{equation} 
which we will denote by $\on{pre-}\!\sF^{\on{geom}}$.

\medskip

By construction, the functor $\on{pre-}\!\sF^{\on{geom}}$ is linear with respect to the actions of $\Sph_G\otimes \Sph_M$
on the two sides. 

\sssec{} \label{sss:F geom and contract}

Consider the composition
\begin{multline} \label{e:F geom and contract}
\Whit^!(G) \underset{\Sph_G}\otimes \on{I}(G,P^-)^{\on{loc}}_{\rho_P(\omega_X)} \overset{\on{Id}\otimes \on{pre-}\!\sF^{\on{geom}}}\longrightarrow \\
\to \Whit^!(G) \underset{\Sph_G}\otimes (\Whit_*(G)\otimes \Whit^!(M))  \to \Whit^!(M),
\end{multline}
where the second arrow, i.e.,
\begin{equation} \label{e:geom and contract 1}
\Whit^!(G) \underset{\Sph_G}\otimes (\Whit_*(G)\otimes \Whit^!(M))  \to \Whit^!(M)
\end{equation}
is induced by the duality pairing
$$\Whit_!(G)\underset{\Sph_G}\otimes \Whit_*(G) \to \Vect.$$
By construction, the functor \eqref{e:F geom and contract} identifies with the functor $J_{\Whit}^{-,\on{enh}}$.

\sssec{} \label{sss:prepre F geom}

Let us view the two sides of \eqref{e:F geom} as acted on by $\Sph_G\otimes \Rep(\cM)$ via
$$\on{Sat}_M^{\on{nv}}:\Rep(\cM)\to \Sph_M.$$

\medskip

The functor $\on{pre-}\!\sF^{\on{geom}}$ is recovered uniquely from the $\Sph_G$-linear functor
$$\on{pre-pre-}\!\sF^{\on{geom}}:\on{I}(G,P^-)^{\on{loc}}_{\rho_P(\omega_X)}\to \Whit_*(G),$$
equal to the composition of $\on{pre-}\!\sF^{\on{geom}}$ followed by the functor
$$\Whit^!(M) \to \Vect,$$
given by taking the !-fiber at the unit point of the affine Grassmannian.

\medskip

Note that the functor $\on{pre-pre-}\!\sF^{\on{geom}}$ equals the composition
\begin{equation} \label{e:prepre F geom}
\on{I}(G,P^-)^{\on{loc}}_{\rho_P(\omega_X)} \overset{\alpha_{\rho_M(\omega_X),\on{taut}}}\longrightarrow  \on{I}(G,P^-)^{\on{loc}}_{\rho(\omega_X)} \to 
\Dmod_{\frac{1}{2}}(\Gr_{G,\rho(\omega_X)})\to \Whit_*(G).
\end{equation} 

\sssec{}

Let 
$$\Omega(R)\in \on{FactAlg}^{\on{untl}}(\Whit_*(G)\otimes \Whit^!(M))$$
be the image of the factorization unit in $\on{I}(G,P^-)^{\on{loc}}_{\rho_P(\omega_X)}$ along the functor $\on{pre-}\!\sF^{\on{geom}}$.

\medskip

Thus, we obtain that the functor $\on{pre-}\!\sF^{\on{geom}}$ lifts to a functor
$$\on{I}(G,P^-)^{\on{loc}}_{\rho_P(\omega_X)}\to \Omega(R)\mod^{\on{fact}}(\Whit_*(G)\otimes \Whit^!(M)),$$
which we denote by $\sF^{\on{geom}}$.

\sssec{}

As in \secref{sss:action of Sph on Omega mod 1}, the $\Sph_G\otimes \Sph_M$-action on $\Whit_*(G)\otimes \Whit^!(M)$ canonically lifts to an action of
$\Sph_G\otimes \Sph_M$ on $\Omega(R)\mod^{\on{fact}}(\Whit_*(G)\otimes \Whit^!(M))$. 
 
\medskip

Furthermore,  as in \secref{sss:action of Sph on Omega mod 2}, 
the $\Sph_G\otimes \Sph_M$-linear structure on $\on{pre-}\!\sF^{\on{geom}}$ canonically lifts to a 
$\Sph_G\otimes \Sph_M$-linear structure on $\sF^{\on{geom}}$.

\sssec{} \label{sss:F geom and contract enh}

Note that the image of
$$\one_{\Whit^!(G)}\otimes \Omega(R)\in \on{FactAlg}^{\on{untl}}\left(\Whit^!(G) \underset{\Sph_G}\otimes (\Whit_*(G)\otimes \Whit^!(M))\right)$$
along the functor \eqref{e:geom and contract 1} identifies canonically with 
$$\Omega\in  \on{FactAlg}^{\on{untl}}(\Whit^!(M)).$$

Hence, the functor \eqref{e:F geom and contract} upgrades to a functor
\begin{equation} \label{e:F geom and contract enh}
\Whit^!(G) \underset{\Sph_G}\otimes \on{I}(G,P^-)^{\on{loc}}_{\rho_P(\omega_X)} \to
\Omega\mod^{\on{fact}}(\Whit^!(M)).
\end{equation}

Unwinding the construction, we obtain that the functor \eqref{e:F geom and contract enh} identifies with the functor
$$\Whit^!(G) \underset{\Sph_G}\otimes \on{I}(G,P^-)^{\on{loc}}_{\rho_P(\omega_X)} \simeq \Whit^!(\fL(G))^{-,\semiinfAccs}
\overset{J_{\Whit,\Omega}^{-,\semiinf}}\longrightarrow \Omega\mod^{\on{fact}}(\Whit^!(M)).$$

\ssec{Connecting the geometric and spectral sides}

We now come to the key point in the proof of \thmref{t:semiinf geom Satake}.
 
\sssec{}

Denote
$$\FLE_{\cG,\infty,\tau}:=\tau_G\circ \FLE_{\cG,\infty}, \quad \Rep(\cG)\to \Whit_*(G).$$

Consider the functor
$$\FLE_{\cG,\infty,\tau}^{-1}\otimes \on{CS}_{M,\tau}:\Whit_*(G)\otimes \Whit^!(M)\to \Rep(\cG)\otimes \Rep(\cM).$$

\medskip

From \thmref{t:ident Omega bis} we obtain:

\begin{cor} \label{c:Omega to Omega spec}
The functor $\FLE_{\cG,\infty,\tau}^{-1}\otimes \on{CS}_{M,\tau}$ sends
$$\Omega(R)\mapsto \Omega(R)^{\on{spec}}.$$
\end{cor}

\sssec{}  \label{sss:prepre F}

Set
$$\on{pre-pre-}\!\sF:=\FLE_{\cG,\infty,\tau}^{-1}\circ \on{pre-pre-}\!\sF^{\on{geom}}, \quad
\on{I}(G,P^-)^{\on{loc}}_{\rho_P(\omega_X)}\to \Rep(\cG),$$
$$\on{pre-}\!\sF:=(\FLE_{\cG,\infty,\tau}^{-1}\otimes \on{CS}_{M,\tau})\circ \on{pre-}\!\sF^{\on{geom}}, \quad
\on{I}(G,P^-)^{\on{loc}}_{\rho_P(\omega_X)}\to \Rep(\cG)\otimes \Rep(\cM).$$

\medskip

From \corref{c:Omega to Omega spec} we obtain that the functor $\on{pre-}\!\sF$ upgrades to a functor
$$\sF:\on{I}(G,P^-)^{\on{loc}}_{\rho_P(\omega_X)}\to \Omega(R)^{\on{spec}}\mod^{\on{fact}}(\Rep(\cG)\otimes \Rep(\cM)).$$

\sssec{}

The functor $\sF$ is $\Sph^{\on{spec}}_\cG\otimes \Sph^{\on{spec}}_\cM$-linear, where the action on $\on{I}(G,P^-)^{\on{loc}}_{\rho_P(\omega_X)}$
is given by
$$\Sph^{\on{spec}}_\cG\otimes \Sph^{\on{spec}}_\cM \overset{\on{Sat}^{-1}_G\otimes \Sat^{-1}_{M,\tau}}\longrightarrow \Sph_G\otimes \Sph_M.$$

\sssec{} \label{sss:F and contract}

Consider the composition
\begin{equation} \label{e:F and contract}
\Whit^!(G)  \underset{\Sph_G}\otimes\on{I}(G,P^-)^{\on{loc}}_{\rho_P(\omega_X)} \overset{\on{CS}_{G,\tau}\otimes \on{pre-}\!\sF}\longrightarrow
\Rep(\cG) \underset{\Sph^{\on{spec}}_\cG}\otimes  (\Rep(\cG)\otimes \Rep(\cM)) \to \Rep(\cM),
\end{equation}
where the second arrow is as in \eqref{e:spec and contract 1}.

\medskip

From \secref{sss:F geom and contract} we obtain that the functor \eqref{e:F and contract} identifies with
$$\Sat_{M,\tau}\circ J_{\Whit}^{-,\on{enh}}.$$

%
%

\ssec{Construction of the functor}

In this subsection we will construct the functor $\Sat^{-,\semiinf}$. We will reduce the assertion that it is an equivalence to
a pointwise statement. 

\sssec{}

Recall that the category $\on{I}(G,P^-)^{\on{loc}}_{\rho_P(\omega_X)}$ is compactly generated\footnote{I.e., this holds when evaluated
on any affine scheme $S$ mapping to $\Ran$.} by objects of the form
$$\CF\underset{M}\star \one_{\on{I}(G,P^-)^{\on{loc}}_{\rho_P(\omega_X)}}, \quad \CF\in (\Sph_M)^c.$$

\medskip

Similarly, the category $\on{I}(\cG,\cP^-)^{\on{spec,loc}}$ is compactly generated by objects of the form
$$\CF\underset{\cM}\star \one_{\on{I}(\cG,\cP^-)^{\on{spec,loc}}}, \quad \CF\in (\Sph^{\on{spec}}_\cM)^c.$$

Moreover, the latter objects are eventually coconnective. In particular, we obtain that 
$$(\on{I}(\cG,\cP^-)^{\on{spec,loc}})^c\subset (\on{I}(\cG,\cP^-)^{\on{spec,loc}})^{>-\infty}.$$

\sssec{} \label{sss:IGP via Omega} 

Since the functors $\sF$ and $\sF^{\on{spec}}$ are $\Sph^{\on{spec}}_\cM$-linear, we obtain that the essential image
under $\sF$ of compact objects of $\on{I}(G,P^-)^{\on{loc}}_{\rho_P(\omega_X)}$ is contained in the essential
image under $\sF^{\on{spec}}$ of $(\on{I}(\cG,\cP^-)^{\on{spec,loc}})^{>-\infty}$.

\medskip

Using \propref{p:F spec ff on pos}, we obtain that the functor $\sF$ canonically lifts to a functor
$$\Sat^{-,\semiinf}:\on{I}(G,P^-)^{\on{loc}}_{\rho_P(\omega_X)}\to \on{I}(\cG,\cP^-)^{\on{spec,loc}},$$
which makes the diagram
\begin{equation} \label{e:Omega triang orgn}
\vcenter
{\xy
(0,0)*+{\on{I}(G,P^-)^{\on{loc}}_{\rho_P(\omega_X)}}="A";
(75,0)*+{\on{I}(\cG,\cP^-)^{\on{spec,loc}}}="B";
(40,-20)*+{\Omega(R)^{\on{spec}}\mod^{\on{fact}}(\Rep(\cG)\otimes \Rep(\cM))}="C";
{\ar@{->}_{\sF} "A";"C"};
{\ar@{->}^{\sF^{\on{spec}}} "B";"C"};
{\ar@{->}^{\Sat^{-,\semiinf}} "A";"B"};
\endxy}
\end{equation}
commute. 

\medskip

Furthermore, by construction, the functor $\Sat^{-,\semiinf}$ preserves compactness.

\sssec{}

The linearity of the functors $\sF$ and $\sF^{\on{spec}}$ with respect to
$$(\Sph^{\on{spec}}_\cG\otimes \Sph^{\on{spec}}_\cM)^c$$
implies that the functor $\Sat^{-,\semiinf}$ is also linear with respect to 
$(\Sph^{\on{spec}}_\cG\otimes \Sph^{\on{spec}}_\cM)^c$ on compacts.

\medskip

Hence, we obtain that the functor $\Sat^{-,\semiinf}$ is linear with respect to
$\Sph^{\on{spec}}_\cG\otimes \Sph^{\on{spec}}_\cM$.

\sssec{}

Combining Sects. \ref{sss:F spec and contract} and \ref{sss:F and contract}, we obtain that the functor
$\Sat^{-,\semiinf}$ makes the diagram \eqref{e:semiinf geom Satake} commute. 

\sssec{} 

Furthermore, it follows from Sects. \ref{sss:F geom and contract enh} and \ref{sss:F spec and contract enh} that we have a commutative diagram
$$
\CD
\Whit^!(G) \underset{\Sph_G}\otimes \on{I}(G,P^-)^{\on{loc}}_{\rho_P(\omega_X)} @>{\on{CS}_{G,\tau}\otimes \Sat^{-,\semiinf}}>> 
\Rep(\cG) \underset{\Sph^{\on{spec}}_\cG}\otimes \on{I}(\cG,\cP^-)^{\on{spec,loc}} \\
@V{\sim}VV @V{\sim}V{\text{\eqref{e:amb vs IGP spec *}}}V  \\
\Whit^!(\fL(G))^{-,\semiinfAccs} & &   \IndCoh^*((\LS^\Mmf_{\cP^-})^\wedge_\mf) \\
@V{J_{\Whit,\Omega}^{-,\semiinf}}VV @VV{\sQ^\wedge_\mf}V \\
\Omega\mod^{\on{fact}}(\Whit^!(M)) @>{\Sat_{M,\tau}}>> \Omega^{\on{spec}}\mod^{\on{fact}}(\Rep(\cM)).
\endCD
$$

Combining with Theorems \ref{t:Whit vs Omega bis} and \ref{t:LS N vs Omega bis}, we obtain a commutative diagram
\begin{equation} \label{e:Sat vs CS semiinf}
\CD
\Whit^!(G) \underset{\Sph_G}\otimes \on{I}(G,P^-)^{\on{loc}}_{\rho_P(\omega_X)} @>{\on{CS}_{G,\tau}\otimes \Sat^{-,\semiinf}}>> 
\Rep(\cG) \underset{\Sph^{\on{spec}}_\cG}\otimes \on{I}(\cG,\cP^-)^{\on{spec,loc}} \\
@V{\sim}VV @V{\text{\eqref{e:amb vs IGP spec *}}}V{\sim}V  \\
\Whit^!(\fL(G))^{-,\semiinfAccs} @>{\on{CS}^{-,\semiinf}_\tau}>>  \IndCoh^*((\LS^\Mmf_{\cP^-})^\wedge_\mf). 
\endCD
\end{equation} 

\sssec{}

In order to prove \thmref{t:semiinf geom Satake}, it remains to show that the functor $\Sat^{-,\semiinf}$
is an equivalence. 

\medskip

Since it preserves compactness, by \cite[Proposition 6.2.6]{GLC2}, it suffices to show that $\Sat^{-,\semiinf}$
induces an equivalence 
\begin{equation} \label{e:Sat semiinf ptw}
\on{I}(G,P^-)^{\on{loc}}_{\rho_P(\omega_X),\ul{x}}\to \on{I}(\cG,\cP^-)_{\ul{x}}^{\on{spec,loc}}
\end{equation} 
for every geometric point $\ul{x}\in \Ran$.

\medskip

We will prove this in the next subsection. 

\ssec{Proof of the pointwise equivalence}

\sssec{}

The functor \eqref{e:Sat semiinf ptw} preserves compactness, and is essentially surjective on compact generators.
Hence, it suffices to show that it is fully faithful on compacts.

\sssec{}

Recall the paradigm of \cite[Sect. 7.1]{GLC2}. Let
$$(\on{I}(G,P^-)^{\on{loc}}_{\rho_P(\omega_X),\ul{x}})_{\on{temp}} \text{ and }
(\on{I}(\cG,\cP^-)_{\ul{x}}^{\on{spec,loc}})_{\on{temp}}$$
denote the temprered localizations of the two sides, where the temperedness is considered 
with respect to $G$ (and not $M$). 

\medskip

First, we claim that the functor
\begin{equation} \label{e:Sat semiinf temp}
(\on{I}(G,P^-)^{\on{loc}}_{\rho_P(\omega_X),\ul{x}})_{\on{temp}} \to (\on{I}(\cG,\cP^-)_{\ul{x}}^{\on{spec,loc}})_{\on{temp}},
\end{equation}
induced by \eqref{e:Sat semiinf ptw}, is an equivalence. 

\medskip

In order to prove this, by \cite[Corollary 7.1.12]{GLC2}, it suffices to show that the induced functor
\begin{equation} \label{e:Sat semiinf tensor up}
\Rep(\cG)\underset{\Sph^{\on{spec}}_\cG}\otimes \on{I}(G,P^-)^{\on{loc}}_{\rho_P(\omega_X),\ul{x}}\to
\Rep(\cG)\underset{\Sph^{\on{spec}}_\cG}\otimes \on{I}(\cG,\cP^-)_{\ul{x}}^{\on{spec,loc}}
\end{equation} 
is an equivalence. 

\medskip

However, by \eqref{e:Sat vs CS semiinf}, the functor \eqref{e:Sat semiinf tensor up} identifies with the functor
$$\on{CS}^{-,\semiinf}_\tau:\Whit^!(\fL(G))_{\ul{x}}^{-,\semiinfAccs}\to \IndCoh^*((\LS^\Mmf_{\cP^-})^\wedge_\mf)_{\ul{x}},$$
and hence is an equivalence by (the pointwise version of) \thmref{t:semiinf CS}. 

\sssec{}

Given that \eqref{e:Sat semiinf temp} is an equivalence, in order to prove that \eqref{e:Sat semiinf ptw} is fully faithful on compacts,
it suffices to show that both functors
\begin{equation} \label{e:IGP temp proj}
\on{I}(G,P^-)^{\on{loc}}_{\rho_P(\omega_X),\ul{x}}\to (\on{I}(G,P^-)^{\on{loc}}_{\rho_P(\omega_X),\ul{x}})_{\on{temp}}
\end{equation}
and 
\begin{equation} \label{e:IGP spec temp proj}
\on{I}(\cG,\cP^-)_{\ul{x}}^{\on{spec,loc}}\to (\on{I}(\cG,\cP^-)_{\ul{x}}^{\on{spec,loc}})_{\on{temp}}
\end{equation}
are fully faithful on compacts.

\sssec{}

We first prove the assertion concerning \eqref{e:IGP spec temp proj}. 

\medskip

At a fixed $\ul{x}\in \Ran$, the category $\on{I}(\cG,\cP^-)_{\ul{x}}^{\on{spec,loc}}$ is
$$\IndCoh(\LS^\reg_{\cG,\ul{x}}\underset{\LS^\mer_{\cG,\ul{x}}}\times \LS^\mer_{\cP^-,\ul{x}}\underset{\LS^\mer_{\cM,\ul{x}}}\times \LS^\reg_{\cM,\ul{x}}),$$
where 
$\LS^\reg_{\cG,\ul{x}}\underset{\LS^\mer_{\cG,\ul{x}}}\times \LS^\mer_{\cP^-,\ul{x}}\underset{\LS^\mer_{\cM,\ul{x}}}\times \LS^\reg_{\cM,\ul{x}}$
is an algebraic stack locally almost of finite type.

\medskip

We claim that the localization \eqref{e:IGP spec temp proj} identifies with
$$\Psi:\IndCoh(\LS^\reg_{\cG,\ul{x}}\underset{\LS^\mer_{\cG,\ul{x}}}\times \LS^\mer_{\cP^-,\ul{x}}\underset{\LS^\mer_{\cM,\ul{x}}}\times \LS^\reg_{\cM,\ul{x}})\to
\QCoh(\LS^\reg_{\cG,\ul{x}}\underset{\LS^\mer_{\cG,\ul{x}}}\times \LS^\mer_{\cP^-,\ul{x}}\underset{\LS^\mer_{\cM,\ul{x}}}\times \LS^\reg_{\cM,\ul{x}}).$$

This would imply the fully faithfulness on compacts: indeed $\Coh(-)$ maps fully faithfully to $\QCoh(-)$.

\sssec{}

By \cite[Corollary 7.1.8]{GLC2}, we can can identify $(\on{I}(\cG,\cP^-)_{\ul{x}}^{\on{spec,loc}})_{\on{temp}}$ with
$$\QCoh(\LS^\reg_{\cG,\ul{x}})\underset{\QCoh((\LS^\mer_{\cG,\ul{x}})^\wedge_\reg)}\otimes 
\left(\QCoh(\LS^\reg_{\cG,\ul{x}})\underset{\Sph^{\on{spec}}_\cG}\otimes \on{I}(\cG,\cP^-)_{\ul{x}}^{\on{spec,loc}}\right),$$
which we further identity with
$$\QCoh(\LS^\reg_{\cG,\ul{x}})\underset{\QCoh((\LS^\mer_{\cG,\ul{x}})^\wedge_\reg)}\otimes  \IndCoh((\LS^\Mmf_{\cP^-,\ul{x}})^\wedge_\mf).$$

Since $(\LS^\Mmf_{\cP^-,\ul{x}})^\wedge_\mf$ is the formal completion of a smooth algebraic stack, namely, $\LS^\Mmf_{\cP^-,\ul{x}}$ along $\LS^\reg_{\cP^-,\ul{x}}$, and
$$\LS^\reg_{\cG,\ul{x}}\underset{\LS^\mer_{\cG,\ul{x}}}\times \LS^\Mmf_{\cP^-,\ul{x}} \text{ and } \LS^\reg_{\cP^-,\ul{x}}$$
coincide as \emph{subsets} of $\LS^\Mmf_{\cP^-,\ul{x}}$, we obtain that
$$\QCoh(\LS^\reg_{\cG,\ul{x}})\underset{\QCoh((\LS^\mer_{\cG,\ul{x}})^\wedge_\reg)}\otimes  \IndCoh((\LS^\Mmf_{\cP^-,\ul{x}})^\wedge_\mf)\simeq
\QCoh(\LS^\reg_{\cG,\ul{x}}\underset{\LS^\mer_{\cG,\ul{x}}}\times \LS^\Mmf_{\cP^-,\ul{x}}),$$
as required. 

\sssec{}

We now prove that \eqref{e:IGP temp proj} is fully faithful on compacts. Since we are working at a fixed $\ul{x}\in \Ran$, the category 
$\on{I}(G,P^-)^{\on{loc}}_{\rho_P(\omega_X),\ul{x}}$ is equivalent to 
$$\Dmod_{\frac{1}{2}}(\Gr_{G,\rho(\omega_X)\ul{x}})^{I^-_P,\on{ren}}.$$

\medskip

We claim that the functor
$$\Dmod_{\frac{1}{2}}(\Gr_{G,\rho(\omega_X)\ul{x}})^{I^-_P,\on{ren}}\to \left(\Dmod_{\frac{1}{2}}(\Gr_{G,\rho(\omega_X)\ul{x}})^{I^-_P,\on{ren}}\right)_{\on{temp}}$$
is fully faithful on the eventually coconnective subcategory (with respect to the usual t-structure). 

\medskip

This would imply the desired assertion, since compact objects in
$\Dmod_{\frac{1}{2}}(\Gr_{G,\rho(\omega_X)\ul{x}})^{I^-_P,\on{ren}}$ are eventually coconnective. 

\sssec{}

To prove the claim, it suffices to show that the functor 
$$\on{temp}:\Dmod_{\frac{1}{2}}(\Gr_{G,\rho(\omega_X)\ul{x}})\to \left(\Dmod_{\frac{1}{2}}(\Gr_{G,\rho(\omega_X)\ul{x}})\right)_{\on{temp}}$$
is fully faithful on the eventually coconnective subcategory.

\medskip

To show this, it suffices to prove that the for any $\CF\in \Dmod_{\frac{1}{2}}(\Gr_{G,\rho(\omega_X)\ul{x}})^c$, the object
$$\on{coFib}(\on{temp}^L\circ \on{temp}(\CF)\to \CF)$$ 
is infinitely connective. 

\medskip

We have
$$\on{coFib}(\on{temp}^L\circ \on{temp}(\CF)\to \CF)\simeq \CF\star \on{coFib}(\on{temp}^L\circ \on{temp}(\delta_{1,\Gr_G})\to \delta_{1,\Gr_G}).$$

Since $\CF$ is compact, the functor $\CF\star -$ has a bounded cohomological amplitude.
Hence, it suffices to show that 
$$\on{coFib}(\on{temp}^L\circ \on{temp}(\delta_{1,\Gr_G})\to \delta_{1,\Gr_G})$$ 
is infinitely connective. 

\sssec{}

Since $\on{Sat}_G$ is a t-exact equivalence, it suffices to show that 
$$\on{coFib}(\on{temp}^L\circ \on{temp}(\CF')\to \CF')$$ 
is infinitely connective for any $\CF'\in \Sph^{\on{spec}}_\cG$.

\medskip

However, 
$$\on{temp}^L\circ \on{temp}(\CF')=\Xi\circ \Psi(\CF'),$$
for the stack $\on{Hecke}_{\cG,\ul{x}}^{\on{spec,loc}}$, and the assertion follows.

\qed[\thmref{t:semiinf geom Satake}]

\ssec{Identification of the semi-infinite IC sheaf} \label{ss:semiinf Sat IC} 

The goal of this subsection is to prove \propref{p:semiinf Sat IC}. 

\sssec{}

Consider $\on{I}(G,P^-)^{\on{loc}}_{\rho_P(\omega_X)}$ as acted on by $\Rep(\cG)\otimes \Rep(\cM)$ via
$$\Sat^{\on{nv}}_G:\Rep(\cG)\to \Sph_G \text{ and } \Sat^{\on{nv}}_{M,\tau}:\Rep(\cM)\to \Sph_M.$$

Let $\ol{\cG/\cN^-_P}$ be the affine closure of $\cG/\cN^-_P$ and consider the algebra $\CO_{\ol{\cG/\cN^-_P}}$ of global
functions on $\ol{\cG/\cN^-_P}$ as a commutative algebra equipped with an action of $\cG\times \cM$. Denote the
resulting commutative factorization algebra in $\Rep(\cG)\otimes \Rep(\cM)$ by $R_{\ol{\cG/\cN^-_P}}$. 

\medskip

According to \cite[Sect. 5.5.2]{Ga5} (adapted to a general parabolic\footnote{The case of a general parabolic is considered in \cite{FH}.}), 
we have a canonical action of $R_{\ol{\cG/\cN^-_P}}$ on 
$$\Delta^{-,\semiinf}\in \on{FactAlg}\left(\on{I}(G,P^-)^{\on{loc}}_{\rho_P(\omega_X)}\right).$$

\medskip

Furthermore, the object $\IC^{-,\semiinf}_\Ran\in \on{I}(G,P^-)^{\on{loc}}_{\rho_P(\omega_X)}$ identifies with
\begin{equation} \label{e:Dr-Pluck geom}
R_{\cG}\underset{R_{\ol{\cG/\cN^-_P}}}\otimes \Delta^{-,\semiinf},
\end{equation} 
where:

\begin{itemize}

\item We view $R_{\cG}$ as a commutative factorization algebra in $\on{ComAlg}(\Rep(\cG)\otimes \Rep(\cM))$ via restriction
along the homomorphism $\cM\to \cP^-\to \cG$;

\medskip

\item The homomorphism $R_{\ol{\cG/\cN^-_P}}\to R_{\cG}$ is induced by the map 
$$\cG\to \cG/\cN^-_P\to \ol{\cG/\cN^-_P}.$$

\end{itemize}

\sssec{} \label{sss:R bar acts on Delta spec}

We consider $\on{I}(\cG,\cP^-)^{\on{spec,loc}}$ as equipped with an action of $\Rep(\cG)\otimes \Rep(\cM)$ via
$$\on{nv}:\Rep(\cG)\to \Sph^{\on{spec}}_\cG \text{ and } \on{nv}:\Rep(\cM)\to \Sph^{\on{spec}}_\cM.$$

Consider the object
$$\Delta^{-,\on{spec},\semiinf} \simeq \iota^\IndCoh_*(\CO_{\LS^\reg_{\cP^-}}) \in \on{I}(\cG,\cP^-)^{\on{spec,loc}}.$$

It is equipped with an action of $R_{\ol{\cG/\cN^-_P}}$, obtained from the homomorphism
$$\CO_{\ol{\cG/\cN^-_P}}\to \CO_{\cG/\cN^-_P},$$
where we identify $\CO_{\cG/\cN^-_P}\in \Rep(\cG)\otimes \Rep(\cM)$ with the direct image of
$\CO_{\on{pt}/\cP^-}$ along the projection
$$\on{pt}/\cP^-\to \on{pt}/\cG\times \on{pt}/\cM.$$

\medskip

Note that the object $\IC^{-,\on{spec},\semiinf}$ identifies tautologically with 
\begin{equation} \label{e:Dr-Pluck spec}
R_{\cG}\underset{R_{\ol{\cG/\cN^-_P}}}\otimes \Delta^{-,\on{spec},\semiinf}.
\end{equation} 

\sssec{}

Under the equivalence $\on{Sat}^{-,\semiinf}$, the object $\Delta^{-,\semiinf}\simeq \one_{\on{I}(G,P^-)^{\on{loc}}_{\rho_P(\omega_X)}}$ corresponds to
$$\Delta^{-,\on{spec},\semiinf} \simeq \one_{\on{I}(\cG,\cP^-)^{\on{spec,loc}}}.$$

Comparing \eqref{e:Dr-Pluck geom} and \eqref{e:Dr-Pluck spec}, we obtain that in order to prove \propref{p:semiinf Sat IC}, it suffices to show that
under the above identification, the actions of $R_{\ol{\cG/\cN^-_P}}$ on $\Delta^{-,\semiinf}$ and $\Delta^{-,\on{spec},\semiinf}$
correspond to one another, \emph{up to} an automorphism of $R_{\ol{\cG/\cN^-_P}}\in \Rep(\cG)\otimes \Rep(\cM)$, induced by
the action on $\ol{\cG/\cN^-_P}$ of an element $z\in Z_\cM$. 

\medskip

Indeed, the required isomorphism between \eqref{e:Dr-Pluck geom} and \eqref{e:Dr-Pluck spec} will then be induced by the
action of $z$ on $R_\cG$. 

\sssec{}

We note that $\Delta^{-,\on{spec},\semiinf}$ belongs to the heart of the natural t-structure on $\on{I}(\cG,\cP^-)^{\on{spec,loc}}$,
i.e., the object $\Delta^{-,\on{spec},\semiinf}_{X_\dr}[-1]$
belongs to heart of the structure on $\on{I}(\cG,\cP^-)^{\on{spec,loc}}_{X_\dr}$ (see \cite[Sect. B.8.13]{GLC2}).

\medskip

Hence, the question of two actions of a \emph{classical} commutative factorization algebra 
being equal can be checked at the level of fibers at closed points $\ul{x}\in \Ran$. By factorization, we can assume that 
$\ul{x}$ is a singleton $\{x\}$. 

\medskip

I.e., we need to show that the two actions of $\CO_{\ol{\cG/\cN^-_P}}$ on 
$$\CO_{\on{pt}/\cP^-}\in \QCoh(\on{pt}/\cP^-)\simeq \Rep(\cP^-)$$
coincide. 

\sssec{}

Note that for a (classical) commutative algebra $A$ equipped with an action of $\cG\times \cM$, the datum of an action 
of $A$ on $\CO_{\on{pt}/\cP^-}$ is equivalent to the datum of a $\cG\times \cM$-equivariant map
$$\cG/\cN^-_P\to \Spec(A).$$

\medskip

We take $A=\CO_{\ol{\cG/\cN^-_P}}$. We need to show that the map
\begin{equation} \label{e:map to base affine}
\cG/\cN^-_P\to \ol{\cG/\cN^-_P}
\end{equation} 
corresponding to the action of $R_{\ol{\cG/\cN^-_P},x}$ on $\Delta^{-,\semiinf}_x$ is the tautological embedding 
$$\cG/\cN^-_P\hookrightarrow \ol{\cG/\cN^-_P}.$$

\sssec{}

We claim that it is enough to show that the image of \eqref{e:map to base affine} is contained in 
$\cG/\cN^-_P\subset \ol{\cG/\cN^-_P}$. 

\medskip

Indeed, the set of $\cG\times \cM$-equivariant self-maps
of $\cG/\cN^-_P$ is in bijection with $Z_\cM$. 

\sssec{}

Suppose, by contradiction that the image of \eqref{e:map to base affine} is \emph{not} contained in $\cG/\cN^-_P\subset \ol{\cG/\cN^-_P}$. 
Then it is contained in the boundary
$$\ol{\cG/\cN^-_P}\!-\!\cG/\cN^-_P.$$

However, this means that there exists a highest weight $\lambda\in \Lambda^+_G$ such that the corresponding piece of the 
Drinfeld-Pl\"ucker datum on $\Delta^{-,\semiinf}_x$ (see \cite[Sect. 5.4.2]{Ga5})
\begin{equation} \label{e:Dr-Pl}
V_\cG^\lambda\underset{\cG}\star \Delta^{-,\semiinf}_x\to V^\lambda_\cM\underset{\cM}\star \Delta^{-,\semiinf}_x
\end{equation}
(here $V^\lambda_\cG$ and $V^\lambda_\cM$ are irreducible representations of $\cG$ and $\cM$, respectively)
is $0$.

\medskip

However, it follows from the construction of the Drinfeld-Pl\"ucker datum on $\Delta^{-,\semiinf}_x$ in \cite[Sect.6.3]{Ga4}
that the maps \eqref{e:Dr-Pl} are all non-zero. 

\qed[\propref{p:semiinf Sat IC}]

\sssec{Proof of \propref{p:semiinf spec IC dual}} \label{sss:semiinf spec IC dual}


The proof will proceed along the lines parallel to those of \propref{p:semiinf Sat IC}:

\medskip 

The object
$\Delta^{-,\on{spec},\semiinf}_{\on{co}}$ carries a natural action on $R_{\ol{\cG/\cN^-_P}}$, and we have
$$\IC^{-,\on{spec},\semiinf}_{\on{co}}\simeq R_{\cG}\underset{R_{\ol{\cG/\cN^-_P}}}\otimes \Delta^{-,\on{spec},\semiinf}_{\on{co}}.$$

Hence, in order to prove \propref{p:semiinf spec IC dual}, it suffices to show that with respect to the identification 
$$\text{\eqref{e:self-duality I G P spec}}(\Delta^{-,\on{spec},\semiinf}_{\on{co}})\simeq \Delta^{-,\on{spec},\semiinf}$$
of \secref{sss:Delta spec co to Delta spec}, the above action of $R_{\ol{\cG/\cN^-_P}}$ on $\Delta^{-,\on{spec},\semiinf}_{\on{co}}$ 
corresponds to the action of $R_{\ol{\cG/\cN^-_P}}$ on $\Delta^{-,\on{spec},\semiinf}$ from \secref{sss:R bar acts on Delta spec}.

\medskip

However, this follows from the (essential) uniqueness of such an action proved above.

\qed[\propref{p:semiinf Sat IC}]

\section{Proof of \thmref{t:FG1}} \label{s:proof of Wak}

Throughout this section, unless otherwise stated, we work at the pointwise level.

\subsection{Regular Kac-Moody representations}

\subsubsection{} 

Recall that we have a DG category $\hg\mod_{\crit_G,\reg}$ of 
Kac-Moody representations at critical level constructed in \cite{FG2}.

\medskip 

This DG category has a t-structure with heart $\hg\mod_{\crit_G,\reg}^{\heartsuit}$
the subcategory of $\hg\mod_{\crit_G}^{\heartsuit}$ consisting
of Kac-Moody representations on which the center $\fZ$ of the (completed, 
critically-twisted) enveloping algebra acts through its quotient 
$\fz = \on{Image}(\fZ \to \End(\BV_{\fg,\crit_G}))$ as in 
\cite[Sect. 4.6.1]{GLC2}.
There is a canonical action of $\QCoh(\Op_{\cG}^{\reg})$ on 
$\hg\mod_{\crit_G,\reg}$, where we are using the Feigin-Frenkel isomorphism
$\fz \simeq \CO_{\Op_{\cG}^{\reg}}$.

\medskip

According to \cite{RY}, there is an action of $\fL(G)$ on $\hg\mod_{\crit_G,\reg}$
compatible with the forgetful functor to $\hg\mod_{\crit_G}$ and 
commuting with the action of $\QCoh(\Op_{\cG}^{\reg})$.

\begin{rem}

The existence of this action is subtle and is related to the localization
theorem of \cite{RY} (cf. below). However, in practice, we only will be interested
in the bounded below category $\hg\mod_{\crit_G,\reg}^+$, where one can use
a substitute as in \cite{RY}. 

\medskip

Moreover, we also will only be interested in Iwahori invariants in 
$\hg\mod_{\crit_G,\reg}$, where the localization theorem itself is more classical,
giving another way to circumvent these issues.

\medskip

This is to say: we are referencing \cite{RY} for the existence of this 
action because the reference exists and the resulting formalism is 
convenient, not because we actually need the full strength of these results.

\end{rem}

\subsubsection{The localization theorem}

We now recall:

\begin{thm}[\cite{RY}]\label{t:fg loc}

There is a canonical t-exact equivalence
\[
\Dmod_{\crit}(\Gr_G) \underset{\Rep(\cG)}{\otimes} \QCoh(\Op_{\cG}^{\reg}) 
\simeq \hg\mod_{\crit_G,\reg}
\]

\noindent of $\fL(G)$-categories.

\end{thm}

We will only need the following consequence due to \cite{FG3}.

\begin{cor}\label{c:fg loc iwahori}

There is a canonical t-exact equivalence
\[
\Dmod_{\crit}(\Gr_G)^{I_P} \underset{\Rep(\cG)}{\otimes} \QCoh(\Op_{\cG}^{\reg}) 
\simeq \hg\mod_{\crit_G,\reg}^{I_P}
\]

\noindent of categories acted on by the $(\fL(G),I_P)$-Hecke category.

\end{cor}

\subsection{Transverse Kac-Moody representations}

\subsubsection{}

Below, we will define a certain subcategory
\[
(\hg\mod_{\crit_G,\reg}^{I_P})_{\on{trans}} \subset 
\hg\mod_{\crit_G,\reg}^{I_P}.
\]

\noindent We will denote the embedding of this subcategory by 
$\jmath_*$. The functor $\jmath_*$ will admit a left adjoint $\jmath^*$.
Geometrically, one is meant to picture $\jmath$ as the embedding of the open
Bruhat cell into the flag variety.

\medskip

First, we define a related category
\[
\Big(\Dmod_{\crit}(\Gr_G)^{I_P} \underset{\Rep(\cG)}{\otimes} 
\Rep(\cB)\Big)_{\on{trans}} \subset 
\Dmod_{\crit}(\Gr_G)^{I_P} \underset{\Rep(\cG)}{\otimes} 
\Rep(\cB).
\]

\noindent Here $\Rep(\cG)$ acts on $\Dmod_{\crit}(\Gr_G)$ via 
naive Satake.

\sssec{}

Namely, note that the $*$-averaging functor
\[
\Dmod_{\crit}(\Gr_G)^{\fL(N_P^-)\cdot \fL^+(M)} \to \Dmod_{\crit}(\Gr_G)^{I_P} 
\]

\noindent is an equivalence, similar to the proof of \lemref{l:KM to Sph gen}.
At the pointwise level, we can apply \thmref{t:semiinf geom Satake} to identify
\begin{equation}\label{eq:iwahori-gr ren}
\Dmod_{\crit}(\Gr_G)^{\fL(N_P^-)\cdot \fL^+(M),\on{ren}} \simeq 
\on{I}(\cG,\cP^-)^{\on{spec,loc}}
\simeq \IndCoh(\on{pt}/\cG \underset{\check{\fg}/\cG}{\times} 
\cn^-_P/\cP^-).
\end{equation}

\noindent Here we implicitly trivialize $\omega_X$ in a neighborhood of our
point, allowing us to omit the $\rho_P(\omega_X)$-twist on the first term. 

\medskip

As in \cite[Sect. 12]{AG}, we obtain a similar equivalence
\begin{equation}\label{eq:iwahori-gr}
\Dmod_{\crit}(\Gr_G)^{I_P} \simeq 
\IndCoh_{\Nilp}(\on{pt}/\cG \underset{\check{\fg}/\cG}{\times} 
\cn^-_P/\cP^-)
\end{equation}

\noindent from the corresponding equivalence with 
$\fL(N_P^-)\cdot \fL^+(M)$-invariants instead of $I_P$.

\medskip

We then obtain
\begin{equation}\label{eq:iwahori-gr B}
\Dmod_{\crit}(\Gr_G)^{I_P} \underset{\Rep(\cG)}{\otimes} 
\Rep(\cB) \simeq 
\IndCoh_{\Nilp}(\cn^-_P/\cP^-\underset{\check{\fg}/\cG}{\times}\on{pt}/\cB ).
\end{equation}

\sssec{}

There is a canonical map from 
$\cn^-_P/\cP^-\underset{\check{\fg}/\cG}{\times}\on{pt}/\cB $ to $\cP^-\backslash \cG/\cB$.
We let 
$$(\cn^-_P/\cP^-\underset{\check{\fg}/\cG}{\times}\on{pt}/\cB )_{\on{trans}} = 
(0 \underset{\check{\fg}}{\times} \cn^-_P)/\cB_{\cM}
$$
denote the preimage of the open Bruhat cell (where 
$\cB(M) \subset \cM \subset \cP^-$
is the Borel subgroup of $\cM$).

\medskip

We then define
\[
\Big(\Dmod_{\crit}(\Gr_G)^{I_P} \underset{\Rep(\cG)}{\otimes} 
\Rep(\cB)\Big)_{\on{trans}}
\]

\noindent as the subcategory corresponding to 
\[
\IndCoh_{\Nilp}\big((\cn^-_P/\cP^-\underset{\check{\fg}/\cG}{\times}\on{pt}/\cB )_{\on{trans}}\big)
\]
under \eqref{eq:iwahori-gr B}.

\sssec{}

We clearly have the desired functors $\jmath_*$ and $\jmath^*$ in this case --
we will use the same notation in the setting of 
either $$\Dmod_{\crit}(\Gr_G)^{I_P}\underset{\Rep(\cG)}{\otimes} 
\Rep(\cB) \text{ or }\hg\mod_{\crit_G,\reg}.$$

\subsubsection{}

Next, note that the canonical map $\Op_{\cG}^{\reg} \to \on{pt}/\cG$ 
factors canonically through $\on{pt}/\cB$. 
Therefore, we can write
\begin{multline*}
\hg\mod_{\crit_G,\reg}^{I_P} \simeq 
\Dmod_{\crit}(\Gr_G)^{I_P} \underset{\Rep(\cG)}{\otimes} \QCoh(\Op_{\cG}^{\reg}) \simeq \\
\simeq 
\Dmod_{\crit}(\Gr_G)^{I_P} \underset{\Rep(\cG)}{\otimes} 
\Rep(\cB) 
\underset{\Rep(\cB)}{\otimes}
\QCoh(\Op_{\cG}^{\reg}).
\end{multline*}

We then define $(\hg\mod_{\crit_G,\reg}^{I_P})_{\on{trans}}$ as the
subcategory corresponding to
\[
\Big(\Dmod_{\crit}(\Gr_G)^{I_P} \underset{\Rep(\cG)}{\otimes} 
\Rep(\cB)\Big)_{\on{trans}}
\underset{\Rep(\cB)}{\otimes}
\QCoh(\Op_{\cG}^{\reg})
\]

\noindent under the above equivalence.

\subsection{Wakimoto modules as objects of the transverse subcategory}

We now reduce the proof of \thmref{t:FG1} to a series of intermediate results.

\subsubsection{}

In what follows, we let\footnote{Unlike in \secref{sss:BW}, we do not 
distinguish between the object of $\hg\mod_{\crit_G,\reg}^{I_P,\heartsuit}$
and the underlying object in $\Vect^{\heartsuit}$.}
\[
\BW_{\fg,\crit} := \Wak(\on{Vac}(M)_{\crit_M-\rhoch_P}).
\]

\noindent In a similar spirit, we let
\[
\BV_{\fg,\crit} := \on{Vac}(G)_{\crit_G}.
\]

\noindent We remark that $\BW_{\fg,\crit}$ receives a map of factorization algebras
from $\BV_{\fg,\crit}$; it follows formally that $\BW_{\fg,\crit}$ lies in
$\hg\mod_{\crit_G,\reg}^{I_P,\heartsuit} \subset \hg\mod_{\crit_G}^{I_P,\heartsuit}$.

\subsubsection{}

First, we claim:

\begin{thm}\label{t:wak-trans}

$\BW_{\fg,\crit}$ lies in $(\hg\mod_{\crit_G,\reg}^{I_P})_{\on{trans}}$.

\end{thm}

\sssec{}

It follows from \thmref{t:wak-trans} that the canonical map 
\[
\BV_{\fg,\crit} \to \BW_{\fg,\crit}
\]

\noindent induces a map
\begin{equation}\label{eq:vac-wak-j}
\jmath^*(\BV_{\fg,\crit}) \to \BW_{\fg,\crit} \in (\hg\mod_{\crit_G,\reg}^{I_P})_{\on{trans}}
\end{equation}

\noindent or equivalently
\[
\jmath_*\circ \jmath^*(\BV_{\fg,\crit}) \to \BW_{\fg,\crit} \in 
\hg\mod_{\crit_G,\reg}^{I_P}.
\]

\sssec{}

Next, we claim: 

\begin{thm}\label{t:vac-to-wak}

The map \eqref{eq:vac-wak-j} is an isomorphism.

\end{thm}

\begin{rem}

In the Borel case, \thmref{t:vac-to-wak} amounts to the assertion of
\cite[Theorem 4.11]{FG1}. The proof we give below is similar in spirit, 
but has some changes to more easily adapt to the general parabolic case.

\end{rem}

\subsubsection{}\label{sss:fg1-pf}

Let us now deduce \thmref{t:FG1} from \thmref{t:vac-to-wak}.

\medskip

At the pointwise level, there are canonical equivalences
\begin{multline} \label{e:rewrite Op as Spr}
\IndCoh^*(\Op_\cG^{\on{mon-free}}) \underset{\Sph_\cG^{\on{spec}}}\otimes \on{I}(\cG,\cP^-)^{\on{spec,loc}} \simeq 
\IndCoh^*(\Op_{\cG}^{\mer} \underset{\LS_\cG^{\mer}}{\times} 
(\cn^-_P)^\wedge_0/\cP^-) \simeq \\
\simeq \QCoh(\Op_{\cG}^{\mer} \underset{\LS_\cG^{\mer}}{\times} 
(\cn^-_P)^\wedge_0/\cP^-).
\end{multline}

\noindent Here the last isomorphism comes as in \cite[Proposition 3.8.7]{GLC2}.

\medskip

Under the equivalence
\begin{multline*} 
\KL(G)_{\crit_G}\underset{\Sph_G}\otimes \on{I}(G,P^-)^{\on{loc}}\overset{\text{\eqref{e:tensored up Sat dual}}}\simeq 
\IndCoh^*(\Op_\cG^{\on{mon-free}}) \underset{\Sph_\cG^{\on{spec}}}\otimes \on{I}(\cG,\cP^-)^{\on{spec,loc}}
\overset{\text{\eqref{e:rewrite Op as Spr}}}\simeq  \\
\simeq \QCoh(\Op_{\cG}^{\mer} \underset{\LS_\cG^{\mer}}{\times} 
(\cn^-_P)^\wedge_0/\cP^-)
\end{multline*}
the image of the $\BV_{\fg,\crit}$ under
$$\KL(G)_{\crit_G}\to \KL(G)_{\crit_G}\underset{\Sph_G}\otimes \on{I}(G,P^-)^{\on{loc}}$$
maps to the structure sheaf of
\[
\Op_{\cG}^{\reg} \underset{\LS_\cG^{\reg}}{\times} 
\on{pt}/\cP^-.
\]

It follows by functoriality that $\jmath_*\circ \jmath^*(\BV_{\fg,\crit})$ maps to the
structure sheaf of
\[
\Op_{\cG}^{\reg} \underset{\LS_\cG^{\reg}}{\times} 
\on{pt}/\cP^- 
\underset{\cB\backslash \cG/\cP^-}{\times} \on{pt}/\cB(M) = 
\MOp_{\cG,\cP^-}^{\reg}.
\]
This is the same as 
$\on{co}\!J_\Op^{-,\on{spec}}(\CO_{\Op_{\cM,\rhoch_P}^{\reg}})$ by construction.
Therefore, the existence of an isomorphism as in \thmref{t:FG1} 
follows from \thmref{t:vac-to-wak}.

\medskip 

The fact that this isomorphism satisfies (i) from \thmref{t:FG1} is obvious from the
constructions above. It remains to see that this isomorphism satisfies (ii) from \thmref{t:FG1}.
We do so below.

\subsubsection{}

We consider the following two objects of $\Rep(\cM)$.

\medskip

On the one hand, we consider the tautological $\cM$-torsor $\CP_{\cM,\MOp_{\cG,\cP^-}^{\reg}}$
on  $\MOp_{\cG,\cP^-}$, and we let $\CO_{\CP_{\cM,\MOp_{\cG,\cP^-}^{\reg}}} \in \Rep(\cM)$
be the space of functions on it, viewed as an $\cM$-representation.

\medskip

On the other hand, we consider the object $\CO'_{\CP_{\cM,\MOp_{\cG,\cP^-}^{\reg}}}$
equal to the image of $\Wak^{\on{enh}}(\on{Vac}_{\crit_M-\rhoch_P})$ along the composition
$\text{\eqref{e:add a leg}}\circ \text{\eqref{e:tensored up Sat dual}}$. Note that we can 
rewrite $\CO'_{\CP_{\cM,\MOp_{\cG,\cP^-}^{\reg}}}$ as
$$\DS\circ \Wak(\Sat^{\on{nv}}_M(R_{\cM})\underset{M}\star \on{Vac}_{\crit_M-\rhoch_P}).$$

\medskip

We need to compare two isomorphisms between $\CO_{\CP_{\cM,\MOp_{\cG,\cP^-}^{\reg}}}$ 
and $\CO'_{\CP_{\cM,\MOp_{\cG,\cP^-}^{\reg}}}$. 

\medskip 

The first isomorphism $\eta_1$ comes from the commutativity used
in \thmref{t:FG1} (ii), i.e., it traces back to \eqref{e:ident pre F}.

\medskip

The second isomorphism $\eta_2$ comes using the calculation of the image 
of $\BW_{\fg,\crit}$ under \eqref{e:tensored up Sat dual} given in 
\secref{sss:fg1-pf}.

\medskip

It suffices to show that $\oblv_{\cM}(\eta_1)=\oblv_{\cM}(\eta_2)$ as morphism
in $\Vect$, where is the forgetful functor $\oblv_\cM:\Rep(\cM)\to \Vect$. 

\sssec{}

Unwinding the construction, we obtain that the two induced maps
$$\oblv_{\cM}(\CO_{\CP_{\cM,\MOp_{\cG,\cP^-}^{\reg}}})\rightrightarrows 
\oblv_{\cM}(\CO'_{\CP_{\cM,\MOp_{\cG,\cP^-}^{\reg}}})$$
coincide when precomposed with
$$\CO_{\Op^\reg_\cG}\to \oblv_{\cM}(\CO_{\CP_{\cM,\MOp_{\cG,\cP^-}^{\reg}}}).$$

\medskip

In particular, we obtain that $\oblv_{\cM}(\eta_1)$ and $\oblv_{\cM}(\eta_2)$ are
equal evaluated on the vacuum vector $1\in \CO_{\CP_{\cM,\MOp_{\cG,\cP^-}^{\reg}}}$.

\medskip

In addition, both $\oblv_{\cM}(\eta_1)$ and $\oblv_{\cM}(\eta_2)$ are 
$\oblv_\cM(\CO_{\CP_{\cM,\MOp_{\cG,\cP^-}^{\reg}}})$-linear, where:

\begin{itemize}

\item $\CO_{\CP_{\cM,\MOp_{\cG,\cP^-}^{\reg}}}$ acts on itself naturally;

\item $\CO_{\CP_{\cM,\MOp_{\cG,\cP^-}^{\reg}}}$ acts on $\CO'_{\CP_{\cM,\MOp_{\cG,\cP^-}^{\reg}}}$
via:

\begin{itemize}

\item The identification $\MOp_{\cG,\cP^-}^{\reg}\simeq \Op^\reg_{\cM,\rhoch_P}$ as schemes 
mapping to $\on{pt}/\cM$;

\medskip

\item The action of\footnote{Here $\CP_{\cM,\Op_{\cM,\rhoch_P}^{\reg}}$ denotes the tautological $\cM$-torsor
on $\Op^\reg_{\cM,\rhoch_P}$.}
$\CO_{\CP_{\cM,\Op_{\cM,\rhoch_P}^{\reg}}}$ on $\Sat^{\on{nv}}_M(R_{\cM})\star \on{Vac}_{\crit_M-\rhoch_P}$ 
(as an object of $\KL(M)_{\crit_M-\rhoch_P}$, equipped with an action of $\cM$). 

\end{itemize} 

\end{itemize}

\medskip

The above two properties combined imply that $\oblv_{\cM}(\eta_1)=\oblv_{\cM}(\eta_2)$.

\subsection{Proof of \thmref{t:wak-trans}}

\subsubsection{}\label{sss:trans-m}

Note that $\Sph_M$ acts on $\hg\mod_{\crit_G,\reg}^{\fL(N_{P^-})\cdot \fL^+(M)}$,
and hence on $\hg\mod_{\crit_G,\reg}^{I_P}$.
This induces an action of $\Rep(\cM)$ on $\hg\mod_{\crit_G,\reg}^{I_P}$
via naive Satake (for $M$). 
We will denote this action as
\[
(\CM \in \hg\mod_{\crit_G,\reg}^{I_P}, W \in \Rep(\cM)) \mapsto 
\CM \underset{M}\star W.
\]

\subsubsection{}

For $V \in \Rep(\cG)^{\heartsuit}$, let $\CV_{\Op_{\cG}^{\reg}}$
denote the corresponding vector bundle on $\Op_{\cG}^{\reg}$.

\medskip

By \eqref{eq:iwahori-gr}, there is a canonical action of $\Rep(\cP^-)$
on $\Dmod_{\crit}(\Gr_G)^{I_P}$, and hence on 
$\hg\mod_{\crit_G,\reg}^{I_P}$. The underlying action of 
$\Rep(\cG)$ is given by the action of $\QCoh(\Op_{\cG}^{\reg})$, while
the underlying action of $\Rep(\cM)$ is given by \secref{sss:trans-m}.

\medskip

It follows that there are natural maps
\begin{equation}\label{eq:dp-sat-m}
\CV_{\Op_{\cG}^{\reg}}^{\lambda} \underset{\CO_{\Op_{\cG}^{\reg}}}{\otimes} \CM \to 
\CM \underset{M}\star W^{\lambda}, \quad \CM\in \hg\mod_{\crit_G,\reg}^{I_P}.
\end{equation}

\noindent where $V^{\lambda}$ is the highest weight representation of
$\cG$ and $W^{\lambda} = (V^{\lambda})_{\cN_{P^-}}$ is the highest weight
representation of $\cM$ for a given dominant coweight $\lambda$.
The above maps satisfy Pl\"ucker relations in the natural sense.

\subsubsection{}

Since the $\cG$-bundle on $\Op_{\cG}^{\reg}$ is equipped with a reduction to 
$\cB$, there are canonical line bundles 
$\CL_{\Op_{\cG}^{\reg}}^{\lambda}$ ($\lambda \in \Lambda$) on 
$\Op_{\cG}^{\reg}$ with maps
\begin{equation}\label{eq:dp-op}
\CL_{\Op_{\cG}^{\reg}}^{\lambda} \to \CV_{\Op_{\cG}^{\reg}}^{\lambda}, \quad \lambda\in \Lambda^+_G
\end{equation}

\noindent satisfying Pl\"ucker relations.

\subsubsection{} \label{sss:trans-axioms}

Unwinding the definitions, we obtain that an object $\CM \in \hg\mod_{\crit_G,\reg}^{I_P}$
lies in the transverse subcategory if and only if the following conditions hold:

\medskip

\noindent For every $\lambda\in \Lambda_G^+\cap \Lambda_{G,P}$ (so that $W^\lambda$ is 1-dimensional), 
%
%
%
%
%
%
%
%
the composition
\[
\CL_{\Op_{\cG}^{\reg}}^{\lambda} \underset{\CO_{\Op_{\cG}^{\reg}}}{\otimes} \CM
\to \CV_{\Op_{\cG}^{\reg}}^{\lambda} \underset{\CO_{\Op_{\cG}^{\reg}}}{\otimes} \CM
\to \CM \underset{M}\star W^{\lambda}
\]
is an isomorphism, where:

\medskip

\begin{itemize}

\item The first arrow comes from \eqref{eq:dp-op};

\medskip

\item The second arrow comes from \eqref{eq:dp-sat-m}.

\end{itemize}

\subsubsection{}

Recall now that the action of $\CO_{\Op_{\cG}^{\reg}}$ on
$\BW_{\fg,\crit}$ extends to an action of $\CO_{\MOp_{\cG,\cP^-}^{\reg}}$.




For $W \in \Rep(\check{P}^-)^{\heartsuit}$, let $\CW_{\MOp_{\cG,\cP^-}^{\reg}}$ 
denote the induced vector bundle on $\MOp_{\cG,\cP^-}^{\reg}$.

\medskip

By ``birth of opers" for $\check{M}$ (see \cite[Sect. 5.2]{GLC2}), we have canonical isomorphisms
\begin{multline} \label{e:birth M}
\BW_{\fg,\crit} \underset{M}\star W :=
\Wak(\on{Vac}(M)_{\crit_M-\rhoch_P}) \underset{M}\star W:=
\Wak\Big(\on{Vac}(M)_{\crit_M-\rhoch_P} \underset{M}\star \Sat^{\on{nv}}_M(W)\Big) 
\simeq \\
\simeq 
\Wak\Big(\on{Vac}(M)_{\crit_M-\rhoch_P} 
\underset{\CO_{\Op_{\cM}^{\reg}}}{\otimes} \CW_{\Op_{\cM}^{\reg}}\Big) 
\simeq 
\CW_{\MOp_{\cG,\cP^-}^{\reg}} \underset{{\CO_{\MOp_{\cG,\cP^-}^{\reg}}}}{\otimes}  
\BW_{\fg,\crit}, \quad \quad W\in \Rep(\cM). 
\end{multline}

%

\sssec{}

Note that there are natural maps 
\begin{equation}\label{eq:dp-mop}
\CV_{\Op_{\cG}^{\reg}}^{\lambda}|_{\MOp_{\cG,\cP^-}^{\reg}} \to 
\CW_{\MOp_{\cG,\cP^-}^{\reg}}^{\lambda},
\end{equation}

\noindent where $\CW_{\MOp_{\cG,\cP^-}^{\reg}}^{\lambda}$ 
is the vector bundle corresponding to $W^\lambda\in \Rep(\cM)\to \Rep(\cP^-)$, $\lambda\in \Lambda^+_G$. 

\medskip

Now \cite[Proposition 2.6.1]{FG1} asserts that these isomorphisms are compatible
with the natural Drinfeld-Pl\"ucker structures on both sides, i.e., the diagrams
$$
\CD
\CV_{\Op_{\cG}^{\reg}}^{\lambda} \underset{\CO_{\Op_{\cG}^{\reg}}}{\otimes} \BW_{\fg,\crit}  @>{\text{\eqref{eq:dp-sat-m}}}>> 
\BW_{\fg,\crit}  \underset{M}\star W^{\lambda} \\
@V{\sim}VV @V{\sim}V{\text{\eqref{e:birth M}}}V \\
\CV_{\Op_{\cG}^{\reg}}^{\lambda}|_{\MOp_{\cG,\cP^-}^{\reg}}  
\underset{\CO_{\MOp_{\cG,\cP^-}^{\reg}}}\otimes \BW_{\fg,\crit} @>{\text{\eqref{eq:dp-mop}}}>> 
\CW_{\MOp_{\cG,\cP^-}^{\reg}} \underset{{\CO_{\MOp_{\cG,\cP^-}^{\reg}}}}{\otimes}  
\BW_{\fg,\crit}
\endCD
$$
commute.  

\medskip

This makes it manifest that that $\BW_{\fg,\crit}$ satisfies the condition from \secref{sss:trans-axioms}, since the maps
\eqref{eq:dp-mop} have the required property. 

\begin{rem}

Formally, \cite[Proposition 2.6.1]{FG1} only deals with the Borel case. However,
it is shown in \cite[Sect. 2.8]{FG1} that Proposition 2.6.1 implies its parallel
version for the parabolic being the whole group; similar reasoning shows that
the Borel case implies the general parabolic case. 

\end{rem}

\subsection{Proof of \thmref{t:vac-to-wak}}

\subsubsection{}

We will prove the following two results.

\begin{lem}\label{l:jmath-psi}

The canonical map $\jmath^*(\BV_{\fg,\crit}) \to \BW_{\fg,\crit}$
is an isomorphism after applying the functor 
$$(\hg\mod_{\crit_G,\reg}^{I_P})_{\on{trans}} \overset{\jmath_*}\hookrightarrow 
\hg\mod_{\crit_G,\reg}^{I_P}\to \hg\mod_{\crit_G}\overset{\DS}\to \Vect.$$
\end{lem}

\begin{lem}\label{l:iwahori-ds-cons}
The functor $\DS$ is conservative
when restricted to the eventually coconnective 
transverse subcategory
$$(\hg\mod_{\crit_G,\reg}^{I_P})_{\on{trans}}^{>-\infty}:=
(\hg\mod_{\crit_G,\reg}^{I_P})_{\on{trans}} \cap (\hg\mod_{\crit_G,\reg}^{I_P})^{>-\infty}.$$
\end{lem}

\sssec{}

Clearly the combination of these two lemmas yields \thmref{t:vac-to-wak}.

\subsection{Proof of \lemref{l:jmath-psi}}\label{ss:jmath-psi-pf}

\subsubsection{More on the transverse subcategory}

We now need to look closer at the transverse subcategory and the
functor $\jmath^*$.

\medskip 

Note that the open cell $(\cB\backslash \cG/\cP^-)_{\on{trans}}$ is isomorphic
to $\cB(M)$. Therefore, the functor
\[
\jmath^*:\Dmod_{\crit}(\Gr_G)^{I_P} \underset{\Rep(\cG)}{\otimes} 
\Rep(\cB) \to 
\Big(\Dmod_{\crit}(\Gr_G)^{I_P} \underset{\Rep(\cG)}{\otimes} 
\Rep(\cB)\Big)_{\on{trans}}
\]

\noindent factors as
\[
\begin{tikzcd}
\Dmod_{\crit}(\Gr_G)^{I_P} \underset{\Rep(\cG)}{\otimes} 
\Rep(\cB) \arrow[r,"\pi^*"] &
\Dmod_{\crit}(\Gr_G)^{I_P} \underset{\Rep(\cG)}{\otimes} 
\Rep(\cB(M))\arrow[d,"\widetilde{\jmath}^*"]
\\
&
\Big(\Dmod_{\crit}(\Gr_G)^{I_P} \underset{\Rep(\cG)}{\otimes} 
\Rep(\cB)\Big)_{\on{trans}}.
\end{tikzcd}
\]

\noindent Each of the functors $\pi^*$ and $\widetilde{\jmath}^*$ admit
right adjoints, which we denote as $\pi_*$ and $\widetilde{\jmath}_*$ respectively.

\medskip

We obtain a canonically defined natural transformation
\begin{equation}\label{eq:jmath-pi}
\pi_*\circ \pi^* \to \pi_*\circ \widetilde{\jmath}_*\circ \widetilde{\jmath}^*\circ \pi^* = 
\jmath_*\circ \jmath^*.
\end{equation}

\begin{rem}\label{r:pi-j-quotes}

On the spectral side, 
$\pi^*$ and $\widetilde{\jmath}^*$ correspond to the maps
\[
\begin{tikzcd}
\on{pt}/\cB(M) \overset{\wt\jmath}{\to}
\cP^-\backslash \cG/\cB(M) \overset{\pi}{\to}
\cP^-\backslash \cG/\cB.
\end{tikzcd}
\]
\end{rem}

\subsubsection{}

Note that we have a Cartesian square
\begin{equation}
\label{e:Miura and BM}
\CD
\MOp_{\cG,\cP^-}^{\reg} @>>> \on{pt}/\cB(M) \\
@V{\sfp^{\on{Miu}}}VV @VVV \\
\Op_{\cG}^{\reg} @>>> \on{pt}/\cB.
\endCD
\end{equation} 

From here, we obtain that the endofunctor 
$\pi_*\circ \pi^*$ (on either side of \eqref{e:Miura and BM}) corresponds to 
$$(-)\underset{\CO_{\Op_{\cG}^{\reg}}}{\otimes} 
\CO_{\MOp_{\cG,\cP^-}^{\reg}}.$$

\medskip 

Therefore, the natural transformation \eqref{eq:jmath-pi} corresponds to 
a map
\begin{equation}\label{eq:mop-jmath}
\BV_{\fg,\crit} \underset{\CO_{\Op_{\cG}^{\reg}}}{\otimes}
\CO_{\MOp_{\cG,\cP^-}^{\reg}} \to \jmath_*\circ \jmath^*(\BV_{\fg,\crit})
\end{equation}

\noindent whose composition with
\[
\BV_{\fg,\crit} \to 
\BV_{\fg,\crit} \underset{\CO_{\Op_{\cG}^{\reg}}}{\otimes}
\CO_{\MOp_{\cG,\cP^-}^{\reg}} 
\]

\noindent is the unit map for the adjunction $(\jmath^*,\jmath_*)$.

\begin{lem}\label{l:mop-jmath}

The composition
\[
\BV_{\fg,\crit} \underset{\CO_{\Op_{\cG}^{\reg}}}{\otimes}
\CO_{\MOp_{\cG,\cP^-}^{\reg}} \to \jmath_*\circ \jmath^*(\BV_{\fg,\crit}) \to 
\BW_{\fg,\crit}
\]

\noindent is the canonical map coming from the map 
$\BV_{\fg,\crit} \to \BW_{\fg,\crit}$ and the canonical action of 
$\CO_{\MOp_{\cG,\cP^-}^{\reg}} $ on $\BW_{\fg,\crit}$ extending the
action of $\CO_{\Op_{\cG}^{\reg}}$.

\end{lem}

This follows again from \cite[Proposition 2.6.1]{FG1}.

\subsubsection{}

Recall now that the map 
\[
\BV_{\fg,\crit} \underset{\CO_{\Op_{\cG}^{\reg}}}{\otimes}
\CO_{\MOp_{\cG,\cP^-}^{\reg}} \to \BW_{\fg,\crit}
\]

\noindent yields an isomorphism after applying the Drinfeld-Sokolov functor
$\DS$. 

\medskip

Therefore, by \lemref{l:mop-jmath}, it suffices to show that 
\begin{equation} \label{e:jj DS}
\BV_{\fg,\crit} \underset{\CO_{\Op_{\cG}^{\reg}}}{\otimes}
\CO_{\MOp_{\cG,\cP^-}^{\reg}} \to \jmath_*\circ \jmath^*(\BV_{\fg,\crit})
\end{equation} 

\noindent yields an isomorphism after applying $\DS$.

\medskip

To prove the latter, it suffices to show that the image of \eqref{e:jj DS} in 
$\Whit_*(\hg\mod_{\crit_G,\reg})$ becomes an isomorphism. 

\medskip

Hence, by the definition of the $\jmath_*\circ \jmath^*$ functor on 
$\hg\mod^{I_P}_{\crit_G,\reg}$, it suffices to prove the following:  

\begin{prop}\label{p:jmath-psi-gr}
For any $\CF_0 \in \Dmod_{\crit}(\Gr_G)^{\fL^+(G)}$ with image $$\CF \in 
\Dmod_{\crit}(\Gr_G)^{I_P} \underset{\Rep(\cG)}{\otimes} 
\Rep(\cB),$$ the map
\[
\pi_*\circ \pi^*(\CF) \overset{\eqref{eq:jmath-pi}}{\to} \jmath_*\circ \jmath^*(\CF)
\]

\noindent yields an isomorphism after applying the functor
\begin{equation} \label{e:Ip to Whit B}
\Dmod_{\crit}(\Gr_G)^{I_P} \underset{\Rep(\cG)}{\otimes} \Rep(\cB)\to
\Dmod_{\crit}(\Gr_G)\underset{\Rep(\cG)}{\otimes} \Rep(\cB)\to \Whit_*(G)\underset{\Rep(\cG)}{\otimes} \Rep(\cB).
\end{equation} 

\end{prop}

The rest of this subsection is devoted to the proof of \propref{p:jmath-psi-gr}.

\sssec{Initial reduction step}

Since $\Dmod_{\crit}(\Gr_G)^{\fL^+(G)}$ is compactly generated by 
eventually coconnective objects, we are reduced to the case where
$\CF_0$ is in the heart of the t-structure. By naive Satake,
it corresponds to a $\cG$-representation $V$. 

\medskip

By $\Rep(\cG)$-equivariance
of the functors appearing in this statement, we can assume
$V$ is trivial, i.e., that $\CF_0$ is the unit object $\delta_{1,\Gr_G}$.

\sssec{} 

We now translate the statement to the spectral side. Note that by the construction of the equivalence $\Sat^{-,\semiinf}$, 
we have a commutative diagram of factorization categories\footnote{We are omiting the translation by $\rho_P(\omega_X)$.}
\[
\begin{tikzcd}
\on{I}(G,P^-)^{\on{loc}} \arrow[r,"\Sat^{-,\semiinf}"] 
\arrow[dd]
&
\on{I}(\cG,\cP^-)^{\on{spec,loc}}
\arrow[d,"\on{pre-}\!\sF^{\on{spec}}"]
\\
&
\Rep(\cG) \otimes \Rep(\cM)
\arrow[d,"\id \otimes \on{inv}_{\cM}"]
\\
\Whit_*(G) 
\arrow[r,"\FLE_{\cG,\infty,\tau}"]
&
\Rep(\cG)
\end{tikzcd}
\]

\noindent where:

\begin{itemize}

\item The left vertical arrow is the composition
\begin{multline*}
\on{I}(G,P^-)^{\on{loc}} := 
\Dmod_{\crit}(\Gr_G)^{\fL(N_P^-)\cdot \fL^+(M),\on{ren}}
\to \\
\to \Dmod_{\crit}(\Gr_G)^{\fL(N_P^-)\cdot \fL^+(M)}
\to \Dmod_{\crit}(\Gr_G) \to \Whit_*(G);
\end{multline*}

\item The right vertical arrow above is the $(\IndCoh,*)$-pushforward
to $\LS_{\cG}^{\reg}$. 

\end{itemize}

\medskip

Now at the \emph{pointwise} level, observe that the left vertical arrow above factors
as
\begin{multline} \label{e:Ip to Whit}
\on{I}(G,P^-)^{\on{loc}} \simeq \Dmod_{\crit}(\Gr_G)^{I_P,\on{ren}} \to \\
\to \Dmod_{\crit}(\Gr_G)^{I_P} \to \Dmod_{\crit}(\Gr_G) \to \Whit_*(G) 
\end{multline} 

\noindent with first arrow being an equivalence.

\medskip 

Therefore, we obtain a commutative diagram
\[
\begin{tikzcd}
\Dmod_{\crit}(\Gr_G)^{I_P,\on{ren}} 
\arrow[r,"\text{\eqref{eq:iwahori-gr ren}}"]
\arrow[d]
&
\IndCoh(\on{pt}/\cG \underset{\check{\fg}/\cG}{\times} 
\cn^-_P/\cP^-)
\arrow[d,"(p_1)^\IndCoh_*"]
\\
\Whit_*(G)
\arrow[r,"\FLE_{\cG,\infty,\tau}"]
&
\QCoh(\on{pt}/\cG) = \Rep(\cG),
\end{tikzcd}
\]
where:

\begin{itemize}

\item $p_1$ denotes the projection
$\on{pt}/\cG \underset{\check{\fg}/\cG}{\times} \cn^-_P/\cP^- \to \on{pt}/\cG$;

\item The left vertical arrow is the composition of the last three arrows in \eqref{e:Ip to Whit}.

\end{itemize} 

\medskip

In turn, this yields
\begin{equation}\label{eq:iwahori-vs-whittaker}
\begin{tikzcd}
\Dmod_{\crit}(\Gr_G)^{I_P,\on{ren}} \underset{\Rep(\cG)}{\otimes} 
\Rep(\cB)
\arrow[r]
\arrow[d]
&
\IndCoh(\cn^-_P/\cP^-\underset{\check{\fg}/\cG}{\times}\on{pt}/\cB )
\arrow[d,"(p_2)^\IndCoh_*"]
\\
\Whit_*(G) \underset{\Rep(\cG)}{\otimes} 
\Rep(\cB)
\arrow[r,"\FLE_{\cG,\infty,\tau}\otimes \on{Id}"]
&
\QCoh(\on{pt}/\cB) = \Rep(\cB),
\end{tikzcd}
\end{equation}
where:

\begin{itemize}

\item $p_2$ denotes the projection $\cn^-_P/\cP^-\underset{\check{\fg}/\cG}{\times}\on{pt}/\cB\to \on{pt}/\cB$;

\item The left vertical arrow is the functor \eqref{e:Ip to Whit B}.

\end{itemize} 

\sssec{}

The image of $\delta_{1,\Gr_G}\otimes \one_{\Rep(\cB)}$ along the top horizontal arrow in \eqref{eq:iwahori-vs-whittaker} is 
$\iota^\IndCoh_*(\CO_{\cP^-\backslash \cG/\cB})$ for $\iota$ being the map
\[
 \cP^-\backslash \cG/\cB =
\on{pt}/\cP^- \underset{\on{pt}/\cG}{\times} \on{pt}/\cB \to 
\cn^-_P/\cP^-\underset{\check{\fg}/\cG}{\times}\on{pt}/\cB.
\]

\medskip 

Thus, we have reduced the assertion of the proposition to the following:

\begin{lem}
The morphism 
$$\pi_*\circ \pi^*(\CO_{\cP^-\backslash \cG/\cB})\to \jmath_*\circ \jmath^*(\CO_{\cP^-\backslash \cG/\cB})$$
induces an isomorphism after taking direct image along
$$\cP^-\backslash \cG/\cB\to \on{pt}/\cB.$$
\end{lem}

\begin{proof}

It suffices to show that the map
$$\CO_{\cP^-\backslash \cG/\cB(M)}\to \wt\jmath_*\circ \wt\jmath^*(\CO_{\cP^-\backslash \cG/\cB(M)})$$
induces an isomorphism after taking direct image along
$$\cP^-\backslash \cG/\cB(M)\to \on{pt}/\cB(M).$$

However, this follows from the fact that
$$k\to \Gamma(\cP^-\backslash \cG,\CO_{\cP^-\backslash \cG})$$
is an isomorphism. 

\end{proof}

\qed[\propref{p:jmath-psi-gr}]

\subsection{Proof of \lemref{l:iwahori-ds-cons}}

\subsubsection{}

We follow the same train of ideas as in \secref{ss:jmath-psi-pf}.


\subsubsection{}

We have the following commutative diagram.
\[
\begin{tikzcd}[column sep = 1em]
(\hg\mod_{\crit_G,\reg}^{I_P})_{\on{trans}}
\arrow[r,phantom,"\simeq"]
\arrow[d,"\jmath_*"]
&
\Big(\Dmod_{\crit}(\Gr_G)^{I_P} 
\underset{\Rep(\cG)}{\otimes} \QCoh(\Op_{\cG}^{\reg})\Big)_{\on{trans}}
\arrow[r]
\arrow[d,"\jmath_*"]
&
\Big(\Dmod_{\crit}(\Gr_G)^{I_P} 
\underset{\Rep(\cG)}{\otimes} \Rep(\cB)\Big)_{\on{trans}}
\arrow[d,"\jmath_*"]
\\
\hg\mod_{\crit_G,\reg}^{I_P}
\arrow[r,phantom,"\simeq"]
\arrow[d]
&
\Dmod_{\crit}(\Gr_G)^{I_P}
\underset{\Rep(\cG)}{\otimes} \QCoh(\Op_{\cG}^{\reg})
\arrow[r]
\arrow[d]
&
\Dmod_{\crit}(\Gr_G)^{I_P}
\underset{\Rep(\cG)}{\otimes} \Rep(\cB)
\arrow[d]
\\
\Whit_*(\hg\mod_{\crit_G,\reg})
\arrow[r,phantom,"\simeq"]
&
\Whit_*(G)
\underset{\Rep(\cG)}{\otimes} \QCoh(\Op_{\cG}^{\reg})
\arrow[r]
&
\Whit_*(G)\underset{\Rep(\cG)}{\otimes} \Rep(\cB).
\end{tikzcd}
\]

In the above diagram, the right horizontal arrows are conservative
because $\Op_{\cG}^{\reg} \to \on{pt}/\cB$ is affine. Moreover,
the equivalences and rightward arrows are t-exact. 
Therefore, it is enough to prove that the right vertical composition is conservative 
when restricted to
$$\Big(\Dmod_{\crit}(\Gr_G)^{I_P} 
\underset{\Rep(\cG)}{\otimes} \QCoh(\Op_{\cG}^{\reg})\Big)_{\on{trans}} \cap
\Big(\Dmod_{\crit}(\Gr_G)^{I_P} 
\underset{\Rep(\cG)}{\otimes} \QCoh(\Op_{\cG}^{\reg})\Big)^{>-\infty}.$$

\subsubsection{}

We now pass to the spectral side. We record the following
(well-known) assertion.

\begin{lem}\label{l:abg-amp}

The equivalence \eqref{eq:iwahori-gr} has bounded cohomological 
amplitude for the natural t-structures on both sides.

\end{lem}

\begin{proof}

The heart of the t-structure on $\IndCoh_{\Nilp}(\on{pt}/\cG \underset{\check{\fg}/\cG}{\times} 
\cn^-_P/\cP^-)$ is $\Rep(\cP^-)^{\heartsuit}$. Therefore, it suffices
to see that for each $V \in \Rep(\cP^-)^{\heartsuit}$ maps to an object in
degrees $[-N,N]$ in 
$\Dmod_{\crit}(\Gr_G)^{I_P}$ for some $N$ independent of $V$.

\medskip 

First, when $V$ comes from a representation of $\cM_{\on{ab}}$, the
assertion is standard, cf. \cite[Theorem 5]{AB} for a perversity statement 
at the level
of sheaves on $I_P\backslash\fL(G)/I_P$ and then note that further pushforward
to $I_P\backslash\Gr_G$ has bounded amplitude.

\medskip 

Next, we deduce the claim for any $V$ of the form $V_1 \otimes V_2$
for $V_1 \in \Rep(\cG)^{\heartsuit}$ and $V_2 \in \Rep(\cM_{\on{ab}})^{\heartsuit}$
from the t-exactness of the action of $\Rep(\cG)$ on $\Dmod_{\crit}(\Gr_G)$.

\medskip

Finally, every $V \in \Rep(\cP)^{\heartsuit}$ has a bounded resolution
by terms of the form $V_1 \otimes V_2$ as above -- this follows by a 
standard argument from quasi-affineness of $\cG/[\cP,\cP]$.

\end{proof}

\sssec{}

Applying \lemref{l:abg-amp} and using \eqref{eq:iwahori-vs-whittaker}, we obtain that it suffices to show that 
the functor
\[
\IndCoh_{\Nilp}\big((\cn^-_P/\cP^-\underset{\check{\fg}/\cG}{\times}\on{pt}/\cB)_{\on{trans}}\big) 
\overset{(p_2)^\IndCoh_*}\longrightarrow \IndCoh(\on{pt}/\cB)\simeq \Rep(\cB)
\]
\noindent is conservative on the eventually coconnective subcategory.

\medskip 

Recall that 
\[
(\cn^-_P/\cP^-\underset{\check{\fg}/\cG}{\times}\on{pt}/\cB )_{\on{trans}} = 
(0 \underset{\check{\fg}}{\times} \cn^-_P)/\cB_{\cM}
\]

We factor the above map $p_2$ as
$$(0 \underset{\check{\fg}}{\times} \cn^-_P)/\cB_{\cM}\to \on{pt}/\cB(M)\to \on{pt}/\cB.$$

The $(\IndCoh,*)$-pushforward functor with respect to the first arrow is conservative 
on the eventually coconnective subcategory, since this morphism is affine. 

\medskip

Finally, the pushforward functor along $\on{pt}/\cB(M)\to \on{pt}/\cB$ is conservative, since
the above map is affine (indeed, $\cB/\cB_{\cM} \simeq \cN_P$ is affine).

\qed[\lemref{l:iwahori-ds-cons}]


\newpage

\begin{thebibliography}{99}


\bibitem[AB]{AB} S.~Arkhipov and R.~Bezrukavnikov, {\it Perverse sheaves on affine flags and Langlands dual group}, 
Israel Journal of Mathematics 170 (2009): 135-183.

\bibitem[ABG]{ABG} S.~Arkhipov, R.~Bezrukavnikov and V.~Ginzburg, {\it Quantum groups, the loop Grassmannian, and the Springer resolution},
Journal of the American Mathematical Society 17.3 (2004): 595-678.

\bibitem[AG]{AG} D.~Arinkin and D.~Gaitsgory, {\it Singular support of coherent sheaves, and the Geometric Langlands Conjecture}, 
Selecta Math. N.S. {\bf 21} (2015), 1--199.

\bibitem[AGKRRV]{AGKRRV} D.~Arinkin, D.~Gaitsgory, D.~Kazhdan, S.~Raskin, N.~Rozenblyum and Y.~Varshavsky, \newline
{\em The stack of local systems with restricted variation and geometric Langlands theory with nilpotent singular support},
arXiv:2010.01906. 

\bibitem[Bar]{Bar} J.~Barlev,  {\it D-modules on spaces of rational maps}, Compositio Mathematica {\bf 150} (214), 835--876.

\bibitem[BD1]{BD1} A.~Beilinson and V.~Drinfeld, {\it Quantization of Hitchin's integrable system and Hecke eigensheaves}, 
available at: \url{https://people.mpim-bonn.mpg.de/gaitsgde/grad_2009/}

\bibitem[BD2]{BD2} A.~Beilinson and V.~Drinfeld, {\it Chiral algebras}, 
American Mathematical Soc. {\bf 51} (2004).

\bibitem[BFGM]{BFGM} A.~Braverman, M.~Finkelberg, D.~Gaitsgory and I~Mirković, {\it Intersection cohomology of Drinfeld's compactifications},
Selecta Mathematica 8.3 (2002), 381-418.

\bibitem[BG1]{BG1} A.~Braverman and D.~Gaitsgory, {\it Geometric Eisenstein series}, Invent. Math. {\bf 150} (2002), 287--84.

\bibitem[BG2]{BG2} A.~Braverman and D.~Gaitsgory, {\it Deformations of local systems and Eisenstein series},
GAFA {\bf 17} (2008), 1788--1850.


\bibitem[CR]{CR} J.~Campbell and S.~Raskin, {\it Langlands duality on the Beilinson-Drinfeld Grassmannian}, 
arXiv:2310.19734.


\bibitem[Ch1]{Ch1} L.~Chen, {\it Nearby cycles on Drinfeld-Gaitsgory-Vinberg interpolation Grassmannian and long intertwining functor},
Duke Math. Jour. {\bf 172} (2023), 447--553. 

\bibitem[Ch2]{Ch2} L.~Chen, {\it Nearby cycles in dualities in geometric Langlands program}, PhD Thesis (Harvard University 2021), \newline
available at: \url{https://windshower.github.io/linchen/papers/dissertation.pdf}

\bibitem[CF]{CF} L.~Chen and Y.~Fu, {\it An extension of Kazhdan-Lusztig equivalence}, arXiv:2111.14606.

\bibitem[CFGY]{CFGY} L.~Chen, Y.~Fu, D.~Gaitsgory and D.~Yang, {\it Loop group actions via factorization}, 
forthcoming.

\bibitem[Dr]{Dr} V.~Drinfeld, {\it Two-dimensional $\ell$--adic representations of the fundamental group of a curve over a finite field and automorphic forms on $GL(2)$}, 
Amer. Jour. of Math. {\bf 105} (1983), 85--114. 

\bibitem[DG1]{DG1} V.~Drinfeld and D.~Gaitsgory, {\it On some finiteness questions for algebraic stacks},
Geometric and Functional Analysis 23.1 (2013): 149-294.

\bibitem[DG2]{DG2} V.~Drinfeld and D.~Gaitsgory, {\it Compact generation of the category of $\mathrm {D} $-modules on the stack of $ G $-bundles on a curve},
Cambridge Journal of Mathematics 3.1 (2015): 19-125.

\bibitem[DG3]{DG3} V.~Drinfeld and D.~Gaitsgory, {\it Geometric constant term functor(s)},
Selecta Mathematica 22 (2016): 1881-1951.

\bibitem[FH]{FH} J.~Faergeman, A.~Hayash, {\it Parabolic geometric Eisenstein series and constant term functors}, forthcoming.  

\bibitem[FR]{FR} J.~Faergeman and S.~Raskin, {\it Non-vanishing of geometric Whittaker coefficients for reductive groups}, \newline
arXiv:2207.02955.



\bibitem[Fr]{Fr} E.~Frenkel, {\it Lectures on Wakimoto modules, opers and the center at the critical level},
arXiv:math/0210029.

\bibitem[FG1]{FG1} E.~Frenkel and D.~Gaitsgory, {\it Geometric realizations of Wakimoto modules at the critical level}, 
Duke Math. J. {\bf 143} (2008), 117--203. 

\bibitem[FG2]{FG2} E.~Frenkel and D.~Gaitsgory,  {\it D-modules on the affine flag variety and representations of 
affine Kac-Moody algebras}, Represent. Theory {\bf 13} (2009), 470--608.

\bibitem[FG3]{FG3} E.~Frenkel and D.~Gaitsgory, {\it Localization of D-modules on the affine Grassmannian},
Annals of Mathematics (2009): 1339-1381.

\bibitem[FGV]{FGV} E.~Frenkel, D.~Gaitsgory and K.~Vilonen, {\it On the geometric Langlands conjecture},
Jour. Amer. Math. Soc. {\bf 15} (2002), 367--417. 

\bibitem[Ga1]{Ga1} D.~Gaitsgory, {\it Outline of the proof of the Geometric Langlands Conjecture for $GL(2)$}, 
Ast\'erisque {\bf 370} (2015), 1--112. 

\bibitem[Ga2]{Ga2} D.~Gaitsgory, {\it A conjectural extension of the Kazhdan-Lusztig equivalence}, Publications of RIMS
{\bf 57} (2021), 1227--1376. 

\bibitem[Ga3]{Ga3} D.~Gaitsgory, {\it A "strange" functional equation for Eisenstein series and Verdier
duality on the moduli stack of bundles}, Annales Scientifiques de l'ENS {\bf 50} (2017), 1123--1162. 

\bibitem[Ga4]{Ga4} D.~Gaitsgory, {\it The semi-infinite intersection cohomology sheaf}, Adv. in Math. {\bf 327} (2018), 789--868. 

\bibitem[Ga5]{Ga5} D.~Gaitsgory, {\it The semi-infinite intersection cohomology sheaf-II: the Ran space
version}, in: Representation theory and algebraic geometry, Trends in Mathematics, Birkh\"auser (2022). 

\bibitem[Ga6]{Ga6} D.~Gaitsgory, {\it The local and global versions of the Whittaker category}, 
PAMQ {\bf 16}, (2020), 775--904. 

\bibitem[GLC1]{GLC1} D.~Gaitsgory and S.~Raskin, {\it Proof of the geometric Langlands conjecture I: construction of the functor}, 
arXiv:2405.03599

\bibitem[GLC2]{GLC2} D.~Arinkin, D~Beraldo, L.~Chen, J.~Faergeman, D.~Gaitsgory, K.~Lin, S.~Raskin and N.~Rozenblyum, \newline
{\it Proof of the geometric Langlands conjecture II: Kac-Moody localization and the FLE},
arXiv:2405.03648

\bibitem[GLC4]{GLC4} D.~Arinkin, D~Beraldo, L.~Chen, J.~Faergeman, D.~Gaitsgory, K.~Lin, S.~Raskin and N.~Rozenblyum, 
{\it Proof of the geometric Langlands conjecture IV: Ambidexterity}, 
available at: \url{https://people.mpim-bonn.mpg.de/gaitsgde/GLC/}

\bibitem[GLC5]{GLC5} D.~Gaitsgory and S.~Raskin, {\it Proof of the geometric Langlands conjecture V: The multiplicity one theorem},
available at: \url{https://people.mpim-bonn.mpg.de/gaitsgde/GLC/}

\bibitem[GLys]{GLys} D.~Gaitsgory and S.~Lysenko, {\it Parameters and duality for the metaplectic geometric Langlands theory}, 
Selecta Math. New Ser. {\bf 24} (2018), 227--301. Also arXiv: 1608.00284.

\bibitem[GaRo3]{GaRo3} D.~Gaitsgory and N.~Rozenblyum, {\it A study in derived algebraic geometry, Vol. 1: Correspondences and Duality}, 
Mathematical surveys and monographs {\bf 221} (2017), AMS, Providence, RI.

\bibitem[GaRo4]{GaRo4} D.~Gaitsgory and N.~Rozenblyum, {\it A study in derived algebraic geometry, Vol. 2: Deformations, Lie Theory and
Formal Geometry}, 
Mathematical surveys and monographs {\bf 221} (2017), AMS, Providence, RI.

\bibitem[Laum]{Laum} G.~Laumon, {\it Correspondance de Langlands géométrique pour les corps de fonctions},
Duke Math. J. 54.1 (1987): 309-359.

\bibitem[Lin]{Lin} K.~Lin, {\it Poincar\'e series and miraculous duality}, 
arXiv:2211.05282.

\bibitem[Ra1]{Ra1} S.~Raskin, {\it A geometric proof of the Feigin-Frenkel theorem}, Representation Theory {\bf 16} (2012), 489--512. 

\bibitem[Ra2]{Ra2} S.~Raskin, {\it W-algebras and Whittaker categories}, 
Selecta Mathematica {\bf 27} (2021).

\bibitem[Ra3]{Ra3} S.~Raskin, {\it Chiral principal series categories I: Finite dimensional calculations}, 
Adv. in Math. {\bf 388} (2021). 

\bibitem[Ra4]{Ra4} S.~Raskin, {\it Chiral principal series categories II: the factorizable Whittaker category},
available at: \url{https://gauss.math.yale.edu/~sr2532/cpsii.pdf}

\bibitem[RY]{RY} S.~Raskin and D.~Yang, {\it Affine Beilinson-Bernstein localization at the critical level},
Annals of Mathematics 200.2 (2024): 487-527.

\bibitem[Ro]{Ro} N.~Rozenblyum, {\it Connections on moduli spaces and infinitesimal Hecke modifications}, 
arXiv:2108.07745. 

\bibitem[Ya-D]{YaD} D.~Yang, {\it Categorical Moy-Prasad theory}, 
arXiv:2104.12917.

\bibitem[Ya-R]{YaR} R.~Yang, {\it Twisted Whittaker category on affine flags and the category of representations of the mixed quantum group},
Compositio Mathematica {\bf 160} (2024), 1349--1417. 



\end{thebibliography}
\end{document}